\documentclass[12pt]{report}

\usepackage{epsfig}
\usepackage{amsfonts}
\usepackage{amssymb}
\usepackage{amsmath}
\usepackage{amsthm}
\usepackage{latexsym}

\newtheorem{theorem}{Theorem}[section]
\newtheorem{lemma}[theorem]{Lemma}
\newtheorem{corollary}[theorem]{Corollary}
\newtheorem{proposition}[theorem]{Proposition}
\newtheorem{definition}[theorem]{Definition}

\newtheorem{conjecture}[theorem]{Conjecture}

\newtheorem{remark}[theorem]{Remark} 
\newtheorem{example}[theorem]{Example} 

\newcommand{\R}[1]{\mathbb{R}^{#1}}

\newcommand{\lbar}{\underline{L}}

\newcommand{\pd}{\partial}

\newcommand{\bp}{\indent$\bullet$\hspace{.1in}}
\newcommand{\sla}[1]{\displaystyle{\not}{#1}}
\newcommand{\ul}[1]{\underline{#1}}

\newcommand{\delm}{\delta_-}
\newcommand{\delp}{\delta_+}
\newcommand{\delh}{\delta_H}
\newcommand{\tg}{\hat{\Gamma}}
\newcommand{\tsla}[1]{\tilde{\sla{#1}}}
\newcommand{\fa}{{\mathfrak{a}}}
\newcommand{\fb}{{\mathfrak{b}}}
\newcommand{\fc}{{\mathfrak{c}}}
\newcommand{\fd}{{\mathfrak{d}}}

\begin{document}

\title{Nonlinear Wave Dynamics in Black Hole Spacetimes}
\author{John Stogin}
\date{May, 2017}
\maketitle

\abstract{

This thesis details a method for proving global boundedness and decay results for nonlinear wave equations on black hole spacetimes. The method is applied to five example problems of increasing difficulty. The first problem, which addresses the semilinear wave equation on Minkowski space, is quite simple and should be accessible to a reader who is still new to the field of partial differential equations. The final problem, which was posed by Ionescu and Klainerman in \cite{ionescu2014global}, constitutes a step toward proving stability for slowly rotating Kerr black holes. The remaining intermediate problems are: a semilinear wave equation on the Schwarzschild spacetime, a semilinear wave equation on any subextremal Kerr spacetime with the additional assumption of axisymmetry, and a restriction of the final problem to the Schwarzschild case.

The method used in this thesis is based on a few particular developments that may be useful for other related problems. These include: a new method for constructing Morawetz-type estimates that is fairly robust (insofar as it may be successfully applied to all five problems), a strategy based on a decay hierarchy for energy estimates on uniformly spacelike hypersurfaces using, in particular, a notion of weak decay, and a technique for handling certain terms with factors that are singular on an axis of symmetry.

}

\tableofcontents

\chapter{Introduction}

A standard problem in the field of nonlinear hyperbolic partial differential equations (PDEs) is to understand the long-time behavior of solutions to a nonlinear equation which initially resemble a known special solution to the same equation. Due to the wave-like nature of hyperbolic PDEs, these solutions (more precisely, the deviation of these solutions from the known special solution) may be referred to as \textit{nonlinear waves}. Much can be learned of these solutions by studying solutions to a corresponding linearized equation, which is much easier to analyze. Solutions to the linearized equation are referred to as \textit{linear waves}. For semilinear PDEs, the difference between the nonlinear equation and the linear one is comprised of terms which are at least quadratic in first derivatives of the waves. In this case, if the waves are sufficiently small, then the quadratic terms are negligible and the two equations are almost the same. Therefore, understanding nonlinear waves can be reduced to the following two sub-goals.
\begin{itemize}
\item Understand properties of the linear waves, especially how they become small over time as they spread out, and
\item Using the properties of the linear waves, show that the difference between the linear and nonlinear equations is negligable so that the properties also belong to the nonlinear waves.
\end{itemize}
This type of argument is often referred to as a bootstrap argument. If the nonlinear waves become smaller over time, the known special solution is said to be stable.

One of the most impressive nonlinear stability results is the seminal work of Christodoulou and Klainerman \cite{christodoulou1993global} which proves stability of the completely flat spacetime known as Minkowski space as a solution to the very elegant, yet incredibly complicated, Einstein vacuum equations (EVEs). In their work, Christodoulou and Klainerman exploit a specific structure of these equations, called \textit{null structure}, which explains the mechanism by which nonlinear waves decay as they expand over time. The EVEs are not semilinear, but one can approximate them by semilinear equations as a first step, and the null structure is still seen in this approximation.

However, other known solutions to the EVEs exist, including the Schwarzschild and (more general) Kerr black holes. The full nonlinear stability of these black holes is an important open problem.

$$*$$

This thesis presents a method for proving global boundedness and decay for solutions to nonlinear wave equations on $3+1$ dimensional vacuum spacetimes. The method, which is explained further in \S\ref{i:general_method_sec}, is sufficiently robust to address a handful of related problems. In particular, five problems are discussed in this thesis. They are
\begin{enumerate}
\item Semilinear waves with null structure on the Minkowski background \\(Chapter \ref{mink_chap}).
\item Semilinear waves with null structure on the Schwarzschild background \\(Chapter \ref{szd_chap}).
\item Axisymmetric semilinear waves with weak null structure on any subextremal Kerr background (Chapter \ref{kerr_chap}).
\item Nonlinear axisymmetric perturbations of a nontrivial solution to a wave map problem on the Schwarzschild background satisfying a condition on the axis corresponding to conservation of angular momentum (Chapter \ref{wm_szd_chap}).
\item Nonlinear axisymmetric perturbations of a nontrivial solution to a wave map problem on any sufficiently slowly rotating Kerr background satisfying a condition on the axis corresponding to conservation of angular momentum \\(Chapter \ref{wm_kerr_chap}).
\end{enumerate}

Each successive problem in the above list comes with additional difficulties, both artificial and inherent to the problem. The intention for solving each problem is to save the reader from considering too many complications at once. Instead, the reader may compare any given problem that introduces a new complication to the previous problem in order to understand how the new complication is handled. For more information about each of these problems, see \S\ref{i:problems_sec}.

\section{Foundational, prior, and related works}
Modern analysis of wave dynamics on asymptotically flat spacetimes began in the 1980s with the introduction of the powerful \textit{vectorfield method}. This method is comprised of two parts. First, \textit{vectorfield multipliers} are used to derive energy-type estimates for a solution $\psi$ to a wave equation. Then \textit{vectorfield commutators} are used to find additional solutions related to $\psi$, which yield more information about $\psi$ upon application of the same energy-type estimates. The vectorfield commutators are generally most useful when they are related to spacetime symmetries.

In Minkowski space, which has a high number of symmetries, a lot of information is available about decay. To precisely quantify decay rates, we first introduce the null pair ($L$, $\lbar$).
\begin{align*}
L &= \pd_t +\pd_r \\
\lbar &= \pd_t-\pd_r,
\end{align*}
where the partial derivatives $\pd_t$ and $\pd_r$ are with respect to the standard spherical coordinate system for Minkowski space. We also denote by $\sla\nabla$ the restriction of the gradient operator to each sphere of constant radius $r$ and time $t$. Together, $L\psi$, $\lbar\psi$, and $\sla\nabla\psi$ tell us all information about the first derivative of a wavefunction $\psi$ at any point not belonging to the center $\{r=0\}$.

Using the vectorfield method in Minkowski space, one can prove the following decay rates for a solution $\psi$ to the linear wave equation in $n+1$ dimensions.
\begin{align}
|\psi(t,r)| &\lesssim (1+|t-r|)^{-1/2}(1+t+r)^{-(n-1)/2} \label{psi_decay_rate} \\
|L\psi(t,r)| &\lesssim (1+|t-r|)^{-1/2}(1+t+r)^{-(n+1)/2} \label{Lpsi_decay_rate} \\
|\sla\nabla\psi(t,r)| &\lesssim (1+|t-r|)^{-1/2}(1+t+r)^{-(n+1)/2} \\
|\lbar\psi(t,r)| &\lesssim (1+|t-r|)^{-3/2}(1+t+r)^{-(n-1)/2}. \label{Lbarpsi_decay_rate}
\end{align}
This was shown by Klainerman \cite{klainerman1985uniform} \cite{klainerman1987remarks} using energy conservation together with the \textit{global Sobolev inequality} (also known as the \textit{Klainerman--Sobolev inequality}), which exploits the fact that some vectorfield commutators have $t$ weights and also the fact that $L$, $\lbar$, and $\sla\nabla$ can be expressed in terms of these commutators with additional $t$ or $r$ weights. In particular, from these estimates, one sees that the waves are largest in the region where $r\approx t$, and that in this region $\lbar\psi$ decays slower than $L\psi$ and $\sla\nabla\psi$.

In the particularly relevant case $n=3$, the best decay rate for $\lbar\psi$ is $(1+t)^{-1}$, which is not integrable in time. This is the reason for the blowup results in $3+1$ dimensions found by John \cite{john1981blow}. However, if the nonlinear term has at least one so-called ``good derivative''
$$\Box_g\psi=\bar{D}\psi D\psi$$
$$\bar{D}\psi = L\psi\text{ or }\sla\nabla\psi,$$
then it is possible to prove global boundedness and decay. This condition is called the \textit{null condition} and it happens to appear in many equations with physical relevance. \cite{klainerman1986null}

Unfortunately, the use of commutators with $t$ weights complicates the nonlinear problem and doesn't generalize easily to black hole spacetimes. A newer method (often referred to as \textit{the new method} at the time of this writing) due to Dafermos and Rodnianski \cite{dafermos2010new} obtains decay by defining a family of weighted energies $E_p(t)$ on an asymptotically null foliation. By taking later and later slices in this foliation, one leaves behind outgoing null rays along which $\lbar\psi$ may decay poorly. This family of estimates yields a decay hierarchy which is sufficient for the same class of nonlinear problems. For example, see Yang's work. \cite{yang2013global}

More recently, Ionescu and Klainerman \cite{ionescu2014global} found a similar family of estimates on uniformly spacelike hypersurfaces. This was accomplished using an asymptotically null timelike vectorfield multiplier instead of the uniformly timelike one used for asymptotically null foliations by Dafermos and Rodnianski. This has the advantage that it is more suitably adapted to Cauchy problems defined on uniformly spacelike initial hypersurfaces. The work presented in this thesis is the first application of these estimates to a nonlinear problem.

The study of nonlinear waves on black hole spacetimes is an active area of research and there are many recent results. Luk \cite{luk2010null} adapted a conformal killing vectorfield in Minkowski space (which has $t$ weights) to prove a global result for slowly rotating Kerr spacetimes. Dafermos, Rodnianski, and Shlapentokh-Rothman \cite{DRSR} used a Fourier decomposition to prove estimates for the linear scalar wave equation on all subextremal Kerr spacetimes. (In contrast, this thesis only addresses \textit{axisymmetric} waves on all subextremal Kerr spacetimes, albeit using only physical-space methods.) Andersson and Blue \cite{AnBlu}, found a purely physical-space approach to prove a similar result in slowly rotating Kerr spacetimes using the hidden Carter operator. \cite{carter}

There are also results for more complicated linear systems that more closely resemble the EVE system. One of these is the result of Ionescu and Klainerman \cite{ionescu2014global}, which laid the foundation for this thesis. More recently, Dafermos, Holzegel, and Rodnianski \cite{dafermos2016linear} proved stability for the complete linearized EVE system about Schwarzschild, and Hintz and Vasy \cite{hintz2016global} proved nonlinear stability of Kerr-de Sitter black holes.

\section{The general method of this thesis}\label{i:general_method_sec}

\subsection{The main idea}

We take a moment to discuss the main idea behind the general method of this thesis. \textbf{The statements made here are only approximately true, as they have been simplified to communicate the main idea.}

The $r^p$ weighted energy norm
$$E_p(t)=\int_{\Sigma_t}r^p\left[(\bar{D}\psi)^2+r^{-2}(D\psi)^2\right],$$
has been adapted to have better control of good derivatives $\bar{D}\psi$, which are expected to decay faster. (See the decay rates \ref{psi_decay_rate}-\ref{Lbarpsi_decay_rate}.) After a considerable amount of work, one can almost\footnote{More accurately, the norm $E_{p-1}(t)$ is replaced with a norm $B_p(t)$, which has the same $r^{p-1}$ weight, but suffers a loss due to trapping in black hole spacetimes.} obtain estimates of the following form in almost\footnote{We will only take $p\in [\delm,2-\delp]$ so that the constant of the estimate is uniform in $p$.} the range $p\in [0,2]$.
$$E_p(t_2)+\int_{t_1}^{t_2}E_{p-1}(t)dt\lesssim E_p(t_1).$$
This roughly suggests that $E_2(t)$ is bounded and that $E_1(t)$ decays like $t^{-1}$ and furthermore that $E_0(t)$ decays like $t^{-2}$ and $E_{-1}(t)$ decays like $t^{-3}$. By interpolation, this suggests that 
$$E_p(t)\sim t^{p-2}$$
for all $p\in[-1,2]$.

Let $\psi^s$ represent a wavefunction obtained by applying $s$ commutator operators to $\psi$ and let $E_p^s(t)$ be its weighted energy, which satisfies the same type of estimate and decay rates. There is an $L^\infty$ estimate roughly of the form\footnote{The increase in number of derivatives from $s$ to $s+3$ actually depends on the particular problem.}
$$|r\bar{D}\psi^s|^2+|D\psi^s|^2+|\psi^s|^2 \lesssim E_0^{s+3}(t),$$
which implies that
$$|r\bar{D}\psi^s|+|D\psi^s|+|\psi^s| \lesssim t^{-1}.$$
This estimate should be compared to the decay rates \ref{psi_decay_rate}-\ref{Lbarpsi_decay_rate}. For example, this estimate implies that $|rL\psi|\lesssim t^{-1}$. The decay rate \ref{Lpsi_decay_rate} implies that $|rL\psi|\lesssim r(1+|t-r|)^{-1/2}(1+t+r)^{-2}$. Note that $\sup_r r(1+|t-r|)^{-1/2}(1+t+r)^{-2}\approx t^{-1}$. \textbf{The decay rates \ref{psi_decay_rate}-\ref{Lbarpsi_decay_rate} certainly provide more information, but the worst case (where $r\approx t$) is captured by the above estimate.}

In fact more generally, there is an $L^\infty$ estimate of the form
$$|r^{p+1}\bar{D}\psi^s|^2+|r^pD\psi^s|^2+|r^p\psi^s|^2\lesssim E_{2p}^{s+3}(t).$$
This implies that
$$r^p(|r\bar{D}\psi^s|+|D\psi^s|+|\psi^s|)\lesssim (E_{2p}^{s+3}(t))^{1/2}\lesssim t^{p-1}.$$
Put another way, multiplying the left side by $r^p$ corresponds to multiplying the right side by $t^p$.
\textbf{This is consistent with the fact that the quantities are largest in the region $r\approx t$, where $r^p\approx t^p$.} Again, a comparison to the decay rates \ref{psi_decay_rate}-\ref{Lbarpsi_decay_rate} can be made.

From the above weighted estimate, the best achievable decay occurs when $2p=-1$, in which case the above quantities decay like $t^{-3/2}$, which is integrable in time. This integrability is essential for the nonlinear problem. In reality, for low values of $p$, one does not actually have decay for $E_p^{s}(t)$, but instead one has decay in a weaker integrated sense
$$\int_{t}^\infty E_p^s(\tau)d\tau \lesssim t^{(p-2)+1}.$$
But this weaker sense is still sufficiently strong to solve the problems of this thesis.

\subsection{The steps involved in the general method}\label{i:steps_sec}

Although there will be variations between each of the five problems, there is a general seven step procedure that will be followed for each problem. The steps are outlined below.

\textbf{Step 1:} Determine an energy quantity that is conserved by solutions to the linearized problem. This is possible using the killing vectorfield $\pd_t$ as a vectorfield multiplier, since all of the spacetimes under consideration are stationary. (However, the way to use this multiplier for the final two problems is not trivial and was an import observation in \cite{ionescu2014global}.) Based on this conserved quantity, prove an energy estimate for the nonlinear problem. Furthermore, prove a slightly more general estimate based on the vectorfield multiplier $h\pd_t$, where $h$ is a function of $r$ only. The $h\pd_t$ estimate will be used in step 3.

\textbf{Step 2:} Prove a Morawetz estimate. This estimate is quite delicate to construct as it is very sensitive to the geometry of the spacetime. A substantial first step towards proving the Morawetz estimate is to prove a \textit{partial Morawetz estimate}. The strategy used in this thesis to construct the partial Morawetz estimate is new. It was discovered while constructing a Morawetz estimate for axisymmetric waves in subextremal Kerr spacetimes starting with an approach adapted from \cite{AnBlu}. Fortunately, it generalizes to the other problems. From the partial Morawetz estimate and a few additional tweaks, the Morawetz estimate follows.

\textbf{Step 3:} Prove the $r^p$ estimate. This estimate is a combination of the $h\pd_t$ estimate from step 1, the Morawetz estimate from step 2, and a third estimate, called the \textit{incomplete $r^p$ estimate near $i^0$}, which is responsible for contributing the $r^p$ weights to the norms used in the $r^p$ estimate. As its name implies, the incomplete $r^p$ estimate near $i^0$ is related only to the asymptotic flatness of the spacetime. Its proof uses a vectorfield multiplier that is only supported near $i^0$, where deviations from Minkowski space can be treated as small errors.

\textbf{Step 4:} Prove the \textit{dynamic estimates}, which encode all of the necessary information about future dynamics of the problem. The dynamic estimates are simply restatements of the energy estimate from step 1 and the $r^p$ estimate from step 3, however additional work is required to show that they hold for higher order wavefunctions obtained by applying commutator operators to the wavefunction $\psi$. Unfortunately, this work can be quite cumbersome for the more complicated problems. It is worth noting that if the reader is not concerned with the construction of the estimates in steps 1-3, then it is possible to begin at step 4.

\textbf{Step 5:} Prove $L^\infty$ estimates compatible with the energy norms in the dynamic estimates. The $L^\infty$ estimates are based on a simple Sobolev estimate, which is applied to various derivatives of $\psi$ with $r$ weights. This step is strongly tied with step 4, because the commutators used in step 4 must be sufficiently useful to make the $L^\infty$ estimates possible. Generally speaking, the commutators will be the operator $\pd_t$ and the rotation operators $\Omega$. (In the Kerr problems, the rotation operators are replaced by the second order Carter operator $Q$.)

\textbf{Step 6:} Identify structural information about the nonlinear term. Most importantly, observe that some form of the null condition is satisified. For the first three problems, this is done rather quickly at the beginning of the final section in the chapter, since the structure is given by choice. For the final two problems, this is a more serious step, because the nonlinear terms are given by the problem itself, and some effort must be spent ensuring that these terms indeed satisfy a form of the null condition. The final problem will also require a new structural condition on the axis.

\textbf{Step 7:} State and prove the main theorem. The statement of the main theorem has three parts. The first part is a statement about how the weighted energy norms from the dynamic estimates behave over time (whether they are bounded or decay, and if they decay, at what rate). The second part is a restatement of the $L^\infty$ estimates. The third part is a direct consequence of the first two parts--it states how derivatives of $\psi$ with various $r$ weights decay in time. The proof of the main theorem uses a standard bootstrap argument.

\section{A description of each problem in this thesis}\label{i:problems_sec}

We take a moment to discuss each of the five problems in more detail and how they relate to each other.

\subsection{Problem 1: Waves in the Minkowski spacetime}

The first problem solves the equation
$$\Box_g\psi=D\psi\bar{D}\psi,$$
where $g$ is the Minkowski metric, $D$ represents any normalized derivative and $\bar{D}$ represents any of the good derivatives $L$ or $\sla\nabla$.

This problem is the simplest of all five problems and is included to present the general method of this thesis as clearly as possible.

Unlike in other problems, two versions of the Morawetz estimate are proved for this problem. The first is a more familiar type of Morawetz estimate with stronger weights and the second is the estimate that must be used by the method of this thesis. The purpose of proving both is to provide comparison to the existing literature and to draw attention to the fact that this method requires a particular type of Morawetz estimate.

Also, some concepts, such as the decay mechanism for the linear problem (Lemma \ref{m:decay_mechanism_lem}), the interpolation between weighted energies (Lemma \ref{m:interpolation_lem}), and the weak decay principle (Lemma \ref{weak_decay_lem}) are only presented for this problem.

There is one slight complication for this problem that is not present in other problems due to the Minkowski space having a center. (For all other problems, only the black hole exterior is considered.) For this reason, the spatial translation commutators $\pd_{x^i}$ are used in addition to the time translation operator $\pd_t$ and the rotation operators $\Omega_{x^i}$ in the solution of this problem.

\subsection{Problem 2: Waves in the Schwarzschild spacetime}

The second problem solves the equation
$$\Box_g\psi=\mathcal{N},$$
where $g$ is the Schwarzschild metric, and $\mathcal{N}$ is a finite sum of quadratic terms. Each particular term in $\mathcal{N}$ is required to have at least one good factor (essentially $\bar{D}\psi$) and an additional requirement is made ensuring that terms with $\pd_r\psi$ vanish at the event horizon. The purpose of this additional requirement is to avoid commuting with the redshift vectorfield, which significantly complicates the argument. (This assumption is removed in the next problem.) For a precise definition of $\mathcal{N}$, see Definition \ref{s:nonlinear_def}.

This is the simplest of all problems involving a black hole.

The commutators used for this problem are the time translation operator $\pd_t$ and the rotation operators $\Omega$.

\subsection{Problem 3: Axisymmetric waves in subextremal Kerr spacetimes}

The third problem investigates axisymmetric solutions
$$\pd_\phi\psi = 0$$
to the equation
$$\Box_g\psi=\mathcal{N},$$
where $g$ is any subextremal Kerr metric and $\mathcal{N}$ is a finite sum of terms. Each particular term is at least quadratic and must have at least one good factor (ie. $\bar{D}\psi$). There is a weak condition on the remaining factors of higher order (at least cubic) terms ensuring that the factors remain bounded during evolution in time. For a precise definition of $\mathcal{N}$, see Definition \ref{k:nonlinear_def}.

The assumption of axisymmetry serves two purposes. First, it ensures that the energy quantity that is conserved by the linear problem is comprised of purely nonnegative terms. (Thus, the issue of superradiance is avoided in this thesis.) Second, it ensures only a single trapping radius, a condition required for the construction of a Morawetz estimate without significant complications. The assumption of axisymmetry could be relaxed at least in the slowly rotating range by the use of higher order operators \cite{AnBlu} or in the full subextremal range using a Fourier decomposition \cite{DRSR}. However, by using only physical-space methods, this problem serves as a primer for the fifth problem.

The rotation commutators $\Omega$ used in the Schwarzschild problem are replaced with the second order Carter operator $Q$, which coincides with $\Omega_x^2+\Omega_y^2+\Omega_z^2$ in Schwarzschild. The $L^\infty$ estimates depend heavily on the fact that $Q$ is an elliptic operator. The operator notation $\Omega$ is redefined to represent the operator $r\sla\nabla$. Since $Q$ is second order, when it acts on a product $\psi_1\psi_2$, there will be some terms that cannot be expressed as a product of commutators applied to either $\psi_1$ or $\psi_2$. The simplest example is the term $\Omega\psi_1\Omega\psi_2$, which shows up in the expansion of $Q(\psi_1\psi_2)$. This leads to the generalized notation $\Omega^l\psi^{s-l}$ to refer to a factor of order $s$ that can be obtained by applying commutators to a product.

Recall that the previous problem made a simplifying assumption that nonlinear terms with $\pd_r\psi$ vanish on the event horizon. Since this is no longer the case, a new commutator $\tg=1_H\pd_r$ (with $1_H$ a function supported near the event horizon) is used to handle these terms. This commutator is often called the \textit{redshift commutator}.  Since $\tg$ does not commute with the wave operator, it is given a different index $k$, and the proof of the nonlinear problem will be a finite inductive argument in $k$. This leads to the generalized notation $\Omega^l\psi^{s-l,k}$ that appears on the left side of the $L^\infty$ estimates.

\subsection{The model problem of \cite{ionescu2014global}}\label{model_problem_motivation_sec}

The remaining two problems are related to a model problem posed in \cite{ionescu2014global}. We take a moment to introduce this model problem to provide a better context for the final two problems.

To derive the model problem, one starts with the Einstein Vacuum Equations
\begin{equation}\label{eve}
Ric[g] = 0
\end{equation}
and imposes axisymetry, since the Kerr solutions themselves have axisymmetry. As previously discussed, this excludes the issue of superradiance and simplifies the set of trapped null geodesics. Furthermore, it significantly simplifies the Einstein Vacuum system (\ref{eve}) as we shall now discuss.

An axisymmetric spacetime $(\mathcal{M},\mathbf{g})$ can be described by the restriction $g$ of the metric $\mathbf{g}$ to the quotient spacetime under the axisymmetry, together with a new complex scalar quantity $\sigma$. Equation (\ref{eve}) takes the form
\begin{align}
\Box_g\sigma &= N[g,\sigma]\label{eve_pi_eqn} \\
Ric(g)_{ij} &= N[g,\sigma]_{ij}.\label{eve_ric_eqn}
\end{align}
For a derivation of this reduction, see \cite{weinstein}.

The model problem of \cite{ionescu2014global} comes from yet another simplification: The equation for the reduced metric (\ref{eve_ric_eqn}) is ignored, and $g$ is replaced with the reduced Kerr metric $g_{Kerr}$. So only the evolution of the scalar $\sigma$ according to equation (\ref{eve_pi_eqn}) with $g=g_{Kerr}$ is studied. This has the simplifying advantage that the wave dynamics are fixed.

If
$$\sigma=X+iY,$$
then equation (\ref{eve_pi_eqn}) is given by the following system of equations for $X$ and $Y$.
\begin{align}
\Box_gX &= \frac{\pd^\alpha X\pd_\alpha X}{X}-\frac{\pd^\alpha Y\pd_\alpha Y}{X}\label{wm_X_eqn} \\
\Box_gY &= 2\frac{\pd^\alpha X\pd_\alpha Y}{X}.\label{wm_Y_eqn}
\end{align}
Coincidentally, if $X$ and $Y$ are taken to be the standard coordinates for the hyperbolic half plane with ranges $X\in (0,\infty)$ and $Y\in (-\infty,\infty)$, then these equations are precisely the equations that govern wave maps from the Kerr spacetime to the hyperbolic half plane. For this reason, the system (\ref{wm_X_eqn}-\ref{wm_Y_eqn}) is referred to as the \textit{wave map system}.

A particular (nontrivial) solution to this system is given by the scalar $\sigma_0$ corresponding to the Kerr metric itself.
$$\sigma_0=A+iB,$$
\begin{align*}
A &= \frac{(r^2+a^2)^2-a^2\sin^2\theta(r^2-2Mr+a^2)}{r^2+a^2\cos^2\theta}\sin^2\theta \\
B &= -2aM(3\cos\theta-\cos^3\theta)-\frac{2a^3M\sin^4\theta\cos\theta}{r^2+a^2\cos^2\theta}.
\end{align*}
The purpose of the remaining two problems in this thesis is to investigate the stability of the scalar $\sigma_0=A+iB$ as an axisymmetric solution to the wave map system (\ref{wm_X_eqn}-\ref{wm_Y_eqn}).

A general fact for any truly vacuum axisymmetric spacetime (ie. a spacetime that solves both equations (\ref{eve_pi_eqn}-\ref{eve_ric_eqn}) and thus equation (\ref{eve})) is that the imaginary part $\Im(\sigma)=Y$, called the \textit{Ernst potential}, must be constant on each connected segment of the axis. (Segments are disconnected if they lie on opposite ends of a black hole.) Furthermore, the difference between the constant values of $Y$ on opposite sides of a black hole is directly related to the angular momentum of the black hole. (One can easily see this is the case for the Kerr spacetime by evaluating $B$ at $\theta=0$ and $\theta=\pi$.) This general fact is derived in \S\ref{ernst_potential_axis_derivation}. The space of perturbations of $\sigma_0$ that will be studied in this thesis are such that $Y-B$ vanishes on the entire axis. We therefore interpret these perturbations as preserving angular momentum.

\subsection{Problem 4: axisymmetric perturbations of the nontrivial wave map in the Schwarzschild spacetime}

The fourth problem is a special case of the final problem, but it can be solved using a simpler method. It seeks axisymmetric solutions
\begin{align*}
\pd_\phi X &= 0 \\
\pd_\phi Y &= 0
\end{align*}
 to the wave map system (\ref{wm_X_eqn}-\ref{wm_Y_eqn}) that are at least initially close to the particular solution
\begin{align*}
X_0 &= A = r^2\sin^2\theta \\
Y_0 &= B = 0.
\end{align*}
This problem, which was first posed in \cite{ionescu2014global}, more closely resembles the full Einstein Vacuum system (\ref{eve}) with axisymmetry.

By making a linearization
\begin{align*}
X &= A + A\phi \\
Y &= A^2\psi,
\end{align*}
(with $\phi$ a dynamic quantity not to be confused with the azimuthal coordinate) it is possible to represent the system as a set of wave equations decoupled at the linear level.
\begin{align*}
\Box_g\phi &= \mathcal{N}_\phi[\phi,\psi] \\
\Box_{\tilde{g}}\psi &= \mathcal{N}_\psi[\phi,\psi].
\end{align*}
The metric $\tilde{g}$ corresponds to a $7+1$ dimensional analogue of the Schwarzschild spacetime. Due to the similarity between the two equations, one might initially hope that the techniques used to solve $\Box_g\phi=0$ can be easily adapted to solve $\Box_{\tilde{g}}\psi=0$. This is indeed the case.

It must be noted that the assumption $Y=A^2\psi$ ensures that if $\psi$ is regular on the axis, then $Y$ vanishes at least to fourth order on the axis. This is related to a geometric interpretation that the perturbation does not add angular momentum to the system. See \S\ref{ernst_potential_axis_derivation}.

The above linearization leads to a few nonlinear terms that are singular on the axis. This greatly complicates the process of applying angular derivatives as commutators. However, if one instead makes the special linearization
\begin{align*}
X &= A +A\phi \\
Y &= X^2\psi,
\end{align*}
then the singular terms disappear. This provides a solution to the fourth problem that is simpler than the ``$a=0$ version'' of the solution to the fifth and final problem.

\subsection{Problem 5: axisymmetric perturbations of the nontrivial wave map in slowly rotating Kerr spacetimes}

The fifth and final problem seeks axisymmetric solutions
\begin{align*}
\pd_\phi X &= 0 \\
\pd_\phi Y &= 0
\end{align*}
 to the wave map system  (\ref{wm_X_eqn}-\ref{wm_Y_eqn}) that are at least initially close to the particular solution
\begin{align*}
X_0 &= A =  \frac{(r^2+a^2)^2-a^2\sin^2\theta(r^2-2Mr+a^2)}{r^2+a^2\cos^2\theta}\sin^2\theta \\
Y_0 &= B =  -2aM(3\cos\theta-\cos^3\theta)-\frac{2a^3M\sin^4\theta\cos\theta}{r^2+a^2\cos^2\theta}
\end{align*}
for slowly rotating ($|a|\ll M$) Kerr spacetimes. This problem, which was first posed in \cite{ionescu2014global}, more closely resembles the full Einstein Vacuum system (\ref{eve}) with axisymmetry.

This problem draws on techniques developed in Problems 3 and 4, since Problem 3 deals with wave equations in Kerr spacetimes, and Problem 4 identifies a good $(\phi,\psi)$ linearization for the wave map problem. However, there are two serious complications for the general Kerr case.

The first serious complication is the fact that if $a\ne 0$, the system is coupled at the linear level. For this reason, it is not even clear a priori that a conserved energy quantity exists, much less a Morawetz estimate. However, there is a special bundle formalism for the wave map problem that makes both estimates possible (at least for $|a|\ll M$). This was a substantial contribution of Ionescu and Klainerman in \cite{ionescu2014global}, who went as far as to prove the $0$th order dynamic estimates. For more information about this bundle, see \S\ref{wave_map_bundle_sec}.

The second serious complication is the presence of terms that appear to be singular on the axis. Unlike in the case of Problem 4, it seems unlikely that a more clever linearization can remove these singular terms. To address this issue, a new technique, using the gothic operators $\fa$, $\fb$, $\fc^l$, and $\fd^l$ has been developed. (See Appendix \ref{regularity_sec}.) This technique requires a new special structural condition on the axis, which is expected to be satisfied by any truly geometric problem that produces these apparently singular terms when expressed in a coordinate system adapted to the axis. It is likely to be useful for many other problems posed in a degenerate coordinate system.

\section{An outline of this thesis}

We complete the introduction by outlining the remainder of this thesis.

Chapter \ref{preliminaries_chap} reviews a number of basic facts for the benefit of a reader who is relatively unacquainted with the subject of this thesis. It begins with some geometric calculations that are used (often implicitly) throughout this thesis and then discusses general ingredients involved in constructing spacetime estimates.

Chapter \ref{mink_chap} discusses the first problem of this thesis. In terms of the seven step procedure outlined in \S\ref{i:steps_sec}, step 1 is presented in \S\ref{m:energy_sec}, step 2 is presented in \S\ref{m:morawetzII_sec}, step 3 is presented in \S\ref{m:rp_sec}, step 4 is presented \S\ref{m:dynamic_sec}, step 5 is presented in \S\ref{m:pointwise_sec}, and steps 6 and 7 are presented in \S\ref{m:main_sec}. In addition to the steps required to solve the problem are \S\ref{m:morawetzI_sec}, which reviews a more well-known Morawetz estimate with strong weights, \S\ref{m:special_sec}, which reviews some estimates that are special to Minkowski space, and \S\ref{m:decay_sec}, which explains the decay mechanism that is used for all the problems in this thesis.

Chapter \ref{szd_chap} discusses the second problem of this thesis. In terms of the seven step procedure outlined in \S\ref{i:steps_sec}, step 1 is presented in \S\ref{s:energy_sec}, step 2 is presented in \S\ref{s:morawetz_sec}, step 3 is presented in \S\ref{s:rp_sec}, step 4 is presented \S\ref{s:dynamic_sec}, step 5 is presented in \S\ref{s:pointwise_sec}, and steps 6 and 7 are presented in \S\ref{s:main_sec}.

Chapter \ref{kerr_chap} discusses the third problem of this thesis. In terms of the seven step procedure outlined in \S\ref{i:steps_sec}, step 1 is presented in \S\ref{k:energy_sec}, step 2 is presented in \S\ref{k:morawetz_sec}, step 3 is presented in \S\ref{k:rp_sec}, step 4 is presented \S\ref{k:dynamic_sec}, step 5 is presented in \S\ref{k:pointwise_sec}, and steps 6 and 7 are presented in \S\ref{k:main_sec}.

Chapter \ref{wm_szd_chap} discusses the fourth problem of this thesis. In terms of the seven step procedure outlined in \S\ref{i:steps_sec}, steps 1-3 are presented in \S\ref{ws:spacetime_estimates_sec}, step 4 is presented in \S\ref{ws:dynamic_estimates_sec}, step 5 is presented in \S\ref{ws:pointwise_sec}, step 6 is presented in \S\ref{ws:structure_sec}, and step 7 is presented in \S\ref{ws:main_thm_sec}.

Chapter \ref{wm_kerr_chap} discusses the fifth and final problem of this thesis. In terms of the seven step procedure outlined in \S\ref{i:steps_sec}, steps 1-2 are presented in \S\ref{wk:xi_estimates_sec} and then again in \S\ref{wk:phi_psi_estimates_sec}, step 3 is presented in \S\ref{wk:phi_psi_estimates_sec}, step 4 is presented in \S\ref{wk:dynamic_estimates_sec}, step 5 is presented in \S\ref{wk:pointwise_sec}, step 6 is presented in \S\ref{wk:structure_sec}, and step 7 is presented in \S\ref{wk:main_thm_sec}.

In addition to the above chapters are a few appendices. Appendix \ref{wave_map_appendix} provides additional general information about the wave map problem. Appendix \ref{elliptic_sec} proves standard elliptic estimates that are used throughout this thesis. Appendix \ref{regularity_sec} outlines a general formalism that is necessary to solve the final problem and should be useful for many geometric problems with a degenerate coordinate system.

\chapter{Preliminaries}\label{preliminaries_chap}

This chapter is for the benefit of a new student.

\section{Geometric calculations in coordinates}

This thesis studies equations with geometric structure, but many calculations will require a choice of coordinates. We begin by reviewing a few computational tricks that will be useful throughout this thesis.

\subsection{The metric and its inverse under a change of coordinates}

We begin by reviewing how the metric and its inverse transform under a change of coordinates.

In practice, a change of coordinates is given by a relation
\begin{align*}
x^1&=x^1(y^1,...,y^n) \\
&\vdots \\
x^n&=x^n(y^1,...,y^n).
\end{align*}
Under such a change of coordinates, one can use the formula
$$g_{\mu'\nu'}=\frac{\pd x^\mu}{\pd y^{\mu'}}\frac{\pd x^\nu}{\pd y^{\nu'}}g_{\mu\nu},$$
or alternatively, one can transform the quantity
$$g_{\mu\nu}dx^\mu dx^\nu\rightarrow g_{\mu'\nu'}dy^{\mu'}dy^{\nu'},$$
computing $dx^\mu$ by using the formula for the exterior derivative of a function
$$df(y^1,...,y^n)=\pd_{y^1}f dy^1+...+\pd_{y^n}f dy^n.$$
Although both approaches are clearly equivalent, in practice the latter approach is often easier to do by hand than the former. An example of translating between cartesian and polar coordinates in the plane is given below.

\begin{example}
Consider the change between cartesian and polar coordinates in the 2-dimensional euclidean plane given by
\begin{align*}
x&=r\cos\theta \\
y&=r\sin\theta.
\end{align*}
Given the metric in cartesian coordinates, we compute the metric in polar coordinates using the exterior derivative approach.
$$g_{\mu\nu}dx^\mu dx^\nu= dx^2+dy^2=(\cos\theta dr-r\sin\theta d\theta)^2+(\sin\theta dr+\cos\theta d\theta)^2=dr^2+r^2d\theta^2.$$
From the last expression, we conclude that $g_{\mu'\nu'}$ is diagonal and that $g_{rr}=1$ and $g_{\theta\theta}=r^2$.
\end{example}

To transform the inverse metric under a change of coordinates, one can transform the metric and compute its inverse, or one can use the formula
$$g^{\mu\nu}=\frac{\pd x^\mu}{\pd y^{\mu'}}\frac{\pd x^\nu}{\pd y^{\nu'}}g^{\mu'\nu'},$$
or alternatively, one can transform the quantity
$$g^{\mu'\nu'}\pd_{\mu'}\psi\pd_{\nu'}\psi \rightarrow g^{\mu\nu}\pd_\mu\psi\pd_\nu\psi,$$
computing $\pd_{\mu'}\psi$ by using the chain rule
$$\pd_{\mu'}\psi=\frac{\pd x^\mu}{\pd y^{\mu'}}\pd_{\mu'}\psi.$$
(The function $\psi$ is an arbitrary differentiable function.) Although all three approaches are equivalent, in practice the third approach is often easier to do by hand than either of the other two. An example of translating between cartesian and polar coordinates in the plane is given below.

\begin{example}
Consider the change of coordinates given in the previous example. For an arbitrary differentiable function $\psi$, by the chain rule,
$$\pd_r\psi=\cos\theta\pd_x\psi+\sin\theta\pd_y\psi$$
$$\pd_\theta\psi =-r\sin\theta\pd_x\psi+r\cos\theta\pd_y\psi.$$
Therefore,
\begin{align*}
g^{\mu'\nu'}\pd_{\mu'}\psi\pd_{\nu'}\psi&=(\pd_r\psi)^2+r^{-2}(\pd_\theta\psi)^2 \\
&=(\cos\theta\pd_x\psi+\sin\theta\pd_y\psi)^2+(-r\sin\theta\pd_x\psi+r\cos\theta\pd_y\psi)^2 \\
&=(\pd_x\psi)^2+(\pd_y\psi)^2.
\end{align*}
From the last expression, we conclude that $g^{\mu\nu}$ is diagonal and that $g^{xx}=1$ and $g^{yy}=1$.
\end{example}

\subsection{Formulas for the divergence and laplacian (or wave operator)}

Let $V=V^\mu\pd_\mu$ be a vectorfield. To compute the divergence of $V$, one may use the formula
$$div V=\nabla_\mu V^\mu =\pd_\mu V^\mu+\Gamma_{\mu\nu}^\mu V^\nu$$
and compute (or have handy a table of) the relevant Christoffel symbols using the formula
$$\Gamma_{\mu\nu}^\lambda = \frac12g^{\lambda \alpha}(\pd_\mu g_{\alpha\nu}+\pd_\nu g_{\mu\alpha}-\pd_\alpha g_{\mu\nu}),$$
or alternatively, one can use the formula
\begin{equation}\label{divergence_eqn}
div V=\frac1{\sqrt{g}}\pd_\mu\left(\sqrt{g}V^\mu\right),
\end{equation}
where $\sqrt{g}=\sqrt{det g}$ is the volume form for the metric.

This latter approach is often more practical as it avoids explicit reference to the Christoffel symbols. The motivation for the formula given in equation (\ref{divergence_eqn}) is that it is consistent with the divergence theorem. See \S\ref{int_and_div_thm_sec}.

\begin{example}
Let $V=f(r)\pd_r$ be an arbitrary radial vectorfield in 3-dimensional euclidean space. Then using spherical coordinates,
$$div V=\frac1{\sqrt{g}}\pd_\mu(\sqrt{g}V^\mu)=\frac1{r^2\sin\theta}\pd_r(r^2\sin\theta V^r)=\frac1{r^2}\pd_r\left(r^2 f(r)\right)=f'(r)+\frac{2}rf(r).$$
\end{example}

It is also possible to compute the laplacian (or wave operator, depending on the signature of the metric) using equation (\ref{divergence_eqn}) and setting $V^\mu=\pd^\mu u=g^{\mu\nu}\pd_\nu u$. That is,
$$\triangle u=div(\pd u)=\frac1{\sqrt{g}}\pd_\mu\left(\sqrt{g}g^{\mu\nu}\pd_\nu u\right).$$

\begin{example}
Consider the 2-dimensional unit sphere with standard spherical coordinates $(\theta,\phi)$. The laplacian of an axisymmetric function $u=u(\theta)$ is given by
$$\sla{\triangle} u=\frac1{\sin\theta}\pd_\mu\left(\sin\theta g^{\mu\nu}\pd_\nu u\right)=\frac{1}{\sin\theta}\pd_\theta\left(\sin\theta \pd_\theta u\right)=\pd_\theta^2 u+\cot\theta\pd_\theta u.$$
\end{example}

\subsection{Integration and the divergence theorem in coordinates}\label{int_and_div_thm_sec}

The divergence theorem states that 
\begin{equation}\label{divergence_thm_eqn}
\int\int_{\Omega}div V =\int_{\pd\Omega}V\cdot n.
\end{equation}

In many applications, the domain $\Omega$ can be expressed as a cube in some set of coordinates.
$$\Omega=[x^1_{min},x^1_{max}]\times...\times [x^d_{min},x^d_{max}].$$
In this situation, the divergence theorem can be seen to be consistent with the forumla (\ref{divergence_eqn}) and a useful formula for the flux of $V$ becomes apparent.
\begin{proposition}
Equation (\ref{divergence_thm_eqn}) is consistent with equation (\ref{divergence_eqn}).
\end{proposition}
\begin{proof}
We begin by computing the left side of equation (\ref{divergence_thm_eqn}) in an adapted coordinate system, using equation (\ref{divergence_eqn}).
\begin{align*}
\int_{x^1_{min}}^{x^1_{max}}...\int_{x^d_{min}}^{x^d_{max}}div V \sqrt{g}dx^d...dx^1 &= \int_{x^1_{min}}^{x^1_{max}}...\int_{x^d_{min}}^{x^d_{max}}\frac1{\sqrt{g}}\pd_\mu(\sqrt{g}V^\mu)\sqrt{g}dx^d...dx^1 \\
&=\int_{x^1_{min}}^{x^1_{max}}...\int_{x^d_{min}}^{x^d_{max}}\pd_\mu(\sqrt{g}V^\mu)dx^d...dx^1.
\end{align*}
Then by the fundamental theorem of calculus, the last expression can be evaluated as a sum of integrals on hypersurfaces of constant $x^i$ ($i$ ranging from $1$ to $d$).
\begin{multline*}
\int_{x^1_{min}}^{x^1_{max}}...\int_{x^d_{min}}^{x^d_{max}}\pd_\mu(\sqrt{g}V^\mu)dx^d...dx^1 \\
=\left.\int_{x^2_{min}}^{x^2_{max}}...\int_{x^d_{min}}^{x^d_{max}}\sqrt{g}V^{x^1}dx^d...dx^2\right|_{x^1_{min}}^{x^1_{max}} +...+\left.\int_{x^1_{min}}^{x^1_{max}}...\int_{x^{d-1}_{min}}^{x^{d-1}_{max}}\sqrt{g}V^{x^d}dx^{d-1}...dx^1\right|_{x^d_{min}}^{x^d_{max}}
\end{multline*}
This sum is the flux of $V$ across the boundary of the coordinate cube $\Omega$, which corresponds to the right side of equation (\ref{divergence_thm_eqn}).
\end{proof}

\begin{example}
Let $V=r^{-2}\pd_r$ be the electric field corresponding to a point charge located at $r=0$ in 3-dimensional Euclidean space. Let $\Omega$ be the region between the sphere of radius 1 and the sphere of radius 2. In spherical coordinates, $\Omega=[1,2]\times[0,\pi]\times[0,2\pi]$.
\begin{align*}
\int\int_{\Omega} div V &= \int_1^2\int_0^\pi\int_0^{2\pi}\frac{1}{r^2\sin\theta}\pd_\mu(r^2\sin\theta V^\mu)r^2\sin\theta d\phi d\theta dr \\
&= \int_1^2\int_0^\pi\int_0^{2\pi}\pd_r(r^2\sin\theta V^r)d\phi d\theta dr \\
&= \left.\int_0^\pi\int_0^{2\pi} r^2\sin\theta V^r(r) d\phi d\theta\right|_{r=1}^{r=2} \\
&= 4\pi r^2 r^{-2}|_{r=1}^{r=2} \\
&= 4\pi - 4\pi \\
&= 0.
\end{align*}
This is to be expected, since $\Omega$ does not contain any charge.
\end{example}

\section{Ingredients for energy estimates}

An essential part of the analysis of wave equations is the construction of integrated estimates based on the divergence theorem. The general procedure is outlined here.

\subsection{Spacetime domains and boundary components}

Each of the spacetimes $(\mathcal{M},g)$ considered in this thesis will have coordinates $(t,r,\theta,\phi)$ similar to the spherical coordinates in Minkowski space. Let $\Omega\subset\mathcal{M}$ be a domain of the form
$$\Omega=[t_1,t_2]\times [r_H,\infty)\times S^2.$$
For the Minkowski spacetime, $r_H=0$, but for a black hole spacetime, $r_H>0$ corresponds to the radius of the event horizon. Let $\Sigma_t$ denote the constant-time hypersurface
$$\Sigma_t=\{t\}\times [r_H,\infty)\times S^2,$$
and for a black hole spacetime, let $H_{t_1}^{t_2}$ denote the part of the event horizon bounded between times $t_1$ and $t_2$
$$H_{t_1}^{t_2}=[t_1,t_2]\times\{r_H\}\times S^2.$$
Then for a black hole spacetime,
$$\pd\Omega = H_{t_1}^{t_2}\cup\Sigma_{t_1}\cup\Sigma_{t_2}.$$
(We will assume that all quantities vanish as $r\rightarrow\infty$, so we ignore the component of the boundary at $r=\infty$.)

\subsection{An integrated identity based on the divergence theorem}

Let $J$ be a vectorfield, called the \textit{current}. Then by the divergence theorem, if $J$ vanishes as $r\rightarrow\infty$,
$$\int_{t_1}^{t_2}\int_{\Sigma_t}div J=\int_{\Sigma_{t_2}}J^t+\int_{H_{t_1}^{t_2}}-J^r+\int_{\Sigma_{t_1}}-J^t.$$

By simply rearranging these terms, we establish the following proposition.
\begin{proposition}\label{general_divergence_estimate_prop}
Let $J$ be a current vanishing at an appropriate rate as $r\rightarrow\infty$. Then
$$\int_{H_{t_1}^{t_2}}J^r+\int_{\Sigma_{t_2}}-J^t+\int_{t_1}^{t_2}\int_{\Sigma_t}div J=\int_{\Sigma_{t_1}}-J^t.$$
\end{proposition}

To construct an integrated estimate, one must find a current $J$ (depending on the wavefunction $\psi$) such that $div J$, $-J^t$, and $J^r$ are nonnegative. \textbf{This is the foundation for all of the $L^2$-type estimates in this thesis.}

\subsection{The energy-momentum tensor}\label{energy_momentum_tensor_sec}

Now, we discuss the construction of useful currents $J$. We start with the energy-momentum tensor.

\begin{definition}\label{em_tensor_def}
For a wavefunction $\psi$, define the \textit{energy-momentum tensor} to be
$$T_{\mu\nu}=2\pd_\mu\psi\pd_\nu\psi-g_{\mu\nu}\pd^\lambda\psi\pd_\lambda\psi.$$
\end{definition}

The energy-momentum tensor is fundamental to the analysis of linear waves for the following reason.
\begin{proposition}\label{divT_prop}
If $\Box_g\psi=0$, then 
$$\nabla^\mu T_{\mu\nu}=0.$$
More generally, for any regular wavefunction $\psi$,
$$\nabla^\mu T_{\mu\nu}=2\Box_g\psi\pd_\nu\psi .$$
\end{proposition}
\begin{proof}
Since the connection $\nabla$ is compatible with the metric, $\nabla^{\mu}g_{\alpha\beta}=0$. Thus,
\begin{align*}
\nabla^\mu T_{\mu\nu} &= \nabla^\mu(2\pd_\mu\psi\pd_\nu\psi-g_{\mu\nu}\pd^\lambda\psi\pd_\lambda\psi)\\
&= 2\nabla^\mu\pd_\mu\psi \pd_\nu\psi+2\pd_\mu\psi\nabla^\mu\pd_\nu\psi - g_{\mu\nu}\nabla^\mu(\pd^\lambda\psi\pd_\lambda\psi) \\
&= 2\Box_g\psi\pd_\nu\psi+2g^{\alpha\beta}\pd_\alpha\psi \nabla_\beta\pd_\nu\psi-g^{\alpha\beta}\nabla_\nu\pd_\alpha\psi\pd_\beta\psi-g^{\alpha\beta}\pd_\alpha\psi\nabla_\nu\pd_\beta\psi \\
&= 2\Box_g\psi\pd_\nu\psi.
\end{align*}
\end{proof}

The fact that the energy-momentum tensor is divergence-free for solutions to the wave equation leads to conservation laws in spacetimes with symmetry (see \S\ref{general_isometries_sec}). This is the underlying reason for the conservation of energy in stationary spacetimes.

\subsection{The current templates $J[X,w]$ and $J[X,w,m]$}\label{current_template_sec}

\begin{definition}\label{current_def}
For a wavefunction $\psi$, a vectorfield $X$, and a function $w$, define the ($X,w$)-current by
$$J[X,w]_\mu=T_{\mu\nu}X^\nu+w\psi\pd_\mu\psi-\frac12\psi^2\pd_\mu w,$$
where
$$T_{\mu\nu}=2\pd_\mu\psi\pd_\nu\psi-g_{\mu\nu}\pd^\lambda\psi\pd_\lambda\psi$$
is the energy-momentum tensor from Definition \ref{em_tensor_def}.
\end{definition}

\begin{remark}
The vectorfield $X$ in the above definition is called a \textit{vectorfield multiplier}.
\end{remark}

\begin{lemma}\label{divJ_lem} (Identity for divJ)
For any regular wavefunction $\psi$,
$$div J[X,w]=K^{\mu\nu}\pd_\mu\psi\pd_\nu\psi+K\psi^2+(2X(\psi)+w\psi)\Box_g\psi,$$
where
\begin{align*}
K^{\mu\nu}&=2\nabla^{(\mu}X^{\nu)}+(w-div X)g^{\mu\nu} \\
&=g^{\mu\lambda}\pd_\lambda X^\nu+g^{\nu\lambda}\pd_\lambda X^\mu-X^\lambda\pd_\lambda g^{\mu\nu}+(w-div X) g^{\mu\nu},
\end{align*}
and
$$K=-\frac12\Box_g w.$$
\end{lemma}

\begin{proof}
Recall that
$$J[X,w]_\mu=T_{\mu\nu}X^\nu+w\psi\pd_\mu\psi-\frac12\psi^2\pd_\mu w,$$
where $T_{\mu\nu}$ is the energy-momentum tensor. The goal is to compute $\nabla^\mu J[X,w]_\mu$.

By Proposition \ref{divT_prop} and the symmetry of $T_{\mu\nu}$ in its indices,
\begin{align}
\nabla^\mu(T_{\mu\nu}X^\nu)&=T_{\mu\nu}\nabla^\mu X^\nu+2\Box_g\psi X^\nu\pd_\nu\psi \nonumber \\
&=T_{\mu\nu}\nabla^{(\mu}X^{\nu)}+2X(\psi)\Box_g\psi \nonumber \\
&= \nabla^{(\mu}X^{\nu)}2\pd_\mu\psi\pd_\nu\psi -divX g^{\mu\nu}\pd_\mu\psi\pd_\nu\psi+2X(\psi)\Box_g\psi. \label{div_TX_identity_eqn}
\end{align}
A similar identity holds for the $w$ terms.
\begin{align*}
\nabla^\mu\left(w\psi\pd_\mu\psi-\frac12\psi^2\pd_\mu w\right) &= \pd^\mu w \psi\pd_\mu\psi+w\pd^\mu\psi\pd_\mu\psi+w\psi\nabla^\mu\pd_\mu\psi-\psi\pd^\mu\psi\pd_\mu w-\frac12\psi^2\nabla^\mu\pd_\mu w \\
&=w\pd^\mu\psi\pd_\mu\psi+w\psi\nabla^\mu\pd_\mu\psi-\frac12\psi^2\nabla^\mu\pd_\mu w \\
&=w g^{\mu\nu}\pd_\mu\psi\pd_\nu\psi-\frac12 \psi^2\Box_g w +w\psi\Box_g\psi.
\end{align*}
Collecting all terms in the divergence of $J$, we obtain
\begin{align*}
\nabla^\mu J[X,w]_\mu &= \nabla^\mu(T_{\mu\nu}X^\nu)+\nabla^\mu\left(w\psi\pd_\mu\psi-\frac12\psi^2\pd_\mu w\right) \\
&=\left(2\nabla^{(\mu}X^{\nu)}+(w-divX)g^{\mu\nu}\right)\pd_\mu\psi\pd_\nu\psi-\frac12\Box_g w \psi^2+(2X(\psi)+w\psi)\Box_g\psi.
\end{align*}
Thus, we define the tensor
$$K^{\mu\nu}=2\nabla^{(\mu}X^{\nu)}+(w-divX)g^{\mu\nu}$$
and scalar
$$K=-\frac12\Box_gw.$$

In practice, it helps to use a more coordinate-dependent formula for $K^{\mu\nu}$ that can be computed directly from the components $g^{\mu\nu}$, $X^\mu$, and $w$ and their derivatives. To obtain this formula, we compute the divergence of $T_{\mu\nu}X^\nu$ a different way.
\begin{align*}
\nabla^\mu(T_{\mu\nu}X^\nu) =& \nabla_\lambda(g^{\lambda\mu}X^\nu T_{\mu\nu})\\
=& \frac1{\sqrt{-g}}\pd_\lambda\left(\sqrt{-g}g^{\lambda\mu}X^\nu(2\pd_\mu\psi\pd_\nu\psi-g_{\mu\nu}g^{\alpha\beta}\pd_\alpha\psi\pd_\beta\psi)\right) \\
=& 2\frac1{\sqrt{-g}}\pd_\lambda\left(\sqrt{-g}g^{\lambda\mu}\pd_\mu\psi\right)X^\nu\pd_\nu\psi-\frac1{\sqrt{-g}}\pd_\lambda\left(\sqrt{-g} g^{\lambda\mu}g_{\mu\nu}X^\nu\right)g^{\alpha\beta}\pd_\alpha\psi\pd_\beta\psi \\
&+ 2\pd_\lambda\left(X^\nu\pd_\nu\psi\right)g^{\lambda\mu}\pd_\mu\psi - g^{\lambda\mu}g_{\mu\nu}X^\nu\pd_\lambda\left(g^{\alpha\beta}\pd_\alpha\psi\pd_\beta\psi\right).
\end{align*}
Since
\begin{multline*}
2\frac1{\sqrt{-g}}\pd_\lambda\left(\sqrt{-g}g^{\lambda\mu}\pd_\mu\psi\right)X^\nu\pd_\nu\psi-\frac1{\sqrt{-g}}\pd_\lambda\left(\sqrt{-g} g^{\lambda\mu}g_{\mu\nu}X^\nu\right)g^{\alpha\beta}\pd_\alpha\psi\pd_\beta\psi \\
= 2\Box_g\psi X(\psi)- divX g^{\alpha\beta}\pd_\alpha\psi\pd_\beta\psi,
\end{multline*}
it follows that
\begin{multline*}
\nabla^\mu(T_{\mu\nu}X^\nu)=2\Box_g\psi X(\psi)-divXg^{\alpha\beta}\pd_\alpha\psi\pd_\beta\psi \\
+2\pd_\lambda\left(X^\nu\pd_\nu\psi\right)g^{\lambda\mu}\pd_\mu\psi - g^{\lambda\mu}g_{\mu\nu}X^\nu\pd_\lambda\left(g^{\alpha\beta}\pd_\alpha\psi\pd_\beta\psi\right).
\end{multline*}
Comparing this to (\ref{div_TX_identity_eqn}), we conclude that
\begin{equation}
\nabla^{(\mu}X^{\nu)}\pd_\mu\psi\pd_\nu\psi=2\pd_\lambda\left(X^\nu\pd_\nu\psi\right)g^{\lambda\mu}\pd_\mu\psi - g^{\lambda\mu}g_{\mu\nu}X^\nu\pd_\lambda\left(g^{\alpha\beta}\pd_\alpha\psi\pd_\beta\psi\right). \label{compare_DT_terms_eqn}
\end{equation}
By direct computation,
\begin{multline*}
2\pd_\lambda\left(X^\nu\pd_\nu\psi\right)g^{\lambda\mu}\pd_\mu\psi - g^{\lambda\mu}g_{\mu\nu}X^\nu\pd_\lambda\left(g^{\alpha\beta}\pd_\alpha\psi\pd_\beta\psi\right)  \\
=2X^\nu g^{\lambda\mu}\pd_\lambda\pd_\nu\psi\pd_\mu\psi+2\pd_\lambda X^\nu g^{\lambda\mu}\pd_\mu\psi\pd_\nu\psi -X^\lambda\pd_\lambda(g^{\alpha\beta})\pd_\alpha\psi\pd_\beta\psi -2X^\lambda g^{\alpha\beta}\pd_\lambda\pd_\alpha\psi\pd_\beta\psi \\
=2\pd_\lambda X^\nu g^{\lambda\mu}\pd_\mu\psi\pd_\nu\psi -X^\lambda\pd_\lambda(g^{\alpha\beta})\pd_\alpha\psi\pd_\beta\psi \\
=\left(g^{\mu\lambda}\pd_\lambda X^\nu+g^{\nu\lambda}\pd_\lambda X^\mu-X^\lambda \pd_\lambda (g^{\mu\nu})\right)\pd_\mu\psi\pd_\nu\psi.
\end{multline*}
Since $\psi$ is arbitrary, we conclude from (\ref{compare_DT_terms_eqn}) that
$$\nabla^{(\mu}X^{\nu)}=g^{\mu\lambda}\pd_\lambda X^\nu+g^{\nu\lambda}\pd_\lambda X^\mu-X^\lambda \pd_\lambda (g^{\mu\nu}).$$
Therefore,
$$K^{\mu\nu}=g^{\mu\lambda}\pd_\lambda X^\nu+g^{\nu\lambda}\pd_\lambda X^\mu-X^\lambda\pd_\lambda g^{\mu\nu}+(w-div X) g^{\mu\nu}.$$
\end{proof}

Sometimes, it will be helpful to use a slightly more general current template defined as follows.
\begin{definition}
For a wavefunction $\psi$, a vectorfield $X$, a function $w$, and a second vectorfield $m$, define the $(X,w,m)$-current by
$$J[X,w,m]_\mu=J[X,w]_\mu+m_\mu\psi^2.$$
\end{definition}
Note that while the wave factor $\Box_g\psi$ showed up in $divJ[X,w]$ coupled to both $X$ and $w$, it does not couple to $m$ in $divJ[X,w,m]$. The additional term $m_\mu\psi^2$ is simply a convenient way to translate between divergence terms and boundary terms.

\subsection{Isometries and conformal transformations}\label{general_isometries_sec}

Given a vectorfield $X$ on a spacetime $(\mathcal{M},g)$, one can define a one-parameter family of transformations ${}^{(X)}\Phi_s:\mathcal{M}\rightarrow \mathcal{M}$ by
$${}^{(X)}\Phi_s(p)=\exp_p sX,$$
where $\exp_p$ is the exponential map. It follows that ${}^{(X)}\Phi_0=Id$ and
$$X=\left.\frac{d}{ds}\right|_{s=0}{}^{(X)}\Phi_s.$$
For a rank $(0,k)$ tensor $U_{\mu_1...\mu_k}$, we define the Lie derivative
$$\mathcal{L}_XU_{\mu_1...\mu_k}=\left.\frac{d}{ds}\right|_{s=0}\left({}^{(X)}\Phi_s\right)^*U_{\mu_1...\mu_k},$$
where $\left({}^{(X)}\Phi_s\right)^*$ is the pullback induced by the map ${}^{(X)}\Phi_s$. By using also the pushforward $\left({}^{(X)}\Phi_{-s}\right)_*$, one can extend the Lie derivative operator $\mathcal{L}_X$ to act on tensors of arbitrary rank. There is an abundance of literature on the Lie derivative, so it will not be discussed any further here, except to mention the important formula for the Lie derivative of the metric.
\begin{equation}\label{lie_derivative_of_metric_eqn}
\mathcal{L}_Xg_{\mu\nu}=2\nabla_{(\mu}X_{\nu)}.
\end{equation}

If $\mathcal{L}_Xg_{\mu\nu}=0$, then the maps ${}^{(X)}\Phi_s$ are isometries. More generally, if $\mathcal{L}_Xg_{\mu\nu}=fg_{\mu\nu}$ for some function $f$, then the maps ${}^{(X)}\Phi_s$ are conformal transformations and $f$ is the derivative with respect to $s$ (at $s=0$) of the conformal factor.

\begin{definition}
We say \textit{$X$ generates an isometry} if $\mathcal{L}_Xg_{\mu\nu}=0$ and that \textit{$X$ generates a conformal transformation} if $\mathcal{L}_Xg_{\mu\nu}=f g_{\mu\nu}$.
\end{definition}

\begin{proposition}
If $\Box_g\psi=0$ and $X$ generates an isometry, then
$$divJ[X,0]=0.$$
Note that $J[X,0]_\mu=T_{\mu\nu}X^\nu$.
\end{proposition}

We will actually prove something more general.

\begin{proposition}\label{conformal_current_prop}
On an ($n+1$)-dimensional spacetime, if $\Box_g\psi=0$ and $X$ generates a conformal transformation so that $\mathcal{L}_Xg_{\mu\nu}=fg_{\mu\nu}$, then
$$divJ\left[X,\frac{n-1}{2}f\right]=-\frac{n-1}4(\Box_gf) \psi^2.$$
\end{proposition}
\begin{proof}
By equation (\ref{lie_derivative_of_metric_eqn}), if $X$ generates a conformal transformation, then
$$2\nabla_{(\mu}X_{\nu)}=fg_{\mu\nu}.$$
It follows that
$$divX=g^{\mu\nu}\nabla_{(\mu}X_{\nu)}=g^{\mu\nu}\frac12f g_{\mu\nu}=\frac{n+1}{2}f.$$
Then from Lemma \ref{divJ_lem},
$$divJ[X,(n-1)f]=K^{\mu\nu}\pd_\mu\psi\pd_\nu\psi+K\psi^2,$$
where
\begin{align*}
K^{\mu\nu} &= 2\nabla^{(\mu}X^{\nu)}+(w-divX)g^{\mu\nu} \\
&=fg^{\mu\nu}+\left(\frac{n-1}2f-\frac{n+1}2f\right)g^{\mu\nu} \\
&=0
\end{align*}
and
$$K=-\frac12\Box_gw =-\frac{n-1}4\Box_gf.$$
\end{proof}

\begin{lemma}\label{conformal_current_lem}
Let $X$ generate a conformal transformation, and let $f$ be the factor in the relation
$$\mathcal{L}_Xg_{\mu\nu}=fg_{\mu\nu}.$$
Then
$$\Box_gf=-\frac{2}{n}\nabla^\lambda(Ric_\lambda{}^\sigma X_\sigma).$$
\end{lemma}
\begin{proof}
By a simple calculation,
$$\nabla^\lambda X_\lambda=g^{\mu\nu}\nabla_{(\mu}X_{\nu)}=g^{\mu\nu}\frac12fg_{\mu\nu}=\frac{n+1}{2}f.$$
Therefore,
$$\Box_g(\nabla^\lambda X_\lambda)=\frac{n+1}{2}\Box_gf.$$
We will prove the lemma by showing that
$$\Box_g(\nabla^\lambda X_\lambda)=\frac{1-n}{2}\Box_g f-2\nabla^\lambda(Ric_\lambda{}^\sigma X_\sigma).$$
We have the following calculation
\begin{align*}
\Box_g(\nabla^\lambda X_\lambda) &= g^{\alpha\beta}g^{\lambda\mu}\nabla_\alpha\nabla_\beta\nabla_\lambda X_\mu \\
&= g^{\alpha\beta}g^{\lambda\mu}\left(\nabla_\alpha\nabla_\lambda\nabla_\beta X_\mu +\nabla_\alpha(R_{\beta\lambda\mu}{}^\sigma X_\sigma)\right) \\
&= g^{\alpha\beta}g^{\lambda\mu}\left(\nabla_\lambda\nabla_\alpha\nabla_\beta X_\mu +R_{\alpha\lambda\beta}{}^\sigma\nabla_\sigma X_\mu+R_{\alpha\lambda\mu}{}^\sigma\nabla_\beta X_\sigma +\nabla_\alpha(R_{\beta\lambda\mu}{}^\sigma X_\sigma)\right) \\
&= \nabla^\lambda(\Box_g X_\lambda) + Ric^{\mu\sigma}\nabla_\sigma X_\mu-Ric^{\beta\sigma}\nabla_\beta X_\sigma -\nabla^\beta(Ric_\beta{}^\sigma X_\sigma) \\
&= \nabla^\lambda(\Box_g X_\lambda)  -\nabla^\beta(Ric_\beta{}^\sigma X_\sigma) \\
&= \nabla^\lambda(\Box_g X_\lambda+Ric_\lambda{}^\sigma X_\sigma) - 2\nabla^\lambda(Ric_\lambda{}^\sigma X_\sigma).
\end{align*}
Therefore, it suffices to show that
$$\nabla^\lambda(\Box_g X_\lambda+Ric_\lambda{}^\sigma X_\sigma)=\frac{1-n}2\Box_gf.$$

We start by defining the tensor
$$\Gamma_{\alpha\beta\lambda}=\nabla_\alpha \nabla_{(\beta}X_{\lambda)}+\nabla_\beta\nabla_{(\alpha}X_{\lambda)}-\nabla_\lambda\nabla_{(\alpha}X_{\beta)}.$$
On one hand, since $\nabla_{(\mu}X_{\nu)}=\frac12fg_{\mu\nu}$, we have that
\begin{align*}
2g^{\alpha\beta}\Gamma_{\alpha\beta\lambda} &= g^{\alpha\beta}\left(\nabla_\alpha(f g_{\beta\lambda})+\nabla_\beta(f g_{\alpha\lambda})-\nabla_\lambda(f g_{\alpha\beta})\right) \\
&= \nabla_\lambda f+\nabla_\lambda f-g^{\alpha\beta}g_{\alpha\beta}\nabla_\lambda f \\
&= 2\nabla_\lambda f-(n+1)\nabla_\lambda f \\
&= (1-n)\nabla_\lambda f.
\end{align*}
On the other hand, since $2\nabla_{(\mu}X_{\nu)}=\nabla_\mu X_\nu +\nabla_\nu X_\mu$, we have that
\begin{align*}
2g^{\alpha\beta}\Gamma_{\alpha\beta\lambda} &= g^{\alpha\beta}\left(\nabla_\alpha\nabla_\beta X_\lambda+\nabla_\alpha\nabla_\lambda X_\beta+\nabla_\beta\nabla_\alpha X_\lambda+\nabla_\beta\nabla_\lambda X_\alpha-\nabla_\lambda\nabla_\alpha X_\beta-\nabla_\lambda\nabla_\beta X_\alpha\right) \\
&= 2\Box_g X_\lambda + g^{\alpha\beta} R_{\alpha\lambda\beta}{}^\sigma X_\sigma + g^{\alpha\beta}R_{\beta\lambda\alpha}{}^\sigma X_\sigma \\
&= 2\Box_g X_\lambda +2Ric_\lambda{}^\sigma X_\sigma
\end{align*}
Thus,
$$\nabla^\lambda\left(\Box_g X_\lambda+Ric_\lambda{}^\sigma X_\sigma\right)=\nabla^\lambda(g^{\alpha\beta}\Gamma_{\alpha\beta\lambda})=\frac{1-n}{2}\nabla^\lambda\nabla_\lambda f = \frac{1-n}{2}\Box_g f.$$
\end{proof}

Therefore, we have the following corollary.
\begin{corollary}\label{conformal_current_cor}
On an $(n+1)$-dimensional vacuum spacetime, if $\Box_g\psi=0$ and $X$ generates a conformal transformation so that $\mathcal{L}_Xg_{\mu\nu}=fg_{\mu\nu}$, then
$$div J\left[X,\frac{n-1}2f\right]=0.$$
\end{corollary}
\begin{proof}
By Proposition \ref{conformal_current_prop} and Lemma \ref{conformal_current_lem}, for any spacetime,
$$div J\left[X,\frac{n-1}2f\right] = -\frac{n-1}4(\Box_gf)\psi^2= \frac{n-1}{2n}\nabla^\lambda(Ric_\lambda{}^\sigma X_\sigma) \psi^2.$$
Thus, for a vacuum ($Ric_{\mu\nu}=0$) spacetime,
$$div J\left[X,\frac{n-1}2f\right]=0.$$
\end{proof}

\chapter{Waves in the Minkowski spacetime}\label{mink_chap}

The first problem of this thesis investigates semilinear waves on the Minkowski background. Minkowski space is the simplest solution to the Einstein Vacuum Equations and corresponds to a completely empty spacetime. It is the principal object of study in Einstein's special theory of relativity.

Minkowski space differs from the other spacetimes studied in this thesis, because it does not correspond to a black hole. However, it is still important to consider for the following reasons: All other problems deal with asymptotically flat spacetimes, which is another way of saying they approximate Minkowski space far from the black hole, and the Minkowski metric is the simplest of all metrics, so many distracting complications can be omitted in the first presentation of the method of this thesis. Thus, Minkowski space is ideal for the first problem of this thesis.

The Minkowski metric is most well-known in Cartesian coordinates.
$$g_{\mu\nu}dx^\mu dx^\nu =-dt^2+dx^2+dy^2+dz^2.$$
However, we will generally use spherical coordinates, in which the metric takes the form
$$g_{\mu\nu}dx^\mu dx^\nu = -dt^2+dr^2+r^2(d\theta^2+\sin^2\theta d\phi^2).$$
Its volume form is
$$\mu = r^2\sin\theta$$
and the Lagrangian for the linear wave equation is
\begin{align*}
\mathcal{L} &=\mu g^{\alpha\beta}\pd_\alpha\psi\pd_\beta\psi \\
&= r^2\sin\theta\left[-(\pd_t\psi)^2+(\pd_r\psi)^2+\sla{g}^{\alpha\beta}\pd_\alpha\psi\pd_\beta\psi\right] \\
&= r^2\sin\theta\left[-(\pd_t\psi)^2+(\pd_r\psi)^2+r^{-2}(\pd_\theta\psi)^2+r^{-2}\sin^{-2}\theta(\pd_\phi\psi)^2\right]
\end{align*}
We briefly observe for some future calculations that the angular part of the inverse metric $\sla{g}^{\alpha\beta}$ satisfies $\pd_r(r^2\sla{g}^{\alpha\beta})=0$.

Many estimates will be derived with respect to a foliation of hypersurfaces $\Sigma_t$, each corresponding to a level set of the time coordinate $t$. For estimates involving $\Sigma_{t_1}$ and $\Sigma_{t_2}$, it should be understood that $\Sigma_{t_2}$ is in the future of $\Sigma_{t_1}$, that is, $t_2>t_1$. There is no event horizon in Minkowski space, so the component $H_{t_1}^{t_2}$ will not show up in this chapter.

\section{The energy estimate and the $h\pd_t$ estimate}\label{m:energy_sec}

The most fundamental of all estimates is the energy estimate. We begin by proving this estimate and then prove the slightly more general $h\pd_t$ estimate.

\subsection{The energy estimate}

Since the vectorfield $\pd_t$ is a Killing vectorfield, we have the following identity.
\begin{lemma}\label{m:ee_identity_lem}(Energy identity for Minkowski)
\begin{equation*}
\int_{\Sigma_{t_2}}-J^t[\pd_t] = \int_{\Sigma_{t_1}}-J^t[\pd_t]+\int_{t_1}^{t_2}\int_{\Sigma_t}-2\pd_t\psi\Box_g\psi
\end{equation*}
In particular, if $\Box_g\psi=0$, then the quantity
$$E(t)=\int_{\Sigma_t}-J^t[\pd_t]$$
is conserved.
\end{lemma}
\begin{proof}
By Proposition \ref{general_divergence_estimate_prop}, since there is no event horizon in Minkowski, we have
$$\int_{\Sigma_{t_2}}-J^t[\pd_t]+\int_{t_1}^{t_2}\int_{\Sigma_t}divJ=\int_{\Sigma_{t_1}}-J^t[\pd_t].$$
By Lemma \ref{divJ_lem},
$$divJ[X] = (g^{\mu\lambda}\pd_\lambda X^\nu+g^{\nu\lambda}\pd_\lambda X^\mu-X^\lambda\pd_\lambda g^{\mu\nu}-divX g^{\mu\nu})\pd_\mu\psi\pd_\nu\psi+2X(\psi)\Box_g\psi.$$
For the particular case $X=\pd_t$, the coefficients of $X$ are constant, the metric does not depend on $t$, and by a simple calculation (see equation (\ref{divergence_eqn})), $X$ is divergence free. Therefore,
$$divJ[\pd_t]=2\pd_t\psi\Box_g\psi.$$
We conclude that
$$\int_{\Sigma_{t_2}}-J^t[\pd_t]+\int_{t_1}^{t_2}\int_{\Sigma_t}2\pd_t\psi\Box_g\psi=\int_{\Sigma_{t_1}}-J^t[\pd_t].$$
The statement of the lemma now follows.
\end{proof}

Now, we calculate the flux terms in the previous lemma.

\begin{lemma}\label{m:ee_bndry_lem}
On each constant-time hypersurface $\Sigma_t$,
$$-J^t[\pd_t] = (\pd_t\psi)^2+(\pd_r\psi)^2+|\sla\nabla\psi|^2.$$
\end{lemma}
\begin{proof}
According to Definition \ref{current_def},
$$J^\mu[X]=2g^{\mu\lambda}\pd_\lambda\psi X^\nu\pd_\nu\psi-X^\mu\pd^\lambda\psi\pd_\lambda\psi.$$
Therefore,
\begin{align*}
J^t[\pd_t] &= 2g^{t\lambda}\pd_\lambda\psi \pd_t\psi-(\pd_t)^t\pd^\lambda\psi\pd_\lambda\psi \\
&= 2g^{tt}(\pd_t\psi)^2-\left(g^{tt}(\pd_t\psi)^2+g^{rr}(\pd_r\psi)^2+|\sla\nabla\psi|^2\right) \\
&= g^{tt}(\pd_t\psi)^2-g^{rr}(\pd_r\psi)^2-|\sla\nabla\psi|^2 \\
&= -(\pd_t\psi)^2-(\pd_r\psi)^2-|\sla\nabla\psi|^2.
\end{align*}
The statement of the lemma now follows.
\end{proof}

From the previous two lemmas, we conclude the following energy estimate.
\begin{proposition}\label{m:classic_ee_prop}(Energy estimate for Minkowski)
\begin{equation*}
\int_{\Sigma_{t_2}}\left[(\pd_r\psi)^2+(\pd_t\psi)^2+|\sla\nabla\psi|^2\right] 
 \lesssim \int_{\Sigma_{t_1}}\left[(\pd_r\psi)^2+(\pd_t\psi)^2+|\sla\nabla\psi|^2\right] + Err_\Box,
\end{equation*}
where 
$$Err_\Box=\int_{t_1}^{t_2}\int_{\Sigma_t}|\pd_t\psi\Box_g\psi|.$$
\end{proposition}
\begin{proof}
This follows directly from Lemmas \ref{m:ee_identity_lem} and \ref{m:ee_bndry_lem}.
\end{proof}

\subsection{The $h\pd_t$ energy estimate}\label{m:hdt_sec}

The next estimate, which uses a vectorfield $X$ of the form $h(r)\pd_t$, is a slight generalization of the energy estimate. The function $h$ introduces weights to the energy in different regions. There is a nice interpretation for this estimate. According Lemma \ref{m:divJhdt_lem} below,
$$divJ[h\pd_t]=\frac{h'}2\left[(L\psi)^2-(\lbar\psi)^2\right].$$
Given that $L$ measures incoming waves and $\lbar$ measures outgoing ones, this says that an energy with a larger weight far away (so that $h'>0$) will grow if there are more outgoing waves than incoming waves and vice versa. In our application of this estimate later on, we will use large weights near the center of the spacetime.

\begin{lemma}\label{m:divJhdt_lem}
If $\Box_g\psi=0$, then
$$divJ[h\pd_t]=\frac{h'}{2}\left[(L\psi)^2-(\lbar\psi)^2\right].$$
\end{lemma}
\begin{proof}
Recall from Lemma \ref{divJ_lem} that
$$divJ[X]=K^{\mu\nu}\pd_\mu\psi\pd_\nu\psi,$$
where
$$K^{\mu\nu}=g^{\mu\lambda}\pd_\lambda X^\nu+g^{\nu\lambda}\pd_\lambda X^\mu-X^\lambda\pd_\lambda(g^{\mu\nu})-divXg^{\mu\nu}.$$
Since $\pd_t$ is killing, $\pd_t(g^{\mu\nu})=0$ and $div(\pd_t)=0$. Also, the only component of $X$ is the $t$ component and the only nonzero derivative of that component is the $\pd_r$ derivative, so
$$g^{\mu\lambda}\pd_\lambda X^\nu+g^{\nu\lambda}\pd_\lambda X^\mu=g^{\mu r}\pd_r X^\nu+g^{\nu r}\pd_r X^\mu.$$
It follows that the only possible nonzero $K^{\mu\nu}$ components are
$$K^{tr}+K^{rt} = 2g^{rr}h'.$$
We conclude that
\begin{align*}
divJ[h\pd_t]&=2 h'\pd_r\psi\pd_t\psi \\
&=\frac{h'}2\left[(\pd_r\psi+\pd_t\psi)^2-(\pd_r\psi-\pd_t\psi)^2\right] \\
&=\frac{h'}2\left[(L\psi)^2-(\lbar\psi)^2\right].
\end{align*}
\end{proof}

Taking $h$ to be a positive function decreasing to zero as $r\rightarrow\infty$ at a particular rate, we obtain the $h\pd_t$ estimate.
\begin{proposition}\label{m:hdt_prop}($h\pd_t$ estimate for Minkowski)
Let $R>0$ be any given radius. Then for all $\epsilon>0$ and $p<2$, there is a small constant $c_\epsilon$ and a large constant $C_\epsilon$, such that
\begin{multline*}
\int_{\Sigma_{t_2}}(1+r)^{p-2}\left[(\pd_r\psi)^2+(\pd_t\psi)^2+|\sla\nabla\psi|^2\right] \\
+\int_{t_1}^{t_2}\int_{\Sigma_t\cap\{R+1<r\}}c_\epsilon (1+r)^{p-3}(\lbar\psi)^2 \\
\lesssim \int_{\Sigma_{t_1}}(1+r)^{p-2}\left[(\pd_r\psi)^2+(\pd_t\psi)^2+|\sla\nabla\psi|^2\right]+Err,
\end{multline*}
where
\begin{align*}
Err&=Err_1+Err_\Box \\
Err_1&=\int_{t_1}^{t_2}\int_{\Sigma_t\cap\{R<r\}}\epsilon (1+r)^{-1}(L\psi)^2 \\
Err_\Box&=\int_{t_1}^{t_2}\int_{\Sigma_t}C_\epsilon (1+r)^{p-2}|\pd_t\psi\Box_g\psi|.
\end{align*}
\end{proposition}
\begin{remark}
Notice the weight in $Err_1$. For the estimate to be uniform in $\epsilon$, it must be $(1+r)^{-1}$. However, for a fixed $\epsilon$, the $r$ weight could be replaced by $(1+r)^{p-3}$. Since $p<2$, $(1+r)^{p-3}\lesssim (1+r)^{-1}$.
\end{remark}
\begin{proof}
Observe that for all $p<2$ and $\epsilon>0$, there is a function $h$ satisfying the conditions
$$h=1\text{ for }r<R,$$
$$-h'\ge 0,$$
$$-h'\le 2\epsilon (1+r)^{-1},$$
$$\exists c_\epsilon>0\text{ such that }2c_\epsilon (1+r)^{p-3}\le -h'\text{ for }r\ge R+1,$$
and
$$h=O(r^{p-2})\text{ for large }r.$$
(To construct the function $h$, it is perhaps easier to construct a positive function $-h'$ supported on the interval $r\in[R,\infty)$ and satisfying $\int_R^\infty-h'(r)dr=1$. In particular, this is only possible because $\int_R^\infty r^{-1}dr$ diverges.)

The boundary terms follow from Lemma \ref{m:ee_bndry_lem} and the fact that $J^\mu[h\pd_t]=h J^\mu[\pd_t]$. We do not worry about the size of $\epsilon$ on the boundary terms, because for a fixed $r$, $h$ increases and is bounded above by $1$ as $\epsilon\rightarrow 0$.
\end{proof}

\section{The Morawetz estimate, version I}\label{m:morawetzI_sec}

Morawetz-type estimates are concerned with the divergence term in Proposition \ref{general_divergence_estimate_prop}. They rely very heavily on the specific geometry of the spacetime, so the derivations in this section will not easily generalize to black hole spacetimes (although they will indeed generalize--this is a key part of this thesis).

In this section, we prove the first of two Morawetz-type estimates for Minkowski space. It is well-known, but not quite compatible with the strategy employed in the remainder of this chapter. The fundamental difference between the two estimates lies in the asymptotics as $r\rightarrow\infty$ and the level of care taken with the flux terms.

\subsection{The partial Morawetz estimate}

The construction of the Morawetz estimate is quite complicated. Before attempting a proof of the entire estimate, it helps to make the following intermediate objective:
\textit{Find a momentum $J[X_0,w_0]$ with the property that all terms in $div J[X_0,w_0]$ are non-negative.}
When attempting to prove a Morawetz estimate, this objective constitutes a substantial first step as there is good reason to believe that after a few additional tweaks, one can derive a truly useful estimate. For this reason, we refer to this as a \textit{partial Morawetz estimate}.

\begin{proposition}\label{m:partial_morawetz_I_prop}(Partial Morawetz estimate in Minkowski (version I)) Let
$$X_0=\frac{r}{1+r}\pd_r$$
and
$$w_0=\frac{2r+3}{(1+r)^2}.$$
Then the current $J[X_0,w_0]$ from Definition \ref{current_def} satisfies
$$ divJ[X_0,w_0] = \frac{2}{(1+r)^2}(\pd_r\psi)^2+\frac{2}{1+r}|\sla\nabla\psi|^2+\frac{r+4}{r(1+r)^4}\psi^2.$$
\end{proposition}

\begin{proof}
For the proof, we will use a strategy that generalizes to the remaining problems in this thesis. However, this particular proof will be significantly simpler. In this proof, we drop the subscript $0$ from $X_0$ and $w_0$ until the very end.

\begin{lemma}\label{m:pme_initial_lem}
With the choices $X=X^r(r)\pd_r$ and $w=w(r)$,
\begin{multline*}
divJ[X,w]= 
(wg^{tt}-r^{-2}\pd_r(r^2X^rg^{tt}))(\pd_t\psi)^2+(w-\pd_rX^r)|\sla\nabla\psi|^2 \\
+(wg^{rr}+\pd_rX^r g^{rr}-2r^{-1}X^rg^{rr}-X^r\pd_rg^{rr})(\pd_r\psi)^2 -\frac12r^{-2}\pd_r(r^2\pd_rw)\psi^2.
\end{multline*}
\end{lemma}
\begin{proof}
By Lemma \ref{divJ_lem}, we have that for general $X$ and $w$,
$$divJ[X,w]=K^{\mu\nu}\pd_\mu\psi\pd_\nu\psi+K\psi^2,$$
where
$$K^{\mu\nu}=g^{\mu\lambda}\pd_\lambda X^\nu+g^{\nu\lambda}\pd_\lambda X^\mu-X^\lambda\pd_\lambda g^{\mu\nu}+(w-divX)g^{\mu\nu}$$
and
$$K=-\frac12\Box_g w=-\frac12r^{-2}\pd_r(r^2\pd_rw).$$
A quick inspection reveals that the nonzero components of $K^{\alpha\beta}$ are $K^{rr}$, $K^{tt}$, and $\sla{K}^{\alpha\beta}$ (the components for the angular derivatives). Using the fact that $divX=r^{-2}\pd_r(r^2 X^r)$, we compute
\begin{align*}
K^{rr} &= wg^{rr}+g^{rr}\pd_r X^r+g^{rr}\pd_rX^r-X^r\pd_rg^{rr}-r^{-2}\pd_r(r^2X^r)g^{rr} \\
&= wg^{rr}+g^{rr}\pd_rX^r-2r^{-1}X^rg^{rr}-X^r\pd_rg^{rr}
\end{align*}
\begin{align*}
K^{tt} &= wg^{tt}-X^r\pd_rg^{tt}-r^{-2}\pd_r(r^2X^r)g^{tt} \\
&= wg^{tt}-r^{-2}\pd_r(r^2X^rg^{tt})
\end{align*}
\begin{align*}
\sla{K}^{\alpha\beta} &= w\sla{g}^{\alpha\beta}-X^r\pd_r\sla{g}^{\alpha\beta}-r^{-2}\pd_r(r^2X^r)\sla{g}^{\alpha\beta} \\
&= w\sla{g}^{\alpha\beta}-r^{-2}\pd_r(X^r(r^2\sla{g}^{\alpha\beta})) \\
&= w\sla{g}^{\alpha\beta}-\pd_rX^r\sla{g}^{\alpha\beta}
\end{align*}
(In the last line, we used the fact that $\pd_r(r^2\sla{g}^{\alpha\beta})=0$.)
\end{proof}

It is difficult to make sense of the coefficients in the above lemma, because most coefficients are sums of multiple terms. By making the choices $X^r(r)=u(r)v(r)$ and $w(r)=v(r)\pd_ru(r)$, each coefficient can be written as a single term.
\begin{lemma}
With the choices $X=u(r)v(r)\pd_r$ and $w=v(r)\pd_ru(r)$,
\begin{multline*}
divJ[X,w]= 
-ur^{-2}\pd_r(r^2vg^{tt})(\pd_t\psi)^2-u\pd_rv |\sla\nabla\psi|^2 \\
+r^2(g^{rr})^2u^{-1}\pd_r\left(\frac{u^2v}{r^2g^{rr}}\right)(\pd_r\psi)^2-\frac12r^{-2}\pd_r(r^2\pd_rw)\psi^2.
\end{multline*}
\end{lemma}
\begin{proof}
This follows directly by substituting $X^r=uv$ and $w=v\pd_ru$ into the expression in the previous lemma and combining terms.
\end{proof}

At this point, it is standard to choose a relation between $X$ and $w$ so that the coefficient of $(\pd_t\psi)^2$ vanishes. (It will be recovered later--see \S\ref{m:wdt_correction_I_sec}.) In terms of the functions $u$ and $v$, it suffices to set $v=r^{-2}$.
\begin{remark}
The function $v=r^{-2}$ is related to the geodesic potential for null trajectories in Minkowski space. This fact will be expounded upon in greater detail in Chapter \ref{szd_chap}.
\end{remark}

\begin{lemma}\label{m:divJ_in_terms_of_u_lem}
With the choice $v=r^{-2}$, the coefficient of the $(\pd_t\psi)^2$ term vanishes and furthermore,
$$divJ[X,w]=\frac{2u}{r^3}|\sla\nabla\psi|^2 
+2\pd_r\left(\frac{u}{r^2}\right)(\pd_r\psi)^2-\frac12r^{-2}\pd_r(r^2\pd_r(v\pd_ru))\psi^2.$$
To ensure that all of these coefficients are positive, it is necessary to find a positive function $u$ growing sufficiently fast that $\pd_r\left(\frac{u}{r^2}\right)$ is also positive and furthermore that the expression
$$-\frac12r^{-2}\pd_r(r^2\pd_r(v\pd_ru))=-\frac12r^{-2}\pd_r(r^2\pd_r(r^{-2}\pd_ru))$$
is also nonnegative.
\end{lemma}
If finding such a function $u$ seems tricky now, it only gets more complicated in the Schwarzschild and Kerr settings. A well-known choice for $u$ is given in the following lemma.
\begin{lemma}\label{m:choice_for_u_I_lem}
The function $u=\frac{r^3}{1+r}$, which is clearly positive, satisfies
$$2\pd_r\left(\frac{u}{r^2}\right)=\frac{2}{(1+r)^2}$$
and
$$-\frac12r^{-2}\pd_r(r^2\pd_r(r^{-2}\pd_ru))=\frac{r+4}{r(1+r)^4}.$$
\end{lemma}
\begin{proof}
First,
$$2\pd_r\left(\frac{u}{r^2}\right)=2\pd_r\left(\frac{r}{1+r}\right)=\frac{2}{(1+r)^2}.$$
This verifies the first identity.

Then we have the following intermediate calculations, which verify the second identity.
$$\pd_ru=\pd_r\left(\frac{r^3}{1+r}\right)=\frac{3r^2(1+r)-r^3}{(1+r)^2}=\frac{2r^3+3r^2}{(1+r)^2}$$
$$\pd_r(r^{-2}\pd_ru)=\pd_r\left(\frac{2r+3}{(1+r)^2}\right)=\frac{2(1+r)-2(2r+3)}{(1+r)^3}=-2\frac{r+2}{(1+r)^3}$$
\begin{multline*}
-\frac12r^{-2}\pd_r(r^2\pd_r(r^{-2}\pd_ru)) = r^{-2}\pd_r\left(r^2\frac{r+2}{(1+r)^3}\right) = r^{-2}\frac{(3r^2+4r)(1+r)-3(r^3+2r^2)}{(1+r)^4} \\
=r^{-2}\frac{r^2+4r}{(1+r)^4}=\frac{r+4}{r(r+1)^4}.
\end{multline*}
\end{proof}

Finally, we compute
$$X_0=uv\pd_r=\frac{r^3}{1+r}r^{-2}\pd_r=\frac{r}{1+r}\pd_r$$
$$w_0=v\pd_ru=r^{-2}\frac{2r^3+3r^2}{(1+r)^2}=\frac{2r+3}{(1+r)^2}.$$
This completes the proof of Proposition \ref{m:partial_morawetz_I_prop}.
\end{proof}

\subsection{Remaining issues}

The current $J[X_0,w_0]$ from Proposition \ref{m:partial_morawetz_I_prop} has a non-negative divergence, but there are a few things that must be done to establish a proper Morawetz estimate. \\
\bp The coefficient of $(\pd_t\psi)^2$ vanishes entirely in the partial Morawetz estimate. This is easily fixed by adding to $w$ the term $\epsilon_{\pd_t}w_{\pd_t}$. (See \S\ref{m:wdt_correction_I_sec}.) \\
\bp Some attention must be paid to the boundary terms, which have so far been ignored. The current $J[X_0,w_0]$ does not yield good boundary terms, but this is easily fixed by adding to $X$ the term $C\pd_t$, where $C$ is a sufficiently large constant. (See \S\ref{m:Cdt_correction_I_sec}.)

Thus, the final current $J$ will be of the form
\begin{align*}
J &= J[X,w], \\
X &= X_0 + C\pd_t, \\
w &= w_0 + \epsilon_{\pd_t}w_{\pd_t}.
\end{align*}

\subsection{The correction $\epsilon_{\pd_t}w_{\pd_t}$}\label{m:wdt_correction_I_sec}
By adding a function $w_{\pd_t}$, we can control $(\pd_t\psi)^2$ at the price of an error term that we can already control.
\begin{lemma}\label{m:wdt_correction_I_lem}
The function $w_{\pd_t}=-(1+r)^{-2}$ satisfies
$$(1+r)^{-2}(\pd_t\psi)^2 \lesssim divJ[0,w_{\pd_t}] + (1+r)^{-2}\left[(\pd_r\psi)^2+|\sla\nabla\psi|^2+r^{-1}(1+r)^{-1}\psi^2\right]$$
\end{lemma}
\begin{proof}
By Lemma \ref{divJ_lem}, we have that
\begin{align*}
divJ[0,w_{\pd_t}] &= w_{\pd_t}g^{\alpha\beta}\pd_\alpha\psi\pd_\beta\psi-\frac12\Box_gw_{\pd_t} \psi^2 \\
&= -w_{\pd_t}(\pd_t\psi)^2 +w_{\pd_t}\left[(\pd_r\psi)^2+|\sla\nabla\psi|^2\right]-\frac12\Box_gw_{\pd_t}\psi^2.
\end{align*}
It suffices to choose $w_{\pd_t}$ satisfying the following conditions.
$$(1+r)^{-2}\lesssim -w_{\pd_t},$$
$$|w_{\pd_t}| \lesssim (1+r)^{-2},$$
$$\left|\frac12\Box_gw_{\pd_t}\right| \lesssim r^{-1}(1+r)^{-3}.$$
With the choice $w_{\pd_t}=-(1+r)^{-2}$, the first two conditions are trivially satisfied, and
\begin{align*}
\frac12\Box_gw_{\pd_t} &= -\frac12r^{-2}\pd_r(r^2\pd_r(1+r)^{-2}) \\
&= r^{-2}\pd_r(r^2(1+r)^{-3}) \\
&= 2r^{-1}(1+r)^{-3}-3(1+r)^{-4} \\
&\lesssim r^{-1}(1+r)^{-3}.
\end{align*}
Thus, $w_{\pd_t}$ satisfies the remaining condition.
\end{proof}

\begin{corollary}\label{m:wdt_correction_I_cor}
There exists a small $\epsilon_{\pd_t}>0$ so that
$$(1+r)^{-2}\left[(\pd_t\psi)^2+(\pd_r\psi)^2\right]+(1+r)^{-1}|\sla\nabla\psi|^2+r^{-1}(1+r)^{-3}\psi^2\lesssim divJ[X_0,w_0+\epsilon_{\pd_t}w_{\pd_t}].$$
\end{corollary}
\begin{proof}
Recall from Proposition \ref{m:partial_morawetz_I_prop} that
$$(1+r)^{-2}(\pd_r\psi)^2+(1+r)^{-1}|\sla\nabla\psi|^2+r^{-1}(1+r)^{-3}\psi^2 \lesssim div J[X_0,w_0].$$
The proof follows from this fact and the inequality in Lemma \ref{m:wdt_correction_I_lem}.
\end{proof}

\subsection{The correction $C\pd_t$ and the flux terms}\label{m:Cdt_correction_I_sec}
\begin{lemma}\label{m:Cdt_correction_I_lem}
$$|J^t[X_0,w_0+\epsilon_{\pd_t}w_{\pd_t}]| \lesssim (\pd_t\psi)^2+(\pd_r\psi)^2+(1+r)^{-2}\psi^2$$
\end{lemma}
\begin{proof}
From Lemma \ref{divJ_lem}, for any $X$ and $w$,
\begin{align*}
J^\mu[X,w] &= 2g^{\mu\lambda}\pd_\lambda\psi X^\nu\pd_\nu\psi-X^\mu\pd^\lambda\psi\pd_\lambda\psi+w\psi g^{\mu\lambda}\pd_\lambda\psi-\frac12\psi^2g^{\mu\lambda}\pd_\lambda w.
\end{align*}
It follows that if $X^t=0$ and $\pd_tw=0$, then
\begin{align*}
-J^t[X,w] &= 2\pd_t\psi X^r\pd_r\psi +w\psi\pd_t\psi.
\end{align*}
Note that $X_0^t=0$ and $\pd_t(w+\epsilon_{\pd_t}w_{\pd_t})=0$. Furthermore, $|X_0^r|\le 1$ and $|w_0+\epsilon_{\pd_t}w_{\pd_t}|\lesssim (1+r)^{-1}$, so
\begin{align*}
|J^t[X_0,w_0+\epsilon_{\pd_t}w_{\pd_t}]| &\lesssim 2|\pd_t\psi\pd_r\psi|+(1+r)^{-1}|\psi\pd_t\psi|\\
&\lesssim (\pd_t\psi)^2+(\pd_r\psi)^2+(1+r)^{-2}\psi^2.
\end{align*}
\end{proof}

\begin{corollary}\label{m:Cdt_correction_I_cor}
If $C$ is a sufficiently large constant, then
\begin{multline*}
(1+r)^{-2}\left[(\pd_t\psi)^2+(\pd_r\psi)^2\right]+(1+r)^{-1}|\sla\nabla\psi|^2+r^{-1}(1+r)^{-3}\psi^2 \\
\lesssim divJ[X_0+C\pd_t,w_0+\epsilon_{\pd_t}w_{\pd_t}],
\end{multline*}
and
$$\int_{\Sigma_t}-J^t[X_0+C\pd_t,w_0+\epsilon_{\pd_t}w_{\pd_t}] \approx \int_{\Sigma_t} (\pd_t\psi)^2+(\pd_r\psi)^2+|\sla\nabla\psi|^2.$$
\end{corollary}
\begin{proof}
Note that $divJ[\pd_t]=0$, which is the fundamental reason for Proposition \ref{m:classic_ee_prop}. Thus, the first estimate reduces to Corollary \ref{m:wdt_correction_I_cor}.

The second estimate follows from the fact that
$$-J^t[\pd_t] = (\pd_t\psi)^2+(\pd_r\psi)^2+|\sla\nabla\psi|^2,$$
which was established in Lemma \ref{m:ee_bndry_lem}, together with Lemma \ref{m:Cdt_correction_I_lem} and the Hardy estimate,
$$\int_{\Sigma_t}r^{-2}\psi^2 \lesssim \int_{\Sigma_t}(\pd_r\psi)^2.$$
\end{proof}

\subsection{The Morawetz estimate}\label{m:morawetz_I_final_subsec}

\begin{theorem}(Morawetz estimate in Minkowski (version I))
The following inequality holds for any sufficiently regular function $\psi$.
\begin{multline*}
\int_{\Sigma_{t_2}}(\pd_t\psi)^2+(\pd_r\psi)^2+|\sla\nabla\psi|^2 \\
+\int_{t_1}^{t_2}\int_{\Sigma_t}(1+r)^{-2}(\pd_t\psi)^2+(1+r)^{-2}(\pd_r\psi)^2+(1+r)^{-1}|\sla\nabla\psi|^2+r^{-1}(1+r)^{-3}\psi^2 \\
\lesssim \int_{\Sigma_{t_1}}(\pd_t\psi)^2+(\pd_r\psi)^2+|\sla\nabla\psi|^2 + Err_\Box,
\end{multline*}
where
\begin{equation*}
Err_\Box=\int_{t_1}^{t_2}\int_{\Sigma_t}|(2X(\psi)+w\psi)\Box_g\psi|.
\end{equation*}
\end{theorem}
\begin{proof}
Apply Proposition \ref{general_divergence_estimate_prop} to the current
$$J[X_0+C\pd_t,w_0+\epsilon_{\pd_t}w_{\pd_t}]$$
and invoke Corollary \ref{m:Cdt_correction_I_cor}.
\end{proof}

The flux terms in the above theorem have strong weights as $r\rightarrow\infty$. This is due to the fact that we used a large constant $C$ in Corollary \ref{m:Cdt_correction_I_cor}. Had we instead taken $C=1$, then the flux terms would still have a strong weight for $(L\psi)^2$ and $|\sla\nabla\psi|^2$, but would have weak error terms with factors decaying like $r^{-1}$. This actually has an advantage, because it means that the right side of the estimate can be weak. The reason the error terms have a factor decaying like $r^{-1}$ is that the vectorfield $X$ can be rewritten as
$$X=X_0+\pd_t=\left(1-\frac{1}{1+r}\right)\pd_r+\pd_t=L+O(r^{-1})\pd_r,$$
so it is a future causal vectorfield up to a term of order $O(r^{-1})$.

\textbf{It turns out that the error terms with factors decaying like $r^{-1}$ are still too large. Instead it will be necessary to use a vectorfield $X$ which is future causal up to a term of order $O(r^{-2})$. This is the motivation for a second version of the Morawetz estimate.}

\section{The Morawetz estimate, version II}\label{m:morawetzII_sec}

Bearing in mind the discussion at the end of \S\ref{m:morawetz_I_final_subsec}, it is necessary to derive a slightly different Morawetz-type estimate with different asymptotics as $r\rightarrow\infty$ and a careful analysis of the flux terms.

\subsection{The partial Morawetz estimate}\label{m:partial_morawetz_II_sec}

The vectorfield $X_0^{(I)}$ defined in Proposition \ref{m:partial_morawetz_I_prop} can be rewritten as
$$X_0^{(I)}=\left(1-\frac{1}{1+r}\right)\pd_r.$$
In the second version of the Morawetz estimate, we instead use the following vectorfield
$$X_0^{(II)}=\left(1-\frac{1}{(1+r)^2}\right)\pd_r.$$

\begin{proposition}\label{m:partial_morawetz_II_prop}
(Partial Morawetz estimate in Minkowski (version II)) Let
$$X_0=\frac{r^2+2r}{(1+r)^2}\pd_r$$
and
$$w_0=\frac{2r^2+12r+12}{(1+r)^3}.$$
Then the current $J[X_0,w_0]$ from Definition \ref{current_def} satisfies
$$ divJ[X_0,w_0] = \frac{4}{(1+r)^3}(\pd_r\psi)^2+\frac{r+2}{(1+r)^2}|\sla\nabla\psi|^2+\frac{12}{r(1+r)^5}\psi^2.$$
\end{proposition}
\begin{proof}
We proceed as in the proof of version I (Proposition \ref{m:partial_morawetz_I_prop}) up to the choice for $u$. From Lemma \ref{m:divJ_in_terms_of_u_lem}, we have
$$divJ[X,w]=\frac{2u}{r^3}|\sla\nabla\psi|^2 
+2\pd_r\left(\frac{u}{r^2}\right)(\pd_r\psi)^2-\frac12r^{-2}\pd_r(r^2\pd_r(v\pd_ru))\psi^2.$$
Now, we pay careful attention to the asymptotics. In version I, we made the choice for $u^{(I)}$ so that
$$u^{(I)}v=1-\frac{1}{1+r}.$$
This time, we choose $u^{(II)}$ so that
$$u^{(II)}v=1-\frac{1}{(1+r)^2}.$$
That is, we choose
$$u=u^{(II)}=r^2\left(1-\frac{1}{(1+r)^2}\right)=\frac{r^4+2r^3}{(1+r)^2}.$$
\begin{lemma}
The function $u=\frac{r^4+2r^3}{(1+r)^2}$, which is clearly positive, satisfies
$$2\pd_r\left(\frac{u}{r^2}\right)=\frac{4}{(1+r)^3}$$
and
$$-\frac12r^{-2}\pd_r(r^2\pd_r(r^{-2}\pd_ru))=\frac{12}{r(1+r)^5}.$$
\end{lemma}
\begin{proof}
First,
$$2\pd_r\left(\frac{u}{r^2}\right)=2\pd_r\left(1-\frac{1}{(1+r)^2}\right)=\frac{4}{(1+r)^3}.$$
This verifies the first identity.

Then we have the following intermediate calculations, which verify the second identity.
$$\pd_ru=\pd_r\left(r^2\left(1-\frac{1}{(1+r)^2}\right)\right) = 2r\left(1-\frac{1}{(1+r)^2}\right)+r^2\frac{2}{(1+r)^3}$$
$$r^{-2}\pd_ru=\frac{2}{r}\left(1-\frac{1}{(1+r)^2}\right)+\frac{2}{(1+r)^3}$$
$$\pd_r(r^{-2}\pd_ru)=-\frac{2}{r^2}\left(1-\frac{1}{(1+r)^2}\right)+\frac{2}{r}\frac{2}{(1+r)^3}-\frac{6}{(1+r)^4}.$$
$$r^2\pd_r(r^{-2}\pd_ru)=-2+\frac{2}{(1+r)^2}+\frac{4r}{(1+r)^3}-\frac{6r^2}{(1+r)^4}.$$
\begin{multline*}
\pd_r(r^2\pd_r(r^{-2}\pd_ru))=-\frac{4}{(1+r)^3}+\frac{4}{(1+r)^3}-\frac{12r}{(1+r)^4}-\frac{12r}{(1+r)^4}+\frac{24r^2}{(1+r)^5} \\
=\frac{-24r(1+r)+24r^2}{(1+r)^5}=\frac{-24r}{(1+r)^5}
\end{multline*}
$$-\frac12r^{-2}\pd_r(r^2\pd_r(r^{-2}\pd_ru))=\frac{12}{r(1+r)^5}.$$
\end{proof}
Finally, we compute
$$X_0=uv\pd_r=\frac{r^4+2r^3}{(1+r)^2}r^{-2}\pd_r=\frac{r^2+2r}{(1+r)^2}\pd_r$$
\begin{multline*}
w_0=v\pd_ru=r^{-2}\pd_ru=\frac2r\left(1-\frac{1}{(1+r)^2}\right)+\frac{2}{(1+r)^3} \\
=\frac{2(1+r)^3-2(1+r)+2r}{r(1+r)^3}=\frac{2r^3+6r^2+6r}{r(1+r)^3}=\frac{2r^2+6r+6}{(1+r)^3}.
\end{multline*}
This completes the proof of Proposition \ref{m:partial_morawetz_II_prop}.
\end{proof}

\subsection{Remaining issues}

Just as after the proof of Proposition \ref{m:partial_morawetz_I_prop}, there are a still few modifications that must be made to the current $J[X_0,w_0]$ to establish a proper Morawetz estimate. \\
\bp To recover the term $(\pd_t\psi)^2$, we add to $w$ the term $\epsilon_{\pd_t}w_{\pd_t}$. This is very similar to the procedure in \S\ref{m:wdt_correction_I_sec}, except that the function $w_{\pd_t}$ used for vesion II must decay at a faster rate as $r\rightarrow \infty$. (See \S\ref{m:wdt_correction_II_sec}.) \\
\bp We must pay careful attention to the boundary terms. Unlike in the first version of the Morawetz estiamte, where we rather carelessly added a large constant times $\pd_t$, this time, we only add $\pd_t$ and observe that there are only small errors on the boundary. (See \S\ref{m:Cdt_correction_II_sec}.)

Thus, the final current $J$ will be of the form
\begin{align*}
J &= J[X,w] \\
X &= X_0 + \pd_t \\
w &= w_0 + \epsilon_{\pd_t}w_{\pd_t}.
\end{align*}

\subsection{The correction $\epsilon_{\pd_t}w_{\pd_t}$}\label{m:wdt_correction_II_sec}
Just as in \S\ref{m:wdt_correction_I_sec}, we add a function $w_{\pd_t}$, to recover control of $(\pd_t\psi)^2$ at the price of an error term that we can already control. However, this time, the function $w_{\pd_t}$ must decay at a faster rate if we want to control the new error terms with what we already have.
\begin{lemma}\label{m:wdt_correction_II_lem}
The function $w_{\pd_t}=-(1+r)^{-4}$ satisfies
$$(1+r)^{-4}(\pd_t\psi)^2 \lesssim divJ[0,w_{\pd_t}] + (1+r)^{-4}\left[(\pd_r\psi)^2+|\sla\nabla\psi|^2+r^{-1}(1+r)^{-1}\psi^2\right]$$
\end{lemma}
\begin{proof}
By Lemma \ref{divJ_lem}, we have that
\begin{align*}
divJ[0,w_{\pd_t}] &= w_{\pd_t}g^{\alpha\beta}\pd_\alpha\psi\pd_\beta\psi-\frac12\Box_gw_{\pd_t} \psi^2 \\
&= -w_{\pd_t}(\pd_t\psi)^2 +w_{\pd_t}\left[(\pd_r\psi)^2+|\sla\nabla\psi|^2\right]-\frac12\Box_gw_{\pd_t}\psi^2.
\end{align*}
It suffices to choose $w_{\pd_t}$ satisfying the following conditions.
$$(1+r)^{-4}\lesssim -w_{\pd_t},$$
$$|w_{\pd_t}| \lesssim (1+r)^{-4},$$
$$\left|\frac12\Box_gw_{\pd_t}\right| \lesssim r^{-1}(1+r)^{-5}.$$
With the choice $w_{\pd_t}=-(1+r)^{-4}$, the first two conditions are trivially satisfied, and
\begin{align*}
\frac12\Box_gw_{\pd_t} &= -\frac12\frac{1}{r^2}\pd_r(r^2\pd_r(1+r)^{-4}) \\
&= \frac{2}{r^2}\pd_r(r^2(1+r)^{-5}) \\
&= \frac{4}{r}(1+r)^{-5}-5(1+r)^{-4} \\
&\lesssim r^{-1}(1+r)^{-5}.
\end{align*}
Thus, $w_{\pd_t}$ satisfies the remaining condition.
\end{proof}

\begin{corollary}\label{m:wdt_correction_II_cor}
There exists a small $\epsilon_{\pd_t}>0$ so that
\begin{multline*}
(1+r)^{-4}(\pd_t\psi)^2+(1+r)^{-3}(\pd_r\psi)^2+(1+r)^{-1}|\sla\nabla\psi|^2+r^{-1}(1+r)^{-5}\psi^2\\
\lesssim divJ[X_0,w_0+\epsilon_{\pd_t}w_{\pd_t}].
\end{multline*}
\end{corollary}
\begin{proof}
Recall from Proposition \ref{m:partial_morawetz_II_prop} that
$$(1+r)^{-3}(\pd_r\psi)^2+(1+r)^{-1}|\sla\nabla\psi|^2+r^{-1}(1+r)^{-5}\psi^2 \lesssim div J[X_0,w_0].$$
The proof follows from this fact and the inequality in Lemma \ref{m:wdt_correction_II_lem}.
\end{proof}

\subsection{The correction $\pd_t$ and the flux terms}\label{m:Cdt_correction_II_sec}

Recall that
$$X_0+\pd_t=\left(1+\frac{1}{(1+r)^2}\right)\pd_r+\pd_t=L+\frac{1}{(1+r)^2}\pd_r.$$
We first prove a lemma about boundary terms, while approximating $X_0+\pd_t$ by $L$.
\begin{lemma}\label{m:dt_correction_IIa_lem}
$$-J^t[L,w]=(L\psi)^2+|\sla\nabla\psi|^2+r^{-2}\pd_r\left(\frac{r^2}{2}w\right)\psi^2-r^{-2}\pd_r\left(\frac{r^2}{2}w\psi^2\right)+Err,$$
where
$$Err\lesssim (1+r)^{-1}|\psi L\psi|.$$
Furthermore, provided $\epsilon_{\pd_t}$ is sufficiently small, then
$$r^{-1}(1+r)^{-1}\psi^2\approx r^{-2}\pd_r\left(\frac{r^2}{2}w\right)\psi^2.$$
\end{lemma}
\begin{proof}
From Lemma \ref{divJ_lem}, for any $X$ and $w$,
\begin{align*}
J^\mu[X,w] &= 2g^{\mu\lambda}\pd_\lambda\psi X^\nu\pd_\nu\psi-X^\mu\pd^\lambda\psi\pd_\lambda\psi+w\psi g^{\mu\lambda}\pd_\lambda\psi-\frac12\psi^2g^{\mu\lambda}\pd_\lambda w.
\end{align*}
Therefore,
\begin{align*}
-J^t[L] &= 2\pd_t\psi L\psi+\pd^\lambda\psi\pd_\lambda\psi \\
&=2\pd_t\psi (\pd_t\psi+\pd_r\psi)+\left[-(\pd_t\psi)^2+(\pd_r\psi)^2+|\sla\nabla\psi|^2\right] \\
&=(\pd_t\psi)^2+2\pd_t\psi\pd_r\psi+(\pd_r\psi)^2+|\sla\nabla\psi|^2 \\
&= (L\psi)^2+|\sla\nabla\psi|^2.
\end{align*}
and
\begin{align*}
-J^t[0,w]&=w\psi\pd_t\psi \\
&=w\psi L\psi-w\psi\pd_r\psi \\
&=w\psi L\psi - r^{-2}\pd_r\left(\frac{r^2}{2}w\psi^2\right) +r^{-2}\pd_r\left(\frac{r^2}{2}w\right)\psi^2.
\end{align*}
Adding both of these together, we conclude
$$-J^t[L,w]=(L\psi)^2+|\sla\nabla\psi|^2 - r^{-2}\pd_r\left(\frac{r^2}{2}w\psi^2\right) +r^{-2}\pd_r\left(\frac{r^2}{2}w\right)\psi^2 +w\psi L\psi,$$
from which the first part of the lemma follows.

Now,
\begin{multline*}
r^{-2}\pd_r\left(\frac{r^2}{2}w_0\right)=r^{-2}\pd_r\left(\frac{2r^4+12r^3+12r^2}{(1+r)^3}\right) \\
=r^{-2}\frac{(8r^3+36r^2+24r)(1+r)-3(2r^4+12r^3+12r^2)}{(1+r)^4} \\
=\frac{2r^4+8r^3+24r^2+24r}{r^2(1+r)^4},
\end{multline*}
so
$$r^{-1}(1+r)^{-1}\approx r^{-2}\pd_r\left(\frac{r^2}{2}w_0\right).$$
Also,
$$
r^{-2}\pd_r\left(\frac{r^2}{2}\epsilon_{\pd_t}(1+r)^{-4}\right)
= \frac{\epsilon_{\pd_t}}{r}(1+r)^{-4}-4\epsilon_{\pd_t}(1+r)^{-5} \approx \epsilon_{\pd_t}r^{-1}(1+r)^{-4}.
$$
It follows that if $\epsilon_{\pd_t}$ is sufficiently small, then
$$r^{-1}(1+r)^{-1}\psi^2\approx r^{-2}\pd_r\left(\frac{r^2}{2}w\right)\psi^2.$$
\end{proof}

Now, we prove a lemma that estimates the remaining boundary terms.
\begin{lemma}\label{m:dt_correction_IIb_lem}
$$|J^t[X_0+\pd_t-L]|\lesssim (1+r)^{-2}\left[(\pd_t\psi)^2+(\pd_r\psi)^2\right].$$
\end{lemma}
\begin{proof}
Note that $X_0+\pd_t-L=(1+r)^{-2}\pd_r$. From Definition \ref{current_def}, for any $X$,
$$J^\mu[X]=2g^{\mu\lambda}\pd_\lambda\psi X^\nu\pd_\nu\psi-X^\mu\pd^\lambda\psi\pd_\lambda\psi.$$
So
$$J^t[(1+r)^{-2}\pd_r]=-2(1+r)^{-2}\pd_t\psi\pd_r\psi,$$
whence
$$|J^t[(1+r)^{-2}\pd_r]|\lesssim (1+r)^{-2}\left[(\pd_t\psi)^2+(\pd_r\psi)^2\right].$$
\end{proof}

\begin{corollary}\label{m:dt_correction_II_cor}
If $\epsilon_{\pd_t}$ is sufficiently small, then
\begin{multline*}
(1+r)^{-4}(\pd_t\psi)^2+(1+r)^{-3}(\pd_r\psi)^2+(1+r)^{-1}|\sla\nabla\psi|^2+r^{-1}(1+r)^{-5}\psi^2\\
\lesssim divJ[X_0+\pd_t,w_0+\epsilon_{\pd_t}w_{\pd_t}].
\end{multline*}
and
$$\int_{\Sigma_t}-J^t[X_0+\pd_t,w_0+\epsilon_{\pd_t}w_{\pd_t}] \approx \int_{\Sigma_t}(L\psi)^2+|\sla\nabla\psi|^2+r^{-1}(1+r)^{-1}\psi^2 +Err,$$
where
$$Err=\int_{\Sigma_t}(1+r)^{-1}|\psi L\psi|+(1+r)^{-2}\left[(\pd_t\psi)^2+(\pd_r\psi)^2\right].$$
\end{corollary}
\begin{proof}
Note that $divJ[\pd_t]=0$, so the first estimate reduces to Corollary \ref{m:wdt_correction_II_cor}. The second estimate follows from Lemmas \ref{m:dt_correction_IIa_lem} and \ref{m:dt_correction_IIb_lem} and the fact that
$$\int_{\Sigma_t}-\frac1{r^2}\pd_r\left(\frac{r^2}{2}w\psi^2\right)=0,$$
because it is an integral of a total derivative of a quantity that vanishes at both ends.
\end{proof}

\subsection{The Morawetz estimate}

Finally, we conclude with the second version of the Morawetz estimate, which will be used later in this chapter.
\begin{theorem}\label{m:morawetz_II_thm}(Morawetz estimate in Minkowski (version II))
\begin{multline*}
\int_{\Sigma_{t_2}}(L\psi)^2+|\sla\nabla\psi|^2+r^{-1}(1+r)^{-1}\psi^2 \\
+\int_{t_1}^{t_2}\int_{\Sigma_t}(1+r)^{-4}(\pd_t\psi)^2+(1+r)^{-3}(\pd_r\psi)^2+(1+r)^{-1}|\sla\nabla\psi|^2+r^{-1}(1+r)^{-5}\psi^2 \\
\lesssim \int_{\Sigma_{t_1}}(L\psi)^2+|\sla\nabla\psi|^2+r^{-1}(1+r)^{-1}\psi^2 + Err,
\end{multline*}
where
\begin{align*}
Err &= Err_1+Err_2+Err_\Box \\
Err_1 &= \int_{\Sigma_{t_2}}(1+r)^{-1}|\psi L\psi| \\
Err_2 &= \int_{\Sigma_{t_2}}(1+r)^{-2}\left[(\pd_r\psi)^2+(\pd_t\psi)^2\right] \\
Err_\Box &=\int_{t_1}^{t_2}\int_{\Sigma_t}|(2X(\psi)+w\psi)\Box_g\psi|.
\end{align*}
\end{theorem}
\begin{proof}
Apply Proposition \ref{general_divergence_estimate_prop} to the current
$$J[X_0+C\pd_t,w_0+\epsilon_{\pd_t}w_{\pd_t}]$$
and invoke Corollary \ref{m:dt_correction_II_cor}.
\end{proof}

\section{The $r^p$ estimate}\label{m:rp_sec}

\subsection{The approach of \cite{dafermos2010new}}\label{m:daf_rod_sec}

The $r^p$ estimate is most simply presented using the approach of Dafermos and Rodnianski. \cite{dafermos2010new} To provide context, we discuss this approach first in its simplest form. \textbf{The reader should keep in mind that no notation defined here extends to the rest of this thesis.}

We fix a constant $v_0$ and pick an incoming cone $\ul{C}_{v_0}=\{t+r=v_0\}$. We define a foliation of outgoing cones $C_u=\{t-r=u\}\cap\{t+r>v_0\}$. We require that  $u<v_0$ so that $C_u$ extends from $\ul{C}_{v_0}$. We also define a spacetime domain $D_{u_1}^{u_2}=\cup_{u\in [u_1,u_2]}C_u$. Finally, we define $I_{u_1}^{u_2}$ to be the portion of future null infinity that borders $D_{u_1}^{u_2}$ in a compactified Penrose diagram.

\begin{proposition}
If $\Box_g\psi=0$, then
\begin{multline*}
\int_{C_{u_2}}\frac{r^p}{4}(L\psi+r^{-1}\psi)^2+\int_{I_{u_1}^{u_2}}r^p|\sla\nabla\psi|^2 \\
+\int\int_{D_{u_1}^{u_2}}r^{p-1}\left[\frac{p}{4}(L\psi+r^{-1}\psi)^2+(2-p)|\sla\nabla\psi|^2\right] \\
= \int_{C_{u_1}}\frac{r^p}{4}(L\psi+r^{-1}\psi)^2 +\int_{\ul{C}_{v_0}}r^p|\sla\nabla\psi|^2.
\end{multline*}
\end{proposition}

\begin{remark}
Due to the factors $p$ and $2-p$ showing up in the bulk integral, this identity is only useful if both $0\le p$ and $0\le 2-p$, equivalently, $0\le p\le 2$. If we wish to have control of both $(L\psi+r^{-1}\psi)^2$ and $|\sla\nabla\psi|^2$, then we actually need $0<p<2$. Finally, if we want uniform control, then we must choose two small constants $\delm,\delp>0$ and require that $\delm\le p\le 2-\delp$.
\end{remark}

\begin{proof}
We start by defining a new coordinate system
\begin{align*}
u&=t-r \\
v&=t+r.
\end{align*}
One can check that
$$\pd_u=\frac12(\pd_t-\pd_r)$$
$$\pd_v=\frac12(\pd_t+\pd_r)$$
and
$$\sqrt{-g}=\frac{r^2}{4}d\omega dvdu.$$

Note that if $\Psi=r\psi$, and $\Box_g\psi=0$, then
\begin{align*}
0&=-\pd_t^2\Psi+\pd_r^2\Psi+\sla\triangle\Psi \\
&=-\pd_u\pd_v\Psi+\sla\triangle\Psi.
\end{align*}
So
$$0=\int\int_{D_{u_1}^{u_2}}8r^{p-2}\pd_v\Psi\left[-\pd_u\pd_v\Psi+\sla\triangle\Psi\right].$$
From this fact, we prove the proposition by integrating by parts. While integrating by parts, we include the part of the volume form depending on $r$ explicitly to simplify the calculation.

On one hand,
\begin{align*}
\int\int_{D_{u_1}^{u_2}}8r^{p-2}\pd_v\Psi (-2\pd_u\pd_v\Psi) &= \int_{u_1}^{u_2}\int_{v_0}^\infty\int_{S^2}-2r^{p-2}\pd_v\Psi \pd_u\pd_v\Psi r^2d\omega dvdu \\
&= \int_{u_1}^{u_2}\int_{v_0}^\infty\int_{S^2}-r^{p} \pd_u((\pd_v\Psi)^2) d\omega dvdu \\
&= \left.-\int_{v_0}^\infty\int_{S^2}r^p(\pd_v\Psi)^2d\omega dv \right|_{u_1}^{u_2} \\
&\hspace{1.3in}+\int_{u_1}^{u_2}\int_{v_0}^\infty\int_{S^2}pr^{p-1} (\pd_ur)(\pd_v\Psi)^2 d\omega dvdu \\
&= \int_{C_{u_1}}2r^{p-2}(\pd_v\Psi)^2 -\int_{C_{u_2}}2r^{p-2}(\pd_v\Psi)^2 +\int_{D_{u_1}^{u_2}}-2pr^{p-3}(\pd_v\Psi)^2.
\end{align*}
On the other hand, using the notation $\Omega=r\sla\nabla$ and observing that $[\Omega,\pd_v]=0$,
\begin{align*}
\int\int_{D_{u_1}^{u_2}}8r^{p-2}\pd_v\Psi\sla\triangle\psi &= \int_{u_1}^{u_2}\int_{v_0}^\infty\int_{S^2}2r^{p-2}\pd_v\Psi \sla\triangle\Psi r^2d\omega dvdu \\
&= \int_{u_1}^{u_2}\int_{v_0}^\infty\int_{S^2}2r^{p-2}\pd_v\Psi \Omega\cdot\Omega\Psi d\omega dvdu  \\
&= -\int_{u_1}^{u_2}\int_{v_0}^\infty\int_{S^2}2r^{p-2}\pd_v\Omega\Psi \cdot\Omega\Psi d\omega dvdu  \\
&= -\int_{u_1}^{u_2}\int_{v_0}^\infty\int_{S^2}r^{p-2}\pd_v|\Omega\Psi|^2 d\omega dvdu \\
&= \left.-\int_{u_1}^{u_2}\int_{S^2}r^{p-2}|\Omega\Psi|^2 d\omega du\right|_{v_0}^\infty \\
&\hspace{1in} +\int_{u_1}^{u_2}\int_{v_0}^\infty\int_{S^2}(p-2)r^{p-3}(\pd_v r)|\Omega\Psi|^2 d\omega dvdu  \\
&= \int_{\ul{C}_{v_0}}2r^{p-2}|\sla\nabla\Psi|^2-\int_{I_{u_1}^{u_2}}2r^{p-2}|\sla\nabla\Psi|^2 +\int_{D_{u_1}^{u_2}}2(p-2)r^{p-3}|\sla\nabla\Psi|^2.
\end{align*}
Combining both calculations, we obtain
\begin{multline*}
0=\int_{C_{u_1}}2r^{p-2}(\pd_v\Psi)^2 -\int_{C_{u_2}}2r^{p-2}(\pd_v\Psi)^2 +\int\int_{D_{u_1}^{u_2}}-2pr^{p-3}(\pd_v\Psi)^2 \\
+\int_{\ul{C}_{v_0}}2r^{p-2}|\sla\nabla\Psi|^2-\int_{I_{u_1}^{u_2}}2r^{p-2}|\sla\nabla\Psi|^2 +\int_{D_{u_1}^{u_2}}2(p-2)r^{p-3}|\sla\nabla\Psi|^2,
\end{multline*}
or equivalently,
\begin{multline*}
\int_{C_{u_2}}r^{p-2}(\pd_v\Psi)^2+\int_{I_{u_1}^{u_2}}r^{p-2}|\sla\nabla\Psi|^2+\int\int_{D_{u_1}^{u_2}}pr^{p-3}(\pd_v\Psi)^2+(2-p)r^{p-3}|\sla\nabla\Psi|^2 \\
= \int_{C_{u_1}}r^{p-2}(\pd_v\Psi)^2 +\int_{\ul{C}_{v_0}}r^{p-2}|\sla\nabla\Psi|^2.
\end{multline*}
Finally, observe that $\pd_v\Psi=\frac12L(r\psi)=\frac{r}2(L\psi+r^{-1}\psi)$, so
\begin{multline*}
\int_{C_{u_2}}\frac{r^p}{4}(L\psi+r^{-1}\psi)^2+\int_{I_{u_1}^{u_2}}r^p|\sla\nabla\psi|^2 \\
+\int\int_{D_{u_1}^{u_2}}r^{p-1}\left[\frac{p}{4}(L\psi+r^{-1}\psi)^2+(2-p)|\sla\nabla\psi|^2\right] \\
= \int_{C_{u_1}}\frac{r^p}{4}(L\psi+r^{-1}\psi)^2 +\int_{\ul{C}_{v_0}}r^p|\sla\nabla\psi|^2.
\end{multline*}
\end{proof}

\subsection{The pre-$r^p$ identity}

We derive the following identity, which will be used for the proof of Proposition \ref{m:incomplete_p_estimate_prop}.
\begin{lemma}\label{m:p_ee_identity_lem}
For any function $f=f(r)$, the following identity holds.
\begin{align*}
&\int_{\Sigma_{t_2}}\left[f(L\psi+r^{-1}\psi)^2+f|\sla\nabla\psi|^2+\epsilon r^{-1}f'\psi^2-r^{-2}\pd_r(r^{2}r^{-1}f\psi^2)\right] \\
&\hspace{.2in}+\int_{t_1}^{t_2}\int_{\Sigma_t}\left[(2r^{-1}f-f')|\sla\nabla\psi|^2+ f'\left(L\psi+\frac{1-\epsilon}{r}\psi\right)^2+\epsilon((1-\epsilon)f'-rf'')r^{-2}\psi^2\right] \\
&=\int_{\Sigma_{t_1}}\left[f|\sla\nabla\psi|^2+f(L\psi+r^{-1}\psi)^2+\epsilon r^{-1}f'\psi^2-r^{-2}\pd_r(r^{2}r^{-1}f\psi^2)\right] \\
&\hspace{.2in}+\int_{t_1}^{t_2}\int_{\Sigma_t}-(2fL\psi+2r^{-1}f\psi)\Box_g\psi.
\end{align*}
\end{lemma}

\begin{proof}
We will use Proposition \ref{general_divergence_estimate_prop} together with the following current template.
$$J[X,w,m]_\mu = T_{\mu\nu} X^\nu +w\psi\pd_\mu\psi-\frac12\psi^2\pd_\mu w+m_\mu \psi^2,$$
$$T_{\mu\nu}=2\pd_\mu\psi\pd_\nu\psi-g_{\mu\nu}\pd^\lambda\psi\pd_\lambda\psi.$$

Assume for now that $\Box_g\psi=0$.

\begin{lemma}\label{m:divJphiX_lem}
$$divJ[fL]=f'(L\psi)^2-2r^{-1}f\left((\pd_r\psi)^2-(\pd_t\psi)^2\right)-f'|\sla\nabla\psi|^2.$$
\end{lemma}
\begin{proof}
Note that
$$div J[X] = K^{\mu\nu}\pd_\mu\psi\pd_\nu\psi,$$
where
$$K^{\mu\nu}=2g^{\mu\lambda}\pd_\lambda X^\nu-X^\lambda \pd_\lambda(g^{\mu\nu})-div X g^{\mu\nu}.$$

Set $X=fL=f(\pd_r+\pd_t)$. From the above formula,
$$K^{tr}+K^{rt}=2g^{rr}\pd_r X^t=2f'.$$
Thus, the expression for $divJ[fL]$ will have a mixed term of the form 
$$2f'\pd_r\psi\pd_t\psi.$$
Note that
\begin{align*}
f'(L\psi)^2 &= f'(\pd_r\psi+\pd_t\psi)^2 \\
&= f'(\pd_r\psi)^2+2f'\pd_r\psi\pd_t\psi + f'(\pd_t\psi)^2.
\end{align*}
We now compute the $(\pd_r\psi)^2$ and $(\pd_t\psi)^2$ components, subtracting the part that will be grouped with the $(L\psi)^2$ term.
\begin{align*}
K^{rr}-f' &= \left[2g^{rr}\pd_rX^r-X^r\pd_r g^{rr}-r^{-2}\pd_r(r^2X^r)g^{rr}\right]-f' \\
&= 2\pd_r f -r^{-2}\pd_r\left(r^2 f\right)-f' \\
&= -2r^{-1}  f
\end{align*}
and
\begin{align*}
K^{tt}-f' &= \left[-X^r\pd_r g^{tt}-r^{-2}\pd_r(r^2X^r)g^{tt}\right]-f' \\
&= r^{-2}\pd_r(r^2f) -f' \\
&= 2r^{-1}f.
\end{align*}
Finally,
\begin{align*}
\sla{K}^{\alpha\beta} &= -X^r\pd_r \sla{g}^{\alpha\beta}-r^{-2}\pd_r\left(r^2X^r\right)\sla{g}^{\alpha\beta} \\
&= -r^{-2}\pd_r\left(r^2\sla{g}^{\alpha\beta}X^r\right) \\
&= -r^{-2}\pd_r((r^2\sla{g}^{\alpha\beta})f) \\
&= -f' \sla{g}^{\alpha\beta}.
\end{align*}
(We used in the last line that $\pd_r(r^2\sla{g}^{\alpha\beta})=0$.) Combining all these terms gives the identity stated in the lemma. 
\end{proof}

Next, we choose $w=2r^{-1}f$ to directly cancel the middle term in the above lemma.
\begin{lemma}
$$divJ[fL,2r^{-1}f]=f'(L\psi)^2+(2r^{-1}f-f')|\sla\nabla\psi|^2-\frac12\Box_g(2r^{-1}f)\psi^2.$$
\end{lemma}
\begin{proof}
Note that
$$divJ[0,w]=wg^{\mu\nu}\pd_\mu\psi\pd_\nu\psi-\frac12\Box_gw \psi^2.$$
We compute the new terms only.
\begin{align*}
divJ\left[0,2r^{-1}f\right] &= 2r^{-1}fg^{\alpha\beta}\pd_\alpha\psi\pd_\beta\psi -\frac12\Box_g\left(2r^{-1}f\right)\psi^2 \\
&= 2r^{-1}f\left((\pd_r\psi)^2-(\pd_t\psi)^2\right) +2r^{-1}f|\sla\nabla\psi| -\frac12\Box_g\left(2r^{-1}f\right)\psi^2.
\end{align*}
When adding these terms to the expression in Lemma \ref{m:divJphiX_lem}, the $(\pd_r\psi)^2-(\pd_t\psi)^2$ terms cancel (this was the reason for the choice of $w=2r^{-1}f$) and the result is as desired. 
\end{proof}

Note that
\begin{align*}
-\frac12\Box_g(2r^{-1}f) &= -\frac12r^{-2}\pd_r(r^2\pd_r(2r^{-1}f)) \\
&= r^{-2}\pd_r(f-rf') \\
&= r^{-2}(f'-f'-rf'') \\
&= -r^{-1}f''.
\end{align*}
Later, when $f\sim r^p$, this will have a sign $-p(p-1)$. The sign will be negative if $p>1$, which is bad. So we include a divergence term to fix it. (But in doing so, we almost lose some other good terms--this is why we need a small parameter $\epsilon$.) This is the point of the following lemma.

\begin{lemma}
\begin{multline*}
f'(L\psi)^2-\frac12\Box_g(2r^{-1}f)\psi^2+(1-\epsilon)div(\psi^2 r^{-1}f'L) \\
=f' \left(L\psi+\frac{1-\epsilon}{r}\psi\right)^2 +\epsilon  \left((1-\epsilon)f'-rf''\right)r^{-2}\psi^2.
\end{multline*}
\end{lemma}
\begin{proof}
We just recently determined that
$$-\frac12\Box_g(2r^{-1}f)=-r^{-1}f''.$$
We also calculate
\begin{align*}
div\left(\psi^2r^{-1}f'L\right) &= r^{-2}\pd_\alpha\left(\psi^2rf' L^\alpha\right) \\
&=r^{-1}f'2\psi L\psi +r^{-2}\pd_r(rf')\psi^2 \\
&=r^{-1} f' 2\psi L\psi +r^{-2}f'\psi^2 +r^{-1} f''\psi^2.
\end{align*}
The first two terms in the last line complete a square with the term $f' (L\psi)^2$. The third term cancels with the first term from the previous calculation. However, it will be beneficial to introduce the factor $1-\epsilon$ that appears in the lemma, so that a good term appears with an $\epsilon$ factor. This is summarized by the following two calculations.
\begin{multline*}
f'(L\psi)^2+(1-\epsilon)\left(r^{-1}f'2\psi L\psi+ r^{-2} f'\psi^2\right) \\
= f' \left(L\psi+\frac{1-\epsilon}{r}\psi\right)^2 -(1-\epsilon)^2r^{-2} f'\psi^2 +(1-\epsilon)r^{-2} f'\psi^2 \\
=  f' \left(L\psi+\frac{1-\epsilon}{r}\psi\right)^2 +\epsilon (1-\epsilon)r^{-2} f'\psi^2 
\end{multline*}
and
$$- r^{-1} f''\psi^2 +(1-\epsilon) r^{-1} f''\psi^2=-\epsilon r^{-1} f''\psi^2.$$
Adding these terms together yields
$$ f' \left(L\psi+\frac{1-\epsilon}{r}\psi\right)^2 +\epsilon  \left((1-\epsilon)f'-rf''\right)r^{-2}\psi^2.$$
\end{proof}

Thus, we have shown that if $\Box_g\psi=0$, then
\begin{multline*}
divJ[fL,2r^{-1}f, (1-\epsilon)r^{-1}f'L] \\
=(2r^{-1}f-f')|\sla\nabla\psi|^2+ f' \left(L\psi+\frac{1-\epsilon}r\psi\right)^2+\epsilon\left((1-\epsilon)f'-rf''\right)r^{-2}\psi^2.
\end{multline*}
If we remove the assumption $\Box_g\psi=0$, there is an additional term
$$(2X(\psi)+w\psi)\Box_g\psi =(2fL\psi+2r^{-1}f\psi)\Box_g\psi.$$

Finally, we turn to the boundary terms.
\begin{lemma}
\begin{multline*}
-J^t\left[fL,2r^{-1}f,(1-\epsilon)r^{-1}f'L\right] \\
=f|\sla\nabla\psi|^2+f(L\psi+r^{-1}\psi)^2+\epsilon r^{-1}f'\psi^2-r^{-2}\pd_r(rf\psi^2)
\end{multline*}
\end{lemma}
\begin{proof}
We have
\begin{align*}
-J^t[fL] &= -2\pd^t\psi fL\psi+fL^t\pd^\lambda\psi\pd_\lambda\psi \\
&=-2 fg^{tt}\pd_t\psi L\psi +f\left(g^{tt}(\pd_t\psi)^2+g^{rr}(\pd_r\psi)^2+|\sla\nabla\psi|^2\right) \\
&=-fg^{tt}(\pd_t\psi)^2-2fg^{tt}\pd_t\psi\pd_r\psi+fg^{rr}(\pd_r\psi)^2 +f|\sla\nabla\psi|^2 \\
&=f(\pd_t\psi)^2+2f\pd_t\psi\pd_r\psi +f(\pd_r\psi)^2+f|\sla\nabla\psi|^2 \\
&=f(L\psi)^2+f|\sla\nabla\psi|^2.
\end{align*}
Also,
\begin{align*}
-J^t_{(\psi)}\left[0,2r^{-1}f\right] &= -2r^{-1}f\psi\pd^t\psi \\
&= -2r^{-1}fg^{tt}\psi\pd_t\psi \\
&= 2r^{-1}f\psi\pd_t\psi \\
&= 2r^{-1}f\psi L\psi -2r^{-1}f\psi\pd_r\psi \\
&= 2r^{-1}f\psi L\psi -r^{-2}\pd_r(rf\psi^2)+(r^{-2}f+r^{-1}f')\psi^2\\
&= \left(2r^{-1}f\psi L\psi+r^{-2}f\psi^2\right) +r^{-1}f'\psi^2-r^{-2}\pd_r(rf\psi^2).
\end{align*}
Now, observe that
$$f(L\psi)^2+2r^{-1}f\psi L\psi+r^{-2}f\psi^2 = f\left(L\psi+r^{-1}\psi\right)^2.$$
Thus,
$$-J^t\left[fL,2r^{-1}f\right] 
=f\left(L\psi+r^{-1}\psi\right)^2 
+f|\sla\nabla\psi|^2+r^{-1}f'\psi^2-r^{-2}\pd_r(rf\psi^2).$$
Also,
$$-J^t\left[0,0,(1-\epsilon)r^{-1}f'L\right] = -(1-\epsilon)r^{-1}f'\psi^2 L^t =-(1-\epsilon)r^{-1}f'\psi^2.$$
Adding these two expressions together yields the result. 
\end{proof}

This completes the proof of Lemma \ref{m:p_ee_identity_lem}.
\end{proof}

\subsection{The incomplete $r^p$ estimate near $i^0$}

Now we use Lemma \ref{m:p_ee_identity_lem} and make a particular choice for $f$ (so that $f=r^p$ for large $r$) to prove the following proposition.

\begin{proposition}\label{m:incomplete_p_estimate_prop}
Fix $\delm,\delp>0$. Let $R$ be a large radius. Then for all $p\in [\delm,2-\delp]$, the following estimate holds if $\psi$ decays sufficiently fast as $r\rightarrow\infty$.
\begin{multline*}
\int_{\Sigma_{t_2}\cap\{r>2R\}}r^p\left[(L\psi)^2+|\sla\nabla\psi|^2+r^{-2}\psi^2\right] \\
+\int_{t_1}^{t_2}\int_{\Sigma_t\cap\{r>2R\}}r^{p-1}\left[(L\psi)^2+|\sla\nabla\psi|^2+r^{-2}\psi^2\right] \\
\lesssim \int_{\Sigma_{t_1}\cap\{r>2R\}}r^p\left[(L\psi)^2+|\sla\nabla\psi|^2+r^{-2}\psi^2\right] +Err,
\end{multline*}
where
\begin{align*}
Err &= Err_1+Err_\Box \\
Err_1 &= \int_{\Sigma_t\cap\{R<r<2R\}} (L\psi)^2+|\sla\nabla\psi|^2+\psi^2 \\
Err_\Box &= \int_{t_1}^{t_2}\int_{\Sigma_t\cap\{R<r\}}r^p(|L\psi|+r^{-1}|\psi|)|\Box_g\psi|.
\end{align*}
\end{proposition}
\begin{proof}
The estimate follows from the identity given in Lemma \ref{m:p_ee_identity_lem} and a particular choice for the function $f$.
$$f(r)=\rho^p,$$
where
$$
\rho = \left\{
\begin{array}{ll}
0 & r\le R \\
smooth, \rho'>0 & r\in[R,2R] \\
r & 2R< r.
\end{array}
\right.
$$
With this choice, we have 
$$f\ge 0$$
$$f'\ge 0$$
and for $r>2R$,
$$f = r^p$$
$$f'=pr^{p-1}.$$
Furthermore, for $r>2R$,
$$2r^{-1}f-f' = 2r^{p-1}-pr^{p-1} = (2-p)r^{p-1}.$$
It follows that if $p\le 2-\delp$, then for $r>2R$,
$$r^{p-1}\lesssim 2r^{-1}f-f'.$$
Also, for $r>2R$,
$$\epsilon ((1-\epsilon)f'-rf'') = \epsilon ((1-\epsilon)pr^{p-1}-p(p-1)r^{p-1}) = \epsilon p (2-\epsilon -p)r^{p-1}.$$
If $\delm\le p\le 2-\delp$ and $\epsilon\le \delp/2$, then
$$\epsilon r^{p-1} \lesssim \epsilon ((1-\epsilon)f'-rf'').$$

Finally, we observe that if $\psi$ vanishes sufficiently fast as $r\rightarrow\infty$, then
$$\int_{\Sigma_t}-r^{-2}\pd_r(rf\psi^2) =0.$$
With these facts having been established, it is straightforward to check that the estimate follows from Lemma \ref{m:p_ee_identity_lem}.
\end{proof}

\subsection{The $r^p$ estimate} 

We conclude this section by proving the $r^p$ estimate. This is a combination of the $h\pd_t$ estimate (Proposition \ref{m:hdt_prop}), the Morawetz estimate (Theorem \ref{m:morawetz_II_thm}), and the incomplete $r^p$ estimate near $i^0$ (Proposition \ref{m:incomplete_p_estimate_prop}).

\begin{proposition}\label{m:rp_prop}
Fix $\delm,\delp>0$ and let $p\in[\delm,2-\delp]$. Then if
$\psi$ decays sufficiently fast as $r\rightarrow\infty$, the following estimate holds.
\begin{multline*}
\int_{\Sigma_{t_2}}(1+r)^p\left[(L\psi)^2+|\sla\nabla\psi|^2+r^{-1}(1+r)^{-1}\psi^2 + (1+r)^{-2}(\pd_r\psi)^2\right] \\
+ \int_{t_1}^{t_2}\int_{\Sigma_t}(1+r)^{p-1}\left[(L\psi)^2+|\sla\nabla\psi|^2+r^{-1}(1+r)^{-1}\psi^2+(1+r)^{-2}(\pd_r\psi)^2\right] \\
\lesssim \int_{\Sigma_{t_1}}(1+r)^p\left[(L\psi)^2+|\sla\nabla\psi|^2+r^{-1}(1+r)^{-1}\psi^2+(1+r)^{-2}(\pd_r\psi)^2\right] 
+Err_\Box,
\end{multline*}
where
\begin{align*}
Err_\Box &= \int_{t_1}^{t_2}\int_{\Sigma_t}|(2X(\psi)+w\psi)\Box_g\psi|,
\end{align*}
and the vectorfield $X$ and function $w$ satisfy the following properties. \\
\bp $X$ is everywhere timelike, but asymptotically null at the rate $X=O(r^p)L+O(r^{p-2})\pd_t$. \\
\bp $w =O(r^{p-1})$ for large $r$.
\end{proposition}
\begin{proof}
We start with the Morawetz estimate (Theorem \ref{m:morawetz_II_thm}) and add a small constant times the incomplete $r^p$ estimates (Proposition \ref{m:incomplete_p_estimate_prop}). The small constant can be chosen so that the bulk error term $Err_1$ from Proposition \ref{m:incomplete_p_estimate_prop} can be absorbed into the bulk in the Morawetz estimate. The result is the following estimate.
\begin{multline*}
\int_{\Sigma_{t_2}}(1+r)^p\left[(L\psi)^2+|\sla\nabla\psi|^2+r^{-1}(1+r)^{-1}\psi^2\right]+(1+r)^{-2}(\pd_r\psi)^2 \\
\hspace{1in}+\int_{t_1}^{t_2}\int_{\Sigma_t} (1+r)^{p-1}\left[(L\psi)^2+|\sla\nabla\psi|^2+r^{-1}(1+r)^{-1}\psi^2\right]+(1+r)^{-3}(\pd_r\psi)^2 \\
\lesssim \int_{\Sigma_{t_1}}(1+r)^p\left[(L\psi)^2+|\sla\nabla\psi|^2+r^{-1}(1+r)^{-1}\psi^2\right]+(1+r)^{-2}(\pd_r\psi)^2 + Err'
\end{multline*}
where
\begin{align*}
Err' &= Err'_1+Err'_2+Err'_\Box \\
Err'_1 &= \int_{\Sigma_{t_2}}(1+r)^{-1}|\psi L\psi| \\
Err'_2 &= \int_{\Sigma_{t_2}}(1+r)^{-2}\left[(\pd_r\psi)^2+(\pd_t\psi)^2\right] \\
Err'_\Box &= \int_{t_1}^{t_2}\int_{\Sigma_t\cap\{R<r\}}(1+r)^p(|L\psi|+r^{-1}|\psi|)|\Box_g\psi| \\
&\hspace{1in} +\int_{t_1}^{t_2}\int_{\Sigma_t}|(2X'(\psi)+w'\psi)|\Box_{g}\psi|,
\end{align*}
and $X'$ and $w'$ are the vectorfield and function defined in the Morawetz estimate (Theorem \ref{m:morawetz_II_thm}).

The error term $Err'_1$ can in fact be removed due to the following argument.
\begin{align*}
Err'_1 &\lesssim \int_{\Sigma_{t_2}}\epsilon (1+r)^p(L\psi)^2+\epsilon^{-1}(1+r)^{-p}(1+r)^{-2}\psi^2 \\
&\lesssim \int_{\Sigma_{t_2}}\epsilon (1+r)^p[(L\psi)^2+(1+r)^{-2}\psi^2] +\int_{\Sigma_{t_2}\cap\{r\le R_\epsilon\}}\epsilon^{-1}\psi^2 .
\end{align*}
The radius $R_\epsilon$ should be chosen sufficiently large so that $\epsilon^{-1}(1+r)^{-p}\le \epsilon (1+r)^p$ whenever $r>R_\epsilon$. This critically depends on the fact that $p\ge\delm>0$. Now, the parameter $\epsilon$ can be taken sufficiently small so as to absorb the first two terms into the left side of the main estimate and the last two terms can be included with the term $Err'_2$ after applying a Hardy estimate.

We return to the main estimate. Notice that most terms have improved weights near $i^0$ and a few error terms remain on $\Sigma_{t_2}$. The next step is to use the $h\pd_t$ estimate (Proposition \ref{m:hdt_prop}) to eliminate these error terms and improve the weights near $i^0$ for the remaining $\pd_r\psi$ terms. The result is the following estimate.
\begin{multline*}
\int_{\Sigma_{t_2}}(1+r)^p\left[(L\psi)^2+|\sla\nabla\psi|^2+r^{-1}(1+r)^{-1}\psi^2+(1+r)^{-2}(\pd_r\psi)^2\right] \\
+\int_{t_1}^{t_2}\int_{\Sigma_t} (1+r)^{p-1}\left[(L\psi)^2+|\sla\nabla\psi|^2+r^{-1}(1+r)^{-1}\psi^2+c_\epsilon (1+r)^{-2}(\pd_r\psi)^2\right] \\
\lesssim \int_{\Sigma_{t_1}}(1+r)^p\left[(L\psi)^2+|\sla\nabla\psi|^2+r^{-1}(1+r)^{-1}\psi^2+(1+r)^2(\pd_r\psi)^2\right] +Err'',
\end{multline*}
where
\begin{align*}
Err'' &= Err''_1+Err''_\Box \\
Err''_1 &= \int_{t_1}^{t_2}\int_{\Sigma_t\cap\{R<r\}} \epsilon (1+r)^{-1}(L\psi)^2 \\
Err''_\Box &= \int_{t_1}^{t_2}\int_{\Sigma_t\cap\{R<r\}}(1+r)^p(|L\psi|+(1+r)^{-1}|\psi|)|\Box_g\psi| \\
&\hspace{.5in}+ \int_{t_1}^{t_2}\int_{\Sigma_t}|(2X'(\psi)+w'\psi)\Box_g\psi|
+ \int_{t_1}^{t_2}\int_{\Sigma_t}C_\epsilon (1+r)^{p-2}|\pd_t\psi\Box_g\psi|.
\end{align*}
By taking $\epsilon$ sufficiently small, since $p>0$, the error term $Err''_1$ can be absorbed into the left side. 

The resulting vectorfield $X$ and function $w$ can be understood by combining the vectorfields and scalar functions used to construct the three estimates that were used in this proof.
\end{proof}

\section{Aside: Special estimates in Minkowski space}\label{m:special_sec}

The $r^p$ estimate was quite complicated to construct. In Minkowski space, there are other special estimates that give similar information and can be constructed using comparatively simple vectorfield multipliers. In this standalone section we derive two estimates, which we call ``$p=1$ type'' and ``$p=2$ type,'' because they have similar asymptotic behavior to the $r^p$ estimate with $p=1$ and $p=2$ respectively. (In a way, the energy estimate itself could be considered as a special ``$p=0$ type'' estimate, thus completing the full range from $p=0$ to $p=2$.) Both of the estimates presented in this section are based on conformal vectorfields in Minkowski space, so they do not generalize as well to black hole spacetimes. (However, that does not mean they cannot be generalized at all--for example, see \cite{luk2010null}.)

This is a standalone section. The results established in this section will not be used later.

\subsection{A $p=1$ type estimate}

The $p=1$ type estimate is based on the scaling vectorfield $S$ defined in the following lemma.
\begin{lemma}\label{S_is_conformal_lem}
The scaling vectorfield 
$$S=r\pd_r+t\pd_t$$
 is a conformal Killing vectorfield. More precisely,
$$\mathcal{L}_S g_{\mu\nu} = 2g_{\mu\nu}.$$
\end{lemma}
\begin{proof}
From equation (\ref{lie_derivative_of_metric_eqn}),
$$\mathcal{L}_Sg_{\mu\nu}=2\nabla_{(\mu}S_{\nu)}.$$
It suffices to show that
$$2\nabla^{(\mu}S^{\nu)} = 2g^{\mu\nu}.$$
From Lemma \ref{divJ_lem},
\begin{align*}
2\nabla^{(\mu}S^{\nu)} = g^{\mu\lambda}\pd_\lambda S^\nu+g^{\nu\lambda}\pd_\lambda S^\mu-S^\lambda\pd_\lambda g^{\mu\nu}.
\end{align*}
Observe that
$$\pd_r \sla{g}^{\mu\nu}=\pd_r(r^{-2}r^2\sla{g}^{\mu\nu})=(r^2\sla{g}^{\mu\nu})\pd_r r^{-2} =-\frac{2}r\sla{g}^{\mu\nu},$$
while the remaining metric components are constant with respect to $r$, and all metric components are constant with respect to $t$. Thus,
$$-S^\lambda\pd_\lambda g^{\mu\nu} = -S^r\pd_rg^{\mu\nu} = -r\left(-\frac{2}{r}\sla{g}^{\mu\nu}\right) = 2\sla{g}^{\mu\nu}.$$
Observe also that $\pd_r S^r=1$ and $\pd_t S^t=1$, while all other components of $\pd_\lambda S^\mu$ are zero. Thus,
$$2\nabla^{(\mu}S^{\nu)} = 2g^{\mu\nu}.$$
\end{proof}

\begin{lemma}\label{m:special_S_2_lem}
If $\Box_g\psi=0$, then the momentum $J[S,2]$ satisfies
$$-J^t[S,2]=r(L\psi+r^{-1}\psi)^2+(t-r)\left[(\pd_t\psi)^2+(\pd_r\psi)^2\right]+t|\sla\nabla\psi|^2+r^{-1}\psi^2-r^{-2}\pd_r(r^2\psi^2).$$
and
$$div J[S,2] = 0.$$
This is useful if $r<t$.
\end{lemma}
\begin{proof}
The fact that
$$div J[S,2]=0$$
follows immediately from Lemma \ref{S_is_conformal_lem} and Corollary \ref{conformal_current_cor}, since Minkowski space is Ricci flat.

Now, we compute $-J^t[S,2]$. Recall that
$$J^\mu[X,w] = 2\pd^\mu \psi X(\psi)-X^\mu \pd^\lambda\psi\pd_\lambda\psi +w\psi\pd^\mu \psi-\frac12\psi^2\pd^\mu w.$$
Thus,
\begin{align*}
-J^t[S,2] &= -2g^{tt}\pd_t\psi S(\psi)+S^t \pd^\lambda\psi\pd_\lambda\psi-2g^{tt}\psi\pd_t\psi \\
&= 2t(\pd_t\psi)^2+2r\pd_t\psi\pd_r\psi+t\left[-(\pd_t\psi)^2+(\pd_r\psi)^2+|\sla\nabla\psi|^2\right]+2\psi\pd_t\psi \\
&= t(\pd_t\psi)^2+2r\pd_t\psi\pd_r\psi+t(\pd_r\psi)^2+t|\sla\nabla\psi|^2+2\psi\pd_t\psi \\
&= r(L\psi)^2+(t-r)\left[(\pd_t\psi)^2+(\pd_r\psi)^2\right]+t|\sla\nabla\psi|^2+2\psi L\psi-2\psi\pd_r\psi \\
&= r(L\psi+r^{-1}\psi)^2+(t-r)\left[(\pd_t\psi)^2+(\pd_r\psi)^2\right]+t|\sla\nabla\psi|^2 -r^{-1}\psi^2-2\psi\pd_r\psi \\
&= r(L\psi+r^{-1}\psi)^2+(t-r)\left[(\pd_t\psi)^2+(\pd_r\psi)^2\right]+t|\sla\nabla\psi|^2 +r^{-1}\psi^2-r^{-2}\pd_r(r^2\psi^2).
\end{align*}
\end{proof}

Due to the term with the factor $t-r$ in the above lemma (which has to do with the fact that $S$ is spacelike for $r>t$) it is necessary to transition to a different vectorfield. We will use the vectorfield $rL$ instead of $S$ for $r>t$, noting that $S=rL$ when $r=t$.

\begin{lemma}\label{m:special_rL_2_lem}
If $\Box_g\psi=0$, then the momentum $J[rL,2]$ satisfies
$$-J^t[rL,2]=r(L\psi+r^{-1}\psi)^2+r|\sla\nabla\psi|^2+r^{-1}\psi^2-r^{-2}\pd_r(r^2\psi^2)$$
and
$$divJ[rL,2] = (L\psi)^2+|\sla\nabla\psi|^2.$$
\end{lemma}
\begin{proof}
Note that since $rL = S+(r-t)\pd_t$,
$$J[rL,2] = J[S,2]+(r-t)J[\pd_t,0].$$
Thus, the formula for $-J^t[rL,2]$ follows from Lemma \ref{m:special_S_2_lem} and the fact that
$$-J^t[\pd_t,0] = (\pd_t\psi)^2+(\pd_r\psi)^2+|\sla\nabla\psi|^2.$$

By the same argument, since $divJ[S,2]=0$, we have
$$divJ[rL,2] = div J[(r-t)\pd_t,0].$$
Let $X=(r-t)\pd_t$. By Lemma \ref{divJ_lem},
$$div J[X,0] = K^{\mu\nu}\pd_\mu\psi\pd_\nu\psi,$$
where
\begin{align*}
K^{\mu\nu} &= g^{\mu\lambda}\pd_\lambda X^\nu +g^{\nu\lambda}\pd_\lambda X^\mu-X^\lambda\pd_\lambda g^{\mu\nu} - div Xg^{\mu\nu} \\
&= g^{\mu\lambda}\pd_\lambda(r-t)\delta_t^\nu +g^{\nu\lambda}\pd_\lambda(r-t)\delta_t^\mu +g^{\mu\nu}.
\end{align*}
By inspecting each case for $\mu$ and $\nu$, we conclude that
\begin{align*}
divJ[(r-t)\pd_t,0] &= 2(\pd_t\psi)^2+2\pd_t\psi\pd_r\psi +\left[-(\pd_t\psi)^2+(\pd_r\psi)^2+|\sla\nabla\psi|^2\right] \\
&= (L\psi)^2+|\sla\nabla\psi|^2.
\end{align*}
\end{proof}

As a consequence of Lemma \ref{m:special_rL_2_lem} only, we have the following identity.
\begin{proposition}
If $\Box_g\psi=0$ and $\psi$ decays sufficiently fast as $r\rightarrow\infty$, then
\begin{multline*}
\int_{\Sigma_{t_2}}r\left[(L\psi+r^{-1}\psi)^2+r|\sla\nabla\psi|^2+r^{-2}\psi^2\right] \\
+\int_{t_1}^{t_2}\int_{\Sigma_t} (L\psi)^2+|\sla\nabla\psi|^2 \\
= \int_{\Sigma_{t_1}}r\left[(L\psi+r^{-1}\psi)^2+r|\sla\nabla\psi|^2+r^{-2}\psi^2\right].
\end{multline*}
\end{proposition}
\begin{proof}
This follows by applying Proposition \ref{general_divergence_estimate_prop} to the current $J[rL,2]$ and observing Lemma \ref{m:special_rL_2_lem} and the fact that
$$\int_{\Sigma_t}-r^{-2}\pd_r(r^2\psi^2) = 0,$$
since it is the integral of a total divergence term.
\end{proof}

As a consequence of both Lemma \ref{m:special_S_2_lem} and Lemma \ref{m:special_rL_2_lem}, we have the following identity and inequality.
\begin{proposition}
If $\Box_g\psi=0$ and $\psi$ decays sufficiently fast as $r\rightarrow\infty$, then
\begin{multline*}
\int_{\Sigma_t\cap\{r<t\}} r(L\psi+r^{-1}\psi)^2+(t-r)\left[(\pd_t\psi)^2+(\pd_r\psi)^2\right]+t|\sla\nabla\psi|^2+r^{-1}\psi^2 \\
+\int_{\Sigma_t\cap\{r>t\}} r(L\psi+r^{-1}\psi)^2+r|\sla\nabla\psi|^2+r^{-1}\psi^2 \\
+\int_0^t\int_{\Sigma_\tau\cap\{r>\tau\}} (L\psi)^2+|\sla\nabla\psi|^2 \\
= \int_{\Sigma_0} r(L\psi+r^{-1}\psi)^2+r|\sla\nabla\psi|^2+r^{-1}\psi^2.
\end{multline*}
In particular,
\begin{multline*}
\int_{\Sigma_t\cap\{r<t\}}(1-r/t)\left[(\pd_t\psi)^2+(\pd_r\psi)^2\right]+|\sla\nabla\psi|^2 \\
\le t^{-1}\int_{\Sigma_0} r\left[(L\psi+r^{-1}\psi)^2+|\sla\nabla\psi|^2+r^{-2}\psi^2\right].
\end{multline*}
\end{proposition}
\begin{proof}
The identity follows by applying Proposition \ref{general_divergence_estimate_prop} to the current $J[X,2]$, where
$$
X = \left\{
\begin{array}{ll}
S & r\le t \\
rL & t\le r,
\end{array}
\right.
$$
and observing Lemma \ref{m:special_S_2_lem} in the region $\{r<t\}$, Lemma \ref{m:special_rL_2_lem} in the region $\{t<r\}$, and the fact that
$$\int_{\Sigma_t}-r^{-2}\pd_r(r^2\psi^2) = 0,$$
since it is the integral of a total divergence term.

In particular, the identity implies
\begin{multline*}
\int_{\Sigma_t\cap\{r<t\}}(t-r)\left[(\pd_t\psi)^2+(\pd_r\psi)^2\right]+t|\sla\nabla\psi|^2 \\
\le \int_{\Sigma_0} r\left[(L\psi+r^{-1}\psi)^2+|\sla\nabla\psi|^2+r^{-2}\psi^2\right].
\end{multline*}
Dividing by $t$ yields the inequality stated in the proposition.
\end{proof}

\subsection{A $p=2$ type estimate}

The $p=2$ type estimate is based on the vectorfield $K_0$ defined in the following lemma.
\begin{lemma}\label{K0_is_conformal_lem}
The vectorfield 
$$K_0=(t^2+r^2)\pd_t+2tr\pd_r$$
 is a conformal Killing vectorfield. More precisely,
$$\mathcal{L}_{K_0} g_{\mu\nu} = 4tg_{\mu\nu}.$$
\end{lemma}
\begin{proof}
From equation (\ref{lie_derivative_of_metric_eqn}),
$$\mathcal{L}_{K_0}g_{\mu\nu}=2\nabla_{(\mu}{K_0}_{\nu)}.$$
It suffices to show that
$$2\nabla^{(\mu}{K_0}^{\nu)} = 4t g^{\mu\nu}.$$
From Lemma \ref{divJ_lem},
\begin{align*}
2\nabla^{(\mu}{K_0}^{\nu)} = g^{\mu\lambda}\pd_\lambda {K_0}^\nu+g^{\nu\lambda}\pd_\lambda {K_0}^\mu-{K_0}^\lambda\pd_\lambda g^{\mu\nu}.
\end{align*}
Since $\pd_tg^{\mu\nu}=0$, recalling the proof of Lemma \ref{S_is_conformal_lem},
$$-{K_0}^\lambda\pd_\lambda g^{\mu\nu} =2t(- S^\lambda \pd_\lambda g^{\mu\nu}) = 4t \sla{g}^{\mu\nu}.$$
Now, observe that
\begin{align*}
g^{t\lambda}\pd_\lambda {K_0}^t &= g^{tt}\pd_t(t^2+r^2) = -2t \\
g^{t\lambda}\pd_\lambda {K_0}^r &= g^{tt}\pd_t(2tr) = -2r \\
g^{r\lambda}\pd_\lambda {K_0}^t &= g^{rr}\pd_r(t^2+r^2) = 2r \\
g^{r\lambda}\pd_\lambda {K_0}^r &= g^{rr}\pd_r(2tr) = 2t.
\end{align*}
At this point, it easy to verify (by checking each component)  that
$$2\nabla^{(\mu}{K_0}^{\nu)} = 4tg^{\mu\nu}.$$
\end{proof}

\begin{lemma}\label{m:special_K0_lem}
If $\Box_g\psi=0$, then the momentum $J[K_0,4t]$ satisfies
\begin{multline*}
-J^t[K_0,4t] = (t^2+r^2)|\sla\nabla\psi|^2+\frac12(t+r)^2(L\psi+r^{-1}\psi)^2+\frac12(t-r)^2(\lbar\psi-r^{-1}\psi)^2 \\
-r^{-2}\pd_r(r(t^2+r^2)\psi^2)
\end{multline*}
and
$$divJ[K_0,4t]=0.$$
\end{lemma}
\begin{proof}
The fact that
$$div J[K_0,4t]=0$$
follows immediately from Lemma \ref{K0_is_conformal_lem} and Corollary \ref{conformal_current_cor}, since Minkowski space is Ricci flat.

Now, we compute $-J^t[K_0,4t]$. Recall that
$$J^\mu[X,w] = 2\pd^\mu \psi X(\psi)-X^\mu \pd^\lambda\psi\pd_\lambda\psi +w\psi\pd^\mu \psi-\frac12\psi^2\pd^\mu w.$$
Thus,
\begin{align*}
-J^t[K_0,4t] &= -2g^{tt}\pd_t\psi K_0(\psi)+{K_0}^t\pd^\lambda\psi\pd_\lambda\psi-4tg^{tt}\psi\pd_t\psi+\frac12\psi^2g^{tt}\pd_t(4t) \\
&= 2(t^2+r^2)(\pd_t\psi)^2+4tr\pd_t\psi\pd_r\psi +(t^2+r^2)\pd^\lambda\psi\pd_\lambda\psi+4t\psi\pd_t\psi-2\psi^2 \\
&= (t^2+r^2)\left[(\pd_t\psi)^2+(\pd_r\psi)^2+|\sla\nabla\psi|^2\right]+4tr\pd_t\psi\pd_r\psi +4t\psi\pd_t\psi-2\psi^2 \\
&= (t^2+r^2)|\sla\nabla\psi|^2 + (t-r)^2\left[(\pd_t\psi)^2+(\pd_r\psi)^2\right]+2tr(L\psi)^2 \\
&\hspace{3.5in}+4t\psi\pd_t\psi-2\psi^2.
\end{align*}
Now, observe that
\begin{align*}
2tr(L\psi)^2&+4t\psi\pd_t\psi-2\psi^2 \\
&= 2tr(L\psi)^2+4t\psi L\psi -4t\psi\pd_r\psi-2\psi^2 \\
&= 2tr(L\psi+r^{-1}\psi)^2-2tr^{-1}\psi^2-4t\psi\pd_r\psi-2\psi^2 \\
&= 2tr(L\psi+r^{-1}\psi)^2-2r^{-2}\pd_r(r^2t\psi^2)+2tr^{-1}\psi-2\psi^2 \\
&= 2tr(L\psi+r^{-1}\psi)^2-2r^{-2}\pd_r(r^2t\psi^2)+2\frac{t-r}{r}\psi^2,
\end{align*}
so
\begin{multline*}
-J^t[K_0,4t] = (t^2+r^2)|\sla\nabla\psi|^2+(t-r)^2\left[(\pd_t\psi)^2+(\pd_r\psi)^2\right] +2tr(L\psi+r^{-1}\psi)^2 \\
+2\frac{t-r}r\psi^2-2r^{-2}\pd_r(r^2t\psi^2).
\end{multline*}
Furthermore,
\begin{align*}
2\frac{t-r}{r}\psi^2 &= - r^{-2}\pd_r\left((t-r)^2r\psi^2\right) +\frac{(t-r)^2}{r^2}\psi^2 + 2\frac{(t-r)^2}{r}\psi\pd_r\psi \\
&= - r^{-2}\pd_r\left((t-r)^2r\psi^2\right) +(t-r)^2(\pd_r\psi+r^{-1}\psi)^2-(t-r)^2(\pd_r\psi)^2,
\end{align*}
so
\begin{multline*}
-J^t[K_0,4t] = (t^2+r^2)|\sla\nabla\psi|^2+(t-r)^2\left[(\pd_t\psi)^2+(\pd_r\psi+r^{-1}\psi)^2\right] +2tr(L\psi+r^{-1}\psi)^2 \\
-2r^{-2}\pd_r(r^2t\psi^2)- r^{-2}\pd_r\left((t-r)^2r\psi^2\right).
\end{multline*}
Finally, since
\begin{multline*}
(t-r)^2\left[(\pd_t\psi)^2+(\pd_r\psi+r^{-1}\psi)^2\right] +2tr(L\psi+r^{-1}\psi)^2 \\
=\frac12(t-r)^2\left[(L\psi+r^{-1}\psi)^2+(\lbar\psi-r^{-1}\psi)^2\right]+2tr(L\psi+r^{-1}\psi)^2 \\
= \frac12(t+r)^2(L\psi+r^{-1}\psi)^2+\frac12(t-r)^2(\lbar\psi-r^{-1}\psi)^2
\end{multline*}
and
$$-2r^{-2}\pd_r(r^2t\psi^2)- r^{-2}\pd_r\left((t-r)^2r\psi^2\right) = -r^{-2}\pd_r(r(t^2+r^2)\psi^2),$$
we conclude that
\begin{multline*}
-J^t[K_0,4t] = (t^2+r^2)|\sla\nabla\psi|^2+\frac12(t+r)^2(L\psi+r^{-1}\psi)^2+\frac12(t-r)^2(\lbar\psi-r^{-1}\psi)^2 \\
-r^{-2}\pd_r(r(t^2+r^2)\psi^2).
\end{multline*}
\end{proof}

As a consequence of Lemma \ref{m:special_K0_lem}, we have the following identity and inequality.
\begin{proposition}
If $\Box_g\psi=0$ and $\psi$ decays sufficiently fast as $r\rightarrow\infty$, then
\begin{multline*}
\int_{\Sigma_{t_2}} (t^2+r^2)|\sla\nabla\psi|^2+\frac12(t+r)^2(L\psi+r^{-1}\psi)^2+\frac12(t-r)^2(\lbar\psi-r^{-1}\psi)^2 \\
= \int_{\Sigma_{t_1}} (t^2+r^2)|\sla\nabla\psi|^2+\frac12(t+r)^2(L\psi+r^{-1}\psi)^2+\frac12(t-r)^2(\lbar\psi-r^{-1}\psi)^2.
\end{multline*}
In particular,
\begin{multline*}
\int_{\Sigma_t} |\sla\nabla\psi|^2+\frac12(L\psi+r^{-1}\psi)^2 \\
\le t^{-2} \int_{\Sigma_0} r^2|\sla\nabla\psi|^2+\frac12r^2(L\psi+r^{-1}\psi)^2+\frac12r^2(\lbar\psi-r^{-1}\psi)^2.
\end{multline*}
\end{proposition}
\begin{proof}
The identity follows by applying Proposition \ref{general_divergence_estimate_prop} to the current $J[K_0,4t]$ and observing Lemma \ref{m:special_K0_lem} and the fact that
$$\int_{\Sigma_t}-r^{-2}\pd_r(r(t^2+r^2)\psi^2)=0,$$
since it is the integral of a total divergence term.

In particular, the identity implies
\begin{multline*}
\int_{\Sigma_t} t^2|\sla\nabla\psi|^2+\frac12 t^2(L\psi+r^{-1}\psi)^2 \\
= \int_{\Sigma_0} r^2|\sla\nabla\psi|^2+\frac12r^2(L\psi+r^{-1}\psi)^2+\frac12r^2(\lbar\psi-r^{-1}\psi)^2.
\end{multline*}
Dividing by $t^2$ yields the inequality stated in the proposition.
\end{proof}

\section{The dynamic estimates}\label{m:dynamic_sec}

In this section, we prove the dynamic estimates, which provide all of the necessary information related to the future dynamics of the wavefunction $\psi$. These estimates are simply restatements of the energy estimate (Proposition \ref{m:classic_ee_prop}) and the $r^p$ estimate (Proposition \ref{m:rp_prop}) applied to either $\psi$ itself or a wavefunction $\psi^s$ derived from $\psi$.

The dynamic estimates take the form
$$E(t_2)\lesssim E(t_1)+\int_{t_1}^{t_2}N(t)dt$$
$$E_p(t_2)+\int_{t_1}^{t_2}B_p(t)dt\lesssim E_p(t_1)+\int_{t_1}^{t_2}N_p(t)dt,$$
where $p$ ranges from $\delm$ to $2-\delp$ for arbitrarily small $\delm,\delp>0$. The norm $E(t)$ is the energy norm and the norm $E_p(t)$ is the $r^p$ weighted energy norm with a weight of $r^p$ near $i^0$. The norm $B_p(t)$ is equivalent to $E_{p-1}(t)$ in Minkowski space, but in the remaining problems that involve black holes, the norm $B_p(t)$ will have a degeneracy at the \textit{photon sphere}. The norm $B_p(t)$ has a weight of $r^{p-1}$ near $i^0$. The norms $N(t)$ and $N_p(t)$ are both nonlinear error norms which can be ignored in the linear problem. The estimates are put in this form for the convenience of the reader, who is strongly encouraged to have this form memorized.

\subsection{The dynamic estimates for $\psi$ ($s=0$)}

We begin with the dynamic estimates for the wavefunction $\psi$ only.
\begin{proposition}\label{m:dynamic_estimates_0_prop}
Fix $\delm,\delp>0$ and let $p\in[\delm,2-\delp]$. Then
$$E(t_2)\lesssim E(t_1)+\int_{t_1}^{t_2}N(t)dt$$
$$E_p(t_2)+\int_{t_1}^{t_2}B_p(t)dt \lesssim E_p(t_1)+\int_{t_1}^{t_2}N_p(t)dt,$$
where
$$E(t)=\int_{\Sigma_t}(\pd_t\psi)^2+(\pd_r\psi)^2+|\sla\nabla\psi|^2,$$
$$E_p(t)=\int_{\Sigma_t}(1+r)^p\left[(L\psi)^2+|\sla\nabla\psi|^2+r^{-1}(1+r)^{-1}\psi^2+(1+r)^{-2}(\pd_r\psi)^2\right],$$
$$B_p(t)=\int_{\Sigma_t}(1+r)^{p-1}\left[(L\psi)^2+|\sla\nabla\psi|^2+r^{-1}(1+r)^{-1}\psi^2+(1+r)^{-2}(\pd_r\psi)^2\right],$$
$$N(t)=\int_{\Sigma_t}|\pd_t\psi\Box_g\psi|,$$
$$N_p(t)=\int_{\Sigma_t}(1+r)^p(|L\psi|+(1+r)^{-1}|\psi|+(1+r)^{-2}|\pd_r\psi|)|\Box_g\psi|.$$
\end{proposition}
\begin{proof}
The first estimate is a restatement of the energy estimate (Proposition \ref{m:classic_ee_prop}) and the second is a restatement of the $r^p$ estimate (Proposition \ref{m:rp_prop}).
\end{proof}

We note that one of the estimates can be simplified slightly by absorbing part of the nonlinear norm $N_p(t)$ into the bulk $B_p(t)$ on the left side.
\begin{corollary}\label{m:dynamic_estimates_0_cor}
Fix $\delm,\delp>0$ and let $p\in[\delm,2-\delp]$. Then
$$E(t_2)\lesssim E(t_1)+\int_{t_1}^{t_2}N(t)dt$$
$$E_p(t_2)+\int_{t_1}^{t_2}B_p(t)dt \lesssim E_p(t_1)+\int_{t_1}^{t_2}N_p(t)dt,$$
where $E(t)$, $E_p(t)$, and $B_p(t)$ are as defined in Proposition \ref{m:dynamic_estimates_0_prop} and
$$N(t)=(E(t))^{1/2}||\Box_g\psi||_{L^2(\Sigma_t)}$$
$$N_p(t)=\int_{\Sigma_t}(1+r)^{p+1}(\Box_g\psi)^2.$$
\end{corollary}
\begin{proof}
The first estimate is due to Proposition \ref{m:dynamic_estimates_0_prop} and the fact that
$$\int_{\Sigma_t}|\pd_t\psi\Box_g\psi|\lesssim ||\pd_t\psi||_{L^2(\Sigma_t)}||\Box_g\psi||_{L^2(\Sigma_t)} \lesssim (E(t))^{1/2}||\Box_g\psi||_{L^2(\Sigma_t)}.$$
We turn to the second estimate.

According to Proposition \ref{m:dynamic_estimates_0_prop}, there is a constant $C$ such that
$$E_p(t_2)+\int_{t_1}^{t_2}B_p(t)dt \le C E_p(t_1)+ C\int_{t_1}^{t_2}N_p'(t)dt,$$
where
$$N_p'(t)=\int_{\Sigma_t}(1+r)^p(|L\psi|+(1+r)^{-1}|\psi|+(1+r)^{-2}|\pd_r\psi|)|\Box_g\psi|.$$
Now,
\begin{multline*}
N_p'(t) \\
\lesssim \epsilon\int_{\Sigma_t}(1+r)^{p-1}((L\psi)^2+(1+r)^{-2}\psi^2+(1+r)^{-4}(\pd_r\psi)^2)+\epsilon^{-1}\int_{\Sigma_t}(1+r)^{p+1}(\Box_g\psi)^2 \\
\lesssim \epsilon B_p(t)+\epsilon^{-1}N_p(t).
\end{multline*}
Thus,
$$E_p(t_2)+(1-C\epsilon)\int_{t_1}^{t_2}B_p(t)dt\le C E_p(t_1)+C\epsilon^{-1}\int_{t_1}^{t_2}N_p(t)dt.$$
It suffices to choose $\epsilon = (2C)^{-1}$. This process is called \textit{absorbing the term $\epsilon\int_{t_1}^{t_2}B_p(t)dt$ into the left side}.
\end{proof}

\subsection{Commutators with $\Box_g$}

Minkowski space has many symmetries. \cite{klainerman1987remarks} These symmetries are generated by the translations
$$\mathbb{T}_\mu=\pd_\mu=\{\pd_x,\pd_y,\pd_z,\pd_t\},$$
the rotations about each axis
$$\Omega_{ij}=x^i\pd_j-x^j\pd_i=\{\Omega_x,\Omega_y,\Omega_z\},$$
and the Lorentz boosts
$$L_{0i}=t\pd_i+x^i\pd_t=\{L_x,L_y,L_z\},$$
which correspond to rotations involving one spatial dimension and the time dimension. These symmetry operators are useful, because they commute with the wave operator $\Box_g$, so, for example, the dynamics of $\psi$ will be similar to the dynamics of $\Omega_x\psi$.

In addition, there are also the conformal transformations, which are generated by the vectorfields
$$\mathbb{S}=x^\mu\pd_\mu$$
and
$$\mathbb{K}_\mu=-2x_\mu \mathbb{S}+\mathbb{S}^\lambda\mathbb{S}_\lambda \pd_\mu,$$
of which $\mathbb{S}$ and $\mathbb{K}_0$ were used as multipliers in \S\ref{m:special_sec}. These operators can also be used as commutators, but they will not be used in this thesis.

Most of the above symmetries do not generalize well to the black hole spacetimes in later chapters, with the exceptions being the time translation symmetry $\pd_t$ and (in varying senses) the rotation symmetries $\Omega$. For this reason, we will try to only use these symmetries going forward. However, due to a slight complication near $r=0$, we will also use the spatial translations. These will not be necessary in black hole spacetimes, because the exteriors of such spacetimes have no center.

\begin{definition}
Denote by $\Gamma$ any of the following operators.
$$\Gamma \in\{\pd_t,\Omega_x,\Omega_y,\Omega_z\}\cup\{\pd_x,\pd_y,\pd_z\}$$
(The union is to emphasize that the spatial translation operators should be considered less important as they will not be used in black hole spacetimes.) Furthermore, denote by $\Gamma^s$ any composition of $s$ of these operators.
\end{definition}

We also define the $s$-order wavefunctions derived from the wavefunction $\psi$ by applying symmetry operators.
\begin{definition}
$$\psi^s=\Gamma^s\psi$$
\end{definition}

As previously discussed, these generalized wavefunctions behave very similarly to $\psi$. This is the statement of the following lemma.
\begin{lemma}\label{m:psi_s_lem}
If $\Box_g\psi=0$, then
$$\Box_g\psi^s=0.$$
More generally,
$$\Box_g\psi^s=\Gamma^s(\Box_g\psi).$$
\end{lemma}
\begin{proof}
We only prove the general case.
$$\Box_g\psi^s=\Box_g\Gamma^s\psi=\Gamma^s(\Box_g\psi)+[\Box_g,\Gamma^s]\psi=\Gamma^s(\Box_g\psi).$$
\end{proof}

\subsection{Higher order dynamic estimates}

Now we can write the dynamic estimates in their most general form.
\begin{proposition}\label{m:dynamic_estimates_s_prop}
Fix $\delm,\delp>0$. The following estimates hold for $p\in[\delm,2-\delp]$.
$$E^s(t_2)\lesssim E^s(t_1)+\int_{t_1}^{t_2}N^s(t)dt$$
$$E_p^s(t_2)+\int_{t_1}^{t_2}B_p^s(t)dt \lesssim E_p^s(t_1)+\int_{t_1}^{t_2}N_p^s(t)dt,$$
where
$$E^s(t)=\sum_{s'\le s} E[\psi^{s'}](t)$$
$$E_p^s(t)=\sum_{s'\le s} E_p[\psi^{s'}](t)$$
$$B_p^s(t)=\sum_{s'\le s} B_p[\psi^{s'}](t)$$
$$N^s(t)=(E^s(t))^{1/2}\sum_{s'\le s}||\Gamma^{s'}(\Box_g\psi)||_{L^2(\Sigma_t)}$$
$$N_p^s(t)=\sum_{s'\le s}\int_{\Sigma_t}(1+r)^{p+1}(\Gamma^{s'}(\Box_g\psi))^2,$$
and the norms $E(t)$, $E_p(t)$, $B_p(t)$, are as defined in Proposition \ref{m:dynamic_estimates_0_prop}.
\end{proposition}

\begin{proof}
The proof is a direct application of Corollary \ref{m:dynamic_estimates_0_cor} by making the substitutions
$$\psi\mapsto\psi^{s'}$$
for all values of $s'$ (and all commutators represented by $\Gamma^{s'}$) where $s'\le s$ and observing Lemma \ref{m:psi_s_lem}.
\end{proof}

\section{The $L^\infty$ estimates}\label{m:pointwise_sec}

In this section, we prove $L^\infty$ estimates for certain derivatives of $\psi$ that will appear in the nonlinear terms belonging to $\Gamma^s(\Box_g\psi)$ when the equation for $\Box_g\psi$ is used. We begin with Lemma \ref{m:sobolev_I_lem}, which estimates an arbitrary function in $L^\infty(\Sigma_t\cap\{r\ge r_0\})$ away from the center $\{r=0\}$. This lemma will be generalized to the remaining problems of this thesis. For this problem only, we also need to prove Lemma \ref{m:sobolev_II_lem}, which estimates an arbitrary function in $L^\infty(\Sigma_t\cap\{r\le r_0\})$ near the center $\{r=0\}$. These two lemmas are combined in Proposition \ref{m:infty_start_prop}.

In \S\ref{m:pointwise_lemmas_sec}, Proposition \ref{m:infty_start_prop} is then repeatedly applied to single derivatives of $\psi^s$, resulting in Sobolev norms that can be estimated by the energy norms. (See Lemmas \ref{m:psi_pointwise_lem}-\ref{m:Lpsi_pointwise_lem}.)  These estimates are all summarized at the end in Proposition \ref{m:infty_prop}.

It is important to pay special attention to the $r$ weights in the lemmas that follow. In the main theorem (Theorem \ref{m:main_thm}), we will see that the weighted energy norms $E_p^{s}(t)$ behave like $t^{p-2+\delp}$ for late times and that the energy $E^{s}(t)$ will remaind bounded in time. In this section, we will see that multiplying a derivative of $\psi^s$ by $r$ changes which energy norm can be used to estimate the Sobolev norm provided by Proposition \ref{m:infty_start_prop}. Some derivatives of $\psi^s$ can have more $r$ factors than others. (See Definition \ref{m:DbarD_def}.) These derivatives will eventually be shown to decay better in time. (See the statement of Theorem \ref{m:main_thm}.)

\subsection{A Sobolev-type estimate}

First, we prove the following lemma, which is a Sobolev-type estimate.
\begin{lemma}\label{m:sobolev_I_lem}
If $u$ decays sufficiently fast as $r\rightarrow\infty$, then
$$||u||_{L^\infty(\Sigma_t\cap\{r\ge r_0\})}^2\lesssim \int_{\Sigma_t\cap\{r\ge r_0\}}r^{-2}[(\pd_r\Omega^{\le 2}u)^2+(\Omega^{\le 2}u)^2].$$
\end{lemma}
\begin{proof}
For a fixed $r$, denote by $\bar{u}:S^2(1)\rightarrow\R{}$ the pullback of the function $u:S^2(r)\rightarrow\R{}$ via the canonical map from $S^2(1)$ to $S^2(r)$. Also, denote by $d\omega$ the measure on $S^2(1)$. Then
\begin{multline*}
||u||_{L^\infty(S^2(r))}^2=||\bar{u}||_{L^\infty(S^2(1))}^2\lesssim \int_{S^2(1)}|\sla\nabla^{\le 2}\bar{u}|^2d\omega = \int_{S^2(1)}|\overline{\Omega^{\le 2}u}|^2d\omega \\
=\int_{S^2(r)}(\Omega^{\le 2}u)^2d\omega =r^{-2}\int_{S^2(r)}(\Omega^{\le 2}u)^2r^2d\omega.
\end{multline*}
Now, set $f(r)=\int_{S^2(r)}(\Omega^{\le 2}u)^2d\omega$. Note that
\begin{multline*}
|f'(r)|\lesssim \int_{S^2(r)}|\Omega^{\le 2}u\pd_r \Omega^{\le 2}u|d\omega \lesssim \int_{S^2(r)}[(\pd_r \Omega^{\le 2}u)^2+(\Omega^{\le 2}u)^2]d\omega \\
\lesssim \int_{S^2(r)}r^{-2}[(\pd_r \Omega^{\le 2}u)^2+(\Omega^{\le 2}u)^2]r^2d\omega.
\end{multline*}
Then, assuming $\lim_{r_\rightarrow\infty}f(r)=0$,
\begin{multline*}
|u(r_0)|^2\lesssim f(r_0)\lesssim \int_{r_0}^\infty |f'(r)|dr 
\lesssim \int_{r_0}^{\infty}\int_{S^2(r)}r^{-2}[(\pd_r\Omega^{\le 2}u)^2+(\Omega^{\le 2}u)]r^2d\omega dr \\
= \int_{\Sigma_t\cap\{r\ge r_0\}}r^{-2}[(\pd_r \Omega^{\le 2}u)^2+(\Omega^{\le 2}u)^2].
\end{multline*}
\end{proof}

Unlike in black hole exterior spacetimes, the Minkowski spacetime has a center $\{r=0\}$. The previous lemma requires strong weights near $r=0$, so we also use the classic Sobolev-type estimate, which uses only spatial translation operators. This strategy is only used in this chapter.
\begin{lemma}\label{m:sobolev_II_lem}
$$||u||^2_{L^\infty(\Sigma_t\cap\{r<r_0\})}\lesssim \int_{\Sigma_t\cap\{r<2r_0\}}(\pd_{x^i}^{\le 2}u)^2.$$
\end{lemma}
\begin{proof}
Let $\chi$ be a smooth cutoff function such that $\chi=1$ for $r\le r_0$ and $\chi=0$ for $r\ge 2r_0$. Then
$$||u||_{L^\infty(\Sigma_t\cap\{r<r_0\})}^2\le ||\chi u||_{L^\infty(\Sigma_t\cap\{r<2r_0\})}^2\lesssim \int_{\Sigma_t\cap\{r<2r_0\}}(\pd_{x^i}^{\le 2}(\chi u))^2\lesssim \int_{\Sigma_t\cap\{r<2r_0\}}(\pd_{x^i}^{\le 2}u)^2.$$
\end{proof}

Combining both of these lemmas, we have an estimate with better weights.
\begin{proposition}\label{m:infty_start_prop}
If $u$ decays sufficiently fast as $r\rightarrow\infty$, then
$$||u||^2_{L^\infty(\Sigma_t)}\lesssim \int_{\Sigma_t}(1+r)^{-2}\left[(\pd_r\Gamma^{\le 2}u)^2+(\Gamma^{\le 2}u)^2\right].$$
\end{proposition}
\begin{proof}
We have
$$||u||^2_{L^\infty(\Sigma_t)}\le ||u||^2_{L^\infty(\Sigma_t\cap\{r>1\})}+||u||^2_{L^\infty(\Sigma_t\cap\{r<1\})}.$$
We estimate the first term on the right side using Lemma \ref{m:sobolev_I_lem} and the second term using Lemma \ref{m:sobolev_II_lem}.

From Lemma \ref{m:sobolev_I_lem}, taking $r_0=1$, we have
\begin{multline*}
||u||^2_{L^\infty(\Sigma_t\cap\{r>1\})}\lesssim \int_{\Sigma_t\cap\{r>1\}}r^{-2}\left[(\pd_r\Gamma^{\le 2}u)^2+(\Gamma^{\le 2}u)^2\right] \\
\lesssim \int_{\Sigma_t}(1+r)^{-2}\left[(\pd_r\Gamma^{\le 2}u)^2+(\Gamma^{\le 2}u)^2\right].
\end{multline*}
From Lemma \ref{m:sobolev_II_lem}, taking $r_0=1$, we have
$$||u||^2_{L^\infty(\Sigma_t\cap\{r<1\})}\lesssim \int_{\Sigma_t\cap\{r<2\}}(\Gamma^{\le 2}u)^2 \lesssim \int_{\Sigma_t}(1+r)^{-2}(\Gamma^{\le 2}u)^2.$$
\end{proof}

\subsection{Estimating derivatives using the Sobolev-type estimate}\label{m:pointwise_lemmas_sec}

Now we repeatedly apply Proposition \ref{m:infty_start_prop} to estimate various derivatives of $\psi$ with $r$ weights. We will assume that $\psi$ decays sufficiently fast as $r\rightarrow\infty$.

The following lemma estimates $\psi$ and the higher order analogues $\psi^s$.
\begin{lemma}\label{m:psi_pointwise_lem}
$$|r^p\psi^s|^2\lesssim E^{s+2}_{2p}(t)$$
$$|\psi^s|^2\lesssim E^{s+2}(t).$$
\end{lemma}
\begin{proof}
First, we apply Proposition \ref{m:infty_start_prop} with $u=r^p\psi^s$.
\begin{align*}
|r^p\psi^s|^2 &\lesssim \int_{\Sigma_t} (1+r)^{-2}\left[(\pd_r\Gamma^{\le 2}(r^p\psi^s))^2+(\Gamma^{\le 2}(r^p\psi^s))^2\right] \\
&\lesssim \int_{\Sigma_t}(1+r)^{2p-2}\left[(\pd_r\Gamma^{\le 2}\psi^s)^2+(\Gamma^{\le 2}\psi^s)^2\right] \\
&\lesssim E_{2p}^{s+2}(t).
\end{align*}
This verifies the first estimate. The second estimate follows from the exact same argument in the special case $p=0$, and using $E^{s+2}(t)$ instead of $E_0^{s+2}(t)$.
\end{proof}

The following lemma estimates $\pd_t\psi$ and the higher order analogues $\pd_t\psi^s$.
\begin{lemma}
$$|r^p\pd_t\psi^s|^2\lesssim E_{2p}^{s+3}(t)$$
$$|r\pd_t\psi^s|^2\lesssim E^{s+3}(t)$$
\end{lemma}
\begin{proof}
The first estimate reduces to Lemma \ref{m:psi_pointwise_lem} by observing that $\pd_t\psi^s=\psi^{s+1}$. We now prove the second estimate.

First, we apply Proposition \ref{m:infty_start_prop} with $u=r\pd_t\psi^s$.
\begin{align*}
|r\pd_t\psi^s|^2 &\lesssim \int_{\Sigma_t}(1+r)^{-2}\left[(\pd_r\Gamma^{\le 2}(r\pd_t\psi^s))^2+(\Gamma^{\le 2}(r\pd_t\psi^s))^2\right] \\
&\lesssim \int_{\Sigma_t}\left[(\pd_r\psi^{s+3})^2+(\pd_t\psi^{s+2})^2\right] \\
&\lesssim E^{s+3}(t).
\end{align*}
The point is that $\pd_t\psi^s$ can either be treated as $\psi^{s+1}$ or as a derivative of $\psi^s$. In the latter case, the energy norm $E^s(t)$ has stronger control, because $||\pd_t\psi^s||_{L^2(\Sigma_t)}^2\le E^s(t)$, while $||r^{-1}\psi^{s+1}||_{L^2(\Sigma_t)}^2\le E^{s+1}(t)$.
\end{proof}

The following lemma estimates $\sla\nabla\psi$ and the higher order analogues $\sla\nabla\psi^s$.
\begin{lemma}
$$|r^{p+1}\sla\nabla\psi^s|^2\lesssim E^{s+3}_{2p}(t)$$
$$|r\sla\nabla\psi^s|^2\lesssim E^{s+3}(t)$$
\end{lemma}
\begin{proof}
This lemma reduces to Lemma \ref{m:psi_pointwise_lem} by observing that $r\sla\nabla\psi^s=\Omega\psi^s=\psi^{s+1}$.
\end{proof}

The following lemma estimates $L\psi$ and the higher order analogues $L\psi^s$.
\begin{lemma}\label{m:Lpsi_pointwise_lem}
$$|r^{p+1}L\psi^s|\lesssim E_{2p}^{s+3}(t)+\int_{\Sigma_t}r^{2p}(\Box_g\psi^{s+2})^2$$
$$|rL\psi^s|\lesssim E^{s+3}(t)+\int_{\Sigma_t}(\Box_g\psi^{s+2})^2$$
\end{lemma}
\begin{proof}
Before beginning the estimates stated by the lemma, it is important to establish the following estimate, which allows us to replace the operator $\pd_rL$ with a sum of other operators, including the wave operator.
$$(\pd_r Lu)^2\lesssim (\Box_gu)^2+(L\pd_tu)^2+r^{-2}(\pd_ru)^2+r^{-2}|\sla\nabla\Omega u|^2$$
To verify this, we expand
\begin{align*}
\Box_gu &= -\pd_t^2u+\pd_r^2u+2r^{-1}\pd_ru+\sla\triangle u \\
&= -\pd_t^2u+\pd_rLu-\pd_r\pd_tu+2r^{-1}\pd_ru+r^{-1}\sla\nabla\Omega u \\
&= -L\pd_tu+\pd_rLu+2r^{-1}\pd_ru+r^{-1}\sla\nabla\Omega u.
\end{align*}

Now we start the usual way, applying Proposition \ref{m:infty_start_prop} with $u=r^{p+1}L\psi^s$.
\begin{align*}
|r^{p+1}L\psi^s|^2 &\lesssim \int_{\Sigma_t}(1+r)^{-2}\left[(\pd_r(\Gamma^{\le 2}r^{p+1}L\psi^s))^2+(\Gamma^{\le 2}r^{p+1}L\psi^s)^2\right] \\
&\lesssim \int_{\Sigma_t}(1+r)^{2p}\left[(\pd_rL\psi^{s+2})^2+(L\psi^{s+2})^2\right]
\end{align*}
At this point, we cannot estimate the term $(\pd_rL\psi^{s+2})^2$ by any of the energy norms, because it has two derivatives that are not commutators. Instead, we replace it with other operators including the wave operator.
\begin{align*}
|r^{p+1}L\psi^s|^2 &\lesssim \int_{\Sigma_t}(1+r)^{2p}\left[(\pd_rL\psi^{s+2})^2+(L\psi^{s+2})^2\right] \\
&\lesssim \int_{\Sigma_t}(1+r)^{2p}\left[(\Box_g\psi^{s+2})^2+(L\pd_t\psi^{s+2})^2+r^{-2}(\pd_r\psi^{s+2})^2+r^{-2}|\sla\nabla\Omega\psi^{s+2}|^2\right. \\
&\hspace{4.4in}\left.+(L\psi^{s+2})^2\right] \\
&\lesssim E_{2p}^{s+3}(t)+\int_{\Sigma_t}(1+r)^{2p}(\Box_g\psi^{s+2})^2.
\end{align*}
The second estimate follows from the same argument in the special case $p=0$ and using $E^{s+3}(t)$ in place of $E_0^{s+3}(t)$.
\end{proof}

\subsection{Summarizing the $L^\infty$ estimates}\label{m:pointwise_summary_sec}

The lemmas from \S\ref{m:pointwise_lemmas_sec} show that some derivatives can have extra $r$ weights. This is due to the fact that the energy norm $E_p(t)$ has different weights for different derivatives. We use the notation $\bar{D}$ to refer to any derivative that can have an extra $r$ weight in the lemmas from \S\ref{m:pointwise_lemmas_sec}. More specifically,

\begin{definition}\label{m:DbarD_def}
We define two families of derivatives.
$$\bar{D}=\{L,\sla\nabla\}$$
$$D=\{L,\pd_t,\sla\nabla\}$$
Derivatives in the family $\bar{D}$ are called ``good derivatives'' and the family $D$ spans the entire tangent space at each point.
\end{definition}

To conclude this section, we summarize the lemmas from \S\ref{m:pointwise_lemmas_sec} in a single proposition.
\begin{proposition}\label{m:infty_prop}
$$|r^{p+1}\bar{D}\psi^s|^2+|r^pD\psi^s|^2+|r^p\psi^s|^2 
\lesssim E_{2p}^{s+3}(t)+\int_{\Sigma_t}(1+r)^{2p}(\Box_g\psi^{s+2})^2$$
$$|rD\psi^s|^2 \lesssim E^{s+3}(t)+\int_{\Sigma_t}(\Box_g\psi^{s+2})^2$$
\end{proposition}
\begin{proof}
Each particular case has been covered in one of the Lemmas \ref{m:psi_pointwise_lem}-\ref{m:Lpsi_pointwise_lem}.
\end{proof}

\section{The decay mechanism based on the dynamic estimates}\label{m:decay_sec}

In this section, we see how the dynamic estimates imply time decay of some energy norms. Since the energy norms can be used to estimate derivatives of $\psi$ in $L^\infty$ (see Proposition \ref{m:infty_prop}), this translates into time decay of $\psi$ and its derivatives.

\subsection{The correspondence between $r$ and $t$}

The following lemma is the basis for the decay mechanism based on the dynamic estimates. It allows one to obtain a factor of $t^{-1}$ by decreasing $p$ by $1$.

\begin{lemma}\label{m:decay_mechanism_lem}
Let $\delp,\delm>0$ be fixed, and suppose $p+1,p\in [\delm,2-\delp]$. Suppose furthermore that
$$\Box_g\psi=0.$$
Then since
$$E_{p+1}^s(t_2)+\int_{t_1}^{t_2}B_{p+1}^s(t)dt\lesssim E_{p+1}^s(t_1),$$
$$E_p^s(t_2)+\int_{t_1}^{t_2}B_p^s(t)dt\lesssim E_p^s(t_1),$$
and
$$E_p^s(t)\lesssim B_{p+1}^s(t),$$
it follows that
$$E_p^s(t)\lesssim (1+t)^{-1}E_{p+1}^s(t/2).$$
\end{lemma}
\begin{remark}
Observe that the increase in $p$ is due to the fact that $B_p(t)$ has a weight of $(1+r)^{p-1}$ and not $(1+r)^p$.
\end{remark}
\begin{remark}
In black hole spacetimes, it will not be true that $E_p^s(t)\lesssim B_{p+1}^s(t)$. This due to the trapping phenomenon, and it will require a loss of derivatives (that is, $s$ will increase as well as $p$).
\end{remark}
\begin{proof}
Following the approach in \cite{ionescu2014global}, we use the mean value theorem.\footnote{A popular approach is to use the so-called ``pigeon-hole principle,'' but the mean value theorem slightly simplifies the argument.} For a given $t$, let $t'\in[t/2,t]$ be the value for which
$$B_{p+1}^s(t')=\frac2t\int_{t/2}^tB_{p+1}^s(\tau)d\tau.$$
Then,
$$E_p^s(t)\lesssim E_p^s(t')\lesssim B_{p+1}^s(t') = \frac2t\int_{t/2}^tB_{p+1}^s(\tau)d\tau \lesssim t^{-1}E_{p+1}^s(t/2).$$
Also,
$$E_p^s(t)\lesssim E_p^s(t/2)\lesssim E_{p+1}^s(t/2).$$
Thus, we conclude
$$E_p^s(t)\lesssim (1+t)^{-1}E_{p+1}^s(t/2).$$
\end{proof}

For an example application of this lemma see \S\ref{m:recover_decay_sec}.

\subsection{Interpolating between weighted energy norms}

The decay mechanism (Lemma \ref{m:decay_mechanism_lem}) allows us to obtain a decay factor of $t^{-1}$ when decreasing $p$ by $1$. By interpolation, it is possible to recover decay for all values of $p$ in-between. This is the purpose of the following lemma.

\begin{lemma}\label{m:interpolation_lem}
Suppose $p_0<p_1$ and
$$E_{p_0}^s(t)\lesssim (1+t)^{p_0-2+\delp},$$
$$E_{p_1}^s(t)\lesssim (1+t)^{p_1-2+\delp}.$$
Then for all $p\in[p_0,p_1]$,
$$E_p^s(t)\lesssim (1+t)^{p-2+\delp}.$$
\end{lemma}
\begin{proof}
Fix $p\in (p_0,p_1)$ and let $\lambda\in (0,1)$ be such that $p=\lambda p_0+(1-\lambda)p_1$. Then by Holder's inequality, and letting $(...)$ represent the derivatives of $\psi$ appearing in the weighted energies,
\begin{multline*}
\int_{\Sigma_t}(1+r)^p(...) = \int_{\Sigma_t}(1+r)^{\lambda p_0+(1-\lambda)p_1}(...) \lesssim \left(\int_{\Sigma_t}(1+r)^{p_0}(...)\right)^\lambda\left(\int_{\Sigma_t}(1+r)^{p_1}(...)\right)^{1-\lambda} \\
\lesssim ((1+t)^{p_0-2+\delp})^\lambda ((1+t)^{p_1-2+\delp})^{1-\lambda} = (1+t)^{\lambda p_0+(1-\lambda)p_1-2+\delp} = (1+t)^{p-2+\delp}.
\end{multline*}
\end{proof}

For an example application of this lemma, see \S\ref{m:recover_decay_sec}.

\subsection{A weaker notion of decay}\label{weak_decay_sec}

In problems involving nonlinear equations, it is necessary to obtain a decay rate for some quantity which is globally integrable in time. This appears to be a problem, because the best estimate (obtained by the energy $E_{\delm}^s(t)$) is $(1+t)^{(\delm-2+\delp)/2}$, which is worse than $(1+t)^{-1}$. However, there is a weaker notion of decay (achieved using the energy $E_{\delm-1}^s(t)$) which essentially allows us to pretend we can get $(1+t)^{(\delm-3+\delp)/2}$ as long as we only use it in an integrated sense. Fortunately, this weaker notion of decay is still strong enough to solve all of the problems in this thesis. (For example, see \S\ref{m:recover_energy_sec}.)

\begin{definition}
Let $T=1+t$. A function $f(t)$ decays like $T^{-a}$ weakly if
$$\int_{T}^{\infty}f(t)dt\lesssim T^{-a+1}.$$
\end{definition}
We will often use the following fact for functions that decay weakly. It says that we can treat a function decaying weakly as if it decays strongly as long as we are integrating in time.
\begin{lemma}\label{weak_decay_lem}
If $f(t)$ decays like $T^{-a}$ weakly, then
$$\int_{t_1}^{t_2}t^b(f(t))^{1/2}dt \lesssim \int_{t_1}^{t_2}t^{b-a/2}dt.$$
\end{lemma}
\begin{proof}
We consider three separate cases based on the sign of $b-a/2+1$. (That is, whether the integral on the right side grows polynomially in $t_2$, logarithmically in $t_2$, or decays polynomially in $t_1$.) In all three cases, the main idea of the proof is the same.

For any $\alpha>1$, we have
\begin{multline*}
\int_{\alpha^k}^{\alpha^{k+1}}t^b(f(t))^{1/2}dt\lesssim \left(\int_{\alpha^k}^{\alpha^{k+1}}t^{2b}dt\right)^{1/2}\left(\int_{\alpha^k}^{\alpha^{k+1}}f(t)dt\right)^{1/2} \\
\lesssim \left((\alpha^{k+1})^{2b+1}\right)^{1/2}\left((\alpha^k)^{-a+1}\right)^{1/2}\lesssim (\alpha^k)^{b-a/2+1}=(\alpha^{b-a/2+1})^k.
\end{multline*}
If $b-a/2+1>0$, choose $\alpha>1$ so that $\alpha^{b-a/2+1}=2$ and choose $n$ so that $t_2\approx \alpha^n$.
Then
\begin{multline*}
\int_1^{t_2}t^b(f(t))^{1/2}dt\lesssim \sum_{k=0}^n\int_{\alpha^k}^{\alpha^{k+1}}t^b(f(t))^{1/2}dt\lesssim \sum_{k=0}^n(\alpha^{b-a/2+1})^k = \sum_{k=0}^n 2^k\lesssim 2^n \\
\approx t_2^{b-a/2+1}.
\end{multline*}
If $b-a/2+1<0$, choose $\alpha>1$ so that $\alpha^{b-a/2+1}=2^{-1}$ and choose $n$ so that $t_1\approx\alpha^n$.
Then
\begin{multline*}
\int_{t_1}^{\infty}t^b(f(t))^{1/2}dt\lesssim \sum_{k=n}^\infty\int_{\alpha^k}^{\alpha^{k+1}}t^b(f(t))^{1/2}dt\lesssim \sum_{k=n}^\infty(\alpha^{b-a/2+1})^k = \sum_{k=n}^\infty 2^{-k}\lesssim 2^{-n} \\
\approx t_1^{b-a/2+1}.
\end{multline*}
If $b-a/2+1=0$, then choose $n$ so that $t_2\approx 2^n$. By a similar technique, one obtains
\begin{multline*}
\int_1^{t_2}t^b(f(t))^{1/2}dt\lesssim \sum_{k=1}^n\int_{2^k}^{2^{k+1}}t^b(f(t))^{1/2}dt\lesssim \sum_{k=1}^n\left((2^{k+1})^{2b+1}\right)^{1/2}((2^k)^{-a+1})^{1/2} \\
\lesssim\sum_{k=1}^n1=n\approx\log(t_2).
\end{multline*}
\end{proof}

\section{Theorem: Global boundedness and decay for solutions to the semilinear wave equation with null structure on the Minkowski background}\label{m:main_sec}

We conclude this chapter with a theorem about wave dynamics for a simple nonlinear problem in Minkowski space.

\subsection{The null condition}

One key facet of the problem studied in Theorem \ref{m:main_thm} is a particular structural requirement for the nonlinear term, which is often referred to as the \textit{null condition}. Without this condition, in 3+1 dimensions, the problem is actually borderline unstable. \cite{john1981blow}

\begin{definition}
Suppose $\Box_g\psi=\mathcal{N}$, where $\mathcal{N}$ is nonlinear in first derivatives of $\psi$. We say that $\mathcal{N}$ satisfies the \textit{null condition} if it can be written as
$$\mathcal{N}=D\psi \bar{D}\psi,$$
where the operator families $D$ and $\bar{D}$ are as defined in Definition \ref{m:DbarD_def}.
\end{definition}

Since we will also be concerned with equations for the higher order wavefunctions $\psi^s$, we will need the following lemma, which explains a consequence of the null condition.

\begin{lemma}\label{m:null_cond_lem}
Recall that if $\Box_g\psi=D\psi\bar{D}\psi$, then $\Box_g\psi^s=\Gamma^s(D\psi\bar{D}\psi)$. We have the following estimate for this higher order nonlinear term.
$$|\Gamma^s(D\psi \bar{D}\psi)|\lesssim |D\psi^s\bar{D}\psi^{s/2}|+|\bar{D}\psi^s D\psi^{s/2}|.$$
\end{lemma}
\begin{proof}
Using the Leibniz rule for differential operators,
$$\Gamma^s(D\psi\bar{D}\psi) \approx \sum_{i+j=s} \Gamma^i(D\psi)\Gamma^j(\bar{D}\psi) \approx \sum_{i+j= s} D\psi^i\bar{D}\psi^j.$$
Choose a particular term in the sum. If $i\le j$, then
$$|D\psi^i\bar{D}\psi^j|\lesssim |\bar{D}\psi^s D\psi^{s/2}|,$$
while if $i\ge j$, then
$$|D\psi^i\bar{D}\psi^j|\lesssim |D\psi^s \bar{D}\psi^{s/2}|.$$
\end{proof}
\begin{remark}
It should be clear from the above proof that if $s$ is odd, then the quantity $s/2$ can be treated as the smaller integer $\lfloor s/2\rfloor$.
\end{remark}

\subsection{The main theorem}

Finally, we are prepared to state and prove the main theorem of this chapter.
\begin{theorem}\label{m:main_thm}
Let $g$ be the Minkowski metric, and suppose a function $\psi$ solves an equation of the form
$$\Box_g\psi=D\psi \bar{D}\psi.$$
Then for $\delp,\delm>0$ sufficiently small, if the initial data on $\Sigma_0$ decay sufficiently fast as $r\rightarrow\infty$ and have size
$$I_0=E^{6}(0)+E_{2-\delp}^{6}(0)$$
sufficiently small, then the following estimates hold for $t\ge 0$ (with $T=1+t$).

I) The energies satisfy
$$E^{6}(t)\lesssim I_0,$$
$$E_{p\in[\delm,2-\delp]}^{6}(t)\lesssim T^{p-2+\delp}I_0,$$
$$\int_t^\infty E_{p\in [\delm-1,\delm]}^{6}(\tau)d\tau\lesssim T^{p-2+\delm+1}I_0.$$

II) The following $L^\infty$ estimates hold for the operator families $\bar{D}$ and $D$ defined in \S\ref{m:pointwise_summary_sec} (with $\rho=1+r$) provided $s\le 6$.
$$|\rho^{p+1}\bar{D}\psi^s|^2+|\rho^pD\psi^s|^2+|\rho^p\psi^s|^2\lesssim E_{2p}^{s+3}(t),$$
$$|\rho D\psi^s|^2\lesssim E^{s+3}(t).$$
(These estimates hold for all $p$, but in the context of the theorem, they are only useful for $p\in[\delm-1,2-\delp]$.)

III) Together, (I) and (II) imply that for all $s\le 6$, for $p\in[\delm/2,(2-\delp)/2]$,
$$|\rho^{p+1}\bar{D}\psi^s|+|\rho^pD\psi^s|+|\rho^p\psi^s|\lesssim T^{(2p-2+\delp)/2}I_0^{1/2},$$
and additionally for $p\in[(\delm-1)/2,\delm/2]$,
$$\int_t^\infty |\rho^{p+1}\bar{D}\psi^s|+|\rho^pD\psi^s|+|\rho^p\psi^s|\lesssim T^{(2p-2+\delp)/2+1}I_0^{1/2}.$$
The final estimate should be interpreted as saying that $|\rho^{(\delm+1)/2}\bar{D}\psi^s|$, $|\rho^{(\delm-1)/2}D\psi^s|$, and $|\rho^{(\delm-1)/2}\psi^s|$ decay like $T^{(\delm-3+\delp)/2}$ in a weak sense.
\end{theorem}

The remainder of this section is devoted to proving Theorem \ref{m:main_thm}.

\subsection{Boostrap assumptions}\label{m:ba_sec}

We begin the proof of Theorem \ref{m:main_thm} by making the following bootstrap assumptions.
\begin{align*}
E^6(t) &\le C_bI_0, \\
E_{p\in[\delm,2-\delp]}^6(t) &\le C_bT^{p-2+\delp}I_0, \\
\int_t^\infty E_{p\in[\delm-1,\delm]}^6(\tau)d\tau &\le C_bT^{p-2+\delp+1}I_0.
\end{align*}
Note that these bootstrap assumptions are consistent with the general principle that $E_p^s(t)\sim T^{p-2+\delp}$, which the reader should keep in mind throughout the proof.

\subsection{Improved $L^\infty$ estimates}

The $L^\infty$ estimates from Proposition \ref{m:infty_prop} are essential to the proof of the main theorem. But for better clarity, we first remove the nonlinear (error) terms from these estimates and summarize them in the following lemma.
\begin{lemma}\label{m:improved_infty_lem}
In the context of the bootstrap assumptions provided in \S\ref{m:ba_sec}, the following $L^\infty$ estimates hold for $s\le 6$ and all $p$ in any bounded range.
$$|\rho^{p+1}\bar{D}\psi^s|^2+|\rho^pD\psi^s|^2+|\rho^p\psi^s|^2 \lesssim E_{2p}^{s+3}(t)$$
$$|\rho D\psi^s|^2 \lesssim E^{s+3}(t)$$
These are the same as the estimates as from Proposition \ref{m:infty_prop}, except that the nonlinear (error) terms have been removed.
\end{lemma}

\begin{proof}
We employ an inner bootstrap to remove nonlinear (error) terms from the estimates in Proposition \ref{m:infty_prop}. For this part of the proof, there is no need to appeal to null structure. The estimates that follow in this proof are unusually easy compared to the later estimates in the proof of the main theorem that address $N_p^s(t)$ and $N^s(t)$.

Suppose the following (inner) bootstrap assumptions for all $s\le 6$ and all $p$ in a bounded range.
\begin{align*}
|\rho^{p+1} D\psi^s|^2+|\rho^p D\psi^s|^2+|\rho^p\psi^s|^2 &\le C_b E_{2p}^{s+3}(t), \\
|\rho D\psi^s|^2 &\le C_b E^{s+3}(t).
\end{align*}
Then
\begin{align*}
\int_{\Sigma_t}\rho^{2p}(\Box_g\psi^{s+2})^2 &\lesssim \int_{\Sigma_t}\rho^{2p-2}| D\psi^{s+2}|^2|\rho D\psi^{s/2}|^2 \\
&\lesssim \int_{\Sigma_t}\rho^{2p-2}| D\psi^{s+2}|^2C_b E^{s/2+3}(t) \\
&\lesssim E_{2p}^{s+2}(t)C_bE^{s/2+3}(t) \\
&\lesssim C_b^2I_0E_{2p}^{s+2}(t),
\end{align*}
and
\begin{align*}
\int_{\Sigma_t}(\Box_g\psi^{s+2})^2 &\lesssim \int_{\Sigma_t}| D\psi^{s+2}|^2| D\psi^{s/2}|^2 \\
&\lesssim \int_{\Sigma_t}| D\psi^{s+2}|^2C_bE^{s/2+3}(t) \\
&\lesssim E^{s+2}(t) C_bE^{s/2+3}(t) \\
&\lesssim C_b^2I_0E^{s+2}(t).
\end{align*}
From Proposition \ref{m:infty_prop} and the above estimates, for all $s\le 6$,
\begin{align*}
|\rho^{p+1} D\psi^s|^2+|\rho^p D\psi^s|^2+|\rho^p\psi^s|^2 &\lesssim E_{2p}^{s+3}+\int_{\Sigma_t}\rho^{2p}(\Box_g\psi^{s+2})^2 \lesssim (1+C_b^2I_0)E_{2p}^{s+3}(t), \\
|\rho D\psi^s|^2 &\lesssim E^{s+3}+\int_{\Sigma_t}(\Box_g\psi^{s+2})^2\lesssim (1+C_b^2I_0)E^{s+3}(t).
\end{align*}
To see how to close the inner bootstrap assumptions, let us examine the first estimate, which says that
$$|\rho^{p+1} D\psi^s|^2+|\rho^p D\psi^s|^2+|\rho^p\psi^s|^2 \le C(1+C_b^2I_0)E_{2p}^{s+3}(t),$$
where $C$ is a universal constant (not depending on $I_0$ or $C_b$). We see that the first inner bootstrap assumption will be improved if
$$C(1+C_b^2I_0)\le C_b/2.$$
This is achievable if $C_b$ is sufficiently large and $I_0$ is sufficiently small. An identical argument applies for the second estimate.
\end{proof}
The conclusion of this lemma is the same as the statement of part (II) of the main theorem.

\begin{remark}
Lemma \ref{m:improved_infty_lem}, which is a simplified version of Proposition \ref{m:infty_prop} in the sense that there are no nonlinear (error) terms on the right side of any of the estimates in Lemma \ref{m:improved_infty_lem}, will be used regularly throughout the remainder of the proof of the main theorem as a replacement for Proposition \ref{m:infty_prop}.
\end{remark}

\subsection{Refined estimates for $N^s(t)$ and $N_p^s(t)$ (key step)}

The $L^\infty$ estimates in Lemma \ref{m:improved_infty_lem} allow us to provide refined estimates for the nonlinear error terms. \textbf{This is the key step of the proof.}

\begin{lemma}\label{m:refined_nl_lem}
In the context of the bootstrap assumptions provided in \S\ref{m:ba_sec}, if $s\le 6$, then
\begin{align*}
N^s(t)&\lesssim (E^s(t))^{1/2}\left((E^s(t))^{1/2}(E^{s/2+3}_{\delm-1}(t))^{1/2}+(E^s_{1-\delm}(t))^{1/2}(E^{s/2+3}_{\delm-1}(t))^{1/2}\right), \\
N^s_p(t)&\lesssim E^s(t)B^{s/2+3}_p(t)+B^s_p(t)E^{s/2+3}(t).
\end{align*}
\end{lemma}
\begin{proof}
First, we recall the definitions of $N^s(t)$ and $N_p^s(t)$ from Proposition \ref{m:dynamic_estimates_s_prop}.
$$N^s(t)=\sum_{s'\le s}(E^s(t))^{1/2}||\Gamma^{s'}(\Box_g\psi)||_{L^2(\Sigma_t)}$$
$$N_p^s(t)=\sum_{s'\le s}\int_{\Sigma_t}\rho^{p+1}(\Gamma^{s'}(\Box_g\psi))^2,$$
Therefore, it suffices to prove the following two estimates.
$$||\Gamma^s(\Box_g\psi)||_{L^2(\Sigma_t)}\lesssim (E^s(t))^{1/2}(E_{\delm-1}^{s/2+3}(t))^{1/2}+(E_{1-\delm}^s(t))^{1/2}(E_{\delm-1}^{s/2+3}(t))^{1/2},$$
$$\int_{\Sigma_t}\rho^{p+1}(\Gamma^s(\Box_g\psi))^2 \lesssim E^s(t)B^{s/2+3}_p(t)+B^s_p(t)E^{s/2+3}(t).$$
Recall from Lemma \ref{m:null_cond_lem} that
$$|\Gamma^s( D\psi\bar{D}\psi)|\lesssim |D\psi^s\bar{D}\psi^{s/2}|+|\bar{D}\psi^s D\psi^{s/2}|.$$
The following is a proof of the first estimate.
\begin{align*}
||\Gamma^{\le s}\Box_g\psi||_{L^2(\Sigma_t)}^2 &\lesssim \int_{\Sigma_t}(D\psi^s)^2(\bar{D}\psi^{s/2})^2+\int_{\Sigma_t}(\bar{D}\psi^s)^2(D\psi^{s/2})^2 \\
&\lesssim \int_{\Sigma_t}(D\psi^s)^2(\bar{D}\psi^{s/2})^2+\int_{\Sigma_t}\rho^{1-\delm}(\bar{D}\psi^s)^2(\rho^{(\delm-1)/2}D\psi^{s/2})^2 \\
&\lesssim \int_{\Sigma_t}(D\psi^s)^2||\bar{D}\psi^{s/2}||^2_{L^\infty(\Sigma_t)}+\int_{\Sigma_t}\rho^{1-\delm}(\bar{D}\psi^s)^2||\rho^{(\delm-1)/2}D\psi^{s/2}||_{L^\infty(\Sigma_t)}^2 \\
&\lesssim E^s(t)E_{\delm-1}^{s/2+3}(t)+E_{1-\delm}^s(t)E_{\delm-1}^{s/2+3}(t).
\end{align*}
The following is a proof of the second estimate.
\begin{align*}
\int_{\Sigma_t}\rho^{p+1}(\Box_g\psi^s)^2 &\lesssim \int_{\Sigma_t}\rho^{p+1}(D\psi^s)(\bar{D}\psi^{s/2})^2+\int_{\Sigma_t}\rho^{p+1}(\bar{D}\psi^s)(D\psi^{s/2})^2 \\
&\lesssim \int_{\Sigma_t}(D\psi^s)(\rho^{(p+1)/2}\bar{D}\psi^{s/2})^2+\int_{\Sigma_t}\rho^{p-1}(\bar{D}\psi^s)(\rho D\psi^{s/2})^2 \\
&\lesssim \int_{\Sigma_t}(D\psi^s)||\rho^{(p+1)/2}\bar{D}\psi^{s/2}||_{L^\infty(\Sigma_t)}^2 \\
&\hspace{1.4in}+\int_{\Sigma_t}\rho^{p-1}(\bar{D}\psi^s)||\rho D\psi^{s/2}||_{L^\infty(\Sigma_t)}^2 \\
&\lesssim E^s(t)B^{s/2+3}_p(t)+B^s_p(t)E^{s/2+3}(t).
\end{align*}
Note that in the final step, we used $B_p(t)$ in place of $E_{p-1}(t)$, since these norms are equivalent.
\end{proof}

\begin{corollary}\label{m:NL_absorb_bulk}
In the context of the bootstrap assumptions provided in \S\ref{m:ba_sec}, Proposition \ref{m:dynamic_estimates_s_prop} and Lemma \ref{m:refined_nl_lem} imply that if $C_bI_0$ is sufficiently small, then
$$E_p^s(t_2)+\int_{t_1}^{t_2}B_p^s(t)dt\lesssim E_p^s(t_1)$$
whenever $s\le 6$ and $s/2+3\le s$ (which means $6\le s$).
\end{corollary}
\begin{proof}
By Proposition \ref{m:dynamic_estimates_s_prop},
$$E_p^s(t_2)+\int_{t_1}^{t_2}B_p^s(t)dt\lesssim E_p^s(t_1)+\int_{t_1}^{t_2}N_p^s(t)dt.$$
By Lemma \ref{m:refined_nl_lem} and the bootstrap assumptions, if $s\le 6$ and $s/2+3\le s$,
\begin{align*}
\int_{t_1}^{t_2}N_p^s(t) &\lesssim \int_{t_1}^{t_2}E^s(t)B_p^{s/2+3}(t)+B_p^s(t)E^{s/2+3}(t)dt \\
&\lesssim \int_{t_1}^{t_2}E^s(t)B_p^s(t) \\
&\lesssim C_bI_0\int_{t_1}^{t_2}B_p^s(t)dt.
\end{align*}
It follows that if $C_bI_0$ is sufficiently small, then this term can be absorbed into the bulk term on the left side.
\end{proof}

\subsection{Recovering boundedness of $E^6(t)$}\label{m:recover_energy_sec}

Now, we recover the first bootstrap assumption. Let $s=6$. By Proposition \ref{m:dynamic_estimates_s_prop} and Lemma \ref{m:refined_nl_lem},
\begin{align*}
E^s(t) &\lesssim E^s(0)+\int_0^t N^s(\tau)d\tau \\
&\lesssim E^s(0)+\int_0^t(E^s(t))^{1/2}\left((E^s(t))^{1/2}(E^{s/2+3}_{\delm-1}(t))^{1/2}+(E^s_{1-\delm}(t))^{1/2}(E^{s/2+3}_{\delm-1}(t))^{1/2}\right)d\tau \\
&\lesssim I_0+\int_0^t(C_bI_0)^{1/2}(C_bI_0)^{1/2}(C_bT^{\delm-3+\delp}I_0)^{1/2}d\tau \\
&\lesssim I_0+(C_bI_0)^{3/2} \\
&\lesssim (1+C_b^{3/2}I_0^{1/2})I_0
\end{align*}
\textbf{In particular, we used the weak decay principle (Lemma \ref{weak_decay_lem}) in the third step.}

It follows that if $C_b^{3/2}I_0^{1/2}$ is sufficiently small,
$$E^s(t)\lesssim I_0.$$
This recovers the first bootstrap assumption.

\subsection{Recovering boundedness/decay for $E_{p\in[\delm,2-\delp]}^6(t)$}\label{m:recover_decay_sec}

Now, we prove the decay estimates for the energy norms $E_p^6(t)$, for all $p\in[\delm,2-\delp]$. Set $s=6$. We start by observing that a direct consequence of Corollary \ref{m:NL_absorb_bulk} is
$$E_{2-\delp}^s(t)\lesssim E_{2-\delp}^s(0)\lesssim I_0.$$

Then, we invoke Lemma \ref{m:decay_mechanism_lem} to conclude that
$$E_{1-\delp}^s(t) \lesssim T^{-1}E_{2-\delp}^s(t/2)\lesssim T^{(1-\delp)-2+\delp}I_0.$$
Then, we interpolate (Lemma \ref{m:interpolation_lem}) to conclude for all $p\in [1-\delp,2-\delp]$ that
$$E_{p\in [1-\delp,2-\delp]}^s(t)\lesssim T^{p-2+\delp}I_0.$$
Then, we invoke Lemma \ref{m:decay_mechanism_lem} for each $p\in [\delm,\delm+1]$ to conclude for the full range $p\in [\delm,2-\delp]$ that
$$E_{p\in [\delm,2-\delp]}^s(t)\lesssim T^{p-2+\delp}I_0.$$
In particular, this recovers the second bootstrap assumption.

\subsection{Recovering weak decay for $E_{p\in [\delm-1,\delm]}^6(t)$}

Finally, we observe that since $E_p^s(t)\lesssim B_{p+1}^s(t)$ and as a direct consequence of Corollary \ref{m:NL_absorb_bulk}, for $s=6$,
\begin{align*}
\int_t^\infty E_p^s(\tau)d\tau &\lesssim \int_t^\infty B_{p+1}^s(\tau)d\tau \\
&\lesssim E_{p+1}^s(t) \\
&\lesssim T^{(p+1)-2+\delm}I_0 \\
&\lesssim T^{p-2+\delm+1}I_0.
\end{align*}
This recovers the final bootstrap assumption and completes part I of the statement of Theorem \ref{m:main_thm}.

\noindent This completes the proof of Theorem \ref{m:main_thm}. \qed

\chapter{Waves in the Schwarzschild spacetime}\label{szd_chap}

The second problem of this thesis investigates semilinear waves on the exterior of a Schwarzschild black hole. A Schwarzschild black hole is the simplest black hole solution to the Einstein Vacuum Equations. It corresponds to a non-rotating black hole and belongs to the more general Kerr black hole family, which also includes rotating black holes.

The metric for a Schwarzschild black hole with mass $M$ is most well-known in Boyer-Lindquist coordinates.
$$g_{\mu\nu}dx^\mu dx^\nu = -\left(1-\frac{2M}r\right)dt^2+\left(1-\frac{2M}r\right)^{-1}dr^2+r^2(d\theta^2+\sin^2\theta d\phi^2).$$
Its volume form is
$$\mu = r^2\sin\theta$$
and the Lagrangian for the linear wave equation in Boyer-Lindquist coordinates is
\begin{align*}
\mathcal{L} &=\mu g^{\alpha\beta}\pd_\alpha\psi\pd_\beta\psi \\
&= r^2\sin\theta\left[-\left(1-\frac{2M}r\right)^{-1}(\pd_t\psi)^2+\left(1-\frac{2M}r\right)(\pd_r\psi)^2+\sla{g}^{\alpha\beta}\pd_\alpha\psi\pd_\beta\psi\right] \\
&= r^2\sin\theta\left[-\left(1-\frac{2M}r\right)^{-1}(\pd_t\psi)^2+\left(1-\frac{2M}r\right)(\pd_r\psi)^2+r^{-2}(\pd_\theta\psi)^2\right. \\
&\hspace{4in}\left.\vphantom{\left(1-\frac{2M}r\right)^{-1}}+r^{-2}\sin^{-2}\theta(\pd_\phi\psi)^2\right].
\end{align*}

The event horizon (the interface between the interior and exterior of the black hole) is located at $r=2M$. At this radius, the Boyer-Lindquist coordinates break down and a different coordinate system is needed. The new coordinate system consists of a modified time coordinate that differs from the Boyer-Lindquist time coordinate only in the neighborhood $r\in [2M,2M+\delh]$ for some small constant $\delh>0$. In the new coordinate system, the Lagrangian takes the form
$$\mathcal{L}= r^2\sin\theta\left[ g^{tt}(\pd_t\psi)^2+2g^{tr}\pd_t\psi\pd_r\psi+\left(1-\frac{2M}r\right)(\pd_r\psi)^2+\sla{g}^{\alpha\beta}\pd_\alpha\psi\pd_\beta\psi\right],$$
where $g^{tt}$ is a bounded, negative function coinciding with $-\left(1-\frac{2M}r\right)^{-1}$ for $r>r+\delh$ and $g^{tr}$ is a positive function supported on the interval $r\in [2M,2M+\delh]$. For more information about this coordinate system, see Lemma \ref{s:foliation_lem}.

\begin{figure}
\centering
\includegraphics[scale=0.7]{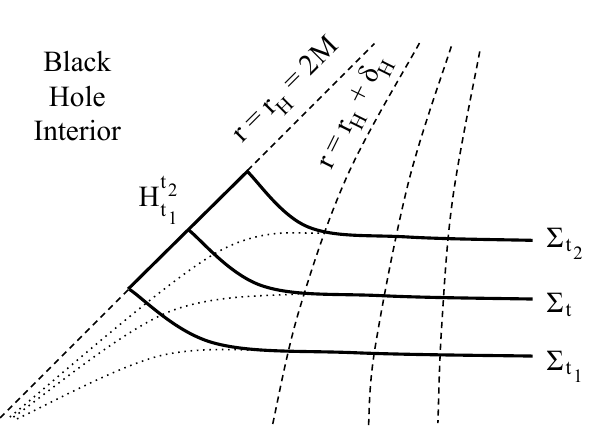}
\caption{The modified time foliation $\Sigma_t$ enters the event horizon and coincides with the Boyer-Lindquist foliation for $r>r_H+\delh$.}\label{s:foliation_fig}
\end{figure}

While hypersurfaces of constant Boyer-Lindquist time all approach the \textit{bifurcation sphere}, the hypersurfaces $\Sigma_t$ of constant time in the new coordinate system meet the event horizon. See Figure \ref{s:foliation_fig}.

We will often use the quantity $\alpha=1-\frac{2M}r$, and in particular, we define the vectorfield $L$ to be
$$L=\alpha\pd_r+\pd_t.$$

\section{The energy estimate and the $h\pd_t$ estimate}\label{s:energy_sec}

The most fundamental of all estimates is the energy estimate. We begin by proving this estimate and then prove the slightly more general $h\pd_t$ estimate.

\subsection{The energy estimate}

Since the vectorfield $\pd_t$ is a Killing vectorfield, we have the following identity.
\begin{lemma}\label{s:ee_identity_lem}(Energy identity for Schwarzschild)
\begin{equation*}
\int_{H_{t_1}^{t_2}}J^r[\pd_t] +\int_{\Sigma_{t_2}}-J^t[\pd_t] = \int_{\Sigma_{t_1}}-J^t[\pd_t]+\int_{t_1}^{t_2}\int_{\Sigma_t}-2\pd_t\psi\Box_g\psi
\end{equation*}
In particular, if $\Box_g\psi=0$, then the quantity
$$E(t)=\int_{H_{t_0}^t}J^r[\pd_t]+\int_{\Sigma_t}-J^t[\pd_t]$$
is conserved.
\end{lemma}
\begin{proof}
By Proposition \ref{general_divergence_estimate_prop}, we have
$$\int_{H_{t_1}^{t_2}}J^r[\pd_t]+\int_{\Sigma_{t_2}}-J^t[\pd_t]+\int_{t_1}^{t_2}\int_{\Sigma_t}divJ[\pd_t]=\int_{\Sigma_{t_1}}-J^t[\pd_t].$$
By Lemma \ref{divJ_lem},
$$divJ[X] = (g^{\mu\lambda}\pd_\lambda X^\nu+g^{\nu\lambda}\pd_\lambda X^\mu-X^\lambda\pd_\lambda g^{\mu\nu}-divX g^{\mu\nu})\pd_\mu\psi\pd_\nu\psi+2X(\psi)\Box_g\psi.$$
For the particular case $X=\pd_t$, the coefficients of $X$ are constant, the metric does not depend on $t$, and by a simple calculation (see equation (\ref{divergence_eqn})), $X$ is divergence free. Therefore,
$$divJ[\pd_t]=2\pd_t\psi\Box_g\psi.$$
We conclude that
$$\int_{H_{t_1}^{t_2}}J^r[\pd_t]+\int_{\Sigma_{t_2}}-J^t[\pd_t]+\int_{t_1}^{t_2}\int_{\Sigma_t}2\pd_t\psi\Box_g\psi=\int_{\Sigma_{t_1}}-J^t[\pd_t].$$
The statement of the lemma now follows.
\end{proof}

Now, we calculate the flux terms in the previous lemma.

\begin{lemma}\label{s:ee_bndry_lem}
On the event horizon $H_{t_1}^{t_2}$,
$$J^r[\pd_t]\approx (\pd_t\psi)^2$$
and on a constant-time hypersurface $\Sigma_t$,
$$-J^t[\pd_t]\approx \chi_H(\pd_r\psi)^2+(\pd_t\psi)^2+|\sla\nabla\psi|^2,$$
where $\chi_H=1-\frac{2M}r$.
\end{lemma}
\begin{proof}
According to Definition \ref{current_def},
$$J^\mu[X]=2g^{\mu\lambda}\pd_\lambda\psi X^\nu\pd_\nu\psi-X^\mu\pd^\lambda\psi\pd_\lambda\psi.$$
Therefore, on the event horizon $H_{t_1}^{t_2}$, since $g^{rr}$ vanishes,
$$J^r[\pd_t]=2g^{rt}(\pd_t\psi)^2.$$
Since $g^{rt}$ is constant and positive on the event horizon,
$$J^r[\pd_t]\approx (\pd_t\psi)^2.$$
Also, on a constant-time hypersurface $\Sigma_t$,
\begin{align*}
J^t[\pd_t] &= 2g^{tt}(\pd_t\psi)^2+2g^{tr}\pd_r\psi\pd_t\psi
-\left(g^{tt}(\pd_t\psi)^2+2g^{tr}\pd_r\psi\pd_t\psi+g^{rr}(\pd_r\psi)^2+|\sla\nabla\psi|^2\right) \\
&= g^{tt}(\pd_t\psi)^2-g^{rr}(\pd_r\psi)^2-|\sla\nabla\psi|^2.
\end{align*}
Since $-g^{tt}\approx 1$ and $g^{rr}=1-\frac{2M}r=\chi_H$,
$$-J^t[\pd_t]\approx\chi_H(\pd_r\psi)^2+(\pd_t\psi)^2+|\sla\nabla\psi|^2.$$
This concludes the proof of the lemma.
\end{proof}

From the previous two lemmas, we conclude the following energy estimate.
\begin{proposition}\label{s:classic_ee_prop}(Energy estimate for Schwarzschild)
\begin{multline*}
\int_{H_{t_1}^{t_2}}(\pd_t\psi)^2+\int_{\Sigma_{t_2}}\left[\chi_H(\pd_r\psi)^2+(\pd_t\psi)^2+|\sla\nabla\psi|^2\right] \\
 \lesssim \int_{\Sigma_{t_1}}\left[\chi_H(\pd_r\psi)^2+(\pd_t\psi)^2+|\sla\nabla\psi|^2\right] + Err_\Box,
\end{multline*}
where $\chi_H=1-\frac{2M}r$ and
$$Err_\Box=\int_{t_1}^{t_2}\int_{\Sigma_t}|\pd_t\psi\Box_g\psi|.$$
\end{proposition}
\begin{proof}
This follows directly from Lemmas \ref{s:ee_identity_lem} and \ref{s:ee_bndry_lem}.
\end{proof}

\subsection{The $h\pd_t$ energy estimate}\label{s:hdt_sec}

We proceed as in \S\ref{m:hdt_sec} to prove the $h\pd_t$ estimate.
\begin{lemma}\label{s:divJhdt_lem}
If $h=h(r)$ is constant in the interval $r\in[r_H,r_H+\delh]$, and $\Box_g\psi=0$, then
$$divJ[h\pd_t]=\frac{h'}{2}\left[(L\psi)^2-(\lbar\psi)^2\right].$$
\end{lemma}
\begin{proof}
Recall from Lemma \ref{divJ_lem} that
$$divJ[X]=K^{\mu\nu}\pd_\mu\psi\pd_\nu\psi,$$
where
$$K^{\mu\nu}=g^{\mu\lambda}\pd_\lambda X^\nu+g^{\nu\lambda}\pd_\lambda X^\mu-X^\lambda\pd_\lambda(g^{\mu\nu})-divXg^{\mu\nu}.$$
Since $\pd_t$ is killing, $\pd_t(g^{\mu\nu})=0$ and $div(\pd_t)=0$. Also, the only component of $X$ is the $t$ component and the only nonzero derivative of that component is the $\pd_r$ derivative, so
$$g^{\mu\lambda}\pd_\lambda X^\nu+g^{\nu\lambda}\pd_\lambda X^\mu=g^{\mu r}\pd_r X^\nu+g^{\nu r}\pd_r X^\mu.$$
It follows that the only possible nonzero $K^{\mu\nu}$ components are
\begin{align*}
K^{tt} &= 2g^{tr}h' \\
K^{tr}+K^{rt} &= 2g^{rr}h'.
\end{align*}
Since $h'=0$ in the region $r\in [r_H,r_H+\delh]$, the first of these actually vanishes. We conclude that
\begin{align*}
divJ[h\pd_t]&=2\alpha h'\pd_r\psi\pd_t\psi \\
&=\frac{h'}2\left[(\alpha\pd_r\psi+\pd_t\psi)^2-(\alpha\pd_r\psi-\pd_t\psi)^2\right] \\
&=\frac{h'}2\left[(L\psi)^2-(\lbar\psi)^2\right].
\end{align*}
\end{proof}

Taking $h$ to be a positive function decreasing to zero as $r\rightarrow\infty$ at a particular rate, we obtain the $h\pd_t$ estimate.
\begin{proposition}\label{s:hdt_prop}($h\pd_t$ estimate for Schwarzschild)
Let $R>r_H+\delh$ be any given radius. Then for all $\epsilon>0$ and $p<2$, there is a small constant $c_\epsilon$ and a large constant $C_\epsilon$, such that
\begin{multline*}
\int_{H_{t_1}^{t_2}}(\pd_t\psi)^2+\int_{\Sigma_{t_2}}r^{p-2}\left[\chi_H(\pd_r\psi)^2+(\pd_t\psi)^2+|\sla\nabla\psi|^2\right] \\
+\int_{t_1}^{t_2}\int_{\Sigma_t\cap\{R+M<r\}}c_\epsilon r^{p-3}(\lbar\psi)^2 \\
\lesssim \int_{\Sigma_{t_1}}r^{p-2}\left[\chi_H(\pd_r\psi)^2+(\pd_t\psi)^2+|\sla\nabla\psi|^2\right]+Err,
\end{multline*}
where $\chi_H=1-\frac{2M}r$ and
\begin{align*}
Err&=Err_1+Err_\Box \\
Err_1&=\int_{t_1}^{t_2}\int_{\Sigma_t\cap\{R<r\}}\epsilon r^{-1}(L\psi)^2 \\
Err_\Box&=\int_{t_1}^{t_2}\int_{\Sigma_t}C_\epsilon r^{p-2}|\pd_t\psi\Box_g\psi|.
\end{align*}
\end{proposition}
\begin{proof}
See the proof of Proposition \ref{m:hdt_prop}, noting Lemmas \ref{s:divJhdt_lem} and \ref{s:ee_bndry_lem}.
\end{proof}

\section{The Morawetz estimate}\label{s:morawetz_sec}

In this section, we prove the Morawetz estimate for Schwarzschild. It is similar to version II of the Morawetz estimate for Minkowski, but its proof is significantly more complicated. In particular, there will be a serious and unavoidable complication related to the presence of trapped null geodesics.

\subsection{Null geodesics in Schwarzschild}\label{s:geodesics_sec}

As will become evident in the Morawetz estimate, the dynamics of waves in Schwarzschild are related to the dynamics of null geodesics. In particular, the existence of trapped null geodesics (null geodesics that do not leave a sphere of constant radius) at the radius $r_{trap}=3M$ is an obstacle to a Morawetz estimate without a loss of certain derivatives at the trapped radius. \cite{Ralston}\cite{Sbierski} Here, we briefly review the dynamics of null geodesics in Schwarzschild.

Due to the number of symmetries in Schwarzschild, which are given by the time translation $\pd_t$ as well as the rotations $\Omega$, the geodesic equation is integrable. That is, there are four independent constants of motion for any geodesic.

\begin{proposition}\label{geodesic_constants_prop}
If
$$\gamma(\lambda)=(t(\lambda),r(\lambda),\theta(\lambda),\phi(\lambda))$$
 is a null geodesic on Schwarzschild, then the following quantities are constant along the trajectory of $\gamma$.
$$E=\left(1-\frac{2M}r\right)\dot{t}$$
$$L_z=r^2\sin^2\theta \dot{\phi}$$
$$L^2=r^4(\dot\theta^2+\sin^2\theta \dot{\phi}^2)$$
And, since $\gamma$ is null,
\begin{equation}\label{geodesic_E_cons_eqn}
E^2=\dot{r}^2+r^{-2}\left(1-\frac{2M}r\right)L^2.
\end{equation}
\end{proposition}
\begin{proof}
If $\vec{V}$ is a vector that generates an isometry of the spacetime, then the quantity $g(\vec{V},\dot{\gamma})$ is conserved. This is because the deformation tensor $\nabla_{(\mu}V_{\nu)}$ vanishes and
\begin{align*}
\frac{d}{d\lambda}g(\vec{V},\dot{\gamma})&=\dot\gamma^\mu\nabla_\mu(V_\nu\dot{\gamma}^\nu) \\
&=\dot\gamma^\nu\dot\gamma^\mu\nabla_{\mu}V_\nu+V_\nu\dot\gamma^\mu\nabla_\mu\dot\gamma^\nu \\
&=\dot\gamma^\nu\dot\gamma^\mu\nabla_{(\mu}V_{\nu)}+0 \\
&=0.
\end{align*}
Since the Schwarzschild metric does not depend on the coordinates $t$ or $\phi$, it follows that the corresponding vectorfields $\pd_t$ and $\pd_\phi$ generate isometries. Therefore, the quantities
$$E=-(\pd_t)^\mu g_{\mu\nu}\dot\gamma^\nu=-g_{tt}\dot{t}=\left(1-\frac{2M}r\right)\dot{t}$$
and
$$L_z=(\pd_\phi)^\mu g_{\mu\nu}\dot\gamma^\nu=g_{\phi\phi}\dot{\phi}=r^2\sin^2\theta\dot{\phi}$$
are conserved.

The vectorfield $\pd_\phi$ is the same as the rotation symmetry operator $\Omega_z$. There are also the rotation symmetry operators $\Omega_x$ and $\Omega_y$, which are less trivially expressed in the spherical coordinate system, but nevertheless also yield conserved quantities
$$L_x=(\Omega_x)^\mu g_{\mu\nu}\dot{\gamma}^\nu,$$
$$L_y=(\Omega_y)^\mu g_{\mu\nu}\dot{\gamma}^\nu.$$
Together with $L_z$, these quantites satsfy the relation
$$L^2=(L_x)^2+(L_y)^2+(L_z)^2=r^4(\dot{\theta}^2+\sin^2\theta\dot{\phi}^2).$$

Finally, all null geodesics satisfy
\begin{align*}
0&=g_{\mu\nu}\dot\gamma^\mu\dot\gamma^\nu \\
&=-\left(1-\frac{2M}r\right)\dot{t}^2+\left(1-\frac{2M}r\right)^{-1}\dot{r}^2+r^2(\dot\theta^2+\sin^2\theta\dot\phi^2) \\
&=-\left(1-\frac{2M}r\right)^{-1}E^2+\left(1-\frac{2M}r\right)^{-1}\dot{r}^2+r^{-2}L^2.
\end{align*}
Therefore,
$$E^2=\dot{r}^2+r^{-2}\left(1-\frac{2M}r\right)L^2.$$
\end{proof}

Differentiating equation (\ref{geodesic_E_cons_eqn}) with respect to the parameter $\lambda$ yields
\begin{equation}\label{geodesic_F_motivation_eqn}
0=2\dot{r}\left[\ddot{r}-r^{-3}\left(1-\frac{3M}r\right)L^2\right],
\end{equation}
clearly suggesting, but not quite proving, the following proposition.
\begin{proposition}
If
$$\gamma(\lambda)=(t(\lambda),r(\lambda),\theta(\lambda),\phi(\lambda))$$
is a null geodesic in Schwarzschild with the conserved squared angular momentum $L^2$ defined in Proposition \ref{geodesic_constants_prop}, then
\begin{equation}\label{geodesic_F_eqn}
\ddot{r}=r^{-3}\left(1-\frac{3M}r\right)L^2.
\end{equation}
\end{proposition}
\begin{proof}
If $\dot{r}\ne 0$, then equation (\ref{geodesic_F_motivation_eqn}) constitutes a valid proof. But for the important case $\dot{r}=0$, we must use a different argument.

A geodesic extremizes the following Lagrangian.
$$\mathcal{L}=g_{\mu\nu}\dot\gamma^\mu\dot\gamma^\nu=-\left(1-\frac{2M}r\right)\dot{t}^2+\left(1-\frac{2M}r\right)^{-1}\dot{r}^2+r^2(\dot{\theta}^2+\sin^2\theta\dot{\phi}^2)$$
According to the Euler-Lagrange equation,
$$\frac{d}{d\lambda}\frac{\pd\mathcal{L}}{\pd\dot{x}^\mu}-\frac{\pd\mathcal{L}}{\pd x^\mu}=0.$$
We substitute $r$ for $x^\mu$ and compute both terms.
\begin{align*}
\frac{d}{d\lambda}\frac{\pd\mathcal{L}}{\pd\dot{r}} &= \frac{d}{d\lambda}\left[2\dot{r}\left(1-\frac{2M}r\right)^{-1}\right] \\
&=2\ddot{r}\left(1-\frac{2M}r\right)^{-1}-2\dot{r}^2\left(1-\frac{2M}r\right)^{-2}\frac{2M}{r^2}.
\end{align*}
\begin{align*}
\frac{\pd\mathcal{L}}{\pd r}=-\frac{2M}{r^2}\dot{t}^2-\left(1-\frac{2M}r\right)^{-2}\frac{2M}{r^2}\dot{r}^2+2r(\dot{\theta}^2+\sin^2\theta\dot{\phi}^2)
\end{align*}
Equating both terms and multiplying by $\frac12\left(1-\frac{2M}r\right)$, we have
\begin{multline*}
\ddot{r}-\left(1-\frac{2M}r\right)^{-1}\frac{2M}{r^2}\dot{r}^2 \\
=-\frac{M}{r^2}\left(1-\frac{2M}r\right)\dot{t}^2-\left(1-\frac{2M}r\right)^{-1}\frac{M}{r^2}\dot{r}^2+r\left(1-\frac{2M}r\right)(\dot{\theta}^2+\sin^2\theta\dot{\phi}^2),
\end{multline*}
which simplifies to
$$\ddot{r}=\frac{M}{r^2}\left[-\left(1-\frac{2M}r\right)\dot{t}^2+\left(1-\frac{2M}r\right)^{-1}\dot{r}^2\right]+r\left(1-\frac{2M}r\right)(\dot{\theta}^2+\sin^2\theta\dot{\phi}^2).$$
The quantity in square brackets appears in the expression $g_{\mu\nu}\dot{\gamma}^\mu\dot{\gamma}^\nu$. Since $g_{\mu\nu}\dot{\gamma}^\mu\dot{\gamma}^\nu=0$ for a null geodesic, the quantity in square brackets can be replaced by $-r^2(\dot{\theta}^2+\sin^2\theta\dot{\phi}^2)$, yielding the equation
$$\ddot{r}=\frac{M}{r^2}\left[-r^2(\dot{\theta}^2+\sin^2\theta\dot{\phi}^2)\right]+r\left(1-\frac{2M}r\right)(\dot{\theta}^2+\sin^2\theta\dot{\phi}^2).$$
Finally, we substitute $\dot{\theta}^2+\sin^2\theta\dot{\phi}^2=r^{-4}L^2$ and simplify to obtain
\begin{align*}
\ddot{r} &= \frac{M}{r^2}(-r^{-2}L^2)+r\left(1-\frac{2M}r\right)(r^{-4}L^2) \\
&= r^{-3}\left(1-\frac{3M}r\right)L^2.
\end{align*}
\end{proof}

\begin{remark}
Although equation (\ref{geodesic_F_motivation_eqn}) does not quite constitute a proof of the formula (\ref{geodesic_F_eqn}), it provides a far more direct explanation for the reason why the quantity on the right side of equation (\ref{geodesic_F_eqn}) is minus one half the derivative with respect to $r$ of the quantity $r^{-2}\left(1-\frac{2M}r\right)L^2$. That is,
$$-\frac12 \frac{d}{dr}\left[r^{-2}\left(1-\frac{2M}r\right)L^2\right]=r^{-3}\left(1-\frac{3M}r\right)L^2.$$
\end{remark}

\begin{corollary}
If
$$\gamma(\lambda)=(t(\lambda),r(\lambda),\theta(\lambda),\phi(\lambda))$$
is a null geodesic in Schwarzschild with the conserved energy $E$ and squared angular momentum $L^2$ defined in Proposition \ref{geodesic_constants_prop}, then
$$E^2=\dot{r}^2+V(r)L^2$$
and
$$\ddot{r}=-\frac12V'(r)L^2,$$
where
$$V(r)=r^{-2}\left(1-\frac{2M}r\right).$$
\end{corollary}

The analogy between the above corollary and a classical particle in a potential $V(r)L^2$ motivates the following definition.
\begin{definition}
The quantity
$$V(r)=r^{-2}\left(1-\frac{2M}r\right)$$
is called the \textit{geodesic potential}.
\end{definition}

\begin{corollary}
There exist null geodesics that are confined to a hypersurface with constant radius $r=3M$. We call these geodesics \textit{trapped geodesics} and we denote $r_{trap}=3M$.
\end{corollary}

\subsection{The partial Morawetz estimate}\label{s:partial_morawetz_sec}

We now proceed as in \S\ref{m:morawetzII_sec} to construct the Morawetz estimate, beginning with the partial Morawetz estimate as in \S\ref{m:partial_morawetz_II_sec}.

\begin{proposition}\label{s:partial_morawetz_prop}(Partial Morawetz estimate)
There exist a vectorfield $X_0$ and a function $w_0$ such that the current $J[X_0,w_0]$ satisfies
$$\left[\frac{M^2}{r^3}\left(1-\frac{r_H}r\right)^2(\pd_r\psi)^2+\frac1r\left(1-\frac{r_{trap}}r\right)^2|\sla\nabla\psi|^2+\frac{M}{r^4}1_{r>4M}\psi^2\right]\lesssim divJ[X_0,w_0]$$
in Boyer-Lindquist coordinates.
\end{proposition}

\begin{proof}
For the proof of this proposition, we drop the $0$ subscript from $X_0$ and $w_0$. We begin by restating Lemma \ref{m:pme_initial_lem} for Schwarzschild.
\begin{lemma}
With the choices $X=X^r(r)\pd_r$ and $w=w(r)$,
\begin{multline*}
divJ[X,w]= 
(wg^{tt}-r^{-2}\pd_r(r^2X^rg^{tt}))(\pd_t\psi)^2+(w-\pd_rX^r)|\sla\nabla\psi|^2 \\
+(wg^{rr}+\pd_rX^r g^{rr}-2r^{-1}X^rg^{rr}-X^r\pd_rg^{rr})(\pd_r\psi)^2 -\frac12r^{-2}\pd_r((r^2-2Mr)\pd_rw)\psi^2.
\end{multline*}
\end{lemma}
\begin{proof}
See the proof of Lemma \ref{m:pme_initial_lem}.
\end{proof}

It is difficult to make sense of the coefficients in the above lemma, because most coefficients are sums of multiple terms. By making the choices $X^r(r)=u(r)v(r)$ and $w(r)=v(r)\pd_ru(r)$, each coefficient can be written as a single term.

\begin{lemma}
With the choices $X=u(r)v(r)\pd_r$ and $w=v(r)\pd_ru(r)$,
\begin{multline*}
divJ[X,w]= 
-ur^{-2}\pd_r(r^2vg^{tt})(\pd_t\psi)^2-u\pd_rv |\sla\nabla\psi|^2 \\
+r^2(g^{rr})^2u^{-1}\pd_r\left(\frac{u^2v}{r^2g^{rr}}\right)(\pd_r\psi)^2-\frac12r^{-2}\pd_r((r^2-2Mr)\pd_rw)\psi^2.
\end{multline*}
\end{lemma}
\begin{proof}
This follows directly by substituting $X^r=uv$ and $w=v\pd_ru$ into the expression in the previous lemma and combining terms.
\end{proof}

At this point, it is standard to choose a relation between $X$ and $w$ so that the coefficient of $(\pd_t\psi)^2$ vanishes. (It will be recovered later--see \S\ref{s:dt_correction_sec}.) In terms of the functions $u$ and $v$, it suffices to set $v=r^{-2}(1-\frac{2M}r)$. 

\begin{remark}
The function $v=r^{-2}(1-\frac{2M}r)$ coincides with the geodesic potential derived in \S\ref{s:geodesics_sec}, which explains why the trapping radius $r_{trap}=3M$ will show up in the calculations that follow.
\end{remark}

\begin{lemma}
With the choice $v=r^{-2}(1-\frac{2M}r)$, the coefficient of the $(\pd_t\psi)^2$ term vanishes and furthermore,
\begin{multline*}
divJ[X,w]=\frac{2u}{r^3}\left(1-\frac{3M}r\right)|\sla\nabla\psi|^2 \\
+2\pd_r\left(\frac{u}{r^2}\right)\left(1-\frac{2M}r\right)^2(\pd_r\psi)^2-\frac12r^{-2}\pd_r((r^2-2Mr)\pd_rw)\psi^2.
\end{multline*}
Therefore, it is necessary to choose $u$ and $w$ so that the following conditions are satisfied.
\begin{eqnarray}
u\left(1-\frac{3M}r\right)\ge 0, \label{s:u_cond_1_eqn}\\
\pd_r\left(\frac{u}{r^2}\right)\ge 0, \label{s:u_cond_2_eqn}\\
\pd_r((r^2-2Mr)\pd_rw)\le 0, \label{s:u_cond_3_eqn}
\end{eqnarray}
Furthermore, $w$ and $u$ must be related by the following constraint.
\begin{equation}
w=r^{-2}\left(1-\frac{2M}r\right)\pd_ru. \label{s:u_constr_eqn}
\end{equation}
\end{lemma}

\begin{lemma}\label{s:choice_for_u_and_w_lem}
It is possible to choose $u$ and $w$ such that all three conditions (\ref{s:u_cond_1_eqn}-\ref{s:u_cond_3_eqn}) and the constraint (\ref{s:u_constr_eqn}) are satisfied. One particular choice is given by the following approach. (See Figure \ref{s:w_fig}.)

i) Require that $u(3M)=0$ and $\pd_ru=r^2\left(1-\frac{2M}r\right)^{-1}w$. This will specify $u$ completely in terms of $w$ and satisfy the constraint (\ref{s:u_constr_eqn}).

ii) Require that $w$ be positive. Then $u$ will be an increasing function and the condition (\ref{s:u_cond_1_eqn}) will be satisfied.

iii) Determine a way to understand $\pd_r\left(\frac{u}{r^2}\right)$ in terms of $w$. Here is one way. Let $\tilde{K}^{rr}=\frac12 r^{3}\pd_r\left(\frac{u}{r^2}\right)$. This particular quantity has the same sign as $\pd_r\left(\frac{u}{r^2}\right)$ and has the property that its derivative can be written as a function of $w$. That is, $\pd_r\tilde{K}^{rr}=r^2\pd_r\left(\frac{w}{2rv}\right)$, where again, $v=r^{-2}(1-\frac{2M}r)$. The condition (\ref{s:u_cond_2_eqn}) now reduces to showing that $\tilde{K}^{rr}\ge 0$.

iv) For sufficiently large $r$, impose the condition that $\pd_r\tilde{K}^{rr}=0$. In particular, this means choosing $w=2rv$ for sufficiently large $r$. The quantity $2rv$ has a maximum at $r_*=4M$. (Again, see Figure \ref{s:w_fig}.) The choice $w=2rv$ for $r\ge r_*$ will satisfy both the remaining conditions (\ref{s:u_cond_2_eqn}-\ref{s:u_cond_3_eqn}) for $r>r_*$.

v) Observe that since $\pd_r(2rv)=0$ at $r_*=4M$, the quantity $(r^2-2Mr)\pd_rw=0$ will necessarily vanish at $r=2M$ and $r=4M$ if $w=2rv$ for $r\ge r_*$. By the mean value theorem, condition (\ref{s:u_cond_3_eqn}) necessitates that $\pd_rw=0$ entirely for $2M\le r\le 4M$. Thus, take $w=w(r_*)$ for $r\le r_*$.

iv) Since $w$ has now been chosen for all $r$, check that $\tilde{K}^{rr}\ge 0$. In particular, with the given choice of $w$, $\tilde{K}^{rr}$ will be decreasing for $r<r_*$ and then remain constant for $r\ge r_*$. Thus, it suffices to compute that $\tilde{K}^{rr}(r_*)>0$.
\end{lemma}

\begin{figure}
\centering
\includegraphics[scale=0.7]{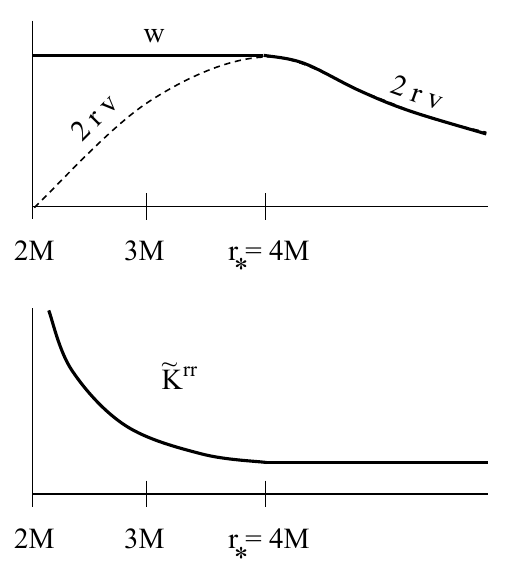}
\caption{Plots of $w$ and $\tilde{K}^{rr}=\frac12r^3\pd_r\left(\frac{u}{r^2}\right)$ for $r\ge 2M$. In the region $r>r_*$, the function $w$ is chosen so that $\tilde{K}^{rr}$ is constant. In the region $r<r_*$, $w$ is chosen to be constant (this is required by condition (\ref{s:u_cond_3_eqn})). The value $r_*=4M$ is chosen so that $w$ remains $C^1$.}\label{s:w_fig}
\end{figure}

\begin{proof}
The lemma consists of multiple statements. Statements (i) and (ii) require no further justification. We turn to statement (iii). Note that
\begin{align*}
\pd_r\tilde{K}^{rr} &= \pd_r\left(\frac12r^3\pd_r\left(\frac{u}{r^2}\right)\right) \\
&= \pd_r\left(\frac12r^3\left(\frac{u'}{r^2}-\frac{2u}{r^3}\right)\right) \\
&= \pd_r\left(\frac{ru'}{2}-u\right) \\
&= \frac{ru''}{2}-\frac{u'}2 \\
&= r^2\pd_r\left(\frac{u'}{2r}\right).
\end{align*}
Recall the constraint (\ref{s:u_constr_eqn}) between $w$ and $u$, which stipulates $w=vu'$. Therefore, we substitute $u'=\frac{w}{v}$ and conclude
$$\pd_r\tilde{K}^{rr} = r^2\pd_r\left(\frac{w}{2rv}\right).$$
This verifies statement (iii). For statement (iv), we must find the value $r_*$ that maximizes the quantity $2rv$. We calculate
\begin{align*}
\pd_r(2rv) &=2\pd_r\left(r^{-1}\left(1-\frac{2M}r\right)\right) \\
&= -\frac{2}{r^2}\left(1-\frac{2M}r\right) +\frac{2}{r}\frac{2M}{r^2} \\
&= -\frac{2}{r^2}\left(1-\frac{4M}r\right).
\end{align*}
Setting $\pd_r(2rv)=0$, we deduce $r_*=4M$. This indeed maximizes the quantity $2rv$, because the quantity $2rv$ vanishes at the event horizon and decays back to zero as $r\rightarrow\infty$.

We also must check that $-\frac12r^{-2}\pd_r\left((r^2-2Mr)\pd_r(2rv)\right)\ge 0$ for $r\ge r_*=4M$. Indeed,
\begin{align*}
-\frac12 r^{-2}\pd_r\left((r^2-2Mr)\pd_r(2rv)\right) &= r^{-2}\pd_r\left(\frac{r^2-2Mr}{r^2}\left(1-\frac{4M}r\right)\right) \\
&= r^{-2}\pd_r\left(\left(1-\frac{2M}r\right)\left(1-\frac{4M}r\right)\right).
\end{align*}
The last quantity is clearly positive for $r\ge 4M$. This verifies statement (iv). Statement (v) requires no further justification. Finally, we turn to statement (vi). First, we observe that for $r<r_*$, since $w=w(r_*)$ is constant, then
$$\pd_r\tilde{K}^{rr}=r^2\pd_r\left(\frac{w(r_*)}{2rv}\right)=-\frac{w(r_*)r^2}{(2rv)^2}\pd_r(2rv).$$
Since we have determined that $2rv$ has a maximum at $r_*=4M$, it follows that $\tilde{K}^{rr}$ decreases to a minimum at $r_*=4M$. Recall further that by the choice of $w$ for $r\ge r_*$, $\tilde{K}^{rr}$ is constant for $r\ge r_*$. It suffices to check the value of $\tilde{K}^{rr}$ at $r=r_*$. (That is, if $\tilde{K}^{rr}(r_*)>0$, then necessarily $\tilde{K}^{rr}>0$ and thus $K^{rr}\ge 0$.) Indeed,
\begin{multline*}
\left.\frac{r^2 v}{w}\pd_r\left(\frac{u}{r^2}\right)\right|_{r_*}
=\left.\frac{r^2v\pd_ru}{wr^2}-\frac{2r^2uv}{wr^3}\right|_{r_*}
=1-\frac{2v(r_*)}{r_*w(r_*)}\int_{r_{trap}}^{r_*}\pd_r u 
=1-\frac{2v(r_*)}{r_*w(r_*)}\int_{r_{trap}}^{r_*}\frac{w}{v} \\
=1-\frac{2v(r_*)}{r_*}\int_{r_{trap}}^{r_*}\frac1v>1-\frac{2(r_*-r_{trap})}{r_*}=1-\frac{2M}{4M}=\frac12,
\end{multline*}
where we used the fact that $w$ is constant for $r\in [r_{trap},r_*]$ and the inequality
$$v(r_*)\int_{r_{trap}}^{r_*}\frac1v<r_*-r_{trap},$$
which follows from the fact that $v$ is decreasing from $r_{trap}$ to $r_*$. 
This concludes the proof of Lemma \ref{s:choice_for_u_and_w_lem}.
\end{proof}

We conclude the proof of Proposition \ref{s:partial_morawetz_prop} by analyzing the asymptotics for large $r$, in particular for $r>r_*$.
\begin{align*}
v &= r^{-2}\left(1-\frac{2M}r\right)=O(r^{-2}) \\
\pd_rv &= O(r^{-3}) \\
w_0 &= 2rv = O(r^{-1}) \\
u &= r^2-c^2 \text{ since }u'=2r \\
X_0 &= uv\pd_r = \left(1-\frac{2M}r\right)\pd_r+O\left(\frac{M^2}{r^2}\right)\pd_r \\
\sla{K}^{\alpha\beta} &= -u\pd_rv\sla{g}^{\alpha\beta}=O(r^{-1})\sla{g}^{\alpha\beta} \\
K^{rr} &=O(r^{-3})\tilde{K}^{rr}=O\left(\frac{M^2}{r^3}\right) \\
K &=r^{-2}\pd_r\left(\left(1-\frac{2M}r\right)\left(1-\frac{4M}r\right)\right) = O\left(\frac{M}{r^4}\right).
\end{align*}
This completes the proof of Proposition \ref{s:partial_morawetz_prop}.
\end{proof}

\subsection{Remaining issues}

So far, we have a current $J[X_0,w_0]$ with a non-negative divergence. However, there are still a few modifications that must be made to the current $J[X_0,w_0]$ to establish a proper Morawetz estimate. \\
\bp The vectorfield $X_0$ behaves like $\log(1-\frac{2M}r)\pd_r$ near the event horizon. So we temper $X$ (and also $w$) by choosing from a family $X_{\epsilon_{temper}}$, $w_{\epsilon_{temper}}$ where the parameter $\epsilon_{temper}>0$, which represents a deviation from $X_0$ and $w_0$, will be small. This will unfortunately introduce an error term of the form $-r^{-2}V_{\epsilon_{temper}}\psi^2$ to the divergence of $J$. (See \S\ref{s:temper_correction_sec}.) \\
\bp The coefficient $K^{rr}$ of $(\pd_r\psi)^2$ vanishes at the event horizon. This can be fixed by using the standard redshift argument, which necessitates a change of foliation near the event horizon and the addition of a redshift vectorfield $\epsilon_{redshift}Y$. This new foliation will be used for the remainder of the chapter. (See \S\ref{s:redshift_correction_sec}.) \\
\bp The coefficient $K$ of $\psi^2$ vanishes for $r\in [r_H,r_*]$. In fact, given the error term introduced by the tempering of $X$ and $w$, it is even slightly negative, with a smallness parameter $\epsilon_{temper}$. This will be fixed by observing a local Hardy estimate. (See \S\ref{s:hardy_correction_sec}.) \\
\bp The coefficient $K^{tt}$ of $(\pd_t\psi)^2$ vanishes entirely. This is easily fixed by adding to $w$ the term $\epsilon_{\pd_t}w_{\pd_t}$ just as in the Minkowski setting. (See \S\ref{s:dt_correction_sec}.) \\
\bp The purpose of the current $J$ is to apply the divergence theorem in the form of Proposition \ref{general_divergence_estimate_prop}. So far, we have paid careful attention to the divergence of $J$, but the boundary terms also must be shown to have good properties. We add $\pd_t$ to $X$ and observe that there are only small errors on the boundary, just as in the Minkowski setting. (See \S\ref{s:boundary_correction_sec}.)

After all the above modifications take place, the final current $J$ will be of the following form.
\begin{align*}
J&=J[X,w]  \\
X&=X_{\epsilon_{temper}}+\epsilon_{redshift}Y+\pd_t\\ 
w&=w_{\epsilon_{temper}}+\epsilon_{\pd_t}w_{\pd_t}
\end{align*}
And the current $J$ will satisfy
\begin{multline*}
\int_{\Sigma_t}\frac{M^2}{r^3}(\pd_r\psi)^2+\frac1r\left(1-\frac{3M}r\right)^2\left[|\sla\nabla\psi|^2+\frac{M^2}{r^2}(\pd_t\psi)^2\right]+\frac{M}{r^4}\psi^2 \\
\lesssim \int_{\Sigma_t}div J[X,w]-(2X(\psi)+w\psi)\Box_g\psi.
\end{multline*}

\subsection{The tempered vectorfield $X_{\epsilon_{temper}}$ and function $w_{\epsilon_{temper}}$}\label{s:temper_correction_sec}

The first issue with the partial Morawetz estimate that we address is the behavior of $X_0$ near the event horizon. Recall that
$$X_0=uv\pd_r,$$
$$w_0=v\pd_ru.$$
Since $w_0$ was constant in a neighborhood of the event horizon and $v$ vanishes to first order at the event horizon, this implies that
$$\pd_ru\sim (r-r_H)^{-1},$$
so
$$u\sim \log(r-r_H).$$
This is problematic when one considers the boundary terms on $\Sigma_{t_1}$ and $\Sigma_{t_2}$ in the complete estimate, and furthermore when one considers the nonlinear problem as the error term $X_0(\psi)\Box_g\psi$ has an unbounded weight.

To solve this problem, we define a family $X_{\epsilon_{temper}},w_{\epsilon_{temper}}$ parametrized by $\epsilon_{temper}$ for which $X_{\epsilon_{temper}}$ is regular if $\epsilon_{temper}>0$. An added benefit is that $w_{\epsilon_{temper}}$ will be zero on a neighborhood of the event horizon. This will simplify calculations of boundary terms near the event horizon.

For the sake of simplicity, we continue to use the Boyer-Lindquist coordinates for Schwarzschild.

The following Lemma is established using the approach from \cite{Ma}.

\begin{lemma}\label{s:temper_lemma}
Suppose $\Box_g\psi=0$. For all $\epsilon_{temper}>0$ sufficiently small, there exist a modified vectorfield $X_{\epsilon_{temper}}$ and a modified function $w_{\epsilon_{temper}}$ agreeing, respectively, with the vectorfield $X_0$ and function $w_0$ from Proposition \ref{s:partial_morawetz_prop} outside of a small neighborhood of the event horizon and satisfying
\begin{multline*}
\left[\frac{M^2}{r^3}\left(1-\frac{r_H}r\right)^2(\pd_r\psi)^2+\frac1r\left(1-\frac{r_{trap}}r\right)^2|\sla\nabla\psi|^2+\frac{M}{r^4}1_{r>r_*}\psi^2\right]-r^{-2}V_{\epsilon_{temper}}\psi^2 \\
\lesssim divJ[X_{\epsilon_{temper}},w_{\epsilon_{temper}}],
\end{multline*}
with a constant independent of $\epsilon_{temper}$. Moreover, the function $V_{\epsilon_{temper}}$ is supported in a neighborhood of the event horizon and satisfies
\begin{equation}\label{s:V_temper_bound_eqn}
||V_{\epsilon_{temper}}||_{L^1(r)}\le \epsilon_{temper}.
\end{equation}
\textbf{This new estimate is in fact worse than the estimate from Proposition \ref{s:partial_morawetz_prop}, but the advantage gained in this lemma is that, unlike the vectorfield $X_0$, the vectorfield $X_{\epsilon_{temper}}$ is regular up to the event horizon.}
\end{lemma}

\begin{remark}
It will be observed during the proof that the function $w_{\epsilon_{temper}}$ is supported outside a neighborhood of the event horizon and $\pd_r(r^2w_{\epsilon_{temper}})$ is nonnegative. These facts will be used when computing boundary terms.
\end{remark}

\begin{proof}
We obtain $X_{\epsilon_{temper}}$ and $w_{\epsilon_{temper}}$ by changing $u$ so that it is bounded near the event horizon. The approach used here was found in \cite{Ma}.

Let $F:\mathbb{R}\rightarrow\mathbb{R}$ be a smooth function satisfying $F(x)=x$ for $x\le1$ and $F(x)=2$ for $x\ge 3$. We use $F(x)$ to temper $u$ by defining
$$\tilde{u}=-\epsilon^{-1}F(-\epsilon u),$$
where $\epsilon>0$ is a small constant to be chosen later. Then we define
$$\tilde{X}=\tilde{u}v\pd_r,$$
$$\tilde{w}=v\pd_r\tilde{u}.$$
Our first claim is that $\tilde{X}$ and $\tilde{w}$ agree with $X_0$ and $w_0$ outside a neighborhood of the event horizon. Given the similarities between the defintions of $\tilde{X}$, $\tilde{w}$ and the identities $X_0=uv\pd_r$, $w_0=v\pd_r u$, it suffices to show that $\tilde{u}$ only differs from $u$ in a neighborhood of the event horizon.

Recall that $u$ is negative and unbounded near the event horizon. Note that when $-\epsilon^{-1}< u$, then $\tilde{u}$ agrees with $u$ as
$$-\epsilon^{-1}<u \Rightarrow (-\epsilon u)< 1,$$
$$\tilde{u}=-\epsilon^{-1}F(-\epsilon u)=(-\epsilon^{-1})(-\epsilon u)=u.$$
So $\tilde{u}$ only differs from $u$ in the region where $u<-\epsilon^{-1}$, which is a small neighborhood of the event horizon.

Since $F(x)\le 2$, it follows that $-2\epsilon^{-1}\le \tilde{u}$, which means that for all $\epsilon>0$, the vectorfield $\tilde{X}$ will be regular up to the event horizon. Furthermore,
\begin{equation}\label{s:tilde_w_in_terms_of_w_0_eqn}
\tilde{w}=v(-\epsilon^{-1}F'(-\epsilon u))(-\epsilon\pd_ru)=v\pd_ruF'(-\epsilon u)=w_0F'(-\epsilon u),
\end{equation}
so the function $\tilde{w}$ will be zero in a neighborhood of the event horizon and $\pd_r(r^2\tilde{w})\ge 0$. This verifies the remark following the lemma.

Recall that
\begin{equation}\label{s:temper_modified_div_eqn}
divJ[\tilde{X},\tilde{w}]=2\pd_r\left(\frac{\tilde{u}}{r^2}\right)\left(1-\frac{2M}r\right)^2(\pd_r\psi)^2+\tilde{u}\pd_rv|\sla\nabla\psi|^2-\frac12\Box_g\tilde{w}\psi^2.
\end{equation}
We claim that the $(\pd_r\psi)^2$ and $|\sla\nabla\psi|^2$ terms in (\ref{s:temper_modified_div_eqn}) satisfy the bounds given in the statement of the lemma. In the region where $\tilde{u}$ differs from $u$,
$$\pd_r\left(\frac{\tilde{u}}{r^2}\right)=-\frac{\tilde{u}}{r^3}+\frac{\pd_r\tilde{u}}{r^2}\ge -\frac{\tilde{u}}{r^2}.$$
It follows that the bound
$$\frac{M^2}{r^3}\lesssim\pd_r\left(\frac{\tilde{u}}{r^2}\right)$$
holds globally with a constant independent of $\epsilon$. The same can be said for the bound
$$\frac1r\left(1-\frac{3M}r\right)^2\lesssim \tilde{u}\pd_rv.$$

Now, we examine the $\psi^2$ term in (\ref{s:temper_modified_div_eqn}). The effect of replacing $w_0$ with $\tilde{w}$ on the $\psi^2$ term will be the indroduction of the quantity
$$-\frac12(\Box_g\tilde{w}-\Box_g w_0)\psi^2.$$
Since $\tilde{w}=w_0$ outside a small neighborhood of the event horizon and since $\Box_gw_0=0$ for $r<r_*=4M$, it follows that
$$-\frac12(\Box_g\tilde{w}-\Box_g w_0)\psi^2=-\frac12\Box_g\tilde{w}1_{r<r_*}\psi^2.$$
Define
$$V_{\epsilon_{temper}}=\frac{r^2}2\Box_g\tilde{w}1_{r<r_*}.$$
To complete the proof of the lemma, it suffices to show that $\epsilon$ can be chosen sufficiently small so that (\ref{s:V_temper_bound_eqn}) holds. From (\ref{s:tilde_w_in_terms_of_w_0_eqn}) and the fact that $\pd_r w_0=0$ in a neighborhood of the event horizon, we calculate
\begin{align*}
\pd_r\tilde{w}&=\pd_r(w_0F'(-\epsilon u)) \\
&=-\epsilon w_0 F''(-\epsilon u)\pd_ru.
\end{align*}
It follows that
\begin{align*}
r^2\Box_g\tilde{w}&=\pd_r\left((r^2-2Mr)\pd_r\tilde{w}\right) \\
&= \pd_r\left((r^2-2Mr)(-\epsilon w_0F''(-\epsilon u)\pd_ru)\right) \\
&=-\epsilon w_0 F''(-\epsilon u) \pd_r\left((r^2-2Mr)\pd_ru\right)+\epsilon^2w_0 F'''(-\epsilon u)(r^2-2Mr)(\pd_ru)^2.
\end{align*}
In the last line, the terms have been grouped based on whether they contiain $F''$ or $F'''$. We estimate each of these terms seperately, observing that $\pd_ru\sim x^{-1}$ near $r=2M$.
\begin{align*}
|-\epsilon w_0 F''(-\epsilon u)\pd_r\left((r^2-2Mr)\pd_ru\right)| &\le \epsilon w_0 |F''(-\epsilon u)| \pd_r(r(r-2M)\pd_ru) \\
&\lesssim \epsilon |F''(-\epsilon u)|\pd_r(r(r-2M)(r-2M)^{-1}) \\
&\lesssim \epsilon |F''(-\epsilon u)|.
\end{align*}
\begin{align*}
|\epsilon^2w_0 F'''(-\epsilon u)(r^2-2Mr)(\pd_ru)^2| &\le \epsilon^2 w_0|F'''(-\epsilon u)|r(r-2M)(\pd_ru)^2 \\
&\lesssim \epsilon^2|rF'''(-\epsilon u)|(r-2M)(r-2M)^{-2} \\
&\lesssim \epsilon^2|rF'''(-\epsilon u)|(r-2M)^{-1}.
\end{align*}
Let $x=r-2M$. From the above bounds, it follows that
$$|r^2\Box_g\tilde{w}|\lesssim \epsilon\chi_\epsilon +\epsilon^2x^{-1}\chi_\epsilon,$$
where $\chi_\epsilon$ is some bounded function supported where $-\epsilon u\in [1,3]$. Since $u\sim \log x$, then $\chi_\epsilon$ is supported in an interval of the form $[e^{-c_2/\epsilon},e^{-c_1/\epsilon}]$. Thus,
\begin{align*}
||V_{temper}||_{L^1(r)} &\lesssim \int \epsilon\chi_\epsilon +\epsilon^2\chi_\epsilon x^{-1}dx \\
 &\lesssim \int_{e^{-c_2/\epsilon}}^{e^{-c_1/\epsilon}}(\epsilon+\epsilon^2x^{-1})dx=\epsilon\left[x+\epsilon \ln|x|\right]_{e^{-c_2/\epsilon}}^{e^{-c_1/\epsilon}}\le\epsilon(e^{-c_1/\epsilon}+c_2-c_1).
\end{align*}
This allows us to conclude that for any $\epsilon_{temper}>0$, it is possible to choose $\epsilon>0$ sufficiently small so that
$$||V_{\epsilon_{temper}}||<\epsilon_{temper}.$$
\end{proof}

\subsection{The redshift vectorfield $\epsilon_{redshift}Y$}\label{s:redshift_correction_sec}

Now, we begin to use the coordinates corresponding to the $\delh$-foliation that is described by the following Lemma.

\begin{lemma}\label{s:foliation_lem}
For all sufficiently small $\delta>0$ (usually denoted by $\delh$ outside of this lemma) there exists a time foliation agreeing with the time foliation of the Boyer-Lindquist coordinates for $r>r_H+\delta$ such that the metric under this foliation satisfies
$$g^{tr}\ge 0,$$
$$r^2g^{tr}=r_H^2+O((r-r_H)^2) \text{ near }r=r_H,$$
and
$$\pd_r(r^2g^{tt})-2\frac{g^{tr}}{g^{rr}}\pd_r(r^2g^{tr})>0$$
for all $r\in [0,\delta]$.
\end{lemma}

\begin{proof}
Denote the Boyer-Lindquist coordinates by $(\bar{t},\bar{r},\bar{\theta},\bar{\phi})$. Then we define a new set of coordinates by
\begin{align*}
t&=\bar{t}+\int_{r_H+\delta}^{\bar{r}}\frac{T(x)dx}{x^2-2Mx} \\
r&=\bar{r} \\
\theta&=\bar{\theta} \\
\phi&=\bar{\phi}
\end{align*}
for some function $T(\bar{r})$ supported for $\bar{r}\le r_H+\delta$. The time foliation is given by the new coordinate $t$, and the metric can be found by the change of coordinates. Under this change of coordinates, since $r=\bar{r}$ everywhere, we have
\begin{align*}
\pd_t&=\pd_{\bar{t}} \\
\pd_r&=\pd_{\bar{r}}-\frac{T}{r^2-2Mr}\pd_{\bar{t}} \\
\pd_\theta&=\pd_{\bar{\theta}} \\
\pd_\phi&=\pd_{\bar{\phi}}
\end{align*}
The inverse metric in these coordinates can be found by computing the following coordinate-invariant quantity.
\begin{align*}
\pd^\lambda\psi\pd_\lambda\psi &= -\left(1-\frac{2M}{r}\right)^{-1}(\pd_{\bar{t}}\psi)^2+\left(1-\frac{2M}r\right)(\pd_{\bar{r}}\psi)^2+\sla{g}^{\bar{\alpha}\bar{\beta}}\pd_{\bar{\alpha}}\psi\pd_{\bar{\beta}}\psi \\
&=-\left(1-\frac{2M}{r}\right)^{-1}(\pd_t\psi)^2+\left(1-\frac{2M}{r}\right)\left(\pd_r\psi+\frac{T}{r^2-2Mr}\pd_t\psi\right)^2+\sla{g}^{\alpha\beta}\pd_\alpha\psi\pd_\beta\psi \\
&=\frac{T^2-r^4}{r^4-2Mr^3}(\pd_t\psi)^2+2\frac{T}{r^2}\pd_t\psi\pd_r\psi+\left(1-\frac{2M}{r}\right)(\pd_r\psi)^2+\sla{g}^{\alpha\beta}\pd_\alpha\psi\pd_\beta\psi
\end{align*}
It follows that in these coordinates,
\begin{align*}
r^2g^{tt}&=\frac{T^2-r^4}{r^2-2Mr} \\
r^2g^{tr}&= T \\
r^2g^{rr}&= r^2-2Mr \\
r^2\sla{g}^{\alpha\beta}&=r^2\sla{g}^{\bar{\alpha}\bar{\beta}}
\end{align*}
If we impose the condition
$$T=r_H^2+O((r-r_H)^2)$$
then all of these metric components are regular and the first two conditions of the lemma are satisfied. We now examine the final condition. Multiplying by $r^2g^{rr}$, we obtain
\begin{align*}
r^2g^{rr}\pd_r(r^2g^{tt})-2r^2g^{tr}\pd_r(r^2g^{tr}) &= \pd_r\left(r^2g^{rr}r^2g^{tt}-(r^2g^{tr})^2\right)-r^2g^{tt}\pd_r(r^2g^{rr}) \\
&=-\pd_r(r^4)-\frac{T^2-r^4}{r^2-2Mr}\pd_r(r^2-2Mr)
\end{align*}
Multiplying again by $r^2g^{rr}=r^2-2Mr$, we obtain
$$(r^4-T^2)\pd_r(r^2-2Mr)-(r^2-2Mr)\pd_r(r^4).$$
Since $T(r_H)=r_H^2$, this quantity vanishes at $r_H$. Since $T'(r_H)=0$, the first derivative of this quantity also vanishes. If $T''(r_H)$ is sufficiently negative, then this quantity can be made positive in a neighborhood of $r_H$. Since we obtained this quantity by multiplying the original quantity by $(r^2-2Mr)^2$, then the original quantity is strictly positive on the support of $T$.
\end{proof}

 The change this has on the partial Morawetz estimate is minimal--it only means that the term ``$(\pd_r\psi)^2$'' really becomes $(\pd_r\psi+f(r)\pd_t\psi)^2$ for some function $f$ supported where $r<r_H+\delh$. Since this term still has the same sign, it may simply be ignored in that region. (It wasn't very useful near the event horizon anyway--it degenerated like $O((r-r_H)^2)$.)

We now show that by considering the effects of the $\delh$-foliation, we can define a redshift vectorfield $Y$ that will give uniform control of $(\pd_r\psi)^2$ and $(\pd_t\psi)^2$ in a neighborhood of the event horizon larger than $\delh$. This is summarized in the following lemma.

\begin{lemma}\label{s:redshift_lem}
On a $\delh$-foliation guaranteed by Lemma \ref{s:foliation_lem}, there exists a vectorfield $Y$ supported on the region $r\in[r_H,r_H+2\delh]$ satisfying
$$1_{r-r_H\in [0,3\delh/2]}\left[(\pd_r\psi)^2+(\pd_t\psi)^2\right]\lesssim divJ[Y]+1_{r-r_H\in (\delh,2\delh)}|\sla\nabla\psi|^2.$$
Moreover, $Y=Y^r(r)\pd_r+Y^t(r)\pd_t$, and $Y^r(r_H)< 0$.
\end{lemma}
\begin{proof}
Let $Y=Y^r(r)\pd_r+Y^t(r)\pd_t$. Then by Lemma \ref{divJ_lem},
$$divJ[Y]=K_Y^{\mu\nu}\pd_\mu\psi\pd_\nu\psi,$$
where
\begin{align*}
K_Y^{\mu\nu}&=g^{\mu\lambda}\pd_\lambda Y^\nu+g^{\nu\lambda}\pd_r Y^\mu-divYg^{\mu\nu}-Y^\lambda\pd_\lambda g^{\mu\nu} \\
&=g^{\mu r}\pd_r Y^\nu+g^{\nu r}\pd_r Y^\mu-\frac{1}{r^2}\pd_r(r^2 Y^r)g^{\mu\nu}-Y^r\pd_r g^{\mu\nu} \\
&=g^{\mu r}\pd_r Y^\nu+g^{\nu r}\pd_r Y^\mu-\frac{1}{r^2}\pd_r(Y^r r^2g^{\mu\nu}).
\end{align*}
It follows that
\begin{align*}
r^2K_Y^{rr}&=(r^2-2Mr) \pd_rY^r-Y^r\pd_r(r^2-2Mr) \\
2r^2K_Y^{(tr)}&=2(r^2g^{rr}\pd_r Y^t-\pd_r(r^2g^{tr})Y^r) \\
r^2K_Y^{tt}&=2r^2g^{tr}\pd_rY^t-\pd_r(Y^rr^2g^{tt}) \\
r^2\sla{K}_Y^{\alpha\beta}&=-\pd_r Y^r \sla{g}^{\alpha\beta}
\end{align*}
We will take $Y^r$ to be a negative constant\footnote{It is traditional to take $\pd_rY^r<0$ near the event horizon to also gain control of the angular derivatives, but this is not necessary here since the current $J[X_{\epsilon_{temper}},w_{\epsilon_{temper}}]$ already has good control of these derivatives and such a choice complicates the proof.} on $r-r_H\in [0,\delh]$ and smoothly increasing to $0$ on $r-r_H\in (\delh,2\delh)$.

With this choice, $K_Y^{rr}$ has the right sign for $r-r_H\in[0,2\delh]$. In particular, this is due to the fact that the term $-Y^r\pd_r(r^2-2Mr)$ is strictly positive at $r=2M$. This fact is related to the surface gravity of the Schwarzschild black hole.

However, $\sla{K}_Y^{\alpha\beta}$ has the wrong sign when $\pd_rY^r\ne 0$. This explains the error term on the right side of the estimate in the lemma.

We choose $Y^t$ to have support coinciding with $g^{tr}$ and such that $K_Y^{(tr)}=0$. That is, $Y^t$ is determined by the relation
$$\pd_rY^t=\frac{1}{g^{rr}}\pd_r(r^2g^{tr})Y^r$$
and the condition $Y^t(\delh)=0$.
Note that although $g^{rr}$ vanishes at the event horizon, the quantity $\frac{1}{g^{rr}}\pd_r(r^2g^{tr})$ remains bounded, because $\pd_r(r^2g^{tr})$ also vanishes to first order as $r\rightarrow r_H$. (This was a condition specified in Lemma \ref{s:foliation_lem}.)

Having made the above choice for $Y^t$, we have in the region $r-r_H\in[0,\delh]$ that
\begin{align*}
r^2K_Y^{tt}&=2\frac{g^{tr}}{g^{rr}}\pd_r(r^2g^{tr})Y^r-\pd_r(Y^rr^2g^{tt}) \\
&= (-Y^r)\left[\pd_r(r^2g^{tt})-2\frac{g^{tr}}{g^{rr}}\pd_r(r^2g^{tr})\right].
\end{align*}
By Lemma \ref{s:foliation_lem}, this is bounded below by a positive constant on $r-r_H\in [0,\delh]$.
 
In the region $r-r_H\in (\delh,2\delh)$,
$$r^2K_Y^{tt}=-\pd_r(Y^rr^2g^{tt})$$
With an appropriate choice of $Y^r$, this term will be positive on $r-r_H\in [\delh,2\delh)$, as $Y^rr^2g^{tt}$ is positive and must go to zero as $Y^r$ vanishes.
\end{proof}

The previous lemma leads to the following corollary.
\begin{corollary}\label{s:redshift_cor}
There exists an $\epsilon_{redshift}$ sufficiently small so that for all $\epsilon_{temper}$ sufficiently small,
\begin{multline*}
\left[\frac{M^2}{r^3}(\pd_r\psi)^2+1_{r-r_H\in[0,3\delh/2]}(\pd_t\psi)^2+\frac1r\left(1-\frac{r_{trap}}r\right)^2|\sla\nabla\psi|^2+\frac{M}{r^4}1_{r>r_*}\psi^2\right]-r^{-2}V_{\epsilon_{temper}}\psi^2 \\
\lesssim divJ[X_{\epsilon_{temper}}+\epsilon_{redshift}Y,w_{\epsilon_{temper}}].
\end{multline*}
\end{corollary}
\begin{proof}
By Lemma \ref{s:temper_lemma}, we have that for some function $f$ supported on the interval $r\in [r_H,r_H+\delh]$,
\begin{multline*}
\left[\frac{M^2}{r^3}\left(1-\frac{r_H}r\right)^2(\pd_r\psi+f\pd_t\psi)^2+\frac1r\left(1-\frac{r_{trap}}r\right)^2|\sla\nabla\psi|^2+\frac{M}{r^4}1_{r>r_*}\psi^2\right]-r^{-2}V_{\epsilon_{temper}}\psi^2 \\
\lesssim divJ[X_{\epsilon_{temper}},w_{\epsilon_{temper}}],
\end{multline*}
Thus,
\begin{multline*}
\left[\frac{M^2}{r^3}1_{r>r_H+\delh}(\pd_r\psi)^2+\frac1r\left(1-\frac{r_{trap}}r\right)^2|\sla\nabla\psi|^2+\frac{M}{r^4}1_{r>r_*}\psi^2\right]-r^{-2}V_{\epsilon_{temper}}\psi^2 \\
\lesssim divJ[X_{\epsilon_{temper}},w_{\epsilon_{temper}}],
\end{multline*}
From Lemma \ref{s:redshift_lem},
$$1_{r-r_H\in [0,3\delh/2]}\left[(\pd_r\psi)^2+(\pd_t\psi)^2\right]\lesssim divJ[Y]+1_{r-r_H\in (\delh,2\delh)}|\sla\nabla\psi|^2.$$
By taking $\epsilon_{redshift}$ sufficiently small, these two estimates can be combined to yield the estimate in the corollary for $divJ[X_{\epsilon_{temper}}+\epsilon_{redshift}Y,w_{\epsilon_{temper}}]$.
\end{proof}

\subsection{The local Hardy estimate}\label{s:hardy_correction_sec}

We now address the $\psi^2$ term in the region $r\in[r_H,r_*]$. Since we finally have good control of $(\pd_r\psi)^2$, we may essentally apply the fundamental theorem of Calculus to ``borrow'' from the $\psi^2$ term in the region $r\in[r_*,\infty)$ to control $\psi^2$ in this region at the additional cost of some of the $(\pd_r\psi)^2$ term. This technique, which we refer to as the \textit{local Hardy estimtate}, is made precise by the lemma below. Its application is explained in the corollary that follows.

\begin{lemma}\label{s:hardy_lem}
Let $V_1$ be a non-negative smooth function of $r$ supported on the interval $[r_H,r_*]$ and let $V_2$ be any non-negative smooth function supported on an interval $[r_*,R_*]$ (for any $R_*>r_*$) and satisfying $||V_2||_{L^1(r)}=2||V_1||_{L^1(r)}$. Then there exists a function $A_{12}$ depending on $V_1$ and $V_2$ such that
$$\int V_1\psi^2dr \le \int A_{12}(\pd_r\psi)^2+ V_2\psi^2dr.$$
Furthermore, $A_{12}$ is supported on the interval $[r_H,R_*]$ and satisfies
$$||A_{12}||_{L^\infty}\le 2(R_*-r_H)||V_1||_{L^1(r)}.$$
\end{lemma}

\begin{proof}
We start with the fundamental theorem of calculus.
$$\psi(r)=-\int_r^R\psi'(\tilde{r})d\tilde{r}+\psi(R).$$
Then, assuming $R>r$,
\begin{align*}
\psi^2(r) &= \left(-\int_r^R\psi'(\tilde{r})d\tilde{r}+\psi(R)\right)^2 \\
&\le 2\left(\int_r^R|\psi'(\tilde{r})|d\tilde{r}\right)^2+2(\psi(R))^2 \\
&\le 2(R-r)\int_r^R(\pd_r\psi)^2+2\psi^2(R).
\end{align*}
Letting $R=R(r)$, multiplying by a non-negative function $V_1(r)$, and integrating in $r$, 
$$\int V_1(r)\psi^2(r)dr\le 2\int V_1(r)(R(r)-r)\int_r^{R(r)}(\psi'(\tilde{r}))^2d\tilde{r}dr+2\int V_1(r)\psi^2(R(r))dr.$$
We may choose another non-negative function $V_2(r)$ with the property that $2\int V_1=\int V_2$ and we may choose $R(r)$ so that 
$$2\int_{-\infty}^r V_1(\tilde{r})d\tilde{r}=\int_{-\infty}^{R(r)}V_2(\tilde{R})d\tilde{R}=\int_{-\infty}^r V_2(R(\tilde{r}))R'(\tilde{r})d\tilde{r}.$$
 Then we have the pointwise identity $2V_1(r)=V_2(R(r))R'(r)$. It follows that
$$\int (V_1(r)-V_2(r))\psi^2(r)dr\le 2\int V_1(r)(R(r)-r)\int_r^{R(r)}(\psi'(\tilde{r}))^2d\tilde{r}dr.$$
Exchanging the order of integration on the right side gives the local Hardy estimate
$$\int (V_1-V_2)\psi^2\le \int A_{12}(\pd_r\psi)^2,$$
where, in the particular case\footnote{A more general case where $r_*\ne R_H$ is very similar.} $r\in [r_H,r_*]$ and $R\in [R_H,R_*]$ with $r_*=R_H$, the function $A_{12}$ is given by
\begin{align*}
A_{12}(r)|_{r_H\le r\le r_*}&=\int_{r_H}^r2V_1(\tilde{r})(R(\tilde{r})-\tilde{r})d\tilde{r} \\
A_{12}(R)|_{R_H\le R\le R_*}&=\int_{r(R)}^{r_*}2V_1(\tilde{r})(R(\tilde{r})-\tilde{r})d\tilde{r},
\end{align*}
and satisfies the following inequality
$$||A_{12}||_{L^\infty}\le \max_{r\in [r_H,r_*]}2(R(r)-r)||V_1||_{L^1}.$$
This completes the proof of Lemma \ref{s:hardy_lem}.
\end{proof}

The previous lemma leads to the following corollary.
\begin{corollary}\label{s:hardy_cor}
There exists $\epsilon_{temper}$ sufficiently small so that
\begin{multline*}
\int_{\Sigma_t}\left[\frac{M^2}{r^3}(\pd_r\psi)^2+1_{r-r_H\in[0,3\delh/2]}(\pd_t\psi)^2+\frac1r\left(1-\frac{r_{trap}}r\right)^2|\sla\nabla\psi|^2+\frac{M}{r^4}\psi^2\right] \\
\lesssim \int_{\Sigma_t}divJ[X_{\epsilon_{temper}}+\epsilon_{redshift}Y,w_{\epsilon_{temper}}].
\end{multline*}
\end{corollary}

\begin{proof}
Fix an $R_*>r_*$. Then choose $\epsilon_{Hardy}$ sufficiently small that
\begin{equation}\label{s:epsilon_hardy_condition_1_eqn}
2\epsilon_{Hardy}\le \frac12\int_{r_*}^{R_*}\frac{M}{r^2}dr 
\end{equation}
and
\begin{equation}\label{s:epsilon_hardy_condition_2_eqn}
2(R_*-r_H)\epsilon_{Hardy}\le \frac12\frac{M^2}{r}\text{ on the interval }r\in[r_H,R_*].
\end{equation}
Recall that $V_{\epsilon_{temper}}$ is supported near the event horizon and satisfies $||V_{\epsilon_{temper}}||_{L^1(r)}\le \epsilon_{temper}$. Choose $\epsilon_{temper}=\frac12\epsilon_{Hardy}$. Then it is possible to choose $V_1$ supported on $r\in[r_H,r_*]$ satisfying $||V_1||_{L^1(r)}=\epsilon_{Hardy}$ and
$$1_{r\in [r_H,r_*]}\lesssim V_1-V_{\epsilon_{temper}}.$$
Given condition (\ref{s:epsilon_hardy_condition_1_eqn}) it is also possible to choose $V_2$ supported on $r\in[r_*,R_*]$ satisfying $||V_2||_{L^1(r)}=2||V_1||_{L^1(r)}=2\epsilon_{Hardy}$ and
$$\frac12\frac{M}{r^2}1_{r>r_*}\le \frac{M}{r^2}1_{r-r*}-V_2.$$
From these two inequalities, we conclude
$$\frac{M}{r^2} \lesssim \frac{M}{r^2}1_{r-r*}-V_{\epsilon_{temper}}+V_1-V_2.$$
By Lemma \ref{s:hardy_lem} and (\ref{s:epsilon_hardy_condition_2_eqn}), since $||V_1||_{L^1(r)}=\epsilon_{Hardy}$,
$$A_{12}[V_1,V_2]\le 2(R_*-r_H)||V_1||_{L^1}= 2(R_*-r_H)\epsilon_{Hardy} \le \frac12\frac{M^2}r$$
on the interval $[r_H,R_*]$. Since $A_{12}[V_1,V_2]$ is supported on the interval $[r_H,R_*]$, it follows that
$$\frac12\frac{M^2}r\le \frac{M^2}r-A_{12}[V_1,V_2].$$

So far we have established
\begin{multline*}
\int_{r_H}^\infty \left[\frac{M^2}{r^3}(\pd_r\psi)^2+\frac{M}{r^4}\psi^2\right]r^2dr \\
\lesssim \int_{r_H}^\infty \left(\frac{M^2}{r}-A_{12}\right)(\pd_r\psi)^2+\left(\frac{M}{r^2}1_{r>r_*}-V_{\epsilon_{temper}}+V_1-V_2\right)\psi^2 dr.
\end{multline*}
From this and the local Hardy estimate (Lemma \ref{s:hardy_lem}) it follows that
\begin{align*}
\int_{r_H}^\infty \left[\frac{M^2}{r^3}(\pd_r\psi)^2+\frac{M}{r^4}\psi^2\right]r^2dr
&\lesssim \int_{r_H}^\infty \frac{M^2}{r}(\pd_r\psi)^2+\left(\frac{M}{r^2}1_{r>r_*}-V_{\epsilon_{temper}}\right)\psi^2 dr \\
&= \int_{r_H}^\infty \left[\frac{M^2}{r^3}(\pd_r\psi)^2+\left(\frac{M}{r^4}1_{r>r_*}-r^{-2}V_{\epsilon_{temper}}\right)\psi^2 \right]r^2dr \\
&\lesssim \int_{r_H}^\infty divJ[X_{\epsilon_{temper}}+\epsilon_{redshift}Y,w_{\epsilon_{temper}}]r^2dr.
\end{align*}
\end{proof}

\subsection{The correction $\epsilon_{\pd_t}w_{\pd_t}$}\label{s:dt_correction_sec}

By including an additional function $w_{\pd_t}$, we complete the term with $(\pd_t\psi)^2$ in our estimate.
\begin{lemma}\label{s:dt_lemma}
There exists a function $w_{\pd_t}$ supported on $r\in [r_H+\delh,\infty)$ such that
\begin{multline*}
\frac{M^2}{r^3}\left(1-\frac{r_{trap}}r\right)^21_{r\ge r_H+3\delh/2}(\pd_t\psi)^2 \\
\lesssim divJ[0,w_{\pd_t}]+\frac{M^2}{r^3}\left(1-\frac{r_{trap}}r\right)^21_{r>r_H+\delh}\left[(\pd_r\psi)^2+|\sla\nabla\psi|^2+r^{-2}\psi^2\right].
\end{multline*}
\end{lemma}
\begin{proof}
The proof of this lemma is rather similar to the proof of Lemmas \ref{m:wdt_correction_I_lem} and \ref{m:wdt_correction_II_lem}, except that additional care must be taken near the horizon and trapping radius. Otherwise, the idea is essentially the same.

Let $w_{\pd_t}$ be a smooth, negative function supported on $r>r_H+\delh$, vanishing to second order at $r=r_{trap}$, and decaying like $O(-M^2/r^3)$ for large $r$. Then
$$divJ[0,w_{\pd_t}]=w_{\pd_t}g^{\alpha\beta}\pd_\alpha\psi\pd_\beta\psi-\frac12\Box_gw_{\pd_t}\psi^2.$$
Since $w_{\pd_t}$ is supported in the region $r>r_H+\delh$,
$$divJ[0,w_{\pd_t}]=w_{\pd_t}g^{tt}(\pd_t\psi)^2+w_{\pd_t}g^{rr}(\pd_r\psi)^2+w_{\pd_t}|\sla\nabla\psi|^2-\frac12\Box_gw_{\pd_t}\psi^2.$$
By the choice of $w_{\pd_t}$, especially since both $w_{\pd_t}$ and $g^{tt}$ are negative, one can check that
$$\frac{M^2}{r^3}\left(1-\frac{r_{trap}}r\right)^21_{r\ge r_H+3\delh/2}(\pd_t\psi)^2\lesssim w_{\pd_t}g^{tt},$$
while
$$|w_{\pd_t}|\lesssim \frac{M^2}{r^3}\left(1-\frac{r_{trap}}r\right)^21_{r>r_H+\delh}.$$
A straightforward calculation shows that $\Box_g(M^2/r^3)=O(M^2/r^5)$, so
$$\left|\frac12\Box_gw_{\pd_t}\right|\lesssim \frac{M^2}{r^5}\left(1-\frac{r_{trap}}r\right)^21_{r>r_H+\delh}.$$
 The result follows.
\end{proof}

The previous lemma leads to the following corollary.
\begin{corollary}\label{s:dt_cor}
There exists a small $\epsilon_{\pd_t}>0$ so that
\begin{multline*}
\int_{\Sigma_t}\left[\frac{M^2}{r^3}(\pd_r\psi)^2+\left(1-\frac{r_{trap}}r\right)^2\left[\frac1r|\sla\nabla\psi|^2+\frac{M^2}{r^3}(\pd_t\psi)^2\right]+\frac{M}{r^4}\psi^2\right] \\
\lesssim \int_{\Sigma_t}divJ[X_{\epsilon_{temper}}+\epsilon_{redshift}Y,w_{\epsilon_{temper}}+\epsilon_{\pd_t}w_{\pd_t}].
\end{multline*}
\end{corollary}
\begin{proof}
Recall from Corollary \ref{s:hardy_cor} that
\begin{multline*}
\int_{\Sigma_t}\left[\frac{M^2}{r^3}(\pd_r\psi)^2+1_{r-r_H\in[0,3\delh/2]}(\pd_t\psi)^2+\frac1r\left(1-\frac{r_{trap}}r\right)^2|\sla\nabla\psi|^2+\frac{M}{r^4}\psi^2\right] \\
\lesssim \int_{\Sigma_t}divJ[X_{\epsilon_{temper}}+\epsilon_{redshift}Y,w_{\epsilon_{temper}}].
\end{multline*}
Now, by Lemma \ref{s:dt_lemma},
\begin{multline*}
\frac{M^2}{r^3}\left(1-\frac{r_{trap}}r\right)^21_{r\ge r_H+3\delh/2}(\pd_t\psi)^2 \\
\lesssim divJ[0,w_{\pd_t}]+\frac{M^2}{r^3}\left(1-\frac{r_{trap}}r\right)^21_{r>r_H+\delh}\left[(\pd_r\psi)^2+|\sla\nabla\psi|^2+r^{-2}\psi^2\right].
\end{multline*}
By taking $\epsilon_{\pd_t}$ sufficiently small, these two estimates can be combined to yield the estimate in the corollary for $divJ[X_{\epsilon_{temper}}+\epsilon_{redshift}Y,w_{\epsilon_{temper}}+\epsilon_{\pd_t}w_{\pd_t}]$.
\end{proof}

\subsection{The correction $\pd_t$ and the flux terms}\label{s:boundary_correction_sec}

To investigate the flux terms, we begin with a few lemmas.

\begin{lemma}\label{s:Jmu_lem} (Identities for $J^\mu$)
\begin{align*}
-J^t[\pd_t] &= -g^{tt}(\pd_t\psi)^2+g^{rr}(\pd_r\psi)^2+|\sla\nabla\psi|^2 \\
-J^t[\pd_r]&= -2g^{tt}\pd_t\psi\pd_r\psi-2g^{tr}(\pd_r\psi)^2 \\
J^r[\pd_t]&=2g^{rr}\pd_r\psi\pd_t\psi+2g^{rt}(\pd_t\psi)^2 \\
J^r[-\pd_r]&=-g^{rr}(\pd_r\psi)^2+g^{tt}(\pd_t\psi)^2+|\sla\nabla\psi|^2.
\end{align*}
\end{lemma}
\begin{proof}
These identities follow from the fact that
$$J^\alpha[\pd_\beta]=g^{\alpha\lambda}T_{\lambda\beta}=2g^{\alpha\lambda}\pd_\lambda\psi\pd_\beta\psi-\delta^\alpha{}_\beta(\pd^\lambda\psi\pd_\lambda\psi).$$
\end{proof}

\begin{lemma}\label{s:Jr_horizon_lem}
On the event horizon,
$$|\sla\nabla\psi|^2 = J^r[-\pd_r]-g^{tt}(\pd_t\psi)^2,$$
$$(\pd_t\psi)^2\approx J^r[\pd_t].$$
\end{lemma}
\begin{proof}
From Lemma \ref{s:Jmu_lem}, and the fact that $g^{rr}=0$ on the event horizon,
\begin{align*}
J^r[\pd_t]&=2g^{rt}(\pd_t\psi)^2 \\
J^r[-\pd_r]&=g^{tt}(\pd_t\psi)^2+|\sla\nabla\psi|^2.
\end{align*}
Note that $g^{rt}$ and $-g^{tt}$ are both positive and bounded in a neighborhood including the event horizon. The estimates of the lemma follow.
\end{proof}

\begin{lemma}\label{s:JL_lem}
On a constant-time hypersurface, outside the the region $r\in[r_H,r_H+\delh]$,
$$-J^t[L] = \alpha^{-1}(L\psi)^2+|\sla\nabla\psi|^2.$$
\end{lemma}
\begin{proof}
Outside the region $r\in[r_H,r_H+\delh]$, $g^{tr}=0$, $-g^{tt}=\alpha^{-1}$, and $g^{rr}=\alpha$. From Lemma \ref{s:Jmu_lem},
\begin{align*}
-J^t[L] &= \alpha (-J^t[\pd_r])+(-J^t[\pd_t]) \\
&= \alpha\left(-2g^{tt}\pd_t\psi\pd_r\psi-2g^{tr}(\pd_r\psi)^2\right)+\left(-g^{tt}(\pd_t\psi)^2+g^{rr}(\pd_r\psi)^2+|\sla\nabla\psi|^2\right) \\
&= 2\pd_t\psi\pd_r\psi +\alpha^{-1}(\pd_t\psi)^2+\alpha(\pd_r\psi)^2+|\sla\nabla\psi|^2 \\
&= \alpha^{-1}(\alpha\pd_r\psi+\pd_t\psi)^2+|\sla\nabla\psi|^2 \\
&= \alpha^{-1}(L\psi)^2+|\sla\nabla\psi|^2.
\end{align*}
\end{proof}

\begin{lemma}\label{s:J0w_lem}
On a constant-time hypersurface, outside the region $r\in[r_H,r_H+\delh]$, if $w$ is a function of $r$ only, then
$$-J^t[0,w]=\alpha^{-1}w\psi L\psi -r^{-2}\pd_r\left(\frac{r^2}{2} w\psi^2\right) +r^{-2}\pd_r\left(\frac{r^2}{2}w\right)\psi^2.$$
\end{lemma}
\begin{proof}
Recall that
$$J^\mu[0,w]=g^{\mu\nu}w\psi\pd_\nu\psi-\frac12g^{\mu\nu}\pd_\nu w \psi^2.$$
It follows that
$$-J^t[0,w]=-g^{tt}w\psi\pd_t\psi-g^{tr}w\psi\pd_r\psi-\frac12g^{tr}\pd_rw\psi^2.$$
Outside the region $r\in[r_H,r_H+\delh]$,
\begin{align*}
-J^t[0,w] &= -g^{tt}w\psi\pd_t\psi \\
&= \alpha^{-1}w\psi\pd_t\psi \\
&= \alpha^{-1}w \psi L\psi -w\psi\pd_r\psi  \\
&= \alpha^{-1}w\psi L\psi -r^{-2}\pd_r\left(\frac{r^2}{2} w\psi^2\right) +r^{-2}\pd_r\left(\frac{r^2}{2}w\right)\psi^2.
\end{align*}
\end{proof}

\begin{lemma}\label{s:morawetz_flux_lem}
The current 
$$J[X_{\epsilon_{temper}}+\epsilon_{redshift}Y+\pd_t,w_{\epsilon_{temper}}+\epsilon_{\pd_t}w_{\pd_t}]$$
satisfies the following flux estimates.

On the hypersurface $\Sigma_t$,
\begin{multline*}
\int_{\Sigma_t} (L\psi)^2+|\sla\nabla\psi|^2+r^{-2}\psi^2 \\
\lesssim  \int_{\Sigma_t}-J^t[X_{\epsilon_{temper}}+\epsilon_{redshift}Y+\pd_t,w_{\epsilon_{temper}}+\epsilon_{\pd_t}w_{\pd_t}] + Err_{\Sigma_t},
\end{multline*}
where 
$$Err_{\Sigma_t}=\int_{\Sigma_t} r^{-1}|\psi L\psi|+\frac{M^2}{r^2}\left[\chi_H(\pd_r\psi)^2+(\pd_t\psi)^2\right],$$
and $\chi_H=1-\frac{2M}r$.

On the hypersurface $H_{t_1}^{t_2}$,
$$\int_{H_{t_1}^{t_2}}|\sla\nabla\psi|^2\lesssim  \int_{H_{t_1}^{t_2}}-J^t[X_{\epsilon_{temper}}+\epsilon_{redshift}Y+\pd_t,w_{\epsilon_{temper}}+\epsilon_{\pd_t}w_{\pd_t}] + Err_{H_{t_1}^{t_2}},$$
where
$$Err_{H_{t_1}^{t_2}}=\int_{H_{t_1}^{t_2}}(\pd_t\psi)^2.$$
\end{lemma}
\begin{proof}
Given the asymptotics at the end of the proof of Proposition \ref{s:partial_morawetz_prop}, we have that for large $r$
$$X_{\epsilon_{temper}}+\epsilon_{redshift}Y=X_0=\alpha\pd_r+\frac{c^2}{r^2}\pd_r$$
and
$$w_{\epsilon_{temper}}+\epsilon_{\pd_t}w_{\pd_t}=O(r^{-1}).$$
Therefore,
$$X_{\epsilon_{temper}}+\epsilon_{redshift}Y+\pd_t=\alpha\pd_r+\pd_t+O\left(\frac{M^2}{r^2}\right)\pd_r = L+O\left(\frac{M^2}{r^2}\right)\pd_r,$$
and
$$r^{-2}\pd_r\left(\frac{r^2}{2}(w_{\epsilon_{temper}}+\epsilon_{\pd_t}w_{\pd_t})\right)=O(r^{-2}).$$
The estimate on $\Sigma_t$ now follows from Lemmas \ref{s:JL_lem} and \ref{s:J0w_lem}, the fact that all remaining error terms can be estimated by\footnote{Note in particular that there is one term $-J^t[-\pd_r]$ that provides control of $(\pd_r\psi)^2$ near the event horizon (see Lemma \ref{s:Jmu_lem}). This explains the $\chi_H$ factor. This is important, because the $h\pd_t$ estimate can control $\frac{M^2}{r^2}\chi_H(\pd_r\psi)^2$, but not $\frac{M^2}{r^2}(\pd_r\psi)^2$.} $\frac{M^2}{r^2}\left[\chi_H(\pd_r\psi)^2+(\pd_t\psi)^2\right]$, and the fact that if $\psi$ decays sufficiently fast as $r\rightarrow\infty$,
$$\int_{\Sigma_t}r^{-2}\pd_r\left(\frac{r^2}2w\psi^2\right)=0,$$
because it is an integral of a total divergence term and $w$ vanishes near the event horizon (see the remark following Lemma \ref{s:temper_lemma}).

The estimate on the hypersurface $H_{t_1}^{t_2}$ follows from Lemma \ref{s:Jr_horizon_lem} and the fact that $X_{\epsilon_{temper}}$ and $Y$ both have negative $\pd_r$ components on the event horizon.
\end{proof}

\subsection{The Morawetz estimate}

Finally, we conclude with the Morawetz estimate.
\begin{theorem}\label{s:morawetz_thm} (Morawetz estimate in Schwarzschild)
The following estimate holds for any sufficiently regular function $\psi$ decaying sufficiently fast as $r\rightarrow\infty$.
\begin{multline*}
\int_{H_{t_1}^{t_2}}|\sla\nabla\psi|^2+\int_{\Sigma_{t_2}}(L\psi)^2+|\sla\nabla\psi|^2+r^{-2}\psi^2 \\
+\int_{t_1}^{t_2}\int_{\Sigma_t}\frac{M^2}{r^3}(\pd_r\psi)^2+\frac1r\left(1-\frac{3M}r\right)^2\left[|\sla\nabla\psi|^2+\frac{M^2}{r^2}(\pd_t\psi)^2\right]+\frac{M}{r^4}\psi^2 \\
\lesssim \int_{\Sigma_{t_1}}(L\psi)^2+|\sla\nabla\psi|^2+r^{-2}\psi^2 + Err,
\end{multline*}
where
\begin{align*}
Err &= Err_1+Err_2+Err_\Box \\
Err_1 &= \int_{\Sigma_{t_2}} r^{-1}|\psi L\psi| \\
Err_2 &= \int_{H_{t_1}^{t_2}}(\pd_t\psi)^2+\int_{\Sigma_{t_2}}\frac{M^2}{r^2}\left[\chi_H(\pd_r\psi)^2+(\pd_t\psi)^2\right] \\
Err_\Box&=\int_{t_1}^{t_2}\int_{\Sigma_t}|(2X(\psi)+w\psi)\Box_g\psi|,
\end{align*}
and $\chi_H=1-\frac{2M}r$.
\end{theorem}

\begin{proof}
Apply Proposition \ref{general_divergence_estimate_prop} to the current
$$J[X_{\epsilon_{temper}}+\epsilon_{redshift}Y+\pd_t,w_{\epsilon_{temper}}+\epsilon_{\pd_t}w_{\pd_t}]$$
and invoke Corollary \ref{s:dt_cor} and Lemma \ref{s:morawetz_flux_lem}.
\end{proof}

\section{The $r^p$ estimate}\label{s:rp_sec}

\subsection{The pre-$r^p$ identity}

We prove the following lemma, which has a much simpler, well-known analogue in Minkowski spacetime (see Lemma \ref{m:p_ee_identity_lem}). It will be used for the proof of Proposition \ref{s:incomplete_p_estimate_prop}.
\begin{lemma}\label{s:p_ee_identity_lem} Let $\alpha=1-\frac{2M}r$ and $L=\alpha\pd_r+\pd_t$. For any function $f=f(r)$ supported where $r>r_H+\delh$, the following identity holds.
\begin{align*}
&\int_{\Sigma_{t_2}}\left[f\left(\alpha^{-1}L\psi+r^{-1}\psi\right)^2 +\alpha^{-1}f|\sla\nabla\psi|^2+\epsilon r^{-1}f'\psi^2 -r^{-2}\pd_r(rf\psi^2)\right] \\
&+\int_{t_1}^{t_2}\int_{\Sigma_t}
\left[
  (2r^{-1}f-f')|\sla\nabla\psi|^2
+ \alpha f'\left(\alpha^{-1}L\psi+\frac{1-\epsilon}{r}\psi\right)^2 
\right. \\
&\hspace{.9in}
\left.
\vphantom{
  (2r^{-1}f-f')|\sla\nabla\psi|^2
+ \alpha f'\left(\alpha^{-1}L\psi+\frac{1-\epsilon}{r}\psi\right)^2 
}
+ \epsilon\alpha\left((1-\epsilon)f'-rf''\right)r^{-2}\psi^2 
- \alpha'\left[-f(\alpha^{-1}fL\psi)^2+(f-\epsilon rf')r^{-2}\psi^2\right] 
\right]\\
=&\int_{\Sigma_{t_1}}\left[f\left(\alpha^{-1}L\psi+r^{-1}\psi\right)^2 +\alpha^{-1}f|\sla\nabla\psi|^2+\epsilon r^{-1}f'\psi^2 -r^{-2}\pd_r(rf\psi^2)\right] \\
&+\int_{t_1}^{t_2}\int_{\Sigma_t}-\left(2\alpha^{-1}fL\psi+2r^{-1}f\psi\right)\Box_g\psi.
\end{align*}
\end{lemma}

\begin{proof}
We will use Proposition \ref{general_divergence_estimate_prop} together with the following current template.
$$J[X,w,m]_\mu = T_{\mu\nu} X^\nu +w\psi\pd_\mu\psi-\frac12\psi^2\pd_\mu w+m_\mu \psi^2,$$
$$T_{\mu\nu}=2\pd_\mu\psi\pd_\nu\psi-g_{\mu\nu}\pd^\lambda\psi\pd_\lambda\psi.$$

Assume for now that $\Box_g\psi=0$. Let $\alpha = 1-\frac{2M}r$, and observe that
$$L=\alpha\pd_r+\pd_t,$$
$$g^{rr}=\alpha,$$
$$g^{tt}=-\alpha^{-1}.$$

\begin{lemma}\label{s:divJphiX_lem}
Without appealing directly to the particular expression for $\alpha$, one can deduce the following identity.
$$divJ[\alpha^{-1}fL]
=(\alpha^{-1}f)'(L\psi)^2-2r^{-1}f\left(\alpha(\pd_r\psi)^2-\alpha^{-1}(\pd_t\psi)^2\right)-f'|\sla\nabla\psi|^2.$$
\end{lemma}
\begin{proof}
Note that
$$div J[X] = K^{\mu\nu}\pd_\mu\psi\pd_\nu\psi,$$
where
$$K^{\mu\nu}=2g^{\mu\lambda}\pd_\lambda X^\nu-X^\lambda \pd_\lambda(g^{\mu\nu})-div X g^{\mu\nu}.$$

Set $X=\alpha^{-1}f(\alpha\pd_r+\pd_t)=f\pd_r+\alpha^{-1}f\pd_t$. From the above formula, since $g^{rt}=0$ where $f$ is supported,
$$K^{tr}+K^{rt}=2g^{rr}\pd_r X^t=2\alpha\pd_r(\alpha^{-1} f).$$
Thus, the expression for $divJ[\alpha^{-1}fL]$ will have a mixed term of the form 
$$2\alpha\pd_r(\alpha^{-1}f)\pd_r\psi\pd_t\psi.$$
Note that
\begin{align*}
(\alpha^{-1} f)'(L\psi)^2 &= (\alpha^{-1} f)'(\alpha\pd_r\psi+\pd_t\psi)^2 \\
&= \alpha^2(\alpha^{-1}f)'(\pd_r\psi)^2+2\alpha(\alpha^{-1}f)'\pd_r\psi\pd_t\psi +(\alpha^{-1} f)'(\pd_t\psi)^2.
\end{align*}
We now compute the $(\pd_r\psi)^2$ and $(\pd_t\psi)^2$ components, subtracting the part that will be grouped with the $(L\psi)^2$ term.
\begin{align*}
K^{rr}-\alpha^2(\alpha^{-1}f)' &= \left[2g^{rr}\pd_rX^r-X^r\pd_r g^{rr}-r^{-2}\pd_r(r^2X^r)g^{rr}\right]-\alpha^2(\alpha^{-1}f)' \\
&= \left[2g^{rr}\pd_rX^r-r^{-2}\pd_r(r^2g^{rr}X^r)\right]-\alpha^2(\alpha^{-1}f)' \\
&= 2\alpha\pd_r f -r^{-2}\pd_r\left(r^2\alpha f\right)-\alpha^2(\alpha^{-1}f)' \\
&= -2r^{-1} \alpha f
\end{align*}
and
\begin{align*}
K^{tt}-(\alpha^{-1}f)' &= \left[-X^r\pd_r g^{tt}-r^{-2}\pd_r(r^2X^r)g^{tt}\right]-(\alpha^{-1}f)' \\
&= -r^{-2}\pd_r(r^2 g^{tt} X^r) -(\alpha^{-1}f)' \\
&= -r^{-2}\pd_r\left(r^2(-\alpha^{-1}) f\right)-(\alpha^{-1}f)' \\
&= 2r^{-1}\alpha^{-1}f.
\end{align*}
Finally,
\begin{align*}
\sla{K}^{\alpha\beta} &= -X^r\pd_r \sla{g}^{\alpha\beta}-r^{-2}\pd_r\left(r^2X^r\right)\sla{g}^{\alpha\beta} \\
&= -r^{-2}\pd_r\left(r^2\sla{g}^{\alpha\beta}X^r\right) \\
&= -r^{-2}\pd_r((r^2\sla{g}^{\alpha\beta})f) \\
&= -f' \sla{g}^{\alpha\beta}.
\end{align*}
(We used in the last line that $\pd_r(r^2\sla{g}^{\alpha\beta})=0$.) Combining all these terms gives the identity stated in the lemma. 
\end{proof}

Next, we choose $w=2r^{-1}f$ to directly cancel the middle term in the above lemma.
\begin{lemma}\label{s:divJphiXw_lem}
$$divJ\left[\alpha^{-1}fL,2r^{-1}f\right] \\
= (\alpha^{-1}f)'(L\psi)^2+\left(2r^{-1}f-f'\right)|\sla\nabla\psi|^2-\frac12\Box_g\left(2r^{-1}f\right)\psi^2.$$
\end{lemma}
\begin{proof}
Note that
$$divJ[0,w]=wg^{\mu\nu}\pd_\mu\psi\pd_\nu\psi-\frac12\Box_gw \psi^2.$$
We compute the new terms only.
\begin{align*}
divJ\left[0,2r^{-1}f\right] &= 2r^{-1}fg^{\alpha\beta}\pd_\alpha\psi\pd_\beta\psi -\frac12\Box_g\left(2r^{-1}f\right)\psi^2 \\
&= 2r^{-1}f\left(\alpha (\pd_r\psi)^2-\alpha^{-1}(\pd_t\psi)^2\right) +2r^{-1}f|\sla\nabla\psi| -\frac12\Box_g\left(2r^{-1}f\right)\psi^2.
\end{align*}
When adding these terms to the expression in Lemma \ref{s:divJphiX_lem}, the $\alpha(\pd_r\psi)^2-\alpha^{-1}(\pd_t\psi)^2$ terms cancel (this was the reason for the choice of $w=2r^{-1}f$) and the result is as desired. 
\end{proof}

The term $-\frac12\Box_g\left(2r^{-1}f\right)\psi^2$ is like $-r^{-1}f''\psi^2$. In the future, when $f\sim r^p$, this will have a sign $-p(p-1)$. The sign will be negative if $p>1$, which is bad. So we include a divergence term to fix it. (But in doing so, we almost lose some other good terms--this is why we need a small parameter $\epsilon$.) This is the point of the following Lemma.
\begin{lemma}\label{s:rp_add_term_lem}
\begin{multline*}
\alpha^{-1}f' (L\psi)^2-\frac12 \Box_g\left(2r^{-1}f\right)\psi^2+(1-\epsilon)div\left(\psi^2r^{-1}f'L\right) \\
=\alpha f' \left(\alpha^{-1}L\psi+\frac{1-\epsilon}{r}\psi\right)^2 + \epsilon\alpha \left((1-\epsilon)f'-rf''\right)r^{-2}\psi^2 
+ \alpha'\left(-\epsilon r^{-1}f'+r^{-2}f\right)\psi^2.
\end{multline*}
\end{lemma}
\begin{proof}
First, we calculate
\begin{align*}
-\frac12\Box_g\left(2r^{-1}f\right)\psi^2 &=-r^{-2}\pd_r\left(r^2\alpha\pd_r\left(r^{-1}f\right)\right)\psi^2 \\
&= -\alpha r^{-2}\pd_r\left(r^2\pd_r\left(r^{-1}f\right)\right) \psi^2 -\alpha'\pd_r\left(r^{-1}f\right)\psi^2 \\
&= -\alpha r^{-2}\pd_r(-f+rf')\psi^2-\alpha'\pd_r\left(r^{-1}f\right)\psi^2 \\
&= -\alpha r^{-1} f''\psi^2-\alpha'\pd_r\left(r^{-1}f\right)\psi^2 .
\end{align*}
We also calculate
\begin{align*}
div\left(\psi^2r^{-1}f'L\right) &= r^{-2}\pd_\alpha\left(\psi^2rf' L^\alpha\right) \\
&=r^{-1}f'2\psi L\psi +r^{-2}\pd_r(rf'\alpha)\psi^2 \\
&=r^{-1} f' 2\psi L\psi +r^{-2}\alpha f'\psi^2 +\alpha r^{-1} f''\psi^2+\alpha' r^{-1} f'\psi^2.
\end{align*}
The first two terms in the last line complete a square with the term $\alpha^{-1} f' (L\psi)^2$. The third term cancels with the first term from the previous calculation. However, it will be beneficial to introduce the factor $1-\epsilon$ that appears in the lemma, so that a good term appears with an $\epsilon$ factor. This is summarized by the following two calculations.
\begin{multline*}
\alpha^{-1}f'(L\psi)^2+(1-\epsilon)\left(r^{-1}f'2\psi L\psi+ r^{-2}\alpha f'\psi^2\right) \\
= \alpha f' \left(\alpha^{-1}L\psi+\frac{1-\epsilon}{r}\psi\right)^2 -(1-\epsilon)^2r^{-2}\alpha f'\psi^2 +(1-\epsilon)r^{-2}\alpha f'\psi^2 \\
= \alpha f' \left(\alpha^{-1}L\psi+\frac{1-\epsilon}{r}\psi\right)^2 +\epsilon (1-\epsilon)r^{-2}\alpha f'\psi^2 
\end{multline*}
and
$$-\alpha r^{-1} f''\psi^2 +(1-\epsilon)\alpha r^{-1} f''\psi^2=-\epsilon \alpha r^{-1} f''\psi^2.$$
Adding these terms together yields
$$\alpha f' \left(\alpha^{-1}L\psi+\frac{1-\epsilon}{r}\psi\right)^2 +\epsilon \alpha \left((1-\epsilon)f'-rf''\right)r^{-2}\psi^2.$$
All the remaining terms (which contain a factor of $\alpha'\sim \frac{M}{r^2}$) are
$$
\alpha'\left[-\pd_r\left(r^{-1}f\right)+(1-\epsilon)r^{-1}f'\right]\psi^2 =\alpha'(-\epsilon r^{-1}f+r^{-2}f)\psi^2
$$
Adding both of these yields the result. 
\end{proof}

Thus, we have shown that if $\Box_g\psi=0$, then
\begin{multline*}
divJ\left[\alpha^{-1} f L,2r^{-1}f,(1-\epsilon)r^{-1}f'L\right] \\
= \alpha f' \left(\alpha^{-1}L\psi+\frac{1-\epsilon}{r}\psi\right)^2 + \epsilon\alpha \left((1-\epsilon)f'-rf''\right)r^{-2}\psi^2 +\left(2r^{-1}f-f'\right)|\sla\nabla\psi|^2 \\
- \alpha' \alpha^{-2}f(L\psi)^2
+ \alpha'\left(-\epsilon r^{-1}f'+r^{-2}f\right)\psi^2.
\end{multline*}
If we remove the assumption that $\Box_g\psi=0$, there is an additional term
$$\left(2X(\psi)+w\psi\right)\Box_g\psi = \left(2\alpha^{-1}fL\psi+2r^{-1}f\psi\right)\Box_g\psi.$$

Finally, we turn to the boundary terms. Since we have assumed that $f$ is supported away from the event horizon, it suffices to compute $-J^t$.
\begin{lemma}\label{s:rp_boundary_terms_lem}
\begin{multline*}
-J^t\left[\alpha^{-1}fL,2r^{-1}f,(1-\epsilon)r^{-1}f'L\right] \\
=f\left(\alpha^{-1}L\psi+r^{-1}\psi\right)^2 +\alpha^{-1}f|\sla\nabla\psi|^2+\epsilon r^{-1}f'\psi^2 -\frac1{q^2}\pd_r(rf\psi^2).
\end{multline*}
\end{lemma}
\begin{proof}
We have
\begin{align*}
-J^t[\alpha^{-1}fL] &= -2\pd^t\psi\alpha^{-1}fL\psi+\alpha^{-1}fL^t\pd^\lambda\psi\pd_\lambda\psi \\
&=-2\alpha^{-1}fg^{tt}\pd_t\psi L\psi +\alpha^{-1}f\left(g^{tt}(\pd_t\psi)^2+g^{rr}(\pd_r\psi)^2+|\sla\nabla\psi|^2\right) \\
&=-\alpha^{-1}fg^{tt}(\pd_t\psi)^2-2\alpha^{-1}fg^{tt}\pd_t\psi\alpha\pd_r\psi+\alpha^{-1}fg^{rr}(\pd_r\psi)^2 +\alpha^{-1}f|\sla\nabla\psi|^2 \\
&=\alpha^{-2}f(\pd_t\psi)^2+2\alpha^{-1}f\pd_t\psi\pd_r\psi +f(\pd_r\psi)^2+\alpha^{-1}f|\sla\nabla\psi|^2 \\
&=\alpha^{-2}f(L\psi)^2+\alpha^{-1}f|\sla\nabla\psi|^2.
\end{align*}
Also,
\begin{align*}
-J^t_{(\psi)}\left[0,2r^{-1}f\right] &= -2r^{-1}f\psi\pd^t\psi \\
&= -2r^{-1}fg^{tt}\psi\pd_t\psi \\
&= 2r^{-1}\alpha^{-1}f\psi\pd_t\psi \\
&= 2r^{-1}\alpha^{-1}f\psi L\psi -2r^{-1}f\psi\pd_r\psi \\
&= 2r^{-1}\alpha^{-1}f\psi L\psi -r^{-2}\pd_r(rf\psi^2)+(r^{-2}f+r^{-1}f')\psi^2\\
&= \left(2r^{-1}\alpha^{-1}f\psi L\psi+r^{-2}f\psi^2\right) +r^{-1}f'\psi^2-r^{-2}\pd_r(rf\psi^2).
\end{align*}
Now, observe that
$$\alpha^{-2}f(L\psi)^2+2r^{-1}\alpha^{-1}f\psi L\psi+r^{-2}f\psi^2 = f\left(\alpha^{-1}L\psi+r^{-1}\psi\right)^2.$$
Thus,
$$-J^t\left[\alpha^{-1}fL,2r^{-1}f\right] 
=f\left(\alpha^{-1}L\psi+r^{-1}\psi\right)^2 
+\alpha^{-1}f|\sla\nabla\psi|^2+r^{-1}f'\psi^2-r^{-2}\pd_r(rf\psi^2).$$
Also,
$$-J^t\left[0,0,(1-\epsilon)r^{-1}f'L\right] = -(1-\epsilon)r^{-1}f'\psi^2 L^t =-(1-\epsilon)r^{-1}f'\psi^2.$$
Adding these two expressions together yields the result. 
\end{proof}

This concludes the proof of Lemma \ref{s:p_ee_identity_lem}.
\end{proof}

\subsection{The incomplete $r^p$ estimate near $i^0$}

Now we use Lemma \ref{s:p_ee_identity_lem} and make a particular choice for the funciton $f$ (so that $f=r^p$ for large $r$) to prove the following.
\begin{proposition}\label{s:incomplete_p_estimate_prop}
Fix $\delm,\delp>0$. Let $R$ be a sufficiently large radius. Then for all $p\in[\delm,2-\delp]$, the following estimate holds if $\psi$ decays sufficiently fast as $r\rightarrow\infty$.
\begin{multline*}
\int_{\Sigma_{t_2}\cap\{r>2R\}}r^p\left[(L\psi)^2+|\sla\nabla\psi|^2+r^{-2}\psi^2\right] \\
+ \int_{t_1}^{t_2}\int_{\Sigma_t\cap\{r>2R\}}r^{p-1}\left[(L\psi)^2+|\sla\nabla\psi|^2+r^{-2}\psi^2\right] \\
\lesssim \int_{\Sigma_{t_2}\cap\{r>2R\}}r^p\left[(L\psi)^2+|\sla\nabla\psi|^2+r^{-2}\psi^2\right] + Err,
\end{multline*}
where
\begin{align*}
Err &= Err_1 + Err_\Box \\
Err_1 &= \int_{t_1}^{t_2}\int_{\Sigma_t\cap\{R<r<2R\}}(L\psi)^2+|\sla\nabla\psi|^2+\psi^2  \\
Err_\Box &= \int_{t_1}^{t_2}\int_{\Sigma_t\cap\{R<r\}}r^p(|L\psi|+r^{-1}|\psi|)|\Box_g\psi|.
\end{align*}
\end{proposition}
\begin{proof}
The estimate follows from the identity given in Lemma \ref{s:p_ee_identity_lem} and a particular choice for the function $f$.
$$f(r)=\rho^p,$$
where
$$
\rho = \left\{
\begin{array}{ll}
0 & r\le R \\
smooth, \rho'>0 & r\in[R,2R] \\
r & 2R< r.
\end{array}
\right.
$$
With this choice, we have 
$$f\ge 0$$
$$f'\ge 0$$
and for $r>2R$,
$$f = r^p$$
$$f'=pr^{p-1}.$$
Furthermore, for $r>2R$,
$$2r^{-1}f-f' = 2r^{p-1}-pr^{p-1} = (2-p)r^{p-1}.$$
It follows that if $p\le 2-\delp$, then for $r>2R$,
$$r^{p-1}\lesssim 2r^{-1}f-f'.$$
Also, for $r>2R$,
$$\epsilon\alpha ((1-\epsilon)f'-rf'') = \epsilon\alpha ((1-\epsilon)pr^{p-1}-p(p-1)r^{p-1}) = \epsilon\alpha p (2-\epsilon -p)r^{p-1}.$$
If $R$ is sufficiently large so that $\alpha>3/4$ and $\delm\le p\le 2-\delp$ and $\epsilon\le \delp/2$, then
$$\epsilon r^{p-1} \lesssim \epsilon\alpha ((1-\epsilon)f'-rf'').$$

We also note that there are some error terms that have a factor of $\alpha'$. Each of these terms has a smallness parameter available, since $R$ can be taken to be very large and
$$\alpha'\lesssim \frac{M}{R}r^{-1}.$$

Finally, we observe that if $\psi$ vanishes sufficiently fast as $r\rightarrow\infty$, then since $f$ is supported for $r>R$, we have
$$\int_{\Sigma_t}-r^{-2}\pd_r(rf\psi^2) =0.$$
With these facts having been established, it is straightforward to check that the estimate follows from Lemma \ref{s:p_ee_identity_lem}.
\end{proof}

\subsection{The $r^p$ estimate}

We conclude this section by proving the $r^p$ estimate. This is a combination of the $h\pd_t$ estimate (Proposition \ref{s:hdt_prop}), the Morawetz estimate (Theorem \ref{s:morawetz_thm}), and the incomplete $r^p$ estimate (Proposition \ref{s:incomplete_p_estimate_prop}).
\begin{proposition}\label{s:rp_prop}
Fix $\delm,\delp>0$ and let $p\in[\delm,2-\delp]$. Then if
$\psi$ decays sufficiently fast as $r\rightarrow\infty$, the following estimate holds.
\begin{multline*}
\int_{\Sigma_{t_2}}r^p\left[(L\psi)^2+|\sla\nabla\psi|^2+r^{-2}\psi^2 + r^{-2}(\pd_r\psi)^2\right] \\
+ \int_{t_1}^{t_2}\int_{\Sigma_t}r^{p-1}\left[\chi_{trap}(L\psi)^2+\chi_{trap}|\sla\nabla\psi|^2+r^{-2}\psi^2+r^{-2}(\pd_r\psi)^2\right] \\
\lesssim \int_{\Sigma_{t_1}}r^p\left[(L\psi)^2+|\sla\nabla\psi|^2+r^{-2}\psi^2+r^{-2}(\pd_r\psi)^2\right] 
+Err_\Box,
\end{multline*}
where $\chi_{trap}=\left(1-\frac{r_{trap}}{r}\right)^2$ and
\begin{align*}
Err_\Box &= \int_{t_1}^{t_2}\int_{\Sigma_t}|(2X(\psi)+w\psi)\Box_g\psi|,
\end{align*}
where the vectorfield $X$ and function $w$ satisfy the following properties. \\
\bp $X$ is everywhere timelike, but asymptotically null at the rate $X=O(r^p)L+O(r^{p-2})\pd_t$. \\
\bp $X|_{r=r_H}=-\lambda\pd_r$ for some positive constant $\lambda$. \\
\bp $X|_{r=r_{trap}}=\lambda\pd_t$ for some positive constant $\lambda$. \\
\bp $w =O(r^{p-1})$ for large $r$.
\end{proposition}
\begin{proof}
We start with the Morawetz estimate (Theorem \ref{s:morawetz_thm}) and add a small constant times the incomplete $r^p$ estimates (Proposition \ref{s:incomplete_p_estimate_prop}). The small constant can be chosen so that the bulk error term $Err_1$ from Proposition \ref{s:incomplete_p_estimate_prop} can be absorbed into the bulk in the Morawetz estimate. The result is the following estimate.
\begin{multline*}
\int_{\Sigma_{t_2}}r^p\left[(L\psi)^2+|\sla\nabla\psi|^2+r^{-2}\psi^2\right]+\frac{M^2}{r^2}(\pd_r\psi)^2 \\
\hspace{1in}+\int_{t_1}^{t_2}\int_{\Sigma_t} r^{p-1}\left[\chi_{trap}(L\psi)^2+\chi_{trap}|\sla\nabla\psi|^2+r^{-2}\psi^2\right]+\frac{M^2}{r^3}(\pd_r\psi)^2 \\
\lesssim \int_{\Sigma_{t_1}}r^p\left[(L\psi)^2+|\sla\nabla\psi|^2+r^{-2}\psi^2\right]+\frac{M^2}{r^2}(\pd_r\psi)^2 + Err'
\end{multline*}
where
\begin{align*}
Err' &= Err'_1+Err'_2+Err'_\Box \\
Err'_1 &= \int_{\Sigma_{t_2}}r^{-1}|\psi L\psi| \\
Err'_2 &= \int_{H_{t_1}^{t_2}}(\pd_t\psi)^2 +\int_{\Sigma_{t_2}}\frac{M^2}{r^2}\left[\chi_H(\pd_r\psi)^2+(\pd_t\psi)^2\right] \\
Err'_\Box &= \int_{t_1}^{t_2}\int_{\Sigma_t\cap\{R<r\}}r^p(|L\psi|+r^{-1}|\psi|)|\Box_g\psi| \\
&\hspace{1in} +\int_{t_1}^{t_2}\int_{\Sigma_t}|(2X'(\psi)+w'\psi)|\Box_{g}\psi|,
\end{align*}
and $X'$ and $w'$ are the vectorfield and function defined in the Morawetz estimate (Theorem \ref{s:morawetz_thm}).

The error term $Err'_1$ can in fact be removed due to the following argument.
\begin{align*}
Err'_1 &\lesssim \int_{\Sigma_{t_2}}\epsilon r^p(L\psi)^2+\epsilon^{-1}r^{-p}r^{-2}\psi^2 \\
&\lesssim \int_{\Sigma_{t_2}}\epsilon r^p[(L\psi)^2+r^{-2}\psi^2] +\int_{\Sigma_{t_2}\cap\{r\le R_\epsilon\}}\epsilon^{-1}\psi^2 .
\end{align*}
The radius $R_\epsilon$ should be chosen sufficiently large so that $\epsilon^{-1}r^{-p}\le \epsilon r^p$ whenever $r>R_\epsilon$. This critically depends on the fact that $p\ge\delm>0$. Now, the parameter $\epsilon$ can be taken sufficiently small so as to absorb the first two terms into the left side of the main estimate and the last two terms can be included with the term $Err'_2$ after applying a Hardy estimate.

We return to the main estimate. Notice that most terms have improved weights near $i^0$ and a few error terms remain on $H_{t_1}^{t_2}$ and $\Sigma_{t_2}$. The next step is to use the $h\pd_t$ estimate (Proposition \ref{s:hdt_prop}) to eliminate these error terms and improve the weights near $i^0$ for the remaining $\pd_r\psi$ terms. The result is the following estimate.
\begin{multline*}
\int_{\Sigma_{t_2}}r^p\left[(L\psi)^2+|\sla\nabla\psi|^2+r^{-2}\psi^2+r^{-2}(\pd_r\psi)^2\right] \\
+\int_{t_1}^{t_2}\int_{\Sigma_t} r^{p-1}\left[\chi_{trap}(L\psi)^2+\chi_{trap}|\sla\nabla\psi|^2+r^{-2}\psi^2+c_\epsilon r^{-2}(\pd_r\psi)^2\right] \\
\lesssim \int_{\Sigma_{t_1}}r^p\left[(L\psi)^2+|\sla\nabla\psi|^2+r^{-2}\psi^2+r^2(\pd_r\psi)^2\right] +Err'',
\end{multline*}
where
\begin{align*}
Err'' &= Err''_1+Err''_\Box \\
Err''_1 &= \int_{t_1}^{t_2}\int_{\Sigma_t\cap\{5M<r\}} \epsilon r^{-1}(L\psi)^2 \\
Err''_\Box &= \int_{t_1}^{t_2}\int_{\Sigma_t\cap\{R<r\}}r^p(|L\psi|+r^{-1}|\psi|)|\Box_g\psi| + \int_{t_1}^{t_2}\int_{\Sigma_t}|(2X'(\psi)+w'\psi)(\Box_g\psi)| \\
&\hspace{3.5in}+ \int_{t_1}^{t_2}\int_{\Sigma_t}C_\epsilon r^{p-2}|\pd_t\psi(\Box_g\psi)|.
\end{align*}
By taking $\epsilon$ sufficiently small, the error term $Err''_1$ can be absorbed into the left side. 

The resulting vectorfield $X$ and function $w$ can be understood by combining the vectorfields and scalar functions used to construct the three estimates that were used in this proof.
\end{proof}

\section{The dynamic estimates}\label{s:dynamic_sec}

In this section, we prove the dynamic estimates, which provide all of the necessary information related to the future dynamics of the wavefunction $\psi$. These estimates are simply restatements of the energy estimate (Proposition \ref{s:classic_ee_prop}) and the $r^p$ estimate (Proposition \ref{s:rp_prop}) applied to either $\psi$ itself or another wavefunction $\psi^s$ derived from $\psi$.

\subsection{The dynamic estimates for $\psi$ ($s=0$)}

We begin with the dynamic estimates for the wavefunction $\psi$ only.
\begin{proposition}\label{s:dynamic_estimates_0_prop}
Fix $\delm,\delp>0$ and let $p\in[\delm,2-\delp]$. Then
$$E(t_2)\lesssim E(t_1)+\int_{t_1}^{t_2}N(t)dt$$
$$E_p(t_2)+\int_{t_1}^{t_2}B_p(t)dt \lesssim E_p(t_1)+\int_{t_1}^{t_2}N_p(t)dt,$$
where
$$E(t)=\int_{\Sigma_t}\chi_H(\pd_r\psi)^2+(\pd_t\psi)^2+|\sla\nabla\psi|^2,$$
$$E_p(t)=\int_{\Sigma_t}r^p\left[(L\psi)^2+|\sla\nabla\psi|^2+r^{-2}\psi^2+r^{-2}(\pd_r\psi)^2\right],$$
$$B_p(t)=\int_{\Sigma_t}r^{p-1}\left[\chi_{trap}(L\psi)^2+\chi_{trap}|\sla\nabla\psi|^2+r^{-2}\psi^2+r^{-2}(\pd_r\psi)^2\right],$$
$$N(t)=\int_{\Sigma_t}|\pd_t\psi\Box_g\psi|,$$
$$N_p(t)=\int_{\Sigma_t}|(2X(\psi)+w\psi)\Box_g\psi|,$$
where $\chi_H=1-\frac{2M}r$, $\chi_{trap}=\left(1-\frac{3M}r\right)^2$, and $X$ and $w$ are the vectorfield and function from Proposition \ref{s:rp_prop}.
\end{proposition}
\begin{proof}
The first estimate is a restatement of the energy estimate (Proposition \ref{s:classic_ee_prop}) and the second is a restatement of the $r^p$ estimate (Proposition \ref{s:rp_prop}).
\end{proof}

We slightly simplify the previous estimates by absorbing part of the nonlinear norm $N_p(t)$ into the bulk $B_p(t)$ on the left side. Note the complication due to trapping.
\begin{corollary}\label{s:dynamic_estimates_0_cor}
Fix $\delm,\delp>0$ and let $p\in[\delm,2-\delp]$. Then
$$E(t_2)\lesssim E(t_1)+\int_{t_1}^{t_2}N(t)dt$$
$$E_p(t_2)+\int_{t_1}^{t_2}B_p(t)dt \lesssim E_p(t_1)+\int_{t_1}^{t_2}N_p(t)dt,$$
where $E(t)$, $E_p(t)$, and $B_p(t)$ are as defined in Proposition \ref{s:dynamic_estimates_0_prop} and
$$N(t)=(E(t))^{1/2}||\Box_g\psi||_{L^2(\Sigma_t)},$$
$$N_p(t)=\int_{\Sigma_t}r^{p+1}(\Box_g\psi)^2+\int_{\Sigma_t\cap\{r\approx 3M\}}|\pd_t\psi\Box_g\psi|.$$
\end{corollary}
\begin{proof}
The first estimate is due to Proposition \ref{s:dynamic_estimates_0_prop} and the fact that
$$\int_{\Sigma_t}|\pd_t\psi\Box_g\psi|\lesssim ||\pd_t\psi||_{L^2(\Sigma_t)}||\Box_g\psi||_{L^2(\Sigma_t)} \lesssim (E(t))^{1/2}||\Box_g\psi||_{L^2(\Sigma_t)}.$$
We turn to the second estimate.
According to Proposition \ref{s:dynamic_estimates_0_prop},
$$E_p(t_2)+\int_{t_1}^{t_2}B_p(t)dt \lesssim E_p(t_1)+\int_{t_1}^{t_2}N_p'(t)dt,$$
where
\begin{align*}
N_p'(t) &= \int_{\Sigma_t}|(2X(\psi)+w\psi)\Box_g\psi| \\
&\lesssim \int_{\Sigma_t}r^p(|L\psi|+r^{-1}|\psi|+r^{-2}|\pd_r\psi|)|\Box_g\psi|.
\end{align*}
We would like to claim that
$$N_p'(t)\lesssim \epsilon B_p(t)+\epsilon^{-1}\int_{\Sigma_t}r^{p+1}(\Box_g\psi)^2.$$
\textbf{This would require an estimate of the form}
$$\int_{\Sigma_t}r^{p-1}(|L\psi|+r^{-1}|\psi|+r^{-2}|\pd_r\psi|)^2\lesssim B_p(t),$$
\textbf{which is not completely true due to the loss of $(L\psi)^2$ at the trapping radius}. Instead,
$$\int_{\Sigma_t}r^{p-1}(|L\psi|+r^{-1}|\psi|+r^{-2}|\pd_r\psi|)^2\lesssim B_p(t)+\int_{\Sigma_t\cap\{r\approx 3M\}}|\pd_t\psi\Box_g\psi|,$$
and therefore
$$N_p'(t)\lesssim \epsilon B_p(t)+\epsilon^{-1}\int_{\Sigma_t}r^{p+1}(\Box_g\psi)^2+\epsilon\int_{\Sigma_t\cap\{r\approx 3M\}}|\pd_t\psi\Box_g\psi|.$$
By taking $\epsilon$ sufficiently small, the term $\epsilon \int_{t_1}^{t_2}B_p(t)dt$ can be absorbed into the left side of the second estimate.
\end{proof}

\subsection{Commutators with $\Box_g$}

There are substantially fewer symmetries in Schwarzschild than in Minkowski. They are generated by the time translation vectorfield $\pd_t$ and the rotation operators $\Omega_x$, $\Omega_y$, and $\Omega_z$.

\begin{definition}
Denote by $\Gamma$ the time translation operator or any of the rotation operators.
$$\Gamma\in \{\pd_t,\Omega_x,\Omega_y,\Omega_z\}$$
Furthermore, denote by $\Gamma^s$ any composition of $s$ of these operators.
\end{definition}

We also define the $s$-order wavefunctions derived from the wavefunction $\psi$ by applying symmetry operators.
\begin{definition}
$$\psi^s=\Gamma^s\psi$$
\end{definition}

As previously discussed, these generalized wavefunctions behave very similarly to $\psi$. This is the statement of the following lemma.
\begin{lemma}\label{s:s_comm_lem}
If $\Box_g\psi=0$, then
$$\Box_g\psi^s=0.$$
More generally,
$$\Box_g\psi^s=\Gamma^s(\Box_g\psi).$$
\end{lemma}
\begin{proof}
We only prove the general case.
$$\Box_g\psi^s=\Box_g\Gamma^s\psi=\Gamma^s(\Box_g\psi)+[\Box_g,\Gamma^s]\psi=\Gamma^s(\Box_g\psi).$$
\end{proof}

\subsection{Higher order dynamic estimates}

Now, we can write the dynamic estimates in a more general form.
\begin{corollary}\label{s:s_cor}
Fix $\delm,\delp>0$. The following estimates hold for $p\in [\delm,2-\delp]$ and every integer $s\ge 0$.
$$E^s(t_2)\lesssim E^s(t_1)+\int_{t_1}^{t_2}N^s(t)dt,$$
$$E_p^s(t_2)+\int_{t_1}^{t_2}B_p^s(t)dt\lesssim E_p^s(t_1)+\int_{t_1}^{t_2}N_p^s(t)dt,$$
where
$$E^s(t)=\sum_{s'\le s} E[\psi^{s'}](t),$$
$$E_p^s(t)=\sum_{s'\le s} E_p[\psi^{s'}](t),$$
$$B_p^s(t)=\sum_{s'\le s} B_p[\psi^{s'}](t),$$
$$N^s(t)= (E^s(t))^{1/2}\sum_{s'\le s}||q^{-2}\Gamma^{s'}(q^2\Box_g\psi)||_{L^2(\Sigma_t)},$$
$$N_p^s(t)=\sum_{s'\le s} \int_{\Sigma_t}r^{p+1}(q^{-2}\Gamma^{\le s}(q^2\Box_g\psi))^2+\sum_{s'\le s}\int_{\Sigma_t\cap\{r\approx 3M\}}|\pd_t\psi^{s'}q^{-2}\Gamma^{s'}(q^2\Box_g\psi)|,$$
and the norms $E(t)$, $E_p(t)$, and $B_p(t)$ are as defined in Proposition \ref{s:dynamic_estimates_0_prop}.
\end{corollary}
\begin{proof}
The proof is a direct application of Corollary \ref{s:dynamic_estimates_0_cor} by making the substitutions
$$\psi\mapsto\psi^{s'}$$
for all values of $s'$ (and all commutators represented by $\Gamma^{s'}$) where $s'\le s$ and observing Lemma \ref{s:s_comm_lem}.
\end{proof}

\section{The $L^\infty$ estimates}\label{s:pointwise_sec}

We now prove the $L^\infty$ estimates. These are very similar to the $L^\infty$ estimates in Minkowski, with the main differences related to the presence of an event horizon in Schwarzschild.

\subsection{A sobolev-type estimate}

We restate Lemma \ref{m:sobolev_I_lem}.
\begin{lemma}\label{s:pointwise_lem}
If $u$ decays sufficiently fast as $r\rightarrow\infty$, then
$$||u||^2_{L^\infty(\Sigma_t\cap\{r\ge r_0\})}\lesssim \int_{\Sigma_t\cap\{r\ge r_0\}}r^{-2}\left[(\pd_r \Omega^{\le 2}u)^2+(\Omega^{\le 2}u)^2\right]$$
\end{lemma}

\begin{proof}
See the proof of Lemma \ref{m:sobolev_I_lem}.
\end{proof}

\subsection{Estimating derivatives using the Sobolev-type estimate}

Now, we repeatedly apply Lemma \ref{s:pointwise_lem} to estimate various derivatives with $r$ weights. We will assume that $\psi$ decays sufficiently fast as $r\rightarrow\infty$.

The following lemma estimates $\psi$ and the higher order analogues $\psi^s$.

\begin{lemma}\label{s:psi_pointwise_lem}
For $r\ge r_H$,
$$|r^p\psi^s|^2\lesssim E^{s+2}_{2p}(t)$$
and for $r\ge r_0>r_H$,
$$|\psi^s|^2\lesssim E^{s+2}(t).$$
\end{lemma}
\begin{proof}
First, we apply Lemma \ref{s:pointwise_lem} with $u=r^p\psi^{s}$.
\begin{align*}
|r^p\psi^{s}|^2 &\lesssim \int_{\Sigma_t} r^{-2}\left[(\pd_r\Omega^{\le 2}(r^p\psi^{s}))^2+(\Omega^{\le 2}(r^p\psi^{s}))^2\right] \\
&\lesssim \int_{\Sigma_t}r^{2p-2}\left[(\pd_r\psi^{s+2})^2+(\psi^{s+2})^2\right] \\
&\lesssim E_{2p}^{s+2}(t).
\end{align*}
This verifies the first estimate. The second estimate follows from the exact same argument in the special case $p=0$, and the observation that as long as $r\ge r_0>r_H$, then $E^{s+2}(t)$ can be used in place of of $E_0^{s+2}(t)$.
\end{proof}

The following lemma estimates $\pd_t\psi$ and the higher order analogues $\pd_t\psi^{s}$.
\begin{lemma}
For $r\ge r_H$,
$$|r^p\pd_t\psi^{s}|^2\lesssim E_{2p}^{s+3}(t)$$
and for $r\ge r_0>r_H$,
$$|r\pd_t\psi^{s}|^2\lesssim E^{s+3}(t)$$
\end{lemma}
\begin{proof}
The first estimate reduces to Lemma \ref{s:psi_pointwise_lem} by observing that $\pd_t\psi^{s}=\psi^{s+1}$. We now prove the second estimate.

First, we apply Lemma \ref{s:pointwise_lem} with $u=r\pd_t\psi^{s}$.
\begin{align*}
|r\pd_t\psi^{s}|^2 &\lesssim \int_{\Sigma_t\cap\{r> r_0\}}r^{-2}\left[(\pd_r\Omega^{\le 2}(r\pd_t\psi^{s}))^2+(\Omega^{\le 2}(r\pd_t\psi^{s}))^2\right] \\
&\lesssim \int_{\Sigma_t\cap\{r> r_0\}}\left[(\pd_r\Omega^{\le 2}(\pd_t\psi^{s}))^2+(\Omega^{\le 2}(\pd_t\psi^{s}))^2\right] \\
&\lesssim \int_{\Sigma_t\cap\{r> r_0\}}\left[(\pd_r\psi^{s+3})^2+(\pd_t\psi^{s+2})^2\right] \\
&\lesssim E^{s+2}(t).
\end{align*}
The point is that $\pd_t\psi^{s+2}$ can either be treated as $\psi^{s+3}$ or as a derivative of $\psi^{s+2}$. In the latter case, the energy norm $E^s(t)$ has stronger control, because $||\pd_t\psi^s||_{L^2(\Sigma_t\cap\{r>r_0\})}^2\le E^s(t)$, while $||r^{-1}\psi^{s+1}||_{L^2(\Sigma_t\cap\{r>r_0\})}^2\le E^{s+1}(t)$.
\end{proof}

The following lemma estimates $\sla\nabla\psi$ and the higher order analogues $\sla\nabla\psi^{s}$.
\begin{lemma}
For $r\ge r_H$,
$$|r^{p+1}\sla\nabla\psi^{s}|^2\lesssim E^{s+3}_{2p}(t)$$
and for $r\ge r_0>r_H$,
$$|r\sla\nabla\psi^{s}|^2\lesssim E^{s+3}(t)$$
\end{lemma}
\begin{proof}
This lemma reduces to Lemma \ref{s:psi_pointwise_lem} by observing that $|r\sla\nabla\psi^{s}|\lesssim |\Omega\psi^{s}|=|\psi^{s+1}|$.
\end{proof}

The following lemma estimates $L\psi$ and the higher order analogues $L\psi^{s}$.
\begin{lemma}\label{s:L_pointwise_lem}
Letting $L=\alpha\pd_r+\pd_t$, where $\alpha=1-\frac{2M}r$, we have that for $r\ge r_H$,
$$|r^{p+1}L\psi^{s}|^2\lesssim E_{2p}^{s+3}(t)+\int_{\Sigma_t}r^{2p}(\Box_g\psi^{s+2})^2$$
and for $r\ge r_0>r_H$,
$$|rL\psi^{s}|^2\lesssim E^{s+3}(t)+\int_{\Sigma_t\cap\{r>r_0\}}(\Box_g\psi^{s+2})^2.$$
\end{lemma}
\begin{proof}
Before beginning the estimates stated by the lemma, it is important to establish
$$(\pd_rLu)^2\lesssim (\Box_gu)^2+(L\pd_tu)^2+r^{-2}(\pd_t^2u)^2+r^{-2}(\pd_r\pd_tu)^2+r^{-2}(\pd_ru)^2+(\sla\triangle u)^2,$$
To verify this, we expand
$$\Box_g=g^{tt}\pd_t^2+g^{rr}\pd_r^2+r^{-2}\pd_r(r^2g^{rr})\pd_r+g^{rt}\pd_r\pd_t+r^{-2}\pd_r(r^2g^{rt})\pd_t+\sla\triangle$$
and observe that for $r>r_H+\delh$, $g^{rt}=0$ and
\begin{multline*}
g^{tt}\pd_t^2u+g^{rr}\pd_r^2u+r^{-2}\pd_r(r^2g^{rr})\pd_ru \\
=-\left(1-\frac{2M}r\right)^{-1}\pd_t^2u+\left(1-\frac{2M}r\right)\pd_r^2u+\frac{2(r-M)}{r^2}\pd_ru.
\end{multline*}
We also expand
\begin{align*}
\pd_rLu &= \pd_r\left(\left(1-\frac{2M}r\right)\pd_ru+\pd_tu\right) \\
&= \left(1-\frac{2M}r\right)\pd_r^2u+\frac{2M}{r^2}\pd_ru+\pd_r\pd_tu
\end{align*}
Note the similarities between these two expressions. It follows that for $r\ge r_H+\delh$,
\begin{multline*}
g^{tt}\pd_t^2u+g^{rr}\pd_r^2u+r^{-2}\pd_r(r^2g^{rr})\pd_ru - \pd_rLu \\
=-\left(1-\frac{2M}r\right)^{-1}\pd_t^2u-\pd_r\pd_tu + \frac{2}r\left(1-\frac{2M}r\right)\pd_ru \\
=-\alpha^{-1}\left(\pd_t^2u+\alpha\pd_r\pd_tu\right)+\frac2{r}\alpha\pd_ru \\
=-\alpha^{-1}L\pd_tu+\frac{2}r\alpha\pd_ru.
\end{multline*}
Keeping in mind the additional terms that show up for $r\le r_H+ \delh$, we arrive at the estimate for $(\pd_rLu)^2$ at the beginning of this proof. With this estimate in mind, we begin to prove the estimates stated by the lemma.

We apply Lemma \ref{s:pointwise_lem} with $u=r^{p+1}L\psi^{s}$.
\begin{align*}
|r^{p+1}L\psi^{s}|^2 
&\lesssim \int_{\Sigma_t}r^{-2}\left[(\pd_r\Omega^{\le 2}(r^{p+1}L\psi^{s}))^2+(\Omega^{\le 2}(r^{p+1}L\psi^{s}))^2\right] \\
&\lesssim \int_{\Sigma_t}r^{2p}\left[(\pd_r\Omega^{\le 2}L\psi^{s})^2+(\Omega^{\le 2}L\psi^{s})^2\right] \\
&\lesssim \int_{\Sigma_t}r^{2p}\left[(\pd_rL\psi^{s+2})^2+(L\psi^{s+2})^2\right] \\
&\lesssim E_{2p}^{s+2}(t) +\int_{\Sigma_t}r^{2p}(\pd_rL\psi^{s+2})^2.
\end{align*}
Now, according to the estimate we previously established,
\begin{align*}
\int_{\Sigma_t}r^{2p}&(\pd_rL\psi^{s+2})^2 \\
&\lesssim \int_{\Sigma_t}r^{2p}\left[(\Box_g\psi^{s+2})^2+(L\pd_t\psi^{s+2})^2+r^{-2}(\pd_t^2\psi^{s+2})^2+r^{-2}(\pd_r\pd_t\psi^{s+2})^2\right.\\
&\hspace{3in}\left.+r^{-2}(\pd_r\psi^{s+2})^2+(\sla\triangle\psi^{s+2})^2\right] \\
&\lesssim \int_{\Sigma_t}r^{2p}\left[(\Box_g\psi^{s+2})^2+(L\psi^{s+3})^2+r^{-2}(\psi^{s+3})^2+r^{-2}(\pd_r\psi^{s+3})^2\right].
\end{align*}
It follows that
$$|r^{p+1}L\psi^{s}|^2\lesssim E_{2p}^{s+3}(t)+\int_{\Sigma_t}r^{2p}(\Box_g\psi^{s+2})^2.$$
This is the first estimate of the lemma. The second estimate follows from the same exact argument in the special case $p=0$, and the observation that as long as $r\ge r_0>r_H$, then $E^{s}(t)$ can be used in place of $E_0^{s}(t)$.
\end{proof}

\subsection{Summarizing the $L^\infty$ estimates}

To conclude this section, we summarize the previous lemmas in a single proposition, making use of the following definition.
\begin{definition}
We define two families of operators
$$\bar{D}=\{L,\sla\nabla\}$$
$$D=\{L,\alpha\pd_r,\sla\nabla\}$$
where $\alpha = 1-\frac{2M}r$.
\end{definition}

\begin{proposition}\label{s:infinity_prop}
For $r\ge r_H$,
$$|r^{p+1}\bar{D}\psi^{s}|^2+|r^pD\psi^{s}|^2+|r^p\psi^{s}|^2\lesssim E_{2p}^{s+3}(t)+\int_{\Sigma_t}r^{2p}(\Box_g\psi^{s+2})^2,$$
and for $r\ge r_0>r_H$,
$$|rD\psi^{s}|^2\lesssim E^{s+3}(t)+\int_{\Sigma_t}(\Box_g\psi^{s+2})^2.$$
\end{proposition}

\begin{proof}
With the exception of $\alpha\pd_r$, all of the cases have been proved in Lemmas \ref{s:psi_pointwise_lem}-\ref{s:L_pointwise_lem}. Finally, observe that since $\alpha\pd_r=L-\pd_t$,
$$|r^p\alpha\pd_r\psi^{s}|^2\le |r^pL\psi^{s}|^2+|r^p\pd_t\psi^{s}|^2$$
Thus, even the case of the operator $\alpha\pd_r$ can be reduced to Lemmas \ref{s:psi_pointwise_lem}-\ref{s:L_pointwise_lem}.
\end{proof}

\begin{remark}
The presense of the factor $\alpha=1-\frac{2M}r$ is bad, because it means that we cannot control $|\pd_r\psi^s|$ up to the event horizon. The reason for this problem is that $L$ and $\pd_t$ become parallel at the event horizon. It is possible to control $|\pd_r\psi^s|$ up to the event horizon by using an additional operator which is like $\pd_r$ near the event horizon. But since this operator does not commute with the wave operator, it introduces serious complications to the proof. For the sake of simplicity, in this chapter, we will avoid use of such operator and necessarily make a stronger assumption about the nonlinear term near the event horizon. However, in the next chapter, we will use this operator. It will be called $\tg$, and the index $k$ will be used to refer to the number of times $\tg$ is applied to $\psi^s$.
\end{remark}

\section{Theorem: Global boundedness and decay for solutions to the semilinear wave equation with null structure on the Schwarzschild background}\label{s:main_sec}
\subsection{Structure of the nonlinear term}\label{s:nonlinear_structure_sec}

Now, we define the precise structure of the nonlinear term that will be assumed by the main theorem of this chapter.

\begin{definition}\label{s:nonlinear_def}
Suppose $\Box_g\psi=\mathcal{N}$, where $\mathcal{N}$ is nonlinear in zeroth and first order derivatives of $\psi$. We say $\mathcal{N}$ satisfies the null condition if it can be written as
$$\mathcal{N}=\sum \gamma\beta,$$
where the sum has a finite number of terms and
$$\gamma\in\{L\psi,r^{-1}\alpha\pd_r\psi,\sla\nabla\psi,r^{-1}\psi\},$$
$$\beta\in\{L\psi,\alpha\pd_r\psi,\sla\nabla\psi,r^{-1}\psi\},$$
where $\alpha=1-\frac{2M}r$. The $\gamma$ factors can be thought of as ``good'' factors, because they can be estimated with an extra $r$ factor using the $L^\infty$ estimates (Proposition \ref{s:infinity_prop}) while, in contrast, the $\beta$ factors can be thought of as potentially ``bad''. The null condition requires that at least one factor be a good factor.
\end{definition}

\begin{proposition}\label{s:commutator_prop}
Generalize the terms $\gamma$ and $\beta$ as follows.
$$\gamma^{s}\in\{L\psi^{s},r^{-1}\alpha\pd_r\psi^{s},\sla\nabla\psi^{s},r^{-1}\psi^{s}\},$$
$$\beta^{s}\in\{L\psi^{s},\pd_r\alpha\psi^{s},\sla\nabla\psi^{s},r^{-1}\psi^{s}\}.$$
If $\Box_g\psi=\mathcal{N}$ and $\mathcal{N}$ satisfies the null condition, then
$$|\Gamma^s\Box_g\psi|\lesssim \sum |\gamma^{hi}\beta^{lo}|+\sum |\beta^{hi}\gamma^{lo}|,$$
where $\gamma^{hi}$ and $\beta^{hi}$ represent $\gamma^{s'}$ and $\beta^{s'}$ with $s'\le s$, while $\gamma^{lo}$ and $\beta^{lo}$ represent $\gamma^{s'}$ and $\beta^{s'}$ with $s'\le s/2$.
\end{proposition}

\begin{proof}
See the proof of Lemma \ref{m:null_cond_lem}.
\end{proof}

\subsection{The main theorem}\label{s:main_thm_sec}

\begin{theorem}\label{s:main_thm}
Let $g$ be a Schwarzschild metric, and suppose a function $\psi$ solves an equation of the form
$$\Box_g\psi=\mathcal{N},$$
where $\mathcal{N}$ satisfies the null condition as defined in \S\ref{s:nonlinear_structure_sec}. Then for $\delp,\delm>0$ sufficiently small, if the initial data on $\Sigma_0$ decay sufficiently fast as $r\rightarrow\infty$ and have size
\begin{equation}
I_0=E^{12}(0)+E_{2-\delp}^{12}(0)
\end{equation}
sufficiently small, then the following estimates hold for $t\ge 0$ (with $T=1+t$).

I) The energies satisfy
$$E^{12}(t)\lesssim I_0,$$
$$E^{12}_{p\in[\delm,2-\delp]}(t)\lesssim I_0,$$
$$E^{11}_{p\in[1-\delp,2-\delp]}(t)\lesssim T^{p-2+\delp}I_0,$$
$$E^{10}_{p\in[\delm,2-\delp]}(t)\lesssim T^{p-2+\delp}I_0,$$
$$\int_{t}^{\infty}E_{p\in[\delm-1,\delm]}^{9}(\tau)d\tau\lesssim T^{p-2+\delp+1}I_0.$$

II) The following $L^\infty$ estimates hold, provided $s\le 12$.
$$|r^{p+1}\bar{D} \psi^{s}|^2+|r^pD\psi^{s}|^2+|r^p\psi^{s}|^2\lesssim E^{s+3}_{2p}(t),$$
$$|rD\psi^{s}|^2\lesssim E^{s+3}(t).$$

III) Together, (I) and (II) imply that if $s\le 9$, for all $p\in [\delm/2,(2-\delp)/2]$,
$$|r^{p+1}\bar{D} \psi^{s}|+|r^pD\psi^{s}|+|r^p\psi^{s}|\lesssim T^{(2p-2+\delp)/2}I_0^{1/2},$$
$$|rD\psi^{s}|\lesssim I_0^{1/2},$$
and additionally for $p\in [(\delm-1)/2,\delm/2]$,
$$\int_t^\infty |r^{p+1}\bar{D}\psi^{s}|+|r^pD\psi^{s}|+|r^p\psi^{s}|\lesssim T^{(2p-2+\delp)/2+1}I_0^{1/2}.$$
The final estimate should be interpreted as saying that $|r^{(\delm+1)/2}\bar{D}\psi^{s}|$, $|r^{(\delm-1)/2}D\psi^{s}|$, and $|r^{(\delm-1)/2}\psi^{s}|$ decay like $T^{(\delm-3+\delp)/2}$ in a weak sense.
\end{theorem}

The remainder of this section is devoted to proving Theorem \ref{s:main_thm}.

\subsection{Bootstrap assumptions}\label{s:bootstrap_assumptions_sec}

We begin the proof of Theorem \ref{s:main_thm} by making the following bootstrap assumptions.
\begin{align*}
E^{12}(t)&\le C_bI_0, \\
E^{12}_{p\in[\delm,2-\delp]}(t)&\le C_b I_0, \\
E^{10}_{p\in[\delm,2-\delp]}(t)&\le C_b T^{\delm-2+\delp}I_0, \\
\int_t^\infty E^{9}_{p\in[\delm-1,\delm]}(\tau)d\tau &\le C_b T^{(\delm-1)-2+\delp+1}I_0.
\end{align*}
Note that, with the exception of the highest order energies, these bootstrap assumptions are consistent with the general principle that $E_p^{s}(t)\sim T^{p-2+\delp}$, which the reader should keep in mind throughout the proof of the main theorem.

\subsection{Improved $L^\infty$ estimates}

The $L^\infty$ estimates from Proposition \ref{s:infinity_prop} are essential to the argument of the proof of the main theorem. But for better clarity, we first remove the nonlinear (error) terms from these estimates and summarize them in the following lemma.
\begin{lemma}\label{s:simplified_pointwise_lemma}
In the context of the bootstrap assumptions provided in \S\ref{s:bootstrap_assumptions_sec}, the following $L^\infty$ estimates hold for $s\le 12$ and all $p$ in any bounded range.
\begin{align}
|r^{p+1}\bar{D}\psi^{s}|^2+|r^pD\psi^{s}|^2+|r^p\psi^{s}|^2 &\lesssim  E_{2p}^{s+3}(t), \\
|rD\psi^{s}|^2 &\lesssim E^{s+3}(t). \label{s:classic_energy_pointwise_estimate}
\end{align}
These are the same as the estimates from Proposition \ref{s:infinity_prop}, except that the nonlinear (error) terms have been removed.
\end{lemma}

\begin{proof} 
See the proof of Lemma \ref{m:improved_infty_lem}, keeping in mind Proposition \ref{s:infinity_prop}.
\end{proof}
The conclusion of this lemma is the same as the statement of part (II) of the main theorem.

\begin{remark}
Lemma \ref{s:simplified_pointwise_lemma}, which is a simplified version of Proposition \ref{s:infinity_prop} in the sense that there are no nonlinear (error) terms on the right side of any of the estimates in Lemma \ref{s:simplified_pointwise_lemma}, will be used in the remainder of the proof of the main theorem as a replacement for Proposition \ref{s:infinity_prop}.
\end{remark}

\subsection{Refined estimates for $N^{s}(t)$ and $N_p^{s}(t)$ (key step)}

The $L^\infty$ estimates given in Lemma \ref{s:simplified_pointwise_lemma} allow us to provide refined estimates for the nonlinear error terms. \textbf{This is the crucial step of the proof.}

\begin{lemma}\label{s:refined_nl_E_lem}
In the context of the bootstrap assumptions provided in \S\ref{s:bootstrap_assumptions_sec}, if $s\le 12$, then
$$N^{s}(t)\lesssim (E^{s}(t))^{1/2}\left((E^{s}(t))^{1/2}(E_{\delm-1}^{s/2+3}(t))^{1/2}+(E^{s}_{1-\delm}(t))^{1/2}(E^{s/2+3}_{\delm-1}(t))^{1/2}\right).$$
\end{lemma}
\begin{proof}
We recall the definition of $N^{s}(t)$ from Corollary \ref{s:s_cor}. For all $s$,
$$N^{s}(t)= (E^{s}(t))^{1/2}\sum_{s'\le s}||\Gamma^{s'}\Box_g\psi||_{L^2(\Sigma_t)}.$$
Therefore, it suffices to prove the following estimate.
$$
||\Gamma^s\Box_g\psi||_{L^2(\Sigma_t)} \\
\lesssim (E^{s}(t))^{1/2}(E_{\delm-1}^{s/2+3}(t))^{1/2}+(E^{s}_{1-\delm}(t))^{1/2}(E^{s/2+3}_{\delm-1}(t))^{1/2}.
$$
Note that by Proposition \ref{s:commutator_prop},
\begin{align*}
||\Gamma^s\Box_g\psi||_{L^2(\Sigma_t)} &\lesssim ||\gamma^{hi}\beta^{lo}||_{L^2(\Sigma_t)}+||\beta^{hi}\gamma^{lo}||_{L^2(\Sigma_t)} \\
&\lesssim ||r^{\frac{1-\delm}2}\gamma^{hi}||_{L^2(\Sigma_t)}||r^{\frac{\delm-1}2}\beta^{lo}||_{L^\infty(\Sigma_t)}+||\beta^{hi}||_{L^2(\Sigma_t)}||\gamma^{lo}||_{L^\infty(\Sigma_t)} \\
&\lesssim (E^s_{1-\delm}(t))^{1/2}(E^{s/2+3}_{\delm-1}(t))^{1/2}+(E^s(t))^{1/2}(E_{\delm-1}^{s/2+3}(t))^{1/2}.
\end{align*}
\end{proof}

\begin{lemma}\label{s:refined_nl_Ep_lem}
In the context of the bootstrap assumptions provided in \S\ref{s:bootstrap_assumptions_sec}, if $s\le 12$ and $C_bI_0\le 1$, then
$$
N_p^{s}(t)\lesssim E^{s}(t)B_p^{s/2+3}(t)+B_p^{s}(t)E^{s/2+3}(t) +E^{s}_{p'}(t)(E_{p''}^{s/2+3}(t))^{1/2}.
$$
for arbitrary $p'$ and $p''<2-\delp$.
\end{lemma}

\begin{proof}
From Corollary \ref{s:s_cor},
$$N_p^{s}(t)=\sum_{s'\le s} \int_{\Sigma_t}r^{p+1}(\Gamma^{s'}\Box_g\psi)^2+\sum_{s'\le s}\int_{\Sigma_t\cap\{r\approx r_{trap}\}}|\pd_t\psi^{s'}\Gamma^{s'}\Box_g\psi|,$$
We rewrite for some fixed radius $R>r_{trap}$,
\begin{multline*}
N_p^{s}(t) = \left[\sum_{s'\le s} \int_{\Sigma_t\cap\{r<R\}}(\Gamma^{s'}\Box_g\psi)^2+\sum_{s'\le s}\int_{\Sigma_t\cap\{r\approx r_{trap}\}}|\pd_t\psi^{s'}\Gamma^{s'}\Box_g\psi|\right] \\
+\sum_{s'\le s} \int_{\Sigma_t\cap\{r>R\}}r^{p+1}(\Gamma^{s'}\Box_g\psi)^2.
\end{multline*}
Using Proposition \ref{s:commutator_prop} again, the terms in square brackets can be bounded by $E^{s}_{p'}(t)E_{p''}^{s/2+3}(t)$ and $E^{s}_{p'}(t)(E_{p''}^{s/2+3}(t))^{1/2}$ respectively, with $p'$ and $p''$ arbitrary because the terms are supported on a compact radial interval. Moreover, by the second bootstrap assumption, since $E_{p''}^{s/2+3}(t)\le C_bI_0\le 1$, we have\footnote{In particular, if there were no trapping, we would have additional decay, since $p''$ could be taken less than $2-\delp$.} $E^{s}_{p'}(t)E_{p''}^{s/2+3}(t)\le E^{s}_{p'}(t)(E_{p''}^{s/2+3}(t))^{1/2}$.

 More importantly,
\begin{align*}
\sum_{s'\le s} \int_{\Sigma_t\cap\{r>R\}}r^{p+1}(\Gamma^{s'}\Box_g\psi)^2 &\lesssim \int_{\Sigma_t\cap\{r>R\}}r^{p+1}\left[(\gamma^{hi})^2(\beta^{lo})^2+(\beta^{hi})^2(\gamma^{lo})^2\right] \\
&\lesssim \int_{\Sigma_t\cap\{r>R\}}r^{p-1}(\gamma^{hi})^2(r\beta^{lo})^2+(\beta^{hi})^2(r^{(p+1)/2}\gamma^{lo})^2 \\
&\lesssim \int_{\Sigma_t\cap\{r>R\}}r^{p-1}(\gamma^{hi})^2||r\beta^{lo}||^2_{L^\infty(\Sigma_t)}+(\beta^{hi})^2||r^{(p+1)/2}\gamma^{lo}||^2_{L^\infty(\Sigma_t)} \\
&\lesssim B_p^s(t)||r\beta^{lo}||^2_{L^\infty(\Sigma_t)}+E^s(t)||r^{(p+1)/2}\gamma^{lo}||^2_{L^\infty(\Sigma_t)} \\
&\lesssim B_p^{s}(t)E^{s/2+3}(t)+E^{s}(t)B_p^{s/2+3}(t).
\end{align*}
\end{proof}

\begin{corollary}\label{s:NL_absorb_bulk}
In the context of the bootstrap assumptions provided in \S\ref{s:bootstrap_assumptions_sec}, Corollary \ref{s:s_cor} and Lemma \ref{s:refined_nl_Ep_lem} imply that if $s\le 12$ and $C_bI_0$ is sufficiently small, then
$$
E_p^{s}(t_2)+\int_{t_1}^{t_2}B_p^{s}(t)dt
\lesssim E_p^{s}(t_1)+(C_bI_0)^{1/2}\int_{t_1}^{t_2}E_p^{s}(t)T^{(\delm-3+\delp)/2}dt.
$$
whenever $s/2+3\le s$ (which means $6\le s$).
\end{corollary}
\begin{proof}
By Corollary \ref{s:s_cor} and Lemma \ref{s:refined_nl_Ep_lem},
\begin{multline*}
E_p^{s}(t_2)+\int_{t_1}^{t_2}B_p^{s}(t)dt \\
\lesssim E_p^{s}(t_1)+\int_{t_1}^{t_2}\left[E^{s}(t)B_p^{s/2+3}(t)+B_p^{s}(t)E^{s/2+3}(t)+E^{s}_{p'}(t)(E_{p''}^{s/2+3}(t))^{1/2}\right]dt
\end{multline*}
Since $s/2+3\le s$,
$$E_p^{s}(t_2)+\int_{t_1}^{t_2}B_p^{s}(t)dt \lesssim E_p^{s}(t_1)+\int_{t_1}^{t_2}\left[E^{s}(t)B_p^{s}(t)+E^{s}_{p'}(t)(E_{p''}^{s/2+3}(t))^{1/2}\right]dt.$$
Recall that the bootstrap assumptions imply
$$E^s(t)\lesssim C_bI_0$$
 and 
$$\int_{t_1}^\infty E_{\delm-1}^9(t)dt\lesssim C_bT^{\delm-3+\delp+1}I_0.$$
Using these two estimates and the weak decay principle with the second estimate, we have
$$E_p^{s}(t_2)+\int_{t_1}^{t_2}B_p^{s}(t)dt \lesssim E_p^{s}(t_1)+\int_{t_1}^{t_2}\left[C_bI_0B_p^{s}(t)+E^{s}_{p}(t)(C_bT^{\delm-3+\delp}I_0)^{1/2}\right]dt.$$
If $C_bI_0$ is sufficiently small, the $\int_{t_1}^{t_2}C_bI_0B_p^s(t)dt$ term can be absorbed into the bulk term on the left side. Therefore,
$$
E_p^{s}(t_2)+\int_{t_1}^{t_2}B_p^{s}(t)dt
\lesssim E_p^{s}(t_1)+(C_bI_0)^{1/2}\int_{t_1}^{t_2}E_p^{s}(t)T^{(\delm-3+\delp)/2}dt.
$$
\end{proof}

\subsection{Recovering boundedness of $E^{s}(t)$ ($s=12$)}

Let $s=12$. By Corollary \ref{s:s_cor} and Lemma \ref{s:refined_nl_E_lem},
\begin{align*}
E^{s}(t) \lesssim & E^{s}(0)+\int_0^tN^{s}(\tau)d\tau \\
\lesssim & I_0 +\int_0^t(E^{s}(\tau))^{1/2}\left((E^{s}(\tau))^{1/2}(E_{\delm-1}^{9}(\tau))^{1/2}\right. \\
&\hspace{2.5in}\left.+(E^{s}_{1-\delm}(\tau))^{1/2}(E^{9}_{\delm-1}(\tau))^{1/2}\right)d\tau \\
\lesssim & I_0+\int_0^t(C_bI_0)^{1/2}(C_bI_0)^{1/2}(C_bT^{\delm-3+\delp}I_0)^{1/2}d\tau \\
\lesssim & (1+C_b^{3/2}I_0^{1/2})I_0.
\end{align*}
In particular, we used the weak decay principle in the third step. It follows that if $C_b^{3/2}I_0^{1/2}$ is sufficiently small, then
$$E^{s}(t)\lesssim I_0.$$
This recovers the first bootstrap assumption.

\subsection{Recovering boundedness of $E^{s}_p(t)$ ($s=12$)}\label{s:recover_Ep_sec}

Our first application of Corollary \ref{s:NL_absorb_bulk} is to prove boundedness of $E^{s}_p(t)$. Let $s=12$. Then
\begin{align*}
E_p^{s}(t)+\int_0^tB_p^{s}(\tau)d\tau &\lesssim E_p^{s}(t)+(C_bI_0)^{1/2}\int_0^tE_p^{s}(\tau)T^{(\delm-3+\delp)/2}d\tau \\
&\lesssim I_0+(C_bI_0)^{3/2}T^{\delm-1/2+\delp} \\
&\lesssim (1+C_b^{3/2}I_0^{1/2})I_0.
\end{align*}
It follows that if $C_b^{3/2}I_0^{1/2}$ is sufficiently small, then
$$E_{p\in[\delm,2-\delp]}^{s}(t)\lesssim I_0.$$
This recovers the second bootstrap assumption.

\subsection{Proving decay for $E^{s}_{p\in[1-\delp,2-\delp]}(t)$ ($s=11$) and $E^{s}_{p\in[\delm,2-\delp]}(t)$ ($s=10$)}

At last we are ready to prove decay. This will be accomplished by repeatedly applying the following lemma.

\begin{lemma}\label{s:NL_WE_decay}
Suppose  $p+1,p\in [\delm,2-\delp]$ and
$$E_{p+1}^{s+1}(t)\lesssim T^{(p+1)-2+\delp}I_0,$$
$$E_p^{s}(t)\le C_bT^{p-2+\delp}I_0.$$
Then if $I_0$ is sufficiently small,
$$E_p^{s}(t)\lesssim T^{p-2+\delp}I_0.$$
\end{lemma}

\begin{proof}
Using the mean value theorem, for a given $t$, let $t'\in [t/2,t]$ be the value for which $B_{p+1}^{s+1}(t')=\frac2t\int_{t/2}^tB_{p+1}^{s+1}(\tau)d\tau$. Then using Corollary \ref{s:NL_absorb_bulk},
\begin{align*}
E_p^{s}(t) &\lesssim E_p^{s}(t')+(C_bI_0)^{1/2}\int_{t'}^tE_p^{s}(\tau)T^{(\delm-3+\delp)/2}d\tau \\
&\lesssim E_p^{s}(t')+(C_bI_0)^{1/2}\int_{t'}^t(C_bT^{p-2+\delp}I_0)T^{(\delm-3+\delp)/2}d\tau \\
&\lesssim E_p^{s}(t')+C_b^{3/2}I_0^{1/2}T^{p-2+\delp}I_0.
\end{align*}
Now by the choice of $t'$,
\begin{multline*}
E_p^{s}(t')\lesssim B_{p+1}^{s+1}(t')=\frac2t\int_{t/2}^tB_{p+1}^{s+1}(\tau)d\tau \\
\lesssim t^{-1}E_{p+1}^{s+1}(t/2)+t^{-1}(C_bI_0)^{1/2}\int_{t/2}^tE_{p+1}^{s+1}(\tau)T^{(\delm-3+\delp)/2}d\tau \\
\lesssim T^{p-2+\delp}I_0+(C_bI_0)^{1/2}T^{p-2+\delp}_0
\end{multline*}
Thus,
$$E_p^{s}(t)\lesssim (1 +C_b^{1/2}I_0^{1/2}+C_b^{3/2}I_0^{1/2})T^{p-2+\delp}I_0.$$
Provided $C_b^{1/2}I_0^{1/2}+C_b^{3/2}I_0^{1/2}$ is sufficiently small, the proof is complete.
\end{proof}

By applying Lemma \ref{s:NL_WE_decay} for $p=1-\delp$ and $s=11$, we obtain
$$E_{1-\delp}^{11}(t)\lesssim T^{-1}I_0.$$
By interpolation we obtain decay for all $p\in [1-\delp,2-\delp]$.
$$E_{p\in [1-\delp,2-\delp]}^{11}(t)\lesssim T^{p-2+\delp}I_0.$$
Then applying Lemma \ref{s:NL_WE_decay} again for each $p\in[\delm,1-\delp]$ and $s=10$, we obtain
\begin{equation*}
E_{p\in [\delm,2-\delp]}^{10}(t)\lesssim T^{p-2+\delp}I_0.
\end{equation*}
This recovers the third bootstrap assumption.

\subsection{Recovering weak decay for $E_{p\in[\delm-1,\delm]}^{s}(t)$ ($s=9$)}

Finally, set $s=9$ and observe that for $p\in [\delm-1,\delm]$,
$$\int_t^\infty E_p^{s}(\tau)d\tau \lesssim \int_t^\infty B_{p+1}^{s+1}(\tau)d\tau\lesssim T^{(p+1)-2+\delp}I_0=T^{p-2+\delp+1}I_0.$$
In particular, by taking $p\in[\delm-1,\delm]$, this recovers the final bootstrap assumption.

\noindent This completes the proof of Theorem \ref{s:main_thm}. \qed

\chapter{Axisymmetric waves in subextremal Kerr spacetimes}\label{kerr_chap}

The third problem of this thesis investigates axisymmetric semilinear waves on the exterior of a subextremal Kerr black hole. The Kerr black hole family includes the Schwarzschild black hole as well as rotating black holes.

The metric for a Kerr black hole with mass $M$ and angular momentum $a$ is most well-known in Boyer-Lindquist coordinates.
\begin{multline*}
g_{\mu\nu}dx^\mu dx^\nu = -\left(1-\frac{2Mr}{r^2+a^2\cos^2\theta}\right)dt^2-2\frac{2Mr a\sin^2\theta}{r^2+a^2\cos^2\theta}dtd\phi \\
+\left(r^2+a^2+\frac{2Mra^2\sin^2\theta}{r^2+a^2\cos^2\theta}\right)\sin^2\theta d\phi^2+\frac{r^2+a^2\cos^2\theta}{r^2-2Mr+a^2}dr^2+(r^2+a^2\cos^2\theta)d\theta^2.
\end{multline*}
Its volume form is
$$\mu = (r^2+a^2\sin^2\theta)\sin\theta.$$
Despite the fact that the Kerr metric is significantly more complicated than the Schwarzschild metric, there are a number of similarities between the problem of this chapter and the problem of the previous chapter, which we now discuss.

\textbf{Both (linearized) problems have Lagrangians with similar structure.} The linearized Lagrangian density
$$\mathcal{L}=\mu g^{\alpha\beta}\pd_\alpha\psi\pd_\beta\psi$$
for the second problem was
$$\mathcal{L}_{Schwarzschild}=r^2\sin\theta\left[-\left(1-\frac{2M}r\right)^{-1}(\pd_t\psi)^2+\left(1-\frac{2M}r\right)(\pd_r\psi)^2+|\sla\nabla\psi|^2\right],$$
while the linearized Lagrangian density for this problem is
\begin{multline*}
\mathcal{L}_{Kerr, axisymmetry} \\
=q^2\sin\theta\left[-\left(\frac{(r^2+a^2)^2}{q^2\Delta}-a^2\sin^2\theta\right)(\pd_t\psi)^2+\frac{\Delta}{q^2}(\pd_r\psi)^2+\frac1{q^2}(\pd_\theta\psi)^2\right],
\end{multline*}
where
$$\Delta = r^2-2Mr+a^2,$$
$$q^2 = r^2+a^2\cos^2\theta.$$

At first glance, $\mathcal{L}_{Kerr,axisymmetry}$ appears to be rather complicated compared to $\mathcal{L}_{Schwarzschild}$. Of particular concern are the dependence on $\theta$ in the quantity $q^2$ and the fact that $g^{tt}g^{rr}\ne -1$. However, $\mathcal{L}_{Kerr,axisymmetry}$ can be rewritten as
\begin{multline*}
\mathcal{L}_{Kerr,axisymmetry} \\
=(r^2+a^2)\sin\theta\left[-\frac{r^2+a^2}{\Delta}(\pd_t\psi)^2+\frac{\Delta}{r^2+a^2}(\pd_r\psi)^2+\frac{1}{r^2+a^2}Q^{\alpha\beta}\pd_\alpha\psi\pd_\beta\psi\right],
\end{multline*}
where
$$Q^{\alpha\beta}\pd_\alpha\psi\pd_\beta\psi=(\pd_\theta\psi)^2+a^2\sin^2\theta(\pd_t\psi)^2.$$
In this form, the coefficients for $(\pd_t\psi)^2$ and $(\pd_r\psi)^2$ have a product of $-1$ and all $\theta$ dependence occurs either in the modified volume form or the tensor $Q^{\alpha\beta}$. This motivates a generalization of the quantity $\alpha$ from the last chapter.
$$\alpha =\frac{\Delta}{r^2+a^2}.$$

The tensor $Q^{\alpha\beta}$, known as \textit{Carter's hidden symmetry}, is a generalization of the angular metric $r^2\sla{g}^{\alpha\beta}$, and it corresponds to a second order operator $Q$ that generalizes the operator $r^2\sla\triangle$ and commutes with $\Box_{g_{Kerr}}$. 

\textbf{Both problems feature simple trapping phenomena.} As was discussed in \S\ref{s:geodesics_sec}, trapped null geodesics complicate the analysis of linear waves. In particular, the classical vectorfield method seems to only work for problems where trapping is confined to a single \textit{photon sphere}, as is the case for all problems in this thesis. (The recent results \cite{AnBlu} and \cite{DRSR}, which go further, are discussed below.)


Unlike the Schwarzschild spacetime, which has a single radius $r_{trap}=3M$ where all trapped geodesics occur, any Kerr $a\ne 0$ spacetime has an entire radial interval where trapped geodesics can occur. However, for a fixed geodesic angular momentum, there is again only a single trapping radius for null geodesics with that angular momentum. One can think of the condition $\pd_\phi\psi=0$ as a restriction to the ``waves with zero angular momentum.'' The trapping radius $r_{trap}$ in this problem coincides with the trapping radius for geodesics with zero angular momentum. It is the radius that maximizes the geodesic potential $\frac{\Delta}{(r^2+a^2)^2}$, which coincides with the familiar potential $r^{-2}(1-\frac{2M}r)$ when $a=0$.

Some recent works handle linear waves with multiple trapping radii. In \cite{AnBlu}, it is shown that higher order symmetry operators can be used to generalize the vectorfield method in a way that is useful for slowly rotating Kerr spacetimes without the requirement of axisymmetry. In \cite{DRSR}, the linear wave equation (without the requirement of axisymmetry) in all subextremal Kerr spacetimes is completely understood by using Fourier methods. Due to the issue of trapping as well as the problem of the ergoregion (discussed next), it appears likely that Fourier methods are necessary for treating the general wave equation in Kerr.

\textbf{Both problems have a positive energy that is conserved at the linear level.} Another new challenge in Kerr is the presence of an \textit{ergoregion}, a region near the event horizon where the killing vectorfield $\pd_t$ becomes spacelike. In the analysis of waves, this manifests itself in the fact that there is no conserved energy that is globally positive. The most obvious candidate, based on the vectorfield multiplier $\pd_t$ in Boyer-Lindquist coordinates, has terms with $\pd_\phi\psi$ that have negative coefficients in the ergoregion. By assuming $\pd_\phi\psi=0$, the energy is still positive.

Given the above similarities between the problem of this chapter and the problem of the previous chapter, we now discuss a few differences.

\textbf{Unlike the problem in the previous chapter, this problem has weaker symmetries.} In the previous chapter, the first order symmetry operators $\pd_t$ and $\Omega_i$ were used as commutators. However, the Kerr spacetimes only have one first order angular symmetry operator $\Omega_z=\pd_\phi$, which is not useful when $\pd_\phi\psi=0$. Instead, the symmetries used in Kerr are the Killing vectorfield $\pd_t$ and the second order hidden symmetry operator $Q$ introduced earlier. The operator $Q$ differs from the spherical Laplacian by a term depending on $\pd_t^2$, so $Q$ and $\pd_t^2$ together provide the spherical Laplacian, which has the elliptic properties necessary to control the angular derivatives. See \S\ref{k:symmetry_operators_sec} and Lemma \ref{main_elliptic_lemma}.

\textbf{Compared to the previous chapter, a weaker assumption is made for the nonlinearity.} To avoid a distracting complication in Chapter \ref{szd_chap}, an extra assumption was made so that nonlinear terms with $\pd_r\psi$ vanished on the event horizon. The reason was to avoid commuting with the redshift vectorfield, which is not actually a symmetry operator. Since the method of Chapter \ref{szd_chap} is now assumed to be well understood, the redshift vectorfield commutator will be used in this chapter to eliminate the need for the additional assumption on the nonlinearity. See \S\ref{k:additional_commutator_sec}. This is in preparation for the final two problems of this thesis.

\section{The energy estimate and the $h\pd_t$ estimate}\label{k:energy_sec}
The most fundamental of all estimates is the energy estimate. We begin by proving this estimate and then prove the slightly more general $h\pd_t$ estimate.

\subsection{The energy estimate}

Since the vectorfield $\pd_t$ is a Killing vectorfield, we have the following identity.
\begin{lemma}\label{k:ee_identity_lem}(Energy identity for Kerr)
\begin{equation*}
\int_{H_{t_1}^{t_2}}J^r[\pd_t] +\int_{\Sigma_{t_2}}-J^t[\pd_t] = \int_{\Sigma_{t_1}}-J^t[\pd_t]+\int_{t_1}^{t_2}\int_{\Sigma_t}-2\pd_t\psi\Box_g\psi
\end{equation*}
In particular, if $\Box_g\psi=0$, then the quantity
$$E(t)=\int_{H_{t_0}^t}J^r[\pd_t]+\int_{\Sigma_t}-J^t[\pd_t]$$
is conserved.
\end{lemma}
\begin{proof}
The proof is identical to the proof of Lemma \ref{s:ee_identity_lem}.
\end{proof}

Now, we calculate the flux terms in the previous lemma.

\begin{lemma}\label{k:ee_bndry_lem}
On the event horizon $H_{t_1}^{t_2}$,
$$J^r[\pd_t]\approx (\pd_t\psi)^2$$
and on a constant-time hypersurface $\Sigma_t$,
$$-J^t[\pd_t]\approx \chi_H(\pd_r\psi)^2+(\pd_t\psi)^2+\frac{Q^{\alpha\beta}}{q^2}\pd_\alpha\psi\pd_\beta\psi,$$
where $\chi_H=1-\frac{r_H}r$.
\end{lemma}
\begin{proof}
The proof is similar to the proof of Lemma \ref{s:ee_bndry_lem}.
\end{proof}

From the previous two lemmas, we conclude the following energy estimate.
\begin{proposition}\label{k:classic_ee_prop}(Energy estimate for Kerr, assuming axisymmetry)
Suppose $\pd_\phi\psi=0$. Then
\begin{multline*}
\int_{H_{t_1}^{t_2}}(\pd_t\psi)^2+\int_{\Sigma_{t_2}}\left[\chi_H(\pd_r\psi)^2+(\pd_t\psi)^2+\frac{Q^{\alpha\beta}}{q^2}\pd_\alpha\psi\pd_\beta\psi\right] \\
 \lesssim \int_{\Sigma_{t_1}}\left[\chi_H(\pd_r\psi)^2+(\pd_t\psi)^2+\frac{Q^{\alpha\beta}}{q^2}\pd_\alpha\psi\pd_\beta\psi\right] + Err_\Box,
\end{multline*}
where $\chi_H=1-\frac{r_H}r$ and
$$Err_\Box=\int_{t_1}^{t_2}\int_{\Sigma_t}|\pd_t\psi\Box_g\psi|.$$
\end{proposition}
\begin{proof}
This follows directly from Lemmas \ref{k:ee_identity_lem} and \ref{k:ee_bndry_lem}.
\end{proof}

\subsection{The $h\pd_t$ energy estimate}\label{k:hdt_sec}

We proceed as in \S\ref{m:hdt_sec} and \S\ref{s:hdt_sec} to prove the $h\pd_t$ estimate.
\begin{lemma}\label{k:divJhdt_lem}
If $h=h(r)$ is constant in the interval $r\in[r_H,r_H+\delh]$, and $\Box_g\psi=0$, then
$$\frac{q^2}{r^2+a^2}divJ[h\pd_t]=\frac{h'}{2}\left[(L\psi)^2-(\lbar\psi)^2\right].$$
\end{lemma}
\begin{proof}
Recall from Lemma \ref{divJ_lem} that
$$divJ[X]=K^{\mu\nu}\pd_\mu\psi\pd_\nu\psi,$$
where
$$K^{\mu\nu}=g^{\mu\lambda}\pd_\lambda X^\nu+g^{\nu\lambda}\pd_\lambda X^\mu-X^\lambda\pd_\lambda(g^{\mu\nu})-divXg^{\mu\nu}.$$
Since $\pd_t$ is killing, $\pd_t(g^{\mu\nu})=0$ and $div(\pd_t)=0$. Also, the only component of $X$ is the $t$ component and the only nonzero derivative of that component is the $\pd_r$ derivative, so
$$g^{\mu\lambda}\pd_\lambda X^\nu+g^{\nu\lambda}\pd_\lambda X^\mu=g^{\mu r}\pd_r X^\nu+g^{\nu r}\pd_r X^\mu.$$
It follows that the only possible nonzero $K^{\mu\nu}$ components are
\begin{align*}
K^{tt} &= 2g^{tr}h' \\
K^{tr}+K^{rt} &= 2g^{rr}h'.
\end{align*}
Since $h'=0$ in the region $r\in [r_H,r_H+\delh]$, the first of these actually vanishes. Recalling that $\alpha=\frac{q^2}{r^2+a^2}g^{rr}$, we conclude that
\begin{align*}
\frac{q^2}{r^2+a^2}divJ[h\pd_t]&=2\alpha h'\pd_r\psi\pd_t\psi \\
&=\frac{h'}2\left[(\alpha\pd_r\psi+\pd_t\psi)^2-(\alpha\pd_r\psi-\pd_t\psi)^2\right] \\
&=\frac{h'}2\left[(L\psi)^2-(\lbar\psi)^2\right].
\end{align*}
\end{proof}

Taking $h$ to be a positive function decreasing to zero as $r\rightarrow\infty$ at a particular rate, we obtain the $h\pd_t$ estimate.
\begin{proposition}\label{k:hdt_prop}($h\pd_t$ estimate for Kerr, assuming axisymmetry)
Suppose $\pd_\phi\psi=0$.
Let $R>r_H+\delh$ be any given radius. Then for all $\epsilon>0$ and $p<2$, there is a small constant $c_\epsilon$ and a large constant $C_\epsilon$, such that
\begin{multline*}
\int_{H_{t_1}^{t_2}}(\pd_t\psi)^2+\int_{\Sigma_{t_2}}r^{p-2}\left[\chi_H(\pd_r\psi)^2+(\pd_t\psi)^2+\frac{Q^{\alpha\beta}}{q^2}\pd_\alpha\psi\pd_\beta\psi\right] \\
+\int_{t_1}^{t_2}\int_{\Sigma_t\cap\{R+M<r\}}c_\epsilon r^{p-3}(\lbar\psi)^2 \\
\lesssim \int_{\Sigma_{t_1}}r^{p-2}\left[\chi_H(\pd_r\psi)^2+(\pd_t\psi)^2+\frac{Q^{\alpha\beta}}{q^2}\pd_\alpha\psi\pd_\beta\psi\right]+Err,
\end{multline*}
where $\chi_H=1-\frac{r_H}r$ and
\begin{align*}
Err&=Err_1+Err_\Box \\
Err_1&=\int_{t_1}^{t_2}\int_{\Sigma_t\cap\{R<r\}}\epsilon r^{-1}(L\psi)^2 \\
Err_\Box&=\int_{t_1}^{t_2}\int_{\Sigma_t}C_\epsilon r^{p-2}|\pd_t\psi\Box_g\psi|.
\end{align*}
\end{proposition}
\begin{proof}
See the proof of Proposition \ref{m:hdt_prop}, noting Lemmas \ref{k:divJhdt_lem} and \ref{k:ee_bndry_lem}.
\end{proof}

\section{The Morawetz estimate}\label{k:morawetz_sec}

In this section, we prove the Morawetz estiamte for Kerr, assuming axisymmetry. It is a direct generalization of the method used to prove the Morawetz estimate for Schwarzschild. As such, we will proceed at a faster pace, dropping some calculations that are essentially the same for both spacetimes.

\subsection{The partial Morawetz estimate}

We now proceed as in \S\ref{s:morawetz_sec} to construct the Morawetz estimate, beginning with the partial Morawetz estimate as in \S\ref{s:partial_morawetz_sec}.

\begin{proposition}\label{k:partial_morawetz_prop}(Partial Morawetz estimate)
Suppose $\pd_\phi\psi=0$. Then there exist a vectorfield $X_0$ and a function $w_0$ such that the current $J[X_0,w_0]$ as defined in \S\ref{current_template_sec} satisfies
$$\left[\frac{M^2}{r^3}\left(1-\frac{r_H}r\right)^2(\pd_r\psi)^2+\frac1r\left(1-\frac{r_{trap}}r\right)^2\frac{Q^{\alpha\beta}}{q^2}\pd_\alpha\psi\pd_\beta\psi+\frac{M}{r^4}1_{r>r_*}\psi^2\right]\lesssim divJ[X_0,w_0]$$
in Boyer-Lindquist coordinates, where $r_*$ is the larger root of the polynomial $r^2-4Mr+a^2$.
\end{proposition}

\begin{proof}
For the proof of this proposition, we drop the $0$ subscript from $X_0$ and $w_0$. We begin by restating Lemma \ref{m:pme_initial_lem} for Kerr.
\begin{lemma}
With the choices $X=X^r(r)\pd_r$ and $w=w(r)$,
\begin{multline*}
divJ[X,w]= 
(wg^{tt}-q^{-2}\pd_r(q^2X^rg^{tt}))(\pd_t\psi)^2+(w-\pd_rX^r)\frac{Q^{\alpha\beta}}{q^2}\pd_\alpha\psi\pd_\beta\psi \\
+(wg^{rr}+\pd_rX^r g^{rr}-2rq^{-2}X^rg^{rr}-X^r\pd_rg^{rr})(\pd_r\psi)^2 -\frac12q^{-2}\pd_r(q^2g^{rr}\pd_rw)\psi^2.
\end{multline*}
\end{lemma}
\begin{proof}
See the proof of Lemma \ref{m:pme_initial_lem}.
\end{proof}

It is difficult to make sense of the coefficients in the above lemma, because most coefficients are sums of multiple terms. By making the choices $X^r(r)=u(r)v(r)$ and $w(r)=v(r)\pd_ru(r)$, each coefficient can be written as a single term.

\begin{lemma}
With the choices $X=u(r)v(r)\pd_r$ and $w=v(r)\pd_ru(r)$,
\begin{multline*}
divJ[X,w]= 
-\frac{u}{q^2}\pd_r\left(vq^2g^{tt}\right)(\pd_t\psi)^2-u\pd_rv\frac{Q^{\alpha\beta}}{q^2}\pd_\alpha\psi\pd_\beta\psi \\
+q^{-2}(q^2g^{rr})^2u^{-1}\pd_r\left(\frac{u^2v}{q^2g^{rr}}\right)(\pd_r\psi)^2-\frac12q^{-2}\pd_r(q^2g^{rr}\pd_rw)\psi^2.
\end{multline*}
\end{lemma}
\begin{proof}
This follows directly by substituting $X^r=uv$ and $w=v\pd_ru$ into the expression in the previous lemma and combining terms.
\end{proof}

As in the proofs of previous partial Morawetz estimates, we choose a relation between $X$ and $w$ so that the coefficient of $(\pd_t\psi)^2$ vanishes. In terms of the functions $u$ and $v$, it suffices to set $v=(q^2g^{tt})^{-1}=\frac{\Delta}{(r^2+a^2)^2}=\frac{r^2-2Mr+a^2}{(r^2+a^2)^2}$.

\begin{remark}
The function $v=\frac{r^2-2Mr+a^2}{(r^2+a^2)^2}$ is the Kerr analogue of the geodesic potential derived in \S\ref{s:geodesics_sec} for geodesics with $L_z=0$. (The condition $L_z=0$ is the geodesic analogue of the condition $\pd_\phi\psi=0$.) It has a maximum at the radius $r_{trap}$, which coincides with $3M$ when $a=0$. This explains why the trapping radius $r_{trap}$ will show up in the calculations that follow.
\end{remark}

\begin{lemma}
With the choice $v=\frac{r^2-2Mr+a^2}{(r^2+a^2)^2}$, the coefficient of the $(\pd_t\psi)^2$ term vanishes and furthermore,
\begin{multline*}
divJ[X,w] \\
=-u\pd_rv\frac{Q^{\alpha\beta}}{q^2}\pd_\alpha\psi\pd_\beta\psi+2\frac{\Delta^2}{q^2(r^2+a^2)}\pd_r\left(\frac{u}{r^2+a^2}\right)(\pd_r\psi)^2-\frac12q^{-2}\pd_r(\Delta\pd_rw)\psi^2.
\end{multline*}
Therefore, it is necessary to choose $u$ and $w$ so that the following conditions are satisfied.
\begin{eqnarray}
u\pd_rv\ge 0, \label{k:u_cond_1_eqn}\\
\pd_r\left(\frac{u}{r^2+a^2}\right)\ge 0, \label{k:u_cond_2_eqn}\\
\pd_r(\Delta\pd_rw)\le 0, \label{k:u_cond_3_eqn}
\end{eqnarray}
Furthermore, $w$ and $u$ must be related by the following constraint.
\begin{equation}
w=\frac{\Delta}{(r^2+a^2)^2}\pd_ru. \label{k:u_constr_eqn}
\end{equation}
\end{lemma}

\begin{lemma}\label{k:choice_for_u_and_w_lem}
It is possible to choose $u$ and $w$ such that all three conditions (\ref{k:u_cond_1_eqn}-\ref{k:u_cond_3_eqn}) and the constraint (\ref{k:u_constr_eqn}) are satisfied. One particular choice is given by the following approach. (See Figure \ref{k:w_fig}.)

i) Require that $u(r_{trap})=0$ and $\pd_ru=\frac{(r^2+a^2)^2}{\Delta}w$. This will specify $u$ completely in terms of $w$ and satisfy the constraint (\ref{k:u_constr_eqn}).

ii) Require that $w$ be positive. Then $u$ will be an increasing function and the condition (\ref{k:u_cond_1_eqn}) will be satisfied.

iii) Determine a way to understand $\pd_r\left(\frac{u}{r^2+a^2}\right)$ in terms of $w$. Here is one way. Let $\tilde{K}^{rr}=\frac{(r^2+a^2)^2}{2r}\pd_r\left(\frac{u}{r^2+a^2}\right)$. This particular quantity has the same sign as $\pd_r\left(\frac{u}{r^2+a^2}\right)$ and has the property that its derivative can be written as a function of $w$. That is, $\pd_r\tilde{K}^{rr}=(r^2+a^2)\pd_r\left(\frac{w}{2rv}\right)$, where again, $v=\frac{\Delta}{(r^2+a^2)^2}$. The condition (\ref{k:u_cond_2_eqn}) now reduces to showing that $\tilde{K}^{rr}\ge 0$.

iv) For sufficiently large $r$, impose the condition that $\pd_r\tilde{K}^{rr}=0$. In particular, this means choosing $w=2rv$ for sufficiently large $r$. The quantity $2rv$ has a maximum at $r_*$, the larger root of the polynomial $r^2-4Mr+a^2$. (Again, see Figure \ref{k:w_fig}.) The choice $w=2rv$ for $r\ge r_*$ will satisfy both the remaining conditions (\ref{k:u_cond_2_eqn}-\ref{k:u_cond_3_eqn}) for $r>r_*$.

v) Observe that since $\pd_r(2rv)=0$ at $r_*$, the quantity $\Delta\pd_rw=0$ will necessarily vanish at $r=r_H$ and $r=_*$ if $w=2rv$ for $r\ge r_*$. By the mean value theorem, condition (\ref{k:u_cond_3_eqn}) necessitates that $\pd_rw=0$ entirely for $r_H\le r\le r_*$. Thus, take $w=w(r_*)$ for $r\le r_*$.

iv) Since $w$ has now been chosen for all $r$, check that $\tilde{K}^{rr}\ge 0$. In particular, with the given choice of $w$, $\tilde{K}^{rr}$ will be decreasing for $r<r_*$ and then remain constant for $r\ge r_*$. Thus, it suffices to compute that $\tilde{K}^{rr}(r_*)>0$.
\end{lemma}

\begin{figure}
\centering
\includegraphics[scale=0.7]{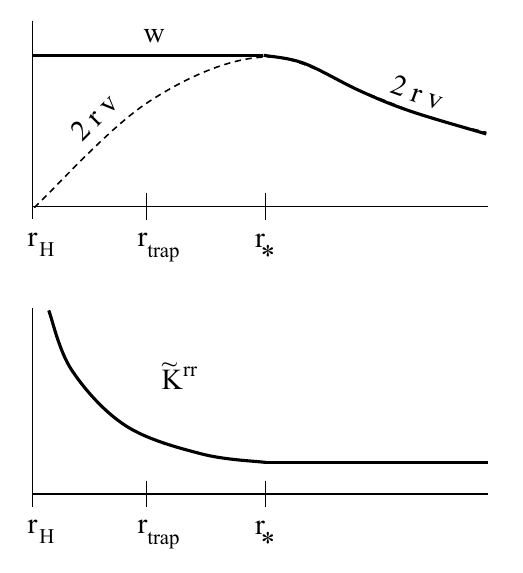}
\caption{Plots of $w$ and $\tilde{K}^{rr}=\frac12r^3\pd_r\left(\frac{u}{r^2+a^2}\right)$ for $r\ge r_H$. In the region $r>r_*$, the function $w$ is chosen so that $\tilde{K}^{rr}$ is constant. In the region $r<r_*$, $w$ is chosen to be constant (this is required by condition (\ref{k:u_cond_3_eqn})). The value $r_*$ is chosen so that $w$ remains $C^1$.}\label{k:w_fig}
\end{figure}

\begin{proof}
The lemma consists of multiple statements. Statements (i) and (ii) require no further justification. We turn to statement (iii). Note that
\begin{align*}
\pd_r\bar{K}^{rr}&=\pd_r\left(\frac{(r^2+a^2)^2}{2r}\pd_r\left(\frac{u}{r^2+a^2}\right)\right) \\
&=\pd_r\left(\frac{r^2+a^2}{2r}u'-u\right) = \frac{r^2+a^2}{2r}\pd_r^2u+\frac12\frac{r^2-a^2}{r^2}u'-u' \\
&=\frac12 (r^2+a^2)(r^{-1}\pd_r^2u-r^{-2}u') \\
&= (r^2+a^2)\pd_r\left(\frac{u'}{2r}\right).
\end{align*}
Recall the constraint (\ref{k:u_constr_eqn}) between $w$ and $u$, which stipulates $w=vu'$. Therefore, we substitute $u'=\frac{w}{v}$ and conclude
$$\pd_r\tilde{K}^{rr} = (r^2+a^2)\pd_r\left(\frac{w}{2rv}\right).$$
This verifies statement (iii). For statement (iv), we must find the value $r_*$ that maximizes the quantity $2rv$. A direct calculation reveals
\begin{align*}
\pd_r(2rv) &=2\pd_r\left(\frac{r(r^2-2Mr+a^2)}{(r^2+a^2)^2}\right) \\
&= -\frac{2}{r^2+a^2}\frac{r^2-a^2}{r^2+a^2}\frac{r^2-4Mr+a^2}{r^2+a^2}.
\end{align*}
Setting $\pd_r(2rv)=0$, we deduce that $r_*$ is the larger root of the polynomial $r^2-4Mr+a^2$. This indeed maximizes the quantity $2rv$, because the quantity $2rv$ vanishes at the event horizon and decays back to zero as $r\rightarrow\infty$.

We also must check that $-\frac12 q^{-2}\pd_r(\Delta\pd_r(2rv))\ge 0$ for $r\ge r_*$. Indeed,
$$-\frac12 q^{-2}\pd_r(\Delta\pd_r w)=q^{-2}\pd_r\left(\frac{\Delta}{r^2+a^2}\frac{r^2-a^2}{r^2+a^2}\frac{r^2-4Mr+a^2}{r^2+a^2}\right).$$
Since the last quantity is a derivative of a product of three functions that are positive (for $r>r_*$) and increasing to $1$, it is necessarily positive. This verifies statement (iv). Statement (v) requires no further justification. Finally, we turn to statement (vi). First, we observe that for $r<r_*$, since $w=w(r_*)$ is constant, then
$$\pd_r\tilde{K}^{rr}=(r^2+a^2)\pd_r\left(\frac{w(r_*)}{2rv}\right)=-\frac{w(r_*)(r^2+a^2)}{(2rv)^2}\pd_r(2rv).$$
Since we have determined that $2rv$ has a maximum at $r_*$, it follows that $\tilde{K}^{rr}$ decreases to a minimum at $r_*$. Recall further that by the choice of $w$ for $r\ge r_*$, $\tilde{K}^{rr}$ is constant for $r\ge r_*$. It suffices to check the value of $\tilde{K}^{rr}$ at $r=r_*$. (That is, if $\tilde{K}^{rr}(r_*)>0$, then necessarily $\tilde{K}^{rr}>0$ and thus $K^{rr}\ge 0$.) First, we observe that 
$$v(r_*)\int_{r_{trap}}^{r_*}\frac1v<r_*-r_{trap},$$
since $v$ is decreasing from $r_{trap}$ to $r_*$. Using this fact, we examine the sign of $\pd_r\left(\frac{u}{r^2+a^2}\right)$ at $r_*$.
\begin{align*}
\left.\frac{(r_*^2+a^2)^2v(r_*)}{w(r_*)}\pd_r\left(\frac{u}{r^2+a^2}\right)\right|_{r_*} 
&= \frac{(r_*^2+a^2)v(r_*)\pd_ru|_{r_*}}{w(r_*)}-\frac{2r_*v(r_*)u(r_*)}{w(r_*)} \\
&= r_*^2+a^2-\frac{2r_*v(r_*)}{w(r_*)}\int_{r_{trap}}^{r_*}\pd_ru \\
&= r_*^2+a^2-\frac{2r_*v(r_*)}{w(r_*)}\int_{r_{trap}}^{r_*}\frac{w}{v} \\
&= r_*^2+a^2-2r_*v(r_*)\int_{r_{trap}}^{r_*}\frac1{v} \\
&> r_*^2+a^2-2r_*(r_*-r_{trap}).
\end{align*}
Since $r_*^2+a^2=4Mr_*$, it suffices to show
$$r_*^2+a^2-2r_*(r_*-r_{trap})=2r_*(2M-(r_*-r_{trap}))\ge 0.$$
This is indeed the case for all $|a|\le M$.
This concludes the proof of Lemma \ref{k:choice_for_u_and_w_lem}.
\end{proof}

We conclude the proof of Proposition \ref{k:partial_morawetz_prop} by analyzing the asymptotics for large $r$, in particular for $r>r_*$.
\begin{align*}
v &= \frac{\Delta}{(r^2+a^2)^2}=O(r^{-2}) \\
\pd_rv &= O(r^{-3}) \\
w_0 &= 2rv = O(r^{-1}) \\
u &= r^2+a^2-c^2 \text{ since }u'=2r \\
X_0 &= uv\pd_r = \frac{\Delta}{r^2+a^2}\pd_r+O\left(\frac{M^2}{r^2}\right)\pd_r \\
{}^{(Q)}K^{\alpha\beta} &= -u\pd_rv\frac{Q^{\alpha\beta}}{q^2}=O(r^{-1})\frac{Q^{\alpha\beta}}{q^2} \\
K^{rr} &=O(r^{-3})\tilde{K}^{rr}=O\left(\frac{M^2}{r^3}\right) \\
K &=q^{-2}\pd_r\left(\frac{\Delta}{r^2+a^2}\frac{r^2-a^2}{r^2+a^2}\frac{r^2-4Mr+a^2}{r^2+a^2}\right) = O\left(\frac{M}{r^4}\right).
\end{align*}
This completes the proof of Proposition \ref{k:partial_morawetz_prop}.
\end{proof}

\subsection{Remaining issues}

So far, we have a current $J[X_0,w_0]$ with a non-negative divergence. However, there are still a few modifications that must be made to the current $J[X_0,w_0]$ to establish a proper Morawetz estimate. \\
\bp The vectorfield $X_0$ behaves like $\log(1-\frac{r_H}r)\pd_r$ near the event horizon. So we temper $X$ (and also $w$) by choosing from a family $X_{\epsilon_{temper}}$, $w_{\epsilon_{temper}}$ where the parameter $\epsilon_{temper}>0$, which represents a deviation from $X_0$ and $w_0$, will be small. This will unfortunately introduce an error term of the form $-q^{-2}V_{\epsilon_{temper}}\psi^2$ to the divergence of $J$. (See \S\ref{k:temper_correction_sec}.) \\
\bp The coefficient $K^{rr}$ of $(\pd_r\psi)^2$ vanishes at the event horizon. This can be fixed by using the standard redshift argument, which necessitates a change of foliation near the event horizon and the addition of a redshift vectorfield $\epsilon_{redshift}Y$. This new foliation will be used for the remainder of the chapter. (See \S\ref{k:redshift_correction_sec}.) \\
\bp The coefficient $K$ of $\psi^2$ vanishes for $r\in [r_H,r_*]$. In fact, given the error term introduced by the tempering of $X$ and $w$, it is even slightly negative, with a smallness parameter $\epsilon_{temper}$. This will be fixed by observing a local Hardy estimate. (See \S\ref{k:hardy_correction_sec}.) \\
\bp The coefficient $K^{tt}$ of $(\pd_t\psi)^2$ vanishes entirely. (Although there is a $(\pd_t\psi)^2$ term corresponding to ${}^{(Q)}K^{tt}$ when $a\ne 0$, we want uniform estimates in $a$.) This is easily fixed by adding to $w$ the term $\epsilon_{\pd_t}w_{\pd_t}$. (See \S\ref{k:dt_correction_sec}.) \\
\bp The purpose of the current $J$ is to apply the divergence theorem in the form of Proposition \ref{general_divergence_estimate_prop}. So far, we have paid careful attention to the divergence of $J$, but the boundary terms also must be shown to have good properties. We add $\pd_t$ to $X$ and observe that there are only small errors on the boundary. (See \S\ref{k:boundary_correction_sec}.)

After all the above modifications take place, the final current $J$ will be of the following form.
\begin{align*}
J&=J[X,w] \\
X&=X_{\epsilon_{temper}}+\epsilon_{redshift}Y+\pd_t\\ 
w&=w_{\epsilon_{temper}}+\epsilon_{\pd_t}w_{\pd_t}
\end{align*}
And the current $J$ will satisfy
\begin{multline*}
\int_{\Sigma_t}\frac{M^2}{r^3}(\pd_r\psi)^2+\frac1r\left(1-\frac{r_{trap}}r\right)^2\left[|\sla\nabla\psi|^2+\frac{M^2}{r^2}(\pd_t\psi)^2\right]+\frac{M}{r^4}\psi^2 \\
\lesssim \int_{\Sigma_t}div J[X,w]-(2X(\psi)+w\psi)\Box_g\psi.
\end{multline*}

\subsection{The tempered vectorfield $X_{\epsilon_{temper}}$ and function $w_{\epsilon_{temper}}$}\label{k:temper_correction_sec}

The tempering of $X_0$ and $w_0$ to produce $X_{\epsilon_{temper}}$ and $w_{\epsilon_{temper}}$ can be done the same way as in \S\ref{s:temper_correction_sec}. We will simply state the result in the context of Kerr without proof.

\begin{lemma}\label{k:temper_lemma}
Suppose $\Box_g\psi=0$ and $\pd_\phi\psi=0$. For all $\epsilon_{temper}>0$ sufficiently small, there exist a modified vectorfield $X_{\epsilon_{temper}}$ and a modified function $w_{\epsilon_{temper}}$ agreeing, respectively, with the vectorfield $X_0$ and function $w_0$ from Proposition \ref{k:partial_morawetz_prop} outside of a small neighborhood of the event horizon and satisfying
\begin{multline*}
\left[\frac{M^2}{r^3}\left(1-\frac{r_H}r\right)^2(\pd_r\psi)^2+\frac1r\left(1-\frac{r_{trap}}r\right)^2\frac{Q^{\alpha\beta}}{q^2}\pd_\alpha\psi\pd_\beta\psi+\frac{M}{r^4}1_{r>r_*}\psi^2\right]-q^{-2}V_{\epsilon_{temper}}\psi^2 \\
\lesssim divJ[X_{\epsilon_{temper}},w_{\epsilon_{temper}}],
\end{multline*}
with a constant independent of $\epsilon_{temper}$. Moreover, the function $V_{\epsilon_{temper}}$ is supported in a neighborhood of the event horizon and satisfies
\begin{equation}\label{k:V_temper_bound_eqn}
||V_{\epsilon_{temper}}||_{L^1(r)}\le \epsilon_{temper}.
\end{equation}
\textbf{This new estimate is in fact worse than the estimate from Proposition \ref{k:partial_morawetz_prop}, but the advantage gained in this lemma is that, unlike the vectorfield $X_0$, the vectorfield $X_{\epsilon_{temper}}$ is regular up to the event horizon.}
\end{lemma}

\begin{remark}
As observed in the proof of Lemma \ref{s:temper_lemma}, the function $w_{\epsilon_{temper}}$ is supported outside a neighborhood of the event horizon and $\pd_r(q^2w_{\epsilon_{temper}})$ is nonnegative. These facts will be used when computing boundary terms.
\end{remark}

\subsection{The redshift vectorfield $\epsilon_{redshift}Y$}\label{k:redshift_correction_sec}

Now, we begin to use the coordinates corresponding to the $\delh$-foliation that is described by the following Lemma.

\begin{lemma}\label{k:foliation_lem}
For all sufficiently small $\delta>0$ (usually denoted by $\delh$ outside of this lemma) there exists a time foliation agreeing with the time foliation of the Boyer-Lindquist coordinates for $r>r_H+\delta$ such that the metric under this foliation satisfies
$$g^{tr}\ge 0,$$
$$q^2g^{tr}=r_H^2+O((r-r_H)^2) \text{ near }r=r_H,$$
and
$$\pd_r(q^2g^{tt})-2\frac{q^{tr}}{g^{rr}}\pd_r(q^2g^{tr})>0$$
for all $r\in [0,\delta]$.
\end{lemma}

\begin{proof}
Denote the Boyer-Lindquist coordinates by $(\bar{t},\bar{r},\bar{\theta},\bar{\phi})$. Then we define a new set of coordinates by\footnote{Usually, one also changes the variable $\phi$, but this is not important for problems with axissymmetry.}
\begin{align*}
t&=\bar{t}+\int_{r_H+\delta}^{\bar{r}}\frac{T(x)dx}{x^2-2Mx+a^2} \\
r&=\bar{r} \\
\theta&=\bar{\theta} \\
\phi&=\bar{\phi}
\end{align*}
for some function $T(\bar{r})$ supported for $\bar{r}\le r_H+\delta$. The time foliation is given by the new coordinate $t$, and the metric can be found by the change of coordinates. Under this change of coordinates, since $r=\bar{r}$ everywhere, we have
\begin{align*}
\pd_t&=\pd_{\bar{t}} \\
\pd_r&=\pd_{\bar{r}}-\frac{T}{\Delta}\pd_{\bar{t}} \\
\pd_\theta&=\pd_{\bar{\theta}} \\
\pd_\phi&=\pd_{\bar{\phi}}
\end{align*}
The inverse metric in these coordinates can be found by computing the following coordinate-invariant quantity.
\begin{align*}
q^2\pd^\lambda\psi\pd_\lambda\psi &= -\frac{(r^2+a^2)^2}{\Delta}(\pd_{\bar{t}}\psi)^2+\Delta(\pd_{\bar{r}}\psi)^2+q^2Q^{\bar{\alpha}\bar{\beta}}\pd_{\bar{\alpha}}\psi\pd_{\bar{\beta}}\psi \\
&=-\frac{(r^2+a^2)^2}{\Delta}(\pd_t\psi)^2+\Delta\left(\pd_r\psi+\frac{T}{\Delta}\pd_t\psi\right)^2+q^2Q^{\alpha\beta}\pd_\alpha\psi\pd_\beta\psi \\
&=\frac{T^2-(r^2+a^2)^2}{\Delta}(\pd_t\psi)^2+2T\pd_t\psi\pd_r\psi+\Delta(\pd_r\psi)^2+q^2Q^{\alpha\beta}\pd_\alpha\psi\pd_\beta\psi
\end{align*}
It follows that in these coordinates,
\begin{align*}
q^2g^{tt}&=\frac{T^2-(r^2+a^2)^2}{\Delta} \\
q^2g^{tr}&= T \\
q^2g^{rr}&= \Delta \\
q^2Q^{\alpha\beta}&=q^2Q^{\bar{\alpha}\bar{\beta}}
\end{align*}
If we impose the condition
$$T=r_H^2+a^2+O((r-r_H)^2)$$
then all of these metric components are regular and the first two conditions of the lemma are satisfied. We now examine the other condition. Multiplying by $q^2g^{rr}$, we obtain
\begin{align*}
q^2g^{rr}\pd_r(q^2g^{tt})-2q^2g^{tr}\pd_r(q^2g^{tr}) &= \pd_r\left(q^2g^{rr}q^2g^{tt}-(q^2g^{tr})^2\right)-q^2g^{tt}\pd_r(q^2g^{rr}) \\
&=-\pd_r\left((r^2+a^2)^2\right)-\frac{T^2-(r^2+a^2)^2}{\Delta}\pd_r\Delta
\end{align*}
Multiplying again by $q^2g^{rr}=\Delta$, we obtain
$$[(r^2+a^2)^2-T^2]\pd_r\Delta-\Delta\pd_r(r^2+a^2)^2.$$
Since $T(r_H)=r_H^2+a^2$, this quantity vanishes at $r_H$. Since $T'(r_H)=0$, the first derivative of this quantity also vanishes. If $T''(r_H)$ is sufficiently negative, then this quantity can be made positive in a neighborhood of $r_H$. Since we obtained this quantity by multiplying the original quantity by $\Delta^2$, then the original quantity is strictly positive on the support of $T$.
\end{proof}

 The change this has on the partial Morawetz estimate is minimal--it only means that the term ``$(\pd_r\psi)^2$'' really becomes $(\pd_r\psi+f(r)\pd_t\psi)^2$ for some function $f$ supported where $r<r_H+\delh$. Since this term still has the same sign, it may simply be ignored in that region. (It wasn't very useful near the event horizon anyway--it degenerated like $O((r-r_H)^2)$.)

We now show that by considering the effects of the $\delh$-foliation, we can define a redshift vectorfield $Y$ that will give uniform control of $(\pd_r\psi)^2$ and $(\pd_t\psi)^2$ in a neighborhood of the event horizon larger than $\delh$. This is summarized in the following lemma.

\begin{lemma}\label{k:redshift_lem}
Suppose $\Box_g\psi=0$ and $\pd_\phi\psi=0$. Let $|a|<M$. On a $\delh$-foliation guaranteed by Lemma \ref{k:foliation_lem}, there exists a vectorfield $Y$ supported on the region $r\in[r_H,r_H+2\delh]$ satisfying
$$1_{r-r_H\in [0,3\delh/2]}\left[(\pd_r\psi)^2+(\pd_t\psi)^2\right]\lesssim divJ[Y]+1_{r-r_H\in (\delh,2\delh)}\frac{Q^{\alpha\beta}}{q^2}\pd_\alpha\psi\pd_\beta\psi.$$
Moreover, $Y=Y^r(r)\pd_r+Y^t(r)\pd_t$, and $Y^r(r_H)< 0$.
\end{lemma}
\begin{proof}
Let $Y=Y^r(r)\pd_r+Y^t(r)\pd_t$. Then by Lemma \ref{divJ_lem},
$$divJ[Y]=K_Y^{\mu\nu}\pd_\mu\psi\pd_\nu\psi,$$
where
\begin{align*}
K_Y^{\mu\nu}&=g^{\mu\lambda}\pd_\lambda Y^\nu+g^{\nu\lambda}\pd_r Y^\mu-divYg^{\mu\nu}-Y^\lambda\pd_\lambda g^{\mu\nu} \\
&=g^{\mu r}\pd_r Y^\nu+g^{\nu r}\pd_r Y^\mu-\frac{1}{q^2}\pd_r(q^2 Y^r)g^{\mu\nu}-Y^r\pd_r g^{\mu\nu} \\
&=g^{\mu r}\pd_r Y^\nu+g^{\nu r}\pd_r Y^\mu-\frac{1}{q^2}\pd_r(Y^r q^2g^{\mu\nu}).
\end{align*}
It follows that
\begin{align*}
q^2K_Y^{rr}&=\Delta \pd_rY^r-Y^r\pd_r\Delta \\
2q^2K_Y^{(tr)}&=2(q^2g^{rr}\pd_r Y^t-\pd_r(q^2g^{tr})Y^r) \\
q^2K_Y^{tt}&=2q^2g^{tr}\pd_rY^t-\pd_r(Y^rq^2g^{tt}) \\
q^2\hspace{3pt}{}^{(Q)}K_Y^{\alpha\beta}&=-\pd_r Y^r Q^{\alpha\beta}
\end{align*}
We will take $Y^r$ to be a negative constant\footnote{It is traditional to take $\pd_rY^r<0$ near the event horizon to also gain control of the angular derivatives, but this is not necessary here since the current $J[X_{\epsilon_{temper}},w_{\epsilon_{temper}}]$ already has good control of these derivatives and such a choice complicates the proof.} on $r-r_H\in [0,\delh]$ and smoothly increasing to $0$ on $r-r_H\in (\delh,2\delh)$.

With this choice, $K_Y^{rr}$ has the right sign for $r-r_H\in[0,2\delh]$. In particular, this is due to the fact that the term $-Y^r\pd_r\Delta$ is strictly positive at $r=r_H$. \textbf{This fact is related to the surface gravity of the Kerr black hole, and is not true in the extremal case $|a|=M$.}

However, ${}^{(Q)}K_Y^{\alpha\beta}$ has the wrong sign when $\pd_rY^r\ne 0$. This explains the error term on the right side of the estimate in the lemma.

We choose $Y^t$ to have support coinciding with $g^{tr}$ and such that $K_Y^{(tr)}=0$. That is, $Y^t$ is determined by the relation
$$\pd_rY^t=\frac{1}{g^{rr}}\pd_r(q^2g^{tr})Y^r$$
and the condition $Y^t(\delh)=0$.
Note that although $g^{rr}$ vanishes at the event horizon, the quantity $\frac{1}{g^{rr}}\pd_r(q^2g^{tr})$ remains bounded, because $\pd_r(q^2g^{tr})$ also vanishes to first order as $r\rightarrow r_H$. (This was a condition specified in Lemma \ref{k:foliation_lem}.)

Having made the above choice for $Y^t$, we have in the region $r-r_H\in[0,\delh]$ that
\begin{align*}
q^2K_Y^{tt}&=2\frac{g^{tr}}{g^{rr}}\pd_r(q^2g^{tr})Y^r-\pd_r(Y^rq^2g^{tt}) \\
&= (-Y^r)\left[\pd_r(q^2g^{tt})-2\frac{g^{tr}}{g^{rr}}\pd_r(q^2g^{tr})\right].
\end{align*}
By Lemma \ref{k:foliation_lem}, this is bounded below by a positive constant on $r-r_H\in [0,\delh]$.
 
In the region $r-r_H\in (\delh,2\delh)$,
$$q^2K_Y^{tt}=-\pd_r(Y^rq^2g^{tt}).$$
With an appropriate choice of $Y^r$, this term will be positive on $r-r_H\in [\delh,2\delh)$, as $Y^rq^2g^{tt}$ is positive and must go to zero as $Y^r$ vanishes.
\end{proof}

The previous lemma leads to the following corollary.
\begin{corollary}\label{k:redshift_cor}
Suppose $\pd_\phi\psi=0$. Then there exists an $\epsilon_{redshift}$ sufficiently small so that for all $\epsilon_{temper}$ sufficiently small,
\begin{multline*}
\left[\frac{M^2}{r^3}(\pd_r\psi)^2+1_{r-r_H\in[0,3\delh/2]}(\pd_t\psi)^2+\frac1r\left(1-\frac{r_{trap}}r\right)^2\frac{Q^{\alpha\beta}}{q^2}\pd_\alpha\psi\pd_\beta\psi+\frac{M}{r^4}1_{r>r_*}\psi^2\right] \\
-r^{-2}V_{\epsilon_{temper}}\psi^2 \\
\lesssim divJ[X_{\epsilon_{temper}}+\epsilon_{redshift}Y,w_{\epsilon_{temper}}].
\end{multline*}
\end{corollary}

\subsection{The local Hardy estimate}\label{k:hardy_correction_sec}

The local Hardy estimate can be applied the same way as in \S\ref{s:hardy_correction_sec}. We will simply state the result in the context of Kerr without proof.
\begin{corollary}\label{k:hardy_cor}
Suppose $\pd_\phi\psi=0$. There exists $\epsilon_{temper}$ sufficiently small so that
\begin{multline*}
\int_{\Sigma_t}\left[\frac{M^2}{r^3}(\pd_r\psi)^2+1_{r-r_H\in[0,3\delh/2]}(\pd_t\psi)^2+\frac1r\left(1-\frac{r_{trap}}r\right)^2\frac{Q^{\alpha\beta}}{q^2}\pd_\alpha\psi\pd_\beta\psi+\frac{M}{r^4}\psi^2\right] \\
\lesssim \int_{\Sigma_t}divJ[X_{\epsilon_{temper}}+\epsilon_{redshift}Y,w_{\epsilon_{temper}}].
\end{multline*}
\end{corollary}

\subsection{The correction $\epsilon_{\pd_t}w_{\pd_t}$}\label{k:dt_correction_sec}

The correction $w_{\pd_t}$ can be used to complete the term with $(\pd_t\psi)^2$ the same way as in \S\ref{s:dt_correction_sec}. We will simply state the result in the context of Kerr without proof.
\begin{corollary}\label{k:dt_cor}
Suppose $\pd_\phi\psi=0$. Then there exist a function $w_{\pd_t}$ and a small $\epsilon_{\pd_t}>0$ so that
\begin{multline*}
\int_{\Sigma_t}\left[\frac{M^2}{r^3}(\pd_r\psi)^2+\left(1-\frac{r_{trap}}r\right)^2\left[\frac1r\frac{Q^{\alpha\beta}}{q^2}\pd_\alpha\psi\pd_\beta\psi+\frac{M^2}{r^3}(\pd_t\psi)^2\right]+\frac{M}{r^4}\psi^2\right] \\
\lesssim \int_{\Sigma_t}divJ[X_{\epsilon_{temper}}+\epsilon_{redshift}Y,w_{\epsilon_{temper}}+\epsilon_{\pd_t}w_{\pd_t}].
\end{multline*}
\end{corollary}

\subsection{The correction $\pd_t$ and the flux terms}\label{k:boundary_correction_sec}

The correction $\pd_t$ can be used to obtain good boundary terms (with only small errors) the same way as in \S\ref{s:boundary_correction_sec}. We will simply state the result in the context of Kerr without proof.

\begin{lemma}\label{k:morawetz_flux_lem}
The current
$$J^t[X_{\epsilon_{temper}}+\epsilon_{redshift}Y+\pd_t,w_{\epsilon_{temper}}+\epsilon_{\pd_t}w_{\pd_t}]$$
satisfies the following flux estimates.

On the hypersurface $\Sigma_t$,
\begin{multline*}
\int_{\Sigma_t} (L\psi)^2+\frac{Q^{\alpha\beta}}{q^2}\pd_\alpha\psi\pd_\beta\psi+r^{-2}\psi^2 \\
\lesssim  \int_{\Sigma_t}-J^t[X_{\epsilon_{temper}}+\epsilon_{redshift}Y+\pd_t,w_{\epsilon_{temper}}+\epsilon_{\pd_t}w_{\pd_t}] + Err_{\Sigma_t},
\end{multline*}
where 
$$Err_{\Sigma_t}=\int_{\Sigma_t} r^{-1}|\psi L\psi|+\frac{M^2}{r^2}\left[\chi_H(\pd_r\psi)^2+(\pd_t\psi)^2\right],$$
and $\chi_H=1-\frac{r_H}r$.

On the hypersurface $H_{t_1}^{t_2}$,
$$\int_{H_{t_1}^{t_2}}\frac{Q^{\alpha\beta}}{q^2}\pd_\alpha\psi\pd_\beta\psi\lesssim  \int_{H_{t_1}^{t_2}}-J^t[X_{\epsilon_{temper}}+\epsilon_{redshift}Y+\pd_t,w_{\epsilon_{temper}}+\epsilon_{\pd_t}w_{\pd_t}] + Err_{H_{t_1}^{t_2}},$$
where
$$Err_{H_{t_1}^{t_2}}=\int_{H_{t_1}^{t_2}}(\pd_t\psi)^2.$$
\end{lemma}

\subsection{The Morawetz estimate}

Finally, we conclude with the Morawetz estimate.
\begin{theorem}\label{k:morawetz_thm} (Morawetz estimate for Kerr $|a|<M$, assuming axisymmetry)
The following estimate holds for any sufficiently regular function $\psi$ satisfying
$$\pd_\phi\psi=0,$$
and decaying sufficiently fast as $r\rightarrow\infty$.
\begin{multline*}
\int_{H_{t_1}^{t_2}}\frac{Q^{\alpha\beta}}{q^2}\pd_\alpha\psi\pd_\beta\psi+\int_{\Sigma_{t_2}}(L\psi)^2+\frac{Q^{\alpha\beta}}{q^2}\pd_\alpha\psi\pd_\beta\psi+r^{-2}\psi^2 \\
+\int_{t_1}^{t_2}\int_{\Sigma_t}\frac{M^2}{r^3}(\pd_r\psi)^2+\frac1r\left(1-\frac{r_{trap}}r\right)^2\left[\frac{Q^{\alpha\beta}}{q^2}\pd_\alpha\psi\pd_\beta\psi+\frac{M^2}{r^2}(\pd_t\psi)^2\right]+\frac{M}{r^4}\psi^2 \\
\lesssim \int_{\Sigma_{t_1}}(L\psi)^2+\frac{Q^{\alpha\beta}}{q^2}\pd_\alpha\psi\pd_\beta\psi+r^{-2}\psi^2 + Err,
\end{multline*}
where
\begin{align*}
Err &= Err_1+Err_2+Err_\Box \\
Err_1 &= \int_{\Sigma_{t_2}} r^{-1}|\psi L\psi| \\
Err_2 &= \int_{H_{t_1}^{t_2}}(\pd_t\psi)^2+\int_{\Sigma_{t_2}}\frac{M^2}{r^2}\left[\chi_H(\pd_r\psi)^2+(\pd_t\psi)^2\right] \\
Err_\Box&=\int_{t_1}^{t_2}\int_{\Sigma_t}(2X(\psi)+w\psi)\Box_g\psi,
\end{align*}
and $\chi_H=1-\frac{r_H}r$.
\end{theorem}

\begin{proof}
Apply Proposition \ref{general_divergence_estimate_prop} to the current
$$J[X_{\epsilon_{temper}}+\epsilon_{redshift}Y+\pd_t,w_{\epsilon_{temper}}+\epsilon_{\pd_t}w_{\pd_t}]$$
and invoke Corollary \ref{k:dt_cor} and Lemma \ref{k:morawetz_flux_lem}.
\end{proof}

\section{The $r^p$ estimate}\label{k:rp_sec}

\subsection{The pre-$r^p$ identity}

We prove the following lemma, which has a much simpler, well-known analogue in Minkowski spacetime (see Lemma \ref{m:p_ee_identity_lem}). It will be used for the proof of Proposition \ref{k:incomplete_p_estimate_prop}.
\begin{lemma}\label{k:p_ee_identity_lem} Suppose $\pd_\phi\psi=0$. Let $\alpha=\frac{\Delta}{r^2+a^2}$ and $L=\alpha\pd_r+\pd_t$. For any function $f=f(r)$ supported where $r>r_H+\delh$, the following identity holds.
\begin{align*}
&\int_{\Sigma_{t_2}}\left[\left(1-\alpha\frac{a^2\sin^2\theta}{r^2+a^2}\right)\frac{r^2+a^2}{q^2}f\left(\alpha^{-1}L\psi+\frac{r}{r^2+a^2}\psi\right)^2 +\frac{\alpha^{-1}f}{q^2}(\pd_\theta\psi)^2+\epsilon\frac{rf'}{q^2}\psi^2\right. \\
&\hspace{1.3in}\left.+\alpha\frac{a^2\sin^2\theta}{q^2}f\left(\pd_r\psi+\frac{r}{r^2+a^2}\psi\right)^2 +\frac{a^2f}{q^2(r^2+a^2)}\psi^2-\frac1{q^2}\pd_r(rf\psi^2)\right] \\
&+\int_{t_1}^{t_2}\int_{\Sigma_t}
\left[
  \left(\frac{2rf}{r^2+a^2}-f'\right)\frac{Q^{\alpha\beta}}{q^2}\pd_\alpha\psi\pd_\beta\psi
+ \alpha f'\frac{r^2+a^2}{q^2}\left(\alpha^{-1}L\psi+\frac{(1-\epsilon)r}{r^2+a^2}\psi\right)^2 
\vphantom{
+ \epsilon\alpha\left((1-\epsilon)f'-rf''\right)\frac{\psi^2}{q^2} 
+ \alpha'\left(\frac{r^2-a^2}{r^2+a^2}f-\epsilon r f'\right)\frac{\psi^2}{q^2}
+ \frac{a^2}{r^2+a^2}\left(\alpha((1-\epsilon)^2-2)f'-\frac{4\alpha rf}{r^2+a^2}\right)\frac{\psi^2}{q^2}
}\right. \\
&\hspace{2.4in} + \epsilon\alpha\left((1-\epsilon)f'-rf''\right)\frac{\psi^2}{q^2} -\alpha'\alpha^{-2}f\frac{r^2+a^2}{q^2}(L\psi)^2 \\
&\hspace{0.9in}
\left.\vphantom{
  \left(\frac{2rf}{r^2+a^2}-f'\right)\frac{Q^{\alpha\beta}}{q^2}\pd_\alpha\psi\pd_\beta\psi
+ \alpha f'\frac{r^2+a^2}{q^2}\left(\alpha^{-1}L\psi+\frac{(1-\epsilon)r}{r^2+a^2}\psi\right)^2 
+ \epsilon\alpha\left((1-\epsilon)f'-rf''\right)\frac{\psi^2}{q^2} 
}
+ \alpha'\left(-\epsilon r f' +\frac{r^2-a^2}{r^2+a^2}f\right)\frac{\psi^2}{q^2}
+ \frac{a^2}{r^2+a^2}\left(-(1+\epsilon)\alpha f'-\frac{4\alpha rf}{r^2+a^2}\right)\frac{\psi^2}{q^2}
\right] \\
=&\int_{\Sigma_{t_1}}\left[\left(1-\alpha\frac{a^2\sin^2\theta}{r^2+a^2}\right)\frac{r^2+a^2}{q^2}f\left(\alpha^{-1}L\psi+\frac{r}{r^2+a^2}\psi\right)^2 +\frac{\alpha^{-1}f}{q^2}(\pd_\theta\psi)^2+\epsilon\frac{rf'}{q^2}\psi^2\right. \\
&\hspace{1.3in}\left.+\alpha\frac{a^2\sin^2\theta}{q^2}f\left(\pd_r\psi+\frac{r}{r^2+a^2}\psi\right)^2 +\frac{a^2f}{q^2(r^2+a^2)}\psi^2-\frac1{q^2}\pd_r(rf\psi^2)\right] \\
&+\int_{t_1}^{t_2}\int_{\Sigma_t}-\left(2\alpha^{-1}fL\psi+\frac{2rf}{r^2+a^2}\psi\right)\Box_g\psi.
\end{align*}
\end{lemma}

\begin{proof}
We will use Proposition \ref{general_divergence_estimate_prop} together with the following current template.
$$J[X,w,m]_\mu = T_{\mu\nu} X^\nu +w\psi\pd_\mu\psi-\frac12\psi^2\pd_\mu w+m_\mu \psi^2,$$
$$T_{\mu\nu}=2\pd_\mu\psi\pd_\nu\psi-g_{\mu\nu}\pd^\lambda\psi\pd_\lambda\psi.$$

Assume for now that $\Box_g\psi=0$. Let $\alpha = \frac{\Delta}{r^2+a^2}$, and observe that
$$L=\alpha\pd_r+\pd_t,$$
$$q^2g^{rr}=(r^2+a^2)\alpha,$$
$$q^2g^{tt}=-(r^2+a^2)\alpha^{-1}.$$

\begin{lemma}\label{divJphiX_lem}
Without appealing directly to the particular expression for $\alpha$, one can deduce the following identity.
\begin{multline*}
\frac{q^2}{r^2+a^2}divJ[\alpha^{-1}fL] \\
=(\alpha^{-1}f)'(L\psi)^2-\frac{2rf}{r^2+a^2}\left(\alpha(\pd_r\psi)^2-\alpha^{-1}(\pd_t\psi)^2\right)-f'\frac{Q^{\alpha\beta}}{r^2+a^2}\pd_\alpha\psi\pd_\beta\psi.
\end{multline*}
\end{lemma}
\begin{proof}
Note that
$$div J[X] = K^{\mu\nu}\pd_\mu\psi\pd_\nu\psi,$$
where
$$K^{\mu\nu}=2g^{\mu\lambda}\pd_\lambda X^\nu-X^\lambda \pd_\lambda(g^{\mu\nu})-div X g^{\mu\nu}.$$

Set $X=\alpha^{-1}f(\alpha\pd_r+\pd_t)=f\pd_r+\alpha^{-1}f\pd_t$. From the above formula, since $g^{rt}=0$,
$$\frac{q^2}{r^2+a^2}(K^{tr}+K^{rt})=2\frac{q^2}{r^2+a^2}g^{rr}\pd_r X^t=2\alpha\pd_r(\alpha^{-1} f).$$
Thus, the expression for $\frac{q^2}{r^2+a^2}divJ[\alpha^{-1}fL]$ will have a mixed term of the form 
$$2\alpha\pd_r(\alpha^{-1}f)\pd_r\psi\pd_t\psi.$$
Note that
\begin{align*}
(\alpha^{-1} f)'(L\psi)^2 &= (\alpha^{-1} f)'(\alpha\pd_r\psi+\pd_t\psi)^2 \\
&= \alpha^2(\alpha^{-1}f)'(\pd_r\psi)^2+2\alpha(\alpha^{-1}f)'\pd_r\psi\pd_t\psi +(\alpha^{-1} f)'(\pd_t\psi)^2.
\end{align*}
We now compute the $(\pd_r\psi)^2$ and $(\pd_t\psi)^2$ components, subtracting the part that will be grouped with the $(L\psi)^2$ term.
\begin{align*}
\frac{q^2}{r^2+a^2}K^{rr}-\alpha^2(\alpha^{-1}f)' &= \frac{q^2}{r^2+a^2}\left[2g^{rr}\pd_rX^r-X^r\pd_r g^{rr}-\frac1{q^2}\pd_r(q^2X^r)g^{rr}\right]-\alpha^2(\alpha^{-1}f)' \\
&= \frac{q^2}{r^2+a^2}\left[2g^{rr}\pd_rX^r-\frac1{q^2}\pd_r(q^2g^{rr}X^r)\right]-\alpha^2(\alpha^{-1}f)' \\
&= 2\alpha\pd_r f -\frac{1}{r^2+a^2}\pd_r\left((r^2+a^2)\alpha f\right)-\alpha^2(\alpha^{-1}f)' \\
&= -\frac{2r \alpha f}{r^2+a^2}
\end{align*}
and
\begin{align*}
\frac{q^2}{r^2+a^2}K^{tt}-(\alpha^{-1}f)' &= \frac{q^2}{r^2+a^2}\left[-X^r\pd_r g^{tt}-\frac{1}{q^2}\pd_r(q^2X^r)g^{tt}\right]-(\alpha^{-1}f)' \\
&= -\frac{q^2}{r^2+a^2}\frac{1}{q^2}\pd_r(q^2 g^{tt} X^r) -(\alpha^{-1}f)' \\
&= -\frac{1}{r^2+a^2}\pd_r\left((r^2+a^2)(-\alpha^{-1}) f\right)-(\alpha^{-1}f)' \\
&= \frac{2r\alpha^{-1}f}{r^2+a^2}.
\end{align*}
Finally,
\begin{align*}
\frac{q^2}{r^2+a^2}{}^{(Q)}K^{\alpha\beta} &= \frac{q^2}{r^2+a^2}\left[-X^r\pd_r {}^{(Q)}g^{\alpha\beta}-\frac1{q^2}\pd_r\left(q^2X^r\right){}^{(Q)}g^{\alpha\beta}\right] \\
&= \frac{q^2}{r^2+a^2}\left[-\frac{1}{q^2}\pd_r\left(q^2{}^{(Q)}g^{\alpha\beta}X^r\right)\right] \\
&= -\frac{1}{r^2+a^2}\pd_r(Q^{\alpha\beta}f) \\
&= -f' \frac{Q^{\alpha\beta}}{r^2+a^2}.
\end{align*}
Combining all these terms gives the identity stated in the lemma. 
\end{proof}

Next, we choose $w=\frac{2rf}{r^2+a^2}$ to directly cancel the middle term in the above lemma.
\begin{lemma}\label{divJphiXw_lem}
\begin{multline*}
\frac{q^2}{r^2+a^2}divJ\left[\alpha^{-1}fL,\frac{2rf}{r^2+a^2}\right] \\
= (\alpha^{-1}f)'(L\psi)^2+\left(\frac{2rf}{r^2+a^2}-f'\right)\frac{Q^{\alpha\beta}}{r^2+a^2}\pd_\alpha\psi\pd_\beta\psi-\frac12\frac{q^2}{r^2+a^2}\Box_g\left(\frac{2rf}{r^2+a^2}\right)\psi^2.
\end{multline*}
\end{lemma}
\begin{proof}
Note that
$$divJ[0,w]=wg^{\mu\nu}\pd_\mu\psi\pd_\nu\psi-\frac12\Box_gw \psi^2.$$
We compute the new terms only.
\begin{multline*}
\frac{q^2}{r^2+a^2}divJ\left[0,\frac{2rf}{r^2+a^2}\right] = \frac{2rf}{r^2+a^2}\frac{q^2 g^{\alpha\beta}}{r^2+a^2}\pd_\alpha\psi\pd_\beta\psi -\frac12\frac{q^2}{r^2+a^2}\Box_g\left(\frac{2rf}{r^2+a^2}\right)\psi^2 \\
= \frac{2rf}{r^2+a^2}\left(\alpha (\pd_r\psi)^2-\alpha^{-1}(\pd_t\psi)^2\right) +\frac{2rf}{r^2+a^2}\frac{Q^{\alpha\beta}}{r^2+a^2}\pd_\alpha\psi\pd_\beta\psi -\frac12\frac{q^2}{r^2+a^2}\Box_g\left(\frac{2rf}{r^2+a^2}\right)\psi^2.
\end{multline*}
When adding these terms to the expression in Lemma \ref{divJphiX_lem}, the $\alpha(\pd_r\psi)^2-\alpha^{-1}(\pd_t\psi)^2$ terms cancel (this was the reason for the choice of $w=\frac{2rf}{r^2+a^2}$) and the result is as desired. 
\end{proof}

The term $-\frac12\frac{q^2}{r^2+a^2}\Box_g\left(\frac{2rf}{r^2+a^2}\right)\psi^2$ is like $-r^{-1}f''\psi^2$. In the future, when $f\sim r^p$, this will have a sign $-p(p-1)$. The sign will be negative if $p>1$, which is bad. So we include a divergence term to fix it. (But in doing so, we almost lose some other good terms--this is why we need a small parameter $\epsilon$.) This is the point of the following Lemma.
\begin{lemma}\label{k:rp_add_term_lem}
\begin{align*}
\alpha^{-1}f'& (L\psi)^2+\frac{q^2}{r^2+a^2}\left[-\frac12 \Box_g\left(\frac{2rf}{r^2+a^2}\right)\psi^2+(1-\epsilon)div\left(\psi^2\frac{r}{q^2}f'L\right)\right] \\
=& \alpha f' \left(\alpha^{-1}L\psi+\frac{(1-\epsilon)r}{r^2+a^2}\psi\right)^2 + \epsilon\alpha \frac{\left((1-\epsilon)f'-rf''\right)}{r^2+a^2}\psi^2 \\
&+ \alpha'\left(-\frac{\epsilon r f'}{r^2+a^2}+\frac{(r^2-a^2)f}{(r^2+a^2)^2}\right)\psi^2 
+ \frac{a^2}{r^2+a^2}\left(\frac{-(1+\epsilon)\alpha f'}{r^2+a^2}-\frac{4\alpha rf}{(r^2+a^2)^2}\right)\psi^2.
\end{align*}
\end{lemma}
\begin{proof}
First, we calculate
\begin{align*}
-\frac{q^2}{r^2+a^2}\frac12\Box_g\left(\frac{2rf}{r^2+a^2}\right)\psi^2 &=-\frac{1}{r^2+a^2}\pd_r\left((r^2+a^2)\alpha\pd_r\left(\frac{rf}{r^2+a^2}\right)\right)\psi^2 \\
&= -\frac{\alpha}{r^2+a^2}\pd_r\left((r^2+a^2)\pd_r\left(\frac{rf}{r^2+a^2}\right)\right) \psi^2 \\
&\hspace{1.8in} -\alpha'\pd_r\left(\frac{rf}{r^2+a^2}\right)\psi^2 \\
&= -\frac{\alpha r f''}{r^2+a^2}\psi^2 -\alpha'\pd_r\left(\frac{rf}{r^2+a^2}\right)\psi^2 \\
&\hspace{1in}-\frac{2\alpha a^2}{(r^2+a^2)^2}\left(f'+\frac{2r}{r^2+a^2}f\right)\psi^2.
\end{align*}
We also calculate
\begin{align*}
\frac{q^2}{r^2+a^2}div\left(\psi^2\frac{r}{q^2}f'L\right) &= \frac{1}{r^2+a^2}\pd_\alpha\left(\psi^2rf' L^\alpha\right) \\
&=\frac{rf'}{r^2+a^2}2\psi L\psi +\frac{\pd_r(rf'\alpha)}{r^2+a^2}\psi^2 \\
&=\frac{r f'}{r^2+a^2}2\psi L\psi +\frac{\alpha f'}{r^2+a^2}\psi^2 +\frac{\alpha r f''}{r^2+a^2}\psi^2+\frac{\alpha' r f'}{r^2+a^2}\psi^2.
\end{align*}
The first two terms in the last line almost complete a square (up to a term on the order of $\frac{a^2}{r^2+a^2}$) with the term $\alpha^{-1} f' (L\psi)^2$. The third term cancels with the first term from the previous calculation. However, it will be beneficial to introduce the factor $1-\epsilon$ that appears in the lemma, so that a good term appears with an $\epsilon$ factor. This is summarized by the following two calculations.
\begin{multline*}
\alpha^{-1}f'(L\psi)^2+(1-\epsilon)\left(\frac{rf'}{r^2+a^2}2\psi L\psi+\frac{\alpha f'}{r^2+a^2}\psi^2\right) \\
= \alpha f' \left(\alpha^{-1}L\psi+\frac{(1-\epsilon)r}{r^2+a^2}\psi\right)^2 -\frac{(1-\epsilon)^2r^2\alpha f'}{(r^2+a^2)^2}\psi^2 +\frac{(1-\epsilon)\alpha f'}{r^2+a^2}\psi^2 \\
= \alpha f' \left(\alpha^{-1}L\psi+\frac{(1-\epsilon)r}{r^2+a^2}\psi\right)^2 +\frac{\epsilon (1-\epsilon)\alpha f'}{r^2+a^2}\psi^2 +\frac{a^2(1-\epsilon)\alpha f'}{(r^2+a^2)^2}\psi^2
\end{multline*}
and
$$-\frac{\alpha r f''}{r^2+a^2}\psi^2 +(1-\epsilon)\frac{\alpha r f''}{r^2+a^2}\psi^2=-\epsilon \frac{\alpha r f''}{r^2+a^2}\psi^2.$$
Adding these terms together and ignoring the term with the $a^2$ factor yields
$$\alpha f' \left(\alpha^{-1}L\psi+\frac{(1-\epsilon)r}{r^2+a^2}\psi\right)^2 +\epsilon \alpha \frac{\left((1-\epsilon)f'-rf''\right)}{r^2+a^2}\psi^2.$$
All the remaining terms (which either contain a factor of $\alpha'\sim \frac{M}{r^2}$ or $\frac{a^2}{r^2+a^2}$) are
\begin{multline*}
\alpha'\left[-\pd_r\left(\frac{rf}{r^2+a^2}\right)+(1-\epsilon)\frac{rf'}{r^2+a^2}\right]\psi^2 \\
+\frac{a^2}{r^2+a^2}\left[-\frac{2\alpha}{r^2+a^2}\left(f'+\frac{2r}{r^2+a^2}f\right)+\frac{(1-\epsilon)\alpha f'}{r^2+a^2}\right]\psi^2
\end{multline*}
Adding both of these yields the result. 
\end{proof}

Thus, we have shown that if $\Box_g\psi=0$, then
\begin{multline*}
\frac{q^2}{r^2+a^2}divJ\left[\alpha^{-1} f L,\frac{2rf}{r^2+a^2},(1-\epsilon)\frac{rf'}{q^2}L\right] \\
= \alpha f' \left(\alpha^{-1}L\psi+\frac{(1-\epsilon)r}{r^2+a^2}\psi\right)^2 + \epsilon\alpha \frac{\left((1-\epsilon)f'-rf''\right)}{r^2+a^2}\psi^2 +\left(\frac{2rf}{r^2+a^2}-f'\right)\frac{Q^{\alpha\beta}}{r^2+a^2}\pd_\alpha\psi\pd_\beta\psi \\
- \alpha' \alpha^{-2}f(L\psi)^2
+ \alpha'\left(-\frac{\epsilon r f'}{r^2+a^2}+\frac{(r^2-a^2)f}{(r^2+a^2)^2}\right)\psi^2  \\
+ \frac{a^2}{r^2+a^2}\left(\frac{-(1+\epsilon)\alpha f'}{r^2+a^2}-\frac{4\alpha rf}{(r^2+a^2)^2}\right)\psi^2.
\end{multline*}
If we remove the assumption that $\Box_g\psi=0$, there is an additional term
$$\left(2X(\psi)+w\psi\right)\Box_g\psi = \left(2\alpha^{-1}fL\psi+\frac{2rf}{r^2+a^2}\psi\right)\Box_g\psi.$$

Finally, we turn to the boundary terms. Since we have assumed that $f$ is supported away from the event horizon, it suffices to compute $-J^t$.
\begin{lemma}\label{k:rp_boundary_terms_lem}
\begin{multline*}
-J^t\left[\alpha^{-1}fL,\frac{2rf}{r^2+a^2},(1-\epsilon)\frac{rf'}{q^2}L\right] \\
=\left(1-\alpha\frac{a^2\sin^2\theta}{r^2+a^2}\right)\frac{r^2+a^2}{q^2}f\left(\alpha^{-1}L\psi+\frac{r}{r^2+a^2}\psi\right)^2 +\frac{\alpha^{-1}f}{q^2}(\pd_\theta\psi)^2+\epsilon\frac{rf'}{q^2}\psi^2 \\
+\alpha\frac{a^2\sin^2\theta}{q^2}f\left(\pd_r\psi+\frac{r}{r^2+a^2}\psi\right)^2 +\frac{a^2f}{q^2(r^2+a^2)}\psi^2-\frac1{q^2}\pd_r(rf\psi^2).
\end{multline*}
\end{lemma}
\begin{proof}
We have
\begin{align*}
-J^t[\alpha^{-1}fL] &= -2\pd^t\psi\alpha^{-1}fL\psi+\alpha^{-1}fL^t\pd^\lambda\psi\pd_\lambda\psi \\
&=-2\alpha^{-1}f(g^{tt}+{}^{(Q)}g^{tt})\pd_t\psi L\psi \\
&\hspace{1.5in}+\alpha^{-1}f\left((g^{tt}+{}^{(Q)}g^{tt})(\pd_t\psi)^2+g^{rr}(\pd_r\psi)^2+{}^{(Q)}g^{\theta\theta}(\pd_\theta\psi)^2\right) \\
&=-\alpha^{-1}f(g^{tt}+{}^{(Q)}g^{tt})(\pd_t\psi)^2-2\alpha^{-1}f(g^{tt}+{}^{(Q)}g^{tt})\pd_t\psi\alpha\pd_r\psi+\alpha^{-1}fg^{rr}(\pd_r\psi)^2 \\
&\hspace{4in}+\alpha^{-1}f{}^{(Q)}g^{\theta\theta}(\pd_\theta\psi)^2 \\
&=\frac{r^2+a^2}{q^2}\left(1-\alpha\frac{a^2\sin^2\theta}{r^2+a^2}\right)\alpha^{-2}f(\pd_t\psi)^2 \\
&\hspace{0.2in}+2\frac{r^2+a^2}{q^2}\left(1-\alpha\frac{a^2\sin^2\theta}{r^2+a^2}\right)\alpha^{-1}f\pd_t\psi\pd_r\psi +\frac{r^2+a^2}{q^2}f(\pd_r\psi)^2+\frac{\alpha^{-1}f}{q^2}(\pd_\theta\psi)^2 \\
&=\frac{r^2+a^2}{q^2}\left(1-\alpha\frac{a^2\sin^2\theta}{r^2+a^2}\right)\alpha^{-2}f(L\psi)^2+\alpha\frac{a^2\sin^2\theta}{q^2}f(\pd_r\psi)^2+\frac{\alpha^{-1}f}{q^2}(\pd_\theta\psi)^2.
\end{align*}
Also,
\begin{align*}
-J^t_{(\psi)}\left[0,\frac{2rf}{r^2+a^2}\right] &= -\frac{2rf}{r^2+a^2}\psi\pd^t\psi \\
&= -\frac{2rf}{r^2+a^2}(g^{tt}+{}^{(Q)}g^{tt})\psi\pd_t\psi \\
&= \frac{2r\alpha^{-1}f}{q^2}\left(1-\alpha\frac{a^2\sin^2\theta}{r^2+a^2}\right)\psi\pd_t\psi \\
&= \frac{2r\alpha^{-1}f}{q^2}\left(1-\alpha\frac{a^2\sin^2\theta}{r^2+a^2}\right)\psi L\psi - \frac{2rf}{q^2}\left(1-\alpha\frac{a^2\sin^2\theta}{r^2+a^2}\right)\psi\pd_r\psi \\
&= \frac{2r\alpha^{-1}f}{q^2}\left(1-\alpha\frac{a^2\sin^2\theta}{r^2+a^2}\right)\psi L\psi -\frac{2rf}{q^2}\psi\pd_r\psi +\frac{2rf}{q^2}\alpha\frac{a^2\sin^2\theta}{r^2+a^2}\psi\pd_r\psi \\
&= \frac{2r\alpha^{-1}f}{q^2}\left(1-\alpha\frac{a^2\sin^2\theta}{r^2+a^2}\right)\psi L\psi +\left(-\frac{1}{q^2}\pd_r(rf\psi^2)+\frac{f+rf'}{q^2}\psi^2\right) \\
&\hspace{3.2in}+\frac{2rf}{q^2}\alpha\frac{a^2\sin^2\theta}{r^2+a^2}\psi\pd_r\psi \\
&= \left(1-\alpha\frac{a^2\sin^2\theta}{r^2+a^2}\right)\left(\frac{2r\alpha^{-1}f}{q^2}\psi L\psi+\frac{f}{q^2}\psi^2\right) \\
&\hspace{0.5in}+\alpha\frac{a^2\sin^2\theta}{r^2+a^2}\left(\frac{2rf}{q^2}\psi\pd_r\psi+\frac{f}{q^2}\psi^2\right) +\frac{rf'}{q^2}\psi^2-\frac{1}{q^2}\pd_r(rf\psi^2).
\end{align*}
Now, observe that
\begin{multline*}
\frac{r^2+a^2}{q^2}\alpha^{-2}f(L\psi)^2+\frac{2r\alpha^{-1}f}{q^2}\psi L\psi+\frac{f}{q^2}\psi^2  \\
= \frac{r^2+a^2}{q^2}f\left(\alpha^{-1}L\psi+\frac{r}{r^2+a^2}\psi\right)^2+\frac{a^2f}{q^2(r^2+a^2)}\psi^2
\end{multline*}
and
$$\frac{r^2+a^2}{q^2}f(\pd_r\psi)^2 +\frac{2rf}{q^2}\psi\pd_r\psi+\frac{f}{q^2}\psi^2 = \frac{r^2+a^2}{q^2}f\left(\pd_r\psi+\frac{r}{r^2+a^2}\psi\right)^2+\frac{a^2f}{q^2(r^2+a^2)}\psi^2.$$
Thus,
\begin{multline*}
-J^t\left[\alpha^{-1}fL,\frac{2rf}{r^2+a^2}\right] \\
=\left(1-\alpha\frac{a^2\sin^2\theta}{r^2+a^2}\right)\frac{r^2+a^2}{q^2}f\left(\alpha^{-1}L\psi+\frac{r}{r^2+a^2}\psi\right)^2 \\
+\alpha\frac{a^2\sin^2\theta}{r^2+a^2}\frac{r^2+a^2}{q^2}f\left(\pd_r\psi+\frac{r}{r^2+a^2}\psi\right)^2 \\
+\frac{a^2f}{q^2(r^2+a^2)}\psi^2+\frac{\alpha^{-1}f}{q^2}(\pd_\theta\psi)^2+\frac{rf'}{q^2}\psi^2-\frac1{q^2}\pd_r(rf\psi^2).
\end{multline*}
Also,
$$-J^t\left[0,0,(1-\epsilon)\frac{rf'}{q^2}L\right] = -(1-\epsilon)\frac{rf'}{q^2}\psi^2 L^t =-(1-\epsilon)\frac{rf'}{q^2}\psi^2.$$
Adding these two expressions together yields the result. 
\end{proof}

\noindent This concludes the proof of Lemma \ref{k:p_ee_identity_lem}.
\end{proof}

\subsection{The incomplete $r^p$ estimate near $i^0$}

Now we use Lemma \ref{k:p_ee_identity_lem} and make a particular choice for the funciton $f$ (so that $f=r^p$ for large $r$) to prove the following.
\begin{proposition}\label{k:incomplete_p_estimate_prop}
Suppose $\pd_\phi\psi=0$. Fix $\delm,\delp>0$. Let $R$ be a sufficiently large radius. Then for all $p\in[\delm,2-\delp]$, the following estimate holds if $\psi$ decays sufficiently fast as $r\rightarrow\infty$.
\begin{multline*}
\int_{\Sigma_{t_2}\cap\{r>2R\}}r^p\left[(L\psi)^2+|\sla\nabla\psi|^2+r^{-2}\psi^2\right] \\
+ \int_{t_1}^{t_2}\int_{\Sigma_t\cap\{r>2R\}}r^{p-1}\left[(L\psi)^2+|\sla\nabla\psi|^2+r^{-2}\psi^2\right] \\
\lesssim \int_{\Sigma_{t_2}\cap\{r>2R\}}r^p\left[(L\psi)^2+|\sla\nabla\psi|^2+r^{-2}\psi^2\right] + Err,
\end{multline*}
where
\begin{align*}
Err &= Err_1 + Err_2 + Err_\Box \\
Err_1 &= \int_{t_1}^{t_2}\int_{\Sigma_t\cap\{R<r<2R\}}(L\psi)^2+|\sla\nabla\psi|^2+\frac{a^2}{M^2}(\pd_t\psi)^2+\psi^2  \\
Err_2 &= \int_{\Sigma_{t_1}\cap\{r>R\}}a^2r^{p-2}(\pd_r\psi)^2 \\
Err_\Box &= \int_{t_1}^{t_2}\int_{\Sigma_t\cap\{R<r\}}r^p(|L\psi|+r^{-1}|\psi|)|\Box_g\psi|.
\end{align*}
\end{proposition}
\begin{proof}
The estimate follows from the identity given in Lemma \ref{k:p_ee_identity_lem} and a particular choice for the function $f$.
$$f(r)=\rho^p,$$
where
$$
\rho = \left\{
\begin{array}{ll}
0 & r\le R \\
smooth & r\in[R,2R] \\
r & 2R< r.
\end{array}
\right.
$$
With this choice, we have 
$$f\ge 0$$
$$f'\ge 0$$
and for $r>2R$,
$$f = r^p$$
$$f'=pr^{p-1}.$$
Furthermore, for $r>2R$,
\begin{multline*}
\frac{2rf}{r^2+a^2}-f' = \frac{2r^{p+1}}{r^2+a^2}-pr^{p-1} = \frac{(2-p)r^{p+1}}{r^2+a^2}-\frac{a^2pr^{p-1}}{r^2+a^2} \\
 \ge \frac{1}{1+a^2/(4R^2)}\left(2-p-\frac{a^2p}{4R^2}\right)r^{p-1}.
\end{multline*}
It follows that if $p\le 2-\delp$ and $R$ is sufficiently large so that $\frac{a^2p}{4R^2}\le \delp/2$, then for $r>2R$,
$$r^{p-1}\lesssim \frac{2rf}{r^2+a^2}-f'.$$
Also, for $r>2R$,
$$\epsilon\alpha ((1-\epsilon)f'-rf'') = \epsilon\alpha ((1-\epsilon)pr^{p-1}-p(p-1)r^{p-1}) = \epsilon\alpha p (2-\epsilon -p)r^{p-1}.$$
If $R$ is sufficiently large so that $\alpha>3/4$ and $\delm\le p\le 2-\delp$ and $\epsilon\le \delp/2$, then
$$\epsilon r^{p-1} \lesssim \epsilon\alpha ((1-\epsilon)f'-rf'').$$

We also note that there are some error terms that either have a factor of $\alpha'$ or $\frac{a^2}{r^2+a^2}$. Each of these terms has a smallness parameter available, since $R$ can be taken to be very large and
$$\alpha'\lesssim \frac{M}{R}r^{-1}$$
and
$$\frac{a^2}{r^2+a^2}\lesssim \frac{M^2}{R^2}.$$

Finally, we observe that if $\psi$ vanishes sufficiently fast as $r\rightarrow\infty$, then since $f$ is supported for $r>R$, we have
$$\int_{\Sigma_t}-\frac1{q^2}\pd_r(rf\psi^2) =0.$$
With these facts having been established, it is straightforward to check that the estimate follows from Lemma \ref{k:p_ee_identity_lem}.

\end{proof}

\subsection{The $r^p$ estimate}

We conclude this section by proving the $r^p$ estimate. This is a combination of the $h\pd_t$ estimate (Proposition \ref{k:hdt_prop}), the Morawetz estimate (Theorem \ref{k:morawetz_thm}), and the incomplete $r^p$ estimate (Proposition \ref{k:incomplete_p_estimate_prop}).
\begin{proposition}\label{k:rp_prop}
Fix $\delm,\delp>0$ and let $p\in[\delm,2-\delp]$. Then if
$$\pd_\phi\psi = 0$$
and $\psi$ decays sufficiently fast as $r\rightarrow\infty$, the following estimate holds.
\begin{multline*}
\int_{\Sigma_{t_2}}r^p\left[(L\psi)^2+|\sla\nabla\psi|^2+r^{-2}\psi^2 + r^{-2}(\pd_r\psi)^2\right] \\
+ \int_{t_1}^{t_2}\int_{\Sigma_t}r^{p-1}\left[\chi_{trap}(L\psi)^2+\chi_{trap}|\sla\nabla\psi|^2+r^{-2}\psi^2+r^{-2}(\pd_r\psi)^2\right] \\
\lesssim \int_{\Sigma_{t_1}}r^p\left[(L\psi)^2+|\sla\nabla\psi|^2+r^{-2}\psi^2+r^{-2}(\pd_r\psi)^2\right] 
+Err,
\end{multline*}
where $\chi_{trap}=\left(1-\frac{r_{trap}}{r}\right)^2$ and
\begin{align*}
Err &= Err_\Box \\
Err_\Box &= \int_{t_1}^{t_2}\int_{\Sigma_t}|(2X(\psi)+w\psi)\Box_g\psi|,
\end{align*}
where the vectorfield $X$ and function $w$ satisfy the following properties. \\
\bp $X$ is everywhere timelike, but asymptotically null at the rate $X=O(r^p)L+O(r^{p-2})\pd_t$. \\
\bp $X|_{r=r_H}=-\lambda\pd_r$ for some positive constant $\lambda$. \\
\bp $X|_{r=r_{trap}}=\lambda\pd_t$ for some positive constant $\lambda$. \\
\bp $w =O(r^{p-1})$ for large $r$.
\end{proposition}
\begin{proof}
We start with the Morawetz estimate (Theorem \ref{k:morawetz_thm}) and add a small constant times the incomplete $r^p$ estimate (Proposition \ref{k:incomplete_p_estimate_prop}). The small constant can be chosen so that the bulk error term $Err_1$ from Proposition \ref{k:incomplete_p_estimate_prop} can be absorbed into the bulk in the Morawetz estimate. The result is the following estimate.
\begin{multline*}
\int_{\Sigma_{t_2}}r^p\left[(L\psi)^2+|\sla\nabla\psi|^2+r^{-2}\psi^2\right]+\frac{M^2}{r^2}(\pd_r\psi)^2 \\
\hspace{1in}+\int_{t_1}^{t_2}\int_{\Sigma_t} r^{p-1}\left[\chi_{trap}(L\psi)^2+\chi_{trap}|\sla\nabla\psi|^2+r^{-2}\psi^2\right]+\frac{M^2}{r^3}(\pd_r\psi)^2 \\
\lesssim \int_{\Sigma_{t_1}}r^p\left[(L\psi)^2+|\sla\nabla\psi|^2+r^{-2}\psi^2\right]+\frac{M^2}{r^2}(\pd_r\psi)^2 + Err'
\end{multline*}
where
\begin{align*}
Err' &= Err'_1+Err'_2+Err'_3+Err'_\Box \\
Err'_1 &= \int_{\Sigma_{t_2}}r^{-1}|\psi L\psi| \\
Err'_2 &= \int_{H_{t_1}^{t_2}}(\pd_t\psi)^2 +\int_{\Sigma_{t_2}}\frac{M^2}{r^2}\left[\chi_H(\pd_r\psi)^2+(\pd_t\psi)^2\right] \\
Err'_3 &= \int_{\Sigma_{t_1}\cap\{r>R\}}a^2r^{p-2}(\pd_r\psi)^2 \\
Err'_\Box &= \int_{t_1}^{t_2}\int_{\Sigma_t\cap\{R<r\}}r^p(|L\psi|+r^{-1}|\psi|)|\Box_g\psi| \\
&\hspace{1in} +\int_{t_1}^{t_2}\int_{\Sigma_t}|(2X'(\psi)+w'\psi)|\Box_{g}\psi|,
\end{align*}
and $X'$ and $w'$ are the vectorfield and function defined in the Morawetz estimate (Theorem \ref{k:morawetz_thm}).

The error term $Err'_1$ can in fact be removed due to the following argument.
\begin{align*}
Err'_1 &\lesssim \int_{\Sigma_{t_2}}\epsilon r^p(L\psi)^2+\epsilon^{-1}r^{-p}r^{-2}\psi^2 \\
&\lesssim \int_{\Sigma_{t_2}}\epsilon r^p[(L\psi)^2+r^{-2}\psi^2] +\int_{\Sigma_{t_2}\cap\{r\le R_\epsilon\}}\epsilon^{-1}\psi^2 .
\end{align*}
The radius $R_\epsilon$ should be chosen sufficiently large so that $\epsilon^{-1}r^{-p}\le \epsilon r^p$ whenever $r>R_\epsilon$. This critically depends on the fact that $p\ge\delm>0$. Now, the parameter $\epsilon$ can be taken sufficiently small so as to absorb the first two terms into the left side of the main estimate and the last two terms can be included with the term $Err'_2$ after applying a Hardy estimate.

We return to the main estimate. Notice that most terms have improved weights near $i^0$ and a few error terms remain on $H_{t_1}^{t_2}$ and $\Sigma_{t_2}$. The next step is to use the $h\pd_t$ estimate (Proposition \ref{k:hdt_prop}) to eliminate these error terms and improve the weights near $i^0$ for the remaining $\pd_r\psi$ terms. The result is the following estimate.
\begin{multline*}
\int_{\Sigma_{t_2}}r^p\left[(L\psi)^2+|\sla\nabla\psi|^2+r^{-2}\psi^2+r^{-2}(\pd_r\psi)^2\right] \\
+\int_{t_1}^{t_2}\int_{\Sigma_t} r^{p-1}\left[\chi_{trap}(L\psi)^2+\chi_{trap}|\sla\nabla\psi|^2+r^{-2}\psi^2+c_\epsilon r^{-2}(\pd_r\psi)^2\right] \\
\lesssim \int_{\Sigma_{t_1}}r^p\left[(L\psi)^2+|\sla\nabla\psi|^2+r^{-2}\psi^2+r^2(\pd_r\psi)^2\right] +Err'',
\end{multline*}
where
\begin{align*}
Err'' &= Err''_1+Err''_\Box \\
Err''_1 &= \int_{t_1}^{t_2}\int_{\Sigma_t\cap\{5M<r\}} \epsilon r^{-1}(L\psi)^2 \\
Err''_\Box &= \int_{t_1}^{t_2}\int_{\Sigma_t\cap\{R<r\}}r^p(|L\psi|+r^{-1}|\psi|)|\Box_g\psi| + \int_{t_1}^{t_2}\int_{\Sigma_t}|(2X'(\psi)+w'\psi)(\Box_g\psi)| \\
&\hspace{3.5in}+ \int_{t_1}^{t_2}\int_{\Sigma_t}C_\epsilon r^{p-2}|\pd_t\psi(\Box_g\psi)|.
\end{align*}
By taking $\epsilon$ sufficiently small, the error term $Err''_1$ can be absorbed into the left side. 

The resulting vectorfield $X$ and function $w$ can be understood by combining the vectorfields and scalar functions used to construct the three estimates that were used in this proof.
\end{proof}

\section{The dynamic estimates}\label{k:dynamic_sec}

In this section, we prove the dynamic estimates, which provide all of the necessary information related to the future dynamics of the wavefunction $\psi$. These estimates are simply restatements of the energy estimate (Proposition \ref{k:classic_ee_prop}) and the $r^p$ estimate (Proposition \ref{k:rp_prop}) applied to either $\psi$ itself or another wavefunction ($\psi^s$ or $\psi^{s,k}$) derived from $\psi$. The $k$ index is new in this problem and corresponds to the operator $\tg$ (see \S\ref{k:additional_commutator_sec}).

\subsection{The dynamic estimates for $\psi$ ($s,k=0$)}

We begin with the dynamic estimates for the wavefunction $\psi$ only.
\begin{proposition}\label{k:dynamic_estimates_0_prop}
Suppose $\pd_\phi\psi=0$.
Fix $\delm,\delp>0$ and let $p\in[\delm,2-\delp]$. Then
$$E(t_2)\lesssim E(t_1)+\int_{t_1}^{t_2}N(t)dt$$
$$E_p(t_2)+\int_{t_1}^{t_2}B_p(t)dt \lesssim E_p(t_1)+\int_{t_1}^{t_2}N_p(t)dt,$$
where
$$E(t)=\int_{\Sigma_t}\chi_H(\pd_r\psi)^2+(\pd_t\psi)^2+|\sla\nabla\psi|^2,$$
$$E_p(t)=\int_{\Sigma_t}r^p\left[(L\psi)^2+|\sla\nabla\psi|^2+r^{-2}\psi^2+r^{-2}(\pd_r\psi)^2\right],$$
$$B_p(t)=\int_{\Sigma_t}r^{p-1}\left[\chi_{trap}(L\psi)^2+\chi_{trap}|\sla\nabla\psi|^2+r^{-2}\psi^2+r^{-2}(\pd_r\psi)^2\right],$$
$$N(t)=\int_{\Sigma_t}|\pd_t\psi\Box_g\psi|,$$
$$N_p(t)=\int_{\Sigma_t}|(2X(\psi)+w\psi)\Box_g\psi|,$$
where $\chi_H=1-\frac{r_H}r$, $\chi_{trap}=\left(1-\frac{r_{trap}\}}r\right)^2$, and $X$ and $w$ are the vectorfield and function from Proposition \ref{k:rp_prop}.
\end{proposition}
\begin{proof}
The first estimate is a restatement of the energy estimate (Proposition \ref{k:classic_ee_prop}) and the second is a restatement of the $r^p$ estimate (Proposition \ref{k:rp_prop}).
\end{proof}

We slightly simplify the previous estimates by absorbing part of the nonlinear norm $N_p(t)$ into the bulk $B_p(t)$ on the left side. Note the complication due to trapping.
\begin{corollary}\label{k:dynamic_estimates_0_cor}
Suppose $\pd_\phi\psi=0$.
Fix $\delm,\delp>0$ and let $p\in[\delm,2-\delp]$. Then
$$E(t_2)\lesssim E(t_1)+\int_{t_1}^{t_2}N(t)dt$$
$$E_p(t_2)+\int_{t_1}^{t_2}B_p(t)dt \lesssim E_p(t_1)+\int_{t_1}^{t_2}N_p(t)dt,$$
where $E(t)$, $E_p(t)$, $B_p(t)$, are as defined in Proposition \ref{k:dynamic_estimates_0_prop} and
$$N(t)=(E(t))^{1/2}||\Box_g\psi||_{L^2(\Sigma_t)},$$
$$N_p(t)=\int_{\Sigma_t}r^{p+1}(\Box_g\psi)^2+\int_{\Sigma_t\cap\{r\approx r_{trap}\}}|\pd_t\psi\Box_g\psi|.$$
\end{corollary}
The following proof, which is essentially the same as the proof of Corollary \ref{s:dynamic_estimates_0_cor}, is provided again for comparison with the proof of Proposition \ref{k:s_k_prop}.
\begin{proof}
The first estimate is due to Proposition \ref{k:dynamic_estimates_0_prop} and the fact that
$$\int_{\Sigma_t}|\pd_t\psi\Box_g\psi|\lesssim ||\pd_t\psi||_{L^2(\Sigma_t)}||\Box_g\psi||_{L^2(\Sigma_t)} \lesssim (E(t))^{1/2}||\Box_g\psi||_{L^2(\Sigma_t)}.$$
We turn to the second estimate.
According to Proposition \ref{k:dynamic_estimates_0_prop},
$$E_p(t_2)+\int_{t_1}^{t_2}B_p(t)dt \lesssim E_p(t_1)+\int_{t_1}^{t_2}N_p'(t)dt,$$
where
\begin{align*}
N_p'(t) &= \int_{\Sigma_t}|(2X(\psi)+w\psi)\Box_g\psi| \\
&\lesssim \int_{\Sigma_t}r^p(|L\psi|+r^{-1}|\psi|+r^{-2}|\pd_r\psi|)|\Box_g\psi|.
\end{align*}
We would like to claim that
$$N_p'(t)\lesssim \epsilon B_p(t)+\epsilon^{-1}\int_{\Sigma_t}r^{p+1}(\Box_g\psi)^2.$$
\textbf{This would require an estimate of the form}
$$\int_{\Sigma_t}r^{p-1}(|L\psi|+r^{-1}|\psi|+r^{-2}|\pd_r\psi|)^2\lesssim B_p(t),$$
\textbf{which is not completely true due to the loss of $(L\psi)^2$ at the trapping radius}. Instead,
$$\int_{\Sigma_t}r^{p-1}(|L\psi|+r^{-1}|\psi|+r^{-2}|\pd_r\psi|)^2\lesssim B_p(t)+\int_{\Sigma_t\cap\{r\approx r_{trap}\}}|\pd_t\psi\Box_g\psi|,$$
and therefore
$$N_p'(t)\lesssim \epsilon B_p(t)+\epsilon^{-1}\int_{\Sigma_t}r^{p+1}(\Box_g\psi)^2+\epsilon\int_{\Sigma_t\cap\{r\approx r_{trap}\}}|\pd_t\psi\Box_g\psi|.$$
By taking $\epsilon$ sufficiently small, the term $\epsilon \int_{t_1}^{t_2}B_p(t)dt$ can be absorbed into the left side of the second estimate.
\end{proof}

\subsection{Commutators with $\Box_g$}\label{k:symmetry_operators_sec}

We will use the operators $\pd_t$ and $Q$ as commutators. We refer to these commutators as symmetry operators, because of the following lemma.
\begin{lemma}
$$[\pd_t,q^2\Box_g]=0,$$
$$[Q,q^2\Box_g]=0.$$
\end{lemma}
\begin{proof}
There is no explicit $t$ dependence in the metric, so $\pd_t$ clearly commutes with $q^2\Box_g$.

Next, observe that
\begin{align*}
q^2\Box_g\psi&=\frac1{\sin\theta}\pd_\alpha\left(\sin\theta q^2 g^{\alpha\beta}\pd_\beta\psi\right) \\
&=-\frac{(r^2+a^2)^2}{\Delta}\pd_t^2\psi+\pd_r(\Delta\pd_r\psi)+Q\psi,
\end{align*}
and since
$$Q\psi = \frac1{\sin\theta}\pd_\alpha\left(\sin\theta Q^{\alpha\beta}\pd_\beta\psi\right),$$
it commutes with the operators $-\frac{(r^2+a^2)^2}{\Delta}\pd_t^2$ and $\pd_r(\Delta\pd_r\cdot)$.
\end{proof}

We define the $s$-order commutators $\Gamma$.
\begin{definition}
$$\Gamma^su = Q^l\pd_t^{s-2l}u$$
where $0\le 2l\le s$.
\end{definition}
We also define the $s$-order wavefunctions derived from the wavefunction $\psi$ by applying symmetry operators.
\begin{definition}
$$\psi^s=\Gamma^s\psi$$
\end{definition}
Thus, we have the following lemma.
\begin{lemma}\label{k:s_comm_lem}
If $\Box_g\psi=0$, then
$$\Box_g\psi^s=0.$$
More generally,
$$\Box_g\psi^s = q^{-2}\Gamma^s(q^2\Box_g\psi).$$
\end{lemma}
\begin{proof}
We only prove the general case.
$$\Box_g\psi^s=\Box_g\Gamma^s\psi=q^{-2}q^2\Box_g\Gamma^s\psi=q^{-2}\Gamma^s(q^2\Box_g\psi).$$
\end{proof}

\subsection{The dynamic estimates for $\psi^s$ ($k=0$)}

Now we can write the dynamic estimates in a more general form.
\begin{corollary}\label{k:s_cor}
Suppose $\pd_\phi\psi=0$.
Fix $\delm,\delp>0$. The following estimates hold for $p\in [\delm,2-\delp]$ and every integer $s\ge 0$.
$$E^s(t_2)\lesssim E^s(t_1)+\int_{t_1}^{t_2}N^s(t)dt,$$
$$E_p^s(t_2)+\int_{t_1}^{t_2}B_p^s(t)dt\lesssim E_p^s(t_1)+\int_{t_1}^{t_2}N_p^s(t)dt,$$
where
$$E^s(t)=\sum_{s'\le s} E[\psi^{s'}](t),$$
$$E_p^s(t)=\sum_{s'\le s} E_p[\psi^{s'}](t),$$
$$B_p^s(t)=\sum_{s'\le s} B_p[\psi^{s'}](t),$$
$$N^s(t)= (E^s(t))^{1/2}\sum_{s'\le s}||q^{-2}\Gamma^{s'}(q^2\Box_g\psi)||_{L^2(\Sigma_t)},$$
$$N_p^s(t)=\sum_{s'\le s} \int_{\Sigma_t}r^{p+1}(q^{-2}\Gamma^{\le s}(q^2\Box_g\psi))^2+\sum_{s'\le s}\int_{\Sigma_t\cap\{r\approx r_{trap}\}}|\pd_t\psi^{s'}q^{-2}\Gamma^{s'}(q^2\Box_g\psi)|,$$
and the norms $E(t)$, $E_p(t)$, and $B_p(t)$ are as defined in Proposition \ref{k:dynamic_estimates_0_prop}.
\end{corollary}
\begin{proof}
The proof is a direct application of Corollary \ref{k:dynamic_estimates_0_cor} by making the substitutions
$$\psi\mapsto\psi^{s'}$$
for all values of $s'$ (and all commutators represented by $\Gamma^{s'}$) where $s'\le s$ and observing Lemma \ref{k:s_comm_lem}.
\end{proof}

\subsection{The additional operator $\tg$ as a commutator}\label{k:additional_commutator_sec}

To handle nonlinear terms with a factor of $\pd_r\psi$ near the event horizon, it will be necessary to use an additional commutator, which we call $\tg$.
\begin{definition}
$$\tg=1_H(r)\pd_r,$$
where $1_H$ is a nonnegative smooth function supported on a neighborhood of the event horizon.
\end{definition}

Unfortunately, $\tg$ does not commute with the wave operator $\Box_g$, but we will see that one particular term in the commutator with $\Box_g$ has an appropriate sign near the event horizon, and that is what allows us to use $\tg$. We now compute the commutator.

\begin{lemma}\label{k:dr_comm_1_lem}
$$[\tg,q^2\Box_g]u=\Delta'\pd_r\tg u -21_H'q^2\Box_gu+\{\pd_r\Gamma^{\le 1}u,\Gamma^{\le 2}u\} $$
where the expression $\{\pd_r\Gamma^{\le 1}u,\Gamma^{\le 2}u\}$ represents any terms of the form $f(r,\theta)\pd_r\Gamma^{\le 1}u$ or $f(r,\theta)\Gamma^{\le 2}u$ for smooth $f$ with compact support.
\end{lemma}
\begin{proof}
We expand
$$q^2\Box_g= q^2g^{tt}\pd_t^2 + 2q^2g^{tr}\pd_r\pd_t +q^2g^{rr}\pd_r^2+\pd_r(q^2g^{tr})\pd_t+\pd_r(q^2g^{rr})\pd_r+\frac1{\sin\theta}\pd_\alpha\left(\sin\theta Q^{\alpha\beta}\pd_\beta\cdot\right).$$
Note that the $Q$ term commutes with $\tg$. We compute terms arising in the commutator by neglecting the highest order derivative terms and the $Q$ terms (both represented by the ellipsis in what follows).
\begin{multline*}
\tg (q^2\Box_g u)= \\
\tg(q^2g^{tt})\pd_t^2u+2\tg(q^2g^{tr})\pd_r\pd_tu+\tg(q^2g^{rr})\pd_r^2u+\tg(\pd_r(q^2g^{tr}))\pd_tu+\tg(\pd_r(q^2g^{rr}))\pd_ru+... \\
=1_H\Delta'\pd_r^2u+\{\pd_t^2u,\pd_r\pd_tu,\pd_tu,\pd_ru\}+...
\end{multline*}
\begin{align*}
q^2\Box_g(\tg u)&=2q^2g^{tr}1_H'\pd_r\pd_tu+2q^2g^{rr}1_H'\pd_r^2u+q^2g^{rr}1_H''\pd_ru+\pd_r(q^2g^{rr})1_H'\pd_ru+... \\
&=2q^2g^{rr}1_H' \pd_r^2u+\{\pd_r\pd_tu,\pd_ru\}+... \\
&=21_H' q^2\Box_g u+\{\pd_t^2u,\pd_r\pd_tu,Qu,\pd_tu,\pd_ru\}+...
\end{align*}
Taking the difference,
\begin{align*}
[\tg,q^2\Box_g]u &= \tg(q^2\Box_g u)-q^2\Box_g(\tg u) \\
&=1_H\Delta'\pd_r^2u - 21_H' q^2\Box_g u+\{\pd_t^2u,\pd_r\pd_tu,Qu,\pd_tu,\pd_ru\} \\
&=\Delta'\pd_r\tg u - 21_H'q^2\Box_g u+\{\pd_t^2u,\pd_r\pd_tu,Qu,\pd_tu,\pd_ru\}.
\end{align*}
\end{proof}

Now we generalize the prevous lemma by commuting with $\tg$ arbitrarily many times.
\begin{lemma}\label{k:dr_comm_2_lem}
$$[\tg^k,q^2\Box_g]u=k\Delta'\pd_r\tg^ku+\{\pd_r\tg^{\le k-1}\Gamma^{\le 1}u,\tg^{\le k-1}\Gamma^{\le 2}u,\tg^{\le k-1}(q^2\Box_gu)\},$$
where the $\{...\}$ notation is the same as in Lemma \ref{k:dr_comm_1_lem}.
\end{lemma}
\begin{proof}
The proof is an induction argument on $k$. The case $k=1$ corresponds to Lemma \ref{k:dr_comm_1_lem}. Assuming that the statement of the lemma holds at the level $k$, we prove the analogous statement at the level $k+1$. We have
\begin{align*}
[\tg^{k+1},q^2\Box_g]u &= \tg^{k+1}(q^2\Box_gu)-q^2\Box_g(\tg^{k+1}u) \\
&= \tg (\tg^k(q^2\Box_gu)-q^2\Box_g(\tg^k u))+(\tg(q^2\Box_g(\tg^ku))-q^2\Box_g(\tg^{k+1}u)) \\
&= \tg[\tg^k,q^2\Box_g]u+[\tg,q^2\Box_g](\tg^k u).
\end{align*}
Now, by the inductive hypothesis,
\begin{align*}
\tg[\tg^k,q^2\Box_g]u &= \tg\left(k\Delta'\pd_r\tg^ku+\{\pd_r\tg^{\le k-1}\Gamma^{\le 1}u,\tg^{\le k-1}\Gamma^{\le 2}u,\tg^{\le k-1}(q^2\Box_gu)\}\right) \\
&= k\Delta'\pd_r\tg^{k+1}u +\{\pd_r\tg^{\le k}\Gamma^{\le 1}u,\tg^{\le k}\Gamma^{\le 2}u,\tg^{\le k}(q^2\Box_gu)\},
\end{align*}
and by the base case,
\begin{align*}
[\tg,q^2\Box_g](\tg^k u) &= \Delta'\pd_r\tg^{k+1}u +\{\pd_r\Gamma^{\le 1}\tg^ku,\Gamma^{\le 2}\tg^ku,q^2\Box_g(\tg^ku)\} \\
&= \Delta'\pd_r\tg^{k+1}u +\{\pd_r\tg^k\Gamma^{\le 1}u,\tg^k\Gamma^{\le 2}u,\tg^{\le k}(q^2\Box_gu)\},
\end{align*}
where in the last step we used the inductive hypothesis a second time.

Summing these two yields
$$[\tg^{k+1},q^2\Box_g]u = (k+1)\Delta'\pd_r\tg^{k+1}u+ \{\pd_r\tg^{\le k}\Gamma^{\le 1}u,\tg^{\le k}\Gamma^{\le 2}u,\tg^{\le k}(q^2\Box_gu)\}.$$
This completes the inductive argument.
\end{proof}

Since $\tg$ does not commute with $\Box_g$, it is treated seperately from the previous commutators. We define the $s,k$-order wavefunctions.
\begin{definition}
$$\psi^{s,k} = \tg^k\Gamma^s\psi.$$
\end{definition}

We have the following generalization of Lemma \ref{k:s_comm_lem}.
\begin{lemma}\label{k:s_k_comm_lem}
$$\Box_g\psi^{s,k}=q^{-2}[q^2\Box_g,\tg^k]\psi^s+q^{-2}\tg^k\Gamma^s(q^2\Box_g\psi).$$
\end{lemma}
\begin{proof}
\begin{align*}
\Box_g\psi^{s,k} &= (\Box_g\psi^{s,k}-q^{-2}\tg^k(q^2\Box_g\psi^s))+q^{-2}\tg^k(q^2\Box_g\psi^s) \\
&=q^{-2}[q^2\Box_g,\tg^k]\psi^s+q^{-2}\tg^k\Gamma^s(q^2\Box_g\psi).
\end{align*}
\end{proof}

\subsection{The dynamic estimates for $\psi^{s,k}$}\label{k:dynamic_estimates_s_k_sec}

Now we have the following generalization of Proposition \ref{k:dynamic_estimates_0_prop}.
\begin{proposition}\label{k:s_k_L_prop}
Suppose $\pd_\phi\psi=0$.
Fix $\delm,\delp>0$. The following estimates hold for $p\in [\delm,2-\delp]$ and for integers $s\ge 0$ and $k\ge 1$.
$$E^{s,k}(t_2)\lesssim E^{s,k}(t_1)+\int_{t_1}^{t_2}L^{s,k}(t)+N^{s,k}(t)dt,$$
$$E_p^{s,k}(t_2)+\int_{t_1}^{t_2}B_p^{s,k}(t)dt\lesssim E_p^{s,k}(t_1)+\int_{t_1}^{t_2}L_p^{s,k}(t)+N_p^{s,k}(t)dt,$$
where
\begin{align*}
E^{s,k}(t) &= \sum_{\substack{s'\le s \\ k'\le k}} E[\psi^{s',k'}](t), \\
E_p^{s,k}(t) &= \sum_{\substack{s'\le s \\ k'\le k}} E_p[\psi^{s',k'}](t), \\
B_p^{s,k}(t) &= \sum_{\substack{s'\le s \\ k'\le k}} B_p[\psi^{s',k'}](t),
\end{align*}
and
$$L^{s,k}(t) =\sum_{\substack{s'\le s \\ k'\le k}}\int_{\Sigma_t}|\pd_t\psi^{s',k'}(q^{-2}[q^2\Box_g,\tg^{k'}]\psi^{s'})|$$
$$L_p^{s,k}(t) = \sum_{\substack{s'\le s \\ k'\le k}}\int_{\Sigma_t}|(2X(\psi^{s',k'})+w\psi^{s',k'})(q^{-2}[q^2\Box_g,\tg^{k'}]\psi^{s'})|$$
$$N^{s,k}(t) =\sum_{\substack{s'\le s \\ k'\le k}}\int_{\Sigma_t}|\pd_t\psi^{s',k'}q^{-2}\tg^{k'}\Gamma^{s'}(q^2\Box_g\psi)|$$
$$N_p^{s,k}(t) = \sum_{\substack{s'\le s \\ k'\le k}}\int_{\Sigma_t}|(2X(\psi^{s',k'})+w\psi^{s',k'})q^{-2}\tg^{k'}\Gamma^{s'}(q^2\Box_g\psi)|$$
and the norms $E(t)$, $E_p(t)$, and $B_p(t)$ are as defined in Proposition \ref{k:dynamic_estimates_0_prop}.
\end{proposition}
\begin{proof}
The proof is a direct application of Proposition \ref{k:dynamic_estimates_0_prop} by making the substitutions
$$\psi\mapsto\psi^{s',k'}$$
for all values of $s'$ and $k'$ (and all commutators represented by $\tg^k\Gamma^{s'}$) where $s'\le s$ and $k'\le k$. By Lemma \ref{k:s_k_comm_lem},
$$|\Box_g\psi^{s',k'}|\le |q^{-2}[q^2\Box_g,\tg^{k'}]\psi^{s'}|+|q^{-2}\tg^{k'}\Gamma^{s'}(q^2\Box_g\psi)|.$$
The resulting error terms have been grouped into the parts that will either be linear or nonlinear when using the nonlinear wave equation.
\end{proof}

Our final task is to simplify Proposition \ref{k:s_k_L_prop} slightly. We begin with the following lemma.
\begin{lemma}\label{k:s_k_L_lem}
Suppose $|a|<M$ and $k\ge 1$. Then the quantities $L^{s,k}(t)$ and $L_p^{s,k}(t)$ defined in Proposition \ref{k:s_k_L_prop} with the absolute values removed satisfy
$$L^{s,k}(t)\lesssim B_{p'}^{s,k}(t)+B_{p'}^{s+2,k-1}(t),$$
$$L_p^{s,k}(t)\lesssim B_{p'}^{s+2,k-1}(t),$$
where $p'$ is arbitrary.
\end{lemma}
\begin{remark}
The reason for the arbitrary $p'$ on the right side is that the norms $L^{s,k}(t)$ and $L_p^{s,k}(t)$ have terms related to commuting with $\tg$. These terms are all supported on a compact radial interval, so the factor $r^{p'-1}$ that appears in $B_{p'}^{s,k}(t)$ can be approximated by a constant.
\end{remark}
\begin{proof}
The first estimate follows directly from Lemma \ref{k:dr_comm_2_lem}, but the second is more interesting. The key observation is to recognize that the term represented by $L_p^{s,k}(t)$ actually has a good sign near the event horizon. That is, according to Lemma \ref{k:dr_comm_2_lem},
\begin{align*}
\int_{\Sigma_t}-X(\tg^k\psi^s)q^{-2}[q^2\Box_g,\tg^k]\psi^s &= \int_{\Sigma_t}-2(-X^r)\pd_r(\tg^k\psi^s)q^{-2}k\Delta'\pd_r(\tg^k\psi^s)+err \\
&=\int_{\Sigma_t}-2(-X^r)q^2\Delta'(\pd_r\tg^k\psi^s)^2+err
\end{align*}
Since $X^r<0$ near the event horizon and since $\Delta'>0$ as long as $|a|<M$, \textbf{the principal term becomes minus a square,} so it can be ignored or used to control small error terms. After this observation, the estimate follows easily.
\end{proof}

We conclude with the following proposition, which applies whenever $k\ge 1$.
\begin{proposition}\label{k:s_k_prop}
Suppose $\pd_\phi\psi=0$.
Fix $\delm,\delp>0$. The following estimates hold for $p\in [\delm,2-\delp]$, $s\ge 0$, $k\ge 1$, and arbitrary $p'$.
$$E^{s,k}(t_2)\lesssim E^{s,k}(t_1)+\int_{t_1}^{t_2}B_{p'}^{s,k}(t)+B_{p'}^{s+2,k-1}(t)+N^{s,k}(t)dt,$$
$$E_p^{s,k}(t_2)+\int_{t_1}^{t_2}B_p^{s,k}(t)dt\lesssim E_p^{s,k}(t_1)+\int_{t_1}^{t_2}B_{p'}^{s+2,k-1}(t)+N_p^{s,k}(t)dt,$$
where
$$N^{s,k}(t)= (E^{s,k}(t))^{1/2}\sum_{\substack{s'\le s \\ k'\le k}}||q^{-2}\tg^{k'}\Gamma^{s'}(q^2\Box_g\psi)||_{L^2(\Sigma_t)},$$
$$N_p^{s,k}(t) = \sum_{\substack{s'\le s \\ k'\le k}}\int_{\Sigma_t}r^{p+1}|q^{-2}\tg^{k'}\Gamma^{s'}(q^2\Box_g\psi)|^2,$$
and the norms $E^{s,k}(t)$, $E^{s,k}_p(t)$, and $B^{s,k}_p(t)$ are as defined in Proposition \ref{k:s_k_L_prop}.
\end{proposition}

\begin{proof}
By Proposition \ref{k:s_k_L_prop} and Lemma \ref{k:s_k_L_lem}, we have
$$E^{s,k}(t_2)\lesssim E^{s,k}(t_1)+\int_{t_1}^{t_2}B_{p'}^{s,k}(t)+B_{p'}^{s+2,k-1}(t)+(N')^{s,k}(t)dt,$$
where $(N')^{s,k}(t)$ is the quantity $N^{s,k}(t)$ from Proposition \ref{k:s_k_L_prop}. By the same simple argument as in the proof of Corollary \ref{k:dynamic_estimates_0_cor},
$$(N')^{s,k}(t)\lesssim N^{s,k}(t).$$

We turn to the second estimate. By Proposition \ref{k:s_k_L_prop} and Lemma \ref{k:s_k_L_lem}, we have
$$E_p^{s,k}(t_2)+\int_{t_1}^{t_2}B_p^{s,k}(t)dt\lesssim E_p^{s,k}(t_1)+\int_{t_1}^{t_2}B_{p'}^{s+2,k-1}(t)+(N')_p^{s,k}(t)dt,$$
where $(N')_p^{s,k}(t)$ is the quantity $N_p^{s,k}(t)$ from Proposition \ref{k:s_k_L_prop}.
\begin{align*}
(N')_p^{s,k}(t) &= \sum_{\substack{s'\le s \\ k'\le k}}\int_{\Sigma_t}|(2X(\psi^{s',k'})+w\psi^{s',k'})q^{-2}\tg^{k'}\Gamma^{s'}(q^2\Box_g\psi)| \\
&\lesssim \sum_{\substack{s'\le s \\ k'\le k}}\int_{\Sigma_t}r^p(|L\psi^{s',k'}|+r^{-1}|\psi^{s',k'}|+r^{-2}|\pd_r\psi^{s',k'}|)|q^{-2}\tg^{k'}\Gamma^{s'}(q^2\Box_g\psi)| \\
&\lesssim \epsilon B_p^{s,k}(t)+ \sum_{\substack{s'\le s \\ k'\le k}}\int_{\Sigma_t\cap\{r\approx r_{trap}\}}|\pd_t\psi^{s',k'}q^{-2}\tg^{k'}\Gamma^{s'}(q^2\Box_g\psi)| \\
&\hspace{2.3in} + \epsilon^{-1}\sum_{\substack{s'\le s \\ k'\le k}}\int_{\Sigma_t}r^{p+1}|q^{-2}\tg^{k'}\Gamma^{s'}(q^2\Box_g\psi)|^2 \\
&\lesssim \epsilon B_p^{s,k}(t)+\epsilon B_{p'}^{s+2,0}(t) + \epsilon^{-1}\sum_{\substack{s'\le s \\ k'\le k}}\int_{\Sigma_t}r^{p+1}|q^{-2}\tg^{k'}\Gamma^{s'}(q^2\Box_g\psi)|^2.
\end{align*}
The above procedure should be compared to the procedure in the proof of Corollary \ref{k:dynamic_estimates_0_cor}. This time, the problematic term near the trapping region can be controlled in part by $B_{p'}^{s+2,0}(t)$, because $\tg$ is supported away from the trapping region. It doesn't hurt to use $B_{p'}^{s+2,0}(t)$, because the term $B_{p'}^{s+2,k-1}(t)$ already appears on the right side of the estimate.

Just as in the proof of Corollary \ref{k:dynamic_estimates_0_cor}, the term $\epsilon B_p^{s,k}(t)$ can be absorbed into the term on the left side. This proves the second estimate.
\end{proof}

\section{The $L^\infty$ estimates}\label{k:pointwise_sec}

In this section, we prove $L^\infty$ estimates for certain derivatives of $\psi$ that will appear in $q^{-2}\tg^k\Gamma^s(q^2\Box_g\psi)$ after using the equation to replace $\Box_g\psi$ with nonlinear terms. Since the commutator $Q$ is a second order differential operator, in the nonlinear problem there will be terms that are not products of derivatives of $\psi^{s,k}$ (see Proposition \ref{k:commutator_prop}). Fortunately, the ellipticity of $Q$ (Lemma \ref{main_elliptic_lemma}) allows us to still control these terms. This fact is encoded in Lemma \ref{k:pointwise_lem}.

First, we introduce the following notation.
\begin{definition}
We define the expression $\Omega^l\psi^{s-l}$ to represent the following coordinate-invariant $r^l$-weighted rank $l$ tensor on $S^2(r)$.
$$\Omega^l\psi^{s-l}=r^l\sla\nabla^l\psi^{s-l}.$$
(It should always be assumed that $l\le s$.)
\end{definition}
Note that the expression $\Omega\psi$ now has a different meaning than in Chapters \ref{mink_chap} and \ref{szd_chap}, but the meanings are still very similar, because the rotation operators used in Chapters \ref{mink_chap} and \ref{szd_chap} behave like $r\sla\nabla$.

\subsection{A sobolev-type estimate}

The following lemma generalizes Lemma \ref{m:sobolev_I_lem}.
\begin{lemma}\label{k:pointwise_lem}
If $u$ decays sufficiently fast as $r\rightarrow\infty$, then
$$||\Omega^l u||^2_{L^\infty(\Sigma_t\cap\{r\ge r_0\})}\lesssim \int_{\Sigma_t\cap\{r\ge r_0\}}r^{-2}\left[(\pd_r \Gamma^{\le l+3}u)^2+(\Gamma^{\le l+3}u)^2\right]$$
If $l$ is even, the same result holds with only $\Gamma^{\le l+2}u$ in the integral on the right side.
\end{lemma}

\begin{proof}
For a fixed $r$, denote by $\bar{u}:S^2(1)\rightarrow\R{}$ the pullback of the function $u:S^2(r)\rightarrow\R{}$ via the canonical map from $S^2(1)$ to $S^2(r)$. Likewise, let $\overline{\Omega^lu}$ be the pullback of the tensor $\Omega^lu$. Denote by $d\omega$ the measure on $S^2(1)$.

\textbf{Let $l$ be even.} By induction, it is straightforward to show that
$$\overline{\Omega^l u}=\sla\nabla^l\bar{u}.$$
Then
\begin{multline*}
||\Omega^lu||^2_{L^\infty(S^2(r))}=||\overline{\Omega^l u}||_{L^\infty(S^2(1))}^2=||\sla\nabla^l\bar{u}||_{L^\infty(S^2(1))}^2 \\
\lesssim \int_{S^2(1)}(\sla\nabla^{\le l+2}\bar{u})^2=\int_{S^2(1)}\left(\overline{\Omega^{\le l+2}u}\right)^2=\int_{S^2(r)}\left(\Omega^{\le l+2}u\right)^2d\omega \lesssim \int_{S^2(r)}\left(\Gamma^{\le l+2}u\right)^2d\omega.
\end{multline*}
The final step was an application of Lemma \ref{main_elliptic_lemma}. 

\textbf{If instead $l$ is odd}, the calculation is the exact same as for the even case, except in the application of Lemma \ref{main_elliptic_lemma} in the final step. That is,
$$||\Omega^lu||_{L^\infty(S^2(r))}^2\lesssim ... \lesssim \int_{S^2(r)}(\Omega^{\le l+2}u)^2\lesssim \int_{S^2(r)}\left(\Gamma^{\le l+3}u\right)^2d\omega.$$

\textbf{Thus, in both cases},
$$||\Omega^lu||_{L^\infty(S^2(r))}^2\lesssim \int_{S^2(r)}\left(\Gamma^{\le l+3}u\right)^2d\omega.$$
Now, set $f(r)=\int_{S^2(r)}(\Gamma^{\le l+3}u)^2d\omega$. Note that
\begin{multline*}
|f'(r)|\lesssim \int_{S^2(r)}|\Gamma^{\le l+3}u\pd_r\Gamma^{\le l+3}u|d\omega\lesssim \int_{S^2(r)}\left[(\pd_r \Gamma^{\le l+3}u)^2+(\Gamma^{\le l+3}u)^2\right]d\omega \\
\lesssim\int_{S^2(r)}r^{-2}\left[(\pd_r \Gamma^{\le l+3}u)^2+(\Gamma^{\le l+3}u)^2\right]q^2d\omega
\end{multline*}
Then, assuming $\lim_{r\rightarrow\infty}f(r)=0$,
\begin{multline*}
|\Omega^l u(r_0)|^2\lesssim f(r_0)\lesssim \int_{r_0}^\infty |f'(r)|dr\lesssim  \int_{r_0}^\infty\int_{S_2(r)}r^{-2}\left[(\pd_r\Gamma^{\le 3}u)^2+(\Gamma^{\le 3}u)^2\right]q^2d\omega dr \\
=\int_{\Sigma_t\cap\{r\ge r_0\}}r^{-2}\left[(\pd_r \Gamma^{\le l+3}u)^2+(\Gamma^{\le l+3}u)^2\right].
\end{multline*}
\end{proof}

\subsection{Estimating derivatives using the Sobolev-type estimate}

Now, we repeatedly apply Lemma \ref{k:pointwise_lem} to estimate various derivatives with $r$ weights. We will assume that $\psi$ decays sufficiently fast as $r\rightarrow\infty$.

The following lemma estimates $\psi$ and the higher order analogues $\Omega^l\psi^{s-l,k}$.

\begin{lemma}\label{k:psi_pointwise_lem}
For $r\ge r_H$,
$$|r^p\psi^{s,k}|^2\lesssim E^{s+3,k}_{2p}(t)$$
and for $r\ge r_0>r_H$,
$$|\psi^{s,k}|^2\lesssim E^{s+3,k}(t).$$
\end{lemma}
\begin{proof}
First, we apply Lemma \ref{k:pointwise_lem} with $u=r^p\psi^{s-l,k}$.
\begin{align*}
|r^p\Omega^l\psi^{s-l,k}|^2 &\lesssim \int_{\Sigma_t} r^{-2}\left[(\pd_r\Gamma^{\le l+3}(r^p\psi^{s-l,k}))^2+(\Gamma^{\le l+3}(r^p\psi^{s-l,k}))^2\right] \\
&\lesssim \int_{\Sigma_t}r^{2p-2}\left[(\pd_r\psi^{s+3,k})^2+(\psi^{s+3,k})^2\right] \\
&\lesssim E_{2p}^{s+3}(t).
\end{align*}
This verifies the first estimate. The second estimate follows from the exact same argument in the special case $p=0$, and the observation that as long as $r\ge r_0>r_H$, then $E^{s+3}(t)$ can be used in place of of $E_0^{s+3}(t)$.
\end{proof}

The following lemma estimates $\pd_t\psi$ and the higher order analogues $\pd_t\Omega^l\psi^{s-l,k}$.
\begin{lemma}
For $r\ge r_H$,
$$|r^p\pd_t\Omega^l\psi^{s-l,k}|^2\lesssim E_{2p}^{s+4}(t)$$
and for $r\ge r_0>r_H$,
$$|r\pd_t\Omega^l\psi^{s-l,k}|^2\lesssim E^{s+4}(t)$$
\end{lemma}
\begin{proof}
The first estimate reduces to Lemma \ref{k:psi_pointwise_lem} by observing that $\pd_t\Omega^l\psi^{s-l,k}=\Omega^l\psi^{s+1-l,k}$. We now prove the second estimate.

First, we apply Lemma \ref{k:pointwise_lem} with $u=r\pd_t\Omega^l\psi^{s-l,k}$.
\begin{align*}
|r\pd_t\Omega^l\psi^{s-l,k}|^2 &\lesssim \int_{\Sigma_t\cap\{r> r_0\}}r^{-2}\left[(\pd_r\Gamma^{\le l+3}(r\pd_t\psi^{s-l,k}))^2+(\Gamma^{\le l+3}(r\pd_t\psi^{s-l,k}))^2\right] \\
&\lesssim \int_{\Sigma_t\cap\{r> r_0\}}\left[(\pd_r\Gamma^{\le l+3}(\pd_t\psi^{s-l,k}))^2+(\Gamma^{\le l+3}(\pd_t\psi^{s-l,k}))^2\right] \\
&\lesssim \int_{\Sigma_t\cap\{r> r_0\}}\left[(\pd_r\psi^{s+4,k})^2+(\pd_t\psi^{s+3,k})^2\right] \\
&\lesssim E^{s+3}(t).
\end{align*}
The point is that $\pd_t\psi^{s+3,k}$ can either be treated as $\psi^{s+4,k}$ or as a derivative of $\psi^{s+3,k}$. In the latter case, the energy norm $E^s(t)$ has stronger control, because $||\pd_t\psi^s||_{L^2(\Sigma_t\cap\{r>r_0\})}^2\le E^s(t)$, while $||r^{-1}\psi^{s+1}||_{L^2(\Sigma_t\cap\{r>r_0\})}^2\le E^{s+1}(t)$.
\end{proof}

The following lemma estimates $\sla\nabla\psi$ and the higher order analogues $\sla\nabla\Omega^l\psi^{s-l,k}$.
\begin{lemma}
For $r\ge r_H$,
$$|r^{p+1}\sla\nabla\Omega^l\psi^{s-l,k}|^2\lesssim E^{s+4}_{2p}(t)$$
and for $r\ge r_0>r_H$,
$$|r\sla\nabla\Omega^l\psi^{s-l,k}|^2\lesssim E^{s+4}(t)$$
\end{lemma}
\begin{proof}
This lemma reduces to Lemma \ref{k:psi_pointwise_lem} by observing that
$$r\sla\nabla\Omega^l\psi^{s-l,k}=\Omega^{l+1}\psi^{s+1-(l+1),k},$$
which is represented by $\Omega^{l+1}\psi^{s+1-l,k}$.
\end{proof}

The following lemma estimates $L\psi$ and the higher order analogues $L\Omega^l\psi^{s-l,k}$.
\begin{lemma}\label{k:L_pointwise_lem}
Letting $L=\alpha\pd_r+\pd_t$, where $\alpha=\frac{\Delta}{r^2+a^2}$, we have that for $r\ge r_H$,
$$|r^{p+1}L\Omega^l\psi^{s-l,k}|^2\lesssim E_{2p}^{s+5,k}(t)+\int_{\Sigma_t}r^{2p}(\Box_g\psi^{s+3,k})^2$$
and for $r\ge r_0>r_H$,
$$|rL\Omega^l\psi^{s-l,k}|^2\lesssim E^{s+5,k}(t)+\int_{\Sigma_t\cap\{r>r_0\}}(\Box_g\psi^{s+3,k})^2.$$
\end{lemma}
\begin{proof}
Before beginning the estimates stated by the lemma, it is important to establish
$$(\pd_rLu)^2\lesssim (\Box_gu)^2+(L\pd_tu)^2+r^{-2}(\pd_t^2u)^2+r^{-2}(\pd_r\pd_tu)^2+r^{-2}(\pd_ru)^2+r^{-2}(Qu)^2,$$
To verify this, we expand
$$q^2\Box_g=q^2g^{tt}\pd_t^2+q^2g^{rr}\pd_r^2+\pd_r(q^2g^{rr})\pd_r+q^2g^{rt}\pd_r\pd_t+\pd_r(q^2g^{rt})\pd_t+Q$$
and observe that for $r>r_H+\delh$, $g^{rt}=0$ and
$$q^2g^{tt}\pd_t^2u+q^2g^{rr}\pd_r^2u+\pd_r(q^2g^{rr})\pd_ru=-\frac{(r^2+a^2)^2}{\Delta}\pd_t^2u+\Delta\pd_r^2u+\Delta'\pd_ru.$$
We also expand
\begin{align*}
(r^2+a^2)\pd_rLu &= (r^2+a^2)\pd_r\left(\frac{\Delta}{r^2+a^2}\pd_ru+\pd_tu\right) \\
&= \Delta\pd_r^2u+\Delta'\pd_ru-\frac{2r\Delta}{r^2+a^2}\pd_ru+(r^2+a^2)\pd_r\pd_tu
\end{align*}
Note that both of these expressions share the terms $\Delta\pd_r^2u$ and $\Delta'\pd_ru$. It follows that for $r\ge r_H+\delh$,
\begin{multline*}
q^2g^{tt}\pd_t^2u+q^2g^{rr}\pd_r^2u+\pd_r(q^2g^{rr})\pd_ru - (r^2+a^2)\pd_rLu \\
=-\frac{(r^2+a^2)^2}{\Delta}\pd_t^2u-(r^2+a^2)\pd_r\pd_tu + \frac{2r\Delta}{r^2+a^2}\pd_ru \\
=-\frac{(r^2+a^2)^2}{\Delta}\left(\pd_t^2u+\frac{\Delta}{r^2+a^2}\pd_r\pd_tu\right)+\frac{2r\Delta}{r^2+a^2}\pd_ru \\
=-\frac{(r^2+a^2)^2}{\Delta}L\pd_tu+\frac{2r\Delta}{r^2+a^2}\pd_ru.
\end{multline*}
Keeping in mind the additional terms that show up for $r\le r_H+ \delh$, we arrive at the estimate for $(\pd_rLu)^2$ at the beginning of this proof. With this estimate in mind, we begin to prove the estimates stated by the lemma.

We apply Lemma \ref{k:pointwise_lem} with $u=r^{p+1}L\psi^{s-l,k}$.
\begin{align*}
|r^{p+1}L\Omega^l\psi^{s-l,k}|^2 &= |\Omega^l(r^{p+1}L\psi^{s-l,k})|^2 \\
&\lesssim \int_{\Sigma_t}r^{-2}\left[(\pd_r\Gamma^{\le l+3}(r^{p+1}L\psi^{s-l,k}))^2+(\Gamma^{\le l+3}(r^{p+1}L\psi^{s-l,k}))^2\right] \\
&\lesssim \int_{\Sigma_t}r^{2p}\left[(\pd_r\Gamma^{\le l+3}L\psi^{s-l,k})^2+(\Gamma^{\le l+3}L\psi^{s-l,k})^2\right] \\
&\lesssim \int_{\Sigma_t}r^{2p}\left[(\pd_rL\psi^{s+3,k})^2+(L\psi^{s+3,k})^2\right] \\
&\lesssim E_{2p}^{s+3,k}(t) +\int_{\Sigma_t}r^{2p}(\pd_rL\psi^{s+3,k})^2.
\end{align*}
Now, according to the estimate we previously established,
\begin{align*}
\int_{\Sigma_t}r^{2p}&(\pd_rL\psi^{s+3,k})^2 \\
&\lesssim \int_{\Sigma_t}r^{2p}\left[(\Box_g\psi^{s+3,k})^2+(L\pd_t\psi^{s+3,k})^2+r^{-2}(\pd_t^2\psi^{s+3,k})^2+r^{-2}(\pd_r\pd_t\psi^{s+3,k})^2\right.\\
&\hspace{3in}\left.+r^{-2}(\pd_r\psi^{s+3,k})^2+r^{-2}(Q\psi^{s+3,k})^2\right] \\
&\lesssim \int_{\Sigma_t}r^{2p}\left[(\Box_g\psi^{s+3,k})^2+(L\psi^{s+4,k})^2+r^{-2}(\psi^{s+5,k})^2+r^{-2}(\pd_r\psi^{s+4,k})^2\right].
\end{align*}
It follows that
$$|r^{p+1}L\Omega^l\psi^{s-l,k}|^2\lesssim E_{2p}^{s+5,k}(t)+\int_{\Sigma_t}r^{2p}(\Box_g\psi^{s+3,k})^2.$$
This is the first estimate of the lemma. The second estimate follows from the same exact argument in the special case $p=0$, and the observation that as long as $r\ge r_0>r_H$, then $E^{s,k}(t)$ can be used in place of $E_0^{s,k}(t)$.
\end{proof}

The following lemma estimates $\tg\psi$ as well as the higher order analogues $\tg\Omega^l\psi^{s-l,k}$. (This is new compared previous chapters.)
\begin{lemma}\label{k:gh_infty_lem}
Keeping in mind that $\tg$ is supported in a neighborhood of the event horizon, for arbitrary $p'$, we have for $r\ge r_H$,
$$|\tg \Omega^l\psi^{s-l,k}|^2 \lesssim E_{p'}^{s+3,k+1}(t)$$
and for $r\ge r_0>r_H$,
$$|\tg \Omega^l\psi^{s-l,k}|^2 \lesssim E^{s+3,k+1}(t).$$
\end{lemma}
\begin{proof}
We apply Lemma \ref{k:pointwise_lem} with $u=\tg\psi^{s-l,k}$, and freely introduce a factor of $r^{p'}$ since $\tg$ is supported on a compact interval in $r$.
\begin{align*}
|\tg\Omega^l\psi^{s-l,k}|^2 &= |\Omega^l\psi^{s-l,k+1}|^2 \\
&\lesssim \int_{\Sigma_t}r^{-2}\left[(\pd_r\Gamma^{\le l+3}\psi^{s-l,k+1})^2+(\Gamma^{\le l+3}\psi^{s-l,k+1})^2\right] \\
&\lesssim \int_{\Sigma_t}r^{p'-2}\left[(\pd_r\Gamma^{\le l+3}\psi^{s-l,k+1})^2+(\Gamma^{\le l+3}\psi^{s-l,k+1})^2\right] \\
&\lesssim E_{p'}^{s+3,k+1}(t).
\end{align*}
This proves the first estimate of the lemma. The second estimate follows from the same argument, and the observation that as long as $r\ge r_0>r_H$, then $E^{s,k}(t)$ can be used in place of $E_{p'}^{s,k}(t)$.
\end{proof}

\subsection{Summarizing the $L^\infty$ estimates}

To conclude this section, we summarize the previous lemmas in a single proposition, making use of the following definition.
\begin{definition}
We define two families of operators
$$\bar{D}=\{L,\sla\nabla\}$$
$$D=\{L,\pd_r,\sla\nabla\}$$
\end{definition}

\begin{proposition}\label{k:infinity_prop}
For $r\ge r_H$,
\begin{multline*}
|r^{p+1}\bar{D}\Omega^l\psi^{s-l,k}|^2+|r^pD\Omega^l\psi^{s-l,k}|^2+|r^p\Omega^l\psi^{s-l,k}|^2 \\
\lesssim E_{2p}^{s+3,k+1}(t)+E_{2p}^{s+5,k}(t)+\int_{\Sigma_t}r^{2p}(\Box_g\psi^{s+3,k})^2,
\end{multline*}
and for $r\ge r_0>r_H$,
$$|rD\Omega^l\psi^{s-l,k}|^2\lesssim E^{s+3,k+1}(t)+E^{s+5,k}(t)+\int_{\Sigma_t}(\Box_g\psi^{s+3,k})^2.$$
\end{proposition}

\begin{proof}
With the exception of $\pd_r$, all of the cases have been proved in Lemmas \ref{k:psi_pointwise_lem}-\ref{k:gh_infty_lem}. Finally, observe that
$$|r^p\pd_r\Omega^l\psi^{s-l,k}|^2\lesssim |r^pL\Omega^l\psi^{s-l,k}|^2+|r^p\pd_t\Omega^l\psi^{s-l,k}|^2+|\tg\Omega^l\psi^{s-l,k}|^2$$
Thus, even the case of the operator $\pd_r$ can be reduced to Lemmas \ref{k:psi_pointwise_lem}-\ref{k:gh_infty_lem}.
\end{proof}

\begin{remark}
In the above proof, we see the reason for the need for the commutator $\tg$. Note that Lemmas \ref{k:psi_pointwise_lem}-\ref{k:L_pointwise_lem} do not have an increase in $k$ on the right side, but Lemma \ref{k:gh_infty_lem} does. Lemma \ref{k:gh_infty_lem} was needed in order to estimate the $\pd_r$ derivative near the event horizon, because $L$ coincides with $\pd_t$ on the event horizon. Excluding this issue, there would be no need to commute with $\tg$ and introduce the $k$ index.
\end{remark}

\section{Theorem: Global boundedness and decay for axisymmetric solutions to the semilinear wave equation with weak null structure on subextremal Kerr spacetimes}\label{k:main_sec}
\subsection{Structure of the nonlinear term}\label{k:nonlinear_structure_sec}

Now, we define the weak null condition for a nonlinear term, which will be assumed by the main theorem of this chapter.

\begin{definition}\label{k:nonlinear_def}
Suppose $\Box_g\psi=\mathcal{N}$, where $\mathcal{N}$ is nonlinear in zeroth and first order derivatives of $\psi$. We say $\mathcal{N}$ satisfies the weak null condition if it can be written as
$$\mathcal{N}=\sum \gamma\beta(r\beta)^j,$$
where the sum has a finite number of terms, and for each term, $0\le j$ and
$$\gamma\in\{L\psi,r^{-1}\pd_r\psi,\sla\nabla\psi,r^{-1}\psi\},$$
$$\beta\in\{L\psi,\pd_r\psi,\sla\nabla\psi,r^{-1}\psi\}.$$
The $\gamma$ factors can be thought of as ``good'' factors, because they can be estimated with an extra $r$ factor using the $L^\infty$ estimates (Proposition \ref{k:infinity_prop}) while, in contrast, the $\beta$ factors can be thought of as potentially ``bad''. The null condition requires that at least one factor be a good factor.
\end{definition}

\begin{proposition}\label{k:commutator_prop}
Generalize the terms $\gamma$ and $\beta$ as follows.
$$\gamma^{s,k}\in\{L\Omega^l\psi^{s-l,k},r^{-1}\pd_r\Omega^l\psi^{s-l,k},\sla\nabla\Omega^l\psi^{s-l,k},r^{-1}\Omega^l\psi^{s-l,k}\},$$
$$\beta^{s,k}\in\{L\Omega^l\psi^{s-l,k},\pd_r\Omega^l\psi^{s-l,k},\sla\nabla\Omega^l\psi^{s-l,k},r^{-1}\Omega^l\psi^{s-l,k}\}.$$
If $\Box_g\psi=\mathcal{N}$ and $\mathcal{N}$ satisfies the weak null condition, then
$$|q^{-2}\tg^k\Gamma^s(q^2\Box_g\psi)|\lesssim \sum |\gamma^{hi}\beta^{lo}(r\beta^{lo})^j|+\sum |\beta^{hi}\gamma^{lo}(r\beta^{lo})^j|,$$
where $\gamma^{hi}$ and $\beta^{hi}$ represent $\gamma^{s',k'}$ and $\beta^{s',k'}$ with $s'\le s$ and $k'\le k$, while $\gamma^{lo}$ and $\beta^{lo}$ represent $\gamma^{s',k'}$ and $\beta^{s',k'}$ with $s'+2k'\le (s+2k)/2$.
\end{proposition}

\begin{proof}
For the proof, we assume the sum consists of just a single term with $j=0$. (It should be clear how to generalize to terms for which $j>0$.) 
$$\Box_g\psi=\gamma\beta.$$
Note that since $\pd_t$ is a first order operator and satisfies the Leibniz rule,
$$|\pd_t^s\Box_g\psi|=|\pd_t^s(\gamma\beta)|\lesssim \sum_{i\le s} |\pd_t^i\gamma\pd_t^{s-i}\beta| \lesssim |\gamma^{hi}\beta^{lo}|+|\beta^{hi}\gamma^{lo}|,$$
The last step follows by observing that for all $i$ in the sum, either $i$ or $s-i$ is less than $s/2$.

We ignore the factors $q^{-2}$ and $q^2$, because whenever an operator acts on $q^2$, the result is a term with fewer derivatives of $\psi$ and an additional $r^{-2}$ weight.

The calculation becomes a bit more complicated when applying $Q$, because $Q$ is a second order operator. In general,
$$Q(\psi_1\psi_2)=Q(\psi_1)\psi_2+\psi_1Q(\psi_2)+2a^2\sin^2\theta \pd_t\psi_1\pd_t\psi_2+2\Omega\psi_1\Omega \psi_2$$
Using the above formula it is straightforward to show that, up to terms with fewer derivatives,
$$q^{-2}\Gamma^s(q^2\psi_1\psi_2)\approx \Omega^l\psi_1^{s_1}\Omega^l\psi_2^{s_2},$$
where $s_1+s_2\le s-2l$, because $l$ counts the number of times the angular derivatives from $Q$ split between the two factors.

Keeping the above observations in mind, it should now be clear that
$$|q^{-2}\tg^k\Gamma^s(q^2\Box_g\psi)|=|q^{-2}\tg^k\Gamma^s(q^2\gamma\beta)|
\lesssim \sum_{\substack{s_1+s_2\le s \\ k_1+k_2\le k}}|\gamma^{s_1,k_1}\beta^{s_2,k_2}|.$$
Now, since $(s_1+2k_1)+(s_2+2k_2)=(s_1+s_2)+2(k_1+k_2)\le s+2k$, it follows that one of $s_1+2k_1$ and $s_2+2k_2$ must be less than $(s+2k)/2$. Choose a particular term in the sum. If $s_1+2k_1\le s_2+2k_2$, then
$$|\gamma^{s_1,k_1}\beta^{s_2,k_2}|\le |\beta^{hi}\gamma^{lo}|,$$
while if $s_1+2k_1\ge s_2+2k_2$, then
$$|\gamma^{s_1,k_1}\beta^{s_2,k_2}|\le |\gamma^{hi}\beta^{lo}|.$$
\end{proof}

\subsection{The main theorem}\label{k:main_thm_sec}

\begin{theorem}\label{k:main_thm}
Let $g$ be any subextremal ($|a|<M$) Kerr metric, and suppose a function $\psi$ is axisymmetric
$$\pd_\phi\psi=0,$$
and solves an equation of the form
$$\Box_g\psi=\mathcal{N},$$
where $\mathcal{N}$ satisfies the weak null condition as defined in \S\ref{k:nonlinear_structure_sec}. Define the energies
$$E^{\ul{n}}(t)=\sum_{s+2k=n}E^{s,k}(t),$$
$$E_p^{\ul{n}}(t)=\sum_{s+2k=n}E_p^{s,k}(t).$$
Then for $\delp,\delm>0$ sufficiently small, if the initial data on $\Sigma_0$ decay sufficiently fast as $r\rightarrow\infty$ and have size
\begin{equation}
I_0=E^{\ul{22}}(0)+E_{2-\delp}^{\ul{22}}(0)
\end{equation}
sufficiently small, then the following estimates hold for $t\ge 0$ (with $T=1+t$).

I) The energies satisfy
$$E^{\ul{22}}(t)\lesssim I_0,$$
$$E^{\ul{22}}_{p\in[\delm,2-\delp]}(t)\lesssim I_0,$$
$$E^{\ul{20}}_{p\in[1-\delp,2-\delp]}(t)\lesssim T^{p-2+\delp}I_0,$$
$$E^{\ul{18}}_{p\in[\delm,2-\delp]}(t)\lesssim T^{p-2+\delp}I_0,$$
$$\int_{t}^{\infty}E_{p\in[\delm-1,\delm]}^{\ul{16}}(\tau)d\tau\lesssim T^{p-2+\delp+1}I_0.$$

II) The following $L^\infty$ estimates hold, provided $s+2k\le 22$.
$$|r^{p+1}\bar{D} \Omega^l\psi^{s-l,k}|^2+|r^pD\Omega^l\psi^{s-l,k}|^2+|r^p\Omega^l\psi^{s-l,k}|^2\lesssim E^{s+5,k}_{2p}(t)+E^{s+3,k+1}_{2p}(t),$$
and for $r\ge r_0>r_H$,
$$|rD\Omega^l\psi^{s-l,k}|^2\lesssim E^{s+5,k}(t)+E^{s+3,k+1}(t).$$

III) Together, (I) and (II) imply that if $s+2k\le 16$, for all $p\in [\delm/2,(2-\delp)/2]$,
$$|r^{p+1}\bar{D} \Omega^l\psi^{s-l,k}|+|r^pD\Omega^l\psi^{s-l,k}|+|r^p\Omega^l\psi^{s-l,k}|\lesssim T^{(2p-2+\delp)/2}I_0^{1/2},$$
$$|rD\Omega^l\psi^{s-l,k}|\lesssim I_0^{1/2},$$
and additionally for $p\in [(\delm-1)/2,\delm/2]$,
$$\int_t^\infty |r^{p+1}\bar{D}\Omega^l\psi^{s-l,k}|+|r^pD\Omega^l\psi^{s-l,k}|+|r^p\Omega^l\psi^{s-l,k}|\lesssim T^{(2p-2+\delp)/2+1}I_0^{1/2}.$$
The final estimate should be interpreted as saying that $|r^{(\delm+1)/2}\bar{D}\Omega^l\psi^{s-l,k}|$, \\
 $|r^{(\delm-1)/2}D\Omega^l\psi^{s-l,k}|$, and $|r^{(\delm-1)/2}\Omega^l\psi^{s-l,k}|$ decay like $T^{(\delm-3+\delp)/2}$ in a weak sense.
\end{theorem}

The remainder of this section is devoted to proving Theorem \ref{k:main_thm}.

\subsection{Bootstrap assumptions}\label{k:bootstrap_assumptions_sec}

We begin the proof of Theorem \ref{k:main_thm} by making the following bootstrap assumptions.
\begin{align*}
E^{\ul{22}}(t)&\le C_bI_0, \\
E^{\ul{22}}_{p\in[\delm,2-\delp]}(t)&\le C_b I_0, \\
E^{\ul{18}}_{p\in[\delm,2-\delp]}(t)&\le C_b T^{p-2+\delp}I_0, \\
\int_t^\infty E^{\ul{16}}_{p\in[\delm-1,\delm]}(\tau)d\tau &\le C_b T^{p-2+\delp+1}I_0.
\end{align*}
Note that, with the exception of the highest order energies, these bootstrap assumptions are consistent with the general principle that $E_p^{s,k}(t)\sim T^{p-2+\delp}$, which the reader should keep in mind throughout the proof of the main theorem.

\subsection{Improved $L^\infty$ estimates}

The $L^\infty$ estimates from Proposition \ref{k:infinity_prop} are essential to the argument of the proof of the main theorem. But for better clarity, we first remove the nonlinear (error) terms from these estimates and summarize them in the following lemma.
\begin{lemma}\label{k:simplified_pointwise_lemma}
In the context of the bootstrap assumptions provided in \S\ref{k:bootstrap_assumptions_sec}, the following $L^\infty$ estimates hold for $s+2k\le 22$ and all $p$ in any bounded range.
\begin{align}
|r^{p+1}\bar{D}\Omega^l\psi^{s-l,k}|^2+|r^pD\Omega^l\psi^{s-l,k}|^2+|r^p\Omega^l\psi^{s-l,k}|^2 &\lesssim  E_{2p}^{s+3,k+1}(t)+E_{2p}^{s+5,k}(t), \\
|rD\Omega^l\psi^{s-l,k}|^2 &\lesssim E^{s+3,k+1}(t)+E^{s+5,k}(t). \label{k:classic_energy_pointwise_estimate}
\end{align}
These are the same as the estimates from Proposition \ref{k:infinity_prop}, except that the nonlinear (error) terms have been removed.
\end{lemma}

\begin{proof} 
See the proof of Lemma \ref{m:improved_infty_lem}, keeping in mind Proposition \ref{k:infinity_prop}.
\end{proof}
The conclusion of this lemma is the same as the statement of part (II) of the main theorem.

\begin{remark}
Lemma \ref{k:simplified_pointwise_lemma}, which is a simplified version of Proposition \ref{k:infinity_prop} in the sense that there are no nonlinear (error) terms on the right side of any of the estimates in Lemma \ref{k:simplified_pointwise_lemma}, will be used in the remainder of the proof of the main theorem as a replacement for Proposition \ref{k:infinity_prop}.
\end{remark}

\subsection{Refined estimates for $N^{s,k}(t)$ and $N_p^{s,k}(t)$ (key step)}\label{k:refined_estimates_sec}

The $L^\infty$ estimates given in Lemma \ref{k:simplified_pointwise_lemma} allow us to provide refined estimates for the nonlinear error terms. \textbf{This is the crucial step of the proof.}

\begin{lemma}\label{k:refined_nl_E_lem}
In the context of the bootstrap assumptions provided in \S\ref{k:bootstrap_assumptions_sec}, if $s+2k\le 22$ and $C_bI_0\le 1$, then
$$N^{s,k}(t)\lesssim (E^{s,k}(t))^{1/2}\left((E^{s,k}(t))^{1/2}(E_{\delm-1}^{\ul{16}}(t))^{1/2}+(E^{s,k}_{1-\delm}(t))^{1/2}(E^{\ul{16}}_{\delm-1}(t))^{1/2}\right).$$
\end{lemma}
\begin{proof}
We recall the definitions of $N^{s,k}(t)$ from Corollary \ref{k:s_cor} ($k=0$) and Proposition \ref{k:s_k_prop} $(k\ge 1)$. For all $s$ and $k$,
$$N^{s,k}(t)= (E^{s,k}(t))^{1/2}\sum_{\substack{s'\le s \\ k'\le k}}||q^{-2}\tg^{k'}\Gamma^{s'}(q^2\Box_g\psi)||_{L^2(\Sigma_t)}.$$
Therefore, it suffices to prove the following estimate.
\begin{multline*}
\sum_{\substack{s'\le s \\ k'\le k}}||q^{-2}\tg^{k'}\Gamma^{s'}(q^2\Box_g\psi)||_{L^2(\Sigma_t)} \\
\lesssim (E^{s,k}(t))^{1/2}(E_{\delm-1}^{\ul{16}}(t))^{1/2}+(E^{s,k}_{1-\delm}(t))^{1/2}(E^{\ul{16}}_{\delm-1}(t))^{1/2}.
\end{multline*}
Note that by Proposition \ref{k:commutator_prop}
\begin{align*}
||q^{-2}\tg^k\Gamma^s(q^2\Box_g\psi)||_{L^2(\Sigma_t)} &\lesssim ||\gamma^{hi}\beta^{lo}(r\beta^{lo})^j||_{L^2(\Sigma_t)}+ ||\beta^{hi}\gamma^{lo}(r\beta^{lo})^j||_{L^2(\Sigma_t)} \\
&\lesssim \left(||r^{\frac{1-\delm}2}\gamma^{hi}||_{L^2(\Sigma_t)}||r^{\frac{\delm-1}2}\beta^{lo}||_{L^\infty(\Sigma_t)}\right. \\
&\hspace{1in}\left.\vphantom{||r^{\frac{1-\delm}2}\gamma^{hi}||_{L^2}||r^{\frac{\delm-1}2}\beta^{lo}||_{L^\infty}}
+||\beta^{hi}||_{L^2(\Sigma_t)}||\gamma^{lo}||_{L^\infty(\Sigma_t)}\right)||r\beta^{lo}||_{L^\infty(\Sigma_t)}^j \\
&\lesssim \left((E^{s,k}_{1-\delm}(t))^{1/2}(E^{\ul{16}}_{\delm-1}(t))^{1/2}\right. \\
&\hspace{1in}\left.+(E^{s,k}(t))^{1/2}(E_{\delm-1}^{\ul{16}}(t))^{1/2}\right)(E^{\ul{16}}(t))^{j/2}
\end{align*}
Noting that
$$(E^{\ul{16}}(t))^{j/2}\le (C_bI_0)^{j/2}\le 1^{j/2}=1,$$
we conclude the proof.
\end{proof}

\begin{lemma}\label{k:refined_nl_Ep_lem}
In the context of the bootstrap assumptions provided in \S\ref{k:bootstrap_assumptions_sec}, if $s+2k\le 22$ and $C_bI_0\le 1$, then if $k=0$,
$$N_p^{s,0}(t)\lesssim E^{s,0}(t)B_p^{s/2+5,0}(t)+B_p^{s,0}(t)E^{s/2+5,0}(t) +E^{s,0}_{p'}(t)(E_{p''}^{\ul{16}}(t))^{1/2},$$
while if $k\ge 1$,
$$N_p^{s,k}(t)\lesssim E^{s,k}_{p'}(t)E_{p''}^{\ul{16}}(t).$$
\end{lemma}

\begin{proof}
If $k=0$, then from Corollary \ref{k:s_cor},
$$N_p^{s,0}(t)=\sum_{s'\le s} \int_{\Sigma_t}r^{p+1}(q^{-2}\Gamma^{s'}(q^2\Box_g\psi))^2+\sum_{s'\le s}\int_{\Sigma_t\cap\{r\approx r_{trap}\}}|\pd_t\psi^{s'}q^{-2}\Gamma^{s'}(q^2\Box_g\psi)|,$$
We rewrite for some fixed radius $R$,
\begin{multline*}
N_p^{s,0}(t) = \left[\sum_{s'\le s} \int_{\Sigma_t\cap\{r<R\}}(q^{-2}\Gamma^{s'}(q^2\Box_g\psi))^2+\sum_{s'\le s}\int_{\Sigma_t\cap\{r\approx r_{trap}\}}|\pd_t\psi^{s'}q^{-2}\Gamma^{s'}(q^2\Box_g\psi)|\right] \\
+\sum_{s'\le s} \int_{\Sigma_t\cap\{r>R\}}r^{p+1}(q^{-2}\Gamma^{s'}(q^2\Box_g\psi))^2.
\end{multline*}
Using Proposition \ref{k:commutator_prop} again, the terms in square brackets can be bounded by
$$E_{p'}^s(t)E_{p''}^{s/2+3}(t)(E^{\ul{16}}(t))^{j/2}$$
and
$$E_{p'}^s(t)(E_{p''}^{s/2+3}(t))^{1/2}(E^{\ul{16}}(t))^{j/2}$$
respectively, with $p'$ and $p''$ arbitrary because the terms are supported on a compact radial interval. Furthermore, by the bootstrap assumptions and the fact that $C_bI_0\le 1$, each of these can be bounded by
$$E_{p'}^s(t)(E_{p''}^{s/2+3}(t))^{1/2}.$$
More importantly,
\begin{align*}
\sum_{s'\le s} \int_{\Sigma_t\cap\{r>R\}}&r^{p+1}(q^{-2}\Gamma^{s'}(q^2\Box_g\psi))^2 \\
&\lesssim \int_{\Sigma_t\cap\{r>R\}}r^{p+1}\left[(\gamma^{hi})^2(\beta^{lo})^2(r\beta^{lo})^{2j}+(\beta^{hi})^2(\gamma^{lo})^2(r\beta^{lo})^{2j}\right] \\
&\lesssim \int_{\Sigma_t\cap\{r>R\}}r^{p-1}(\gamma^{hi})^2(r\beta^{lo})^{2j+2}+(\beta^{hi})^2(r^{(p+1)/2}\gamma^{lo})^2(r\beta^{lo})^{2j} \\
&\lesssim \int_{\Sigma_t\cap\{r>R\}}r^{p-1}(\gamma^{hi})^2||r\beta^{lo}||^{2j+2}_{L^\infty(\Sigma_t)} \\
&\hspace{1in}+(\beta^{hi})^2||r^{(p+1)/2}\gamma^{lo}||^2_{L^\infty(\Sigma_t)}||r\beta^{lo}||^{2j}_{L^\infty(\Sigma_t)} \\
&\lesssim B_p^{s,0}(t)||r\beta^{lo}||^{2j+2}_{L^\infty(\Sigma_t)}+E^{s,0}(t)||r^{(p+1)/2}\gamma^{lo}||^2_{L^\infty(\Sigma_t)}||r\beta^{lo}||^{2j}_{L^\infty(\Sigma_t)} \\
&\lesssim B_p^{s,0}(t)(E^{s/2+5,0}(t))^{1+j}+E^{s,0}(t)B_p^{s/2+5,0}(t)(E^{s/2+5,0}(t))^j.
\end{align*}
Again, since $C_bI_0\le 1$, we may ignore the additional powers of $j$ to conclude
$$\sum_{s'\le s} \int_{\Sigma_t\cap\{r>R\}}r^{p+1}(q^{-2}\Gamma^{s'}(q^2\Box_g\psi))^2 
\lesssim B_p^{s,0}(t)(E^{s/2+5,0}(t))+E^{s,0}(t)B_p^{s/2+5,0}(t).$$
This establishes the $k=0$ estimate of the lemma. If $k\ge 1$, then from Proposition \ref{k:s_k_prop},
$$N_p^{s,k}(t) = \sum_{\substack{s'\le s \\ k'\le k}}\int_{\Sigma_t}r^{p+1}(q^{-2}\tg^{k'}\Gamma^{s'}(q^2\Box_g\psi))^2.$$
Since $\tg$ is only supported near the event horizon, these terms are also supported on a compact radial interval, so $N_p^{s,k\ge 1}(t) \lesssim E_{p'}^{s,k}(t)E_{p''}^{\ul{16}}(t)$.
\end{proof}

\begin{corollary}\label{k:NL_absorb_bulk}
In the context of the bootstrap assumptions provided in \S\ref{k:bootstrap_assumptions_sec}, Corollary \ref{k:s_cor} ($k=0$) and Proposition \ref{k:s_k_prop} ($k\ge 1$) together with Lemma \ref{k:refined_nl_Ep_lem} imply that if $s+2k\le 22$ and $C_bI_0$ is sufficiently small, then
\begin{multline*}
E_p^{s,k}(t_2)+\int_{t_1}^{t_2}B_p^{s,k}(t)dt \\
\lesssim E_p^{s,k}(t_1)+(C_bI_0)^{1/2}\int_{t_1}^{t_2}E_p^{s,k}(t)T^{(\delm-3+\delp)/2}dt+\sum_{\substack{s'+2k'\le s+2k \\ k'<k}}\int_{t_1}^{t_2}B_p^{s',k'}(t)dt.
\end{multline*}
whenever $s/2+5\le s$ (which means $10\le s$).
\end{corollary}
\begin{proof}
By Corollary \ref{k:s_cor} ($k=0$) and Proposition \ref{k:s_k_prop} ($k\ge 1$) together with Lemma \ref{k:refined_nl_Ep_lem},
\begin{multline*}
E_p^{s,k}(t_2)+\int_{t_1}^{t_2}B_p^{s,k}(t)dt \\
\lesssim E_p^{s,k}(t_1)+\int_{t_1}^{t_2}\left[B_{p'}^{s+2,k-1}(t)+E^{s,0}(t)B_p^{s/2+5,0}(t)+B_p^{s,0}(t)E^{s/2+5,0}(t)\right. \\
\left.+E^{s,k}_{p''}(t)(E_{p'''}^{\ul{16}}(t))^{1/2}\right]dt
\end{multline*}
We briefly discuss how to treat each term in the square brackets. The first term, which is only relevant when $k\ge 1$, is estimated as follows
$$B_{p'}^{s+2,k-1}(t) \lesssim \sum_{\substack{s'+2k'\le s+2k \\ k'<k}}B_p^{s',k'}(t),$$
with the value $p$ chosen for $p'$.

The second term is absorbed into the left side of the estimate by a now-standard argument (the bootstrap assumptions ensure that $E^{s,0}(t)\lesssim C_bI_0$, so as long as $C_bI_0$ is sufficiently small, given that $s/2+5\le s$, this can be done). The third term is treated like the second term, but it only requires the weaker assumption that $s/2+5\le 22$, which is always the case.

The fourth term is estimated by choosing $p''=p$ and $p'''=\delm-1$ and appealing to the weak decay principle, which, together with the bootstrap assumptions, implies
$$\int_{t_1}^{t_2}E^{s,k}_{p''}(t)(E_{p'''}^{\ul{16}}(t))^{1/2} \lesssim\int_{t_1}^{t_2}E_p^{s,k}(t)(C_bI_0T^{\delm-3+\delp})^{1/2}dt.$$
\end{proof}

\subsection{Inductive assumptions}\label{k:inductive_assumptions_sec}

The remainder of the proof is a finite induction argument. First, estimates are proved for $k=0$, and then for $k=1$, etc. until $k=11$ (which saturates $s+2k=22$). For each $k$, it will be necessary to use estimates established for $k-1$. These assumptions are listed here.

Either
$$k=0,$$
or $s+2k=22$ and
$$\sum_{\substack{s'+2k'\le 22 \\ k'<k}}\int_t^\infty B_{p\in[\delm,2-\delp]}^{s',k'}(\tau)d\tau \lesssim I_0,$$
or $s+2k=20$ and
$$\sum_{\substack{s'+2k'\le 20 \\ k'<k}}\int_t^\infty B_{p\in[1-\delp,2-\delp]}^{s',k'}(\tau)d\tau \lesssim T^{p-2+\delp}I_0,$$
or $s+2k=18$ and
$$\sum_{\substack{s'+2k'\le 18 \\ k'<k}}\int_t^\infty B_{p\in[\delm,2-\delp]}^{s',k'}(\tau)d\tau \lesssim T^{p-2+\delp}I_0.$$

For the remainder of the proof, $k$ should be considered fixed. These estimates will be used for the fixed $k$, and eventually (in \S\ref{k:establish_inductive_assumptions_sec}) the corresponding estimates obtained by replacing $k$ with $k+1$ will be proved, thus closing the inductive argument.

\subsection{Recovering boundedness of $E^{s,k}_p(t)$ ($s+2k=22$)}\label{k:recover_Ep_sec}

Our first application of Corollary \ref{k:NL_absorb_bulk} is to prove boundedness of $E^{s,k}_p(t)$. Let $s+2k=22$ and $p\in [\delm,2-\delp]$. Then
\begin{align*}
E_p^{s,k}(t)+\int_0^tB_p^{s,k}(\tau)d\tau &\lesssim E_p^{s,k}(t)+(C_bI_0)^{1/2}\int_0^tE_p^{s,k}(\tau)T^{(\delm-3+\delp)/2}d\tau \\
&\hspace{2in}+\sum_{\substack{s'+2k'\le s+2k \\ k'<k}}\int_0^tB_p^{s',k'}(\tau)d\tau \\
&\lesssim I_0+(C_bI_0)^{3/2}T^{\delm-1/2+\delp}+I_0 \\
&\lesssim (1+C_b^{3/2}I_0^{1/2})I_0.
\end{align*}
It follows that if $C_b^{3/2}I_0^{1/2}$ is sufficiently small, then
$$E_p^{s,k}(t)\lesssim I_0.$$
This recovers the second bootstrap assumption at the level $k$.

\subsection{Recovering boundedness of $E^{s,k}(t)$ ($s+2k=22$)}

Let $s+2k=22$. By Corollary \ref{k:s_cor} ($k=0$) and Proposition \ref{k:s_k_prop} ($k\ge 1$) and Lemma \ref{k:refined_nl_E_lem},
\begin{align*}
E^{s,k}(t) \lesssim & E^{s,k}(0)+\int_0^tB_{2-\delp}^{s,k}(\tau)+B_{2-\delp}^{s+2,k-1}(\tau)+N^{s,k}(\tau)d\tau \\
\lesssim & I_0+\int_0^t N^{s,k}(\tau)d\tau \\
\lesssim & I_0 +\int_0^t(E^{s,k}(\tau))^{1/2}\left((E^{s,k}(\tau))^{1/2}(E_{\delm-1}^{\ul{16}}(\tau))^{1/2}\right. \\
&\hspace{2.5in}\left.+(E^{s,k}_{1-\delm}(\tau))^{1/2}(E^{\ul{16}}_{\delm-1}(\tau))^{1/2}\right)d\tau \\
\lesssim & I_0+\int_0^t(C_bI_0)^{1/2}(C_bI_0)^{1/2}(C_bT^{\delm-3+\delp}I_0)^{1/2}d\tau \\
\lesssim & (1+C_b^{3/2}I_0^{1/2})I_0.
\end{align*}
In particular, we used the inductive assumption and the result from \S\ref{k:recover_Ep_sec} in the second step and the weak decay principle in the fourth step.

It follows that if $C_b^{3/2}I_0^{1/2}$ is sufficiently small, then
$$E^{s,k}(t)\lesssim I_0.$$
This recovers the first bootstrap assumption at the level $k$.

\subsection{Proving decay for $E^{s,k}_{p\in[1-\delp,2-\delp]}(t)$ ($s+2k=20$) and $E^{s,k}_{p\in[\delm,2-\delp]}(t)$ ($s+2k=18$)}

At last we are ready to prove decay. This will be accomplished by repeatedly applying the following lemma.

\begin{lemma}\label{k:NL_WE_decay}
Suppose  $p+1,p\in [\delm,2-\delp]$ and
$$E_{p+1}^{s+2,k}(t)\lesssim T^{(p+1)-2+\delp}I_0,$$
$$\sum_{\substack{s'+2k'\le s+2k+2 \\ k'<k}}\int_t^\infty B_{p+1}^{s',k'}(\tau)d\tau\lesssim T^{(p+1)-2+\delp}I_0,$$
$$E_p^{s,k}(t)\le C_bT^{p-2+\delp}I_0,$$
$$\sum_{\substack{s'+2k'\le s+2k \\ k'<k}}\int_t^\infty B_p^{s',k'}(\tau)d\tau\lesssim T^{p-2+\delp}I_0.$$
Then if $I_0$ is sufficiently small,
$$E_p^{s,k}(t)\lesssim T^{p-2+\delp}I_0.$$
\end{lemma}

\begin{proof}
Using the mean value theorem, for a given $t$, let $t'\in [t/2,t]$ be the value for which $B_{p+1}^{s+2,k}(t')=\frac2t\int_{t/2}^tB_{p+1}^{s+2,k}(\tau)d\tau$. Then using Corollary \ref{k:NL_absorb_bulk},
\begin{align*}
E_p^{s,k}(t) &\lesssim E_p^{s,k}(t')+(C_bI_0)^{1/2}\int_{t'}^tE_p^{s,k}(\tau)T^{(\delm-3+\delp)/2}d\tau+\sum_{\substack{s'+2k'\le s+2k \\ k'<k}}\int_{t'}^tB_p^{s',k'}(\tau)d\tau  \\
&\lesssim E_p^{s,k}(t')+(C_bI_0)^{1/2}\int_{t'}^t(C_bT^{p-2+\delp}I_0)T^{(\delm-3+\delp)/2}d\tau+T^{p-2+\delp}I_0 \\
&\lesssim E_p^{s,k}(t')+C_b^{3/2}I_0^{1/2}T^{p-2+\delp}I_0+T^{p-2+\delp}I_0.
\end{align*}
Now by the choice of $t'$,
\begin{multline*}
E_p^{s,k}(t')\lesssim B_{p+1}^{s+2,k}(t')=\frac2t\int_{t/2}^tB_{p+1}^{s+2,k}(\tau)d\tau \\
\lesssim t^{-1}E_{p+1}^{s+2,k}(t/2)+t^{-1}(C_bI_0)^{1/2}\int_{t/2}^tE_{p+1}^{s+2,k}(\tau)T^{(\delm-3+\delp)/2}d\tau \\
+t^{-1}\sum_{\substack{s'+2k'\le s+2k+2 \\ k'<k}}\int_{t/2}^t B_{p+1}^{s',k'}(\tau)d\tau.
\end{multline*}
Thus,
$$E_p^{s,k}(t)\lesssim (1+C_b^{3/2}I_0^{1/2})T^{p-2+\delp}I_0.$$
Provided $C_b^{3/2}I_0^{1/2}$ is sufficiently small, the proof is complete.
\end{proof}

By applying Lemma \ref{k:NL_WE_decay} for $p=1-\delp$ and $s+2k=20$, we obtain
$$E_{1-\delp}^{s,k}(t)\lesssim T^{-1}I_0.$$
By interpolation we obtain decay for all $p\in [1-\delp,2-\delp]$.
$$E_{p\in [1-\delp,2-\delp]}^{s,k}(t)\lesssim T^{p-2+\delp}I_0.$$
Then applying Lemma \ref{k:NL_WE_decay} again for each $p\in[\delm,1-\delp]$ and $s+2k=18$, we obtain
\begin{equation*}
E_{p\in [\delm,2-\delp]}^{s,k}(t)\lesssim T^{p-2+\delp}I_0.
\end{equation*}
This recovers the third bootstrap assumption at the level $k$.

\subsection{Establishing inductive assumptions for $k+1$}\label{k:establish_inductive_assumptions_sec}

By Corollary \ref{k:NL_absorb_bulk},
\begin{multline*}
\int_t^\infty B_p^{s,k}(\tau)d\tau 
\lesssim E_p^{s,k}(t)+(C_bI_0)^{1/2}\int_t^\infty E_p^{s,k}(\tau)T^{(\delm-3+\delp)/2}d\tau \\
+\sum_{\substack{s'+2k'\le s+2k \\ k'<k}}\int_t^\infty B_p^{s',k'}(\tau)d\tau.
\end{multline*}
The quantities $E_p^{s,k}(t)$, and $\int_t^\infty B_p^{s',k'}(\tau)d\tau$ ($s'+2k'\le s+2k$ and $k'<k$) all have the same proven decay rates. For $s+2k=22$ and $p\in[\delm,2-\delp]$, they are bounded in time by $I_0$, for $s+2k=20$ and $p\in[1-\delp,2-\delp]$, they decay at least as fast as $T^{p-2+\delp}I_0$, and for $s+2k=18$ and $p\in[\delm,2-\delp]$, they decay at least as fast as $T^{p-2+\delp}I_0$. Thus,
$$\sum_{\substack{s'+2k'\le 22 \\ k'<k+1}}\int_t^\infty B_{p\in[\delm,2-\delp]}^{s',k'}(\tau)d\tau \lesssim I_0,$$
$$\sum_{\substack{s'+2k'\le 20 \\ k'<k+1}}\int_t^\infty B_{p\in[1-\delp,2-\delp]}^{s',k'}(\tau)d\tau \lesssim T^{p-2+\delp}I_0,$$
$$\sum_{\substack{s'+2k'\le 18 \\ k'<k+1}}\int_t^\infty B_{p\in[\delm,2-\delp]}^{s',k'}(\tau)d\tau \lesssim T^{p-2+\delp}I_0.$$
These are the inductive assumptions at the next level $k+1$.

\subsection{Recovering weak decay for $E_{p\in[\delm-1,\delm]}^{s,k}(t)$ ($s+2k=16$)}

Finally, set $s+2k=16$ and observe that for $p\in [\delm-1,\delm]$,
$$\int_t^\infty E_p^{s,k}(\tau)d\tau \lesssim \int_t^\infty B_{p+1}^{s+2,k}(\tau)d\tau\lesssim T^{(p+1)-2+\delp}I_0=T^{p-2+\delp+1}I_0.$$
This recovers the final bootstrap assumption at the level $k$.

\noindent This completes the proof of Theorem \ref{k:main_thm}. \qed

\chapter{A class of perturbations of the nontrivial wave map solution in the Schwarzschild spacetime}\label{wm_szd_chap}

The following conjecture was posed in \cite{ionescu2014global}.
\begin{conjecture}\label{ws:iokl_conj}
The stationary solution
\begin{align*}
X_0 &= A = \frac{(r^2+a^2)^2-a^2\sin^2\theta(r^2-2Mr+a^2)}{r^2+a^2\cos^2\theta}\sin^2\theta \\
Y_0 &= B = -2aM(3\cos\theta-\cos^3\theta)-\frac{2a^3M\sin^4\theta\cos\theta}{r^2+a^2\cos^2\theta}.
\end{align*}
of the wave map system
\begin{align}
\Box_gX &= \frac{\pd^\alpha X\pd_\alpha X}{X}-\frac{\pd^\alpha Y\pd_\alpha Y}{X}\label{ws:X_eqn} \\
\Box_gY &= 2\frac{\pd^\alpha X\pd_\alpha Y}{X}.\label{ws:Y_eqn}
\end{align}
with $g$ any subextremal ($|a|<M$) Kerr metric is future asymptotically stable in the domain of outer communication of the Kerr spacetime for all smooth, axially symmetric, admissible perturbations, ie. perturbations vanishing in some suitable way on the axis of symmetry.
\end{conjecture}

In this chapter, we resolve the above conjecture for the Schwarzschild case, in which the quantities $A$ and $B$ greatly simplify.
\begin{align*}
A &= r^2\sin^2\theta \\
B &= 0.
\end{align*}
 In the next chapter, we will extend the result of this chapter to slowly rotating ($|a|\ll M$) Kerr spacetimes using a significantly more complicated approach. As discussed in \S\ref{model_problem_motivation_sec}, the motivation for this conjecture is its relation to the stability of Kerr black holes as solutions to the Einstein Vacuum System (\ref{eve}).

\section{Introduction}

There are two signficant differences between the problem of this chapter and the problem in Chapter \ref{szd_chap}. The first significant difference is that the system (\ref{ws:X_eqn}-\ref{ws:Y_eqn}) when linearized about the nontrivial solution yields a different linear equation. This equation is best understood as a scalar wave equation for a modified Schwarzschild spacetime $(\tilde{\mathcal{M}},\tilde{g})$. (For more details, continue reading in \S\ref{ws:phi_psi_sec}.) Fortunately, the techniques in Chapter \ref{szd_chap} easily generalize to this new spacetime. This is established in \S\ref{ws:spacetime_estimates_sec}.

The second significant difference is that the nonlinear terms require substantial more consideration than in Chapter \ref{szd_chap}. In previous chapters, the structure of the nonlinear terms was given by assumption, while in this and the next chapter, the structure is given by the system (\ref{ws:X_eqn}-\ref{ws:Y_eqn}). (This is the purpose of \S\ref{ws:structure_sec}.) As we shall see, there are some challenging nonlinear terms with factors that are singular on the axis of symmetry. In the Schwarzschild case, we observe that by making a very special choice of linearization, these terms can be eliminated. (See Proposition \ref{ws:special_choice_prop}.)

\subsection{The $(\phi,\psi)$ system}\label{ws:phi_psi_sec}

In \cite{ionescu2014global}, the equations (\ref{ws:X_eqn}-\ref{ws:Y_eqn}) are linearized in a certain way to obtain\footnote{The scalar functions $\xi_1$ and $\xi_2$ are named as such to be consistent with the scalar functions $\xi_1$ and $\xi_2$ in Chapter \ref{wm_kerr_chap}. They will not be used any further in this chapter.}
\begin{align*}
\Box_g\xi_1 &= V\xi_1 \\
\Box_g\xi_2 &= 0,
\end{align*}
where the potential $V$ is given by
$$V=\frac4{r^2}\left(1-\frac{2M}r\right)+\frac{4\cot^2\theta}{r^2}.$$
The motivation for such a linearization is that it derives from a bundle formalism. (See \S\ref{wk:xi_a_intro_sec} or \S\ref{wave_map_bundle_sec}.) However, since the potential function $V$ is singular on the axis, this suggests that the scalar $\xi_1$ should vanish on the axis. By making the particular choice
\begin{align*}
\xi_1 &= r^2\sin^2\theta \psi \\
\xi_2 &= -\phi,
\end{align*}
where $\phi$ and $\psi$ are both scalar wavefunctions,\footnote{The scalar wavefunction $\phi$ is not to be confused with the azimuthal coordinate.} then the linearized equations become
\begin{align*}
\Box_g\psi &+\left(1-\frac{2M}r\right)\frac{4}r \pd_r\psi +\frac{4\cot\theta}{r^2}\pd_\theta\psi = 0 \\
\Box_g\phi &= 0.
\end{align*}
The equation for $\phi$ was the main subject of Chapter \ref{szd_chap}, but the equation for $\psi$ is new. There is an interpretation for the additional terms showing up in the equation for $\psi$. The Schwarzschild spacetime can be written as a warped product of two manifolds, one being the sphere $S^2$. If the sphere $S^2$ is instead replaced with the hypersphere $S^6$, yielding a modified $7+1$ dimensional Schwarzschild spacetime, then the additional terms in the equation for $\psi$ are completely accounted for and the linear system reduces to the following simple system.
\begin{align*}
\Box_{\tilde{g}}\psi &= 0 \\
\Box_g\phi &= 0,
\end{align*}
where $\tilde{g}$ is the metric for the modified $7+1$ dmensional Schwarzschild spacetime. This fact suggests that the equations for $\phi$ and $\psi$ are somewhat similar. As we shall see in \S\ref{ws:spacetime_estimates_sec}, the same essential estimates can be proved for this new spacetime.

Perhaps more important than simplifying the equations, the wavefunctions $(\phi,\psi)$ capture the essential behavior of the linear system on the axis. As we shall see in the main theorem of this chapter, the quantity $\psi$ will be regular on the axis, which implies that the quantity $\xi_1=r^2\sin^2\theta\psi$ will vanish to second order on the axis.

\subsection{The modified Schwarzschild spacetime $(\tilde{\mathcal{M}},\tilde{g})$}

It was just stated that the equation for $\psi$ naturally belongs to a modified Schwarzschild spacetime. We take a moment to explan this fact further.

Let $(\mathcal{M},g)$ denote the Schwarzschild spacetime. In the usual Boyer-Lindquist coordinate system, the metric is
$$g_{\alpha\beta}dx^\alpha dx^\beta=-\left(1-\frac{2M}r\right)dt^2+\left(1-\frac{2M}r\right)^{-1}dr^2+r^2d\omega_{S^2}^2,$$
and its volume form is
$$\sqrt{-\det g}=r^2\sin\theta.$$

The equation for $\psi$ is
$$\Box_g\psi +2\frac{\pd^\alpha A}{A}\pd_\alpha \psi = 0.$$
The linear operator in this equation is a wave operator for a different ($7+1$ dimensional) spacetime $(\tilde{\mathcal{M}},\tilde{g})$, whose metric is given by\footnote{The difference between the two metrics $g$ and $\tilde{g}$ is very subtle. For $g$, the last term is $r^2d\omega^2_{S^2}$ and for $\tilde{g}$ the last term is $r^2d\omega^2_{S^6}$.}
$$\tilde{g}_{\alpha\beta}dx^\alpha dx^\beta=-\left(1-\frac{2M}r\right)dt^2+\left(1-\frac{2M}r\right)^{-1}dr^2+r^2d\omega_{S^6}^2,$$
and whose volume form is effectively
$$\sqrt{-\det\tilde{g}}=r^6\sin^5\theta.$$

Let us take a moment to derive this fact. The sphere $S^6$ can be given coordinates $\theta_1,...,\theta_5,\phi$, where $\theta_i\in [0,\pi]$ and $\phi\in [0,2\pi]$. An axisymmetric function $f:\mathcal{M}\rightarrow\R{}$ can be written as $f(t,r,\theta)$. We consider only functions $f:\tilde{\mathcal{M}}\rightarrow\R{}$ of the form $f(t,r,\theta_1)$. We therefore identify $\theta_1$ on $\tilde{\mathcal{M}}$ with $\theta$ on $\mathcal{M}$. In the calculations to follow, $\theta_1$ and $\theta$ will be used interchangeably. The metric for $S^6$ in these coordinates is
$$d\omega^2_{S^6}=d\theta_1^2+\sin^2\theta_1(d\theta_2^2+\sin^2\theta_2(...(d\theta_5^2+\sin^2\theta_5d\phi^2)...)).$$
It follows that
$$\sqrt{-\det\tilde{g}}=r^6\sin^5\theta_1\sin^4\theta_2\sin^3\theta_3\sin^2\theta_4\sin\theta_5.$$

Assume that $\psi:\tilde{\mathcal{M}}\rightarrow\R{}$ satisfies
$$\pd_{\theta_2}\psi=...=\pd_{\theta_5}\psi=\pd_\phi\psi=0.$$
Then since
$$g^{tt}=\tilde{g}^{tt}\text{, }g^{rr}=\tilde{g}^{rr}\text{, and }g^{\theta\theta}=\tilde{g}^{\theta_1\theta_1},$$
it follows that
$$g^{\alpha\beta}\pd_\beta\psi=\tilde{g}^{\alpha\beta}\pd_\beta\psi.$$
Therefore,
\begin{align*}
\Box_{\tilde{g}}\psi &= \frac{1}{\sqrt{-\det\tilde{g}}}\pd_\alpha\left(\sqrt{-\det\tilde{g}}\tilde{g}^{\alpha\beta}\pd_\beta\psi\right) \\
&= \frac{1}{\sqrt{-\det\tilde{g}}}\pd_\alpha\left(\sqrt{-\det\tilde{g}}g^{\alpha\beta}\pd_\beta\psi\right) \\
&= \frac{1}{r^6\sin^5\theta_1}\pd_\beta\left(r^6\sin^5\theta_1g^{\alpha\beta}\pd_\beta\psi\right) \\
&= \frac{1}{r^2\sin\theta}\pd_\alpha\left(r^2\sin\theta g^{\alpha\beta}\pd_\beta\psi\right)+\frac{1}{r^4\sin^4\theta}\pd^\alpha(r^4\sin^4\theta)\pd_\alpha\psi \\
&= \frac{1}{\sqrt{-det g}}\pd_\alpha\left(\sqrt{-\det g}g^{\alpha\beta}\pd_\beta\psi\right)+\frac{2}{r^2\sin^2\theta}\pd^\alpha(r^2\sin^2\theta)\pd_\alpha\psi \\
&= \Box_g\psi +2\frac{\pd^\alpha A}{A}\pd_\alpha\psi.
\end{align*}

\subsection{Conventions for spacetime norms}

Since there are now effectively two different spacetimes, we briefly lay out a few conventions for the remainder of this chapter. \\
\bp All functions will be assumed to depend only on $t$, $r$, and $\theta$. Therefore, any quantity can be treated as a function defined on either spacetime. \\
\bp To avoid ambiguity, the tilde mark ( $\tilde{ }$ ) will be used to denote quantities corresponding to $(\tilde{\mathcal{M}},\tilde{g})$. This includes the effective volume form 
$$\tilde{\mu}=A^2r^2\sin\theta=r^6\sin^5\theta,$$
and the constant-time hypersurface 
$$\tilde{\Sigma}_t=\{t\}\times [r_H,\infty)\times S^6.$$
\bp Integrated expressions will depend on a volume form that is implicitly defined by the manifold of integration. That is,
$$\int_{\Sigma_t}f=\int_{r_H}^\infty\int_0^\pi\int_0^{2\pi}f(t,r,\theta)r^2\sin\theta d\phi d\theta dr.$$
and
\begin{multline*}
\int_{\tilde{\Sigma}_t}f=\int_{r_H}^\infty\int_0^\pi\int_0^\pi\int_0^\pi\int_0^\pi\int_0^\pi\int_0^{2\pi}f(t,r,\theta)r^6\sin^5\theta\sin^4\theta_2\sin^3\theta_3\sin^2\theta_4\sin\theta_5 \\
d\phi d\theta_5 d\theta_4 d\theta_3 d\theta_2 d\theta dr.
\end{multline*}
Since all relevant integral estimates in this chapter are valid up to a constant, one may equivalently define
$$\int_{\Sigma_t}f=\int_{r_H}^\infty\int_0^\pi f(t,r,\theta)r^2\sin\theta d\theta dr$$
and
$$\int_{\tilde{\Sigma}_t}f=\int_{r_H}^\infty\int_0^\pi f(t,r,\theta)r^6\sin^5\theta d\theta dr.$$
The point is that the only difference occurs in the factors that show up in the volume form. \\
\bp We also observe that $L^\infty$ estimates are weaker in the higher-dimensional spacetime. That is,
$$||f||_{L^\infty(S^2)}\lesssim \sum_{i\le 2}||\sla\nabla^if||_{L^2(S^2)},$$
while
$$||f||_{L^\infty(S^6)}\lesssim \sum_{i\le 4}||\tilde{\sla\nabla}^if||_{L^2(S^6)}.$$

\subsection{The nonlinear equations for $(\phi,\psi)$}

In terms of $\phi$ and $\psi$, the nonlinear system is
\begin{align*}
\Box_g\phi &= \mathcal{N}_\phi \\
\Box_{\tilde{g}}\psi &= \mathcal{N}_\psi,
\end{align*}
where the nonlinear terms are
\begin{align*}
(1+\phi)\mathcal{N}_\phi &= \pd^\alpha\phi\pd_\alpha\phi -\mathcal{N} \\
(1+\phi)^2\mathcal{N}_\psi &= 2\psi\mathcal{N}-2(1+\phi)\pd^\alpha\phi\pd_\alpha\psi \\
\mathcal{N} &= A^2(1+\phi)^4 \pd^\alpha\psi\pd_\alpha\psi +4A^2(1+\phi)^3\psi\pd^\alpha\phi\pd_\alpha\psi +4A^2(1+\phi)^2\psi^2\pd^\alpha\phi\pd_\alpha\phi \\
&\hspace{12pt}+ 4A\pd^\alpha A (1+\phi)^4 \psi\pd_\alpha\psi +8 A\pd^\alpha A (1+\phi)^3\psi^2\pd_\alpha\phi +4\pd^\alpha A\pd_\alpha A (1+\phi)^4\psi^2
\end{align*}
and again, the modified wave operator $\Box_{\tilde{g}}$ appearing in the equation for $\psi$ is defined by
$$\Box_{\tilde{g}}=\Box_g+2\frac{\pd^\alpha A}{A}\pd_\alpha.$$
The derivation of this system is given in \S\ref{ws:special_choice_sec}.

\section{Estimates on the spacetime $(\tilde{M},\tilde{g})$}\label{ws:spacetime_estimates_sec}

Since the vectorfield estimates for $\tilde{\mathcal{M}}$ are so similar to the vectorfield estimates for $\mathcal{M}$ derived in Chapter \ref{szd_chap}, and since a more general problem is handled in Chapter \ref{wm_kerr_chap}, we will only prove the most essential vectorfield estimates in this chapter, which are the partial Morawetz estimate and the pre-$r^p$ identity.

\subsection{The energy estimate and $h\pd_t$ estimate for $\tilde{\mathcal{M}}$}

Since the vectorfield $\pd_t$ is a Killing vectorfield for $\tilde{\mathcal{M}}$, we have both the energy estimate and $h\pd_t$ estimate stated below.

\begin{proposition}\label{ws:classic_ee_prop}(Energy estimate for $\tilde{\mathcal{M}}$)
\begin{multline*}
\int_{H_{t_1}^{t_2}}(\pd_t\psi)^2+\int_{\tilde\Sigma_{t_2}}\left[\chi_H(\pd_r\psi)^2+(\pd_t\psi)^2+|\sla\nabla\psi|^2\right] \\
 \lesssim \int_{\tilde\Sigma_{t_1}}\left[\chi_H(\pd_r\psi)^2+(\pd_t\psi)^2+|\sla\nabla\psi|^2\right] + Err_\Box,
\end{multline*}
where $\chi_H=1-\frac{2M}r$ and
$$Err_\Box=\int_{t_1}^{t_2}\int_{\tilde\Sigma_t}|\pd_t\psi\Box_{\tilde{g}}\psi|.$$
\end{proposition}

\begin{proposition}\label{ws:hdt_prop}($h\pd_t$ estimate for $\tilde{\mathcal{M}}$)
Let $R>r_H+\delh$ be any given radius. Then for all $\epsilon>0$ and $p<2$, there is a small constant $c_\epsilon$ and a large constant $C_\epsilon$, such that
\begin{multline*}
\int_{H_{t_1}^{t_2}}(\pd_t\psi)^2+\int_{\tilde{\Sigma}_{t_2}}r^{p-2}\left[\chi_H(\pd_r\psi)^2+(\pd_t\psi)^2+|\sla\nabla\psi|^2\right] \\
+\int_{t_1}^{t_2}\int_{\tilde{\Sigma}_t\cap\{R+M<r\}}c_\epsilon r^{p-3}(\lbar\psi)^2 \\
\lesssim \int_{\tilde{\Sigma}_{t_1}}r^{p-2}\left[\chi_H(\pd_r\psi)^2+(\pd_t\psi)^2+|\sla\nabla\psi|^2\right]+Err,
\end{multline*}
where $\chi_H=1-\frac{2M}r$ and
\begin{align*}
Err&=Err_1+Err_\Box \\
Err_1&=\int_{t_1}^{t_2}\int_{\tilde{\Sigma}_t\cap\{R<r\}}\epsilon r^{-1}(L\psi)^2 \\
Err_\Box&=\int_{t_1}^{t_2}\int_{\tilde{\Sigma}_t}C_\epsilon r^{p-2}|\pd_t\psi\Box_{\tilde{g}}\psi|.
\end{align*}
\end{proposition}

\subsection{The Morawetz estimate for $\tilde{\mathcal{M}}$}

The key ingredient for the Morawetz estimate is the partial Morawetz estimate. Since this estimate is very sensitive to the geometry of the spacetime, we will prove it completely for $\tilde{\mathcal{M}}$.

\begin{proposition}\label{ws:partial_morawetz_prop}(Partial Morawetz estimate for $\tilde{\mathcal{M}}$)
There exist a vectorfield $X_0$ and a function $w_0$ such that the current $J[X_0,w_0]$ satisfies
$$\left[\frac{M^6}{r^7}\left(1-\frac{r_H}r\right)^2(\pd_r\psi)^2+\frac1r\left(1-\frac{r_{trap}}r\right)^2|\sla\nabla\psi|^2+\frac{1}{r^3}1_{r>4M}\psi^2\right]\lesssim divJ[X_0,w_0]$$
in Boyer-Lindquist coordinates.
\end{proposition}

\begin{proof}
For the proof of this proposition, we drop the $0$ subscript from $X_0$ and $w_0$. We begin by restating Lemma \ref{m:pme_initial_lem} for $\tilde{\mathcal{M}}$.
\begin{lemma}
With the choices $X=X^r(r)\pd_r$ and $w=w(r)$,
\begin{multline*}
divJ[X,w]= 
(wg^{tt}-r^{-6}\pd_r(r^6X^rg^{tt}))(\pd_t\psi)^2+(w-\pd_rX^r)|\sla\nabla\psi|^2 \\
+(wg^{rr}+\pd_rX^r g^{rr}-6r^{-1}X^rg^{rr}-X^r\pd_rg^{rr})(\pd_r\psi)^2 -\frac12r^{-6}\pd_r((r^6-2Mr^5)\pd_rw)\psi^2.
\end{multline*}
\end{lemma}
\begin{proof}
See the proof of Lemma \ref{m:pme_initial_lem}.
\end{proof}

It is difficult to make sense of the coefficients in the above lemma, because most coefficients are sums of multiple terms. By making the choices $X^r(r)=u(r)v(r)$ and $w(r)=v(r)\pd_ru(r)$, each coefficient can be written as a single term.

\begin{lemma}
With the choices $X=u(r)v(r)\pd_r$ and $w=v(r)\pd_ru(r)$,
\begin{multline*}
divJ[X,w]= 
-ur^{-6}\pd_r(r^6vg^{tt})(\pd_t\psi)^2-ur^{-4}\pd_r(r^4v) |\sla\nabla\psi|^2 \\
+r^6(g^{rr})^2u^{-1}\pd_r\left(\frac{u^2v}{r^6g^{rr}}\right)(\pd_r\psi)^2-\frac12r^{-6}\pd_r((r^6-2Mr^5)\pd_rw)\psi^2.
\end{multline*}
\end{lemma}
\begin{proof}
This follows directly by substituting $X^r=uv$ and $w=v\pd_ru$ into the expression in the previous lemma and combining terms.
\end{proof}

At this point, we choose a relation between $X$ and $w$ so that the coefficient of $(\pd_t\psi)^2$ vanishes. In terms of the functions $u$ and $v$, we set $v=r^{-6}(1-\frac{2M}r)$. 

\begin{lemma}
With the choice $v=r^{-6}(1-\frac{2M}r)$, the coefficient of the $(\pd_t\psi)^2$ term vanishes and furthermore,
\begin{multline*}
divJ[X,w]=\frac{2u}{r^7}\left(1-\frac{3M}r\right)|\sla\nabla\psi|^2 \\
+2\pd_r\left(\frac{u}{r^6}\right)\left(1-\frac{2M}r\right)^2(\pd_r\psi)^2-\frac12r^{-6}\pd_r((r^6-2Mr^5)\pd_rw)\psi^2.
\end{multline*}
Therefore, it is necessary to choose $u$ and $w$ so that the following conditions are satisfied.
\begin{eqnarray}
u\left(1-\frac{3M}r\right)\ge 0, \label{ws:u_cond_1_eqn}\\
\pd_r\left(\frac{u}{r^6}\right)\ge 0, \label{ws:u_cond_2_eqn}\\
\pd_r((r^6-2Mr^5)\pd_rw)\le 0, \label{ws:u_cond_3_eqn}
\end{eqnarray}
Furthermore, $w$ and $u$ must be related by the following constraint.
\begin{equation}
w=r^{-6}\left(1-\frac{2M}r\right)\pd_ru. \label{ws:u_constr_eqn}
\end{equation}
\end{lemma}

\begin{lemma}\label{ws:choice_for_u_and_w_lem}
It is possible to choose $u$ and $w$ such that all three conditions (\ref{ws:u_cond_1_eqn}-\ref{ws:u_cond_3_eqn}) and the constraint (\ref{ws:u_constr_eqn}) are satisfied. One particular choice is given by the following approach.

i) Require that $u(3M)=0$ and $\pd_ru=r^6\left(1-\frac{2M}r\right)^{-1}w$. This will specify $u$ completely in terms of $w$ and satisfy the constraint (\ref{ws:u_constr_eqn}).

ii) Require that $w$ be positive. Then $u$ will be an increasing function and the condition (\ref{ws:u_cond_1_eqn}) will be satisfied.

iii) Determine a way to understand $\pd_r\left(\frac{u}{r^6}\right)$ in terms of $w$. Here is one way. Let $\tilde{K}^{rr}=\frac16 r^{7}\pd_r\left(\frac{u}{r^6}\right)$. This particular quantity has the same sign as $\pd_r\left(\frac{u}{r^6}\right)$ and has the property that its derivative can be written as a function of $w$. That is, $\pd_r\tilde{K}^{rr}=r^6\pd_r\left(\frac{w}{6r^5v}\right)$, where again, $v=r^{-6}(1-\frac{2M}r)$. The condition (\ref{ws:u_cond_2_eqn}) now reduces to showing that $\tilde{K}^{rr}\ge 0$.

iv) For sufficiently large $r$, impose the condition that $\pd_r\tilde{K}^{rr}=0$. In particular, this means choosing $w=6r^5v$ for sufficiently large $r$. The quantity $6r^5v$ has a maximum at $r_*=4M$. The choice $w=6r^5v$ for $r\ge r_*$ will satisfy both the remaining conditions (\ref{ws:u_cond_2_eqn}-\ref{ws:u_cond_3_eqn}) for $r>r_*$.

v) Observe that since $\pd_r(6r^5v)=0$ at $r_*=4M$, the quantity $(r^6-2Mr^5)\pd_rw=0$ will necessarily vanish at $r=2M$ and $r=4M$ if $w=6r^5v$ for $r\ge r_*$. By the mean value theorem, condition (\ref{ws:u_cond_3_eqn}) necessitates that $\pd_rw=0$ entirely for $2M\le r\le 4M$. Thus, take $w=w(r_*)$ for $r\le r_*$.

iv) Since $w$ has now been chosen for all $r$, check that $\tilde{K}^{rr}\ge 0$. In particular, with the given choice of $w$, $\tilde{K}^{rr}$ will be decreasing for $r<r_*$ and then remain constant for $r\ge r_*$. Thus, it suffices to compute that $\tilde{K}^{rr}(r_*)>0$.
\end{lemma}

\begin{proof}
The lemma consists of multiple statements. Statements (i) and (ii) require no further justification. We turn to statement (iii). Note that
\begin{align*}
\pd_r\tilde{K}^{rr} &= \pd_r\left(\frac16r^7\pd_r\left(\frac{u}{r^6}\right)\right) \\
&= \pd_r\left(\frac16r^7\left(\frac{u'}{r^6}-\frac{6u}{r^7}\right)\right) \\
&= \pd_r\left(\frac{ru'}{6}-u\right) \\
&= \frac{ru''}{6}-\frac{5u'}6 \\
&= r^6\pd_r\left(\frac{u'}{6r^5}\right).
\end{align*}
Recall the constraint (\ref{ws:u_constr_eqn}) between $w$ and $u$, which stipulates $w=vu'$. Therefore, we substitute $u'=\frac{w}{v}$ and conclude
$$\pd_r\tilde{K}^{rr} = r^6\pd_r\left(\frac{w}{6r^5v}\right).$$
This verifies statement (iii). For statement (iv), we must find the value $r_*$ that maximizes the quantity $6r^5v$. Note that
$$6r^5v = 6r^{-1}\left(1-\frac{2M}r\right).$$
In the proof of Lemma \ref{s:choice_for_u_and_w_lem}, we already calculated
$$\pd_r\left(r^{-1}\left(1-\frac{2M}r\right)\right) = -\frac{1}{r^2}\left(1-\frac{4M}r\right).$$
Setting $\pd_r(6r^5v)=0$, we deduce $r_*=4M$. This indeed maximizes the quantity $6r^5v$, because the quantity $6r^5v$ vanishes at the event horizon and decays back to zero as $r\rightarrow\infty$.

We also must check that $-\frac12r^{-6}\pd_r\left((r^6-2Mr^5)\pd_r(6r^5v)\right)\ge 0$ for $r\ge r_*=4M$. Indeed,
\begin{align*}
-\frac12 r^{-6}\pd_r\left((r^6-2Mr^5)\pd_r(6r^5v)\right) &= 3r^{-6}\pd_r\left(\frac{r^6-2Mr^5}{r^2}\left(1-\frac{4M}r\right)\right) \\
&= 3r^{-6}\pd_r\left(r^4\left(1-\frac{2M}r\right)\left(1-\frac{4M}r\right)\right).
\end{align*}
The last quantity is clearly positive for $r\ge 4M$. This verifies statement (iv). Statement (v) requires no further justification. Finally, we turn to statement (vi). First, we observe that for $r<r_*$, since $w=w(r_*)$ is constant, then
$$\pd_r\tilde{K}^{rr}=r^6\pd_r\left(\frac{w(r_*)}{6r^5v}\right)=-\frac{w(r_*)r^6}{(6r^5v)^2}\pd_r(6r^5v).$$
Since we have determined that $6r^5v$ has a maximum at $r_*=4M$, it follows that $\tilde{K}^{rr}$ decreases to a minimum at $r_*=4M$. Recall further that by the choice of $w$ for $r\ge r_*$, $\tilde{K}^{rr}$ is constant for $r\ge r_*$. It suffices to check the value of $\tilde{K}^{rr}$ at $r=r_*$. (That is, if $\tilde{K}^{rr}(r_*)>0$, then necessarily $\tilde{K}^{rr}>0$ and thus $K^{rr}\ge 0$.) Note that
\begin{multline*}
\left.\frac{r^6 v}{w}\pd_r\left(\frac{u}{r^6}\right)\right|_{r_*}
=\left.\frac{r^6v\pd_ru}{wr^6}-\frac{6r^6uv}{wr^7}\right|_{r_*}
=1-\frac{6v(r_*)}{r_*w(r_*)}\int_{r_{trap}}^{r_*}\pd_r u 
=1-\frac{6v(r_*)}{r_*w(r_*)}\int_{r_{trap}}^{r_*}\frac{w}{v} \\
=1-\frac{6v(r_*)}{r_*}\int_{r_{trap}}^{r_*}\frac1v 
=1-\frac{6r_*^4v(r_*)}{r_*^5}\int_{r_{trap}}^{r_*}\frac{r^4}{r^4v} 
>1-\frac{6}{r_*^5}\int_{r_{trap}}^{r_*}r^4
\end{multline*}
where we used the fact that $w$ is constant for $r\in [r_{trap},r_*]$ and in the final step we used the fact that $r^4v$ is decreasing from $r_{trap}$ to $r_*$. Finally, we calculate
$$1-\frac{6}{r_*^5}\int_{r_{trap}}^{r_*}r^4 =1-\frac{6(r_*^5-r_{trap}^5)}{5r_*^5} = 1-\frac65\left(1-\left(\frac{3M}{4M}\right)^5\right) = \frac{217}{2560} > 0.$$
This concludes the proof of Lemma \ref{ws:choice_for_u_and_w_lem}.
\end{proof}

We conclude the proof of Proposition \ref{ws:partial_morawetz_prop} by analyzing the asymptotics for large $r$, in particular for $r>r_*$.
\begin{align*}
v &= r^{-6}\left(1-\frac{2M}r\right)=O(r^{-6}) \\
\pd_rv &= O(r^{-7}) \\
w_0 &= 6r^5v = O(r^{-1}) \\
u &= r^6-c^6 \text{ since }u'=6r^5 \\
X_0 &= uv\pd_r = \left(1-\frac{2M}r\right)\pd_r+O\left(\frac{M^6}{r^6}\right)\pd_r \\
\sla{K}^{\alpha\beta} &= -ur^{-4}\pd_r(r^4v)\sla{g}^{\alpha\beta}=O(r^{-1})\sla{g}^{\alpha\beta} \\
K^{rr} &=O(r^{-7})\tilde{K}^{rr}=O\left(\frac{M^6}{r^7}\right) \\
K &=r^{-6}\pd_r\left(r^4\left(1-\frac{2M}r\right)\left(1-\frac{4M}r\right)\right) = O\left(r^{-3}\right).
\end{align*}
This completes the proof of Proposition \ref{ws:partial_morawetz_prop}.
\end{proof}

With the partial Morawetz estimate having been proved, the Morawetz estimate now follows with no serious additional difficulty.

\begin{theorem}\label{ws:morawetz_thm} (Morawetz estimate for $\tilde{\mathcal{M}}$)
The following estimate holds for any sufficiently regular function $\psi$ decaying sufficiently fast as $r\rightarrow\infty$.
\begin{multline*}
\int_{H_{t_1}^{t_2}}|\sla\nabla\psi|^2+\int_{\tilde{\Sigma}_{t_2}}(L\psi)^2+|\sla\nabla\psi|^2+r^{-2}\psi^2 \\
+\int_{t_1}^{t_2}\int_{\tilde{\Sigma}_t}\frac{M^6}{r^7}(\pd_r\psi)^2+\frac1r\left(1-\frac{3M}r\right)^2\left[|\sla\nabla\psi|^2+\frac{M^6}{r^6}(\pd_t\psi)^2\right]+\frac{1}{r^3}\psi^2 \\
\lesssim \int_{\tilde{\Sigma}_{t_1}}(L\psi)^2+|\sla\nabla\psi|^2+r^{-2}\psi^2 + Err,
\end{multline*}
where
\begin{align*}
Err &= Err_1+Err_2+Err_\Box \\
Err_1 &= \int_{\tilde{\Sigma}_{t_2}} r^{-1}|\psi L\psi| \\
Err_2 &= \int_{H_{t_1}^{t_2}}(\pd_t\psi)^2+\int_{\tilde{\Sigma}_{t_2}}\frac{M^6}{r^6}\left[\chi_H(\pd_r\psi)^2+(\pd_t\psi)^2\right] \\
Err_\Box&=\int_{t_1}^{t_2}\int_{\tilde{\Sigma}_t}|(2X(\psi)+w\psi)\Box_g\psi|,
\end{align*}
and $\chi_H=1-\frac{2M}r$.
\end{theorem}
\begin{proof}
See the proof of Theorem \ref{s:morawetz_thm}, keeping in mind Proposition \ref{ws:partial_morawetz_prop}.
\end{proof}

\subsection{The $r^p$ estimate for $\tilde{\mathcal{M}}$}

The key ingredient for the $r^p$ estimate is the pre-$r^p$ identity. This identity is also sensitive to the geometry of the spacetime, so we will prove it completely for $\tilde{\mathcal{M}}$.

\begin{lemma}\label{ws:p_identity_psi_lem} (Pre-$r^p$ identity for $\tilde{\mathcal{M}}$) Let $\alpha=1-\frac{2M}r$ and $L=\alpha\pd_r+\pd_t$. For any function $f=f(r)$ supported where $r>r_H+\delh$, the following identity holds.
\begin{align*}
&\int_{\tilde{\Sigma}_{t_2}}f\left(\alpha^{-1}L\psi+\frac{3}r\psi\right)^2 +\alpha^{-1}f|\tsla\nabla\psi|^2 +6f\frac{\psi^2}{r^2} -\frac1{A^2r^2}\pd_r(A^23rf\psi^2) \\
&+\int_{t_1}^{t_2}\int_{\tilde{\Sigma}_t}
\left[
  \left(2r^{-1}f-f'\right)|\tsla\nabla\psi|^2
+ \alpha f'\left(\alpha^{-1}L\psi+\frac3r\psi\right)^2 
\vphantom{
+ \epsilon\alpha\left((1-\epsilon)f'-rf''\right)\frac{\phi^2}{q^2} 
+ \alpha'\left(\frac{r^2-a^2}{r^2+a^2}f-\epsilon r f'\right)\frac{\phi^2}{q^2}
+ \frac{a^2}{r^2+a^2}\left(\alpha((1-\epsilon)^2-2)f'-\frac{4\alpha rf}{r^2+a^2}\right)\frac{\phi^2}{q^2}
}\right. \\
&\hspace{1.2in}
\left.\vphantom{
  \left(\frac{2rf}{r^2+a^2}-f'\right)\frac{Q^{\alpha\beta}}{q^2}\pd_\alpha\phi\pd_\beta\phi
+ \alpha f'\frac{r^2+a^2}{q^2}\left(\alpha^{-1}L\phi+\frac{(1-\epsilon)r}{r^2+a^2}\phi\right)^2 
+ \epsilon\alpha\left((1-\epsilon)f'-rf''\right)\frac{\phi^2}{q^2} 
}
+ 6\alpha(2r^{-1}f-f')\frac{\psi^2}{r^2} -\alpha'\alpha^{-2}f(L\psi)^2 + 3\alpha'f\frac{\psi^2}{r^2}
\right] \\
=&\int_{\tilde{\Sigma}_{t_1}} f\left(\alpha^{-1}L\psi+\frac{3}r\psi\right)^2 +\alpha^{-1}f|\tsla\nabla\psi|^2 +6f\frac{\psi^2}{r^2} -\frac1{A^2r^2}\pd_r(A^23rf\psi^2) \\
&+\int_{t_1}^{t_2}\int_{\tilde{\Sigma}_t}-\left(2\alpha^{-1}fL\psi+6r^{-1}f\psi\right)\Box_{\tilde{g}}\psi.
\end{align*}
\end{lemma}

\begin{proof}
The proof is similar to the proof of Lemma \ref{s:p_ee_identity_lem}, however there are a few subtle differences, some of which actually simplify the proof. Again, we will use Proposition \ref{general_divergence_estimate_prop} together with the following current template.
$$J[X,w,m]_\mu = T_{\mu\nu} X^\nu +w\psi\pd_\mu\psi-\frac12\psi^2\pd_\mu w+m_\mu \psi^2,$$
$$T_{\mu\nu}=2\pd_\mu\psi\pd_\nu\psi-g_{\mu\nu}\pd^\lambda\psi\pd_\lambda\psi.$$

Assume for now that $\Box_{\tilde{g}}\psi=0$. Let $\alpha = 1-\frac{2M}r$, and observe that
$$L=\alpha\pd_r+\pd_t,$$
$$g^{rr}=\alpha,$$
$$g^{tt}=-\alpha^{-1}.$$

\begin{lemma}\label{ws:divJpsiX_lem}
Without appealing directly to the particular expression for $\alpha$, one can deduce the following.
\begin{multline*}
divJ[\alpha^{-1}fL]=(\alpha^{-1}f)'(L\psi)^2-6r^{-1}f\left(\alpha(\pd_r\psi)^2-\alpha^{-1}(\pd_t\psi)^2\right) \\
-\left(4r^{-1}f+f'\right)|\tsla\nabla\psi|^2.
\end{multline*}
\end{lemma}
\begin{proof}
Note that
$$div J[X] = K^{\mu\nu}\pd_\mu\psi\pd_\nu\psi,$$
where
$$K^{\mu\nu}=2g^{\mu\lambda}\pd_\lambda X^\nu-X^\lambda \pd_\lambda(g^{\mu\nu})-div X g^{\mu\nu}.$$

Set $X=\alpha^{-1}f(\alpha\pd_r+\pd_t)=f\pd_r+\alpha^{-1}f\pd_t$. From the above formula, since $g^{rt}=0$,
$$K^{tr}+K^{rt}=2g^{rr}\pd_r X^t=2\alpha\pd_r(\alpha^{-1} f).$$
Thus, the expression for $divJ[\alpha^{-1}fL]$ will have a mixed term of the form 
$$2\alpha\pd_r(\alpha^{-1}f)\pd_r\psi\pd_t\psi.$$
Note that
\begin{align*}
(\alpha^{-1} f)'(L\psi)^2 &= (\alpha^{-1} f)'(\alpha\pd_r\psi+\pd_t\psi)^2 \\
&= \alpha^2(\alpha^{-1}f)'(\pd_r\psi)^2+2\alpha(\alpha^{-1}f)'\pd_r\psi\pd_t\psi +(\alpha^{-1} f)'(\pd_t\psi)^2.
\end{align*}
We now compute the $(\pd_r\psi)^2$ and $(\pd_t\psi)^2$ components, subtracting the part that will be grouped with the $(L\psi)^2$ term.
\begin{align*}
K^{rr}-\alpha^2(\alpha^{-1}f)' &=2g^{rr}\pd_rX^r-X^r\pd_r g^{rr}-\frac1{A^2r^2}\pd_r(A^2r^2X^r)g^{rr}-\alpha^2(\alpha^{-1}f)' \\
&= 2g^{rr}\pd_rX^r-\frac1{A^2r^2}\pd_r(A^2r^2g^{rr}X^r)-\alpha^2(\alpha^{-1}f)' \\
&= 2\alpha\pd_r f -\frac{1}{r^6}\pd_r\left(r^6\alpha f\right)-\alpha^2(\alpha^{-1}f)' \\
&= -6\alpha r^{-1}f
\end{align*}
and
\begin{align*}
K^{tt}-(\alpha^{-1}f)' &= -X^r\pd_r g^{tt}-\frac{1}{A^2r^2}\pd_r(A^2r^2X^r)g^{tt}-(\alpha^{-1}f)' \\
&= -\frac{1}{A^2r^2}\pd_r(A^2r^2 g^{tt} X^r) -(\alpha^{-1}f)' \\
&= -\frac{1}{r^6}\pd_r\left(r^6(-\alpha^{-1}) f\right)-(\alpha^{-1}f)' \\
&= 6\alpha^{-1}r^{-1}f.
\end{align*}
Finally,
\begin{align*}
\sla{K}^{\alpha\beta} &= -X^r\pd_r \sla{g}^{\alpha\beta}-\frac1{A^2r^2}\pd_r\left(A^2r^2X^r\right)\sla{g}^{\alpha\beta} \\
&= -\frac{1}{A^2r^2}\pd_r\left(A^2r^2\sla{g}^{\alpha\beta}X^r\right) \\
&= -\frac{1}{r^6}\pd_r(r^6\sla{g}^{\alpha\beta}f) \\
&= -\frac{1}{r^6}\pd_r(r^4f(r^2\sla{g}^{\alpha\beta})) \\
&= -\frac{1}{r^6}\pd_r(r^4f)(r^2\sla{g}^{\alpha\beta}) \\
&= -\left(4r^{-1}f+f'\right) \sla{g}^{\alpha\beta}
\end{align*}
Combining all these terms gives the identity stated in the lemma. 
\end{proof}

Next, we choose $w=6r^{-1}f$ to directly cancel with the middle term in the above lemma.
\begin{lemma}\label{ws:divJpsiXw_lem}
$$divJ\left[\alpha^{-1}fL,6r^{-1}f\right] = (\alpha^{-1}f)'(L\psi)^2+\left(2r^{-1}f-f'\right)|\tsla\nabla\psi|^2-\frac12\Box_{\tilde{g}}\left(6r^{-1}f\right)\psi^2.$$
\end{lemma}
\begin{proof}
Note that
$$divJ[0,w]=wg^{\mu\nu}\pd_\mu\psi\pd_\nu\psi-\frac12\Box_{\tilde{g}}w \psi^2.$$
We compute the new terms only.
\begin{multline*}
divJ\left[0,6r^{-1}f\right] = 6r^{-1}fg^{\alpha\beta}\pd_\alpha\psi\pd_\beta\psi -\frac12\Box_{\tilde{g}}\left(6r^{-1}f\right)\psi^2 \\
= 6r^{-1}f\left(\alpha (\pd_r\psi)^2-\alpha^{-1}(\pd_t\psi)^2\right) +6r^{-1}f|\tsla\nabla\psi|^2 -\frac12\Box_{\tilde{g}}\left(6r^{-1}f\right)\psi^2.
\end{multline*}
When adding these terms to the expression in Lemma \ref{ws:divJpsiX_lem}, the $\alpha(\pd_r\psi)^2-\alpha^{-1}(\pd_t\psi)^2$ terms cancel (this was the reason for the choice of $w=6r^{-1}f$) and the result is as desired. 
\end{proof}

The term $-\frac12\Box_{\tilde{g}}\left(6r^{-1}f\right)\psi^2$ is like $-6r^{-1}f''-24r^{-1}\pd_r(r^{-1}f)\phi^2$. In the future, when $f\sim r^p$, this will have a sign $-p^2-3p+4=-(p+4)(p-1)$. The sign will be negative if $p>1$, which is bad. So we include a divergence term to fix it. (Unlike in the analogous step for the proof of Lemma \ref{s:p_ee_identity_lem}, there will not be a need for a smallness parameter $\epsilon$.) This is the point of the following Lemma.
\begin{lemma}
\begin{multline*}
\alpha^{-1}f'(L\psi)^2+\left[-\frac12\Box_{\tilde{g}}\left(6r^{-1}f\right)\psi^2+div\left(\psi^23r^{-1}f'L\right)\right] \\
= \alpha f'\left(\alpha^{-1}L\psi+\frac{3}{r}\psi\right)^2+6\alpha(2r^{-1}f-f')\frac{\psi^2}{r^2}
+3\alpha' r^{-2}f\psi^2
\end{multline*}
\end{lemma}
\begin{proof}
First, borrowing from a calculation in the proof of Lemma \ref{s:rp_add_term_lem}, we obtain
\begin{align*}
-\frac12\Box_{\tilde{g}}\left(6r^{-1}f\right)\psi^2 
&= -\frac{1}{A^2r^2}\pd_r\left(A^2r^2\alpha\pd_r\left(3r^{-1}f\right)\right)\psi^2 \\
&= -\Box_g\left(3r^{-1}f\right)\psi^2-\frac{\pd_rA^2}{A^2}\alpha\pd_r\left(3r^{-1}f\right)\psi^2 \\
&= -3\alpha r^{-1}f''\psi^2-\frac{\pd_rA^2}{A^2}\alpha\pd_r\left(3r^{-1}f\right)\psi^2 -\alpha'\pd_r\left(3r^{-1}f\right)\psi^2.
\end{align*}
We also calculate
\begin{multline*}
div\left(\psi^23r^{-1}f'L\right) = \frac{1}{A^2r^2}\pd_\alpha(A^2\psi^2 3rf'L^\alpha) \\
=3r^{-1}f'2\psi L\psi + r^{-2}\pd_r(3r f'\alpha)\psi^2+\frac{\pd_rA^2}{A^2}3\alpha r^{-1}f'\psi^2 \\
= 3r^{-1}f'2\psi L\psi +3\alpha r^{-2}f'\psi^2 + 3\alpha r^{-1}f''\psi^2 +\frac{\pd_rA^2}{A^2}3\alpha r^{-1}f'\psi^2 +3\alpha' r^{-1}f'\psi^2  \\
= 3r^{-1}f'2\psi L\psi +9\alpha r^{-2} f'\psi^2 + 3\alpha r^{-1}f''\psi^2 -6\alpha r^{-2}f'\psi^2+\frac{\pd_rA^2}{A^2}3\alpha r^{-1}f'\psi^2 +3\alpha'r^{-1} f'\psi^2.
\end{multline*}
The first two terms in the last line complete a square with the term $\alpha^{-1}f'(L\psi)^2$. The third term cancels with the first term from the previous calculation. Due to the fourth and fifth terms, which did not show in the calculation in the proof of Lemma \ref{s:rp_add_term_lem}, the $\epsilon$ parameter is not needed here, allowing for a slightly simpler calculation.
$$\alpha^{-1}f'(L\psi)^2+3r^{-1}f'2\psi L\psi +9\alpha r^{-2}f'\psi^2=\alpha f'\left(\alpha^{-1}L\psi+\frac{3}{r}\psi\right)^2$$
and
$$-3\alpha r^{-1}f''\psi^2+3\alpha r^{-1}f''\psi^2=0.$$
The new terms are
\begin{align*}
-\frac{\pd_rA^2}{A^2}&\alpha\pd_r\left(3r^{-1}f\right)\psi^2 -6\alpha r^{-2} f'\psi^2+\frac{\pd_rA^2}{A^2}3\alpha r^{-1} f'\psi^2 \\
&= -6\alpha r^{-2} f'\psi^2 -3\frac{\pd_r A^2}{A^2} \pd_r\left(r^{-1}\right)\alpha f \psi^2 \\
&= -6\alpha r^{-2} f'\psi^2 -3\left(4r^{-1}\right)\left(-r^{-2}\right)\alpha f \psi^2 \\
&= -6\alpha r^{-2}f'\psi^2+12\alpha r^{-3}f\psi^2
\end{align*}
Adding these terms together yields
$$\alpha f'\left(\alpha^{-1}L\psi+\frac{3}{r}\psi\right)^2+6\alpha(2r^{-1}f-f')\frac{\psi^2}{r^2}.$$
All the remaining terms (which contain a factor of $\alpha'\sim \frac{M}{r^2}$) are
\begin{align*}
-\alpha'&\pd_r\left(3r^{-1}f\right)\psi^2 +3\alpha' r^{-1} f'\psi^2 \\
&= -3\alpha'\pd_r(r^{-1})f\psi^2 \\
&= 3\alpha' r^{-2}f \psi^2.
\end{align*}
Adding both of these yields the result. 
\end{proof}

Thus, we have shown that if $\Box_{\tilde{g}}\psi = 0$, then
\begin{multline*}
div J\left[\alpha^{-1}fL,6r^{-1}f,3r^{-1}f'L\right] \\
= \alpha f'\left(\alpha^{-1}L\psi+\frac{3}{r}\psi\right)^2+6\alpha(2r^{-1}f-f')\frac{\psi^2}{r^2} +\left(2r^{-1}f-f'\right)|\tsla\nabla\psi|^2 \\
-\alpha'\alpha^{-2}f(L\psi)^2+3\alpha'f\frac{\psi^2}{r^2}.
\end{multline*}
If we remove the assumption that $\Box_{\tilde{g}}\psi=0$, there is an additional term
$$(2X(\psi)+w\psi)\Box_{\tilde{g}}\psi =\left(2\alpha^{-1}fL\psi+6r^{-1}f\psi\right)\Box_{\tilde{g}}\psi$$
appearing in the expression for $div J$.

Finally, we turn to the boundary terms. Since we have assumed that $f$ is supported away from the event horizon, it suffices to compute $-J^t$.
\begin{lemma}
\begin{multline*}
-J^t\left[\alpha^{-1}fL,6r^{-1}f,3r^{-1}f'L\right] \\
=f\left(\alpha^{-1}L\psi+\frac{3}{r}\psi\right)^2 +\alpha^{-1}f|\tsla\nabla\psi|^2 +6f\frac{\psi^2}{r^2} -\frac1{A^2r^2}\pd_r(A^23rf\psi^2).
\end{multline*}
\end{lemma}
\begin{proof}
Borrowing a calculation from the proof of Lemma \ref{s:rp_boundary_terms_lem}, we have
$$-J^t[\alpha^{-1}fL] = \alpha^{-2}f(L\psi)^2+\alpha^{-1}f|\tsla\nabla\psi|^2.$$
Borrowing another calculation from the proof of Lemma \ref{s:rp_boundary_terms_lem}, we have
\begin{align*}
-J^t &\left[0,6r^{-1}f\right] \\
&= 6\alpha^{-1}r^{-1}f\psi L\psi+3f\frac{\psi^2}{r^2} +3r^{-1}f'\psi^2-\frac{3}{r^2}\pd_r(rf\psi^2) \\
&= 6\alpha^{-1}r^{-1}f\psi L\psi+9f\frac{\psi^2}{r^2} -6f\frac{\psi^2}{r^2}+3r^{-1}f'\psi^2+\frac{\pd_rA^2}{A^2}3rf\frac{\psi^2}{r^2}-\frac{1}{A^2r^2}\pd_r(A^23rf\psi^2).
\end{align*}
Following a similar procedure as in the proof of Lemma \ref{s:rp_boundary_terms_lem}, we notice that
$$\alpha^{-2}f(L\psi)^2+6\alpha^{-1}r^{-1}f\psi L\psi+9f\frac{\psi^2}{r^2} =f\left(\alpha^{-1}L\psi+\frac{3}{r}\psi\right)^2.$$
Also, there are two new terms, which we now combine.
$$-6f\frac{\psi^2}{r^2}+\frac{\pd_rA^2}{A^2}3rf\frac{\psi^2}{r^2} = -6f\frac{\psi^2}{r^2} +12f\frac{\psi^2}{r^2} = 6f\frac{\psi^2}{r^2}.$$
Thus,
\begin{multline*}
-J^t\left[\alpha^{-1}fL,6r^{-1}f,0\right] \\
=f\left(\alpha^{-1}L\psi+\frac{3}{r}\psi\right)^2 
+\alpha^{-1}f|\tsla\nabla\psi|^2+6f\frac{\psi^2}{r^2}+3r^{-1}f'\psi^2-\frac1{A^2r^2}\pd_r(A^23rf\psi^2).
\end{multline*}
Finally,
$$-J^t\left[0,0,3r^{-1}f'L\right]=-3r^{-1}f'\psi^2L^t=-3r^{-1}f'\psi^2.$$
Adding these two expressions together yields the result. 
\end{proof}

This concludes the proof of Lemma \ref{ws:p_identity_psi_lem}.
\end{proof}

Now, we briefly observe that the coefficients for each of the terms in Lemma \ref{ws:p_identity_psi_lem} are very similar to the coefficients for their corresponding terms in Lemma \ref{s:p_ee_identity_lem} with the one exception that the terms with a factor of $\epsilon$ in Lemma \ref{s:p_ee_identity_lem} are now replaced with the stronger terms $6f\frac{\psi^2}{r^2}$ and $6\alpha(2r^{-1}f-f')\frac{\psi^2}{r^2}$. From these facts, we conclude the incomplete $r^p$ estimate near $i^0$.

\begin{proposition}\label{ws:incomplete_p_estimate_prop}
Fix $\delm,\delp>0$. Let $R$ be a sufficiently large radius. Then for all $p\in[\delm,2-\delp]$, the following estimate holds if $\psi$ decays sufficiently fast as $r\rightarrow\infty$.
\begin{multline*}
\int_{\tilde\Sigma_{t_2}\cap\{r>2R\}}r^p\left[(L\psi)^2+|\sla\nabla\psi|^2+r^{-2}\psi^2\right] \\
+ \int_{t_1}^{t_2}\int_{\tilde\Sigma_t\cap\{r>2R\}}r^{p-1}\left[(L\psi)^2+|\sla\nabla\psi|^2+r^{-2}\psi^2\right] \\
\lesssim \int_{\tilde\Sigma_{t_2}\cap\{r>2R\}}r^p\left[(L\psi)^2+|\sla\nabla\psi|^2+r^{-2}\psi^2\right] + Err,
\end{multline*}
where
\begin{align*}
Err &= Err_1 + Err_\Box \\
Err_1 &= \int_{t_1}^{t_2}\int_{\tilde\Sigma_t\cap\{R<r<2R\}}(L\psi)^2+|\sla\nabla\psi|^2+\psi^2  \\
Err_\Box &= \int_{t_1}^{t_2}\int_{\tilde\Sigma_t\cap\{R<r\}}r^p(|L\psi|+r^{-1}|\psi|)|\Box_{\tilde{g}}\psi|.
\end{align*}
\end{proposition}
\begin{proof}
See the proof of Proposition \ref{ws:incomplete_p_estimate_prop}, keeping in mind Lemma \ref{ws:p_identity_psi_lem}.
\end{proof}

From the incomplete $r^p$ estimate near $i^0$, the Morawetz estimate, and the $h\pd_t$ estimate, we conclude the $r^p$ estimate.

\begin{proposition}\label{ws:rp_prop}
Fix $\delm,\delp>0$ and let $p\in[\delm,2-\delp]$. Then if
$\psi$ decays sufficiently fast as $r\rightarrow\infty$, the following estimate holds.
\begin{multline*}
\int_{\tilde\Sigma_{t_2}}r^p\left[(L\psi)^2+|\sla\nabla\psi|^2+r^{-2}\psi^2 + r^{-2}(\pd_r\psi)^2\right] \\
+ \int_{t_1}^{t_2}\int_{\tilde\Sigma_t}r^{p-1}\left[\chi_{trap}(L\psi)^2+\chi_{trap}|\sla\nabla\psi|^2+r^{-2}\psi^2+r^{-2}(\pd_r\psi)^2\right] \\
\lesssim \int_{\tilde\Sigma_{t_1}}r^p\left[(L\psi)^2+|\sla\nabla\psi|^2+r^{-2}\psi^2+r^{-2}(\pd_r\psi)^2\right] 
+Err_\Box,
\end{multline*}
where $\chi_{trap}=\left(1-\frac{r_{trap}}{r}\right)^2$ and
\begin{align*}
Err_\Box &= \int_{t_1}^{t_2}\int_{\tilde\Sigma_t}|(2X(\psi)+w\psi)\Box_{\tilde{g}}\psi|,
\end{align*}
where the vectorfield $X$ and function $w$ satisfy the following properties. \\
\bp $X$ is everywhere timelike, but asymptotically null at the rate $X=O(r^p)L+O(r^{p-2})\pd_t$. \\
\bp $X|_{r=r_H}=-\lambda\pd_r$ for some positive constant $\lambda$. \\
\bp $X|_{r=r_{trap}}=\lambda\pd_t$ for some positive constant $\lambda$. \\
\bp $w =O(r^{p-1})$ for large $r$.
\end{proposition}
\begin{proof}
See the proof of Proposition \ref{s:rp_prop}, keeping in mind Proposition \ref{ws:hdt_prop}, Theorem \ref{ws:morawetz_thm}, and Proposition \ref{ws:incomplete_p_estimate_prop}.
\end{proof}

\section{The dynamic estimates}\label{ws:dynamic_estimates_sec}

\subsection{The dynamic estimates for $(\phi,\psi)$ ($s,k=0$)}

We state the dynamic estimates for the wavefunctions $(\phi,\psi)$ only, which follow from the dynamic estimates in Schwarzschild as well as from the estimates derived in the previous section.

\begin{proposition}\label{ws:dynamic_estimates_0_prop}
Fix $\delm,\delp>0$. The following estimates hold for $p\in [\delm,2-\delp]$.
$$E(t_2)\lesssim E(t_1)+\int_{t_1}^{t_2}N(t)dt,$$
$$E_p(t_2)+\int_{t_1}^{t_2}B_p(t)dt\lesssim E_p(t_1)+\int_{t_1}^{t_2}N_p(t)dt,$$
where
\begin{align*}
E(t)=&\int_{\Sigma_t}\chi_H(\pd_r\phi)^2+(\pd_t\phi)^2+|\sla\nabla\phi|^2 \\
&+\int_{\tilde\Sigma_t}\chi_H(\pd_r\psi)^2+(\pd_t\psi)^2+|\tsla\nabla\psi|^2,
\end{align*}
\begin{align*}
E_p(t)=&\int_{\Sigma_t}r^p\left[(L\phi)^2+|\sla\nabla\phi|^2+r^{-2}\phi^2+r^{-2}(\pd_r\phi)^2\right] \\
&+\int_{\tilde\Sigma_t}r^p\left[(L\psi)^2+|\tsla\nabla\psi|^2+r^{-2}\psi^2+r^{-2}(\pd_r\psi)^2\right],
\end{align*}
\begin{align*}
B_p(t)=&\int_{\Sigma_t}r^{p-1}\left[\chi_{trap}(L\phi)^2+\chi_{trap}|\sla\nabla\phi|^2+r^{-2}\phi^2+r^{-2}(\pd_r\phi)^2\right] \\
&+\int_{\tilde\Sigma_t}r^{p-1}\left[\chi_{trap}(L\psi)^2+\chi_{trap}|\tsla\nabla\psi|^2+r^{-2}\psi^2+r^{-2}(\pd_r\psi)^2\right],
\end{align*}
$$N(t)=(E(t))^{1/2}\left(||\Box_g\phi||_{L^2(\Sigma_t)}+||\Box_{\tilde{g}}\psi||_{L^2(\tilde{\Sigma}_t)}\right),$$
\begin{align*}
N_p(t) =& \int_{\Sigma_t}r^{p+1}(\Box_g\phi)^2 + \int_{\tilde\Sigma_t}r^{p+1}(\Box_{\tilde{g}}\psi)^2 \\
&+\int_{\Sigma_t\cap\{r\approx r_{trap}\}}|\pd_t\phi\Box_g\phi|+\int_{\tilde\Sigma_t\cap\{r\approx r_{trap}\}}|\pd_t\psi\Box_{\tilde{g}}\psi|
\end{align*}
and $\chi_H=1-\frac{2M}r$, $\chi_{trap}=\left(1-\frac{r_{trap}}r\right)^2$.
\end{proposition}
\begin{proof}
If $\psi=0$, this is exactly Corollary \ref{s:dynamic_estimates_0_cor}. If $\psi\ne 0$, then by the additional estimates in \S\ref{ws:spacetime_estimates_sec}, this proposition can be proved the same way as Corollary \ref{s:dynamic_estimates_0_cor}.
\end{proof}

\subsection{Commutators with $\Box_g$ and commutators with $\Box_{\tilde{g}}$}

We have a few options for commutators with $\Box_g$ and commutators with $\Box_{\tilde{g}}$. First, there is the operator $\pd_t$, which is a commutator with both $\Box_g$ and $\Box_{\tilde{g}}$. Then there are rotation operators. The three rotation operators from Chapter \ref{szd_chap}, which commute with $\Box_g$, are analogous to seven rotation operators that commute with $\Box_{\tilde{g}}$. However, these rotation operators destroy the axisymmetry that makes it possible to easily translate between both spacetimes. Instead of using the rotation operators, we will commute with the second order operators $r^2\sla\triangle$ and $r^2\tsla\triangle$, which we call $Q$ and $\tilde{Q}$ respectively.

\begin{definition}
\begin{align*}
Q &= r^2\sla\triangle \\
\tilde{Q} &= r^2\tsla\triangle.
\end{align*}
\end{definition}

\begin{remark} When acting on a funciton $f(\theta)$, the operators $Q$ and $\tilde{Q}$ are
$$Q f(\theta) = \pd_\theta^2f +\cot\theta\pd_\theta f$$
$$\tilde{Q}f(\theta) = \pd_\theta^2 f + 5\cot\theta\pd_\theta f$$
\end{remark}

The following lemma states that $Q$ and $\tilde{Q}$ are good commutators for $\Box_g$ and $\Box_{\tilde{g}}$ respectively. Note that there is no need as in the Kerr case to premultiply by $q^2$, since $q^2=r^2$ when $a=0$ and $r^2$ does not depend on $\theta$.
\begin{lemma}\label{ws:s_comm_lem}
\begin{align*}
[\pd_t,\Box_g] &= 0 \\
[\pd_t,\Box_{\tilde{g}}] &= 0 \\
[Q,\Box_g] &= 0 \\
[\tilde{Q},\Box_{\tilde{g}}] &= 0.
\end{align*}
\end{lemma}
\begin{proof}
There is no explicit $t$ dependence in either metric $g$ or $\tilde{g}$, so $\pd_t$ clearly commutes with both $\Box_g$ and $\Box_{\tilde{g}}$.

It is easy to check that the operator $\Box_g-\sla\triangle$ does not depend on $\theta$ and has no $\pd_\theta$ operators. Since $Q$ only depends on $\theta$ and has only $\pd_\theta$ operators, that means
$$0=[Q,\Box_g-\sla\triangle]=[Q,\Box_g]-[Q,r^{-2}Q] =[Q,\Box_g].$$
The final step required that $Q(r^2)=0$. In Kerr, since $Q(q^2)\ne 0$, there is a need to premultiply the operator $\Box_g$ by $q^2$.

The same argument shows that
$$0=[\tilde{Q},\Box_{\tilde{g}}].$$
\end{proof}

We now define the $s$-order commutators $\Gamma$ and $\tilde{\Gamma}$.
\begin{definition}
\begin{align*}
\Gamma^s u &= Q^l\pd_t^{s-2l}u \\
\tilde{\Gamma}^s u &= \tilde{Q}^l\pd_t^{s-2l}u
\end{align*}
where $0\le 2l \le s$.
\end{definition}

We also define the $s$-order wavefunctions derived from the $\phi$ and $\psi$ by applying symmetry operators.
\begin{definition}
\begin{align*}
\phi^s &= \Gamma^s\phi \\
\psi^s &= \tilde{\Gamma}^s\psi.
\end{align*}
\end{definition}

\subsection{The additional operator $\tg$ as a commutator}

We also note that we will need to use the operator $\tg$ defined in \S\ref{k:additional_commutator_sec}. The commutator estimates generalize to $\Box_{\tilde{g}}$ without any additional difficulty.

Thus, we have the following definition for the $s,k$-order wavefunctions.
\begin{definition}
\begin{align*}
\phi^{s,k} &= \tg^k\Gamma^s\phi \\
\psi^{s,k} &= \tg^k\tilde{\Gamma}^s\psi.
\end{align*}
\end{definition}

\subsection{The dynamic estimates for $(\phi^{s,k},\psi^{s,k})$}

By repeating the procedure in \S\ref{k:dynamic_estimates_s_k_sec}, we can prove the following proposition without any additional difficulty.

\begin{proposition}
Fix $\delm,\delp>0$. The following estimates hold for $p\in [\delm,2-\delp]$, $s\ge 0$, $k\ge 1$, and arbitrary $p'$.
$$E^{s,k}(t_2)\lesssim E^{s,k}(t_1)+\int_{t_1}^{t_2}B_{p'}^{s,k}(t)+B_{p'}^{s+2,k-1}(t)+N^{s,k}(t)dt,$$
$$E_p^{s,k}(t_2)+\int_{t_1}^{t_2}B_p^{s,k}(t)dt\lesssim E_p^{s,k}(t_1)+\int_{t_1}^{t_2}B_{p'}^{s+2,k-1}(t)+N_p^{s,k}(t)dt,$$
where
\begin{align*}
E^{s,k}(t) &= \sum_{\substack{s'\le s \\ k'\le k}} E[\psi^{s',k'}](t), \\
E_p^{s,k}(t) &= \sum_{\substack{s'\le s \\ k'\le k}} E_p[\psi^{s',k'}](t), \\
B_p^{s,k}(t) &= \sum_{\substack{s'\le s \\ k'\le k}} B_p[\psi^{s',k'}](t),
\end{align*}
$$N^{s,k}(t)= (E^{s,k}(t))^{1/2}\sum_{\substack{s'\le s \\ k'\le k}}\left(||\tg^{k'}\Gamma^{s'}\Box_g\phi||_{L^2(\Sigma_t)}+||\tg^{k'}\tilde\Gamma^{s'}\Box_g\psi||_{L^2(\tilde{\Sigma}_t)}\right),$$
$$N_p^{s,k}(t) = \sum_{\substack{s'\le s \\ k'\le k}}\int_{\Sigma_t}r^{p+1}|\tg^{k'}\Gamma^{s'}\Box_g\phi|^2 + \sum_{\substack{s'\le s \\ k'\le k}}\int_{\Sigma_t}r^{p+1}|\tg^{k'}\Gamma^{s'}\Box_g\psi|^2,$$
and the norms $E(t)$, $E_p(t)$, and $B_p(t)$ are as defined in Proposition \ref{ws:dynamic_estimates_0_prop}.
\end{proposition}
\begin{proof}
The proof is a direct application of Proposition \ref{ws:dynamic_estimates_0_prop} by making the substitutions
\begin{align*}
\phi &\mapsto\phi^{s'} \\
\psi &\mapsto\psi^{s'}
\end{align*}
for all values of $s'$ (and all commutators represented by $\Gamma^{s'}$ and $\tilde{\Gamma}^{s'}$) where $s'\le s$ and observing Lemma \ref{ws:s_comm_lem}.
\end{proof}

\section{The $L^\infty$ estimates}\label{ws:pointwise_sec}

In this section, we prove $L^\infty$ estimates for certain derivatives of $\phi$ and $\psi$ that will appear in the quantities $\tg^k\Gamma^s\Box_g\phi$ and $\tg^k\tilde{\Gamma}^s\Box_{\tilde{g}}\psi$ that become nonlinear when using the equations for $\phi$ and $\psi$.

We will use the following notation, which is similar to the notation introduced in Chapter \ref{kerr_chap}.
\begin{definition}
We define the expressions $\Omega^l\phi^{s-l}$ and $\Omega^l\psi^{s-l}$ to represent the following coordinate-invariant $r^l$-weighted rank $l$ tensors on $S^2(r)$ and $S^6(r)$.
$$\Omega^l\phi^{s-l}=r^l\sla\nabla^l\phi^{s-l}\text{ or }r^l\tsla\nabla^l\phi^{s-l}$$
$$\Omega^l\psi^{s-l}=r^l\sla\nabla^l\psi^{s-l}\text{ or }r^l\tsla\nabla^l\psi^{s-l}.$$
(It should always be assumed that $l\le s$.)
\end{definition}
In Chapter \ref{kerr_chap}, this notation was necessitated by the fact that the operator $Q$ is second order, and therefore does not satisfy the Leibniz property when acting on products. This fact is still true for both $Q$ and $\tilde{Q}$ in this problem, but there is also a new issue due to the fact that $\psi$ appears in the nonlinear terms in the equation for $\Box_g\phi$ and $\phi$ appears in the nonlinear terms in the equation for $\Box_{\tilde{g}}\psi$. So, for example, when applying $\tilde{Q}$ to the equation for $\Box_{\tilde{g}}\psi$, even the terms where $\tilde{Q}$ does not split between factors may not be expressable as products of $\phi^{s,k}$ and $\psi^{s,k}$, since the factor $\tilde{Q}\phi$ is not of the form $\phi^{s,k}$ for any $s$ and $k$. Instead, the factor $\tilde{Q}\phi$ must be treated as $r^2\tsla\nabla^2\phi$, which, since it is a coordinate-invariant tensor, is similar to $r^2\sla\nabla^2\phi$.

\subsection{Sobolev-type estimates for $\Sigma_t$ and $\tilde{\Sigma}_t$}

First, we prove the following lemma, which includes the Sobolev-type estimate for $\Sigma_t$ (Lemma \ref{k:pointwise_lem}) as well as a new Sobolev-type estimate for $\tilde{\Sigma}_t$. In particular, the $r$ weight that is gained depends on the volume form for the associated space, so the estimate on $\Sigma_t$ gains two factors of $r$ and the estimate on $\tilde{\Sigma}_t$ gains six factors of $r$. Also, the fact that the volume form for $\tilde{\Sigma}_t$ has additional factors of $\sin\theta$ means that more derivatives are required in the Sobolev estimate for $\tilde{\Sigma}_t$.

\begin{lemma}\label{ws:infty_base_lem}
If $u$ decays sufficiently fast as $r\rightarrow\infty$, then
$$||\Omega^lu||^2_{L^\infty(\Sigma_t\cap\{r>r_0\})}\lesssim \int_{\Sigma_t\cap\{r>r_0\}}r^{-2}\left[(\pd_r\Gamma^{\le l+3}u)^2+(\Gamma^{\le l+3}\phi)^2\right]$$
and
$$||\Omega^lu||^2_{L^\infty(\tilde\Sigma_t\cap\{r>r_0\})}\lesssim \int_{\tilde{\Sigma}_t\cap\{r>r_0\}}r^{-6}\left[(\pd_r\tilde{\Gamma}^{\le l+5} u)^2+(\tilde{\Gamma}^{\le l+5} u)^2\right].$$
If $l$ is even, the same results hold with only $\Gamma^{\le l+2}$ and $\tilde{\Gamma}^{\le l+4}$ respectively.
\end{lemma}
\begin{proof}
The first estimate is the same as the estimate in Lemma \ref{k:pointwise_lem}. The second estimate is proved in a similar way to the first estimate, with a few differences highlighted below.

To start, we note that in the case where $l$ is even,
\begin{multline*}
||\Omega^l u||_{L^\infty(S^6(r))}^2=||\tsla\nabla^l \bar{u}||_{L^\infty(S^6(1))}^2 \\
\lesssim \int_{S^6(1)}|\tsla\nabla^{l+4}\bar{u}|^2 = \int_{S^6(r)}|\Omega^{\le l+4}u|^2d\tilde{\omega} \lesssim \int_{S^6(r)}(\tilde\Gamma^{\le s+4}u)^2d\tilde{\omega}.
\end{multline*}
The main difference is in the second step, which relies on a classic Sobolev estimate. Since $S^6(1)$ has $6$ dimensions, there is a loss of $4$ derivatives instead of $2$ derivatives for $S^2(1)$.

As with the first estimate, there is an additional loss of one derivative if $l$ is odd. Therefore, for arbitrary $l$,
$$||\Omega^l u||_{L^\infty(S^6(r))}^2 \lesssim \int_{S^6(r)}(\tilde\Gamma^{\le s+5}u)^2d\tilde{\omega}.$$

As with the first estimate, we now set $f(r)=\int_{S^6(r)}(\tilde\Gamma^{\le s+5}u)^2d\tilde{\omega}$ and conclude that
\begin{multline*}
|f'(r)|\lesssim \int_{S^6(r)} \left[(\pd_r\tilde\Gamma^{\le l+5}u)^2+(\tilde\Gamma^{\le l+3}u)^2\right]d\tilde{\omega} \\
\lesssim \int_{S^6(r)} r^{-6} \left[(\pd_r\tilde\Gamma^{\le l+5}u)^2+(\tilde\Gamma^{\le l+3}u)^2\right]r^6d\tilde{\omega}.
\end{multline*}
Finally, assuming $\lim_{r\rightarrow\infty}f(r)=0$, we conclude
\begin{multline*}
|\Omega^lu(r_0)|^2\lesssim f(r_0)\lesssim \int_{r_0}^\infty |f'(r)|dr \lesssim \int_{r_0}^\infty\int_{S^6(r)}r^{-6}\left[(\pd_r\tilde\Gamma^{\le l+5}u)^2+(\tilde\Gamma^{\le l+3}u)^2\right]r^6d\tilde{\omega}dr \\
\lesssim \int_{\tilde\Sigma_t} r^{-6}\left[(\pd_r\tilde\Gamma^{\le l+5}u)^2+(\tilde\Gamma^{\le l+3}u)^2\right]
\end{multline*}
The main difference is that the factor of $r^6$ is needed for the volume form of $\tilde{\Sigma}_t$ instead of $r^2$.
\end{proof}

\subsection{Estimating derivatives using the Sobolev-type estimate}

Now, we repeatedly apply Lemma \ref{ws:infty_base_lem} to estimate various derivatives with $r$ weights. We will assume that $\phi$ and $\psi$ decay sufficiently fast as $r\rightarrow\infty$.

The following lemma estimates $\phi$ and $\psi$, as well as the higher order analogues $\Omega^l\phi^{s-l,k}$ and $\Omega^l\psi^{s-l,k}$.
\begin{lemma}\label{ws:low_order_infty_lem}
For $r\ge r_H$,
$$|r^p\Omega^l\phi^{s-l,k}|^2+|r^{p+2}\Omega^l\psi^{s-l,k}|^2\lesssim E_{2p}^{s+5,k}(t)$$
and for $r\ge r_0>r_H$,
$$|\Omega^l\phi^{s-l,k}|^2+|r^2\Omega^l\psi^{s-l,k}|^2\lesssim E^{s+5,k}(t).$$
\end{lemma}
\begin{proof}
First, we apply Lemma \ref{ws:infty_base_lem} with $u=r^p\phi^{s-l,k}$.
\begin{align*}
|r^p\Omega^l\phi^{s-l,k}|^2 &=|\Omega^l(r^p\phi^{s-l,k})|^2 \\
&\lesssim \int_{\Sigma_t}r^{-2}\left[(\pd_r\Gamma^{\le l+3}(r^p\phi^{s-l,k}))^2+(\Gamma^{\le l+3}(r^p\phi^{s-l,k}))^2\right] \\
&\lesssim \int_{\Sigma_t}r^{2p-2}\left[(\pd_r\Gamma^{\le l+3}\phi^{s-l,k})^2+(\Gamma^{\le l+3}\phi^{s-l,k})^2\right] \\
&\lesssim E_{2p}^{s+3}(t).
\end{align*}
Then, we apply Lemma \ref{ws:infty_base_lem} with $u=r^{p+2}\psi^{s-l,k}$.
\begin{align*}
|r^{p+2}\Omega^l\psi^{s-l,k}|^2 &=|\Omega^l(r^{p+2}\psi^{s-l,k})|^2 \\
&\lesssim \int_{\tilde\Sigma_t}r^{-6}\left[(\pd_r\Gamma^{\le l+5}(r^{p+2}\psi^{s-l,k}))^2+(\Gamma^{\le l+5}(r^{p+2}\psi^{s-l,k}))^2\right] \\
&\lesssim \int_{\tilde\Sigma_t}r^{2p-2}\left[(\pd_r\Gamma^{\le l+5}\psi^{s-l,k})^2+(\Gamma^{\le l+5}\psi^{s-l,k})^2\right] \\
&\lesssim E_{2p}^{s+5}(t).
\end{align*}
Together, these estimates prove the first estimate of the lemma. The second estimate follows from the same exact argument in the special case $p=0$, and the observation that as long as $r\ge r_0>r_H$, then $E^{s,k}(t)$ can be used in place of $E_0^{s,k}(t)$.
\end{proof}

The following lemma estimates $\pd_t\phi$ and $\pd_t\psi$ as well as the higher order analogues $\pd_t\Omega^l\phi^{s-l,k}$ and $\pd_t\Omega^l\psi^{s-l,k}$.
\begin{lemma}\label{ws:pd_t_infty_lem}
For $r\ge r_H$,
$$|r^p\pd_t\Omega^l\phi^{s-l,k}|^2+|r^{p+2}\pd_t\Omega^l\psi^{s-l,k}|^2 \lesssim E_{2p}^{s+6,k}(t)$$
and for $r\ge r_0>r_H$,
$$|r\pd_t\Omega^l\phi^{s-l,k}|^2+|r^3\pd_t\Omega^l\psi^{s-l,k}|^2 \lesssim E^{s+6,k}(t).$$
\end{lemma}
\begin{proof}
The first estimate reduces to Lemma \ref{ws:low_order_infty_lem} by observing that $\pd_t\Omega^l\phi^{s-l,k}=\Omega^l\pd_t\phi^{s-l,k}=\Omega^l\phi^{s+1-l,k}$ and likewise $\pd_t\Omega^l\psi^{s-l,k}=\Omega^l\psi^{s+1-l,k}$. We now prove the second estimate.

First, we apply Lemma \ref{ws:infty_base_lem} with $u=r\pd_t\phi^{s-l,k}$.
\begin{align*}
|r\pd_t\Omega^l\phi^{s-l,k}|^2 &= |\Omega^l(r\pd_t\phi^{s-l,k})|^2 \\
&\lesssim \int_{\Sigma_t\cap\{r>r_0\}}r^{-2}\left[(\pd_r\Gamma^{\le l+3}(r\pd_t\phi^{s-l,k}))^2+(\Gamma^{\le l+3}(r\pd_t\phi^{s-l,k}))^2\right] \\
&\lesssim \int_{\Sigma_t\cap\{r>r_0\}}\left[(\pd_r\Gamma^{\le l+3}\pd_t\phi^{s-l,k})^2+(\Gamma^{\le l+3}\pd_t\phi^{s-l,k})^2\right] \\
&\lesssim \int_{\Sigma_t\cap\{r>r_0\}}\left[(\pd_r\Gamma^{\le l+3}\phi^{s+1-l,k})^2+(\pd_t\Gamma^{\le l+3}\phi^{s-l,k})^2\right] \\
&\lesssim E^{s+4,k}(t).
\end{align*}
Next, by applying Lemma \ref{ws:infty_base_lem} with $u=r^3\pd_t\psi^{s-l,k}$ and repeating the same procedure, we arrive at the following estimate.
$$|r^3\pd_t\Omega^l\psi^{s-l,k}|^2\lesssim E^{s+6,k}(t).$$
Together, these estimates prove the second estimate of the lemma.
\end{proof}

The following lemma estimates $\sla\nabla\phi$ and $\tsla\nabla\psi$ as well as the higher order analogues $\sla\nabla\Omega^l\phi^{s-l,k}$ and $\tsla\nabla\Omega^l\psi^{s-l,k}$.
\begin{lemma}\label{ws:pd_theta_infty_lem}
For $r\ge r_H$,
$$|r^{p+1}\sla\nabla\Omega^l\phi^{s-l,k}|^2+|r^{p+3}\tsla\nabla\Omega^l\psi^{s-l,k}|^2\lesssim E_{2p}^{s+6,k}(t)$$
and for $r\ge r_0>r_H$,
$$|r\sla\nabla\Omega^l\phi^{s-l,k}|^2+|r^3\tsla\nabla\Omega^l\psi^{s-l,k}|^2\lesssim E^{s+6,k}(t).$$
\end{lemma}
\begin{proof}
This lemma reduces to Lemma \ref{ws:low_order_infty_lem} by observing that
$$r^{p+1}\sla\nabla\Omega^l\phi^{s-l,k}=r^p\Omega\Omega^l\phi^{s-l,k}=r^p\Omega^{l+1}\phi^{s-l,k}\subset r^p\Omega^l\phi^{s+1-l,k},$$
and likewise
$$r^{p+3}\tsla\nabla\Omega^l\psi^{s-l,k}\subset r^{p+2}\Omega^l\psi^{s+1-l,k}.$$
\end{proof}

The following lemma estimates $L\phi$ and $L\psi$ as well as the higher order analogues $L\Omega^l\phi^{s-l,k}$ and $L\Omega^l\psi^{s-l,k}$.
\begin{lemma}\label{ws:L_infty_lem}
Letting $L=\alpha\pd_r+\pd_t$, where $\alpha=\frac{\Delta}{r^2+a^2}$, we have that for $r\ge r_H$,
$$|r^{p+1}L\Omega^l\phi^{s-l,k}|^2+|r^{p+3}L\Omega^l\psi^{s-l,k}|^2\lesssim E_{2p}^{s+7,k}(t)+\int_{\Sigma_t}r^{2p}(\Box_g\phi^{s+3,k})^2+\int_{\tilde\Sigma_t}r^{2p}(\Box_g\psi^{s+5,k})^2$$
and for $r\ge r_0>r_H$,
$$|rL\Omega^l\phi^{s-l,k}|^2+|r^3L\Omega^l\psi^{s-l,k}|^2\lesssim E^{s+7,k}(t)+\int_{\Sigma_t\cap\{r>r_0\}}(\Box_g\phi^{s+3,k})^2+\int_{\tilde\Sigma_t\cap\{r>r_0\}}(\Box_g\psi^{s+5,k})^2.$$
\end{lemma}
\begin{proof}
We recall from the proof of Lemma \ref{s:L_pointwise_lem} that
$$(\pd_rLu)^2\lesssim (\Box_gu)^2+(L\pd_tu)^2+r^{-2}(\pd_t^2u)^2+r^{-2}(\pd_r\pd_tu)^2+r^{-2}(\pd_ru)^2+r^{-2}(Qu)^2,$$
By the same technique,
$$(\pd_rLu)^2\lesssim (\Box_{\tilde g}u)^2+(L\pd_tu)^2+r^{-2}(\pd_t^2u)^2+r^{-2}(\pd_r\pd_tu)^2+r^{-2}(\pd_ru)^2+r^{-2}(\tilde{Q}u)^2.$$

We now apply Lemma \ref{ws:infty_base_lem} with $u=r^{p+1}L\phi^{s-l,k}$.
\begin{align*}
|r^{p+1}L\Omega^l\phi^{s-l,k}|^2 &= |\Omega^l(r^{p+1}L\phi^{s-l,k})|^2 \\
&\lesssim \int_{\Sigma_t}r^{-2}\left[(\pd_r\Gamma^{\le l+3}(r^{p+1}L\phi^{s-l,k}))^2+(\Gamma^{\le l+3}(r^{p+1}L\phi^{s-l,k}))^2\right] \\
&\lesssim \int_{\Sigma_t}r^{2p}\left[(\pd_r\Gamma^{\le l+3}L\phi^{s-l,k})^2+(\Gamma^{\le l+3}L\phi^{s-l,k})^2\right] \\
&\lesssim \int_{\Sigma_t}r^{2p}\left[(\pd_rL\phi^{s+3,k})^2+(L\phi^{s+3,k})^2\right] \\
&\lesssim E_{2p}^{s+3,k}(t) +\int_{\Sigma_t}r^{2p}(\pd_rL\phi^{s+3,k})^2.
\end{align*}
Now,
\begin{align*}
\int_{\Sigma_t}r^{2p}&(\pd_rL\phi^{s+3,k})^2 \\
&\lesssim \int_{\Sigma_t}r^{2p}\left[(\Box_g\phi^{s+3,k})^2+(L\pd_t\phi^{s+3,k})^2+r^{-2}(\pd_t^2\phi^{s+3,k})^2+r^{-2}(\pd_r\pd_t\phi^{s+3,k})^2\right.\\
&\hspace{3in}\left.+r^{-2}(\pd_r\phi^{s+3,k})^2+r^{-2}(Q\phi^{s+3,k})^2\right] \\
&\lesssim \int_{\Sigma_t}r^{2p}\left[(\Box_g\phi^{s+3,k})^2+(L\phi^{s+4,k})^2+r^{-2}(\phi^{s+5,k})^2+r^{-2}(\pd_r\phi^{s+4,k})^2\right].
\end{align*}
It follows that
$$|r^{p+1}L\Omega^l\phi^{s-l,k}|^2\lesssim E_{2p}^{s+5,k}(t)+\int_{\Sigma_t}r^{2p}(\Box_g\phi^{s+3,k})^2.$$
By a similar argument,
$$|r^{p+3}L\Omega^l\psi^{s-l,k}|^2\lesssim E_{2p}^{s+7,k}(t)+\int_{\tilde{\Sigma}_t}r^{2p}(\Box_{\tilde{g}}\psi^{s+5,k})^2.$$
These two estimates together establish the first estimate of the lemma. The second estimate follows from the same exact argument in the special case $p=0$, and the observation that as long as $r\ge r_0>r_H$, then $E^{s,k}(t)$ can be used in place of $E_0^{s,k}(t)$.
\end{proof}

The following lemma estimates $\tg\phi$ and $\tg\psi$ as well as the higher order analogues $\tg\Omega^l\phi^{s-l,k}$ and $\tg\Omega^l\psi^{s-l,k}$.
\begin{lemma}\label{ws:gh_infty_lem}
Keeping in mind that $\tg$ is supported in a neighborhood of the event horizon, for arbitrary $p'$, we have for $r\ge r_H$,
$$|\tg \Omega^l\phi^{s-l,k}|^2+|\tg \Omega^l\psi^{s-l,k}|^2 \lesssim E_{p'}^{s+5,k+1}(t)$$
and for $r\ge r_0>r_H$,
$$|\tg \Omega^l\phi^{s-l,k}|^2+|\tg \Omega^l\psi^{s-l,k}|^2 \lesssim E^{s+5,k+1}(t).$$
\end{lemma}
\begin{proof}
We apply Lemma \ref{ws:infty_base_lem} with $u=\tg\phi^{s-l,k}$, and freely introduce a factor of $r^{p'}$ since $\tg$ is supported on a compact interval in $r$.
\begin{align*}
|\tg\Omega^l\phi^{s-l,k}|^2 &= |\Omega^l\phi^{s-l,k+1}|^2 \\
&\lesssim \int_{\Sigma_t}r^{-2}\left[(\pd_r\Gamma^{\le l+3}\phi^{s-l,k+1})^2+(\Gamma^{\le l+3}\phi^{s-l,k+1})^2\right] \\
&\lesssim \int_{\Sigma_t}r^{p'-2}\left[(\pd_r\Gamma^{\le l+3}\phi^{s-l,k+1})^2+(\Gamma^{\le l+3}\phi^{s-l,k+1})^2\right] \\
&\lesssim E_{p'}^{s+3,k+1}(t).
\end{align*}
A similar argument shows that
$$|\tg \Omega^l\psi^{s-l,k}|^2\lesssim E_{p'}^{s+5,k+1}(t).$$
Together, these estimates prove the first estimate of the lemma. The second estimate follows from the same argument, and the observation that as long as $r\ge r_0>r_H$, then $E^{s,k}(t)$ can be used in place of $E_{p'}^{s,k}(t)$.
\end{proof}

\subsection{Summarizing the $L^\infty$ estimates}
To conclude this section, we summarize the previous lemmas in a single proposition.
\begin{definition} We define two families of operators.
$$\bar{D}=\{L,r^{-1}\pd_\theta\},$$
$$D=\{L,\pd_r,r^{-1}\pd_\theta\}.$$
\end{definition}

\begin{proposition}\label{ws:infinity_prop}
For $r\ge r_H$,
\begin{multline*}
|r^{p+1}\bar{D}\Omega^l\phi^{s-l,k}|^2+|r^pD\Omega^l\phi^{s-l,k}|^2+|r^p\Omega^l\phi^{s-l,k}|^2 \\
+|r^{p+3}\bar{D}\Omega^l\psi^{s-l,k}|^2+|r^{p+2}D\Omega^l\psi^{s-l,k}|^2+|r^{p+2}\Omega^l\psi^{s-l,k}|^2 \\
\lesssim E_{2p}^{s+5,k+1}(t)+E_{2p}^{s+7,k}(t)+\int_{\Sigma_t}r^{2p}(\Box_g\phi^{s+5,k})^2+\int_{\tilde\Sigma_t}r^{2p}(\Box_{\tilde{g}}\psi^{s+5,k})^2
\end{multline*}
and for $r\ge r_0>r_H$,
\begin{multline*}
|rD\Omega^l\phi^{s-l,k}|^2+|r^3D\Omega^l\psi^{s-l,k}|^2 \\
\lesssim E^{s+5,k+1}(t)+E^{s+7,k}(t)+\int_{\Sigma_t\cap\{r>r_0\}}(\Box_g\phi^{s+5,k})^2+\int_{\tilde\Sigma_t\cap\{r>r_0\}}(\Box_{\tilde{g}}\psi^{s+5,k})^2.
\end{multline*}
\end{proposition}
\begin{proof}
With the exception of the operator $\pd_r$, all of the cases have been proved in Lemmas \ref{ws:low_order_infty_lem}-\ref{ws:gh_infty_lem}. Finally, observe that
$$|r^p\pd_r\Omega^l\phi^{s-l,k}|^2\lesssim |r^pL\Omega^l\phi^{s-l,k}|^2+|r^p\pd_t\Omega^l\phi^{s-l,k}|^2+|\tg\Omega^l\phi^{s-l,k}|^2$$
and
$$|r^p\pd_r\Omega^l\phi^{s-l,k}|^2\lesssim |r^pL\Omega^l\phi^{s-l,k}|^2+|r^p\pd_t\Omega^l\phi^{s-l,k}|^2+|\tg\Omega^l\phi^{s-l,k}|^2.$$
Thus, even the case of the operator $\pd_r$ can be reduced to Lemmas \ref{ws:low_order_infty_lem}-\ref{ws:gh_infty_lem}.
\end{proof}

\section{The structure of the nonlinear terms}\label{ws:structure_sec}

In this section, we carefully examine the nonlinear terms $\mathcal{N}_\phi$ and $\mathcal{N}_\psi$ as well as their higher order analogues $\tg^k\Gamma^s\mathcal{N}_\phi$ and $\tg^k\tilde\Gamma^s\mathcal{N}_\psi$ and determine a procedure for estimating them in the proof of the main theorem. In particular, we note that these terms satisfy a version of the null condition and (as a consequence of a special linearization) all coefficients of these terms are bounded on the axis, thus allowing us to avoid the formalism in Appendix \ref{regularity_sec} that will be necessary for the final problem in Chapter \ref{wm_kerr_chap}.

We start in \S\ref{ws:special_choice_sec} by deriving the equations for $(\phi,\psi)$. In particular, in the proof of Proposition \ref{ws:special_choice_prop}, we see the benefit for choosing the special linearization. Then in \S\ref{ws:example_terms_sec} we look at a few example terms to illustrate the procedure to be used by the proof of the main theorem (Theorem \ref{ws:main_thm}). In particular, we will see how to estimate terms that contain products of both $\phi$ and $\psi$ as well as the role of the null condition. Then in \S\ref{ws:nl_structure_sec}, we define and prove the precise structures of $\mathcal{N}_\phi$ and $\mathcal{N}_\psi$. These structures will then be used in \S\ref{ws:nl_s_k_structure_sec} to define and prove the precise structures of $\tg^k\Gamma^s\mathcal{N}_\phi$ and $\tg^k\tilde\Gamma^s\mathcal{N}_\psi$. Finally, in \S\ref{ws:nl_strategy_revisited_sec} we prove a proposition that uses these structures to estimate nonlinear quantities that will show up in the proof of the main theorem.

\subsection{Deriving the equations for $(\phi,\psi)$}\label{ws:special_choice_sec}

We start by deriving the equations for $\phi$ and $\psi$. Recall that
$$\frac{\pd_\theta A}{A} = 2\cot\theta.$$
Since $\cot\theta$ is singular on the axis, we generally want to avoid having terms with $\frac{\pd_\alpha A}{A}$ as they could potentially be very bad when we commute the equations with angular commutators.\footnote{In Chapter \ref{wm_kerr_chap}, we will see a more general procedure that can handle some of these terms.} Our particular choice for the linearization
\begin{align*}
X &= A+A\phi \\
Y &= X^2\psi
\end{align*}
is motivated by the desire to eliminate terms with $\frac{\pd_\alpha A}{A}$ from the equations.

\begin{lemma}\label{ws:k_sub_lem}
If one makes the substitutions
\begin{align*}
X &= A + A\phi \\
Y &= A^k\zeta
\end{align*}
and requires that $\phi$ and $\zeta$ are axisymmetric functions, then the wave map system (\ref{ws:X_eqn}-\ref{ws:Y_eqn}) reduces to the following system of equations for $\phi$ and $\zeta$.
\begin{align*}
\Box_g\phi &= \frac{\pd^\alpha\phi\pd_\alpha\phi}{1+\phi} - \frac{A^{2k-2}}{1+\phi}\left(\pd^\alpha\zeta\pd_\alpha\zeta +2k\frac{\pd^\alpha A}{A}\zeta\pd_\alpha\zeta +k^2\frac{\pd^\alpha A\pd_\alpha A}{A^2}\zeta^2\right) \\
\Box_g\zeta &+2(k-1)\frac{\pd^\alpha A}{A}\pd_\alpha \zeta+k(k-2)\frac{\pd^\alpha A\pd_\alpha A}{A^2}\zeta
= 2k\frac{\pd^\alpha A}{A} \frac{\pd_\alpha\phi}{1+\phi}\zeta+2\frac{\pd^\alpha \phi\pd_\alpha\zeta}{1+\phi}
\end{align*}
In particular, if $k\ge 2$, then all of the terms in the first equation have nonsingular factors.
\end{lemma}
\begin{proof}
We start with the substitution
$$X = A+A\phi = A(1+\phi) .$$
In this case,
$$\frac{\pd_\alpha X}{X} =\frac{\pd_\alpha (A(1+\phi))}{A(1+\phi)} = \frac{\pd_\alpha A}{A}+\frac{\pd_\alpha\phi}{1+\phi}$$
and
$$X^{-1}\Box_g X = \frac{\Box_g (A(1+\phi))}{A(1+\phi)} =\frac{\Box_g A}{A}+\frac{\Box_g\phi}{1+\phi}+2\frac{\pd^\alpha A\pd_\alpha\phi}{A(1+\phi)}.$$
Therefore, the equation
$$X^{-1}\Box_g X = \frac{\pd^\alpha X\pd_\alpha X}{X^2}-\frac{\pd^\alpha Y\pd_\alpha Y}{X^2}$$
becomes
$$\frac{\Box_g A}{A}+2\frac{\pd^\alpha A\pd_\alpha\phi}{A(1+\phi)}+\frac{\Box_g\phi}{1+\phi}
= \frac{\pd^\alpha A\pd_\alpha A}{A^2}+2\frac{\pd^\alpha A\pd_\alpha\phi}{A(1+\phi)}+\frac{\pd^\alpha \phi\pd_\alpha \phi}{(1+\phi)^2}-\frac{\pd^\alpha Y\pd_\alpha Y}{A^2(1+\phi)^2}.$$
The first terms on both sides cancel because $A$ itself is a solution to the equation, and the second terms on both sides cancel because they are identitcal. The result, after multiplying by $1+\phi$ is the following equation.
$$\Box_g\phi = \frac{\pd^\alpha \phi\pd_\alpha \phi}{1+\phi}-\frac{\pd^\alpha Y\pd_\alpha Y}{A^2(1+\phi)}.$$
Now, by substituting
$$Y = A^k\zeta,$$
this equation becomes
$$\Box_g\phi = \frac{\pd^\alpha \phi\pd_\alpha \phi}{1+\phi} - \frac{A^{2k-2}}{1+\phi}\left(\pd^\alpha\zeta\pd_\alpha\zeta +2k\frac{\pd^\alpha A}{A}\zeta\pd_\alpha\zeta +k^2\frac{\pd^\alpha A\pd_\alpha A}{A^2}\zeta^2\right).$$
This verifies the first equation.

We now turn to the second equation.
$$\Box_g Y = 2\frac{\pd^\alpha X\pd_\alpha Y}{X}.$$
Given the substitution for $Y$,
\begin{align*}
A^{-k}\Box_g Y &= A^{-k}\Box_g(A^k\zeta) \\
&= \Box_g\zeta +2k\frac{\pd^\alpha A}{A}\pd_\alpha \zeta +\zeta A^{-k}\Box_g(A^k) \\
&= \Box_g\zeta +2k\frac{\pd^\alpha A}{A}\pd_\alpha \zeta +\zeta A^{-k}(kA^{k-1}\Box_gA+k(k-1)A^{k-2}\pd^\alpha A\pd_\alpha A) \\
&= \Box_g\zeta +2k\frac{\pd^\alpha A}{A}\pd_\alpha \zeta +k^2\frac{\pd^\alpha A\pd_\alpha A}{A^2}\zeta,
\end{align*}
where we have used the fact that $\Box_g A = \frac{\pd^\alpha A\pd_\alpha A}{A}$ in the last step.

Given the substitutions for $X$ and $Y$,
\begin{align*}
2\frac{\pd^\alpha X\pd_\alpha Y}{X A^k} &= 2\frac{\pd^\alpha (A(1+\phi))\pd_\alpha (A^k\zeta)}{A(1+\phi)A^k} \\
&= 2k\frac{\pd^\alpha A\pd_\alpha A}{A^2}\zeta +2\frac{\pd^\alpha A}{A}\pd_\alpha \zeta +2k\frac{\pd^\alpha A\pd_\alpha \phi}{A(1+\phi)}\zeta +2\frac{\pd^\alpha\phi\pd_\alpha \zeta}{1+\phi}.
\end{align*}

Equating both expressions, we obtain
\begin{multline*}
\Box_g\zeta +2k\frac{\pd^\alpha A}{A}\pd_\alpha \zeta +k^2\frac{\pd^\alpha A\pd_\alpha A}{A^2}\zeta \\
= 2k\frac{\pd^\alpha A\pd_\alpha A}{A^2}\zeta +2\frac{\pd^\alpha A}{A}\pd_\alpha \zeta +2k\frac{\pd^\alpha A\pd_\alpha \phi}{A(1+\phi)}\zeta +2\frac{\pd^\alpha\phi\pd_\alpha \zeta}{1+\phi}.
\end{multline*}
Thus,
$$\Box_g\zeta +2(k-1)\frac{\pd^\alpha A}{A}\pd_\alpha \zeta+k(k-2)\frac{\pd^\alpha A\pd_\alpha A}{A^2}\zeta
= 2k\frac{\pd^\alpha A\pd_\alpha \phi}{A(1+\phi)}\zeta+2\frac{\pd^\alpha \phi\pd_\alpha\zeta}{1+\phi}.$$
This verifies the second equation.
\end{proof}

\begin{lemma}\label{ws:2_sub_lem}
If one makes the same substitutions as in Lemma \ref{ws:k_sub_lem} and chooses $k=2$,
\begin{align*}
X &= A + A\phi \\
Y &= A^2\zeta
\end{align*}
then the wave map system (\ref{ws:X_eqn}-\ref{ws:Y_eqn}) reduces to the following system of equations for $\phi$ and $\zeta$.
\begin{align*}
\Box_g\phi &= \frac{1}{1+\phi}\left[\pd^\alpha\phi\pd_\alpha\phi - A^2\pd^\alpha\zeta\pd_\alpha\zeta -4 A\pd^\alpha A\zeta\pd_\alpha\zeta-4\pd^\alpha A\pd_\alpha A \zeta^2\right] \\
\Box_g\zeta &+2\frac{\pd^\alpha A}{A}\pd_\alpha \zeta
= 4\frac{\pd^\alpha A}{A} \frac{\pd_\alpha\phi}{1+\phi} \zeta + 2\frac{\pd^\alpha \phi\pd_\alpha\zeta}{1+\phi}
\end{align*}
In particular, the highly singular potential term in the second equation vanished by choosing $k=2$.
\end{lemma}
\begin{proof}
This follows directly from Lemma \ref{ws:k_sub_lem} by choosing $k=2$.
\end{proof}

\begin{proposition}\label{ws:special_choice_prop}
If one makes the substitutions
\begin{align*}
X &= A+A\phi \\
Y &= X^2\psi,
\end{align*}
which agree with the substitutions of Lemma \ref{ws:2_sub_lem} up to terms that are nonlinear in $\phi$ and $\psi$, then the wave map system (\ref{ws:X_eqn}-\ref{ws:Y_eqn}) reduces to the following system of equations for $\phi$ and $\psi$.
\begin{align*}
\Box_g\phi &= \mathcal{N}_\phi \\
\Box_{\tilde{g}}\psi &= \mathcal{N}_\psi,
\end{align*}
where
\begin{align*}
(1+\phi)\mathcal{N}_\phi &= \pd^\alpha\phi\pd_\alpha\phi -\mathcal{N} \\
(1+\phi)^2\mathcal{N}_\psi &= 2\psi\mathcal{N}-2(1+\phi)\pd^\alpha\phi\pd_\alpha\psi \\
\mathcal{N} &= A^2(1+\phi)^4 \pd^\alpha\psi\pd_\alpha\psi +4A^2(1+\phi)^3\psi\pd^\alpha\phi\pd_\alpha\psi +4A^2(1+\phi)^2\psi^2\pd^\alpha\phi\pd_\alpha\phi \\
&\hspace{12pt}+ 4A\pd^\alpha A (1+\phi)^4 \psi\pd_\alpha\psi +8 A\pd^\alpha A (1+\phi)^3\psi^2\pd_\alpha\phi +4\pd^\alpha A\pd_\alpha A (1+\phi)^4\psi^2
\end{align*}
In particular, no terms in $\mathcal{N}_\phi$ or $\mathcal{N}_\psi$ have singular factors.
\end{proposition}

\begin{proof}
\textbf{We begin the proof by motivating the choice $Y = X^2\psi$ as an alternative to the choice $Y=A^2\zeta$.} Recall the second equation from Lemma \ref{ws:2_sub_lem}
$$\Box_g\zeta +2\frac{\pd^\alpha A}{A}\pd_\alpha \zeta
= 4\frac{\pd^\alpha A}{A} \frac{\pd_\alpha\phi}{1+\phi} \zeta + 2\frac{\pd^\alpha \phi\pd_\alpha\zeta}{1+\phi}$$
There are still two terms in this equation with $\frac{\pd_\alpha A}{A}$ as a factor. As observed in the beginning of this chapter, the entire left side can be reinterpreted as $\Box_{\tilde{g}}\zeta$, with $\tilde{g}$ the metric for a higher dimensional spacetime. Essentially all of the work done in this chapter so far has prepared us to handle the term on the left side. But the term on the right side still causes concern.

Due to the similarity between the two terms, it is possible to combine them, yielding the following equation.
$$\Box_g\zeta +2\frac{\pd^\alpha A}{A}\left[\pd_\alpha \zeta -2\frac{\pd_\alpha\phi}{1+\phi}\zeta\right] =  2\frac{\pd^\alpha \phi\pd_\alpha\zeta}{1+\phi}$$
\textbf{Now, we make the key observation that if we choose}
$$\zeta = (1+\phi)^2\psi,$$
\textbf{then}
$$\pd_\alpha \zeta-2\frac{\pd_\alpha\phi}{1+\phi}\zeta = (1+\phi)^2\pd_\alpha\psi.$$
\textbf{As a consequence,}
$$Y = A^2\zeta = A^2(1+\phi)^2\psi = X^2\psi.$$

Continuing the substitution $\zeta = (1+\phi)^2\psi$ for the rest of the equation, we have
$$\Box_g((1+\phi)^2\psi) +2(1+\phi)^2\frac{\pd^\alpha A}{A}\pd_\alpha\psi =  2\frac{\pd^\alpha \phi\pd_\alpha((1+\phi)^2\psi)}{1+\phi}.$$
Expanding and dividing by $(1+\phi)^2$,
$$\Box_g\psi +4\frac{\pd^\alpha\phi\pd_\alpha\psi}{1+\phi}+2\frac{\Box_g\phi}{1+\phi}\psi+2\frac{\pd^\alpha\phi\pd_\alpha\phi}{(1+\phi)^2}\psi +2\frac{\pd^\alpha A}{A}\pd_\alpha\psi = 4\frac{\pd^\alpha\phi\pd_\alpha\phi}{(1+\phi)^2}\psi +2\frac{\pd^\alpha\phi\pd_\alpha\psi}{1+\phi}.$$
By rearranging terms, the second equation of Lemma \ref{ws:2_sub_lem} becomes
$$\Box_g\psi+2\frac{\pd^\alpha A}{A}\pd_\alpha\psi = 2\psi\left[\frac{\pd^\alpha\phi\pd_\alpha\phi}{(1+\phi)^2}-\frac{\Box_g\phi}{1+\phi}\right] -2\frac{\pd^\alpha\phi\pd_\alpha\psi}{1+\phi}.$$
If we define
$$\mathcal{N} = \pd^\alpha\phi\pd_\alpha\phi -(1+\phi)\Box_g\phi,$$
then the second equation from Lemma \ref{ws:2_sub_lem} reduces to
$$(1+\phi)^2\Box_{\tilde{g}}\psi = 2\psi\mathcal{N}-2(1+\phi)\pd^\alpha\phi\pd_\alpha\psi.$$

Recall also the first equation from Lemma \ref{ws:2_sub_lem}.
$$\Box_g\phi = \frac{1}{1+\phi}\left[\pd^\alpha\phi\pd_\alpha\phi - A^2\pd^\alpha\zeta\pd_\alpha\zeta -4 A\pd^\alpha A\zeta\pd_\alpha\zeta-4\pd^\alpha A\pd_\alpha A \zeta^2\right]$$
From this equation, we determine that
\begin{align*}
\mathcal{N} &= \pd^\alpha\phi\pd_\alpha\phi-(1+\phi)\Box_g\phi \\
&= A^2\pd^\alpha\zeta\pd_\alpha\zeta +4 A\pd^\alpha A\zeta\pd_\alpha\zeta + 4\pd^\alpha A\pd_\alpha A \zeta^2.
\end{align*}
Finally, we substitute $\zeta = (1+\phi)^2\psi$ in the expression for $\mathcal{N}$ and expand to conclude
\begin{align*}
\mathcal{N} &= A^2\pd^\alpha\zeta\pd_\alpha\zeta +4 A\pd^\alpha A\zeta\pd_\alpha\zeta + 4\pd^\alpha A\pd_\alpha A \zeta^2
 \\
&= A^2\pd^\alpha ((1+\phi)^2\psi)\pd_\alpha ((1+\phi)^2\psi) +4A\pd^\alpha A (1+\phi)^2\psi\pd_\alpha ((1+\phi)^2\psi) \\
&\hspace{2in}+4\pd^\alpha A\pd_\alpha A (1+\phi)^4\psi^2 \\
&= A^2(1+\phi)^4 \pd^\alpha\psi\pd_\alpha\psi +4A^2(1+\phi)^3\psi\pd^\alpha\phi\pd_\alpha\psi +4A^2(1+\phi)^2\psi^2\pd^\alpha\phi\pd_\alpha\phi \\
&\hspace{12pt}+ 4A\pd^\alpha A (1+\phi)^4 \psi\pd_\alpha\psi +8 A\pd^\alpha A (1+\phi)^3\psi^2\pd_\alpha\phi +4\pd^\alpha A\pd_\alpha A (1+\phi)^4\psi^2.
\end{align*}
\end{proof}

\subsection{The general strategy illustrated by two example terms}\label{ws:example_terms_sec}

There are many terms in $\mathcal{N}_\phi$ and $\mathcal{N}_\psi$. Let us illustrate how to handle them by looking at two particular examples.

First, we examine the term $r^4\sin^4\theta L\psi\lbar\psi$, which arises from the term $A^2\pd^\alpha\psi \pd_\alpha\psi$, which shows up in $\mathcal{N}_\phi$.
\begin{lemma}\label{ws:example_Nphi_term_lem} (Example $\mathcal{N}_\phi$ term)
Assume the following estimates.
$$||r^{(p-1)/2+3}L\Omega^l\psi^{s-l}||_{L^\infty(t)}^2\lesssim B_p^{s+7}(t),$$
$$||r^3\lbar\Omega^l\psi^{s-l}||_{L^\infty(t)}^2\lesssim E^{s+7}(t).$$
Then
$$\int_{\Sigma_t\cap\{r>R\}}r^{p+1}\left(\Gamma^s\left(r^4\sin^4\theta L\psi\lbar\psi\right)\right)^2\lesssim B_p^s(t)E^{s/2+7}(t)+E^s(t)B_p^{s/2+7}(t).$$
\end{lemma}
\begin{proof}
First, we compute
\begin{align*}
\int_{\Sigma\cap\{r>R\}}r^{p+1}r^8\sin^8\theta (L\psi^s)^2(\lbar\psi^{s/2})^2
&\lesssim \int_{\Sigma_t\cap\{r>R\}}r^{p-1}(L\psi^s)^2r^4\sin^4\theta (r^3\lbar\psi^{s/2})^2 \\
&\lesssim \int_{\Sigma_t\cap\{r>R\}}r^{p-1}(L\psi^s)^2r^4\sin^4\theta ||r^3\lbar\psi^{s/2}||_{L^\infty(t)}^2 \\
&\lesssim \int_{\tilde{\Sigma}_t\cap\{r>R\}}r^{p-1}(L\psi^s)^2||r^3\lbar\psi^{s/2}||_{L^\infty(t)}^2 \\
&\lesssim B_p^s(t)E^{s/2+7}(t).
\end{align*}
One important step in this calculation is the gain of $r^4$ by passing to the volume form on $\tilde{\Sigma}_t$. However, this same step also loses a factor of $\sin^4\theta$.

Likewise,
\begin{align*}
\int_{\Sigma_t\cap\{r>R\}}r^{p+1}r^8\sin^8\theta (\lbar\psi^s)^2(L\psi^{s/2})^2
&\lesssim \int_{\Sigma_t\cap\{r>R\}}(\lbar\psi^s)^2r^4\sin^4\theta (r^{(p-1)/2+3}L\psi^{s/2})^2 \\
&\lesssim \int_{\Sigma_t\cap\{r>R\}}(\lbar\psi^s)^2r^4\sin^4\theta ||r^{(p-1)/2+3}L\psi^{s/2}||_{L^\infty(t)}^2 \\
&\lesssim \int_{\tilde{\Sigma}_t\cap\{r>R\}}(\lbar\psi^s)^2 ||r^{(p-1)/2+3}L\psi^{s/2}||_{L^\infty(t)}^2 \\
&\lesssim E^s(t)B_p^{s/2+7}(t).
\end{align*}

So far, we have assumed that each of the $\Gamma$ operators in the expression $\Gamma^{\le s}(r^4\sin^4\theta L\psi\lbar\psi)$ acted on either $L\psi$ or $\lbar\psi$. But more generally,
$$\Gamma^s(r^4\sin^4\theta L\psi\lbar\psi) \approx \sum_{i+j+k+s_1+s_2\le s} r^4\Omega^i(\sin^4\theta)L\Omega^j\psi^{s_1}\lbar\Omega^k\psi^{s_2}.$$
At first glance, a few of these terms may be concerning, because a factor of $\sin^2\theta$ was required for the conversion from an integral over $\Sigma_t$ to an integral over $\tilde{\Sigma}_t$. Note, however, that $\max(j,k)\le s-i$ so that Lemma \ref{gain_sin_lem} can be applied to the high order term. This is illustrated in the example below.
\begin{align*}
\int_{\Sigma_t\cap\{r>R\}}&r^{p+1}r^8(L\Omega^l\psi^{s-2-l})^2(\lbar\Omega^{l'}\psi^{s/2-{l'}})^2 \\
&\lesssim E^{s/2+7}(t)\int_{\Sigma_t\cap\{r>R\}}r^{p-1}(L\Omega^l\psi^{s-2-l})^2r^4 \\
&\lesssim E^{s/2+7}(t)\int_{\Sigma_t\cap\{r>R\}}r^{p-1}(L\Omega^{l+2}\psi^{s-2-l})^2r^4\sin^4\theta \\
&\lesssim B_p^s(t)E^{s/2+7}(t).
\end{align*}
\end{proof}

Now we examine the terms $L\phi\lbar\psi$ and $\lbar\phi L\psi$, which arise from the term $\pd^\alpha \phi\pd_\alpha \psi$, which shows up in $\mathcal{N}_\psi$.
\begin{lemma} (Example $\mathcal{N}_\psi$ term)
Assume the following estimates.
$$||r^{(p-1)/2+1}L\Omega^l\phi^{s-l}||^2_{L^\infty(t)}+||r^{(p-1)/2+3}L\Omega^l\psi^{s-l}||^2_{L^\infty(t)}\lesssim B_p^{s+7}(t),$$
$$||r\lbar\Omega^l\phi^{s-l}||^2_{L^\infty(t)}+||r^{3}\lbar\Omega^l\psi^{s-l}||^2_{L^\infty(t)}\lesssim E^{s+7}(t).$$
Then
$$\int_{\tilde{\Sigma}_t\cap\{r>R\}}r^{p+1}\left(\Gamma^s(L\phi\lbar\psi+\lbar\phi L\psi)\right)^2\lesssim B_p^s(t)E^{s/2+7}(t)+E^s(t)B_p^{s/2+7}(t).$$
\end{lemma}
\begin{proof}
We examine four separate cases, each one depending on whether the $\phi$ factor or the $\psi$ factor has the highest number of derivatives and also whether $L$ acts on $\phi$ or $\psi$.
\begin{align*}
\int_{\tilde{\Sigma}_t\cap\{r>R\}}r^{p+1}(L\Omega^l\psi^{s-l})^2(\lbar\Omega^l\phi^{s/2-l})^2
&\lesssim \int_{\tilde{\Sigma}_t\cap\{r>R\}}r^{p-1}(L\Omega^l\psi^{s-l})^2(r\lbar\Omega^l\phi^{s/2-l})^2 \\
&\lesssim \int_{\tilde{\Sigma}_t\cap\{r>R\}}r^{p-1}(L\Omega^l\psi^{s-l})^2||r\lbar\Omega^l\phi^{s/2-l}||_{L^\infty(t)}^2 \\
&\lesssim B_p^s(t)E^{s/2+7}(t).
\end{align*}
\begin{align*}
\int_{\tilde{\Sigma}_t\cap\{r>R\}}r^{p+1}(L\Omega^l\phi^{s-l})^2(\lbar\Omega^l\psi^{s/2-l})^2
&\lesssim \int_{\tilde{\Sigma}_t\cap\{r>R\}}r^{p-1}(L\Omega^l\phi^{s-l})^2r^{-4}(r^3\lbar\Omega^l\psi^{s/2-l})^2 \\
&\lesssim \int_{\tilde{\Sigma}_t\cap\{r>R\}}r^{p-1}(L\Omega^l\phi^{s-l})^2r^{-4}||r^3\lbar\Omega^l\psi^{s/2-l}||_{L^\infty(t)}^2 \\
&\lesssim \int_{\Sigma_t\cap\{r>R\}}r^{p-1}(L\Omega^l\phi^{s-l})^2||r^3\lbar\Omega^l\psi^{s/2-l}||_{L^\infty(t)}^2 \\
&\lesssim B_p^s(t)E^{s/2+7}(t).
\end{align*}

\begin{align*}
\int_{\tilde{\Sigma}_t\cap\{r>R\}}r^{p+1}(\lbar\Omega^l\psi^{s-l})^2(L\Omega^l\phi^{s/2-l})^2
&\lesssim \int_{\tilde{\Sigma}_t\cap\{r>R\}}(\lbar\Omega^l\psi^{s-l})^2(r^{(p-1)/2+1}L\Omega^l\phi^{s/2-l})^2 \\
&\lesssim \int_{\tilde{\Sigma}_t\cap\{r>R\}}(\lbar\Omega^l\psi^{s-l})^2||r^{(p-1)/2+1}L\Omega^l\phi^{s/2-l}||_{L^\infty(t)}^2 \\
&\lesssim E^s(t)B_p^{s/2+7}(t).
\end{align*}
\begin{align*}
\int_{\tilde{\Sigma}_t\cap\{r>R\}}r^{p+1}(\lbar\Omega^l\phi^{s-l})^2(L\Omega^l\psi^{s/2-l})^2
&\lesssim \int_{\tilde{\Sigma}_t\cap\{r>R\}}(\lbar\Omega^l\phi^{s-l})^2r^{-4}(r^{(p-1)/2+3}L\Omega^l\psi^{s/2-l})^2 \\
&\lesssim \int_{\tilde{\Sigma}_t\cap\{r>R\}}(\lbar\Omega^l\phi^{s-l})^2r^{-4}||r^{(p-1)/2+3}L\Omega^l\psi^{s/2-l}||_{L^\infty(t)}^2 \\
&\lesssim \int_{\Sigma_t\cap\{r>R\}}(\lbar\Omega^l\phi^{s-l})^2||r^{(p-1)/2+3}L\Omega^l\psi^{s/2-l}||_{L^\infty(t)}^2 \\
&\lesssim E^s(t)B_p^{s/2+7}(t).
\end{align*}
These four cases illustrate the entirety of the proof. 
\end{proof}

\subsection{The structures of $\mathcal{N}_\phi$ and $\mathcal{N}_\psi$}\label{ws:nl_structure_sec}

We now categorize each term in $\mathcal{N}_\phi$ and $\mathcal{N}_\psi$ so the previous examples can be generalized systematically. We begin with the definition of the null condition.

\begin{definition} We define two families of terms.
$$\gamma = \{L\phi,r^{-1}\pd_r\phi,\sla\nabla\phi,r^{-1}\phi,r^2L\psi,r\pd_r\psi,r^2\sla\nabla\psi,r\psi\}$$
$$\beta = \{L\phi,\pd_r\phi,\sla\nabla\phi,r^{-1}\phi,r^2L\psi,r^2\pd_r\psi,r^2\sla\nabla\psi,r\psi\}$$
The null condition states that any nonlinear term must be a product with at least one $\gamma$ factor.
\end{definition}

\begin{lemma}\label{ws:N_phi_structure_lem} (Structure of $\mathcal{N}_\phi$) The nonlinear term $\mathcal{N}_\phi$ can be expressed as a sum of terms of the form
$$\frac{f\gamma\beta(r\beta)^{j_0}}{1+\phi}$$
satisfying the following additional rules. \\
i) The factor $f$ is smooth and bounded. \\
ii) If $\psi$ appears at least once in the term, then $f$ has a factor of $\sin^2\theta$.
\end{lemma}
\begin{proof}
In this proof, we use the sign $\approx$ to emphasize that smooth and bounded factors (including $M/r$) are neglected. Recall that
\begin{multline*}
(1+\phi)\mathcal{N}_\phi = \pd^\alpha\phi\pd_\alpha\phi \\
- \left[A^2(1+\phi)^4 \pd^\alpha\psi\pd_\alpha\psi +4A^2(1+\phi)^3\psi\pd^\alpha\phi\pd_\alpha\psi +4A^2(1+\phi)^2\psi^2\pd^\alpha\phi\pd_\alpha\phi\right. \\
\left.4A\pd^\alpha A (1+\phi)^4 \psi\pd_\alpha\psi +8 A\pd^\alpha A (1+\phi)^3\psi^2\pd_\alpha\phi +4\pd^\alpha A\pd_\alpha A (1+\phi)^4\psi^2\right]
\end{multline*}
and $A=r^2\sin^2\theta$.

We investigate each term in $(1+\phi)\mathcal{N}_\phi$.
\begin{align*}
\pd^\alpha \phi\pd_\alpha \phi &\approx L\phi\lbar \phi + \sla\nabla\phi\cdot \sla\nabla\phi \\
&\approx (L\phi)(\lbar\phi)+(\sla\nabla\phi)(\sla\nabla\phi) \\
&\approx (L\phi)(\pd_r\phi)+(L\phi)(L\phi)+(\sla\nabla\phi)(\sla\nabla\phi).
\end{align*}
Some of the remaining terms will have factors of $\phi$ coming from the expansion of powers of $(1+\phi)$. Since $\phi = r(r^{-1}\phi)=r\beta$, this accounts for some factors of $r\beta$. In the remaining terms, we will ignore factors that are powers of $(1+\phi)$.
\begin{align*}
A^2\pd^\alpha\psi \pd_\alpha \psi &\approx r^4\sin^4\theta (L\psi \lbar\psi+\sla\nabla\psi\cdot \sla\nabla\psi) \\
&\approx \sin^4\theta (r^2L\psi)(r^2\lbar\psi)+\sin^4\theta(r^2\sla\nabla\psi)(r^2\sla\nabla\psi) \\
&\approx \sin^4\theta (r^2L\psi)(r^2\pd_r\psi)+\sin^4\theta (r^2L\psi)(r^2L\psi)+\sin^4\theta(r^2\sla\nabla\psi)(r^2\sla\nabla\psi).
\end{align*}
\begin{align*}
A^2\psi\pd^\alpha\phi\pd_\alpha\psi &\approx r^4\sin^4\theta (\psi L\phi \lbar\psi +\psi L\psi\lbar\phi +\psi \sla\nabla\phi\cdot\sla\nabla\psi) \\
&\approx \sin^4\theta (L\phi) (r^2\lbar \psi)(r^2\psi) +\sin^4\theta (r^2 L\psi)(\lbar\phi)(r^2\psi) \\
&\hspace{.5in}+ \sin^4\theta(\sla\nabla\phi)(r^2\sla\nabla\psi)(r^2\psi) \\
&\approx \sin^4\theta (L\phi)(r^2\pd_r\psi)(r^2\psi) +\sin^4\theta(L\phi)(r^2 L\psi)(r^2\psi) \\
&\hspace{.5in} +\sin^4\theta (r^2L\psi)(\pd_r\phi)(r^2\psi) + \sin^4\theta (r^2L\psi)(L\phi)(r^2\psi) \\
&\hspace{.5in} +\sin^4\theta (\sla\nabla\phi)(r^2\sla\nabla\psi)(r^2\psi)
\end{align*}
\begin{align*}
A^2\psi^2\pd^\alpha\phi\pd_\alpha\phi &\approx r^4\sin^4\theta \psi^2\pd^\alpha\phi\pd_\alpha\phi \\
&\approx \sin^4\theta \pd^\alpha\phi\pd_\alpha\phi (r^2\psi)^2 \\
&\approx \sin^4\theta (L\phi)(\pd_r\phi)(r^2\psi)^2+\sin^4\theta(L\phi)(L\phi)(r^2\psi)^2+\sin^4\theta(\sla\nabla\phi)(\sla\nabla\phi)(r^2\psi)^2
\end{align*}
\begin{align*}
A\pd^\alpha A \psi\pd_\alpha\psi &\approx r^3\sin^4\theta \psi \pd_r\psi + r^4\sin^2\theta \sla\nabla(\sin^2\theta) \psi \sla\nabla\psi \\
&\approx \sin^4\theta (r\psi)(r^2\pd_r\psi) + \sin^2\theta r\sla\nabla(\sin^2\theta) (r^2\sla\nabla \psi)(r\psi)
\end{align*}
\begin{align*}
A\pd^\alpha A \psi^2\pd_\alpha\phi &\approx r^3\sin^4\theta \psi^2\pd_r\phi+r^4\sin^2\theta\sla\nabla(\sin^2\theta)\psi^2\sla\nabla\phi \\
&\approx \sin^4\theta (r\psi)(\pd_r\phi)(r^2\psi) + \sin^2\theta r\sla\nabla(\sin^2\theta)(\sla\nabla\phi)(r\psi)(r^2\psi)
\end{align*}
\begin{align*}
\pd^\alpha A\pd_\alpha A \psi^2 &\approx r^2\sin^4\theta\psi^2 + r^2\sin^2\theta\psi^2 \\
&\approx \sin^4\theta (r\psi)(r\psi) +\sin^2\theta (r\psi)(r\psi)
\end{align*}
It is now straightforward to check for each of the above calculations that the first factor in parentheses is a $\gamma$ term, the second factor is a $\beta$ term, and any additional factors are $r\beta$ terms. It is also straightforward to check that any term containing a $\psi$ factor also has a factor of $\sin^2\theta$.
\end{proof}

\begin{lemma}\label{ws:N_psi_structure_lem} (Structure of $\mathcal{N}_\psi$) The nonlinear term $\mathcal{N}_\psi$ can be expressed as a sum of terms of the form
$$\frac{r^{-2}f\gamma\beta (r\beta)^{j_0}}{(1+\phi)^2},$$
where the factor $f$ is smooth and bounded.
\end{lemma}
\begin{proof}
Just as in the proof of Lemma \ref{ws:N_phi_structure_lem}, we use the sign $\approx$ to emphasize that smooth and bounded factors (including $M/r$) are neglected. Recall that
\begin{align*}
(1+\phi)^2\mathcal{N}_\psi &= 2\psi\left[A^2(1+\phi)^4 \pd^\alpha\psi\pd_\alpha\psi +4A^2(1+\phi)^3\psi\pd^\alpha\phi\pd_\alpha\psi +4A^2(1+\phi)^2\psi^2\pd^\alpha\phi\pd_\alpha\phi\right. \\
&\hspace{.5in}\left.4A\pd^\alpha A (1+\phi)^4 \psi\pd_\alpha\psi +8 A\pd^\alpha A (1+\phi)^3\psi^2\pd_\alpha\phi +4\pd^\alpha A\pd_\alpha A (1+\phi)^4\psi^2\right] \\
&\hspace{.3in}-2(1+\phi)\pd^\alpha\phi\pd_\alpha\psi
\end{align*}
Each of the terms in square brackets was examined in the proof of Lemma \ref{ws:N_phi_structure_lem} and found to have null structure. Since these terms are multiplied by an additional factor of $\psi=r^{-2}(r\beta)$, this accounts for the additional factor of $r^{-2}$ while also contrbuting a factor of $r\beta$. It suffices to check the final term. For the same reason as in the proof of Lemma \ref{ws:N_phi_structure_lem}, we ignore the factor of $1+\phi$.
\begin{align*}
\pd^\alpha\phi\pd_\alpha\psi &\approx L\phi \lbar\psi + L\psi\lbar\phi +\sla\nabla\phi\cdot\sla\nabla\psi \\
&\approx r^{-2}(L\phi)(r^2\lbar\psi) +r^{-2}(r^2L\psi)(\lbar\phi) + r^{-2}(\sla\nabla\phi)(r^2\sla\nabla\psi) \\
&\approx r^{-2}(L\phi)(r^2\pd_r\psi)+r^{-2}(L\phi)(r^2L\psi)+r^{-2}(r^2L\psi)(\pd_r\phi)+r^{-2}(r^2L\psi)(L\phi) \\
&\hspace{.5in}+r^{-2}(\sla\nabla\phi)(r^2\sla\nabla\psi)
\end{align*}
It is now straightforward to check for each of the above terms that the first factor in parentheses is a $\gamma$ term, and the second factor is a $\beta$ term. Also, all of these terms have an additional factor of $r^{-2}$.
\end{proof}

\subsection{The structures of $\tg^k\Gamma^s\mathcal{N}_\phi$ and $\tg^k\tilde\Gamma^s\mathcal{N}_\psi$}\label{ws:nl_s_k_structure_sec}

\begin{definition} We generalize the previous families of terms to the following.
\begin{multline*}
\gamma^{s,k} = \{L\Omega^l\phi^{s-l,k},r^{-1}\pd_r\Omega^l\phi^{s-l,k},r^{-1}\pd_\theta\Omega^l\phi^{s-l,k},r^{-1}\Omega^l\phi^{s-l,k}, \\
r^2L\Omega^l\psi^{s-l,k},r\pd_r\Omega^l\psi^{s-l,k},r\pd_\theta\Omega^l\psi^{s-l,k},r\Omega^l\psi^{s-l,k}\}
\end{multline*}
\begin{multline*}
\beta^{s,k} = \{L\Omega^l\phi^{s-l,k},\pd_r\Omega^l\phi^{s-l,k},r^{-1}\pd_\theta\Omega^l\phi^{s-l,k},r^{-1}\Omega^l\phi^{s-l,k}, \\
r^2L\Omega^l\psi^{s-l,k},r^2\pd_r\Omega^l\psi^{s-l,k},r\pd_\theta\Omega^l\psi^{s-l,k},r\Omega^l\psi^{s-l,k}\}
\end{multline*}
The null condition still states that any nonlinear term must be a product with at least one $\gamma$ factor.
\end{definition}

The following two lemmas (one for $\tg^k\Gamma^s\mathcal{N}_\phi$ and one for $\tg^k\tilde{\Gamma}^s\mathcal{N}_\psi$) generalize the lemmas in \S\ref{ws:nl_structure_sec}.

\begin{lemma}\label{ws:Nphi_s_k_structure_lem}
(Structure of $\tg^k\Gamma^s\mathcal{N}_\phi$) The nonlinear term $\tg^k\Gamma^s\mathcal{N}_\phi$ can be expressed as a sum of terms of the following form. (The integer $j$ represents the number of times a differential operator acts on the denominator.) \\
$$\frac{f\alpha^{s_1,k_1}\beta^{s_2,k_2}(r\beta^{s_3,k_3})...(r\beta^{s_{2+j_0+j},k_{2+j_0+j}})}{(1+\phi)^{j+1}},$$
where $0\le j\le s+k$, $s_1+...+s_{2+j_0+j}\le s$, $k_1+...+k_{2+j_0+j}\le k$, and the following rules apply. \\
\bp The factor $f$ is smooth and bounded. \\
\bp If $\psi$ appears at least once in the term, then $f$ has a factor of $(\sin\theta)^{\max(0,2-\min_i(s-s_i))}$.
\end{lemma}
\begin{proof}
We start with the case $s=k=0$, which was proved in Lemma \ref{ws:N_phi_structure_lem}.  Each term in $\mathcal{N}_\phi$ can be written in the form
$$\frac{f\gamma \beta(r\beta)^{j_0}}{1+\phi}.$$
The effect of applying $\pd_t^{s'}$ to this type of term is to obtain terms of the form
$$\frac{f\pd_t^{s_1}\gamma\pd_t^{s_2}\beta (r\pd_t^{s_3}\beta)...(r\pd_t^{s_{2+j_0+j}}\beta)}{(1+\phi)^{1+j}}$$
where $s_1+...+s_{2+j_0+j}\le s'$ and $j\le s'$. The additional factors of $r\pd_t^{s_i}\beta$ $(i>s+j_0)$ appear each time one of the $\pd_t$ operators acts on the denominator.

Recall that $\tg^k\Gamma^s$ is composed not only of $\pd_t$, but also $Q$ and $\tg$. Since these operators commute with each other, the order in which they are applied is not important. For the sake of simplicity, we first apply the $\pd_t^{s'}$ operators (which has already been done) and then the $\tg^k$ operators. Since these are both first order operators, it should be clear that the resulting terms are of the form
$$\frac{f\pd_t^{s_1,k_1}\gamma_0\pd_t^{s_2,k_2}\beta_0 (r\pd_t^{s_3,k_3}\beta_0)...(r\pd_t^{s_{2+j_0+j},k_{2+j_0+j}}\beta_0)}{(1+\phi)^{1+j}}$$
where $j\le s'+k$, $s_1+...+s_{2+j_0+j}\le s'$, $k_1+...+k_{2+j_0+j}\le k$, and $\gamma_0^{s_1,k_1}=\tg^{k_1}\pd_t^{s_1}\gamma$ and $\beta_0^{s_i,k_i}=\tg^{k_i}\pd_t^{s_i}\beta$.

This is the point where we apply the operator $Q$. Since
$$Q=\Omega\cdot\Omega,$$
it follows that each application of $Q$ will introduce $\Omega$ operators to each of the terms. Since in general,
$$\gamma^{s,k} = \Omega^l \gamma_0^{s-l,k}$$
$$\beta^{s,k} = \Omega^l \beta_0^{s-l,k}$$
the result after applying $Q^{s-s'}$ is to obtain terms of the form
$$\frac{f\alpha^{s_1,k_1}\beta^{s_2,k_2}(r\beta^{s_3,k_3})...(r\beta^{s_{2+j_0+j},k_{2+j_0+j}})}{(1+\phi)^{j+1}},$$
where $j\le s+k$, $s_1+...+s_{2+j_0+j}\le s$, and $k_1+...+k_{2+j_0+j}\le k$.

Finally, we observe that in the $s=k=0$ case, if $\psi$ appears at least once in the term then $f$ has a factor of $\sin^2\theta$. In general, when applying $Q$, it is possible that $\Omega$ could act on this factor, reducing it to $\sin\theta$ or $1$. By counting how many times this happens, it is clear that the resulting factor is at worst $(\sin\theta)^{\max(0,2-\min_i(s-s_i))}$.
\end{proof}

\begin{lemma}\label{ws:Npsi_s_k_structure_lem}
(Structure of $\tg^k\tilde\Gamma^s\mathcal{N}_\psi$) The nonlinear term $\tg^k\tilde\Gamma^s\mathcal{N}_\psi$ can be expressed as a sum of terms of the following form. (The integer $j$ represents the number of times a differential operator acts on the denominator.) \\
$$\frac{r^{-2}f\alpha^{s_1,k_1}\beta^{s_2,k_2}(r\beta^{s_3,k_3})...(r\beta^{s_{2+j_0+j},k_{2+j_0+j}})}{(1+\phi)^{j+2}}$$
where $0\le j\le s+k$, $s_1+...+s_{2+j_0+j}\le s$, $k_1+...+k_{2+j_0+j}\le k$, and the factor $f$ is smooth and bounded.
\end{lemma}
\begin{proof}
The proof is the same as the proof of Lemma \ref{ws:Nphi_s_k_structure_lem}, except that it starts with Lemma \ref{ws:N_psi_structure_lem}.
\end{proof}

\subsection{The general strategy revisited}\label{ws:nl_strategy_revisited_sec}

Now that the structures of $\tg^k\Gamma^s\mathcal{N}_\phi$ and $\tg^k\tilde\Gamma^s\mathcal{N}_\psi$ have been determined, the strategy that was illustrated in \S\ref{ws:example_terms_sec} can be summarized in the following proposition.

\begin{proposition}\label{ws:strategy_revisited_prop}
Suppose the following estimates hold.
$$(1+\phi)^{-1}\lesssim 1,$$
$$\int_R^\infty\int_0^\pi (r^{(p-1)/2}\gamma^s)^2 r^2\sin^5\theta d\theta dr \lesssim E_{p-1}^s(t),$$
$$||r\beta^s||_{L^\infty(\Sigma_t\{r>R\})}^2\lesssim E^{s+7}(t),$$
$$\int_R^\infty\int_0^\pi (\beta^s)^2r^2\sin^5\theta d\theta dr\lesssim E^s(t),$$
$$||r^{(p-1)/2+1}\gamma^s||_{L^\infty(\Sigma_t\cap\{r>R\})}^2\lesssim E_{p-1}^{s+7}(t).$$
Then
\begin{multline*}
\int_{\Sigma_t\cap\{r>R\}}r^{p+1}|\tg^k\Gamma^s\mathcal{N}_\phi|^2 + \int_{\tilde{\Sigma}_t\cap\{r>R\}}r^{p+1}|\tg^k\tilde{\Gamma}^s\mathcal{N}_\psi|^2 \\
\lesssim \left(E_{p-1}^s(t)E^{s/2+7}(t)+E^s(t)E_{p-1}^{s/2+7}(t)\right)\sum_{j\le s+j_0}(E^{s/2+7}(t))^j.
\end{multline*}
\end{proposition}
\begin{remark}
The second and fourth assumptions are automatically true. The first, third, and fifth assumptions will be true in the context of the bootstrap assumptions in the proof of the main theorem (Theorem \ref{ws:main_thm}).
\end{remark}
\begin{proof}
Note that Lemmas \ref{ws:Nphi_s_k_structure_lem} and \ref{ws:Npsi_s_k_structure_lem} both had the requirement $s_1+...+s_{2+j_0+j}\le s$. This means there can be at most one factor (denoted $\gamma^{hi}$ or $\beta^{hi}$) with up to $s$ derivatives and all of the remaining factors (denoted $\gamma^{lo}$ or $\beta^{lo}$) have at most $s/2$ derivatives. Since the integrals are taken over the range $r>R$, where $\tg=0$, we are free to ignore all $k$ indices. 

The strategy in this proof is rather simple. All of the factors that have at most $s/2$ derivatives are estimated in $L^\infty$, while the remaining factor with at most $s$ derivatives is estimated in $L^2$. The procedure then takes one of two possible directions, depending on whether the high derivative term is a $\gamma$ term or a $\beta$ term.

\textbf{First, we estimate the integral over $\Sigma_t$, using the form given in Lemma \ref{ws:Nphi_s_k_structure_lem}.}
\begin{multline*}
\int_{\Sigma_t\cap\{r>R\}}r^{p+1}|\tg^k\Gamma^s\mathcal{N}_\phi|^2 \\
\lesssim \sum_{0\le j\le s+j_0} \int_{\Sigma_t\cap\{r>R\}}r^{p+1}\left(f\gamma^{hi}\beta^{lo}(r\beta^{lo})^j\right)^2 + \sum_{0\le j\le s+j_0} \int_{\Sigma_t\cap\{r>R\}}r^{p+1}\left(f\beta^{hi}\gamma^{lo}(r\beta^{lo})^j\right)^2 .
\end{multline*}
Now,
\begin{align*}
\int_{\Sigma_t\cap\{r>R\}}r^{p+1}\left(f\gamma^{hi}\beta^{lo}(r\beta^{lo})^j\right)^2
&\lesssim \int_{\Sigma_t\cap\{r>R\}}r^{p-1}\left(f\gamma^{hi}(r\beta^{lo})(r\beta^{lo})^j\right)^2 \\
&\lesssim ||r\beta^{lo}||_{L^\infty}^{2(j+1)}\int_{\Sigma_t\cap\{r>R\}}r^{p-1}(f\gamma^{hi})^2 \\
&\lesssim (E^{s/2+7}(t))^{j+1}\int_{\Sigma_t\cap\{r>R\}}r^{p-1}(f\gamma^{hi})^2 \\
&\lesssim E_{p-1}^s(t)(E^{s/2+7}(t))^{j+1}.
\end{align*}
The last step requires further justification, because the $\gamma$ term could possibly be a $\psi$ term. But in that case, Lemma \ref{ws:Nphi_s_k_structure_lem} also states that $f$ has an additional factor of $\sin^2\theta$ (or in the case where $\gamma$ does not have exactly $s$ derivatives, this factor might be either $\sin\theta$ or $1$, but we can apply Lemma \ref{gain_sin_lem}--see the end of the proof of the example Lemma \ref{ws:example_Nphi_term_lem}).

Also,
\begin{align*}
\int_{\Sigma_t\cap\{r>R\}}r^{p+1}\left(f\beta^{hi}\gamma^{lo}(r\beta^{lo})^j\right)^2 
&\lesssim \int_{\Sigma_t\cap\{r>R\}}\left(f\beta^{hi}(r^{(p+1)/2}\gamma^{lo})(r\beta^{lo})^j\right)^2 \\
&\lesssim ||r^{(p+1)/2}\gamma^{lo}||_{L^\infty}^2||r\beta||_{L^\infty}^{2j}\int_{\Sigma_t\cap\{r>R\}}(f\beta^{hi})^2 \\
&\lesssim E_{p-1}^{s/2+7}(t)(E^{s/2+7}(t))^j\int_{\Sigma_t\cap\{r>R\}}(f\beta^{hi})^2 \\
&\lesssim E^s(t)E_{p-1}^{s/2+7}(t)(E^{s/2+7}(t))^j.
\end{align*}
The last step requires the same justification given in the preceeding calculation.

Combining both estimates, we conclude that 
$$\int_{\Sigma_t\cap\{r>R\}}r^{p+1}|\tg^k\Gamma^s\mathcal{N}_\phi|^2 \lesssim \left(E_{p-1}^s(t)E^{s/2+7}(t)+E^s(t)E_{p-1}^{s/2+7}(t)\right)\sum_{j\le s+j_0}(E^{s/2+7}(t))^j.$$

\textbf{Next, we estimate the integral over $\tilde{\Sigma}_t$, using the form given in Lemma \ref{ws:Npsi_s_k_structure_lem}.} The form given in Lemma \ref{ws:Npsi_s_k_structure_lem} is very similar to the form given in Lemma \ref{ws:Nphi_s_k_structure_lem}, so we will not repeat the estimates in detail. The only difference is the presence of an additional factor of $r^{-2}$ that comes along with $f$. This is required in the two final steps.
$$\int_{\tilde{\Sigma}_t\cap\{r>R\}}r^{p-1}(r^{-2}f\gamma^{hi})^2\lesssim E_{p-1}^s(t)$$
or
$$\int_{\tilde{\Sigma}_t\cap\{r>R\}}(r^{-2}f\beta^{hi})^2\lesssim E^s(t).$$
The factor of $r^{-2}$ (inside the square) is needed to offset the additional factor of $r^4$ in the volume form for $\tilde{\Sigma}_t$.

We conclude that
$$\int_{\tilde{\Sigma}_t\cap\{r>R\}}r^{p+1}|\tg^k\tilde{\Gamma}^s\mathcal{N}_\psi|^2 \lesssim \left(E_{p-1}^s(t)E^{s/2+7}(t)+E^s(t)E_{p-1}^{s/2+7}(t)\right)\sum_{j\le s}(E^{s/2+7}(t))^j.$$
This completes the proof. 
\end{proof}

\section{Theorem: Global boundedness and decay for axisymmetric perturbations of the nontrivial solution to the wave map problem from Schwarzschild to the hyperbolic plane preserving angular momentum}\label{ws:main_thm_sec}

In this final section, we state the main theorem of Chapter \ref{wm_szd_chap} and summarize its proof, which is very similar to the proof of the main theorem of Chapter \ref{kerr_chap}.

\subsection{The main theorem}\label{ws:main_thm_statement_sec}

\begin{theorem}\label{ws:main_thm}
Let
\begin{align*}
X&=A+A\phi, \\
Y&=X^2\psi.
\end{align*}
(In particular, the assumption $Y=O(\sin^4\theta)$ near the axis excludes any perturbations corresponding to a change in angular momentum in the sense described in \S\ref{ernst_potential_axis_derivation}.) 

Suppose the pair $(X,Y)$ is axisymmetric and satisfies the wave map system
\begin{align*}
X\Box_g X&=\pd^\alpha X\pd_\alpha X-\pd^\alpha Y\pd_\alpha Y \\
X\Box_g Y&=2\pd^\alpha X\pd_\alpha Y,
\end{align*}
where $g$ is the Schwarzschild metric.

Define the energies
$$E^{\ul{n}}(t)=\sum_{s+2k=n}E^{s,k}(t),$$
$$E_p^{\ul{n}}(t)=\sum_{s+2k=n}E_p^{s,k}(t).$$
Then for $\delp,\delm>0$ sufficiently small, if the initial data for $(\phi,\psi)$ decay sufficiently fast as $r\rightarrow \infty$ and have size
\begin{equation*}
I_0=E^{\ul{26}}(0)+E^{\ul{26}}_{2-\delp}(0)
\end{equation*}
sufficiently small, then the following estimates hold for $t\ge 0$ (with $T=1+t$).

I) The energies satisfy
$$E^{\ul{26}}(t)\lesssim I_0$$
$$E^{\ul{26}}_{p\in[\delm,2-\delp]}(t)\lesssim I_0$$
$$E^{\ul{24}}_{p\in[1-\delp,2-\delp]}(t)\lesssim T^{p-2+\delp}I_0$$
$$E^{\ul{22}}_{p\in[\delm,2-\delp]}(t)\lesssim T^{p-2+\delp}I_0$$
$$\int_t^{\infty}E^{\ul{20}}_{p\in[\delm-1,\delm]}(\tau)d\tau\lesssim T^{p-2+\delp+1}I_0$$

II) For all $s,k$ such that $s+2k\le 26$, the following $L^\infty$ estimates hold.
\begin{multline*}
|r^{p+1}\bar{D}\Omega^l\phi^{s-l,k}|^2+|r^pD\Omega^l\phi^{s-l,k}|^2+|r^p\Omega^l\phi^{s-l,k}|^2 \\
+|r^{p+3}\bar{D}\Omega^l\psi^{s-l,k}|^2+|r^{p+2}D\Omega^l\psi^{s-l,k}|^2+|r^{p+2}\Omega^l\psi^{s-l,k}|^2 \\
\lesssim E_{2p}^{s+5,k+1}(t)+E_{2p}^{s+7,k}(t)
\end{multline*}

III) Together, (I) and (II) imply that if $s+2k\le 13$, for all $p\in [\delm/2,(2-\delp)/2]$,
\begin{multline*}
|r^{p+1}\bar{D}\Omega^l\phi^{s-l,k}|+|r^pD\Omega^l\phi^{s-l,k}|+|r^p\Omega^l\phi^{s-l,k}| \\
+|r^{p+3}\bar{D}\Omega^l\psi^{s-l,k}|+|r^{p+2}D\Omega^l\psi^{s-l,k}|+|r^{p+2}\Omega^l\psi^{s-l,k}| \\
\lesssim T^{(2p-2+\delp)/2}I_0^{1/2}
\end{multline*}
and additionally for $p\in [(\delm-1)/2,\delm/2]$,
\begin{multline*}
\int_t^\infty |r^{p+1}\bar{D}\Omega^l\phi^{s-l,k}|+|r^pD\Omega^l\phi^{s-l,k}|+|r^p\Omega^l\phi^{s-l,k}| \\
+\int_t^\infty |r^{p+3}\bar{D}\Omega^l\psi^{s-l,k}|+|r^{p+2}D\Omega^l\psi^{s-l,k}|+|r^{p+2}\Omega^l\psi^{s-l,k}| \\
\lesssim T^{(2p-2+\delp)/2+1}I_0^{1/2}.
\end{multline*}
The final estimate should be interpreted as saying that $|r^{(\delm-1)/2+1}\bar{D}\Omega^l\phi^{s-l,k}|$, \\
$|r^{(\delm-1)/2}D\Omega^l\phi^{s-l,k}|$, $|r^{(\delm-1)/2}\Omega^l\phi^{s-l,k}|$, $|r^{(\delm-1)/2+3}\bar{D}\Omega^l\psi^{s-l,k}|$, $|r^{(\delm-1)/2+2}D\Omega^l\psi^{s-l,k}|$, and $|r^{(\delm-1)/2+2}\Omega^l\psi^{s-l,k}|$ decay like $T^{(\delm-3+\delp)/2}$ in a weak sense.
\end{theorem}

\subsection{A summary of the proof}

The proof of Theorem \ref{ws:main_thm} is very similar to the proof the main theorem of the previous chapter (Theorem \ref{k:main_thm}) and is not worth repeating entirely. Instead, a summary of the proof is given below.

The proof begins as in \S\ref{k:bootstrap_assumptions_sec} with the following bootstrap assumptions,
\begin{align*}
E^{\ul{26}}(t)&\le C_bI_0, \\
E^{\ul{26}}_{p\in[\delm,2-\delp]}(t)&\le C_b I_0, \\
E^{\ul{22}}_{p\in[\delm,2-\delp]}(t)&\le C_b T^{p-2+\delp}I_0, \\
\int_t^\infty E^{\ul{20}}_{p\in[\delm-1,\delm]}(\tau)d\tau &\le C_b T^{p-2+\delp+1}I_0.
\end{align*}
as well as the additional bootstrap assumption
\begin{equation}\label{ws:additional_bootstrap_assumption}
|\phi|\le \frac12,
\end{equation}
which is needed to ensure that factors of $(1+\phi)^{-1}$ that show up in the nonlinear terms remain bounded.

With these bootstrap assumptions, the $L^\infty$ estimates from Proposition \ref{ws:infinity_prop} are improved by removing the terms with $\Box_g\phi$ and $\Box_{\tilde{g}}\psi$. The improved estimates are given in the following lemma.
\begin{lemma}\label{ws:improved_pointwise_lem}
In the context of the bootstrap assumptions, the following $L^\infty$ estimates hold for $s+2k\le 26$ and all $p$ in any bounded range.

For $r\ge r_H$,
\begin{multline*}
|r^{p+1}\bar{D}\Omega^l\phi^{s-l,k}|^2+|r^pD\Omega^l\phi^{s-l,k}|^2+|r^p\Omega^l\phi^{s-l,k}|^2 \\
+|r^{p+3}\bar{D}\Omega^l\psi^{s-l,k}|^2+|r^{p+2}D\Omega^l\psi^{s-l,k}|^2+|r^{p+2}\Omega^l\psi^{s-l,k}|^2 \\
\lesssim E_{2p}^{s+5,k+1}(t)+E_{2p}^{s+7,k}(t),
\end{multline*}
and for $r\ge r_0>r_H$,
$$|rD\Omega^l\phi^{s-l,k}|^2+|r^3D\Omega^l\psi^{s-l,k}|^2 \lesssim E^{s+5,k+1}(t)+E^{s+7,k}(t).$$
\end{lemma}

The next step, in which the norms $N^{s,k}(t)$ and $N_p^{s,k}(t)$ are refined, is the key step. (See \S\ref{k:refined_estimates_sec}.) The new work that is required to perform this step has already been done in \S\ref{ws:nl_strategy_revisited_sec}. The result is summarized in the following set of estimates.

\begin{lemma}
In the context of the bootstrap assumptions, if $s+2k\le 22$ and $C_bI_0\le 1$, then
$$N^{s,k}(t)\lesssim (E^{s,k}(t))^{1/2}\left((E^{s,k}(t))^{1/2}(E_{\delm-1}^{\ul{16}}(t))^{1/2}+(E^{s,k}_{1-\delm}(t))^{1/2}(E^{\ul{16}}_{\delm-1}(t))^{1/2}\right)$$
and
\begin{multline*}
E_p^{s,k}(t_2)+\int_{t_1}^{t_2}B_p^{s,k}(t)dt \\
\lesssim E_p^{s,k}(t_1)+(C_bI_0)^{1/2}\int_{t_1}^{t_2}E_p^{s,k}(t)T^{(\delm-3+\delp)/2}dt+\sum_{\substack{s'+2k'\le s+2k \\ k'<k}}\int_{t_1}^{t_2}B_p^{s',k'}(t)dt.
\end{multline*}
\end{lemma}
This lemma contains the same information that was given by Lemma \ref{k:refined_nl_E_lem} and Corollary \ref{k:NL_absorb_bulk} in Chapter \ref{kerr_chap}.

At this point, the proof proceeds the exact same way as in the proof of Theorem \ref{k:main_thm} starting at \S\ref{k:inductive_assumptions_sec}, since the remaining estimates rely only on the information contained in the above lemma.

All that remains is to recover the additional bootstrap assumption (\ref{ws:additional_bootstrap_assumption}) by proving that
$$|\phi|\le \frac14.$$
This is easily done using the following estimate from Lemma \ref{ws:improved_pointwise_lem}
$$||\phi||_{L^\infty(\Sigma_t)}^2 \lesssim E_0^{5,1}(t)+E_0^{7,0}(t) \lesssim I_0,$$
and choosing $I_0$ sufficiently small.

\noindent This completes the proof of Theorem \ref{ws:main_thm}. \qed

\chapter{A class of perturbations of the nontrivial wave map solution in slowly rotating Kerr spacetimes}\label{wm_kerr_chap}

The following conjecture was posed in \cite{ionescu2014global}.
\begin{conjecture}\label{wk:iokl_conj}
The stationary solution
\begin{align*}
X_0 &= A = \frac{(r^2+a^2)^2-a^2\sin^2\theta(r^2-2Mr+a^2)}{r^2+a^2\cos^2\theta}\sin^2\theta \\
Y_0 &= B = -2aM(3\cos\theta-\cos^3\theta)-\frac{2a^3M\sin^4\theta\cos\theta}{r^2+a^2\cos^2\theta}.
\end{align*}
of the wave map system
\begin{align}
\Box_gX &= \frac{\pd^\alpha X\pd_\alpha X}{X}-\frac{\pd^\alpha Y\pd_\alpha Y}{X}\label{wk:X_eqn} \\
\Box_gY &= 2\frac{\pd^\alpha X\pd_\alpha Y}{X}.\label{wk:Y_eqn}
\end{align}
with $g$ any subextremal ($|a|<M$) Kerr metric is future asymptotically stable in the domain of outer communication of the Kerr spacetime for all smooth, axially symmetric, admissible perturbations, ie. perturbations vanishing in some suitable way on the axis of symmetry.
\end{conjecture}

In this final chapter, we resolve the above conjecture for the slowly rotating ($|a|\ll M$) Kerr spacetimes. As discussed in \S\ref{model_problem_motivation_sec}, the motivation for this conjecture is its relation to the stability of Kerr black holes as solutions to the Einstein Vacuum System (\ref{eve}).

\section{Introduction}

When introducing the model problem of Conjecture \ref{wk:iokl_conj}, Ionescu and Klainerman outlined three main difficulties that must be overcome to find a solution.
\begin{enumerate}
\item \textit{Strong linear stability} The main result of \cite{ionescu2014global} is equivalent to the dynamic estimates for the linearized wave map system before applying commutators. These estimates are reproven in this chapter in \S\ref{wk:xi_estimates_sec}-\S\ref{wk:dynamic_estimates_sec} using the methods developed so far in this thesis and are summarized most concisely in Theorem \ref{p_L_thm}. It was perhaps not known at the time  \cite{ionescu2014global} was published that these decay estimates are sufficient to handle the nonlinear problem. (The missing ingredient is the weak decay principle.)
\item \textit{Nonlinear stability} The nonlinear terms must have special structure, including the null condition, which is compatible with the linear decay estimates. In addition to this condition, there is a new structural condition on the axis. See \S\ref{intro_regularity_sec} for a brief description and \S\ref{nl_structure_sec} for a detailed examination.
\item \textit{Degeneracy on the axis} There are difficulties associated with the axisymmetric reduction of the equations, which are manifest in terms with factors that are singular on the axis. To overcome these difficulties, a new formalism detailed in Appendix \ref{regularity_sec} is essential in this chapter and likely will be useful in future works on problems with axisymmetry. Again, see \S\ref{intro_regularity_sec} for a brief description.
\end{enumerate}

\subsection{The $\xi_a$ system of \cite{ionescu2014global}}\label{wk:xi_a_intro_sec}

The model problem (\ref{wk:X_eqn}-\ref{wk:Y_eqn}) represents a wave map from Kerr to the hyperbolic plane. (For more information, see \S\ref{wave_map_general_theory_sec}.) This fact naturally gives rise to a bundle $\mathcal{B}$, of which a section $\xi_a$ describes a linearized perturbation about a given solution. (For an explanation of this fact, see \S\ref{wave_map_bundle_sec}.) One of the main important discoveries in \cite{ionescu2014global} is that this bundle reveals a conserved energy quantity for the linearized system that has good positivity properties, and it also suggests a good grouping of terms when deriving a Morawetz, estimate especially near the trapping radius.

 We take a moment here to define the bundle $\mathcal{B}$.

\begin{definition}
Let $\mathcal{B}$ be a bundle spanned by two basis elements at each point in the Kerr spacetime
$$\mathcal{B}=span\{e^1,e^2\}$$
endowed with a metric $h$ given by
$$h(e^i,e^j)=\delta^{ij},$$
where $\delta^{ij}$ is the Kroneker delta, and a connection $D$ given by
\begin{align*}
D_\mu e^1 &= -\frac{\pd_\mu B}{A} e^2 \\
D_\mu e^2 &= \frac{\pd_\mu B}{A} e^1.
\end{align*}
It is easy to check that the connection $D$ is compatible with the metric $h$.
\end{definition}

Given two scalar functions $\xi_1$ and $\xi_2$, we can write a section of the bundle $\mathcal{B}$ as
$$\xi = \xi_1 e^1+\xi_2 e^2.$$
Introducing the tensor index $a$, we often write
$$\xi_a = \xi_1 (e^1)_a+\xi_2 (e^2)_a.$$

We extend the wave operator $\Box_g$ to act on sections of $\mathcal{B}$
$$\Box_g = D^\mu D_\mu$$
and we define the potential tensor
\begin{align*}
V &= g^{\mu\nu}d\Phi_\mu d\Phi_\nu \\
&= g^{\mu\nu}\left(\frac{\pd_\mu A}{A}e_1+\frac{\pd_\mu B}{A}e_2\right)\left(\frac{\pd_\nu A}{A}e_1+\frac{\pd_\nu B}{A}e_2\right),
\end{align*}
where $e_1$ and $e_2$ are dual to $e^1$ and $e^2$ with respect to the metric $h$.

The equation for the wave map system (\ref{wk:X_eqn}-\ref{wk:Y_eqn}) linearized about the solution $(X_0,Y_0)=(A,B)$ is
$$\Box_g\xi_a - V_a{}^b\xi_b =0.$$
This fact is proven indirectly later in this chapter by first relating the $\xi_a$ bundle to the pair of functions $(\phi,\psi)$ in Lemma \ref{translate_nl_lem}, and then relating $(\phi,\psi)$ to the system (\ref{wk:X_eqn}-\ref{wk:Y_eqn}) in Proposition \ref{N_identities_prop}. See also \S\ref{three_systems_equivalent_sec}.

\subsection{The modified Kerr spacetime  $(\tilde{\mathcal{M}},\tilde{g})$}

Most of this chapter (from \S\ref{wk:phi_psi_estimates_sec} and on) will deal with the pair $(\phi,\psi)$ instead of the section $\xi_a$. The simplification in the beginning of \S\ref{ws:phi_psi_sec} serves as motivation for this pair. Unfortunately, the equations for $(\phi,\psi)$ are not as simple as in the Schwarzschild case, but more importantly both $\phi$ and $\psi$ are still regular on the axis. The equations for $(\phi,\psi)$ will be discussed soon, but first it is necessary to (re)introduce the modified spacetime $(\tilde{M},\tilde{g})$.

In the general Kerr $a\ne 0$ case, the operator
$$\Box_g+\frac{\pd^\alpha A}{A}\pd_\alpha$$
has fewer commutators, because the function
$$A=\frac{(r^2+a^2)^2-a^2\sin^2\theta(r^2-2Mr+a^2)}{r^2+a^2\cos^2\theta}\sin^2\theta$$
depends nontrivially on $r$ and $\theta$. For this reason, we write
\begin{align*}
A &= A_1A_2 \\
A_1 &= (r^2+a^2)\sin^2\theta \\
A_2 &= \left(1+\frac{a^2\sin^2\theta}{q^2}\right)\left(1-a^2\sin^2\theta v\right) \\
v &= \frac{r^2-2Mr+a^2}{(r^2+a^2)^2},
\end{align*}
and we generalize the spacetime metric $\tilde{g}$ from Chapter \ref{wm_szd_chap} by replacing $A$ with $A_1$. (Note that $A_1$ reduces to $A$ in the Schwarzschild case--generally speaking, $A_1$ behaves very similarly to $A$ in the Schwarzschild case, while terms depending on $A_2$ are treated as error terms in very much the same way that terms depending on $B$ are treated.)

That is, we extend the calculation for Schwarzschild by defining
$$\Box_{\tilde{g}} = \Box_g\psi +\frac{\pd^\alpha A_1}{A_1}\pd_\alpha$$
and
$$\int_{\tilde{\Sigma}_t}f = \int_{\Sigma_t}f A_1^2 =\int_{r_H}^\infty\int_0^\pi f A_1^2 q^2\sin\theta d\theta dr.$$

\subsection{Conventions for spacetime norms}

We repeat here the conventions used in Chapter \ref{wm_szd_chap}, but generalized to Kerr spacetimes. \\
\bp All functions will be assumed to depend only on $t$, $r$, and $\theta$. Therefore, any quantity can be treated as a function defined on either spacetime $(\mathcal{M},g)$ or $(\mathcal{\tilde{M}},\tilde{g})$. \\
\bp When necessary, the tilde mark ( $\tilde{ }$ ) will be used to denote quantities corresponding to $(\tilde{\mathcal{M}},\tilde{g})$. This includes the effective volume form 
$$\tilde{\mu}=A_1^2q^2\sin\theta=q^2(r^2+a^2)^2\sin^5\theta,$$
and the constant-time hypersurface 
$$\tilde{\Sigma}_t=\{t\}\times [r_H,\infty)\times S^6.$$
\bp Integrated expressions will depend on a volume form that is implicitly defined by the manifold of integration. That is,
$$\int_{\Sigma_t}f=\int_{r_H}^\infty\int_0^\pi\int_0^{2\pi}f(t,r,\theta)q^2\sin\theta d\phi d\theta dr.$$
and
\begin{multline*}
\int_{\tilde{\Sigma}_t}f=\int_{r_H}^\infty\int_0^\pi\int_0^\pi\int_0^\pi\int_0^\pi\int_0^\pi\int_0^{2\pi}f(t,r,\theta)q^2(r^2+a^2)^2\sin^5\theta\sin^4\theta_2\sin^3\theta_3\sin^2\theta_4\sin\theta_5 \\
d\phi d\theta_5 d\theta_4 d\theta_3 d\theta_2 d\theta dr.
\end{multline*}
Since all relevant integral estimates in this chapter are valid up to a constant, one may equivalently define
$$\int_{\Sigma_t}f=\int_{r_H}^\infty\int_0^\pi f(t,r,\theta)q^2\sin\theta d\theta dr$$
and
$$\int_{\tilde{\Sigma}_t}f=\int_{r_H}^\infty\int_0^\pi f(t,r,\theta)q^2(r^2+a^2)^2\sin^5\theta d\theta dr.$$
The point is that the only difference occurs in the factors that show up in the volume form. \\
\bp We also observe that $L^\infty$ estimates are weaker in the higher-dimensional spacetime. That is,
$$||f||_{L^\infty(S^2)}\lesssim \sum_{i\le 2}||\sla\nabla^if||_{L^2(S^2)},$$
while
$$||f||_{L^\infty(S^6)}\lesssim \sum_{i\le 4}||\tilde{\sla\nabla}^if||_{L^2(S^6)}.$$

\subsection{The $(\phi,\psi)$ system}\label{wk:phi_psi_equations_intro_sec}

The linearized wave map system is given by the equations
\begin{align*}
\Box_g\phi &= \mathcal{L}_\phi \\
\Box_{\tilde{g}}\psi &= \mathcal{L}_\psi,
\end{align*}
where
$$\mathcal{L}_\phi=-2\frac{\pd^\alpha B}{A}A\pd_\alpha \psi + 2\frac{\pd^\alpha B\pd_\alpha B}{A^2}\phi-4\frac{\pd^\alpha A\pd_\alpha B}{A^2} A\psi$$
$$\mathcal{L}_\psi=-2\frac{\pd^\alpha A_2}{A_2}\pd_\alpha\psi+2\frac{\pd^\alpha B\pd_\alpha B}{A^2}\psi + 2A^{-1}\frac{\pd^\alpha B}{A}\pd_\alpha\phi,$$
and again, the modified wave operator $\Box_{\tilde{g}}$ appearing in the equation for $\psi$ is defined by
$$\Box_{\tilde{g}}=\Box_g+2\frac{\pd^\alpha A_1}{A_1}\pd_\alpha.$$

The fully nonlinear system is given by
\begin{align}
\Box_g\phi &= \mathcal{L}_\phi+\mathcal{N}_\phi, \label{phi_nonlinear_system_eqn}\\
\Box_{\tilde{g}}\psi &= \mathcal{L}_\psi+\mathcal{N}_\psi, \label{psi_nonlinear_system_eqn}
\end{align}
where the nonlinear terms are
\begin{align*}
(1+\phi)\mathcal{N}_\phi &= \pd^\alpha \phi \pd_\alpha \phi - A\pd^\alpha \psi A\pd_\alpha \psi +2\frac{\pd^\alpha B}{A} \phi A\pd_\alpha \psi -4\frac{\pd^\alpha A}{A}A\psi A\pd_\alpha \psi \\
&\hspace{.5in}-\frac{\pd^\alpha B\pd_\alpha B}{A^2}\phi^2+4\frac{\pd^\alpha A\pd_\alpha B}{A^2}\phi A\psi -4\frac{\pd^\alpha A\pd_\alpha A}{A^2}(A\psi)^2,
\end{align*}
$$(1+\phi)\mathcal{N}_\psi = 2\pd^\alpha \phi \pd_\alpha\psi  +4\frac{\pd^\alpha A}{A}\psi\pd_\alpha \phi -2\frac{\pd^\alpha B}{A}A^{-1}\phi\pd_\alpha\phi.$$
\textbf{These terms are not the same as the terms $\mathcal{N}_\phi$ and $\mathcal{N}_\psi$ defined in Chapter \ref{wm_szd_chap}.} Unlike the terms $\mathcal{N}_\phi$ and $\mathcal{N}_\psi$ defined in Chapter \ref{wm_szd_chap}, these terms have factors that are singular on the axis even when $a=0$. In Chapter \ref{wm_szd_chap}, a special linearization $Y=X^2\psi$ was used to eliminate all terms with singular factors, while in this chapter a linearization corresponding to $Y=A^2\psi$ in Schwarzschild is used. In the general $a\ne 0$ case, it appears to be impossible to remove all terms with singular factors by choosing a special linearization. Instead, to solve the problem in this chapter, the techniques in Appendix \ref{regularity_sec} are required.

The above system will be further studied starting in \S\ref{wk:phi_psi_estimates_sec} and throughout the remainder of this chapter.

\subsection{The energy norms in terms of $\phi$ and $\psi$}

The starting point for the decay argument is a family of spacetime estimates roughly of the form
$$E_p(t_2)+\int_{t_1}^{t_2}B_p(t)dt \lesssim E_p(t_1)+\int_{t_1}^{t_2}N_p(t)dt,$$
for values of $p$ ranging almost from $0$ to $2$. (See Theorem \ref{p_thm}.) These estimates and their higher order analogues are developed in \S\ref{wk:dynamic_estimates_sec}.

For a typical scalar wave problem, the weighted energy $E_p(t)$ is given by
$$E_p(t)=\int_{\Sigma_t}r^p\left[(L\phi)^2+|\sla\nabla\phi|^2+r^{-2}\phi^2+r^{-2}(\lbar\phi)^2\right].$$
For the particular problem studied in this chapter, the weighted energy is actually given by
\begin{multline*}
E_p(t)= \int_{\Sigma_t}r^p\left[(L\phi)^2+|\sla\nabla\phi|^2+r^{-2}\phi^2+r^{-2}(\lbar\phi)^2\right] \\
+\int_{\tilde{\Sigma}_t}r^p\left[(L\psi)^2+|\tsla\nabla\psi|^2+r^{-2}\psi^2+r^{-2}(\lbar\psi)^2\right].
\end{multline*}
Note that the second integral is over $\tilde{\Sigma}_t$, implying that there is an additional weight of $r^4\sin^4\theta$. This will generally be the case for integrals of quantities depending only on $\psi$, and has to do with the fact that the equation for $\psi$ is a wave equation naturally belonging to the modified spacetime $(\tilde{\mathcal{M}},\tilde{g})$.

The bulk quantity $B_p(t)$ is similar in weight to the weighted energy $E_{p-1}(t)$, except that it also has a degeneracy on the photon sphere.
\begin{multline*}
B_p(t)= \int_{\Sigma_t}r^{p-1}\left[\chi_{trap}(L\phi)^2+\chi_{trap}|\sla\nabla\phi|^2+r^{-2}\phi^2+r^{-2}(\pd_r\phi)^2\right] \\
+\int_{\tilde{\Sigma}_t}r^{p-1}\left[\chi_{trap}(L\psi)^2+\chi_{trap}|\tsla\nabla\psi|^2+r^{-2}\psi^2+r^{-2}(\pd_r\psi)^2\right].
\end{multline*}
The function $\chi_{trap}$ vanishes to second order at the trapping radius $r_{trap}$, which defines the photon sphere. $r_{trap}$ is the radius that maximizes the geodesic potential 
$$v=\frac{r^2+2Mr+a^2}{(r^2+a^2)^2}$$
and corresponds to the radius at which null geodesics with zero angular momentum orbit the black hole at a constant radius. It coincides with $3M$ in the Schwarzschild case.

There are $L^\infty$ estimates that pair with the energy norms in the dynamic estimates. They are proved in \S\ref{wk:pointwise_sec}, but a brief summary is given here.

 Let $\phi^s$ and $\psi^s$ denote quantities obtained by applying up to $s$ commutators to either $\phi$ or $\psi$, and let $E_p^s(t)$ be the generalized $s$-order energy corresponding to $\phi^s$ and $\psi^s$. One particular estimate roughly states
$$|rL\phi^s|+|r\sla\nabla\phi^s|+|\phi^s|+|\lbar\phi^s|+|r^3L\phi^s|+|r^3\sla\nabla\phi^s|+|r^2\phi^s|+|r^2\lbar\phi^s| \lesssim (E_0^{s+5}(t))^{1/2}.$$
A more general estimate is roughly the following.
\begin{multline*}|r^{p+1}L\phi^s|+|r^{p+1}\sla\nabla\phi^s|+|r^p\phi^s|+|r^p\lbar\phi^s|+|r^{p+3}L\phi^s|+|r^{p+3}\sla\nabla\phi^s|+|r^{p+2}\phi^s|+|r^{p+2}\lbar\phi^s| \\
\lesssim (E_{2p}(t))^{1/2}.
\end{multline*}
Note the additional factors of $r$ accompanying the $\psi$ terms. These factors are due to the additional factors present in the volume form for $\tilde{\Sigma}_t$.

\subsection{Regularity on the axis and a new structural condition}\label{intro_regularity_sec}

One of the most significant of the new complications of the problem of this chapter is the presence of apparently singular terms, ie. terms which have factors of $\sec\theta$ or $\csc\theta$. It is not a priori clear whether these challenges are related to the true dynamics of the problem or whether they are due to the fact that the coordinate system is degenerate on the axis. As we shall see, in the $(\phi,\psi)$ picture these terms are merely due to the degeneracy of the coordinate system. To handle these terms, we use a new formalism which is developed in detail in Appendix \ref{regularity_sec}, although a brief summary is given here.

Consider the two functions $\cos\theta$ and $\sin\theta$. While these are both smooth functions of $\theta$, if they are treated as functions on the spacetime, one of them is actually much less regular than the other. The problem is that the function $\sin\theta$ behaves like $|\theta|$ in a neighborhood of the half axis $\theta=0$, which means that $\sin\theta$ is not twice differentiable. In contrast, the function $\cos\theta$ is everywhere smooth. Informally, we will say that $\cos\theta$ belongs to a space of functions that are regular on the axis, but $\sin\theta$ does not.

The fact that $\sin\theta$ is not regular on the axis is not clear when applying two coordinate derivatives, because $\pd_\theta^2\sin\theta=-\sin\theta$ appears to be bounded. In order to measure the singular nature of $\sin\theta$ on the axis, it is necessary to use a second operator $\cot\theta\pd_\theta$, which by no coincidence appears in the spherical laplacian. Since $\cot\theta\pd_\theta\sin\theta=\cos^2\theta\csc\theta$, it is now clear that something goes wrong on the axis. Since the two operators $\pd_\theta^2$ and $\cot\theta\pd_\theta$ will often be used, they are given names.
\begin{align*}
\fa &= \pd_\theta^2 \\
\fb &= \cot\theta\pd_\theta.
\end{align*}

The operator $\fb$ itself seems to be in some sense singular on the axis, because it has a factor of $\cot\theta$. However, for any twice-differentiable axisymmetric function $f$, the first derivative $\pd_\theta f$ should vanish at least to first order on the axis. So for twice-differentiable axisymmetric functions, there is a cancellation effect with the factor $\cot\theta$. Essentially, \textit{the operators $\fa$ and $\fb$ preserve the space of regular functions on the axis.}

In the nonlinear problem, we will commute the fully nonlinear equations (\ref{phi_nonlinear_system_eqn}-\ref{psi_nonlinear_system_eqn}) with the following angular operators.
\begin{align*}
Q &= \fa+\fb+a^2\sin^2\theta\pd_t^2 \\
\tilde{Q} &= \fa+5\fb+a^2\sin^2\theta\pd_t^2.
\end{align*}
Let us ignore for a moment the $a^2\sin^2\theta\pd_t^2$ part of these operators. We will need to estimate the terms belonging to
$$(\fa+\fb)^l\mathcal{N}_\phi\hspace{.5in}\text{and}\hspace{.5in}(\fa+5\fb)^l\mathcal{N}_\psi.$$
We introduce the operator family $\fc^l$ to represent any term in the expansion of $(\fa+\fb)^l$. So for example,
$$\fc^2=\{\fa^2,\fa\fb,\fb\fa,\fb^2\}.$$
We will estimate $\fc^l(\mathcal{N}_\phi)$ and $\fc^l(\mathcal{N}_\psi)$, but must do so carefully to ensure that we stay in the space of regular functions on the axis.

There are a few ways things could go wrong if we are not careful. First, if at any point we expand an operator such as $\fa\fb$ or $\fb^2$, then we get terms that are truly singular on the axis. As an example, consider the following calculation.
$$\fb^2f =\cot\theta\pd_\theta(\cot\theta\pd_\theta f) =\cot^2\theta\pd_\theta^2f-\cot\theta\csc^2\theta\pd_\theta f.$$
If $f$ is regular on the axis, then $\pd_\theta f$ will vanish at least to first order on the axis, but that is not enough to ensure that each of the two terms on the right side remain regular--in general they do not. For this reason, it is necessary to treat the operator $\fb$ as an atomic operator.

The second thing that can go wrong happens when applying $\fa$ or $\fb$ to products of functions. As an example, consider the following calculation.
$$\fa(fg)=\pd_\theta^2(fg)=\fa f g+2\pd_\theta f \pd_\theta g+f\fa g.$$
If $f$ and $g$ are regular on the axis, then so are $\fa f$ and $\fa g$, but the factors $\pd_\theta f$ and $\pd_\theta g$ are not--for the same reason that $\cos\theta$ is regular on the axis, but $-\sin\theta=\pd_\theta\cos\theta$ is not. Here, it is important to observe that \textit{the product $\pd_\theta f\pd_\theta g$ is indeed regular on the axis}, because it is a product of two functions behaving like $|\theta|$, which will behave like $\theta^2$. Without being careful, it is possible to expand $\fb^2\fa(fg)$ to get a term of the form $\fb^2(\pd_\theta f)\pd_\theta g$. This term would not be regular on the axis. The correct way to handle such an expansion is illustrated by the following intermediate calculation.
$$\fb(\pd_\theta f\pd_\theta g)=\cot\theta \pd_\theta^2f\pd_\theta g+\cot\theta\pd_\theta f\pd_\theta^2 g = \fa f \fb g+\fb f\fa g.$$
Then one can apply any of the additional $\fc$ operators.

To ensure that products are regular on the axis (especially after applying $\fa$ or $\fb$), the most general form of products that can be permitted in the nonlinear term must be something like 
$$\pd_\theta \fc^{l_1}f_1\pd_\theta \fc^{l_2}f_2...\pd_\theta \fc^{l_{2k}}f_{2k}\fc^{l_{2k+1}}f_{2k+1}...\fc^{l_{2k+k'}}f_{2k+k'}.$$
That is, \textit{it is essential that an even number of single $\pd_\theta$ derivatives show up and that each $\pd_\theta$ be applied only after the $\fc$ operators are applied.} To write this more compactly, we define yet another family of operators.
$$
\fd^l = \left\{\begin{array}{cc}
\fc^{l/2} & l\in 2\mathbb{Z} \\
\pd_\theta \fc^{(l-1)/2} & l\not\in 2\mathbb{Z}.
\end{array}\right.
$$
If $f_1,...,f_k$ are each regular functions, then regular product terms will be of the form
$$\fd^{i_1}f_1...\fd^{i_k}f_k$$
where $i_1+...+i_k=2n$. We call these products \textit{terms of degree $n$}. An important fact is that if $\fa$ or $\fb$ is applied to a term of degree $n$, the result can be expressed as a sum of terms of degree $n+1$. \textbf{A new important structural condition for the nonlinear terms arising in the wave map problem is that they be of this form.}

\section{Estimates for the $\xi_a$ system}\label{wk:xi_estimates_sec}
The geometric nature of the wave map problem suggests a vector bundle formalism with which one can describe a linearization of the wave map system (\ref{wk:X_eqn}-\ref{wk:Y_eqn}). For more information about this bundle, see \S\ref{wave_map_bundle_sec}. In this section, the vector bundle formalism is used to derive the most delicate spacetime estimates, including the $h\pd_t$ estimate and the Morawetz estimate.

\subsection{A current template for the bundle $\mathcal{B}$}

To start, we generalize the current template from Lemma \ref{divJ_lem} to the bundle $\mathcal{B}$.

\begin{lemma}\label{vectorized_current_template_lem}
For a vectorfield $X$, a function $w$, and a $\mathcal{B}'\otimes\mathcal{B}'$-valued one-form $m_\mu^{ab}$, define
$$J[X,w,m]_\mu=T_{\mu\nu}X^\nu+w\xi\cdot D_\mu\xi-\frac12|\xi|^2\pd_\mu w+m_\mu^{ab}\xi_a\xi_b,$$
where
$$T_{\mu\nu}=2D_\mu\xi^a D_\nu\xi_a-g_{\mu\nu}D^\lambda\xi^a D_\lambda\xi_a-g_{\mu\nu}V^{ab}\xi_a\xi_b.$$
Then
\begin{multline*}
div J = K^{\mu\nu}D_\mu\xi\cdot D_\nu\xi-\frac12\Box_g w|\xi|^2+((w-divX)V^{ab}-D_XV^{ab})\xi_a\xi_b \\
+D^\mu m_\mu^{ab}\xi_a\xi_b+2m_\mu^{ab}\xi_aD^\mu\xi_b  +2R_{\mu\nu a b}X^\mu \xi^a D^\nu \xi^b \\
+(\Box_g\xi_a-V_a{}^b\xi_b)(2D_X\xi^a+w\xi^a).
\end{multline*}
where
$$K^{\mu\nu}=2\nabla^{(\mu}X^{\nu)}+(w-divX)g^{\mu\nu}$$
is the same tensor as defined in Definition \ref{em_tensor_def}.
\end{lemma}
\begin{remark}
The formula for $divJ$ in the lemma seems rather complicated, but it can be broken down into the following parts. (i) The part $K^{\mu\nu}D_\mu\xi\cdot D_\nu\xi-\frac12\Box_gw|\xi|^2$ is directly analogous to the scalar wave equation with no potential. (ii) The part $((w-divX)V^{ab}-D_XV^{ab})\xi_a\xi_b$ is new, because it depends on the potential $V^{ab}$ that is introduced in \S\ref{wk:xi_a_intro_sec}. (If the scalar equation also had a potential, this part would be analogous to an additional part for that scalar equation.) (iii) The part $D^\mu m_\mu^{ab}\xi_a\xi_b+2m_\mu^{ab}\xi_aD^\mu \xi_b$ is also new, but only because it will be helpful, since it could be excluded by choosing $m_\mu^{ab}=0$. (iv) The part $2R_{\mu\nu a b}X^\mu \xi^aD^\nu \xi^b$ is also new, but this time it is purely due to the use of a bundle instead of a scalar. It will be considered an error term because the curvature tensor vanishes in the Schwarzschild case (see Lemma \ref{bundle_R_calculation_lem} below). (v) Finally, the part $(\Box_g\xi_a-V_a{}^b\xi_b)(2D_X\xi^a+w\xi^a)$ is a nonlinear error term, because it vanishes when the linear equation is satisfied.
\end{remark}
\begin{proof}
To start, we compute the divergence of the energy-momentum tensor.
\begin{align*}
\nabla^\nu T_{\mu\nu} &= 2D_\mu\xi^a D^\nu D_\nu\xi_a + 2D^\nu D_\mu\xi^a D_\nu\xi_a - 2D_\mu D^\lambda\xi^a D_\lambda\xi_a -D_\mu V^{ab}\xi_a\xi_b -2V^{ab}\xi_aD_\mu\xi_b \\
&= 2D_\mu\xi^a\Box_g\xi_a+2(D_\nu D_\mu\xi_a-D_\mu D_\nu\xi_a)D^\nu\xi^a - D_\mu V^{ab}\xi_a\xi_b -2V_a{}^b\xi_bD_\mu \xi^a \\
&= 2(\Box_g\xi_a-V_a{}^b\xi_b)D_\mu\xi^a+2R_{\nu\mu ab}\xi^b D^\nu\xi^a -D_\mu V^{ab}\xi_a\xi_b \\
&= 2(\Box_g\xi_a-V_a{}^b\xi_b)D_\mu\xi^a+2R_{\mu\nu ab}\xi^a D^\nu\xi^b -D_\mu V^{ab}\xi_a\xi_b.
\end{align*}
It follows that
\begin{align*}
\nabla^\nu(T_{\mu\nu}X^\mu) &= \nabla^\nu T_{\mu\nu}X^\mu + T_{\mu\nu}\nabla^{(\mu}X^{\nu)} \\
&= 2(\Box_g\xi_a-V_a{}^b\xi_b)D_X\xi^a+2R_{\mu\nu ab}X^\mu \xi^aD^\nu\xi^b-D_XV^{ab}\xi_a\xi_b \\
&\hspace{.2in}+(2\nabla^{(\mu}X^{\nu)}-divX g^{\mu\nu})D_\mu\xi^aD_\nu\xi_a-(divX) V^{ab}\xi_a\xi_b.
\end{align*}
Also,
\begin{align*}
\nabla^\mu &(w\xi^a D_\mu\xi_a-\frac12\xi^a\xi_a\pd_\mu w) \\
&= w\xi^a D^\mu D_\mu\xi_a+wD^\mu\xi^a D_\mu\xi_a +\pd^\mu w\xi^a D_\mu\xi_a -\pd_\mu w\xi^aD^\mu\xi_a-\frac12\xi^a\xi_a\nabla^\mu \pd_\mu w \\
&= w\xi^a \Box_g\xi_a+wD^\mu\xi^a D_\mu\xi_a -\frac12\Box_g w \xi^a\xi_a \\
&= (\Box_g\xi_a-V_a{}^b\xi_b)w\xi^a + wV^{ab}\xi_a\xi_b+wD^\mu\xi^a D_\mu\xi_a -\frac12\Box_g w \xi^a\xi_a.
\end{align*}
And
$$\nabla^\mu(m_\mu^{ab}\xi_a\xi_b) = D^\mu m_\mu^{ab}\xi_a\xi_b + 2m_\mu^{ab}\xi_aD^\mu\xi_b.$$
Summing all these terms proves the statement of the lemma.
\end{proof}

The next lemma condenses the part of $divJ$ having to do with the potential.
\begin{lemma}\label{muvVab_lem}
If $X=uv\pd_r$ and $w=v\pd_ru$, then
$$(w-divX)V^{ab}-D_XV^{ab}=-\frac{u}{\mu}D_r(\mu v V^{ab}).$$
\end{lemma}
\begin{proof}
\begin{align*}
(w-divX)V^{ab}-X^\mu D_\mu V^{ab} &= (v\pd_ru - \frac1\mu\pd_r(\mu u v))V^{ab}-uvD_rV^{ab} \\
&= -\frac{u}{\mu}\pd_r(\mu v)V^{ab}-uv D_rV^{ab} \\
&= -\frac{u}{\mu}D_r(\mu v V^{ab}).
\end{align*}
\end{proof}

The final lemma gives a formula for the curvature tensor $R_{\mu\nu ab}$ for the bundle $\mathcal{B}$. Note in particular that in the Schwarzschild case, since $B=0$, the curvature tensor vanishes.
\begin{lemma}\label{bundle_R_calculation_lem}
The curvature tensor is given by
$$R_{\mu\nu a b} = -2\frac{\pd_{[\mu} A\pd_{\nu]} B}{A^2}\epsilon_{ab},$$
where $\epsilon_{ab}$ is the volume element for the bundle $\mathcal{B}$.
\end{lemma}
\begin{proof}
We calculate $R_{\mu\nu a b}$ according to the formula
$$D_\mu D_\nu (e^i)_a -D_\nu D_\mu (e^i)_a = R_{\mu\nu a b}(e^i)^b.$$
We start with
\begin{align*}
D_\mu D_\nu e^1 &= D_\mu\left(-\frac{\pd_\nu B}{A} e^2\right) \\
&= -\frac{\nabla_\mu \pd_\nu B}{A}e^2+\frac{\pd_\mu A\pd_\nu B}{A^2}e^2-\frac{\pd_\mu B\pd_\nu B}{A^2}e^1.
\end{align*}
Therefore,
$$D_\mu D_\nu e^1 -D_\nu D_\mu e^1 = 2\frac{\pd_{[\mu }A\pd_{\nu]}B}{A^2}e^2.$$
Likewise, we calculate
\begin{align*}
D_\mu D_\nu e^2 &= D_\mu\left(\frac{\pd_\nu B}{A} e^1\right) \\
&= \frac{\nabla_\mu \pd_\nu B}{A}e^1-\frac{\pd_\mu A\pd_\nu B}{A^2}e^1-\frac{\pd_\mu B\pd_\nu B}{A^2}e^2.
\end{align*}
Therefore,
$$D_\mu D_\nu e^2 -D_\nu D_\mu e^2 = -2\frac{\pd_{[\mu }A\pd_{\nu]}B}{A^2}e^1.$$
The formula is verified by the following two relations.
$$(D_\mu D_\nu e^1 -D_\nu D_\mu e^1)_a = 2\frac{\pd_{[\mu}A\pd_{\nu]}B}{A^2}(e^2)_a = -2\frac{\pd_{[\mu}A\pd_{\nu]}B}{A^2}\epsilon_{ab}(e^1)^b,$$
$$(D_\mu D_\nu e^2 -D_\nu D_\mu e^2)_a = -2\frac{\pd_{[\mu}A\pd_{\nu]}B}{A^2}(e^1)_a = -2\frac{\pd_{[\mu}A\pd_{\nu]}B}{A^2}\epsilon_{ab}(e^2)^b.$$
This completes the proof.
\end{proof}

\subsection{The energy estimate for the $\xi_a$ system}

The most fundamental of all estimates is the energy estimate. We begin by proving this estimate. The slightly more general $h\pd_t$ estimate will be proved in \S\ref{wk:xi_hdt_sec}.

One of the key advantages to using the vector bundle formalism is that it reveals a conserved quantity based on the multiplier $\pd_t$.

\begin{lemma}\label{wk:energy_identity_lem}
(Energy identity for the $\xi_a$ system)
$$\int_{H_{t_1}^{t_2}}J^r[\pd_t]+\int_{\Sigma_{t_2}}-J^t[\pd_t]=\int_{\Sigma_{t_1}}-J^t[\pd_t]+\int_{t_1}^{t_2}\int_{\Sigma_t}-2D_t\xi^a(\Box_g\xi_a-V_a{}^b\xi_b)$$
In particular, if $\Box_g\xi_a-V_a{}^b\xi_b=0$, then the quantity
$$E(t)=\int_{H_{t_0}^t}J^r[\pd_t]+\int_{\Sigma_t}-J^t[\pd_t]$$
is conserved.
\end{lemma}
\begin{proof}
By Lemma \ref{vectorized_current_template_lem},
\begin{multline*}
div J[X] = K_X^{\mu\nu}D_\mu\xi\cdot D_\nu\xi+(-divX V^{ab}-D_XV^{ab})\xi_a\xi_b 
+2R_{\mu\nu a b}X^\mu \xi^a D^\nu \xi^b \\
+(\Box_g\xi_a-V_a{}^b\xi_b)(2D_X\xi^a).
\end{multline*}
For the particular case $X=\pd_t$, it is already known from the scalar problem that $K_X^{\mu\nu}=0$ and $div X=0$. Next, since $V^{ab}$ does not depend on $t$ and since $D_{\pd_t}e_i=0$, it follows that $D_XV^{ab}=0$. Also, by Lemma \ref{bundle_R_calculation_lem} and the fact that $A$ and $B$ do not depend on $t$, $R_{\mu\nu a b}X^\mu=0$. Therefore,
$$div J[\pd_t] = 2 D_t\xi^a (\Box_g\xi_a-V_a{}^b\xi_b).$$
We conclude that
$$\int_{H_{t_1}^{t_2}}J^r[\pd_t]+\int_{\Sigma_{t_2}}-J^t[\pd_t]+\int_{t_1}^{t_2}\int_{\Sigma_t}2D_t\xi^a(\Box_g\xi_a-V_a{}^b\xi_b)=\int_{\Sigma_{t_1}}-J^t[\pd_t].$$
The statement of the lemma now follows.
\end{proof}

Now we calculate the flux terms in the previous lemma.

\begin{lemma}\label{h_dt_J_components_lem}
On the event horizon $H_{t_1}^{t_2}$,
$$J^r[\pd_t] \approx |D_t\xi|^2$$
and on a constant-time hypersurface $\Sigma_t$,
$$-J^t[\pd_t] \approx \chi_H|D_r\xi|^2+|D_t\xi|^2+q^{-2}|D_\theta\xi|^2+V^{ab}\xi_a\xi_b,$$
where $\chi_H=1-\frac{r_H}r$.
\end{lemma}
\begin{proof}
From Lemma \ref{vectorized_current_template_lem},
$$J^\mu[\pd_t] = 2g^{\mu\lambda}D_\lambda\xi\cdot D_t\xi -\delta^\mu{}_t D^\lambda\xi\cdot D_\lambda\xi-\delta^\mu{}_tV^{ab}\xi_a\xi_b.$$
On the event horizon $H_{t_1}^{t_2}$, since $g^{rr}$ vanishes,
$$J^r[\pd_t] = 2g^{rt}|D_t\xi|^2.$$
Since $g^{rt}$ is constant and positive on the event horizon,
$$J^r[\pd_t] \approx |D_t\xi|^2.$$
Also, on a constant-time hypersurface $\Sigma_t$,
\begin{multline*}
-J^t[\pd_t] = -2g^{tt}|D_t\xi|^2-2g^{tr}D_r\xi\cdot D_t\xi \\
+\left(g^{tt}|D_t\xi|^2+2g^{tr}D_r\xi\cdot D_t\xi+g^{rr}|D_r\xi|^2+q^{-2}|D_\theta\xi|^2+\frac{a^2\sin^2\theta}{q^2}|D_t\xi|^2\right) \\
+V^{ab}\xi_a\xi_b.
\end{multline*}
Thus,
$$-J^t[\pd_t] = \left(-g^{tt}+\frac{a^2\sin^2\theta}{q^2}\right)|D_t\xi|^2+g^{rr}|D_r\xi|^2+q^{-2}|D_\theta\xi|^2+V^{ab}\xi_a\xi_b.$$
The quantity $-g^{tt}+\frac{a^2\sin^2\theta}{q^2}$ is uniformly bounded below by a positive constant for all $|a|<M$.
\end{proof}

From the previous two lemmas, we conclude the following energy estimate.
\begin{proposition}\label{xi_energy_estimate_prop} (Energy estimate for the $\xi_a$ system)
\begin{multline*}
\int_{H_{t_1}^{t_2}}|D_t\xi|^2+\int_{\Sigma_{t_2}}\left[\chi_H|D_r\xi|^2+|D_t\xi|^2+q^{-2}|D_\theta\xi|^2+V^{ab}\xi_a\xi_b\right] \\
\lesssim \int_{\Sigma_{t_1}}\left[\chi_H|D_r\xi|^2+|D_t\xi|^2+q^{-2}|D_\theta\xi|^2+V^{ab}\xi_a\xi_b\right] + Err_{nl},
\end{multline*}
where $\chi_H=1-\frac{r_H}{r}$ and
$$Err_{nl} = \int_{t_1}^{t_2}\int_{\Sigma_t} |D_t\xi^a(\Box_g\xi_a-V_a{}^b\xi_b)|.$$
\end{proposition}
\begin{proof}
This follows directly from Lemmas \ref{wk:energy_identity_lem} and \ref{h_dt_J_components_lem}.
\end{proof}

\subsection{The $h\pd_t$ estimate for the $\xi_a$ system}\label{wk:xi_hdt_sec}

\begin{lemma}\label{h_dt_divJ_lem}
If $h=h(r)$ is constant in the interval $r\in[r_H,r_H+\delh]$ and $\Box_g\xi_a-V_a{}^b\xi_b=0$, then
$$\frac{q^2}{r^2+a^2}divJ[h\pd_t]= \frac{h'}{2}\left[|D_L\xi|^2-|D_{\lbar}\xi|^2\right].$$
\end{lemma}
\begin{proof}
Recall from Lemma \ref{vectorized_current_template_lem} that
$$divJ[X] = K^{\mu\nu}D_\mu\xi\cdot D_\nu\xi - divXV^{ab}-D_XV^{ab}+2R_{\mu\nu a b}X^\mu\xi^aD^\nu\xi^b.$$
If $X=h\pd_t$, then since
$$div(h\pd_t)=0$$
and
$$D_{\pd_t}V^{ab}=0$$
and
$$2R_{\mu\nu a b}(\pd_t)^\mu=0,$$
the only nonzero term is $K^{\mu\nu}D_\mu\xi\cdot D_\nu\xi$. 

Recall also from Lemma \ref{vectorized_current_template_lem} that
\begin{align*}
K^{\mu\nu}&=2\nabla^{(\mu}X^{\nu)}+(w-divX)g^{\mu\nu}.
\end{align*}
Since $divX=0$ and $w=0$, we have
$$K^{\mu\nu}=2\nabla^{(\mu}X^{\nu)}= g^{\mu\lambda}\pd_\lambda X^\nu+g^{\nu\lambda}\pd_\lambda X^\mu - X^\lambda\pd_\lambda(g^{\mu\nu}).$$
Since $\pd_t$ is killing, $\pd_t(g^{\mu\nu})=0$. Also, the only component of $X$ is the $t$ component and the only nonzero derivative of that component is the $\pd_r$ derivative, so
\begin{align*}
g^{\mu\lambda}\pd_\lambda X^\nu +g^{\nu\lambda}\pd_\lambda X^\mu &=
g^{\mu r}\pd_r X^\nu +g^{\nu r}\pd_r X^\mu.
\end{align*}
It follows that the only possible nonzero $K^{\mu\nu}$ components are
\begin{align*}
K^{tt} &= 2g^{tr}h' \\
K^{tr}+K^{rt} &= 2g^{rr}h'.
\end{align*}
Since $h'=0$ in the region $r\in[r_H,r_H+\delh]$, the first of these actually vanishes. We conclude that
\begin{align*}
\frac{q^2}{r^2+a^2}divJ[h\pd_t]&=2\frac{q^2}{r^2+a^2}g^{rr}h'D_r\xi\cdot D_t\xi \\
&= 2\alpha h'D_r\xi\cdot D_t\xi \\
&= \frac{h'}2 \left[|\alpha D_r\xi+D_t\xi|^2-|\alpha D_r\xi-D_t\xi|^2\right] \\
&= \frac{h'}2 \left[|D_L\xi|^2-|D_{\lbar}\xi|^2\right]
\end{align*}
\end{proof}

Taking $h$ to be a positive function decreasing to zero as $r\rightarrow\infty$ at a particular rate, we obtain the $h\pd_t$ estimate.
\begin{proposition}\label{xi_h_dt_prop}($h\pd_t$ estimate for the $\xi_a$ system)
Fix $\delp>0$ and let $p\le 2-\delp$. Let $R>r_H+\delh$ be any given radius. Then for all $\epsilon >0$, there is a small constant $c_\epsilon$ and a large constant $C_\epsilon$, such that
\begin{multline*}
\int_{H_{t_1}^{t_2}}|D_t\xi|^2+\int_{\Sigma_{t_2}}r^{p-2}\left[\chi_H|D_r\xi|^2+|D_t\xi|^2+q^{-2}|D_\theta\xi|^2 +V^{ab}\xi_a\xi_b\right] \\
+\int_{t_1}^{t_2}\int_{\Sigma_t\cap\{R+M<r\}}c_\epsilon r^{p-3}|D_{\lbar}\xi|^2 \\
\lesssim \int_{\Sigma_{t_1}}C_\epsilon r^{p-2}\left[\chi_H|D_r\xi|^2+|D_t\xi|^2+q^{-2}|D_\theta\xi|^2+V^{ab}\xi_a\xi_b\right] + Err,
\end{multline*}
where $\chi_H=1-\frac{r_H}{r}$ and
\begin{align*}
Err&=Err_1+Err_{nl} \\
Err_1 &= \int_{t_1}^{t_2}\int_{\Sigma_t\cap\{R<r\}}\epsilon r^{-1}|D_L\xi|^2 \\
Err_{nl} &= \int_{t_1}^{t_2}\int_{\Sigma_t}C_\epsilon r^{p-2}|D_t\xi^a(\Box_g\xi_a-V_a{}^b\xi_b)|.
\end{align*}
\end{proposition}
\begin{proof}
See the proof of Proposition \ref{m:hdt_prop}, noting Lemmas \ref{h_dt_J_components_lem} and \ref{h_dt_divJ_lem}.
\end{proof}

\subsection{A result from the scalar wave equation}

We restate a result from \S\ref{k:morawetz_sec}.

\begin{lemma}\label{Kmunu_lem}
Let $X$ be a vectorfield and $w$ a scalar function. Define
$$K^{\mu\nu}=2\nabla^{(\mu}X^{\nu)}+(w-divX)g^{\mu\nu}.$$
For all $\epsilon_{temper}>0$ sufficiently small, there exists a function $u(r)$ satisfying the following conditions.

If $X=uv\pd_r$ and $w=v\pd_r u$, where
$$v=\frac{\Delta}{(r^2+a^2)^2}$$
is the geodesic potential for Kerr, then for any axisymmetric function $\psi$,
$$K^{\mu\nu}\pd_\mu\psi\pd_\nu\psi = 2\left(u'-\frac{2ru}{r^2+a^2}\right)\frac{\Delta^2}{q^2(r^2+a^2)}(\pd_r\psi)^2-\frac{u\pd_r v}{q^2}Q^{\mu\nu}\pd_\mu\psi\pd_\nu\psi.$$
Furthermore, we define
\begin{align*}
K^{tt}&=0 \\
K^{rr}&=2\left(u'-\frac{2ru}{r^2+a^2}\right)\frac{\Delta^2}{q^2(r^2+a^2)} \\
K_{Q}^{\mu\nu} &= -u\pd_r v \frac{Q^{\mu\nu}}{q^2}.
\end{align*}
Then
\begin{align*}
K^{rr} &\approx \frac{M^2}{r^3}\left(1-\frac{r_H}r\right)^2 \\
K_{Q}^{\theta\theta} &\approx \frac1{r^3}\left(1-\frac{r_{trap}}r\right) \\
K_{Q}^{tt} &\approx \frac{a^2\sin^2\theta}{r^3}\left(1-\frac{r_{trap}}r\right).
\end{align*}
Also,
$$\frac{M}{r^4}1_{r\ge r_*}-q^{-2}V_{\epsilon_{temper}} \le -\frac12\Box_gw,$$
where $V_{\epsilon_{temper}}$ is a positive function supported near the event horizon and satisfying 
$$||V_{\epsilon_{temper}}||_{L^1(r)} < \epsilon_{temper}.$$

Let $r_{trap}$ be the radius where $\pd_rv$ changes sign from positive to negative, and let $r_*$ be the radius where $\pd_r(2rv)$ changes sign from positive to negative. (In Schwarzschild, $r_{trap}=3M$ and $r_*=4M$.) Then $u$ changes sign from negative to positive at $r_{trap}$ and $\pd_ru=2r$ for $r>r_*$.
\end{lemma}
\begin{proof}
This was established in \S\ref{k:morawetz_sec}.
\end{proof}

\subsection{The partial Morawetz estimate for the $\xi_a$ system in the Schwarzschild case}\label{partial_morawetz_szd_sec}

Here, we prove the partial Morawetz estimate for the Schwarzschild case in a way that it can be adapted to Kerr for $|a|\ll M$. Since the proof of the estimate in Kerr is complicated, it can help to understand the simpler proof for the Schwarzschild case first.

\begin{proposition}\label{partial_morawetz_szd_prop} (Partial Morawetz estimate for the $\xi_a$ system in the Schwarzschild case) Suppose $a=0$, and suppose $\xi_a$ satisfies the equation $\Box_g\xi_a - V_a{}^b\xi_b=0$. Then there exists a current $J$ such that
\begin{multline*}
\frac{M^2}{r^3}\left(1-\frac{2M}r\right)^2(\pd_r\xi_1)^2+\frac{\chi_{trap}}{r}|\sla\nabla\xi_1|^2+\frac{1}{r^3}(\xi_1)^2 + \chi_{trap}\frac{\cot^2\theta}{r^3}(\xi_1)^2 \\
+ \frac{M^2}{r^3}\left(1-\frac{2M}r\right)^2(\pd_r\xi_2)^2+\frac{\chi_{trap}}{r}|\sla\nabla\xi_2|^2+\frac{M}{r^4}1_{r\ge 4M}(\xi_2)^2 \\
-r^{-2}V_{\epsilon_{temper}}((\xi_1)^2+(\xi_2)^2) \\
\lesssim div J,
\end{multline*}
where $\chi_{trap}=\left(1-\frac{3M}r\right)^2$, and $V_{\epsilon_{temper}}$ is the potential defined in Lemma \ref{Kmunu_lem}.
\end{proposition}

\begin{proof}
We choose
$$J_\mu=J[X,w,m]_\mu,$$
where $J[X,w,m]$ is the current template defined in Lemma \ref{vectorized_current_template_lem}.

We use the same vectorfield
$$X=uv\pd_r$$
and scalar function
$$w=v\pd_ru$$
that are used for the scalar wave equation (see Lemma \ref{Kmunu_lem}).

However, we now choose an additional one-form $m_\mu^{ab}$ with components
\begin{align*}
m_\mu^{ab}&=m_\mu^{ij} (e_i)^a(e_j)^b \\
m^{11}_\mu dx^\mu&=\frac{4\chi v}{r^2}dr+(2-\epsilon)u\pd_rv\cot\theta d\theta \\
m^{12}_\mu dx^\mu&=m^{21}=m^{22}=0,
\end{align*}
where the function $\chi$ will be defined in Lemma \ref{0_K_r_lem}. It is not common to include an $m$ term  with a nonzero $d\theta$ component. The reason will become clear in the proof of Lemma \ref{0_K_theta_lem}. (See the remark following the proof of Lemma \ref{0_K_theta_lem}.)

By Lemma \ref{vectorized_current_template_lem}, we have that
\begin{multline*}
div J = K^{\mu\nu}D_\mu\xi\cdot D_\nu\xi-\frac12\Box_g w|\xi|^2+((w-divX)V^{ab}-D_XV^{ab})\xi_a\xi_b \\
+D^\mu m_\mu^{ab}\xi_a\xi_b+2m_\mu^{ab}\xi_aD^\mu\xi_b  +2R_{\mu\nu a b}X^\mu \xi^a D^\nu \xi^b \\
+(\Box_g\xi_a-V_a{}^b\xi_b)(2D_X\xi^a+w\xi^a).
\end{multline*}
We rearrange these terms according to the following lemma.
\begin{lemma}\label{0_divJ_rearranged_lem}
$$divJ =\mathcal{K}-\frac12\Box_g w |\xi|^2+(\Box_g\xi_a-V_a{}^b\xi_b)(2D_X\xi^a+w\xi^a),$$
where
$$\mathcal{K}= \mathcal{K}_{(r)}+\mathcal{K}_{(\theta)}+\mathcal{K}_{(2)}$$
$$\mathcal{K}_{(r)}=K^{rr}(\pd_r\xi_1)^2+((w-divX)V_{(r)}^{11}-X^\mu\pd_\mu V_{(r)}^{11})(\xi_1)^2+\nabla^rm_r^{11}(\xi_1)^2+2m_r^{11}\xi_1\pd^r\xi_1$$
$$\mathcal{K}_{(\theta)}=K^{\theta\theta}(\pd_\theta\xi_1)^2+((w-divX)V_{(\theta)}^{11}-X^\mu\pd_\mu V_{(\theta)}^{11})(\xi_1)^2+\nabla^\theta m_\theta^{11}(\xi_1)^2+2m_\theta^{11}\xi_1\pd^\theta\xi_1$$
$$\mathcal{K}_{(2)}=K^{rr}(\pd_r\xi_2)^2+K^{\theta\theta}(\pd_\theta\xi_2)^2,$$
and
\begin{align*}
V_{(r)}^{11} &= g^{rr}A^{-2}(\pd_rA)^2 = \frac4{r^2}\left(1-\frac{2M}r\right) \\
V_{(\theta)}^{11} &= g^{\theta\theta}A^{-2}(\pd_\theta A)^2 = \frac{4}{r^2}\cot^2\theta.
\end{align*}
\end{lemma}
\begin{proof}
Since $R_{\mu\nu ab}=0$, it suffices to show that
$$\mathcal{K}=K^{\mu\nu}D_\mu\xi\cdot D_\nu\xi +((w-divX)V^{ab}-D_XV^{ab})\xi_a\xi_b+D^\mu m_\mu^{ab}\xi_a\xi_b+2m_\mu^{ab}\xi_aD^\mu \xi_b.$$
Since $K^{tt}=0$, this can be directly verified.
\end{proof}

In what follows, we examine the positivity properties of each of the terms $\mathcal{K}_{(\theta)}$ (Lemma \ref{0_K_theta_lem}), $\mathcal{K}_{(r)}$ (Lemma \ref{0_K_r_lem}), and $\mathcal{K}_{(2)}$ (Lemma \ref{0_K_2_lem}).

First, we show that even after subtracting a good term, the quantity $\mathcal{K}_{(\theta)}$ controls two angular terms. See the remark following the proof of this lemma for an interpretation of the subtracted term.
\begin{lemma}\label{0_K_theta_lem}
For all $\epsilon>0$,
$$\frac{\chi_{trap}}{r}|\sla\nabla\xi_1|^2+\chi_{trap}\frac{\cot^2\theta}{r^3}(\xi_1)^2\lesssim_\epsilon \mathcal{K}_{(\theta)}-(2-\epsilon)\left(-\frac{u\pd_r v}{r^2}\right)(\xi_1)^2.$$
\end{lemma}
\begin{proof}
Recall that
$$\mathcal{K}_{(\theta)}=K^{\theta\theta}(\pd_\theta\xi_1)^2+((w-divX)V_{(\theta)}^{11}-X^\mu\pd_\mu V_{(\theta)}^{11})(\xi_1)^2+D^\theta m_\theta^{11}(\xi_1)^2+2m_\theta^{11}\xi_1\pd^\theta\xi_1.$$
We calculate (using Lemma \ref{Kmunu_lem} in the first line and Lemma \ref{muvVab_lem} in the second line)
\begin{align*}
K^{\theta\theta}(\pd_\theta\xi_1)^2&=-\frac{u\pd_rv}{r^2}(\pd_\theta\xi_1)^2 \\
((w-divX)V_{(\theta)}^{11}-X^r\pd_rV_{(\theta)}^{11})(\xi_1)^2&=-\frac{u\pd_rv}{r^2}4\cot^2\theta(\xi_1)^2 \\
\nabla^\theta m_\theta^{11}(\xi_1)^2&=\frac1{r^2\sin\theta}\pd_\theta(\sin\theta (2-\epsilon)u\pd_rv\cot\theta)(\xi_1)^2 \\
&=-\frac{u\pd_rv}{r^2}(2-\epsilon)(\xi_1)^2 \\
2m_\theta^{11}\xi_1\pd^\theta\xi_1&=2(2-\epsilon)u\pd_rv\cot\theta \xi_1r^{-2}\pd_\theta\xi_1.
\end{align*}
Summing these, we obtain
$$\mathcal{K}_{(\theta)}=-\frac{u\pd_rv}{r^2}\left[(\pd_\theta\xi_1)^2+4\cot^2\theta(\xi_1)^2+(2-\epsilon)(\xi_1)^2-2(2-\epsilon)\cot\theta\xi_1\pd_\theta\xi_1\right].$$
Now,
\begin{align*}
-2(2-\epsilon)\cot\theta\xi_1\pd_\theta\xi_1 &=  \left(1-\frac{\epsilon}2\right)\left(-2\cot\theta\xi_1\pd_\theta\xi_1\right) \\
&=\left(1-\frac{\epsilon}2\right)\left( (\pd_\theta\xi_1-2\cot\theta\xi_1)^2-(\pd_\theta\xi_1)^2-4\cot^2\theta(\xi_1)^2\right).
\end{align*}
Therefore,
$$\mathcal{K}_{(\theta)}=-\frac{u\pd_rv}{r^2}\left[\frac{\epsilon}{2}(\pd_\theta\xi_1)^2+\frac{\epsilon}{2}4\cot^2\theta(\xi_1)^2+\left(1-\frac{\epsilon}2\right)(\pd_\theta\xi_1-2\cot\theta\xi_1)^2+(2-\epsilon)(\xi_1)^2\right].$$
Put another way,
\begin{multline*}
\mathcal{K}_{(\theta)}-(2-\epsilon)\left(-\frac{u\pd_rv}{r^2}\right)(\xi_1)^2 \\
=-\frac{u\pd_rv}{r^2}\left[\frac{\epsilon}{2}(\pd_\theta\xi_1)^2+\frac{\epsilon}{2}4\cot^2\theta(\xi_1)^2+\left(1-\frac{\epsilon}2\right)(\pd_\theta\xi_1-2\cot\theta\xi_1)^2\right].
\end{multline*}
Since $-u\pd_rv\approx \frac{\chi_{trap}}{r}$, the bound stated in Lemma \ref{0_K_theta_lem} follows.
\end{proof}
\begin{remark}
The case $\epsilon=0$ corresponds directly to a calculation for the $(\phi,\psi)$ system of the next section. Indeed,
$$\pd_\theta\xi_1-2\cot\theta\xi_1=A\pd_\theta\psi.$$
This is precisely the need for the $m_\theta^{11}$ component. It accounts for the fact that integrated quantities for $\xi_a$ will differ from their $(\phi,\psi)$ counterparts up to the divergence of a term with a $\theta$ component. For another example, see Lemma \ref{xi_le_phi_psi_2_lem}.
\end{remark}

Next, we show that if the term subtracted from $\mathcal{K}_{(\theta)}$ is added to the term $\mathcal{K}_{(r)}$, then the result is a positive quantity.
\begin{lemma}\label{0_K_r_lem}
If $\epsilon>0$ is sufficiently small, then there exists a function $\chi$ (determining the component $m_r^{11}$) such that
$$\frac{M^2}{r^3}\left(1-\frac{2M}r\right)^2(\pd_r\xi_1)^2+\frac{1}{r^3}(\xi_1)^2 \lesssim \mathcal{K}_{(r)}+(2-\epsilon)\left(-\frac{u\pd_r v}{r^2}\right)(\xi_1)^2.$$
\end{lemma}
\begin{proof}
Recall that
$$\mathcal{K}_{(r)}=K^{rr}(\pd_r\xi_1)^2+((w-divX)V_{(r)}^{11}-X^\mu\pd_\mu V_{(r)}^{11})(\xi_1)^2+D^rm_r^{11}(\xi_1)^2+2m_r^{11}\xi_1\pd^r\xi_1.$$
We calculate  (using Lemma \ref{Kmunu_lem} in the first line and Lemma \ref{muvVab_lem} in the second line)
\begin{align*}
K^{rr}(\pd_r\xi_1)^2 &= 2\left(\frac{u'}{r^2}-\frac{2u}{r^3}\right)\left(1-\frac{2M}r\right)^2(\pd_r\xi_1)^2 \\
((w-divX)V_{(r)}^{11}-X^\mu\pd_\mu V_{(r)}^{11})(\xi_1)^2 &= -\frac{u}{r^2}\pd_r(r^2 vV_{(r)}^{11}) (\xi_1)^2\\
&= -u\pd_r\left(\frac{4}{r^2}\left(1-\frac{2M}r\right)^2\right)(r^{-1}\xi_1)^2 \\
\nabla^rm_r^{11}(\xi_1)^2 &= \frac1{r^2}\pd_r(r^2 g^{rr}m_r^{11}) (\xi_1)^2\\
&= \pd_r\left(\frac{4\chi}{r^2}\left(1-\frac{2M}r\right)^2\right)(r^{-1}\xi_1)^2 \\
2m_r^{11}\xi_1\pd^r\xi_1 &= \frac{4\chi}{r^3}\left(1-\frac{2M}r\right)^22(r^{-1}\xi_1)\pd_r\xi_1.
\end{align*}
Summing these, we obtain
\begin{multline*}
\mathcal{K}_{(r)} =
2\left(\frac{u'}{r^2}-\frac{2u}{r^3}\right)\left(1-\frac{2M}r\right)^2(\pd_r\xi_1)^2
+(\chi-u)\pd_r\left(\frac{4}{r^2}\left(1-\frac{2M}r\right)^2\right)(r^{-1}\xi_1)^2 \\
+\chi'\frac{4}{r^2}\left(1-\frac{2M}r\right)^2(r^{-1}\xi_1)^2 + \frac{4\chi}{r^3}\left(1-\frac{2M}r\right)^22(r^{-1}\xi_1)\pd_r\xi_1.
\end{multline*}
Now,
\begin{multline*}
\frac{4\chi}{r^3}\left(1-\frac{2M}r\right)^22(r^{-1}\xi_1)\pd_r\xi_1 \\
= \frac{4\chi}{r^3}\left(1-\frac{2M}r\right)^2(\pd_r\xi_1+r^{-1}\xi_1)^2 - \frac{4\chi}{r^3}\left(1-\frac{2M}r\right)^2(\pd_r\xi_1)^2 -\frac{4\chi}{r^3}\left(1-\frac{2M}r\right)^2(r^{-1}\xi_1)^2.
\end{multline*}
Therefore,
\begin{multline*}
\mathcal{K}_{(r)} =
2\left(\frac{u'}{r^2}-\frac{2u}{r^3}-\frac{2\chi}{r^3}\right)\left(1-\frac{2M}r\right)^2(\pd_r\xi_1)^2
+(\chi-u)\pd_r\left(\frac{4}{r^2}\left(1-\frac{2M}r\right)^2\right)(r^{-1}\xi_1)^2 \\
+4\left(\frac{\chi'}{r^2}-\frac{\chi}{r^3}\right)\left(1-\frac{2M}r\right)^2(r^{-1}\xi_1)^2 + \frac{4\chi}{r^3}\left(1-\frac{2M}r\right)^2(\pd_r\xi_1+r^{-1}\xi_1)^2.
\end{multline*}
We make an initial choice for $\chi$.
$$\chi=\left\{
\begin{array}{ll}
  0 & r<3M \\
  u & r\in[3M,4M]. \\
  u(4M) & 4M\le r
\end{array}
\right.$$
Note that since $u(3M)=0$, this piecewise function is continuous. With this choice for $\chi$, we derive estimates for the terms appearing in the above expression for $\mathcal{K}_{(r)}$. These estimates are given in Lemmas \ref{0_K_r_1_lem}-\ref{0_K_r_3_lem}.
\begin{lemma}\label{0_K_r_1_lem}
For all $r>r_H$,
$$\frac{M^2}{r^3}\left(1-\frac{2M}r\right)^2(\pd_r\xi_1)^2 \lesssim 2\left(\frac{u'}{r^2}-\frac{2u}{r^3}-\frac{2\chi}{r^3}\right)\left(1-\frac{2M}r\right)^2(\pd_r\xi_1)^2.$$
\end{lemma}
\begin{proof}
First, we will show that
\begin{equation}\label{up_4u_eqn}
\left(\frac{u'}{r^2}-\frac{4u}{r^3}\right)1_{r\le 4M}>0.
\end{equation}
Note that $u'\ge 0$ everywhere and that $u<0$ for $r< 3M$. The challenge is to prove the inequality for $r\in [3M,4M]$. Since $u$ vanishes at $3M$,
$$\frac{u'}{r^2}-\frac{4u}{r^3}=\frac{u'}{r^2}-\frac{4}{r^3}\int_{3M}^ru'.$$
We now divide by $w=vu'$, which is constant in the interval $[3M,4M]$.
\begin{align*}
\frac1w\left(\frac{u'}{r^2}-\frac{u}{r^3}\right) &= \frac1{r^2v}-\frac{4}{r^3}\int_{3M}^r \frac1{v}.
\end{align*}
Next, we multiply by $r^2v$ so that
\begin{align*}
\frac{r^2v}w\left(\frac{u'}{r^2}-\frac{u}{r^3}\right) &= 1-\frac{4v}{r}\int_{3M}^r \frac1{v}.
\end{align*}
Now, observe that the quantity $v^{-1}$ is increasing starting at $r=3M$, where $v$ has a minimum. Therefore,
$$1-\frac{4v}{r}\int_{3M}^r\frac1v >1-\frac{4v}{r}\left(\frac1v(r-3M)\right)=1-4\left(1-\frac{3M}r\right)=-3\left(1-\frac{4M}r\right).$$
Given the strict inequality in the first step, this establishes (\ref{up_4u_eqn}).

It remains to examine the interval $[4M,\infty)$. In this range, $u=r^2-c_1^2$ and $\chi=c_2^2$ for constants $c_1$ and $c_2$. Thus,
$$\left(\frac{u'}{r^2}-\frac{2u}{r^3}-\frac{2\chi}{r^3}\right)1_{r\ge 4M} = \left(\frac{2r}{r^2}-\frac{2(r^2-c_1^2)}{r^3}-\frac{2c_2^2}{r^3}\right)1_{r\ge 4M} = \frac{2(c_1^2-c_2^2)}{r^3}1_{r>4M}.$$
It follows that $\frac{u'}{r^2}-\frac{2u}{r^3}-\frac{2\chi}{r^3}$ does not change sign and behaves like $O(r^{-3})$ for large $r$.
\end{proof}

\begin{lemma}\label{0_K_r_2_lem}
Since $\chi-u$ has the same sign as $-\left(1-\frac{4M}r\right)$ and $\chi-u=O(-r^2)$ for large $r$,
\begin{multline*}
-\frac1{M^3}\left(1-\frac{2M}r\right)\left(1-\frac{3M}r\right)1_{r\le 3M}(\xi_1)^2+\frac1{r^3} \left(1-\frac{4M}r\right)^21_{r\ge 4M}(\xi_1)^2 \\
\lesssim (\chi-u)\pd_r\left(\frac{4}{r^2}\left(1-\frac{2M}r\right)^2\right)(r^{-1}\xi_1)^2.
\end{multline*}
\end{lemma}
\begin{proof}
We compute
$$(\chi-u)\pd_r\left(\frac{4}{r^2}\left(1-\frac{2M}r\right)^2\right)=-\frac{4(\chi-u)}{r^3}\left(1-\frac{2M}r\right)\left(1-\frac{4M}r\right).$$
By the choice of $\chi$, this quantity vanishes on the interval $[3M,4M]$. Then, $\chi-u>0$ for $r< 3M$ and $\chi-u<0$ for $r>4M$. The fact that $\chi-u$ vanishes linearly as $r\nearrow 3M$ and as $r\searrow 4M$ accounts for the factor $1-\frac{3M}r$ and the additional factor $1-\frac{4M}r$ respectively. Finally, since $\chi-u=O(-r^2)$, this accounts for the appropriate $r^{-3}$ weight for large $r$.
\end{proof}

\begin{lemma}\label{0_K_r_3_lem}
For all $r>r_H$,
\begin{multline*}
\frac{1}{M^3}\left(1-\frac{3M}r\right)^21_{r\le 3M}(\xi_1)^2+\frac{1}{M^3}1_{3M\le r\le 4M}(\xi_1)^2+\frac{1}{r^3}\left(1-\frac{4M}r\right)1_{r\ge 4M}(\xi_1)^2 \\
\lesssim 4\left(\frac{\chi'}{r^2}-\frac{\chi}{r^3}\right)\left(1-\frac{2M}r\right)^2(r^{-1}\xi_1)^2-2\frac{u\pd_rv}{r^2}(\xi_1)^2.
\end{multline*}
\end{lemma}
\begin{proof}
For the region $r\le 3M$, we have that $\chi=0$, so the inequality reduces to the fact that $\left(1-\frac{3M}r\right)^2\lesssim -u\pd_rv$ for $r\le 3M$. For the region $3M\le r\le 4M$, we have that $\chi=u$. We ignore the $-u\pd_rv$ term, because it has the right sign. Then the inequality reduces to the fact that 
$$\frac{\chi'}{r^2}-\frac{\chi}{r^3}=\frac{u'}{r^3}-\frac{u}{r^4}\ge \frac{u'}{r^3}-\frac{4u}{r^4}\ge 0.$$
The last inequality follows from (\ref{up_4u_eqn}).

It is in the region $r\ge 4M$ where the term $-2\frac{u\pd_rv}{r^2}(\xi_1)^2$ is essential. In this region, $\chi'=0$, so
\begin{multline*}
4\left(\frac{\chi'}{r^2}-\frac{\chi}{r^3}\right)\left(1-\frac{2M}r\right)^2(r^{-1}\xi_1)^2-2\frac{u\pd_rv}{r^2}(\xi_1)^2 \\
= -\frac{4\chi}{r^3}\left(1-\frac{2M}r\right)^2(r^{-1}\xi_1)^2 -2\chi\pd_rv(r^{-1}\xi_1)^2-2(u-\chi)\pd_rv (r^{-2}\xi_1)^2.
\end{multline*}
Now,
$$\frac1{r^3}\left(1-\frac{4M}r\right)(\xi_1)^2 \lesssim -2(u-\chi)\pd_rv(r^{-1}\xi_1)^2.$$
It remains to show that
$$0\le -\frac{4\chi}{r^3}\left(1-\frac{2M}r\right)^2(r^{-1}\xi_1)^2-2\chi\pd_rv(r^{-1}\xi_1)^2.$$
Indeed,
\begin{align*}
-\frac{4\chi}{r^3}\left(1-\frac{2M}r\right)^2-2\chi\pd_rv &= -\frac{4\chi}{r^3}\left(1-\frac{2M}r\right)^2+\frac{4\chi}{r^3}\left(1-\frac{3M}r\right) \\
&= \frac{4\chi}{r^3}\left[-1+\frac{4M}r-\frac{4M^2}{r^2} + 1 -\frac{3M}r\right] \\
&= \frac{4\chi}{r^3}\left[\frac{M}{r}-\frac{4M^2}{r^2}\right] \\
&= \frac{4\chi}{r^3}\frac{M}{r}\left(1-\frac{4M}r\right).
\end{align*}
This completes the proof of Lemma \ref{0_K_r_3_lem}.
\end{proof}

Combining Lemmas \ref{0_K_r_1_lem}-\ref{0_K_r_3_lem}, we have the following estimate.
$$\frac{M^2}{r^3}\left(1-\frac{2M}r\right)^2(\pd_r\xi_1)^2+\frac{1}{r^3}\left|1-\frac{3M}r\right|\left|1-\frac{4M}r\right|(\xi_1)^2 \lesssim \mathcal{K}_{(r)}-2\frac{u\pd_r v}{r^2}(\xi_1)^2.$$
By slightly modifying $\chi$, the weak degeneracies at $r=3M$ and $r=4M$ can be removed and the following estimate can be established.
$$\frac{M^2}{r^3}\left(1-\frac{2M}r\right)^2(\pd_r\xi_1)^2+\frac{1}{r^3}(\xi_1)^2 \lesssim \mathcal{K}_{(r)}-2\frac{u\pd_r v}{r^2}(\xi_1)^2.$$
It follows that if $\epsilon>0$ is sufficiently small, then 
$$\frac{M^2}{r^3}\left(1-\frac{2M}r\right)^2(\pd_r\xi_1)^2+\frac{1}{r^3}(\xi_1)^2 \lesssim \mathcal{K}_{(r)}-(2-\epsilon)\frac{u\pd_r v}{r^2}(\xi_1)^2.$$
This completes the proof of Lemma \ref{0_K_r_lem}.
\end{proof}

\begin{lemma}\label{0_K_2_lem}
$$\frac{M^2}{r^3}\left(1-\frac{2M}r\right)^2(\pd_r\xi_2)^2+\frac{\chi_{trap}}{r}|\sla\nabla\xi_2|^2 \lesssim \mathcal{K}_{(2)}$$
\end{lemma}
\begin{proof}
Recall that
$$\mathcal{K}_{(2)}=K^{rr}(\pd_r\xi_2)^2+K^{\theta\theta}(\pd_\theta\xi_2)^2.$$
This result now follows directly from Lemma \ref{Kmunu_lem}.
\end{proof}

Finally, we are prepared to complete the proof of the partial Morawetz estimate in Schwarzschild. From Lemma \ref{0_divJ_rearranged_lem}, since $\Box_g\xi_a-V_a{}^b\xi_b=0$,
\begin{align*}
divJ &= \mathcal{K}_{(\theta)}+\mathcal{K}_{(r)}+\mathcal{K}_{(2)}-\frac12\Box_g w |\xi|^2 \\
&= \left(\mathcal{K}_{(\theta)}-(2-\epsilon)\left(-\frac{u\pd_r v}{r^2}\right)(\xi_1)^2\right) +\left(\mathcal{K}_{(r)}+(2-\epsilon)\left(-\frac{u\pd_r v}{r^2}\right)(\xi_1)^2\right) \\
&\hspace{4in}+\mathcal{K}_{(2)}-\frac12\Box_g w|\xi|^2.
\end{align*}
From Lemma \ref{0_K_theta_lem},
$$\frac{\chi_{trap}}{r}|\sla\nabla\xi_1|^2+\chi_{trap}\frac{\cot^2\theta}{r^3}(\xi_1)^2\lesssim_\epsilon \mathcal{K}_{(\theta)}-(2-\epsilon)\left(-\frac{u\pd_r v}{r^2}\right)(\xi_1)^2.$$
From Lemma \ref{0_K_r_lem},
$$\frac{M^2}{r^3}\left(1-\frac{2M}r\right)^2(\pd_r\xi_1)^2+\frac{1}{r^3}(\xi_1)^2 \lesssim \mathcal{K}_{(r)}+(2-\epsilon)\left(-\frac{u\pd_r v}{r^2}\right)(\xi_1)^2.$$
From Lemma \ref{0_K_2_lem},
$$\frac{M^2}{r^3}\left(1-\frac{2M}r\right)^2(\pd_r\xi_2)^2+\frac{\chi_{trap}}{r}|\sla\nabla\xi_2|^2 \lesssim \mathcal{K}_{(2)}.$$
And from Lemma \ref{Kmunu_lem},
$$\frac{M}{r^4}1_{r\ge 4M}((\xi_1)^2+(\xi_2)^2)-r^{-2}V_{\epsilon_{temper}}((\xi_1)^2+(\xi_2)^2)
\lesssim -\frac12\Box_g w|\xi|^2.$$
The bound stated in the partial Morawetz estimate is obtained by summing each of these.
\end{proof}

\subsection{The partial Morawetz estimate for the $\xi_a$ system in slowly rotating ($|a|\ll M$) Kerr spacetimes}\label{partial_morawetz_kerr_sec}

Here, we prove the partial Morawetz estimate for slowly rotating Kerr spacetimes by slightly generalizing the much simpler proof for the Schwarzschild case presented previously.

First, we begin with a calculation of quantites that arise in Kerr.
\begin{lemma}\label{small_a_quantities_lem}
Recall that
\begin{align*}
A &= A_1A_2 \\
A_1 &= (r^2+a^2)\sin^2\theta \\
A_2 &= \left(1+\frac{a^2\sin^2\theta}{q^2}\right)\left(1-a^2\sin^2\theta v\right) \\
v &= \frac{r^2-2Mr+a^2}{(r^2+a^2)^2} \\
B &= -2aM(3\cos\theta-\cos^3\theta)-\frac{2a^3M\sin^4\theta\cos\theta}{r^2+a^2\cos^2\theta}.
\end{align*}
The following identities hold.
$$\frac{\pd_\mu A}{A} = \frac{\pd_\mu A_1}{A_1}+\frac{\pd_\mu A_2}{A_2},$$
\begin{align*}
\frac{\pd_r A_1}{A_1} &= \frac{2r}{r^2+a^2} \\
\frac{\pd_\theta A_1}{A_1} &= 2\cot\theta,
\end{align*}
\begin{align*}
\frac{\pd_rA_2}{A_2} &= a^2\sin^2\theta\left(\frac{-2r}{q^2(r^2+a^2)}+\frac{-\pd_rv}{1-a^2\sin^2\theta v}\right) \\
\frac{\pd_\theta A_2}{A_2} &= a^2\sin\theta\left(\frac{4Mr\cos\theta}{q^2(r^2+a^2)(1-a^2\sin^2\theta v)}\right) \\
\frac{\pd_rB}{A} &= a^3\sin^2\theta\left(\frac{4Mr\cos\theta}{q^2(r^2+a^2)^2(1-a^2\sin^2\theta v)}\right) \\
\frac{\pd_\theta B}{A} &= a\sin\theta\left(\frac{2M(2r^2(r^2+a^2)+(r^2-a^2)q^2)}{q^2(r^2+a^2)^2(1-a^2\sin^2\theta v)}\right).
\end{align*}
In particular,
\begin{align*}
\left|\frac{\pd_r A_2}{A_2}\right| &\lesssim \frac{a^2}{r^3} \\
\left|\frac{\pd_\theta A_2}{A_2}\right| &\lesssim \frac{a^2M}{r^3} \\
\left|\frac{\pd_r B}{A}\right| &\lesssim \frac{|a|^3M}{r^5} \\
\left|\frac{\pd_\theta B}{A}\right| &\lesssim \frac{|a|M}{r^2}.
\end{align*}
\end{lemma}
\begin{proof}
These follow from direct calculation.
\end{proof}

We also prove the following lemma, which allows us to estimate pure partial derivatives.
\begin{lemma}\label{Dxi_lem}
$$\left(1-\left|\frac{r\pd_r B}{A}\right|\right)\left[(\pd_r\xi_1)^2+(\pd_r\xi_2)^2\right] \le |D_r\xi|^2 +\left|\frac{r\pd_r B}{A}\right|\frac{|\xi|^2}{r^2}$$
$$\left(1-\left|\frac{\pd_{\theta} B}{A}\right|\right)\left[(\pd_{\theta}\xi_1)^2+(\pd_{\theta}\xi_2)^2\right] \le |D_{\theta}\xi|^2+\left|\frac{\pd_{\theta} B}{A}\right||\xi|^2.$$
\end{lemma}
\begin{proof}
For arbitrary quantities $x$ and $y$ and some positive function $f$, we have
$$(x\pm fy)^2+f (x^2+y^2) = x^2 \pm 2fxy+f^2y^2+fx^2+fy^2 = x^2+f(x \pm y)^2+f^2y^2.$$
Therefore,
$$x^2\le (x \pm fy)^2+ f(x^2+y^2),$$
whence
$$(1-f)x^2 \le (x \pm fy)^2+fy^2.$$

Now, we have that
$$|D_r\xi|^2 = \left(\pd_r\xi_1 +\frac{r\pd_rB}{A}\frac{\xi_2}{r}\right)^2 + \left(\pd_r\xi_2 -\frac{r\pd_rB}{A}\frac{\xi_1}{r}\right)^2$$
and
$$|D_{\theta}\xi|^2 = \left(\pd_{\theta}\xi_1 +\frac{\pd_{\theta}B}{A}\xi_2\right)^2 + \left(\pd_{\theta}\xi_2 -\frac{\pd_{\theta}B}{A}\xi_1\right)^2.$$
At this point, both estimates in the lemma can be easily deduced.
\end{proof}

Finally, we turn to the partial Morawetz estimate for slowly rotating Kerr spacetimes. 
\begin{proposition}\label{partial_morawetz_kerr_prop} (Partial Morawetz estimate for the $\xi_a$ system in slowly rotating Kerr spacetimes) Suppose $|a|/M$ is sufficiently small, and suppose $\xi_a$ satisfies the equation $\Box_g\xi_a - V_a{}^b\xi_b=0$.

Then there exists a current $J$ such that
\begin{multline*}
\frac{M^2}{r^3}\left(1-\frac{r_H}r\right)^2(\pd_r\xi_1)^2+\frac{\chi_{trap}}{r}|\sla\nabla\xi_1|^2+\frac{1}{r^3}(\xi_1)^2 + \chi_{trap}\frac{\cot^2\theta}{r^3}(\xi_1)^2 \\
+ \frac{M^2}{r^3}\left(1-\frac{r_H}r\right)^2(\pd_r\xi_2)^2+\frac{\chi_{trap}}{r}|\sla\nabla\xi_2|^2+\frac{M}{r^4}1_{r\ge r_*}(\xi_2)^2 \\
-q^{-2}V_{\epsilon_{temper}}((\xi_1)^2+(\xi_2)^2) -q^{-2}V_{\epsilon_a}(\xi_2)^2\\
\lesssim div J,
\end{multline*}
where $\chi_{trap}=\left(1-\frac{r_{trap}}r\right)^2$, $V_{\epsilon_{temper}}$ is the potential defined in Lemma \ref{Kmunu_lem}, and $V_{\epsilon_a}$ is a positive function supported on $r\in [r_H,r_*]$ and satisfying $||V_{\epsilon_a}||_{L^1(r)}\le \epsilon_a$ when $|a|/M$ is chosen sufficiently small.
\end{proposition}
The following proof is motivated by the proof of Proposition \ref{partial_morawetz_szd_prop}.
\begin{proof}
We again choose the current
$$J_\mu = J[X,w,m]_\mu,$$
where $J[X,w,m]$ is defined in Lemma \ref{vectorized_current_template_lem}.

We use the same vectorfield
$$X=uv\pd_r$$
and scalar function
$$w=v\pd_r u$$
that are used in the scalar wave equation for Kerr (see Lemma \ref{Kmunu_lem}).

As in the simpler proof for Schwarzschild, we also use a one-form $m_\mu^{ab}$ with components
\begin{align*}
m_\mu^{ab}&=m_\mu^{ij} (e_i)^a(e_j)^b \\
m^{11}_\mu dx^\mu&=\frac{4r^2\chi v}{(r^2+a^2)^2}dr+(2-\epsilon)u\pd_rv\cot\theta d\theta \\
m^{12}_\mu dx^\mu&=m^{21}=m^{22}=0,
\end{align*}
where the function $\chi$ will be defined in Lemma \ref{a_K_r_lem}.

By Lemma \ref{vectorized_current_template_lem}, we have that
\begin{multline*}
div J = K^{\mu\nu}D_\mu\xi\cdot D_\nu\xi-\frac12\Box_g w|\xi|^2+((w-divX)V^{ab}-D_XV^{ab})\xi_a\xi_b \\
+D^\mu m_\mu^{ab}\xi_a\xi_b+2m_\mu^{ab}\xi_aD^\mu\xi_b  +2R_{\mu\nu a b}X^\mu \xi^a D^\nu \xi^b \\
+(\Box_g\xi_a-V_a{}^b\xi_b)(2D_X\xi^a+w\xi^a).
\end{multline*}
We rearrange these terms according to the following lemma.
\begin{lemma}\label{a_divJ_rearranged_lem}
$$divJ = \mathcal{K}-\frac12\Box_g w|\xi|^2+(\Box_g\xi_a-V_a{}^b\xi_b)(2D_X\xi^a+w\xi^a),$$
where
$$\mathcal{K} = \mathcal{K}_{(r)}+\mathcal{K}_{(\theta)}+\mathcal{K}_{(2)}+\mathcal{K}_{(t)}+\mathcal{K}_{(a)},$$
\begin{multline*}
\mathcal{K}_{(r)} = K^{rr}\left(\pd_r\xi_1+\frac{\pd_r B}{A}\xi_2\right)^2+((w-divX)V_{(r)}^{11}+X^\mu\pd_\mu V_{(r)}^{11})(\xi_1)^2 \\ +\nabla^r (m_r^{11})(\xi_1)^2+2m_r^{11}\xi_1\pd^r\xi_1
\end{multline*}
\begin{multline*}
\mathcal{K}_{(\theta)} = K_{Q}^{\theta\theta}\left(\pd_\theta\xi_1+\frac{\pd_\theta B}{A}\xi_2\right)^2+((w-divX)V_{(\theta)}^{11}+X^\mu\pd_\mu V_{(\theta)}^{11})(\xi_1)^2 \\ +\nabla^\theta (m_\theta^{11})(\xi_1)^2+2m_\theta^{11}\xi_1\pd^\theta\xi_1
\end{multline*}
$$\mathcal{K}_{(2)}=K^{rr}\left(\pd_r\xi_2-\frac{\pd_r B}{A}\xi_1\right)^2+K_{Q}^{\theta\theta}\left(\pd_\theta\xi_2-\frac{\pd_\theta B}{A}\xi_1\right)^2$$
$$\mathcal{K}_{(t)}=K_Q^{tt}\left((\pd_t\xi_1)^2+(\pd_t\xi_2)^2\right)$$
\begin{multline*}
\mathcal{K}_{(a)} = ((w-divX)V_{(a)}^{ij}+X^\mu \pd_\mu V_{(a)}^{ij})\xi_i\xi_j +X^\mu V^{ij}D_\mu(e_i^a e_j^b)\xi_a\xi_b +m_\mu^{ij}D^\mu(e_i^a e_j^b)\xi_i\xi_j \\
+2m_\mu^{ij} e_i^aD^\mu (e_j^b)\xi_a\xi_b +2R_{\mu\nu ab}X^\mu\xi^aD^\nu\xi^b
\end{multline*}
and
$$V_{(r)}^{11}=\frac{\Delta}{q^2}\left(\frac{2r}{r^2+a^2}\right)^2$$
$$V_{(\theta)}^{11}=\frac1{q^2}\left(2\cot\theta\right)^2$$
$$V_{(a)}^{ij}=V^{ij}-(V_{(r)}^{11}+V_{(\theta)}^{11})\delta_1^i\delta_1^j.$$
\end{lemma}
\begin{proof}
The proof is the same as the proof of Lemma \ref{0_divJ_rearranged_lem}, except that there are new terms which are grouped into the quantities $\mathcal{K}_{(t)}$ and $\mathcal{K}_{(a)}$.
\end{proof}
In what follows, we examine the positivity properties of each of the terms $\mathcal{K}_{(\theta)}$ (Lemma \ref{a_K_theta_lem}), $\mathcal{K}_{(r)}$ (Lemma \ref{a_K_r_lem}), and $\mathcal{K}_{(2)}$ (Lemma \ref{a_K_2_lem}). We also estimate the error term $\mathcal{K}_{(a)}$ (Lemma \ref{a_K_a_lem}). The term $\mathcal{K}_{(t)}$ has a good sign, but also a factor of $a^2/M^2$, so it is generally ignored.

First, we show that even after subtracting a good term, the quantity $\mathcal{K}_{(\theta)}$ almost controls two angular terms.
\begin{lemma}\label{a_K_theta_lem}
For all $\epsilon>0$, if $|a|/M$ is sufficiently small compared to $\epsilon$, then
$$\frac{\chi_{trap}}{r}|\sla\nabla\xi_1|^2+\chi_{trap}\frac{\cot^2\theta}{r^3}(\xi_1)^2\lesssim_\epsilon \mathcal{K}_{(\theta)}-(2-\epsilon)\left(-\frac{u\pd_r v}{q^2}\right)(\xi_1)^2+\chi_{trap}\frac{|a|M}{r^5}(\xi_2)^2.$$
\end{lemma}
\begin{proof}
Recall that
\begin{multline*}
\mathcal{K}_{(\theta)} = K^{\theta\theta}\left(\pd_\theta\xi_1+\frac{\pd_\theta B}{A}\xi_2\right)^2+((w-divX)V_{(\theta)}^{11}+X^\mu\pd_\mu V_{(\theta)}^{11})(\xi_1)^2 \\+\nabla^\theta (m_\theta^{11})(\xi_1)^2+2m_\theta^{11}\xi_1\pd^\theta\xi_1.
\end{multline*}
By Lemma \ref{Dxi_lem}, since $K^{\theta\theta}\ge 0$,
$$\left(1-\left|\frac{\pd_\theta B}{A}\right|\right)K^{\theta\theta}(\pd_\theta\xi_1)^2 \le K^{\theta\theta}\left(\pd_\theta\xi_1+\frac{\pd_\theta B}{A}\xi_2\right)^2 +K^{\theta\theta}\left|\frac{\pd_\theta B}{A}\right|(\xi_2)^2.$$
We calculate  (using Lemma \ref{Kmunu_lem} in the first line and Lemma \ref{muvVab_lem} in the second line)
\begin{align*}
\left(1-\left|\frac{\pd_\theta B}{A}\right|\right)K^{\theta\theta}(\pd_\theta\xi_1)^2 &= \left(1-\left|\frac{\pd_\theta B}{A}\right|\right)\left(-\frac{u\pd_r v}{q^2}\right)(\pd_\theta\xi_1)^2 \\
((w-div X)V^{11}_{(\theta)}-X^\mu \pd_\mu V_{(\theta)}^{11})(\xi_1)^2 &= -\frac{u\pd_r v}{q^2}4\cot^2\theta (\xi_1)^2 \\
\nabla^\theta m_\theta^{11}(\xi_1)^2 &= \frac{1}{q^2\sin\theta}\pd_\theta\left(q^2\sin\theta \frac{u\pd_rv}{q^2}\cot\theta\right)(\xi_1)^2 \\
&= -\frac{u\pd_r v}{q^2}(2-\epsilon)(\xi_1)^2 \\
2m_\theta^{11}\xi_1\pd^\theta \xi_1 &= 2(2-\epsilon)u\pd_r v\cot\theta \xi_1 q^{-2}\pd_\theta\xi_1
\end{align*}
Proceeding exactly as in the Schwarzschild case, we conclude that 
\begin{multline*}
-\frac{u\pd_r v}{q^2}\left[\left(\frac{\epsilon}{2}-\left|\frac{\pd_\theta B}{A}\right|\right)(\pd_\theta\xi_1)^2+\frac{\epsilon}{2}4\cot^2\theta(\xi_1)^2+\left(1-\frac{\epsilon}{2}\right)(\pd_\theta\xi_1-2\cot\theta\xi_1)^2\right] \\
\le \mathcal{K}_{(\theta)}-(2-\epsilon)\left(-\frac{u\pd_r v}{q^2}\right)(\xi_1)^2+\left|\frac{u\pd_r v}{q^2}\frac{\pd_\theta B}{A}\right|(\xi_2)^2.
\end{multline*}
Lastly, since $-\frac{u\pd_r v}{q^2}\sim \frac{\chi_{trap}}{r^3}$ and $\left|\frac{\pd_\theta B}{A}\right|\lesssim \frac{|a|M}{r^2}$,
$$\left|\frac{u\pd_r v}{q^2}\frac{\pd_\theta B}{A}\right|(\xi_2)^2 \lesssim \chi_{trap}\frac{aM}{r^5}(\xi_2)^2.$$
Thus, the bound stated in Lemma \ref{a_K_theta_lem} follows.
\end{proof}

Next, we show that if the term subtracted from $\mathcal{K}_{(\theta)}$ is added to the term $\mathcal{K}_{(r)}$, then the result is almost a coercive quantity.
\begin{lemma}\label{a_K_r_lem}
If $\epsilon>0$ and $|a|/M$ are both sufficiently small, then there exists a function $\chi$ (determining the component $m_r^{11}$) such that 
$$\frac{M^2}{r^3}\left(1-\frac{r_H}r\right)^2(\pd_r\xi_1)^2+\frac1{r^3}(\xi_1)^2\lesssim \mathcal{K}_{(r)}+(2-\epsilon)\left(-\frac{u\pd_rv}{q^2}\right)(\xi_1)^2+\frac{|a|^3M^3}{r^9}\left(1-\frac{r_H}r\right)^2(\xi_2)^2.$$
\end{lemma}
\begin{proof}
Recall that
\begin{multline*}
\mathcal{K}_{(r)} = K^{rr}\left(\pd_r\xi_1+\frac{\pd_r B}{A}\xi_2\right)^2+((w-divX)V_{(r)}^{11}+X^\mu\pd_\mu V_{(r)}^{11})(\xi_1)^2 \\+\nabla^r (m_r^{11})(\xi_1)^2+2m_r^{11}\xi_1\pd^r\xi_1.
\end{multline*}
By Lemma \ref{Dxi_lem}, since $K^{rr}\ge 0$,
$$\left(1-\left|\frac{r\pd_r B}{A}\right|\right)K^{rr}(\pd_r\xi_1)^2 \le K^{rr}\left(\pd_\theta\xi_1+\frac{\pd_rB}{A}\xi_2\right)^2+K^{rr}\left|\frac{r\pd_rB}{A}\right|\frac{(\xi_2)^2}{r^2}.$$
We calculate (using Lemma \ref{Kmunu_lem} in the first line and Lemma \ref{muvVab_lem} in the second line)
$$\left(1-\left|\frac{r\pd_r B}{A}\right|\right)K^{rr}(\pd_r\xi_1)^2 = 2\left(1-\left|\frac{r\pd_r B}{A}\right|\right)\left(\frac{u'}{r^2+a^2}-\frac{2ru}{(r^2+a^2)^2}\right)\frac{\Delta^2}{q^2(r^2+a^2)}(\pd_r\xi_1)^2$$
\begin{align*}
((w-divX)V_{(r)}^{11}+X^\mu\pd_\mu V_{(r)}^{11})(\xi_1)^2 &= -\frac{u}{q^2}\left(q^2v V_{(r)}^{11}\right) (\xi_1)^2 \\
&=-\frac{u}{q^2}\left(q^2\frac{\Delta}{(r^2+a^2)^2}\frac{\Delta}{q^2}\frac{4r^2}{(r^2+a^2)^2}\right)(\xi_1)^2 \\
&=-\frac{u}{q^2}\left(\frac{4r^2\Delta^2}{(r^2+a^2)^4}\right)(\xi_1)^2
\end{align*}
\begin{align*}
\nabla^r m_r^{11} (\xi_1)^2 &= \frac1{q^2}\pd_r\left(q^2 g^{rr}m_r^{11}\right) (\xi_1)^2 \\
&= \frac1{q^2}\pd_r\left(q^2\frac{\Delta}{q^2}\frac{4r^2\chi v}{(r^2+a^2)^2}\right)(\xi_1)^2 \\
&= \frac1{q^2}\pd_r\left(\frac{4r^2\chi\Delta^2}{(r^2+a^2)^4}\right) (\xi_1)^2
\end{align*}
$$2m_r^{11}\xi_1\pd^r\xi_1 = \frac{4r^2\chi v}{(r^2+a^2)^2}\frac{\Delta}{q^2}  2\xi_1\pd_r\xi_1.$$
Summing these, we obtain
\begin{multline*}
2\left(1-\left|\frac{r\pd_r B}{A}\right|\right)\left(\frac{u'}{r^2+a^2}-\frac{2ru}{(r^2+a^2)^2}\right) \frac{\Delta^2}{q^2(r^2+a^2)}(\pd_r\xi_1)^2  \\
+\frac{(\chi-u)}{q^2}\pd_r\left(\frac{4r^2\Delta^2}{(r^2+a^2)^4}\right)(\xi_1)^2 
+\frac{\chi'}{q^2}\frac{4r^2\Delta^2}{(r^2+a^2)^4}(\xi_1)^2 + \frac{4r^2\chi\Delta^2}{q^2(r^2+a^2)^4}(2\xi_1\pd_r\xi_1)
\\ \le \mathcal{K}_{(r)} + 2\frac{\Delta^2}{q^2(r^2+a^2)}\pd_r\left(\frac{u}{r^2+a^2}\right)\frac{r\pd_r B}{A} r^{-2}(\xi_2)^2.
\end{multline*}
Now,
\begin{multline*}
\frac{4r^2\chi\Delta^2}{q^2(r^2+a^2)^4}2\xi_1\pd_r\xi_1  \\
= \frac{4r\chi\Delta^2}{q^2(r^2+a^2)^3}\left(\pd_r\xi_1+\frac{r}{r^2+a^2}\xi_1\right)^2 -\frac{4r\chi\Delta^2}{q^2(r^2+a^2)^3}(\pd_r\xi_1)^2-\frac{4r^3\chi \Delta^2}{q^2(r^2+a^2)^5}(\xi_1)^2.
\end{multline*}
Therefore,
\begin{multline*}
2\left(\left(1-\left|\frac{r\pd_r B}{A}\right|\right)\left(\frac{u'}{r^2+a^2}-\frac{2ru}{(r^2+a^2)^2}\right)-\frac{2r\chi}{(r^2+a^2)^2} \right) \frac{\Delta^2}{q^2(r^2+a^2)}(\pd_r\xi_1)^2  \\
+\frac{(\chi-u)}{q^2}\pd_r\left(\frac{4r^2\Delta^2}{(r^2+a^2)^4}\right)(\xi_1)^2 
+\left(\frac{\chi'}{r^2+a^2}-\frac{r\chi}{(r^2+a^2)^2}\right)\frac{4r^2\Delta^2}{q^2(r^2+a^2)^3}(\xi_1)^2 \\
+ \frac{4r\chi\Delta^2}{q^2(r^2+a^2)^3}\left(\pd_r\xi_1+\frac{r}{r^2+a^2}\xi_1\right)^2
\\ \le \mathcal{K}_{(r)} + 2\frac{\Delta^2}{q^2(r^2+a^2)}\pd_r\left(\frac{u}{r^2+a^2}\right)\left|\frac{r\pd_r B}{A}\right| r^{-2}(\xi_2)^2.
\end{multline*}
The error term can be estimated as
\begin{align*}
2\frac{\Delta^2}{q^2(r^2+a^2)}\pd_r\left(\frac{u}{r^2+a^2}\right)\left|\frac{r\pd_r B}{A}\right| r^{-2}(\xi_2)^2 &\lesssim \left(1-\frac{r_H}r\right)^2\frac{M^2}{r^3}\frac{|a|^3M}{r^4}\frac{(\xi_2)^2}{r^2}(\xi_2)^2 \\
&\lesssim \frac{|a|^3M^3}{r^9}\left(1-\frac{r_H}r\right)^2(\xi_2)^2.
\end{align*}
We make an initial choice for $\chi$.
$$\chi=\left\{
\begin{array}{ll}
  0 & r<r_{trap} \\
  u & r\in[r_{trap},r_*]. \\
  u(r_*) & r_*\le r
\end{array}
\right.$$
Note that since $u(r_{trap})=0$, this piecewise function is continuous. With this choice for $\chi$, we derive estimates for the terms appearing in the expression for $\mathcal{K}_{(r)}$ above. These estimates are given in Lemmas \ref{a_K_r_1_lem}-\ref{a_K_r_3_lem}.
\begin{lemma}\label{a_K_r_1_lem}
If $|a|/M$ is sufficiently small, then for all $r>r_H$,
\begin{multline*}
\frac{M^2}{r^3}\left(1-\frac{r_H}r\right)^2(\pd_r\xi_1) \\
\lesssim 2\left(\left(1-\left|\frac{r\pd_r B}{A}\right|\right)\left(\frac{u'}{r^2+a^2}-\frac{2ru}{(r^2+a^2)^2}\right)-\frac{2r\chi}{(r^2+a^2)^2} \right) \frac{\Delta^2}{q^2(r^2+a^2)}(\pd_r\xi_1)^2
\end{multline*}
\end{lemma}
\begin{proof}
The case $a=0$ reduces to Lemma \ref{0_K_r_1_lem}. One can check that the inequality is not affected by taking $|a|/M$ to be small.
\end{proof}

\begin{lemma}\label{a_K_r_2_lem}
For all  $|a|<M$,
\begin{multline*}
-\frac1{M^3}\left(1-\frac{r_H}r\right)\left(1-\frac{r_{trap}}r\right)1_{r\le r_{trap}}(\xi_1)^2 +\frac1{r^3}\left(1-\frac{r_*}r\right)^21_{r\ge r_*}(\xi_1)^2 \\
\lesssim \frac{(\chi-u)}{q^2}\pd_r\left(\frac{4r^2\Delta^2}{(r^2+a^2)^4}\right)(\xi_1)^2.
\end{multline*}
\end{lemma}
\begin{proof}
The case $a=0$ reduces to Lemma \ref{0_K_r_2_lem}. Recall that $r_*$ is by definition the radius at which the function $\frac{2r\Delta}{(r^2+a^2)^2}$ has a maximum value. Therefore the function $\frac{4r^2\Delta^2}{(r^2+a^2)^4}$ also has a maximum value at $r_*$. With this fact in mind, one can check that the argument in the proof of Lemma \ref{0_K_r_2_lem} applies for the general case $|a|<M$ after replacing $3M$ and $4M$ with $r_{trap}$ and $r_*$ respectively.
\end{proof}

\begin{lemma}\label{a_K_r_3_lem}
If $|a|/M$ is sufficiently small, then for all $r>r_H$, there exists a constant $c$ such that
\begin{multline*}
\frac1{M^3}\left(1-\frac{r_{trap}}r\right)^21_{r\le r_{trap}}(\xi_1)^2+\frac1{M^3}1_{r_{trap}\le r\le r_*}(\xi_1)^2+\frac1{r^3}\left(1-\frac{r_*}r-c\frac{|a|}{M}\right)1_{r\ge r_*}(\xi_1)^2 \\
\lesssim \left(\frac{\chi'}{r^2+a^2}-\frac{r\chi}{(r^2+a^2)^2}\right)\frac{4r^2\Delta^2}{q^2(r^2+a^2)^3}(\xi_1)^2 +2\left(-\frac{u\pd_r v}{q^2} \right)(\xi_1)^2.
\end{multline*}
\end{lemma}
\begin{proof}
The case $a=0$ reduces to Lemma \ref{0_K_r_3_lem}. One can check that the inequality is not affected by taking $|a|/M$ to be small except for a small neighborhood of $r_*$, which will have an error term for $r>r_*$. This term is accounted for by the introduction of a constant $c$ so that the term
$$\frac{1}{r^3}\left(1-\frac{r_*}{r}-c\frac{|a|}{M}\right)1_{r\ge r_*}(\xi_1)^2,$$
which is slightly negative near $r_*$, is still bounded by the right side.
\end{proof}
Combining Lemmas \ref{a_K_r_1_lem}-\ref{a_K_r_3_lem}, we have the following estimate.
\begin{multline*}
\frac{M^2}{r^3}\left(1-\frac{r_H}r\right)^2(\pd_r\xi_1)^2+\frac{1}{r^3}\left|1-\frac{r_{trap}}r\right|\left|1-\frac{r_*}r\right|(\xi_1)^2 -c\frac{|a|}{M}r^{-3}1_{r\approx r_*}(\xi_1)^2 \\
\lesssim \mathcal{K}_{(r)}-2\frac{u\pd_r v}{q^2}(\xi_1)^2+\frac{|a|^3M^3}{r^9}\left(1-\frac{r_H}r\right)^2(\xi_2)^2.
\end{multline*}
If $|a|/M$ is sufficiently small, then by slightly modifying $\chi$, the weak degeneracies at $r_{trap}$ and $r_*$ can be removed as well as the small error term $-c\frac{|a|}{M}r^{-3}1_{r\approx r_*}(\xi_1)^2$. Then the following estimate can be established.
$$\frac{M^2}{r^3}\left(1-\frac{r_H}r\right)^2(\pd_r\xi_1)^2+\frac{1}{r^3}(\xi_1)^2 
\lesssim \mathcal{K}_{(r)}-2\frac{u\pd_r v}{q^2}(\xi_1)^2+\frac{|a|^3M^3}{r^9}\left(1-\frac{r_H}r\right)^2(\xi_2)^2.$$
It follows that if $\epsilon>0$ is sufficiently small, 
$$\frac{M^2}{r^3}\left(1-\frac{r_H}r\right)^2(\pd_r\xi_1)^2+\frac{1}{r^3}(\xi_1)^2 
\lesssim \mathcal{K}_{(r)}-(2-\epsilon)\frac{u\pd_r v}{q^2}(\xi_1)^2+\frac{|a|^3M^3}{r^9}\left(1-\frac{r_H}r\right)^2(\xi_2)^2.$$
This completes the proof of Lemma \ref{a_K_r_lem}.
\end{proof}

\begin{lemma}\label{a_K_2_lem}
If $|a|/M$ is sufficiently small, then
$$\frac{M^2}{r^3}\left(1-\frac{r_H}r\right)^2(\pd_r\xi_2)^2+\frac{\chi_{trap}}{r}|\sla\nabla\xi_2|^2 \lesssim \mathcal{K}_{(2)} +\frac{|a|M}{r^5}(\xi_1)^2.$$
\end{lemma}
\begin{proof}
Recall that
\begin{align*}
\mathcal{K}_{(2)}&=K^{rr}\left(\pd_r\xi_2-\frac{\pd_r B}{A}\xi_1\right)^2+K_Q^{\theta\theta}\left(\pd_\theta\xi_2-\frac{\pd_\theta B}{A}\xi_1\right)^2 \\
&=K^{rr}\left(\pd_r\xi_2-\frac{r\pd_r B}{A}r^{-1}\xi_1\right)^2+K_Q^{\theta\theta}\left(\pd_\theta\xi_2-\frac{\pd_\theta B}{A}\xi_1\right)^2.
\end{align*}
By applying the procedure in the proof of Lemma \ref{Dxi_lem}, we conclude that
$$\left(1-\left|\frac{r\pd_r B}{A}\right|\right)K^{rr}(\pd_r\xi_2)^2 \le K^{rr}\left(\pd_r\xi_2-\frac{r\pd_r B}{A}r^{-1}\xi_1\right)^2 +\left|\frac{r\pd_r B}{A}\right|K^{rr}r^{-2}(\xi_1)^2.$$
Therefore, if $|a|/M$ is sufficiently small,
$$K^{rr}(\pd_r\xi_2)^2 \lesssim K^{rr}\left(\pd_r\xi_2-\frac{\pd_r B}{A}\xi_1\right)^2 +\frac{|a|^3M}{r^4}K^{rr}r^{-2}(\xi_1)^2.$$
Again by applying the procedure in the proof of Lemma \ref{Dxi_lem}, we also conclude that
$$\left(1-\left|\frac{\pd_\theta B}{A}\right|\right)K_Q^{\theta\theta}(\pd_\theta\xi_2)^2 \le K_Q^{\theta\theta}\left(\pd_\theta\xi_2-\frac{\pd_\theta B}{A}\xi_1\right)^2+\left|\frac{\pd_\theta B}{A}\right|K_Q^{\theta\theta}(\xi_1)^2.$$
Therefore, if $|a|/M$ is sufficiently small,
$$K_Q^{\theta\theta}(\pd_\theta\xi_2)^2 \lesssim K_Q^{\theta\theta}\left(\pd_\theta\xi_2-\frac{\pd_\theta B}{A}\xi_1\right)^2+\frac{|a|M}{r^2}K_Q^{\theta\theta}(\xi_1)^2.$$
From both estimates, since $K^{rr}=O(M^2/r^3)$ and $K_Q^{\theta\theta}=O(r^{-3})$, we have
\begin{multline*}
K^{rr}(\pd_r\xi_2)^2+K_Q^{\theta\theta}(\pd_\theta\xi_2)^2 \\
\lesssim K^{rr}\left(\pd_r\xi_2-\frac{r\pd_r B}{A}\xi_1\right)^2+ K_Q^{\theta\theta}\left(\pd_\theta\xi_2-\frac{\pd_\theta B}{A}\xi_1\right)^2+\frac{|a|M}{r^5}(\xi_1)^2 \\
\lesssim \mathcal{K}_{(2)}+\frac{|a|M}{r^5}(\xi_1)^2.
\end{multline*}
The lemma follows from the lower bound estimates for $K^{rr}$ and $K_Q^{\theta\theta}$.
\end{proof}

Now we estimate the error term $\mathcal{K}_{(a)}$, which did not appear in the Schwarzschild case.
\begin{lemma}\label{a_K_a_lem}
The following bound holds.
$$|\mathcal{K}_{(a)}| \lesssim \chi_{trap}\frac{|a|M}{r^3}\left(|\sla\nabla\xi_1|^2+|\sla\nabla\xi_2|^2\right)+ \frac{|a|M}{r^5}|\xi|^2.$$
\end{lemma}
\begin{proof}
Recall that
\begin{multline*}
\mathcal{K}_{(a)} = ((w-divX)V_{(a)}^{ij}+X^\mu \pd_\mu V_{(a)}^{ij})\xi_i\xi_j +X^\mu V^{ij}D_\mu(e_i^a e_j^b)\xi_a\xi_b +m_\mu^{ij}D^\mu(e_i^a e_j^b)\xi_i\xi_j \\
+2m_\mu^{ij} e_i^aD^\mu (e_j^b)\xi_a\xi_b +2R_{\mu\nu ab}X^\mu\xi^aD^\nu\xi^b,
\end{multline*}
where
$$V_{(a)}^{ij}=V^{ij}-(V_{(r)}^{11}+V_{(\theta)}^{11})\delta_1^i\delta_1^j.$$

There is one term (the last term) in the expression for $\mathcal{K}_{(a)}$ that contains a derivative of $\xi$. We calculate it explicitly using Lemma \ref{bundle_R_calculation_lem} and the fact that $X$ has only $t$ and $r$ components.
\begin{align*}
2R_{\mu\nu a b}X^\mu\xi^aD^\nu\xi^b &= -4\frac{\pd_{[\mu}A\pd_{\nu]}B}{A^2}\epsilon_{ab}X^\mu \xi^aD^\nu\xi^b \\
&= -4\frac{\pd_{[r}A\pd_{\theta]}B}{A^2}\epsilon_{ab}X^r\xi^aD^\theta\xi^b \\
&= -4\frac{\pd_{[r}A\pd_{\theta]}B}{A^2}X^rg^{\theta\theta}\epsilon^{ab}\xi_aD_\theta\xi_b.
\end{align*}
Now,
\begin{align*}
\epsilon^{ab}\xi_aD_\theta\xi_b &= \epsilon^{ab}\xi_aD_\theta(\xi_j(e^j)_b) \\
&= \epsilon^{ij}\xi_i\pd_\theta\xi_j +\epsilon^{ab}(e^i)_a(D_\theta e^j)_b\xi_i\xi_j
\end{align*}
Thus, the one term that contains a derivative of $\xi_i$ is
\begin{multline*}
-4\frac{\pd_{[r}A\pd_{\theta]}B}{A^2}X^rg^{\theta\theta}\epsilon^{ij}\xi_i\pd_\theta\xi_j \\
=2\left(\frac{\pd_\theta A\pd_r B}{A^2}-\frac{\pd_r A\pd_\theta B}{A^2}\right)\frac{uv}{q^2}(\xi_1\pd_\theta\xi_2-\xi_2\pd_\theta\xi_1) \\
\lesssim \frac{uv}{q^2}\left(\frac{|a|^3M}{r^5}+\frac1r\frac{|a|M}{r^2}\right)(|\xi_1\pd_\theta\xi_2|+|\xi_2\pd_\theta\xi_1|) \\
\lesssim \frac{u^2v^2}{q^2}\frac{|a|M}{r^3}((\pd_\theta\xi_1)^2+(\pd_\theta\xi_2)^2)+\frac1{q^2}\frac{|a|M}{r^3}((\xi_1)^2+(\xi_2)^2) \\
\lesssim \chi_{trap}\frac{|a|M}{r^3}\left(|\sla\nabla\xi_1|^2+|\sla\nabla\xi_2|^2\right) + \frac{|a|M}{r^5}|\xi|^2.
\end{multline*}
The remaining terms to estimate are
\begin{multline*}
((w-divX)V_{(a)}^{ij}+X^\mu \pd_\mu V_{(a)}^{ij})\xi_i\xi_j +X^\mu V^{ij}D_\mu(e_i^a e_j^b)\xi_a\xi_b +m_\mu^{ij}D^\mu(e_i^a e_j^b)\xi_i\xi_j \\
+2m_\mu^{ij} e_i^aD^\mu (e_j^b)\xi_a\xi_b -4\frac{\pd_{[r}A\pd_{\theta]}B}{A^2}X^rg^{\theta\theta}\epsilon^{ab}(e^i)_a(D_\theta e^j)_b\xi_i\xi_j.
\end{multline*}
Each of these terms can be estimated by $\frac{|a|M}{r^5}|\xi|^2$.
\end{proof}

Finally, we are prepared to complete the proof of the partial Morawetz estimate in slowly rotating Kerr spacetimes. From Lemma \ref{a_divJ_rearranged_lem}, since $\Box_g\xi_a-V_a{}^b\xi_b=0$,
\begin{align*}
divJ &= \mathcal{K}_{(\theta)}+\mathcal{K}_{(r)}+\mathcal{K}_{(2)}+\mathcal{K}_{(t)}+\mathcal{K}_{(a)}-\frac12\Box_g w |\xi|^2 \\
&= \left(\mathcal{K}_{(\theta)}-(2-\epsilon)\left(-\frac{u\pd_r v}{r^2}\right)(\xi_1)^2\right) 
+\left(\mathcal{K}_{(r)}+(2-\epsilon)\left(-\frac{u\pd_r v}{r^2}\right)(\xi_1)^2\right) \\
&\hspace{3in}+\mathcal{K}_{(2)}+\mathcal{K}_{(t)}+\mathcal{K}_{(a)}-\frac12\Box_g w|\xi|^2.
\end{align*}
From Lemma \ref{a_K_theta_lem},
$$\frac{\chi_{trap}}{r}|\sla\nabla\xi_1|^2+\chi_{trap}\frac{\cot^2\theta}{r^3}(\xi_1)^2\lesssim_\epsilon \mathcal{K}_{(\theta)}-(2-\epsilon)\left(-\frac{u\pd_r v}{q^2}\right)(\xi_1)^2+\chi_{trap}\frac{|a|M}{r^5}(\xi_2)^2.$$
From Lemma \ref{a_K_r_lem},
$$\frac{M^2}{r^3}\left(1-\frac{r_H}r\right)^2(\pd_r\xi_1)^2+\frac1{r^3}(\xi_1)^2\lesssim \mathcal{K}_{(r)}+(2-\epsilon)\left(-\frac{u\pd_rv}{q^2}\right)(\xi_1)^2+\frac{|a|^3M^3}{r^9}\left(1-\frac{r_H}r\right)^2(\xi_2)^2.$$
From Lemma \ref{a_K_2_lem},
$$\frac{M^2}{r^3}\left(1-\frac{r_H}r\right)^2(\pd_r\xi_2)^2+\frac{\chi_{trap}}{r}|\sla\nabla\xi_2|^2 \lesssim \mathcal{K}_{(2)} +\frac{|a|M}{r^5}(\xi_1)^2.$$
From Lemma \ref{a_K_a_lem},
$$|\mathcal{K}_{(a)}| \lesssim \chi_{trap}\frac{|a|M}{r^3}\left(|\sla\nabla\xi_1|^2+|\sla\nabla\xi_2|^2\right) + \frac{|a|M}{r^5}|\xi|^2.$$
And from Lemma \ref{Kmunu_lem},
$$\frac{M}{r^4}1_{r\ge r_*}((\xi_1)^2+(\xi_2)^2)-q^{-2}V_{\epsilon_{temper}}((\xi_1)^2+(\xi_2)^2)
\lesssim -\frac12\Box_g w|\xi|^2.$$
We ignore the term
$$\mathcal{K}_{(t)}=K_Q^{tt}((\pd_t\xi_1)^2+(\pd_t\xi_2)^2)\approx \chi_{trap}\frac{a^2\sin^2\theta}{r^3}((\pd_t\xi_1)^2+(\pd_t\xi_2)^2),$$
because it has a good sign, but vanishes on the axis and is of order $a^2/M^2$.

Combining each of these estimates, we obtain the following estimate.
\begin{multline*}
\frac{M^2}{r^3}\left(1-\frac{r_H}r\right)^2(\pd_r\xi_1)^2+\frac{\chi_{trap}}{r}|\sla\nabla\xi_1|^2+\frac{1}{r^3}(\xi_1)^2 + \chi_{trap}\frac{\cot^2\theta}{r^3}(\xi_1)^2 \\
+ \frac{M^2}{r^3}\left(1-\frac{r_H}r\right)^2(\pd_r\xi_2)^2+\frac{\chi_{trap}}{r}|\sla\nabla\xi_2|^2+\frac{M}{r^4}1_{r\ge r_*}(\xi_2)^2 \\
-q^{-2}V_{\epsilon_{temper}}((\xi_1)^2+(\xi_2)^2) \\
\lesssim div J +\chi_{trap}\frac{|a|M}{r^3}\left(|\sla\nabla\xi_1|^2+|\sla\nabla\xi_2|^2\right) + \frac{|a|M}{r^5}|\xi|^2.
\end{multline*}
By taking $|a|/M$ sufficiently small, the new terms on the right side can be absorbed into terms on the left, except for the term $\frac{|a|M}{r^5}(\xi_2)^2$, which only can be absorbed for $r\ge r_*$. This is the reason for the new error term $-q^{-2}V_{\epsilon_a}(\xi_2)^2$ on the left side of the bound stated in the partial Morawetz estimate.

This concludes the proof of the partial Morawetz estimate in slowly rotating Kerr spacetimes.
\end{proof}

\subsection{The Morawetz estimate}\label{xi_morawetz_sec}

We now prove the Morawetz estimate, using the partial Morawetz estimate established earlier.
\begin{proposition}\label{xi_morawetz_estimate_prop} (Morawetz estimate for the $\xi_a$ system in slowly rotating Kerr spacetimes)
Suppose $|a|/M$ is sufficiently small. Then
\begin{align*}
&\int_{H_{t_1}^{t_2}}q^{-2}|D_\theta \xi|^2+V^{ab}\xi_a\xi_b
+\int_{\Sigma_{t_2}} |D_L\xi|^2+q^{-2}|D_\theta\xi|^2+V^{ab}\xi_a\xi_b+r^{-2}|\xi|^2 +\frac{M^2}{r^2}|D_r\xi|^2  \\
&+\int_{t_1}^{t_2}\int_{\Sigma_t} \left[\frac{M^2}{r^3}(\pd_r\xi_1)^2+\chi_{trap}\left(\frac{M^2}{r^3}(\pd_t\xi_1)^2+\frac{1}{r}|\sla\nabla\xi_1|^2+\frac{\cot^2\theta}{r^3}(\xi_1)^2\right)+\frac{1}{r^3}(\xi_1)^2\right. \\
&\hspace{1.25in}\left.+ \frac{M^2}{r^3}(\pd_r\xi_2)^2+\chi_{trap}\left(\frac{M^2}{r^3}(\pd_t\xi_2)^2+\frac{1}{r}|\sla\nabla\xi_2|^2\right)+\frac{M}{r^4}(\xi_2)^2\right] \\
&\lesssim \int_{\Sigma_{t_1}} |D_L\xi|^2+q^{-2}|D_\theta\xi|^2+V^{ab}\xi_a\xi_b+r^{-2}|\xi|^2 +\frac{M^2}{r^2}|D_r\xi|^2  + Err,
\end{align*}
where $\chi_{trap}=\left(1-\frac{r_{trap}}{r}\right)^2$ and
\begin{align*}
Err &= Err_1+Err_{nl} \\
Err_1 &= \int_{H_{t_1}^{t_2}}|D_t\xi|^2 + \int_{\Sigma_{t_2}}r^{-1}|\xi\cdot D_L\xi|+\frac{M^2}{r^2}\left[\chi_H|D_r\xi|^2+|D_t\xi|^2+r^{-2}|\xi|^2\right] \\
Err_{nl} &= \int_{t_1}^{t_2}\int_{\Sigma_t}|(2D_X\xi^a+w\xi^a)(\Box_g\xi_a-V_a{}^b\xi_b)|,
\end{align*}
where $\chi_H=1-\frac{r_H}{r}$.
\end{proposition}
\begin{proof}
Let
\begin{align*}
X &= X_{\epsilon_{temper}}+\epsilon_{redshift}Y+\pd_t \\
w &= w_{\epsilon_{temper}}+\epsilon_{\pd_t}w_{\pd_t},
\end{align*}
where $X_{\epsilon_{temper}}$ and $w_{\epsilon_{temper}}$ are the vectorfield and function used in the proof of the partial Morawetz estimate (the dependence on the parameter $\epsilon_{temper}$ is made more explicit for the argument that will follow), $Y$ is the redshift vectorfield introduced in \S\ref{k:redshift_correction_sec}, and $w_{\pd_t}$ is a function defined in the following lemma.

First, we establish an estimate for the bulk term.
\begin{lemma}
If $X$ and $w$ are as defined above and $m$ is the one-form used in the proof of the partial Morawetz estimate, then
\begin{multline*}
\int_{t_1}^{t_2}\int_{\Sigma_t} \left[\frac{M^2}{r^3}(\pd_r\xi_1)^2+\chi_{trap}\left(\frac{M^2}{r^3}(\pd_t\xi_1)^2+\frac{1}{r}|\sla\nabla\xi_1|^2+\frac{\cot^2\theta}{r^3}(\xi_1)^2\right)+\frac{1}{r^3}(\xi_1)^2\right. \\
\hspace{1.25in}\left.+ \frac{M^2}{r^3}(\pd_r\xi_2)^2+\chi_{trap}\left(\frac{M^2}{r^3}(\pd_t\xi_2)^2+\frac{1}{r}|\sla\nabla\xi_2|^2\right)+\frac{M}{r^4}(\xi_2)^2\right] \\
\lesssim \int_{t_1}^{t_2}\int_{\Sigma_t}div J[X,w,m]
\end{multline*}
\end{lemma}
\begin{proof}
To prove this estimate, we start with the partial Morawetz estimate (Proposition \ref{partial_morawetz_kerr_prop}) and make a few slight modifications.  We have
\begin{multline*}
\frac{M^2}{r^3}\left(1-\frac{r_H}r\right)^2(\pd_r\xi_1)^2+\frac{\chi_{trap}}{r}|\sla\nabla\xi_1|^2+\frac{1}{r^3}(\xi_1)^2 + \chi_{trap}\frac{\cot^2\theta}{r^3}(\xi_1)^2 \\
+ \frac{M^2}{r^3}\left(1-\frac{r_H}r\right)^2(\pd_r\xi_2)^2+\frac{\chi_{trap}}{r}|\sla\nabla\xi_2|^2+\frac{M}{r^4}1_{r\ge r_*}(\xi_2)^2 \\
-q^{-2}V_{\epsilon_{temper}}((\xi_1)^2+(\xi_2)^2) -q^{-2}V_{\epsilon_a}(\xi_2)^2\\
\lesssim div J[X_{\epsilon_{temper}},w_{\epsilon_{temper}},m].
\end{multline*}
By applying a small constant times the redshift vectorfield $Y$, the degeneracy of the $(\pd_r\xi_i)^2$ terms near the horizon can be removed without significant consequence. That is
\begin{multline*}
\frac{M^2}{r^3}(\pd_r\xi_1)^2+\frac{\chi_{trap}}{r}|\sla\nabla\xi_1|^2+\frac{1}{r^3}(\xi_1)^2 + \chi_{trap}\frac{\cot^2\theta}{r^3}(\xi_1)^2 \\
+ \frac{M^2}{r^3}(\pd_r\xi_2)^2+\frac{\chi_{trap}}{r}|\sla\nabla\xi_2|^2+\frac{M}{r^4}1_{r\ge r_*}(\xi_2)^2 \\
-q^{-2}V_{\epsilon_{temper}}((\xi_1)^2+(\xi_2)^2) -q^{-2}V_{\epsilon_a}(\xi_2)^2\\
\lesssim div J[X_{\epsilon_{temper}}+\epsilon_{redshift}Y,w_{\epsilon_{temper}},m].
\end{multline*}
It is worth mention that $Y$ is supported near the event horizon and $Y^r<0$.

Next, by choosing $\epsilon_{temper}$ and $\epsilon_a$ sufficiently small and applying a local Hardy estimate, we obtain in an integrated sense
\begin{multline*}
\int_{t_1}^{t_2}\int_{\Sigma_t}\left[\frac{M^2}{r^3}(\pd_r\xi_1)^2+\frac{\chi_{trap}}{r}|\sla\nabla\xi_1|^2+\frac{1}{r^3}(\xi_1)^2 + \chi_{trap}\frac{\cot^2\theta}{r^3}(\xi_1)^2\right. \\
\left.+\frac{M^2}{r^3}(\pd_r\xi_2)^2+\frac{\chi_{trap}}{r}|\sla\nabla\xi_2|^2+\frac{M}{r^4}(\xi_2)^2\right] \\
\lesssim \int_{t_1}^{t_2}\int_{\Sigma_t}div J[X_{\epsilon_{temper}}+\epsilon_{redshift}Y,w_{\epsilon_{temper}},m].
\end{multline*}

Next, we add to $w$ a small constant $\epsilon_{\pd_t}$ times a function $w_{\pd_t}$. Note that according to Lemma \ref{vectorized_current_template_lem}, in the linear case,
$$divJ[0,w_{\pd_t}] = w_{\pd_t}D^\lambda\xi\cdot D_\lambda\xi-\frac12\Box_g w_{\pd_t}+w_{\pd_t}V^{ab}\xi_a\xi_b.$$
One particular term in the above formula is $w_{\pd_t}g^{tt}|D_t\xi|^2=w_{\pd_t}g^{tt}[(\pd_t\xi_1)^2+(\pd_t\xi_2)^2]$. Thus, if $w_{\pd_t}\le 0$ so that $w_{\pd_t}g^{tt}\ge 0$, then the result is added control of the time derivatives. As long as $w_{\pd_t}$ vanishes to second order at the trapping radius and decays like $O(-M^2/r^3)$, by taking $\epsilon_{\pd_t}$ is sufficiently small, the remaining terms have no significant effect.
\begin{multline*}
\int_{t_1}^{t_2}\int_{\Sigma_t} \left[\frac{M^2}{r^3}(\pd_r\xi_1)^2+\chi_{trap}\left(\frac{M^2}{r^3}(\pd_t\xi_1)^2+\frac{1}{r}|\sla\nabla\xi_1|^2+\frac{\cot^2\theta}{r^3}(\xi_1)^2\right)+\frac{1}{r^3}(\xi_1)^2\right. \\
\hspace{1.25in}\left.+ \frac{M^2}{r^3}(\pd_r\xi_2)^2+\chi_{trap}\left(\frac{M^2}{r^3}(\pd_t\xi_2)^2+\frac{1}{r}|\sla\nabla\xi_2|^2\right)+\frac{M}{r^4}(\xi_2)^2\right] \\
\lesssim \int_{t_1}^{t_2}\int_{\Sigma_t}div J[X_{\epsilon_{temper}}+\epsilon_{redshift}Y,w_{\epsilon_{temper}}+\epsilon_{\pd_t}w_{\pd_t},m].
\end{multline*}
Finally, since
$$divJ[\pd_t]=0,$$
the estimate is not affected by adding $\pd_t$ to the vectorfield $X$. (The purpose of $\pd_t$ is for the boundary terms.)

This concludes the proof of the lemma.
\end{proof}

Now, all that remains to be done is to investigate the boundary terms. This is the purpose of the remaining few calculations.

We first approximate the vectorfield $X$ and function $w$ in order to compute the boundary terms. Since $X_{\epsilon_{temper}}=uv\pd_r = \frac{u\Delta}{(r^2+a^2)^2}\pd_r$, and since $\pd_ru=2r$ and $w_{\epsilon_{temper}}=\frac{2r\Delta}{(r^2+a^2)^2}$ for $r>r_*$, it follows that for $r>r_*$,
\begin{align*}
X&=\frac{(r^2+a^2-c^2)\Delta}{(r^2+a^2)^2}\pd_r+\pd_t=L+O(M^2/r^2)\pd_r, \\
w&=\frac{2r\Delta}{(r^2+a^2)^2} +\epsilon_{\pd_t}w_{\pd_t}=\frac{2r\Delta}{(r^2+a^2)^2}+O(M^2/r^3).
\end{align*}
We have the following lemma.
\begin{lemma}
\begin{multline*}
-J^t\left[L,\frac{2r\alpha}{r^2+a^2}\right] 
= \frac{r^2+a^2}{q^2}\left(1-\alpha\frac{a^2\sin^2\theta}{r^2+a^2}\right)\alpha^{-1}|D_L\xi|^2+\frac1{q^2}|D_\theta\xi|^2+V^{ab}\xi_a\xi_b \\
+\frac{\alpha+r\alpha'}{q^2}|\xi|^2-\frac1{q^2}\pd_r(r\alpha|\xi|^2)+Err',
\end{multline*}
where
$$|Err'|\lesssim r^{-1}|\xi\cdot D_L\xi|+\frac{a^2}{r^2}\left[\left(1-\frac{r_H}r\right)^2|D_r\xi|^2+r^{-2}|\xi|^2\right].$$
\end{lemma}
\begin{proof}
By a direct calculation,
\begin{align*}
-J^t[L] &= -2D^t\xi\cdot D_L\xi+L^tD^\lambda\xi\cdot D_\lambda\xi+L^tV^{ab}\xi_a\xi_b \\
&= -2(g^{tt}+{}^{(Q)}g^{tt})D_t\xi\cdot D_L\xi+\left[(g^{tt}+{}^{(Q)}g^{tt})|D_t\xi|^2+g^{rr}|D_r\xi|^2+{}^{(Q)}g^{\theta\theta}|D_\theta\xi|^2\right] \\
&\hspace{4.75in}+V^{ab}\xi_a\xi_b \\
&= -(g^{tt}+{}^{(Q)}g^{tt})|D_t\xi|^2-2\alpha (g^{tt}-{}^{(Q)}g^{tt})D_t\xi\cdot D_r\xi+g^{rr}|D_r\xi|^2+{}^{(Q)}g^{\theta\theta}|D_\theta\xi|^2 \\
&\hspace{4.75in}+V^{ab}\xi_a\xi_b \\
&= \frac{r^2+a^2}{q^2}\left(1-\alpha\frac{a^2\sin^2\theta}{r^2+a^2}\right)\alpha^{-1}|D_t\xi|^2+\frac{r^2+a^2}{q^2}\left(1-\alpha\frac{a^2\sin^2\theta}{r^2+a^2}\right)D_t\xi\cdot D_r\xi  \\
&\hspace{2.65in}+\frac{r^2+a^2}{q^2}\alpha|D_r\xi|^2+\frac1{q^2}|D_\theta\xi|^2+V^{ab}\xi_a\xi_b \\
&= \frac{r^2+a^2}{q^2}\left(1-\alpha\frac{a^2\sin^2\theta}{r^2+a^2}\right)\alpha^{-1}|D_L\xi|^2+\frac{r^2+a^2}{q^2}\alpha^2\frac{a^2\sin^2\theta}{r^2+a^2}|D_r\xi|^2+\frac1{q^2}|D_\theta\xi|^2 \\
&\hspace{4.75in}+V^{ab}\xi_a\xi_b.
\end{align*}
Also,
\begin{align*}
-J^t\left[0,\frac{2r\alpha}{r^2+a^2}\right] &= -\frac{2r\alpha}{r^2+a^2}\xi\cdot D^t\xi \\
&= -\frac{2r\alpha}{r^2+a^2}(g^{tt}+{}^{(Q)}g^{tt})\xi\cdot D_t\xi \\
&= \frac{2r}{q^2}\left(1-\alpha\frac{a^2\sin^2\theta}{r^2+a^2}\right)\xi\cdot D_t\xi \\
&= \frac{2r}{q^2}\left(1-\alpha\frac{a^2\sin^2\theta}{r^2+a^2}\right)\xi\cdot D_L\xi-\frac{2r\alpha}{q^2}\left(1-\alpha\frac{a^2\sin^2\theta}{r^2+a^2}\right)\xi\cdot D_r\xi \\
&= \frac{2r}{q^2}\left(1-\alpha\frac{a^2\sin^2\theta}{r^2+a^2}\right)\xi\cdot D_L\xi-\frac{2r\alpha}{q^2}\xi\cdot D_r\xi +\frac{2r\alpha^2 a^2\sin^2\theta}{q^2(r^2+a^2)}\xi\cdot D_r\xi  \\
&= \frac{2r}{q^2}\left(1-\alpha\frac{a^2\sin^2\theta}{r^2+a^2}\right)\xi\cdot D_L\xi +\left(-\frac{1}{q^2}\pd_r(r\alpha |\xi|^2)+\frac{\alpha+r\alpha'}{q^2}|\xi|^2\right) \\
&\hspace{3.15in}+\frac{2r\alpha^2 a^2\sin^2\theta}{q^2(r^2+a^2)}\xi\cdot D_r\xi.
\end{align*}
Comparing to the identity given for $-J^t\left[L,\frac{2r\alpha}{r^2+a^2}\right]$, we see that
$$Err' = \frac{2r}{q^2}\left(1-\alpha\frac{a^2\sin^2\theta}{r^2+a^2}\right)\xi\cdot D_L\xi+ \alpha^2\frac{a^2\sin^2\theta}{q^2}|D_r\xi|^2 +\frac{2r\alpha^2 a^2\sin^2\theta}{q^2(r^2+a^2)}\xi\cdot D_r\xi.$$
Thus,
$$Err'\lesssim r^{-1}|\xi\cdot D_L\xi| +\left(1-\frac{r_H}r\right)\frac{a^2}{r^2}|D_r\xi|^2 +\frac{a^2}{r^2}|\xi|^2.$$
This concludes the proof of the lemma.
\end{proof}

Given that
\begin{align*}
X-L&=O(M^2/r^2)\pd_r, \\
w-\frac{2r\alpha}{r^2+a^2}&=O(M^2/r^2)r^{-1}.
\end{align*}
We can estimate the remainder
$$\left|J^t[X,w]-J^t\left[L,\frac{2r\alpha}{r^2+a^2}\right]\right| \lesssim \frac{M^2}{r^2}\left[\left(1-\frac{r_H}r\right)|D_r\xi|^2+|D_t\xi|^2+r^{-2}|\xi|^2\right].$$

But it is also important to specifically take into account the effect of the redshift vectorfield $Y$ near the horizon, because it will remove the degeneracy of the term $\chi_H|D_r\xi|^2$ on the spacelike hypersurface $\Sigma_{t_2}$.
\begin{lemma}
On the event horizon $H_{t_1}^{t_2}$,
$$q^{-2}|D_\theta \xi|^2+V^{ab}\xi_a\xi_b \lesssim J^r[Y] + |D_t\xi|^2$$
and on the constant-time hypersurface $\Sigma_{t_2}$,
$$\frac{M^2}{r^2}|D_r\xi|^2\lesssim -J^t[Y] +\frac{M^2}{r^2}\left[\chi_H|D_r\xi|^2+|D_t\xi|^2\right],$$
where $\chi_H=1-\frac{r_H}r$.
\end{lemma}
\begin{proof}
Since $Y\approx -\pd_r+c\pd_t$ near the event horizon, given the estimates for $J^\mu[\pd_t]$ in Lemma \ref{h_dt_J_components_lem}, it suffices to compute the components $J^r[-\pd_r]$ and $-J^t[-\pd_r]$.

From Lemma \ref{vectorized_current_template_lem},
$$J^\mu[-\pd_r] = -2g^{\mu\lambda}D_\lambda\xi\cdot D_r\xi +\delta^\mu{}_r D^\lambda\xi\cdot D_\lambda\xi+\delta^\mu{}_rV^{ab}\xi_a\xi_b.$$
Therefore,
\begin{multline*}
J^r[-\pd_r] = -2g^{rr}|D_r\xi|^2-2g^{rt}D_t\xi\cdot D_r\xi \\
+\left(g^{tt}|D_t\xi|^2+2g^{tr}D_r\xi\cdot D_t\xi+g^{rr}|D_r\xi|^2+q^{-2}|D_\theta\xi|^2+\frac{a^2\sin^2\theta}{q^2}|D_t\xi|^2\right) \\
+V^{ab}\xi_a\xi_b.
\end{multline*}
On the event horizon $H_{t_1}^{t_2}$, since $g^{rr}$ vanishes,
$$J^r[-\pd_r]=q^{-2}|D_\theta\xi|^2+V^{ab}\xi_a\xi_b +\left(g^{tt}+\frac{a^2\sin^2\theta}{q^2}\right)|D_t\xi|^2.$$
This implies the first estimate.

On the constant-time hypersurface $\Sigma_{t_2}$,
$$-J^t[-\pd_r] = 2g^{tr}|D_r\xi|^2+2g^{tt}D_t\xi\cdot D_r\xi.$$
Since $g^{tr}>0$ near the horizon, this implies the second estimate.
\end{proof}
This accounts for all of the boundary terms in the Morawetz estimate.
\end{proof}

\section{Estimates for the $(\phi,\psi)$ System}\label{wk:phi_psi_estimates_sec}

The purpose of this section is to prove the energy estimate (Proposition \ref{translated_energy_estimate_prop}) and $r^p$ estimate (Proposition \ref{p_estimates_prop}) for the $(\phi,\psi)$ system. The energy estimate is simply a translation of Proposition \ref{xi_energy_estimate_prop}, which is in terms of $\xi_a$. The $r^p$ estimate is a combination of three other estimates, the $h\pd_t$ estimate (Proposition \ref{translated_h_dt_prop}, translated from \ref{xi_h_dt_prop}), the Morawetz estimate (Proposition \ref{translated_morawetz_prop}, translated from \ref{xi_morawetz_estimate_prop}), and the incomplete $r^p$ estimate (Proposition \ref{incomplete_p_estimates_prop}).

This section is split into three parts. In \S\ref{phi_psi_translating_sec}-\ref{wk:translate_morawetz_sec}, the estimates from the previous section, which are written in terms of $\xi_a$, are translated to be in terms of $(\phi,\psi)$. In \S\ref{p_identity_M_sec}-\ref{p_identity_combined_sec}, the incomplete $r^p$ estimate is proved. Finally, in \S\ref{phi_psi_p_ee_sec} the $r^p$ estimate is proved.

Recall that the $\xi_a$ system corresponds to the linearization
\begin{align*}
X &= A-A\xi_2 \\
Y &= B+A\xi_1
\end{align*}
of the full nonlinear wave map system (\ref{wk:X_eqn}-\ref{wk:Y_eqn}). To solve the full nonlinear wave map system, we will instead use the linearization
\begin{align*}
X &= A+A\phi \\
Y &= B+A^2\psi.
\end{align*}
The motivation for this second linearization is given in \S\ref{ws:phi_psi_sec}. For this reason, from now on, we assume that
\begin{align*}
\xi_1 &= A\psi \\
\xi_2 &= -\phi.
\end{align*}

\subsection{Translating estimates from the $\xi_a$ system to the $(\phi,\psi)$ system}\label{phi_psi_translating_sec}

To begin, we prove Lemmas \ref{phi_psi_le_xi_lem}-\ref{translate_nl_lem}, which allow us to translate estimates for the $\xi_a$ system into estimates for the $(\phi,\psi)$ system. To prove these lemmas, we will make repeated use of the following calculations.
\begin{align}
|D_\mu\xi|^2 &= \left(\pd_\mu\xi_1+\frac{\pd_\mu B}{A}\xi_2\right)^2+\left(\pd_\mu\xi_2-\frac{\pd_\mu B}{A}\xi_1\right)^2 \nonumber \\
&= \left(A\pd_\mu \psi +\frac{\pd_\mu A}{A}A\psi-\frac{\pd_\mu B}{A}\phi\right)^2+\left(\pd_\mu \phi+\frac{\pd_\mu B}{A}A\psi\right)^2 \nonumber \\
&= \left(A\pd_\mu \psi +\frac{\pd_\mu A_1}{A_1}A\psi+\frac{\pd_\mu A_2}{A_2}A\psi-\frac{\pd_\mu B}{A}\phi\right)^2+\left(\pd_\mu \phi+\frac{\pd_\mu B}{A}A\psi\right)^2.\label{D_mu_xi_eqn}
\end{align}
\begin{align}
V^{ab}\xi_a\xi_b &= g^{\mu\nu}\left(\frac{\pd_\mu A}{A}\xi_1+\frac{\pd_\mu B}{A}\xi_2\right)\left(\frac{\pd_\nu A}{A}\xi_1+\frac{\pd_\nu B}{A}\xi_2\right) \nonumber \\
&= \frac{\Delta}{q^2}\left(\frac{\pd_r A}{A}\xi_1+\frac{\pd_r B}{A}\xi_2\right)^2+\frac{1}{q^2}\left(\frac{\pd_\theta A}{A}\xi_1+\frac{\pd_\theta B}{A}\xi_2\right)^2 \nonumber \\
&= \frac{\Delta}{q^2}\left(\frac{\pd_r A}{A}A\psi-\frac{\pd_r B}{A}\phi\right)^2+\frac{1}{q^2}\left(\frac{\pd_\theta A}{A}A\psi-\frac{\pd_\theta B}{A}\phi\right)^2 \nonumber \\
&= \frac{\Delta}{q^2}\left(\frac{\pd_r A_1}{A_1}A\psi+\frac{\pd_rA_2}{A_2}A\psi-\frac{\pd_r B}{A}\phi\right)^2+\frac{1}{q^2}\left(\frac{\pd_\theta A_1}{A_1}A\psi+\frac{\pd_\theta A_2}{A_2}A\psi-\frac{\pd_\theta B}{A}\phi\right)^2. \label{V_xi_xi_eqn}
\end{align}
\begin{equation}\label{D_r_A1_eqn}
\frac{\pd_r A_1}{A_1}=\frac{2r}{r^2+a^2}
\end{equation}
\begin{equation}\label{D_theta_A1_eqn}
\frac{\pd_\theta A_1}{A_1}=2\cot\theta
\end{equation}
\begin{equation}\label{D_r_A2_B_eqn}
\left(\frac{\pd_r A_2}{A}\right)^2+\left(\frac{\pd_r B}{A}\right)^2\lesssim \frac{a^2}{M^2}r^{-2}
\end{equation}
\begin{equation}\label{D_theta_A2_B_eqn}
\left(\frac{\pd_\theta A_2}{A}\right)^2+\left(\frac{\pd_\theta B}{A}\right)^2\lesssim \frac{a^2}{M^2}
\end{equation}

The first lemma allows us to estimate $(\phi,\psi)$ terms by $\xi_a$ terms.
\begin{lemma}\label{phi_psi_le_xi_lem}
$$(\pd_t\phi)^2+A^2(\pd_t\psi)^2=|D_t\xi|^2,$$
$$(\pd_r\phi)^2+A^2\left(\pd_r\psi+\frac{2r}{r^2+a^2}\psi\right)^2\lesssim |D_r\xi|^2+\frac{a^2}{M^2}r^{-2}(\phi^2+A^2\psi^2),$$
$$(\pd_\theta\phi)^2+A^2\left(\pd_\theta\psi+2\cot\theta\psi\right)^2 \lesssim |D_\theta\xi|^2+\frac{a^2}{M^2}(\phi^2+A^2\psi^2)$$
$$A^2\chi_H\left(\frac{2r}{r^2+a^2}\psi\right)^2+A^2q^{-2}\left(2\cot\theta\psi\right)^2\lesssim V^{ab}\xi_a\xi_b+\frac{a^2}{M^2}r^{-2}(\phi^2+A^2\psi^2),$$
where $\chi_H=1-\frac{r_H}r$.
\end{lemma}
\begin{proof}
From equation (\ref{D_mu_xi_eqn}),
\begin{align*}
|D_\mu\xi|^2 = \left(A\pd_\mu \psi +\frac{\pd_\mu A_1}{A_1}A\psi+\frac{\pd_\mu A_2}{A_2}A\psi-\frac{\pd_\mu B}{A}\phi\right)^2+\left(\pd_\mu \phi+\frac{\pd_\mu B}{A}A\psi\right)^2.
\end{align*}
In the case $\mu=t$, since $\pd_tA_1=\pd_tA_2=\pd_t B=0$, we get the identity
$$(\pd_t\phi)^2+A^2(\pd_t\psi)^2=|D_t\xi|^2.$$
Given equation (\ref{D_r_A1_eqn}) and estimate (\ref{D_r_A2_B_eqn}),
\begin{align*}
(\pd_r\phi)^2+A^2\left(\pd_r\psi+\frac{2r}{r^2+a^2}\psi\right)^2 &= (\pd_r\phi)^2+\left(A\pd_r\psi +\frac{\pd_r A_1}{A_1}A\psi\right)^2 \\
&\lesssim |D_r\xi|^2+\frac{a^2}{M^2}r^{-2}(\phi^2+A^2\psi^2).
\end{align*}
Similarly, given equation (\ref{D_theta_A1_eqn}) and estimate (\ref{D_theta_A2_B_eqn}),
\begin{align*}
(\pd_\theta\phi)^2+A^2\left(\pd_\theta\psi+2\cot\theta\psi\right)^2 &= (\pd_\theta\phi)^2+\left(A\pd_\theta\psi +\frac{\pd_\theta A_1}{A_1}A\psi\right)^2 \\
&\lesssim |D_\theta\xi|^2+\frac{a^2}{M^2}(\phi^2+A^2\psi^2).
\end{align*}
Finally, from equation (\ref{V_xi_xi_eqn}),
\begin{align*}
V^{ab}\xi_a\xi_b &= \frac{\Delta}{q^2}\left(\frac{\pd_r A_1}{A_1}A\psi+\frac{\pd_rA_2}{A_2}A\psi-\frac{\pd_r B}{A}\phi\right)^2+\frac{1}{q^2}\left(\frac{\pd_\theta A_1}{A_1}A\psi+\frac{\pd_\theta A_2}{A_2}A\psi-\frac{\pd_\theta B}{A}\phi\right)^2.
\end{align*}
Using equations (\ref{D_r_A1_eqn}) and (\ref{D_theta_A1_eqn}) together with estimates (\ref{D_r_A2_B_eqn}) and (\ref{D_theta_A2_B_eqn}) again, we conclude
\begin{align*}
A^2\chi_H\left(\frac{2r}{r^2+a^2}\psi\right)^2+A^2q^{-2}\left(2\cot\theta\psi\right)^2 &\lesssim \frac{\Delta}{q^2}\left(\frac{\pd_r A_1}{A_1}A\psi\right)^2+\frac1{q^2}\left(\frac{\pd_\theta A_1}{A_1}A\psi\right)^2\\
&\lesssim V^{ab}\xi_a\xi_b+\frac{a^2}{M^2}r^{-2}(\phi^2+A^2\psi^2).
\end{align*}
 
\end{proof}

The next lemma allows us to estimate certain $\xi_a$ terms by $(\phi,\psi)$ terms.
\begin{lemma}\label{xi_le_phi_psi_1_lem}
$$|D_t\xi|^2 = (\pd_t\phi)^2+A^2(\pd_t\psi)^2$$
$$|D_r\xi|^2 \lesssim (\pd_r\phi)^2+A^2\left(\pd_r\psi+\frac{2r}{r^2+a^2}\psi\right)^2+\frac{a^2}{M^2}r^{-2}(\phi^2+A^2\psi^2).$$
\end{lemma}
\begin{proof}
The identity for $|D_t\xi|^2$ was already proved in the previous lemma, but is simply restated in this lemma for the sake of completeness. 

To prove the estimate for $|D_r\xi|^2$, we use equation (\ref{D_mu_xi_eqn}), then equation (\ref{D_r_A1_eqn}) and finally estimate (\ref{D_r_A2_B_eqn}).
\begin{align*}
|D_r\xi|^2 &= \left(A\pd_r \psi +\frac{\pd_r A_1}{A_1}A\psi+\frac{\pd_r A_2}{A_2}A\psi-\frac{\pd_r B}{A}\phi\right)^2+\left(\pd_r \phi+\frac{\pd_r B}{A}A\psi\right)^2 \\
&= \left(A\pd_r \psi +\frac{2r}{r^2+a^2}A\psi+\frac{\pd_r A_2}{A_2}A\psi-\frac{\pd_r B}{A}\phi\right)^2+\left(\pd_r \phi+\frac{\pd_r B}{A}A\psi\right)^2 \\
&\lesssim (\pd_r\phi)^2+A^2\left(\pd_r\psi+\frac{2r}{r^2+a^2}\psi\right)^2+\frac{a^2}{M^2}r^{-2}(\phi^2+A^2\psi^2).
\end{align*}
 
\end{proof}

Note that the previous lemma did not estimate $|D_\theta\xi|^2$ or $V^{ab}\xi_a\xi_b$. The reason is that both of these terms are singluar on the axis and would require a term like $A^2\cot^2\theta \psi^2$ on the right side of an estimate. It turns out that if these terms are combined in just the right way, then up to a divergence term (with only a $\theta$ component), the singularities cancel. This fact should be compared to the remark following the proof of Lemma \ref{0_K_theta_lem}.
\begin{lemma}\label{xi_le_phi_psi_2_lem}
For an arbitrary function $f(r)$,
\begin{multline*}
\int_{\Sigma_t}f(r)\left[q^{-2}|D_\theta\xi|^2+V^{ab}\xi_a\xi_b\right] \lesssim \int_{\Sigma_t}f(r)\left[q^{-2}(\pd_\theta\phi)^2+r^{-2}\phi^2\right] \\
+\int_{\tilde{\Sigma}_t}f(r)\left[q^{-2}(\pd_\theta\psi)^2+r^{-2}\psi^2\right].
\end{multline*}
\end{lemma}
\begin{proof}
From equation (\ref{D_mu_xi_eqn}),
\begin{align*}
|D_\theta\xi|^2 &= \left(A\pd_\theta \psi +\frac{\pd_\theta A}{A}A\psi-\frac{\pd_\theta B}{A}\phi\right)^2+\left(\pd_\theta \phi+\frac{\pd_\theta B}{A}A\psi\right)^2,
\end{align*}
and from equation (\ref{V_xi_xi_eqn}),
\begin{align*}
V^{ab}\xi_a\xi_b &= \frac{\Delta}{q^2}\left(\frac{\pd_r A}{A}A\psi-\frac{\pd_r B}{A}\phi\right)^2+\frac{1}{q^2}\left(\frac{\pd_\theta A}{A}A\psi-\frac{\pd_\theta B}{A}\phi\right)^2.
\end{align*}
The only challenge arises when dealing with the terms that contain the factor $\frac{\pd_\theta A}{A}$, because this quantity diverges on the axis. We will now show that after an integration by parts on the sphere, the divergent parts cancel. For the sake of simplicity, we take $f(r)=1$.
\begin{multline*}
\int_{\Sigma_t}q^{-2} \left(A\pd_\theta \psi +\frac{\pd_\theta A}{A}A\psi-\frac{\pd_\theta B}{A}\phi\right)^2+q^{-2}\left(\frac{\pd_\theta A}{A}A\psi-\frac{\pd_\theta B}{A}\phi\right)^2 \\
=2\pi\int_{r_H}^{\infty}\int_0^\pi \left(A\pd_\theta \psi +\frac{\pd_\theta A}{A}A\psi-\frac{\pd_\theta B}{A}\phi\right)^2+\left(\frac{\pd_\theta A}{A}A\psi-\frac{\pd_\theta B}{A}\phi\right)^2 \sin\theta d\theta dr.
\end{multline*}
The singular part is contained in the following expression, which we evaluate by expanding the squares and integrating by parts.
\begin{multline*}
\int_0^\pi \left(A\pd_\theta \psi+\frac{\pd_\theta A}{A}A\psi\right)^2+\left(\frac{\pd_\theta A}{A}A\psi\right)^2 \sin\theta d\theta \\
= \int_0^\pi \left[(\pd_\theta\psi)^2+2\frac{\pd_\theta A}{A}\psi\pd_\theta\psi +2\left(\frac{\pd_\theta A}{A}\right)^2\psi^2\right]A^2\sin\theta d\theta \\
= \int_0^\pi \left[(\pd_\theta\psi)^2-\frac{1}{A^2\sin\theta}\pd_\theta\left(A^2\sin\theta \frac{\pd_\theta A}{A}\right)\psi^2 +2\left(\frac{\pd_\theta A}{A}\right)^2\psi^2\right]A^2\sin\theta d\theta.
\end{multline*}
It now suffices to show that the quantity
$$-\frac{1}{A^2\sin\theta}\pd_\theta\left(A^2\sin\theta \frac{\pd_\theta A}{A}\right)+2\left(\frac{\pd_\theta A}{A}\right)^2$$
is regular on the axis. For the Schwarzschild case, where $A=r^2\sin^2\theta$, this expression is rather simple.
\begin{align*}
-\frac{1}{A^2\sin\theta}&\pd_\theta\left(A^2\sin\theta \frac{\pd_\theta A}{A}\right)+2\left(\frac{\pd_\theta A}{A}\right)^2 \\
&=-\frac{1}{r^4\sin^5\theta}\pd_\theta\left(r^4\sin^5\theta (2\cot\theta)\right)+2(2\cot\theta)^2 \\
&= -\frac{1}{\sin^5\theta}\pd_\theta(2\sin^4\theta\cos\theta)+8\cot^2\theta \\
&= -8\cot^2\theta +2 + 8\cot^2\theta \\
&= 2.
\end{align*}
\begin{remark}
Recall that a term was exchanged between Lemmas \ref{0_K_theta_lem} and \ref{0_K_r_lem} with a factor of $2-\epsilon$. Furthermore, there was a remark following the proof of Lemma \ref{0_K_theta_lem} claiming that $\epsilon=0$ corresponds to the $(\phi,\psi)$ system. The above calculation that yields the number $2$ is directly related to the $2-\epsilon$ factor in the exchanged term.
\end{remark}

In the Kerr case, since $A_1=(r^2+a^2)\sin^2\theta$, we also have
$$-\frac{1}{A_1^2\sin\theta}\pd_\theta\left(A_1^2\sin\theta \frac{\pd_\theta A_1}{A_1}\right)+2\left(\frac{\pd_\theta A_1}{A_1}\right)^2=2.$$
This fact, together with the fact that $\frac{\pd_\theta A_2}{A_2}=O(a^2/r^2)\sin\theta$, implies
\begin{align*}
-\frac{1}{A^2\sin\theta}&\pd_\theta\left(A^2\sin\theta \frac{\pd_\theta A}{A}\right)+2\left(\frac{\pd_\theta A}{A}\right)^2 \\
&= -\frac{1}{A_1^2A_2^2\sin\theta}\pd_\theta\left(A_1^2A_2^2\sin\theta \left(\frac{\pd_\theta A_1}{A_1}+\frac{\pd_\theta A_2}{A_2}\right)\right)+2\left(\frac{\pd_\theta A_1}{A_1}+\frac{\pd_\theta A_2}{A_2}\right)^2 \\
&= 2 + O(a^2/r^2).
\end{align*}
This concludes the proof. 
\end{proof}

Finally, we prove an identity that relates the linear equation for $\xi_a$ with the linear equation for $(\phi,\psi)$. In particular, the proof will demonstrate the cancellation of two singular terms in the expression involving $\Box_g\psi$, which is directly related to the choice $\xi_1=A\psi$.
\begin{lemma}\label{translate_nl_lem}
If $(\xi_1,\xi_2)=(A\psi,-\phi)$,
$$\xi_a = \xi_1 (e^1)_a+\xi_2 (e^2)_a = A\psi (e^1)_a-\phi (e^2)_a,$$
then
\begin{align*}
(e_1)^a(\Box_g\xi_a-V_a{}^b\xi_b) &= A\left(\Box_g\psi+2\frac{\pd^\alpha A}{A}\pd_\alpha \psi -2\frac{\pd^\alpha B}{A^2}\pd_\alpha \phi -2\frac{\pd^\alpha B\pd_\alpha B}{A^2}\psi \right) \\
&= A\left(\Box_{\tilde{g}}\psi+2\frac{\pd^\alpha A_2}{A_2}\pd_\alpha \psi -2\frac{\pd^\alpha B}{A^2}\pd_\alpha \phi -2\frac{\pd^\alpha B\pd_\alpha B}{A^2}\psi \right) \\
&= A(\Box_{\tilde{g}}\psi-\mathcal{L}_\psi) \\
-(e_2)^a(\Box_g\xi_a-V_a{}^b\xi_b) &= \Box_g\phi +2\pd^\alpha B\pd_\alpha \psi -2\frac{\pd^\alpha B\pd_\alpha B}{A^2}\phi +4\frac{\pd^\alpha A\pd_\alpha B}{A}\psi \\
&= \Box_g\phi-\mathcal{L}_\phi.
\end{align*}
where $\mathcal{L}_\phi$ and $\mathcal{L}_\psi$, which are both linear in $(\phi,\psi)$, are the quantities defined in \S\ref{wk:phi_psi_equations_intro_sec}.
\end{lemma}
\begin{proof}
Recall the following identities.
$$\xi=\xi_1e^1+\xi_2e^2$$
$$D_\alpha e^1 = -\frac{\pd_\alpha B}{A}e^2$$
$$D_\alpha e^2 = \frac{\pd_\alpha B}{A}e^1$$
Since $(A,B)$ solves the nonlinear wave map system (\ref{wk:X_eqn}-\ref{wk:Y_eqn}), we also have the following two equations.
\begin{align*}
\Box_g A &= \frac{\pd^\alpha A\pd_\alpha A}{A}-\frac{\pd^\alpha B\pd_\alpha B}{A} \\
\Box_g B &= 2\frac{\pd^\alpha A\pd_\alpha B}{A}.
\end{align*}
Using the identity for $\Box_g B$, we compute
\begin{align*}
\Box_g e^1 &= D^\alpha\left(-\frac{\pd_\alpha B}{A}e^2\right) \\
&= -\frac{\Box_g B}{A}e^2+\frac{\pd^\alpha A\pd_\alpha B}{A^2}e^2-\frac{\pd^\alpha B\pd_\alpha B}{A^2}e^1 \\
&= -\frac{\pd^\alpha A\pd_\alpha B}{A^2}e^2-\frac{\pd^\alpha B\pd_\alpha B}{A^2}e^1.
\end{align*}
\begin{align*}
\Box_g e^2 &= D^\alpha\left(\frac{\pd_\alpha B}{A}e^1\right) \\
&= \frac{\Box_g B}{A}e^1-\frac{\pd^\alpha A\pd_\alpha B}{A^2}e^1-\frac{\pd^\alpha B\pd_\alpha B}{A^2}e^2 \\
&=\frac{\pd^\alpha A\pd_\alpha B}{A^2}e^1-\frac{\pd^\alpha B\pd_\alpha B}{A^2}e^2.
\end{align*}
Therefore,
\begin{align*}
\Box_g\xi &= \Box_g(\xi_1e^2)+\Box_g(\xi_2e^2) \\
&= (\Box_g\xi_1 e^1+2\pd^\alpha\xi_1 D_\alpha e^1+\xi_1\Box_g e^1)+(\Box_g\xi_2 e^2+2\pd^\alpha \xi_2 D_\alpha e^2+\xi_2\Box_g e^2) \\
&= \left(\Box_g \xi_1 +2\frac{\pd^\alpha B}{A}\pd_\alpha \xi_2 -\frac{\pd^\alpha B\pd_\alpha B}{A^2}\xi_1+\frac{\pd^\alpha A\pd_\alpha B}{A^2}\xi_2\right)e^1 \\
&\hspace{.5in} +\left(\Box_g\xi_2-2\frac{\pd^\alpha B}{A}\pd_\alpha \xi_1 -\frac{\pd^\alpha A\pd_\alpha B}{A^2}\xi_1-\frac{\pd^\alpha B\pd_\alpha B}{A^2}\xi_2\right)e^2.
\end{align*}

Recall also that
\begin{align*}
V&=g^{\alpha\beta}d\Phi_\alpha d\Phi_\beta \\
&=g^{\alpha\beta}\left(\frac{\pd_\alpha A}{A}e_1+\frac{\pd_\alpha B}{A}e_2\right)\left(\frac{\pd_\beta A}{A}e_1+\frac{\pd_\beta B}{A}e_2\right).
\end{align*}
It follows that
$$V\cdot\xi =\left(\frac{\pd^\alpha A}{A}e^1+\frac{\pd^\alpha B}{A}e^2\right)\left(\frac{\pd_\alpha A}{A}\xi_1+\frac{\pd_\alpha B}{A}\xi_2\right).$$

Combining both identities, we conclude
\begin{align*}
\Box_g\xi-V\cdot\xi &= 
\left(\Box_g \xi_1 +2\frac{\pd^\alpha B}{A}\pd_\alpha \xi_2 -\frac{\pd^\alpha A\pd_\alpha A}{A^2}\xi_1-\frac{\pd^\alpha B\pd_\alpha B}{A^2}\xi_1\right)e^1 \\
&\hspace{.5in} +\left(\Box_g\xi_2-2\frac{\pd^\alpha B}{A}\pd_\alpha \xi_1 -2\frac{\pd^\alpha A\pd_\alpha B}{A^2}\xi_1-2\frac{\pd^\alpha B\pd_\alpha B}{A^2}\xi_2\right)e^2.
\end{align*}
Now, we replace $(\xi_1,\xi_2)=(A\psi,-\phi)$. In the first calculation, we also use the identity for $\Box_gA$, and see that the term with the factor $\Box_gA$ cancels with the term with the factor $\frac{\pd^\alpha A\pd_\alpha A}{A^2}$. Since these terms are singular on the axis, this cancellation is the purpose of the choice $\xi_1=A\psi$.
\begin{align*}
(e_1)^a(\Box_g\xi_a-V_a{}^b\xi_b) &= \Box_g \xi_1 +2\frac{\pd^\alpha B}{A}\pd_\alpha \xi_2 -\frac{\pd^\alpha A\pd_\alpha A}{A^2}\xi_1-\frac{\pd^\alpha B\pd_\alpha B}{A^2}\xi_1 \\
&= \Box_g(A\psi) +2\frac{\pd^\alpha B}{A}\pd_\alpha(-\phi) -\frac{\pd^\alpha A\pd_\alpha A}{A^2}(A\psi)-\frac{\pd^\alpha B\pd_\alpha B}{A^2}(A\psi) \\
&= (A\Box_g \psi +2\pd^\alpha A\pd_\alpha\psi +\psi \Box_gA) -2\frac{\pd^\alpha B}{A}\pd_\alpha\phi -\frac{\pd^\alpha A\pd_\alpha A}{A}\psi \\
&\hspace{3.5in}-\frac{\pd^\alpha B\pd_\alpha B}{A}\psi \\
&= A\Box_g\psi +2\pd^\alpha A\pd_\alpha\psi-2\frac{\pd^\alpha B}{A}\pd_\alpha\phi \\
&\hspace{2in}+\left(\Box_gA-\frac{\pd^\alpha A\pd_\alpha A}{A}-\frac{\pd^\alpha B\pd_\alpha B}{A}\right)\psi \\
&= A\Box_g\psi +2\pd^\alpha A\pd_\alpha\psi-2\frac{\pd^\alpha B}{A}\pd_\alpha\phi -2\frac{\pd^\alpha B\pd_\alpha B}{A}\psi \\
&= A\left(\Box_g\psi +2\frac{\pd^\alpha A}{A}\pd_\alpha \psi -2\frac{\pd^\alpha B}{A^2}\pd_\alpha\phi -2\frac{\pd^\alpha B\pd_\alpha B}{A^2}\psi\right).
\end{align*}
This verifies the first identity of the lemma.
\begin{align*}
-(e_2)^a(\Box_g\xi_a&-V_a{}^b\xi_b) \\
&= -\left(\Box_g\xi_2-2\frac{\pd^\alpha B}{A}\pd_\alpha \xi_1 -2\frac{\pd^\alpha A\pd_\alpha B}{A^2}\xi_1-2\frac{\pd^\alpha B\pd_\alpha B}{A^2}\xi_2\right) \\
&=-\left(\Box_g(-\phi)-2\frac{\pd^\alpha B}{A}\pd_\alpha (A\psi) -2\frac{\pd^\alpha A\pd_\alpha B}{A^2}(A\psi)-2\frac{\pd^\alpha B\pd_\alpha B}{A^2}(-\phi)\right) \\
&= \Box_g\phi+2\pd^\alpha B\pd_\alpha \psi +4\frac{\pd^\alpha A\pd_\alpha B}{A}\psi-2\frac{\pd^\alpha B\pd_\alpha B}{A^2}\phi.
\end{align*}
This verifies the second identity of the lemma. 
\end{proof}

\subsection{Translating the energy estimate from the $\xi_a$ system}

We rewrite the energy estimate (Proposition \ref{xi_energy_estimate_prop}) in terms of $(\phi,\psi)$.
\begin{proposition}\label{translated_energy_estimate_prop}(Energy estimate for the $(\phi,\psi)$ system)
\begin{multline*}
\int_{\Sigma_{t_2}}\chi_H(\pd_r\phi)^2+(\pd_t\phi)^2+|\sla\nabla\phi|^2+r^{-2}\phi^2 +\int_{\tilde\Sigma_{t_2}}\chi_H(\pd_r\psi)^2+(\pd_t\psi)^2+|\tsla\nabla\psi|^2+r^{-2}\psi^2 \\
\lesssim \int_{\Sigma_{t_1}}\chi_H(\pd_r\phi)^2+(\pd_t\phi)^2+|\sla\nabla\phi|^2+r^{-2}\phi^2 +\int_{\tilde\Sigma_{t_1}}\chi_H(\pd_r\psi)^2+(\pd_t\psi)^2+|\tsla\nabla\psi|^2+r^{-2}\psi^2 \\
+Err_{nl},
\end{multline*}
where $\chi_H=1-\frac{r_H}r$ and
$$Err_{nl} = \int_{t_1}^{t_2}\int_{\Sigma_t}|\pd_t\phi(\Box_g\phi-\mathcal{L}_\phi)|
+\int_{t_1}^{t_2}\int_{\tilde{\Sigma}_t}|\pd_t\psi(\Box_{\tilde{g}}\psi-\mathcal{L}_\psi)|.$$
\end{proposition}
\begin{proof}
According to Proposition \ref{xi_energy_estimate_prop},
\begin{multline*}
\int_{H_{t_1}^{t_2}}|D_t\xi|^2+\int_{\Sigma_{t_2}}\chi_H|D_r\xi|^2+|D_t\xi|^2+q^{-2}|D_\theta\xi|^2+V^{ab}\xi_a\xi_b \\
\lesssim \int_{\Sigma_{t_1}}\chi_H|D_r\xi|^2+|D_t\xi|^2+q^{-2}|D_\theta\xi|^2+V^{ab}\xi_a\xi_b + \int_{t_1}^{t_2}\int_{\Sigma_t} |D_t\xi^a(\Box_g\xi_a-V_a{}^b\xi_b)|.
\end{multline*}
We ignore the term on the horizon, which is nonnegative.

By a standard Hardy estimate and then Lemma \ref{phi_psi_le_xi_lem},
\begin{multline*}
\int_{\Sigma_{t_2}}\chi_H(\pd_r\phi)^2+(\pd_t\phi)^2+|\sla\nabla\phi|^2+r^{-2}\phi^2 +\int_{\tilde\Sigma_{t_2}}\chi_H(\pd_r\psi)^2+(\pd_t\psi)^2+|\tsla\nabla\psi|^2+r^{-2}\psi^2 \\
\lesssim \int_{\Sigma_{t_2}}\chi_H(\pd_r\phi)^2+(\pd_t\phi)^2+|\sla\nabla\phi|^2 +\int_{\tilde\Sigma_{t_2}}\chi_H(\pd_r\psi)^2+(\pd_t\psi)^2+|\tsla\nabla\psi|^2 \\
\lesssim \int_{\Sigma_{t_2}}\chi_H|D_r\xi|^2+|D_t\xi|^2+q^{-2}|D_\theta\xi|^2+V^{ab}\xi_a\xi_b +\int_{\Sigma_{t_2}}\frac{|a|}{M}r^{-2}(\phi^2+A^2\psi^2).
\end{multline*}
By taking $|a|/M$ sufficiently small, the final term on the right side can be absorbed into the left side.

By Lemmas \ref{xi_le_phi_psi_1_lem} and \ref{xi_le_phi_psi_2_lem},
\begin{multline*}
\int_{\Sigma_{t_1}}\chi_H|D_r\xi|^2+|D_t\xi|^2+q^{-2}|D_\theta\xi|^2+V^{ab}\xi_a\xi_b \\
\lesssim \int_{\Sigma_{t_1}}\chi_H(\pd_r\phi)^2+(\pd_t\phi)^2+|\sla\nabla\phi|^2+r^{-2}\phi^2 +\int_{\tilde\Sigma_{t_1}}\chi_H(\pd_r\psi)^2+(\pd_t\psi)^2+|\tsla\nabla\psi|^2+r^{-2}\psi^2.
\end{multline*}

Finally, by Lemma \ref{translate_nl_lem},
$$\int_{t_1}^{t_2}\int_{\Sigma_t} |D_t\xi^a(\Box_g\xi_a-V_a{}^b\xi_b)|
\le \int_{t_1}^{t_2}\int_{\Sigma_t}|\pd_t\phi(\Box_g\phi-\mathcal{L}_\phi)|
+\int_{t_1}^{t_2}\int_{\tilde{\Sigma}_t}|\pd_t\psi(\Box_{\tilde{g}}\psi-\mathcal{L}_\psi)|.$$

These estimates together prove the proposition. 
\end{proof}

\subsection{Translating the $h\pd_t$ estimate from the $\xi_a$ system}

We rewrite the $h\pd_t$ estimate (Proposition \ref{xi_h_dt_prop}) in terms of $(\phi,\psi)$.
\begin{proposition}\label{translated_h_dt_prop} ($h\pd_t$ estimate for the $(\phi,\psi)$ system)
Fix $\delp>0$ and let $p\le 2-\delp$. Then for all $\epsilon>0$, there is a small constant $c_\epsilon$ and a large constant $C_\epsilon$, such that
\begin{align*}
&\hspace{.5in}\int_{\Sigma_{t_2}}r^{p-2}\left[\chi_H(\pd_r\phi)^2+(\pd_t\phi)^2+|\sla\nabla\phi|^2\right] +\int_{\tilde\Sigma_{t_2}}r^{p-2}\left[\chi_H(\pd_r\psi)^2+(\pd_t\psi)^2+|\tsla\nabla\psi|^2\right] \\
&+\int_{H_{t_1}^{t_2}}(\pd_t\phi)^2 +\int_{\tilde{H}_{t_1}^{t_2}}(\pd_t\psi)^2 \\
&\hspace{1in}+\int_{t_1}^{t_2}\int_{\Sigma_t\cap \{6M<r\}} c_\epsilon r^{p-3}(\lbar\phi)^2 + \int_{t_1}^{t_2}\int_{\tilde\Sigma_t\cap \{6M<r\}} c_\epsilon r^{p-3}\left(\lbar\psi+\frac{\lbar A}{A}\psi\right)^2 \\
&\lesssim \int_{\Sigma_{t_1}}C_\epsilon r^{p-2}\left[\chi_H(\pd_r\phi)^2+(\pd_t\phi)^2+|\sla\nabla\phi|^2+r^{-2}\phi^2\right] \\
&\hspace{2in}+\int_{\tilde\Sigma_{t_1}}C_\epsilon r^{p-2}\left[\chi_H(\pd_r\psi)^2+(\pd_t\psi)^2+|\tsla\nabla\psi|^2+r^{-2}\psi^2\right]  +Err,
\end{align*}
where $\chi_H=1-\frac{r_H}r$ and
\begin{align*}
Err &= Err_1 + Err_2 + Err_3 + Err_{nl} \\
Err_1 &= \int_{t_1}^{t_2}\int_{\Sigma_t\cap\{5M<r\}} \epsilon r^{-1}((L\phi)^2+r^{-2}\phi^2) + \int_{t_1}^{t_2}\int_{\tilde\Sigma_t\cap\{5M<r\}} \epsilon r^{-1}((L\psi)^2+r^{-2}\psi^2) \\
Err_2 &= \int_{\Sigma_{t_2}}\frac{|a|}{M}r^{p-4}(\phi^2+A^2\psi^2) \\
Err_3 &= \int_{t_1}^{t_2}\int_{\Sigma_t\cap\{6M<r\}} c_\epsilon\frac{|a|}{M}r^{p-5}(\phi^2+A^2\psi^2) \\
Err_{nl} &= \int_{t_1}^{t_2}\int_{\Sigma_t}C_\epsilon r^{p-2}|\pd_t\phi(\Box_g\phi-\mathcal{L}_\phi)| +\int_{t_1}^{t_2}\int_{\tilde{\Sigma}_t}C_\epsilon r^{p-2}|\pd_t\psi(\Box_{\tilde{g}}\psi-\mathcal{L}_\psi)|.
\end{align*}
\end{proposition}

\begin{proof}
By Proposition \ref{xi_h_dt_prop} with $R=5M$,
\begin{multline*}
\int_{H_{t_1}^{t_2}}|D_t\xi|^2+\int_{\Sigma_{t_2}}r^{p-2}\left[\chi_H|D_r\xi|^2+|D_t\xi|^2+q^{-2}|D_\theta\xi|^2 +V^{ab}\xi_a\xi_b\right] \\
+\int_{t_1}^{t_2}\int_{\Sigma_t\cap\{6M<r\}}c_\epsilon r^{p-3}|D_{\lbar}\xi|^2 \\
\lesssim \int_{\Sigma_{t_1}}C_\epsilon r^{p-2}\left[\chi_H|D_r\xi|^2+|D_t\xi|^2+q^{-2}|D_\theta\xi|^2+V^{ab}\xi_a\xi_b\right] + \int_{t_1}^{t_2}\int_{\Sigma_t\cap\{5M<r\}}\epsilon r^{-1}|D_L\xi|^2 \\
+  \int_{t_1}^{t_2}\int_{\Sigma_t}C_\epsilon r^{p-2}|D_t\xi^a(\Box_g\xi_a-V_a{}^b\xi_b)|.
\end{multline*}

By Lemma \ref{phi_psi_le_xi_lem},
\begin{multline*}\int_{\Sigma_{t_2}}r^{p-2}\left[\chi_H(\pd_r\phi)^2+(\pd_t\phi)^2+|\sla\nabla\phi|^2\right] +\int_{\tilde\Sigma_{t_2}}r^{p-2}\left[\chi_H(\pd_r\psi)^2+(\pd_t\psi)^2+|\tsla\nabla\psi|^2\right] \\
\lesssim \int_{\Sigma_{t_2}}r^{p-2}\left[\chi_H|D_r\xi|^2+|D_t\xi|^2+q^{-2}|D_\theta\xi|^2+V^{ab}\xi_a\xi_b\right] +\int_{\Sigma_{t_2}}\frac{|a|}{M}r^{p-4}(\phi^2+A^2\psi^2).
\end{multline*}
and
$$\int_{H_{t_1}^{t_2}}(\pd_t\phi)^2+\int_{\tilde{H}_{t_1}^{t_2}}(\pd_t\psi)^2 = \int_{H_{t_1}^{t_2}}|D_t\xi|^2,$$
and
\begin{multline*}
\int_{t_1}^{t_2}\int_{\Sigma_t\cap \{6M<r\}} c_\epsilon r^{p-3}(\lbar\phi)^2 + \int_{t_1}^{t_2}\int_{\tilde\Sigma_t\cap \{6M<r\}} c_\epsilon r^{p-3}\left(\lbar\psi+\frac{\lbar A}{A}\psi\right)^2 \\
\lesssim \int_{t_1}^{t_2}\int_{\Sigma_t\cap\{6M<r\}}c_\epsilon r^{p-3}|D_{\lbar}\xi|^2 + \int_{t_1}^{t_2}\int_{\Sigma_t\cap\{6M<r\}}c_\epsilon \frac{|a|}{M}r^{p-5}(\phi^2+A^2\psi^2).
\end{multline*}
At this point, we have estimated all the terms on the left side of the main estimate. Now, we must estimate the additional terms:
\begin{multline*}
\int_{\Sigma_{t_2}}\frac{|a|}{M}r^{p-4}(\phi^2+A^2\psi^2) + \int_{t_1}^{t_2}\int_{\Sigma_t\cap\{6M<r\}}c_\epsilon \frac{|a|}{M}r^{p-5}(\phi^2+A^2\psi^2) \\
+  \int_{\Sigma_{t_1}}C_\epsilon r^{p-2}\left[\chi_H|D_r\xi|^2+|D_t\xi|^2+q^{-2}|D_\theta\xi|^2+V^{ab}\xi_a\xi_b\right] \\
+ \int_{t_1}^{t_2}\int_{\Sigma_t\cap\{5M<r\}}\epsilon r^{-1}|D_L\xi|^2 +\int_{t_1}^{t_2}\int_{\Sigma_t}C_\epsilon r^{p-2}|D_t\xi^a(\Box_g\xi_a-V_a{}^b\xi_b)|.
\end{multline*}
By definition,
$$\int_{\Sigma_{t_2}}\frac{|a|}{M}r^{p-4}(\phi^2+A^2\psi^2) + \int_{t_1}^{t_2}\int_{\Sigma_t\cap\{6M<r\}}c_\epsilon\frac{|a|}{M}r^{p-5}(\phi^2+A^2\psi^2) = Err_2+Err_3.$$
Also, by Lemmas \ref{xi_le_phi_psi_1_lem} and \ref{xi_le_phi_psi_2_lem},
\begin{multline*}
\int_{\Sigma_{t_1}}C_\epsilon r^{p-2}\left[\chi_H|D_r\xi|^2+|D_t\xi|^2+q^{-2}|D_\theta\xi|^2+V^{ab}\xi_a\xi_b\right] \\
\lesssim  \int_{\Sigma_{t_1}}C_\epsilon r^{p-2}\left[\chi_H(\pd_r\phi)^2+(\pd_t\phi)^2+|\sla\nabla\phi|^2+r^{-2}\phi^2\right] \\
+\int_{\tilde\Sigma_{t_1}}C_\epsilon r^{p-2}\left[\chi_H(\pd_r\psi)^2+(\pd_t\psi)^2+|\tsla\nabla\psi|^2+r^{-2}\psi^2\right]
\end{multline*}
Also,
\begin{multline*}
\int_{t_1}^{t_2}\int_{\Sigma_t\cap\{5M<r\}}\epsilon r^{-1}|D_L\xi|^2 \\
\lesssim \int_{t_1}^{t_2}\int_{\Sigma_t\cap\{5M<r\}} \epsilon r^{-1}((L\phi)^2+r^{-2}\phi^2) + \int_{t_1}^{t_2}\int_{\tilde\Sigma_t\cap\{5M<r\}} \epsilon r^{-1}((L\psi)^2+r^{-2}\psi^2) \\
= Err_1.
\end{multline*}
Finally, by Lemma \ref{translate_nl_lem},
\begin{multline*}
\int_{t_1}^{t_2}\int_{\Sigma_t}C_\epsilon r^{p-2}|D_t\xi^a(\Box_g\xi_a-V_a{}^b\xi_b)| \\
\le \int_{t_1}^{t_2}\int_{\Sigma_t}C_\epsilon r^{p-2}|\pd_t\phi(\Box_g\phi-\mathcal{L}_\phi)| +\int_{t_1}^{t_2}\int_{\tilde{\Sigma}_t}C_\epsilon r^{p-2}|\pd_t\psi(\Box_{\tilde{g}}\psi-\mathcal{L}_\psi)| \\
\le Err_{nl}.
\end{multline*}
These estimates complete the proof. 
\end{proof}

\subsection{Translating the Morawetz estimate from the $\xi_a$ system}\label{wk:translate_morawetz_sec}

We rewrite the Morawetz estimate (Proposition \ref{xi_morawetz_estimate_prop}) in terms of $(\phi,\psi)$.
\begin{proposition}\label{translated_morawetz_prop} (Morawetz estimate for the $(\phi,\psi)$ system in slowly rotating Kerr spacetimes) Suppose $|a|/M$ is sufficiently small. Then
\begin{multline*}
\int_{\Sigma_{t_2}}(L\phi)^2+|\sla\nabla\phi|^2+r^{-2}\phi^2+\frac{M^2}{r^2}(\pd_r\phi)^2 +\int_{\tilde{\Sigma}_{t_2}}(L\psi)^2+|\tsla\nabla\psi|^2+r^{-2}\psi^2+\frac{M^2}{r^2}(\pd_r\psi)^2 \\
+\int_{t_1}^{t_2}\int_{\Sigma_t}\frac{M^2}{r^3}(\pd_r\phi)^2+\chi_{trap}\left(\frac{M^2}{r^3}(\pd_t\phi)^2+\frac1r|\sla\nabla\phi|^2\right)+\frac{M}{r^4}\phi^2 \\
+\int_{t_1}^{t_2}\int_{\tilde{\Sigma}_t}\frac{M^2}{r^3}(\pd_r\psi)^2+\chi_{trap}\left(\frac{M^2}{r^3}(\pd_t\psi)^2+\frac1r|\tsla\nabla\psi|^2\right)+\frac{1}{r^3}\psi^2 \\
\lesssim \int_{\Sigma_{t_1}}(L\phi)^2+|\sla\nabla\phi|^2+r^{-2}\phi^2+\frac{M^2}{r^2}(\pd_r\phi)^2 +\int_{\tilde{\Sigma}_{t_1}}(L\psi)^2+|\tsla\nabla\psi|^2+r^{-2}\psi^2+\frac{M^2}{r^2}(\pd_r\psi)^2 \\
+Err,
\end{multline*}
where $\chi_{trap}=\left(1-\frac{r_{trap}}{r}\right)^2$ and
\begin{align*}
Err &= Err_1+Err_2+Err_{nl} \\
Err_1 &= \int_{\Sigma_{t_2}}r^{-1}|\phi L\phi|+\int_{\tilde{\Sigma}_{t_2}}r^{-1}|\psi L\psi|+r^{-2}\psi^2 \\
Err_2 &= \int_{H_{t_1}^{t_2}}(\pd_t\phi)^2+\int_{H_{t_1}^{t_2}}(\pd_t\psi)^2 +\int_{\Sigma_{t_2}}\frac{M^2}{r^2}\left[\chi_H(\pd_r\phi)^2+(\pd_t\phi)^2\right] \\
&\hspace{1in} +\int_{\tilde{\Sigma}_{t_2}}\frac{M^2}{r^2}\left[\chi_H(\pd_r\psi)^2+(\pd_t\psi)^2\right] \\
Err_{nl} &= \int_{t_1}^{t_2}\int_{\Sigma_t}|(2X(\phi)+w\phi+w_{(a)}\psi)(\Box_g\phi-\mathcal{L}_\phi)| \\
&\hspace{1in} +\int_{t_1}^{t_2}\int_{\tilde{\Sigma}_t}|(2X(\psi)+\tilde{w}\psi+\tilde{w}_{(a)}\phi)(\Box_{\tilde{g}}\psi-\mathcal{L}_\psi)|,
\end{align*}
where $\chi_H=1-\frac{r_H}r$ and the new functions $\tilde{w}$, $w_{(a)}$, and $\tilde{w}_{(a)}$, which are defined in terms of the original vectorfield $X$ and function $w$ used in the proof of Proposition \ref{xi_morawetz_estimate_prop}, are given by the following relations.
\begin{align*}
2X(\phi)+w\phi+w_{(a)}\psi &= X(\phi)+w\phi+\frac{X(B)}{A}A\psi \\
2X(\psi)+\tilde{w}\psi+\tilde{w}_{(a)}\phi &= 2X(\psi)+\left(\frac{X(A)}{A}+w\right)\psi-\frac{X(B)}{A}A^{-1}\phi.
\end{align*}
\end{proposition}
\begin{proof}
From Proposition \ref{xi_morawetz_estimate_prop}, we have
\begin{align*}
&\int_{H_{t_1}^{t_2}}q^{-2}|D_\theta \xi|^2+V^{ab}\xi_a\xi_b
+\int_{\Sigma_{t_2}} |D_L\xi|^2+q^{-2}|D_\theta\xi|^2+V^{ab}\xi_a\xi_b+r^{-2}|\xi|^2 +\frac{M^2}{r^2}|D_r\xi|^2  \\
&+\int_{t_1}^{t_2}\int_{\Sigma_t} \left[\frac{M^2}{r^3}(\pd_r\xi_1)^2+\chi_{trap}\left(\frac{M^2}{r^3}(\pd_t\xi_1)^2+\frac{1}{r}|\sla\nabla\xi_1|^2+\frac{\cot^2\theta}{r^3}(\xi_1)^2\right)+\frac{1}{r^3}(\xi_1)^2\right. \\
&\hspace{1.25in}\left.+ \frac{M^2}{r^3}(\pd_r\xi_2)^2+\chi_{trap}\left(\frac{M^2}{r^3}(\pd_t\xi_2)^2+\frac{1}{r}|\sla\nabla\xi_2|^2\right)+\frac{M}{r^4}(\xi_2)^2\right] \\
&\lesssim \int_{\Sigma_{t_1}} |D_L\xi|^2+q^{-2}|D_\theta\xi|^2+V^{ab}\xi_a\xi_b+r^{-2}|\xi|^2 +\frac{M^2}{r^2}|D_r\xi|^2  + Err',
\end{align*}
where
\begin{align*}
Err' &= Err'_1+Err'_{nl} \\
Err'_1 &= \int_{H_{t_1}^{t_2}}|D_t\xi|^2 + \int_{\Sigma_{t_2}}r^{-1}|\xi\cdot D_L\xi|+\frac{M^2}{r^2}\left[\chi_H|D_r\xi|^2+|D_t\xi|^2+r^{-2}|\xi|^2\right] \\
Err'_{nl} &= \int_{t_1}^{t_2}\int_{\Sigma_t}|(2D_X\xi^a+w\xi^a)(\Box_g\xi_a-V_a{}^b\xi_b)|.
\end{align*}
We ignore the term on the horizon, which has an appropriate sign.

By Lemma \ref{phi_psi_le_xi_lem},
\begin{multline*}
\int_{\Sigma_{t_2}}(L\phi)^2+|\sla\nabla\phi|^2+r^{-2}\phi^2+\frac{M^2}{r^2}(\pd_r\phi)^2 +\int_{\tilde{\Sigma}_{t_2}}(L\psi)^2+|\tsla\nabla\psi|^2+r^{-2}\psi^2+\frac{M^2}{r^2}(\pd_r\psi)^2 \\
\lesssim \int_{\Sigma_{t_2}} |D_L\xi|^2+q^{-2}|D_\theta\xi|^2+V^{ab}\xi_a\xi_b+r^{-2}|\xi|^2 +\frac{M^2}{r^2}|D_r\xi|^2.
\end{multline*}
and
\begin{multline*}
\int_{t_1}^{t_2}\int_{\Sigma_t}\frac{M^2}{r^3}(\pd_r\phi)^2+\chi_{trap}\left(\frac{M^2}{r^3}(\pd_t\phi)^2+\frac1r|\sla\nabla\phi|^2\right)+\frac{M}{r^4}\phi^2 \\
+\int_{t_1}^{t_2}\int_{\tilde{\Sigma}_t}\frac{M^2}{r^3}(\pd_r\psi)^2+\chi_{trap}\left(\frac{M^2}{r^3}(\pd_t\psi)^2+\frac1r|\tsla\nabla\psi|^2\right)+\frac{1}{r^3}\psi^2 \\
\lesssim \int_{t_1}^{t_2}\int_{\Sigma_t} \left[\frac{M^2}{r^3}(\pd_r\xi_1)^2+\chi_{trap}\left(\frac{M^2}{r^3}(\pd_t\xi_1)^2+\frac{1}{r}|\sla\nabla\xi_1|^2+\frac{\cot^2\theta}{r^3}(\xi_1)^2\right)+\frac{1}{r^3}(\xi_1)^2\right. \\
\hspace{1.25in}\left.+ \frac{M^2}{r^3}(\pd_r\xi_2)^2+\chi_{trap}\left(\frac{M^2}{r^3}(\pd_t\xi_2)^2+\frac{1}{r}|\sla\nabla\xi_2|^2\right)+\frac{M}{r^4}(\xi_2)^2\right].
\end{multline*}
By Lemmas \ref{xi_le_phi_psi_1_lem} and \ref{xi_le_phi_psi_2_lem},
\begin{multline*}
\int_{\Sigma_{t_1}} |D_L\xi|^2+q^{-2}|D_\theta\xi|^2+V^{ab}\xi_a\xi_b+r^{-2}|\xi|^2 +\frac{M^2}{r^2}|D_r\xi|^2  \\
\lesssim \int_{\Sigma_{t_1}}(L\phi)^2+|\sla\nabla\phi|^2+r^{-2}\phi^2+\frac{M^2}{r^2}(\pd_r\phi)^2 +\int_{\tilde{\Sigma}_{t_1}}(L\psi)^2+|\tsla\nabla\psi|^2+r^{-2}\psi^2+\frac{M^2}{r^2}(\pd_r\psi)^2.
\end{multline*}
Since 
$$Err'_1\lesssim Err_1+Err_2,$$
it remains to check the nonlinear terms.

We calculate
\begin{align*}
2D_X\xi+w\xi &= 2D_X(\xi_1e^1+\xi_2e^2)+w(\xi_1e^1+\xi_2e^2) \\
&= \left(2X(\xi_1)+w\xi_1+\frac{X(B)}{A}\xi_2\right)e^1+\left(2X(\xi_2)+w\xi_2-\frac{X(B)}{A}\xi_1\right)e^2 \\
&= A\left(2X(\psi)+\left(\frac{X(A)}{A}+w\right)\psi-\frac{X(B)}{A}A^{-1}\phi\right)e^1 \\
&\hspace{2.3in}-\left(2X(\phi)+w\phi+\frac{X(B)}{A}A\psi\right)e^2.
\end{align*}
This motivates the definition
\begin{align*}
2X(\phi)+w\phi+w_{(a)}\psi &= X(\phi)+w\phi+\frac{X(B)}{A}A\psi \\
2X(\psi)+\tilde{w}\psi+\tilde{w}_{(a)}\phi &= 2X(\psi)+\left(\frac{X(A)}{A}+w\right)\psi-\frac{X(B)}{A}A^{-1}\phi.
\end{align*}
Since according to Lemma \ref{translate_nl_lem}
\begin{align*}
(e_1)^a(\Box_g\xi_a-V_a{}^b\xi_b) &= A(\Box_{\tilde{g}}\psi-\mathcal{L}_\psi) \\
-(e_2)^a(\Box_g\xi_a-V_a{}^b\xi_b) &= \Box_g\phi -\mathcal{L}_\phi,
\end{align*}
we conclude that
\begin{multline*}
(2D_X\xi^a+w\xi^a)(\Box_g\xi_a-V_a{}^b\xi_b) \\
=A^2(2X(\psi)+\tilde{w}\psi+\tilde{w}_{(a)}\phi)(\Box_{\tilde{g}}-\mathcal{L}_\psi)+(2X(\phi)+w\phi+w_{(a)}\psi)(\Box_g\phi-\mathcal{L}_\phi).
\end{multline*}
From this identity, it is clear that $Err'_{nl}\le Err_{nl}$. This completes the proof. 
\end{proof}

\subsection{The pre-$r^p$ identity for $\mathcal{M}$}\label{p_identity_M_sec}

At this point, all of the relevant estimates from \S\ref{wk:xi_estimates_sec} have been rewritten in terms of $(\phi,\psi)$. But in order to prove the $r^p$ estimate, we still must prove the incomplete $r^p$ estimate near $i^0$. Such an estimate was proved already for the equation $\Box_g\phi=0$ in Chapter \ref{kerr_chap}, and is restated below. In \S\ref{p_identity_M_tilde_sec}, we will prove an analogous estimate for the equation $\Box_{\tilde{g}}\psi =0$. Then, in \S\ref{p_identity_combined_sec}, we will combine both estimates into a single estimate for the $(\phi,\psi)$ system, which will treate the terms $\mathcal{L}_\phi$ and $\mathcal{L}_\psi$ as error terms.

The following lemma, which has a much simpler, well-known analogue in Minkowski spacetime, is a restatement of Lemma \ref{k:p_ee_identity_lem}.
\begin{lemma}\label{p_identity_phi_lem} (Pre-$r^p$ identity for $\mathcal{M}$) Let $\alpha=\frac{\Delta}{r^2+a^2}$ and $L=\alpha\pd_r+\pd_t$. For any function $f=f(r)$ supported where $r>r_H+\delh$, the following identity holds.
\begin{align*}
&\int_{\Sigma_{t_2}}\left[\left(1-\alpha\frac{a^2\sin^2\theta}{r^2+a^2}\right)\frac{r^2+a^2}{q^2}f\left(\alpha^{-1}L\phi+\frac{r}{r^2+a^2}\phi\right)^2 +\frac{\alpha^{-1}f}{q^2}(\pd_\theta\phi)^2+\epsilon\frac{rf'}{q^2}\phi^2\right. \\
&\hspace{1.2in}\left.+\alpha\frac{a^2\sin^2\theta}{q^2}f\left(\pd_r\phi+\frac{r}{r^2+a^2}\phi\right)^2 +\frac{a^2f}{q^2(r^2+a^2)}\phi^2-\frac1{q^2}\pd_r(rf\phi^2)\right] \\
&+\int_{t_1}^{t_2}\int_{\Sigma_t}
\left[
  \left(\frac{2rf}{r^2+a^2}-f'\right)\frac{Q^{\alpha\beta}}{q^2}\pd_\alpha\phi\pd_\beta\phi
+ \alpha f'\frac{r^2+a^2}{q^2}\left(\alpha^{-1}L\phi+\frac{(1-\epsilon)r}{r^2+a^2}\phi\right)^2 
\vphantom{
+ \epsilon\alpha\left((1-\epsilon)f'-rf''\right)\frac{\phi^2}{q^2} 
+ \alpha'\left(\frac{r^2-a^2}{r^2+a^2}f-\epsilon r f'\right)\frac{\phi^2}{q^2}
+ \frac{a^2}{r^2+a^2}\left(\alpha((1-\epsilon)^2-2)f'-\frac{4\alpha rf}{r^2+a^2}\right)\frac{\phi^2}{q^2}
}\right. \\
&\hspace{2.2in} + \epsilon\alpha\left((1-\epsilon)f'-rf''\right)\frac{\phi^2}{q^2} -\alpha'\alpha^{-2}f\frac{r^2+a^2}{q^2}(L\phi)^2 \\
&\hspace{0.8in}
\left.\vphantom{
  \left(\frac{2rf}{r^2+a^2}-f'\right)\frac{Q^{\alpha\beta}}{q^2}\pd_\alpha\phi\pd_\beta\phi
+ \alpha f'\frac{r^2+a^2}{q^2}\left(\alpha^{-1}L\phi+\frac{(1-\epsilon)r}{r^2+a^2}\phi\right)^2 
+ \epsilon\alpha\left((1-\epsilon)f'-rf''\right)\frac{\phi^2}{q^2} 
}
+ \alpha'\left(-\epsilon r f' +\frac{r^2-a^2}{r^2+a^2}f\right)\frac{\phi^2}{q^2}
+ \frac{a^2}{r^2+a^2}\left(-(1+\epsilon)\alpha f'-\frac{4\alpha rf}{r^2+a^2}\right)\frac{\phi^2}{q^2}
\right] \\
=&\int_{\Sigma_{t_1}}\left[\left(1-\alpha\frac{a^2\sin^2\theta}{r^2+a^2}\right)\frac{r^2+a^2}{q^2}f\left(\alpha^{-1}L\phi+\frac{r}{r^2+a^2}\phi\right)^2 +\frac{\alpha^{-1}f}{q^2}(\pd_\theta\phi)^2+\epsilon\frac{rf'}{q^2}\phi^2\right. \\
&\hspace{1.2in}\left.+\alpha\frac{a^2\sin^2\theta}{q^2}f\left(\pd_r\phi+\frac{r}{r^2+a^2}\phi\right)^2 +\frac{a^2f}{q^2(r^2+a^2)}\phi^2-\frac1{q^2}\pd_r(rf\phi^2)\right] \\
&+\int_{t_1}^{t_2}\int_{\Sigma_t}-\left(2\alpha^{-1}fL\phi+\frac{2rf}{r^2+a^2}\phi\right)\Box_g\phi.
\end{align*}
\end{lemma}

\subsection{The pre-$r^p$ identity for $\tilde{\mathcal{M}}$}\label{p_identity_M_tilde_sec}
We prove the following lemma for the spacetime $\tilde{\mathcal{M}}$ by following a similar procedure to that in the proof of Lemma \ref{k:p_ee_identity_lem}.
\begin{lemma}\label{p_identity_psi_lem} (Pre-$r^p$ identity for $\tilde{\mathcal{M}}$) Let $\alpha=\frac{\Delta}{r^2+a^2}$ and $L=\alpha\pd_r+\pd_t$. For any function $f=f(r)$ supported where $r>r_H+\delh$, the following identity holds.
\begin{align*}
&\int_{\tilde{\Sigma}_{t_2}}\left[\left(1-\alpha\frac{a^2\sin^2\theta}{r^2+a^2}\right)\frac{r^2+a^2}{q^2}f\left(\alpha^{-1}L\psi+\frac{3r}{r^2+a^2}\psi\right)^2 +\frac{\alpha^{-1}f}{q^2}(\pd_\theta\psi)^2 +6f\frac{\psi^2}{q^2} \right.\\
&\hspace{0.95in}+\left.\frac{a^2\sin^2\theta \alpha f}{q^2}\left(\pd_r\psi+\frac{3r}{r^2+a^2}\psi\right)^2 -\frac{a^23f}{r^2+a^2}\frac{\psi^2}{q^2}-\frac1{A_1^2q^2}\pd_r(A_1^23rf\psi^2)\right] \\
&+\int_{t_1}^{t_2}\int_{\tilde{\Sigma}_t}
\left[
  \left(\frac{2rf}{r^2+a^2}-f'\right)\frac{Q^{\alpha\beta}}{q^2}\pd_\alpha\psi\pd_\beta\psi
+ \alpha f'\frac{r^2+a^2}{q^2}\left(\alpha^{-1}L\psi+\frac{3r}{r^2+a^2}\psi\right)^2 
\vphantom{
+ \epsilon\alpha\left((1-\epsilon)f'-rf''\right)\frac{\phi^2}{q^2} 
+ \alpha'\left(\frac{r^2-a^2}{r^2+a^2}f-\epsilon r f'\right)\frac{\phi^2}{q^2}
+ \frac{a^2}{r^2+a^2}\left(\alpha((1-\epsilon)^2-2)f'-\frac{4\alpha rf}{r^2+a^2}\right)\frac{\phi^2}{q^2}
}\right. \\
&\hspace{2.7in} + \frac{6r\alpha(2f-rf')}{r^2+a^2}\frac{\psi^2}{q^2} -\alpha'\alpha^{-2}f\frac{r^2+a^2}{q^2}(L\psi)^2 \\
&\hspace{1.2in}
\left.\vphantom{
  \left(\frac{2rf}{r^2+a^2}-f'\right)\frac{Q^{\alpha\beta}}{q^2}\pd_\alpha\phi\pd_\beta\phi
+ \alpha f'\frac{r^2+a^2}{q^2}\left(\alpha^{-1}L\phi+\frac{(1-\epsilon)r}{r^2+a^2}\phi\right)^2 
+ \epsilon\alpha\left((1-\epsilon)f'-rf''\right)\frac{\phi^2}{q^2} 
}
- \alpha'(r^2+a^2)\pd_r\left(\frac{3r}{r^2+a^2}\right)f\frac{\psi^2}{q^2}
- \frac{a^2}{r^2+a^2}\left(3\alpha f'+\frac{36\alpha rf}{r^2+a^2}\right)\frac{\psi^2}{q^2}
\right] \\
=&\int_{\tilde{\Sigma}_{t_1}}\left[\left(1-\alpha\frac{a^2\sin^2\theta}{r^2+a^2}\right)\frac{r^2+a^2}{q^2}f\left(\alpha^{-1}L\psi+\frac{3r}{r^2+a^2}\psi\right)^2 +\frac{\alpha^{-1}f}{q^2}(\pd_\theta\psi)^2 +6f\frac{\psi^2}{q^2} \right.\\
&\hspace{0.95in}+\left.\frac{a^2\sin^2\theta \alpha f}{q^2}\left(\pd_r\psi+\frac{3r}{r^2+a^2}\psi\right)^2 -\frac{a^23f}{r^2+a^2}\frac{\psi^2}{q^2}-\frac1{A_1^2q^2}\pd_r(A_1^23rf\psi^2)\right] \\
&+\int_{t_1}^{t_2}\int_{\tilde{\Sigma}_t}-\left(2\alpha^{-1}fL\psi+\frac{6rf}{r^2+a^2}\psi\right)\Box_{\tilde{g}}\psi.
\end{align*}
\end{lemma}

\begin{proof}
The proof is similar to the proof of Lemma \ref{k:p_ee_identity_lem}, however there are a few subtle differences, some of which actually simplify the proof. Again, we will use Proposition \ref{general_divergence_estimate_prop} together with the following current template.
$$J[X,w,m]_\mu = T_{\mu\nu} X^\nu +w\psi\pd_\mu\psi-\frac12\psi^2\pd_\mu w+m_\mu \psi^2,$$
$$T_{\mu\nu}=2\pd_\mu\psi\pd_\nu\psi-g_{\mu\nu}\pd^\lambda\psi\pd_\lambda\psi.$$

Assume for now that $\Box_{\tilde{g}}\psi=0$. Let $\alpha = \frac{\Delta}{r^2+a^2}$, and observe that
$$L=\alpha\pd_r+\pd_t,$$
$$q^2g^{rr}=(r^2+a^2)\alpha,$$
$$q^2g^{tt}=-(r^2+a^2)\alpha^{-1}.$$

\begin{lemma}\label{divJpsiX_lem}
Without appealing directly to the particular expression for $\alpha$, one can deduce the following.
\begin{multline*}
\frac{q^2}{r^2+a^2}divJ[\alpha^{-1}fL]=(\alpha^{-1}f)'(L\psi)^2-\frac{6rf}{r^2+a^2}\left(\alpha(\pd_r\psi)^2-\alpha^{-1}(\pd_t\psi)^2\right) \\
-\left(\frac{4rf}{r^2+a^2}+f'\right)\frac{Q^{\alpha\beta}}{r^2+a^2}\pd_\alpha\psi\pd_\beta\psi.
\end{multline*}
\end{lemma}
\begin{proof}
Note that
$$div J[X] = K^{\mu\nu}\pd_\mu\psi\pd_\nu\psi,$$
where
$$K^{\mu\nu}=2g^{\mu\lambda}\pd_\lambda X^\nu-X^\lambda \pd_\lambda(g^{\mu\nu})-div X g^{\mu\nu}.$$

Set $X=\alpha^{-1}f(\alpha\pd_r+\pd_t)=f\pd_r+\alpha^{-1}f\pd_t$. From the above formula, since $g^{rt}=0$,
$$\frac{q^2}{r^2+a^2}(K^{tr}+K^{rt})=2\frac{q^2}{r^2+a^2}g^{rr}\pd_r X^t=2\alpha\pd_r(\alpha^{-1} f).$$
Thus, the expression for $\frac{q^2}{r^2+a^2}divJ[\alpha^{-1}fL]$ will have a mixed term of the form 
$$2\alpha\pd_r(\alpha^{-1}f)\pd_r\psi\pd_t\psi.$$
Note that
\begin{align*}
(\alpha^{-1} f)'(L\psi)^2 &= (\alpha^{-1} f)'(\alpha\pd_r\psi+\pd_t\psi)^2 \\
&= \alpha^2(\alpha^{-1}f)'(\pd_r\psi)^2+2\alpha(\alpha^{-1}f)'\pd_r\psi\pd_t\psi +(\alpha^{-1} f)'(\pd_t\psi)^2.
\end{align*}
We now compute the $(\pd_r\psi)^2$ and $(\pd_t\psi)^2$ components, subtracting the part that will be grouped with the $(L\psi)^2$ term.
\begin{align*}
\frac{q^2}{r^2+a^2}K^{rr}-\alpha^2(\alpha^{-1}f)' &= \frac{q^2}{r^2+a^2}\left[2g^{rr}\pd_rX^r-X^r\pd_r g^{rr}-\frac1{A_1^2q^2}\pd_r(A_1^2q^2X^r)g^{rr}\right] \\
&\hspace{3.2in}-\alpha^2(\alpha^{-1}f)' \\
&= \frac{q^2}{r^2+a^2}\left[2g^{rr}\pd_rX^r-\frac1{A_1^2q^2}\pd_r(A_1^2q^2g^{rr}X^r)\right]-\alpha^2(\alpha^{-1}f)' \\
&= 2\alpha\pd_r f -\frac{1}{(r^2+a^2)^3}\pd_r\left((r^2+a^2)^3\alpha f\right)-\alpha^2(\alpha^{-1}f)' \\
&= -\frac{6r \alpha f}{r^2+a^2}
\end{align*}
and
\begin{align*}
\frac{q^2}{r^2+a^2}K^{tt}-(\alpha^{-1}f)' &= \frac{q^2}{r^2+a^2}\left[-X^r\pd_r g^{tt}-\frac{1}{A_1^2q^2}\pd_r(A_1^2q^2X^r)g^{tt}\right]-(\alpha^{-1}f)' \\
&= -\frac{q^2}{r^2+a^2}\frac{1}{A_1^2q^2}\pd_r(A_1^2q^2 g^{tt} X^r) -(\alpha^{-1}f)' \\
&= -\frac{1}{(r^2+a^2)^3}\pd_r\left((r^2+a^2)^3(-\alpha^{-1}) f\right)-(\alpha^{-1}f)' \\
&= \frac{6r\alpha^{-1}f}{r^2+a^2}.
\end{align*}
Finally,
\begin{align*}
\frac{q^2}{r^2+a^2}{}^{(Q)}K^{\alpha\beta} &= \frac{q^2}{r^2+a^2}\left[-X^r\pd_r {}^{(Q)}g^{\alpha\beta}-\frac1{A_1^2q^2}\pd_r\left(A_1^2q^2X^r\right){}^{(Q)}g^{\alpha\beta}\right] \\
&= \frac{q^2}{r^2+a^2}\left[-\frac{1}{A_1^2q^2}\pd_r\left(A_1^2q^2{}^{(Q)}g^{\alpha\beta}X^r\right)\right] \\
&= -\frac{1}{(r^2+a^2)^3}\pd_r((r^2+a^2)^2Q^{\alpha\beta}f) \\
&= -\left(\frac{4rf}{r^2+a^2}+f'\right) \frac{Q^{\alpha\beta}}{r^2+a^2}.
\end{align*}
Combining all these terms gives the identity stated in the lemma. 
\end{proof}

Next, we choose $w=\frac{6rf}{r^2+a^2}$ to directly cancel with the middle term in the above lemma.
\begin{lemma}\label{divJpsiXw_lem}
\begin{multline*}
\frac{q^2}{r^2+a^2}divJ\left[\alpha^{-1}fL,\frac{6rf}{r^2+a^2}\right] \\
= (\alpha^{-1}f)'(L\psi)^2+\left(\frac{2rf}{r^2+a^2}-f'\right)\frac{Q^{\alpha\beta}}{r^2+a^2}\pd_\alpha\psi\pd_\beta\psi-\frac12\frac{q^2}{r^2+a^2}\Box_{\tilde{g}}\left(\frac{6rf}{r^2+a^2}\right)\psi^2.
\end{multline*}
\end{lemma}
\begin{proof}
Note that
$$divJ[0,w]=wg^{\mu\nu}\pd_\mu\psi\pd_\nu\psi-\frac12\Box_{\tilde{g}}w \psi^2.$$
We compute the new terms only.
\begin{multline*}
\frac{q^2}{r^2+a^2}divJ\left[0,\frac{6rf}{r^2+a^2}\right] = \frac{6rf}{r^2+a^2}\frac{q^2 g^{\alpha\beta}}{r^2+a^2}\pd_\alpha\psi\pd_\beta\psi -\frac12\frac{q^2}{r^2+a^2}\Box_{\tilde{g}}\left(\frac{6rf}{r^2+a^2}\right)\psi^2 \\
= \frac{6rf}{r^2+a^2}\left(\alpha (\pd_r\psi)^2-\alpha^{-1}(\pd_t\psi)^2\right) +\frac{6rf}{r^2+a^2}\frac{Q^{\alpha\beta}}{r^2+a^2}\pd_\alpha\psi\pd_\beta\psi -\frac12\frac{q^2}{r^2+a^2}\Box_{\tilde{g}}\left(\frac{6rf}{r^2+a^2}\right)\psi^2.
\end{multline*}
When adding these terms to the expression in Lemma \ref{divJpsiX_lem}, the $\alpha(\pd_r\psi)^2-\alpha^{-1}(\pd_t\psi)^2$ terms cancel (this was the reason for the choice of $w=\frac{6rf}{r^2+a^2}$) and the result is as desired. 
\end{proof}

The term $-\frac12\frac{q^2}{r^2+a^2}\Box_{\tilde{g}}\left(\frac{6rf}{r^2+a^2}\right)\psi^2$ is like $-6r^{-1}f''-24r^{-1}\pd_r(r^{-1}f)\phi^2$. In the future, when $f\sim r^p$, this will have a sign $-p^2-3p+4=-(p+4)(p-1)$. The sign will be negative if $p>1$, which is bad. So we include a divergence term to fix it. (Unlike in the analogous step for the proof of Lemma \ref{k:p_ee_identity_lem}, there will not be a need for a smallness parameter $\epsilon$.) This is the point of the following Lemma.
\begin{lemma}
\begin{multline*}
\alpha^{-1}f'(L\psi)^2+\frac{q^2}{r^2+a^2}\left[-\frac12\Box_{\tilde{g}}\left(\frac{6rf}{r^2+a^2}\right)\psi^2+div\left(\psi^2\frac{3rf'}{q^2}L\right)\right] \\
= \alpha f'\left(\alpha^{-1}L\psi+\frac{3r}{r^2+a^2}\psi\right)^2+\frac{6r\alpha(2f-rf')}{(r^2+a^2)^2}\psi^2 \\
-\alpha'\pd_r\left(\frac{3r}{r^2+a^2}\right)f\psi^2 +\frac{a^2}{r^2+a^2}\left(\frac{-3\alpha f'}{r^2+a^2}+\frac{-36\alpha r f}{(r^2+a^2)^2}\right)\psi^2.
\end{multline*}
\end{lemma}
\begin{proof}
First, borrowing from a calculation in the proof of Lemma \ref{k:rp_add_term_lem}, we obtain
\begin{align*}
-\frac{q^2}{r^2+a^2}\frac12\Box_{\tilde{g}}\left(\frac{6rf}{r^2+a^2}\right)\psi^2 
&= -\frac{1}{A_1^2(r^2+a^2)}\pd_r\left(A_1^2(r^2+a^2)\alpha\pd_r\left(\frac{3rf}{r^2+a^2}\right)\right)\psi^2 \\
&= -\frac{q^2}{r^2+a^2}\Box_g\left(\frac{3rf}{r^2+a^2}\right)\psi^2-\frac{\pd_rA_1^2}{A_1^2}\alpha\pd_r\left(\frac{3rf}{r^2+a^2}\right)\psi^2 \\
&= -\frac{3\alpha rf''}{r^2+a^2}\psi^2-\frac{\pd_rA_1^2}{A_1^2}\alpha\pd_r\left(\frac{3rf}{r^2+a^2}\right)\psi^2 \\
&\hspace{.2in}-\alpha'\pd_r\left(\frac{3rf}{r^2+a^2}\right)\psi^2-\frac{6\alpha a^2}{(r^2+a^2)^2}\left(f'+\frac{2r}{r^2+a^2}f\right)\psi^2.
\end{align*}
We also calculate
\begin{multline*}
\frac{q^2}{r^2+a^2}div\left(\psi^2\frac{3r}{q^2}f'L\right) = \frac{1}{A_1^2(r^2+a^2)}\pd_\alpha(A_1^2\psi^2 3rf'L^\alpha) \\
=\frac{3rf'}{r^2+a^2}2\psi L\psi + \frac{\pd_r(3r f'\alpha)}{r^2+a^2}\psi^2+\frac{\pd_rA_1^2}{A_1^2}\frac{3r\alpha f'}{r^2+a^2}\psi^2 \\
= \frac{3rf'}{r^2+a^2}2\psi L\psi +\frac{3\alpha f'}{r^2+a^2}\psi^2 + \frac{3\alpha rf''}{r^2+a^2}\psi^2 +\frac{\pd_rA_1^2}{A_1^2}\frac{3r\alpha f'}{r^2+a^2}\psi^2 +\alpha'\frac{3 r f'}{r^2+a^2}\psi^2  \\
= \frac{3rf'}{r^2+a^2}2\psi L\psi +\frac{9\alpha f'}{r^2+a^2}\psi^2 + \frac{3\alpha rf''}{r^2+a^2}\psi^2 -\frac{6\alpha f'}{r^2+a^2}\psi^2+\frac{\pd_rA_1^2}{A_1^2}\frac{3r\alpha f'}{r^2+a^2}\psi^2 +\alpha'\frac{3 r f'}{r^2+a^2}\psi^2.
\end{multline*}
The first two terms in the last line almost complete a square (up to a term on the order of $\frac{a^2}{r^2+a^2}$) with the term $\alpha^{-1}f'(L\psi)^2$. The third term cancels with the first term from the previous calculation. Due to the fourth and fifth terms, which did not show in the calculation in the proof of Lemma \ref{k:rp_add_term_lem}, the $\epsilon$ parameter is not needed here, allowing for a slightly simpler calculation.
$$\alpha^{-1}f'(L\psi)^2+\frac{3rf'}{r^2+a^2}2\psi L\psi +\frac{9\alpha f'}{r^2+a^2}\psi^2=\alpha f'\left(\alpha^{-1}L\psi+\frac{3r}{r^2+a^2}\psi\right)^2+\frac{a^2 9\alpha f'}{(r^2+a^2)^2}\psi^2$$
and
$$-\frac{3\alpha rf''}{r^2+a^2}\psi^2+\frac{3\alpha rf''}{r^2+a^2}\psi^2=0.$$
The new terms are
\begin{align*}
-\frac{\pd_rA_1^2}{A_1^2}&\alpha\pd_r\left(\frac{3rf}{r^2+a^2}\right)\psi^2 -\frac{6\alpha f'}{r^2+a^2}\psi^2+\frac{\pd_rA_1^2}{A_1^2}\frac{3r\alpha f'}{r^2+a^2}\psi^2 \\
&= -\frac{6\alpha f'}{r^2+a^2}\psi^2 -\frac{\pd_r A_1^2}{A_1^2} \pd_r\left(\frac{3r}{r^2+a^2}\right)\alpha f \psi^2 \\
&= -\frac{6\alpha f'}{r^2+a^2}\psi^2 -\left(\frac{4r}{r^2+a^2}\right)\left(\frac{3}{r^2+a^2}-\frac{6r^2}{(r^2+a^2)^2}\right)\alpha f \psi^2 \\
&= -\frac{6r^2\alpha f'}{(r^2+a^2)^2}\psi^2 -\frac{6a^2\alpha f'}{(r^2+a^2)^2}-\left(\frac{4r}{r^2+a^2}\right)\left(-\frac{3}{r^2+a^2}+\frac{6a^2}{(r^2+a^2)^2}\right)\alpha f\psi^2 \\
&= \frac{-6r^2\alpha f'+12r\alpha f}{(r^2+a^2)^2}\psi^2 +\frac{a^2}{r^2+a^2}\left(\frac{-6(r^2+a^2)\alpha f'-24r\alpha f}{(r^2+a^2)^2}\right)\psi^2
\end{align*}
Adding these terms together and ignoring terms with an $a^2$ factor yields
$$\alpha f'\left(\alpha^{-1}L\psi+\frac{3r}{r^2+a^2}\psi\right)^2+\frac{6r\alpha(2f-rf')}{(r^2+a^2)^2}\psi^2.$$
All the remaining terms (which either contain a factor of $\alpha'\sim \frac{M}{r^2}$ or $\frac{a^2}{r^2+a^2}$) are
\begin{align*}
-\alpha'&\pd_r\left(\frac{3rf}{r^2+a^2}\right)\psi^2-\frac{6\alpha a^2}{(r^2+a^2)^2}\left(f'+\frac{2r}{r^2+a^2}f\right)\psi^2 +\alpha'\frac{3 r f'}{r^2+a^2}\psi^2 +\frac{a^2 9\alpha f'}{(r^2+a^2)^2}\psi^2 \\
&\hspace{2.7in}+\frac{a^2}{r^2+a^2}\left(\frac{-6(r^2+a^2)\alpha f'-24r\alpha f}{(r^2+a^2)^2}\right)\psi^2 \\
&= -\alpha'\pd_r\left(\frac{3r}{r^2+a^2}\right)f\psi^2 +\frac{a^2}{r^2+a^2}\left(\frac{(-6+9-6)\alpha f'}{r^2+a^2}+\frac{(-12-24)\alpha r f}{(r^2+a^2)^2}\right)\psi^2 \\
&= -\alpha'\pd_r\left(\frac{3r}{r^2+a^2}\right)f\psi^2 +\frac{a^2}{r^2+a^2}\left(\frac{-3\alpha f'}{r^2+a^2}+\frac{-36\alpha r f}{(r^2+a^2)^2}\right)\psi^2.
\end{align*}
Adding both of these yields the result. 
\end{proof}

Thus, we have shown that if $\Box_{\tilde{g}}\psi = 0$, then
\begin{multline*}
\frac{q^2}{r^2+a^2}div J\left[\alpha^{-1}fL,\frac{6rf}{r^2+a^2},\frac{3rf'}{q^2}L\right] \\
= \alpha f'\left(\alpha^{-1}L\psi+\frac{3r}{r^2+a^2}\psi\right)^2+\frac{6r\alpha(2f-rf')}{(r^2+a^2)^2}\psi^2 +\left(\frac{2rf}{r^2+a^2}-f'\right)\frac{Q^{\alpha\beta}}{r^2+a^2}\pd_\alpha\psi\pd_\beta\psi \\
-\alpha'\alpha^{-2}f(L\psi)^2-\alpha'\pd_r\left(\frac{3r}{r^2+a^2}\right)f\psi^2 +\frac{a^2}{r^2+a^2}\left(\frac{-3\alpha f'}{r^2+a^2}+\frac{-36\alpha r f}{(r^2+a^2)^2}\right)\psi^2.
\end{multline*}
If we remove the assumption that $\Box_{\tilde{g}}\psi=0$, there is an additional term
$$(2X(\psi)+w\psi)\Box_{\tilde{g}}\psi =\left(2\alpha^{-1}fL\psi+\frac{6rf}{r^2+a^2}\psi\right)\Box_{\tilde{g}}\psi$$
appearing in the expression for $div J$.

Finally, we turn to the boundary terms. Since we have assumed that $f$ is supported away from the event horizon, it suffices to compute $-J^t$.
\begin{lemma}
\begin{multline*}
-J^t\left[\alpha^{-1}fL,\frac{6rf}{r^2+a^2},\frac{3rf'}{q^2}L\right] \\
=\left(1-\alpha\frac{a^2\sin^2\theta}{r^2+a^2}\right)\frac{r^2+a^2}{q^2}f\left(\alpha^{-1}L\psi+\frac{3r}{r^2+a^2}\psi\right)^2 +\frac{\alpha^{-1}f}{q^2}(\pd_\theta\psi)^2 +\frac{6f}{q^2}\psi^2\\
+\alpha\frac{a^2\sin^2\theta}{q^2}f\left(\pd_r\psi+\frac{3r}{r^2+a^2}\psi\right)^2 
-\frac{3a^2f}{q^2(r^2+a^2)}\psi^2-\frac1{A_1^2q^2}\pd_r(A_1^23rf\psi^2).
\end{multline*}
\end{lemma}
\begin{proof}
Borrowing a calculation from the proof of Lemma \ref{k:rp_boundary_terms_lem}, we have
$$-J^t[\alpha^{-1}fL] = \frac{r^2+a^2}{q^2}\left(1-\alpha\frac{a^2\sin^2\theta}{r^2+a^2}\right)\alpha^{-2}f(L\psi)^2+\alpha\frac{a^2\sin^2\theta}{q^2}f(\pd_r\psi)^2+\frac{\alpha^{-1}f}{q^2}(\pd_\theta\psi)^2.$$
Borrowing another calculation from the proof of Lemma \ref{k:rp_boundary_terms_lem}, we have
\begin{align*}
-J^t &\left[0,\frac{6rf}{r^2+a^2}\right] \\
&= \left(1-\alpha\frac{a^2\sin^2\theta}{r^2+a^2}\right)\left(\frac{6r\alpha^{-1}f}{q^2}\psi L\psi+\frac{3f}{q^2}\psi^2\right)+\alpha\frac{a^2\sin^2\theta}{r^2+a^2}\left(\frac{6rf}{q^2}\psi\pd_r\psi+\frac{3f}{q^2}\psi^2\right) \\
&\hspace{4in}+\frac{3rf'}{q^2}\psi^2-\frac{3}{q^2}\pd_r(rf\psi^2) \\
&= \left(1-\alpha\frac{a^2\sin^2\theta}{r^2+a^2}\right)\left(\frac{6r\alpha^{-1}f}{q^2}\psi L\psi+\frac{9f}{q^2}\psi^2\right)+\alpha\frac{a^2\sin^2\theta}{r^2+a^2}\left(\frac{6rf}{q^2}\psi\pd_r\psi+\frac{9f}{q^2}\psi^2\right) \\
&\hspace{2in}-\frac{6f}{q^2}\psi^2+\frac{3rf'}{q^2}\psi^2+\frac{\pd_rA_1^2}{A_1^2}\frac{3rf}{q^2}\psi^2-\frac{1}{A_1^2q^2}\pd_r(A_1^23rf\psi^2).
\end{align*}
Following a similar procedure as in the proof of Lemma \ref{k:rp_boundary_terms_lem}, we notice that
\begin{multline*}
\frac{r^2+a^2}{q^2}\alpha^{-2}f(L\psi)^2+\frac{6r\alpha^{-1}f}{q^2}\psi L\psi+\frac{9f}{q^2}\psi^2 \\
=\frac{r^2+a^2}{q^2}f\left(\alpha^{-1}L\psi+\frac{3r}{r^2+a^2}\psi\right)^2+\frac{9a^2f}{q^2(r^2+a^2)}\psi^2
\end{multline*}
and
$$\frac{r^2+a^2}{q^2}f(\pd_r\psi)^2+\frac{6rf}{q^2}\psi\pd_r\psi+\frac{9f}{q^2}\psi^2 = \frac{r^2+a^2}{q^2}f\left(\pd_r\psi+\frac{3r}{r^2+a^2}\psi\right)^2+\frac{9a^2f}{q^2(r^2+a^2)}\psi^2.$$
Also, there are two new terms, which we now combine.
$$-\frac{6f}{q^2}\psi^2+\frac{\pd_rA_1^2}{A_1^2}\frac{3rf}{q^2}\psi^2 = -\frac{6f}{q^2}\psi^2 +\frac{12r^2f}{q^2(r^2+a^2)}\psi^2 
= \frac{6f}{q^2}\psi^2-\frac{12a^2f}{q^2(r^2+a^2)}\psi^2.$$
Thus,
\begin{multline*}
-J^t\left[\alpha^{-1}fL,\frac{6rf}{r^2+a^2}\right] \\
=\left(1-\alpha\frac{a^2\sin^2\theta}{r^2+a^2}\right)\frac{r^2+a^2}{q^2}f\left(\alpha^{-1}L\psi+\frac{3r}{r^2+a^2}\psi\right)^2 \\
+\alpha\frac{a^2\sin^2\theta}{r^2+a^2}\frac{r^2+a^2}{q^2}f\left(\pd_r\psi+\frac{3r}{r^2+a^2}\psi\right)^2 \\
-\frac{3a^2f}{q^2(r^2+a^2)}\psi^2+\frac{\alpha^{-1}f}{q^2}(\pd_\theta\psi)^2+\frac{6f}{q^2}\psi^2+\frac{3rf'}{q^2}\psi^2-\frac1{A_1^2q^2}\pd_r(A_1^23rf\psi^2).
\end{multline*}
Finally,
$$-J^t\left[0,0,\frac{3rf'}{q^2}L\right]=-\frac{3rf'}{q^2}\psi^2L^t=-\frac{3rf'}{q^2}\psi^2.$$
Adding these two expressions together yields the result. 
\end{proof}

This concludes the proof. 
\end{proof}

\subsection{The incomplete $r^p$ estimate near $i^0$}\label{p_identity_combined_sec}
Now we combine the identities from Lemmas \ref{p_identity_phi_lem} and \ref{p_identity_psi_lem} and make a choice for the function $f$ (so that $f=r^p$ for large $r$) to prove the following proposition.
\begin{proposition}\label{incomplete_p_estimates_prop}
Fix $\delm,\delp>0$. Let $R$ be a sufficiently large radius. Then for all $p\in[\delm,2-\delp]$, the following estimate holds if $\phi$ and $\psi$ decay sufficiently fast as $r\rightarrow\infty$.
\begin{align*}
\int_{\Sigma_{t_2}\cap\{r>2R\}}&r^p\left[(L\phi)^2+|\sla\nabla\phi|^2+r^{-2}\phi^2\right]
+\int_{\tilde{\Sigma}_{t_2}\cap\{r>2R\}}r^p\left[(L\psi)^2+|\sla\nabla\psi|^2+r^{-2}\psi^2\right] \\
&+ \int_{t_1}^{t_2}\int_{\Sigma_t\cap\{r>2R\}}r^{p-1}\left[(L\phi)^2+|\sla\nabla\phi|^2+r^{-2}\phi^2\right] \\
&+ \int_{t_1}^{t_2}\int_{\tilde{\Sigma}_t\cap\{r>2R\}}r^{p-1}\left[(L\psi)^2+|\sla\nabla\psi|^2+r^{-2}\psi^2\right] \\
\lesssim \int_{\Sigma_{t_2}\cap\{r>2R\}}&r^p\left[(L\phi)^2+|\sla\nabla\phi|^2+r^{-2}\phi^2\right]
+\int_{\tilde{\Sigma}_{t_2}\cap\{r>2R\}}r^p\left[(L\psi)^2+|\sla\nabla\psi|^2+r^{-2}\psi^2\right] \\
&\hspace{4.3in}+Err,
\end{align*}
where
\begin{align*}
Err &= Err_1 + Err_2 + Err_\Box \\
Err_1 &= \int_{t_1}^{t_2}\int_{\Sigma_t\cap\{R<r<2R\}}(L\phi)^2+|\sla\nabla\phi|^2+\frac{a^2}{M^2}(\pd_t\phi)^2+\phi^2  \\
&\hspace{1in}+ \int_{t_1}^{t_2}\int_{\tilde{\Sigma}_t\cap\{R<r<2R\}}(L\psi)^2+|\sla\nabla\psi|^2+\frac{a^2}{M^2}(\pd_t\phi)^2+\psi^2 \\
Err_2 &= \int_{\Sigma_{t_1}\cap\{r>R\}}a^2r^{p-2}(\pd_r\phi)^2 + \int_{\tilde{\Sigma}_{t_1}\cap\{r>R\}}a^2r^{p-2}(\pd_r\psi)^2 \\
Err_\Box &= \int_{t_1}^{t_2}\int_{\Sigma_t\cap\{R<r\}}r^p(|L\phi|+r^{-1}|\phi|)|\Box_g\phi| + \int_{t_1}^{t_2}\int_{\tilde{\Sigma}_t\cap\{R<r\}}r^p(|L\psi|+r^{-1}|\psi|)|\Box_{\tilde{g}}\psi|.
\end{align*}
\end{proposition}
\begin{remark}
Unlike most estimates, this estimate could have been separated into two valid estimates--one estimate depending only on $\phi$ and the other depending only on $\psi$. The reason is that the linear terms $\mathcal{L}_\phi$ and $\mathcal{L}_\psi$, which cause the mixing, were ignored. (See $Err_l$ in Proposition \ref{p_estimates_prop}.)
\end{remark}
\begin{proof}
The estimate follows from the identities given in Lemmas \ref{p_identity_phi_lem} and \ref{p_identity_psi_lem}, and a particular choice for the function $f$.
$$f(r)=\rho^p,$$
where
$$
\rho = \left\{
\begin{array}{ll}
0 & r\le R \\
smooth,\rho'\ge 0 & r\in[R,2R] \\
r & 2R< r.
\end{array}
\right.
$$
With this choice, we have 
$$f\ge 0$$
$$f'\ge 0$$
and for $r>2R$,
$$f = r^p$$
$$f'=pr^{p-1}.$$
Furthermore, for $r>2R$,
$$\frac{2rf}{r^2+a^2}-f' = \frac{2r^{p+1}}{r^2+a^2}-pr^{p-1} = \frac{(2-p)r^{p+1}}{r^2+a^2}-\frac{a^2pr^{p-1}}{r^2+a^2} \ge \frac{1}{1+a^2/(4R^2)}\left(2-p-\frac{a^2p}{4R^2}\right)r^{p-1}.$$
It follows that if $p\le 2-\delp$ and $R$ is sufficiently large so that $\frac{a^2p}{4R^2}\le \delp/2$, then for $r>2R$,
$$r^{p-1}\lesssim \frac{2rf}{r^2+a^2}-f'.$$
Also, for $r>2R$,
$$\epsilon\alpha ((1-\epsilon)f'-rf'') = \epsilon\alpha ((1-\epsilon)pr^{p-1}-p(p-1)r^{p-1}) = \epsilon\alpha p (2-\epsilon -p)r^{p-1}.$$
If $R$ is sufficiently large so that $\alpha>3/4$ and $p\le 2-\delp$ and $\epsilon\le \delp/2$, then
$$\epsilon r^{p-1} \lesssim \epsilon\alpha ((1-\epsilon)f'-rf'').$$
Also, for $r>2R$, if $p\le 2-\delp$, then
$$\frac{6r\alpha(2f-rf')}{r^2+a^2} = \frac{6r\alpha (2-p)r^p}{r^2+a^2} \sim (2-p)r^{p-1}.$$

We also note that there are some error terms that either have a factor of $\alpha'$ or $\frac{a^2}{r^2+a^2}$. Each of these terms has a smallness parameter available, since $R$ can be taken to be very large and
$$\alpha'\lesssim \frac{M}{R}r^{-1}$$
and
$$\frac{a^2}{r^2+a^2}\lesssim \frac{M^2}{R^2}.$$

Finally, we observe that if $\phi$ and $\psi$ vanish sufficiently fast as $r\rightarrow\infty$, then since $f$ is supported for $r>R$, we have
$$\int_{\Sigma_t}-\frac1{q^2}\pd_r(rf\phi^2) =0$$
and
$$\int_{\tilde{\Sigma}_t}-\frac{1}{A_1^2q^2}\pd_r(A_1^23rf\psi^2) = 0.$$

With these facts having been established, it is straightforward to check that the estimate follows from Lemmas  \ref{p_identity_phi_lem} and \ref{p_identity_psi_lem}. 
\end{proof}

\subsection{The $r^p$ estimate}\label{phi_psi_p_ee_sec}

We conclude this section by proving the $r^p$ estimate. This is a combination of the $h\pd_t$ estimate (Proposition \ref{translated_h_dt_prop}), the Morawetz estimate (Proposition \ref{translated_morawetz_prop}), and the incomplete $r^p$ estimate near $i^0$ (Proposition \ref{incomplete_p_estimates_prop}).
\begin{proposition}\label{p_estimates_prop}
Suppose $|a|/M$ is sufficiently small. Fix $\delm,\delp>0$ and let $p\in[\delm,2-\delp]$. Then the following estimate holds if $(\phi,\psi)$ decay sufficiently fast as $r\rightarrow\infty$.
\begin{align*}
&\int_{\Sigma_{t_2}}r^p\left[(L\phi)^2+|\sla\nabla\phi|^2+r^{-2}\phi^2 + r^{-2}(\pd_r\phi)^2\right] \\
&\hspace{2in}+\int_{\tilde\Sigma_{t_2}}r^p\left[(L\psi)^2+|\tsla\nabla\psi|^2+r^{-2}\psi^2+r^{-2}(\pd_r\psi)^2\right] \\
&\hspace{.8in}+ \int_{t_1}^{t_2}\int_{\Sigma_t}r^{p-1}\left[\chi_{trap}(L\phi)^2+\chi_{trap}|\sla\nabla\phi|^2+r^{-2}\phi^2+r^{-2}(\pd_r\phi)^2\right] \\
&\hspace{.8in}+ \int_{t_1}^{t_2}\int_{\tilde\Sigma_t}r^{p-1}\left[\chi_{trap}(L\psi)^2+\chi_{trap}|\tsla\nabla\psi|^2+r^{-2}\psi^2 +r^{-2}(\pd_r\psi)^2\right] \\
&\lesssim \int_{\Sigma_{t_1}}r^p\left[(L\phi)^2+|\sla\nabla\phi|^2+r^{-2}\phi^2+r^{-2}(\pd_r\phi)^2\right] \\
&\hspace{2in}+\int_{\tilde\Sigma_{t_1}}r^p\left[(L\psi)^2+|\tsla\nabla\psi|^2+r^{-2}\psi^2+r^{-2}(\pd_r\psi)^2\right] 
+Err,
\end{align*}
where $\chi_{trap}=\left(1-\frac{r_{trap}}{r}\right)^2$ and
\begin{align*}
Err &= Err_l+Err_{nl} \\
Err_l &= \int_{t_1}^{t_2}\int_{\Sigma_t\cap\{R<r\}}r^p(|L\phi|+r^{-1}|\phi|)|\mathcal{L}_\phi| + \int_{t_1}^{t_2}\int_{\tilde{\Sigma}_t\cap\{R<r\}}r^p(|L\psi|+r^{-1}|\psi|)|\mathcal{L}_\psi| \\
Err_{nl} &= \int_{t_1}^{t_2}\int_{\Sigma_t}|(2X(\phi)+w\phi+w_{(a)}\psi)(\Box_g\phi-\mathcal{L}_\phi)| \\
&\hspace{2in}+\int_{t_1}^{t_2}\int_{\tilde{\Sigma}_t}|(2X(\psi)+\tilde{w}\psi+\tilde{w}_{(a)}\phi)(\Box_{\tilde{g}}\psi-\mathcal{L}_\psi)|,
\end{align*}
where the vectorfield $X$ and functions $w$, $w_{(a)}$, $\tilde{w}$, and $\tilde{w}_{(a)}$ satisfy the following properties. \\
\bp $X$ is everywhere timelike, but asymptotically null at the rate $X=O(r^p)L+O(r^{p-2})\pd_t$. \\
\bp $X|_{r=r_H}=-\lambda\pd_r$ for some positive constant $\lambda$. \\
\bp $X|_{r=r_{trap}}=\lambda\pd_t$ for some positive constant $\lambda$. \\
\bp $w$ and $\tilde{w}$ are both $O(r^{p-1})$. \\
\bp $w_{(a)}$ and $\tilde{w}_{(a)}$ are the same functions as defined in Proposition \ref{translated_morawetz_prop}. In particular,
$$|w_{(a)}\psi|\lesssim \frac{|a|^3M\sin^2\theta}{r^5}|A \psi|,$$
$$|\tilde{w}_{(a)}\phi|\lesssim \frac{|a|^3M\sin^2\theta}{r^5}|A^{-1}\phi|.$$
\end{proposition}
\begin{proof}
We start with the Morawetz estimate (Proposition \ref{translated_morawetz_prop}) and add a small constant times the incomplete $r^p$ estimates (Proposition \ref{incomplete_p_estimates_prop}). The small constant can be chosen so that the bulk error term $Err_1$ from Proposition \ref{incomplete_p_estimates_prop} can be absorbed into the bulk in the Morawetz estimate. The result is the following estimate.
\begin{align*}
&\int_{\Sigma_{t_2}}r^p\left[(L\phi)^2+|\sla\nabla\phi|^2+r^{-2}\phi^2\right]+\frac{M^2}{r^2}(\pd_r\phi)^2 \\
&\hspace{2.5in}+\int_{\tilde{\Sigma}_{t_2}}r^p\left[(L\psi)^2+|\tsla\nabla\psi|^2+r^{-2}\psi^2\right]+\frac{M^2}{r^2}(\pd_r\psi)^2 \\
&\hspace{1in}+\int_{t_1}^{t_2}\int_{\Sigma_t} r^{p-1}\left[\chi_{trap}(L\phi)^2+\chi_{trap}|\sla\nabla\phi|^2+r^{-2}\phi^2\right]+\frac{M^2}{r^3}(\pd_r\phi)^2 \\
&\hspace{1in}+\int_{t_1}^{t_2}\int_{\tilde{\Sigma}_t} r^{p-1}\left[\chi_{trap}(L\psi)^2+\chi_{trap}|\tsla\nabla\psi|^2+r^{-2}\psi^2\right]+\frac{M^2}{r^3}(\pd_r\psi)^2 \\
&\lesssim \int_{\Sigma_{t_1}}r^p\left[(L\phi)^2+|\sla\nabla\phi|^2+r^{-2}\phi^2\right]+\frac{M^2}{r^2}(\pd_r\phi)^2 \\
&\hspace{2in}+\int_{\tilde{\Sigma}_{t_1}}r^p\left[(L\psi)^2+|\tsla\nabla\psi|^2+r^{-2}\psi^2\right]+\frac{M^2}{r^2}(\pd_r\psi)^2 +Err',
\end{align*}
where
\begin{align*}
Err' &= Err'_1+Err'_2+Err'_3+Err'_\Box+Err'_{nl} \\
Err'_1 &= \int_{\Sigma_{t_2}}r^{-1}|\phi L\phi|+\int_{\tilde{\Sigma}_{t_2}}r^{-1}|\psi L\psi|+r^{-2}\psi^2 \\
Err'_2 &= \int_{H_{t_1}^{t_2}}(\pd_t\phi)^2+\int_{H_{t_1}^{t_2}}(\pd_t\psi)^2 +\int_{\Sigma_{t_2}}\frac{M^2}{r^2}\left[\chi_H(\pd_r\phi)^2+(\pd_t\phi)^2\right] \\
&\hspace{1in} +\int_{\tilde{\Sigma}_{t_2}}\frac{M^2}{r^2}\left[\chi_H(\pd_r\psi)^2+(\pd_t\psi)^2\right] \\ 
Err'_3 &= \int_{\Sigma_{t_1}\cap\{r>R\}}a^2r^{p-2}(\pd_r\phi)^2 + \int_{\tilde{\Sigma}_{t_1}\cap\{r>R\}}a^2r^{p-2}(\pd_r\psi)^2 \\
Err'_\Box &= \int_{t_1}^{t_2}\int_{\Sigma_t\cap\{R<r\}}r^p(|L\phi|+r^{-1}|\phi|)|\Box_g\phi| + \int_{t_1}^{t_2}\int_{\tilde{\Sigma}_t\cap\{R<r\}}r^p(|L\psi|+r^{-1}|\psi|)|\Box_{\tilde{g}}\psi| \\
Err'_{nl} &= \int_{t_1}^{t_2}\int_{\Sigma_t}|(2X'(\phi)+w'\phi+w'_{(a)}\psi)(\Box_g\phi-\mathcal{L}_\phi)| \\
&\hspace{1in} +\int_{t_1}^{t_2}\int_{\tilde{\Sigma}_t}|(2X'(\psi)+\tilde{w}'\psi+\tilde{w}'_{(a)}\phi)(\Box_{\tilde{g}}\psi-\mathcal{L}_\psi)|,
\end{align*}
and $X'$, $w'$, $\tilde{w}'$, $w'_{(a)}$, and $\tilde{w}'_{(a)}$ are the vectorfield and functions defined in the Morawetz estimate (Proposition \ref{translated_morawetz_prop}).

The error term $Err'_1$ can in fact be removed due to the following argument.
\begin{align*}
Err'_1 &\lesssim \int_{\Sigma_{t_2}}\epsilon r^p(L\phi)^2+\epsilon^{-1}r^{-p}r^{-2}\phi^2 + \int_{\tilde\Sigma_{t_2}} \epsilon r^p(L\psi)^2+ (\epsilon^{-1}r^{-p}+1)r^{-2}\psi^2 \\
&\lesssim \int_{\Sigma_{t_2}}\epsilon r^p[(L\phi)^2+r^{-2}\phi^2] + \int_{\tilde\Sigma_{t_2}} \epsilon r^p[(L\psi)^2+r^{-2}\psi^2] \\
&\hspace{2in}+\int_{\Sigma_{t_2}\cap\{r\le R_\epsilon\}}\epsilon^{-1}\phi^2 + \int_{\tilde{\Sigma}_{t_2}\cap\{r\le R_\epsilon \}}\epsilon^{-1}\psi^2.
\end{align*}
The radius $R_\epsilon$ should be chosen sufficiently large so that $\epsilon^{-1}r^{-p}\le \epsilon r^p$ and $\epsilon^{-1}r^{-p}+1\le \epsilon r^p$ whenever $r>R_\epsilon$. This critically depends on the fact that $p\ge\delm>0$. Now, the parameter $\epsilon$ can be taken sufficiently small so as to absorb the first two terms into the left side of the main estimate and the last two terms can be included with the term $Err'_2$ after applying a Hardy estimate.

We return to the main estimate. Notice that most terms have improved weights near $i^0$ and a few error terms remain on $H_{t_1}^{t_2}$ and $\Sigma_{t_2}$. The next step is to use the $h\pd_t$ estimate (Proposition \ref{translated_h_dt_prop}) to eliminate these error terms and improve the weights near $i^0$ for the remaining $\pd_r\phi$ and $\pd_r\psi$ terms. The result is the following estimate.
\begin{align*}
&\int_{\Sigma_{t_2}}r^p\left[(L\phi)^2+|\sla\nabla\phi|^2+r^{-2}\phi^2 + r^{-2}(\pd_r\phi)^2\right] \\
&\hspace{2in}+\int_{\tilde\Sigma_{t_2}}r^p\left[(L\psi)^2+|\tsla\nabla\psi|^2+r^{-2}\psi^2+r^{-2}(\pd_r\psi)^2\right] \\
&\hspace{.8in}+ \int_{t_1}^{t_2}\int_{\Sigma_t}r^{p-1}\left[\chi_{trap}(L\phi)^2+\chi_{trap}|\sla\nabla\phi|^2+r^{-2}\phi^2+r^{-2}(\pd_r\phi)^2\right] \\
&\hspace{.8in}+ \int_{t_1}^{t_2}\int_{\tilde\Sigma_t}r^{p-1}\left[\chi_{trap}(L\psi)^2+\chi_{trap}|\tsla\nabla\psi|^2+r^{-2}\psi^2 +r^{-2}(\pd_r\psi)^2\right] \\
&\lesssim \int_{\Sigma_{t_1}}r^p\left[(L\phi)^2+|\sla\nabla\phi|^2+r^{-2}\phi^2+r^{-2}(\pd_r\phi)^2\right] \\
&\hspace{2in}+\int_{\tilde\Sigma_{t_1}}r^p\left[(L\psi)^2+|\tsla\nabla\psi|^2+r^{-2}\psi^2+r^{-2}(\pd_r\psi)^2\right] 
+Err'',
\end{align*}
where
\begin{align*}
Err'' &= Err''_1+Err''_2+Err''_3+Err''_\Box+Err''_{nl} \\
Err''_1 &= \int_{t_1}^{t_2}\int_{\Sigma_t\cap\{5M<r\}} \epsilon r^{-1}((L\phi)^2+r^{-2}\phi^2) + \int_{t_1}^{t_2}\int_{\tilde\Sigma_t\cap\{5M<r\}} \epsilon r^{-1}((L\psi)^2+r^{-2}\psi^2) \\
Err''_2 &= \int_{\Sigma_{t_2}}\frac{|a|}{M}r^{p-4}(\phi^2+A^2\psi^2) \\
Err''_3 &= \int_{t_1}^{t_2}\int_{\Sigma_t\cap\{6M<r\}} c_\epsilon\frac{|a|}{M}r^{p-5}(\phi^2+A^2\psi^2) \\
Err''_\Box &= \int_{t_1}^{t_2}\int_{\Sigma_t\cap\{R<r\}}r^p(|L\phi|+r^{-1}|\phi|)|\Box_g\phi| + \int_{t_1}^{t_2}\int_{\tilde{\Sigma}_t\cap\{R<r\}}r^p(|L\psi|+r^{-1}|\psi|)|\Box_{\tilde{g}}\psi| \\
Err''_{nl} &= \int_{t_1}^{t_2}\int_{\Sigma_t}|(2X'(\phi)+w'\phi+w'_{(a)}\psi)(\Box_g\phi-\mathcal{L}_\phi)| \\
&\hspace{1in} +\int_{t_1}^{t_2}\int_{\tilde{\Sigma}_t}|(2X'(\psi)+\tilde{w}'\psi+\tilde{w}'_{(a)}\phi)(\Box_{\tilde{g}}\psi-\mathcal{L}_\psi)|, \\
&\hspace{1in}+ \int_{t_1}^{t_2}\int_{\Sigma_t}C_\epsilon r^{p-2}|\pd_t\phi(\Box_g\phi-\mathcal{L}_\phi)| +\int_{t_1}^{t_2}\int_{\tilde{\Sigma}_t}C_\epsilon r^{p-2}|\pd_t\psi(\Box_{\tilde{g}}\psi-\mathcal{L}_\psi)|.
\end{align*}
Note that each of the error terms $Err''_1$, $Err''_2$, and $Err''_3$ comes with a smallness parameter. By taking $\epsilon$ and $|a|/M$ sufficiently small, these error terms can be absorbed into the left side. 

We are left only with the error terms $Err''_\Box$ and $Err''_{nl}$, which we now combine by replacing $|\Box_g\phi|$ and $|\Box_{\tilde{g}}\psi|$ in $Err''_\Box$ with $|\mathcal{L}_\phi|+|\Box_g\phi-\mathcal{L}_\phi|$ and $|\mathcal{L}_\psi|+|\Box_{\tilde{g}}\psi-\mathcal{L}_\psi|$ respectively. The $|\mathcal{L}_\phi|$ and $|\mathcal{L}_\psi|$ terms are collected into the linear error term $Err_l$, and the $|\Box_g\phi-\mathcal{L}_\phi|$ and $|\Box_{\tilde{g}}\psi-\mathcal{L}_\psi|$ terms are combined with the terms in $Err''_{nl}$ to form the nonlinear error term $Err_{nl}$.
\begin{multline*}
Err''_\Box + Err''_{nl} \\
\lesssim \int_{t_1}^{t_2}\int_{\Sigma_t\cap\{R<r\}}r^p(|L\phi|+r^{-1}|\phi|)|\mathcal{L}_\phi| + \int_{t_1}^{t_2}\int_{\tilde{\Sigma}_t\cap\{R<r\}}r^p(|L\psi|+r^{-1}|\psi|)|\mathcal{L}_\psi| \\
+\int_{t_1}^{t_2}\int_{\Sigma_t}|(2X(\phi)+w\phi+w_{(a)}\psi)(\Box_g\phi-\mathcal{L}_\phi)| \\
+\int_{t_1}^{t_2}\int_{\tilde{\Sigma}_t}|(2X(\psi)+\tilde{w}\psi+\tilde{w}_{(a)}\phi)(\Box_{\tilde{g}}\psi-\mathcal{L}_\psi)| \\
\lesssim Err_l+Err_{nl},
\end{multline*}
where
\begin{align*}
X &= X'+O(r^{p-2})\pd_t+O(r^p)L \\
w &= w'+O(r^{p-1}) \\
\tilde{w} &= \tilde{w}'+O(r^{p-1}) \\
w_{(a)} &= w'_{(a)} \\
\tilde{w}_{(a)} &= \tilde{w}'_{(a)}.
\end{align*}
This concludes the proof. 
\end{proof}

\section{The dynamic estimates}\label{wk:dynamic_estimates_sec}

In this section, we prove the dynamic estimates, which provide all of the necessary information related to the future dynamics of the system. These estimates are simply restatements of the energy estimate (Proposition \ref{translated_energy_estimate_prop}) and the $r^p$ estimate (Proposition \ref{p_estimates_prop}) applied to either $(\phi,\psi)$ itself or higher order pairs $(\pd_t^s\phi,\pd_t^s\psi)$, $(\phi^s,\psi^s)$, or $(\phi^{s,k},\psi^{s,k})$.

By applying certain operators to the wave map system, the dynamic estimates are generalized to higher derivatives of the quantities $\phi$ and $\psi$. This will generalize, for example, the classic energy norm $E(t)$ to the homogenous classic energy norm $\mathring{E}^s(t)$ (obtained by applying the only true commutators of the linear system, $\pd_t^s$), as well as to the norms $E^s(t)$ and $E^{s,k}(t)$. The other norms generalize the same way.

The various versions of the dynamic estimates are given in Theorems \ref{p_L_thm}, \ref{p_thm}, \ref{p_o_thm}, \ref{p_s_L_thm}, \ref{p_s_thm}, \ref{p_s_k_L_thm}, and \ref{p_s_k_thm}. Theorem \ref{p_L_thm} is an immediate consequence of Propositions \ref{translated_energy_estimate_prop} and \ref{p_estimates_prop}. Theorems \ref{p_s_L_thm} and \ref{p_s_k_L_thm} follow from Theorem \ref{p_L_thm} by applying operators to the wave map system. These theorems are then improved by handling various linear error terms (which mainly arise because the operators do not commute with the entire linear system) and applying the equation. More specifically, Theorem \ref{p_thm} follows from Theorem \ref{p_L_thm}, Theorem \ref{p_s_thm} follows from Theorem \ref{p_s_L_thm}, and Theorem \ref{p_s_k_thm} follows from Theorem \ref{p_s_k_L_thm}. Additionally, Theorem \ref{p_o_thm} follows directly from Theorem \ref{p_thm}, since the operators $\pd_t^s$ do in fact commute with the entire linear system.

\subsection{The dynamic estimates for $(\phi,\psi)$}

We start by combining Propositions \ref{translated_energy_estimate_prop} and \ref{p_estimates_prop} to obtain Theorem \ref{p_L_thm}. But first, we define the following two expressions to simplify the nonlinear term.
\begin{definition} Let $X$, $w$, $w_{(a)}$, $\tilde{w}$, and $\tilde{w}_{(a)}$ be as defined in Proposition \ref{p_estimates_prop}. Then define
\begin{align*}
\mathcal{X}_\phi &= 2X(\phi)+w\phi+w_{(a)}\psi \\
\mathcal{X}_\psi &= 2X(\psi)+\tilde{w}\psi+\tilde{w}_{(a)}\phi.
\end{align*}
\end{definition}
With this definition, we now have the following theorem.
\begin{theorem}\label{p_L_thm}
Suppose $|a|/M$ is sufficiently small and fix $\delm,\delp>0$. The following estimates hold for $p\in [\delm,2-\delp]$.
$$E(t_2)\lesssim E(t_1)+\int_{t_1}^{t_2}N(t)dt,$$
$$E_p(t_2)+\int_{t_1}^{t_2}B_p(t)dt\lesssim E_p(t_1)+\int_{t_1}^{t_2}(L_1)_p(t)+N_p(t)dt,$$
where
\begin{align*}
E(t)=&\int_{\Sigma_t}\chi_H(\pd_r\phi)^2+(\pd_t\phi)^2+|\sla\nabla\phi|^2+r^{-2}\phi^2 \\
&+\int_{\tilde\Sigma_t}\chi_H(\pd_r\psi)^2+(\pd_t\psi)^2+|\tsla\nabla\psi|^2+r^{-2}\psi^2,
\end{align*}
\begin{align*}
E_p(t)=&\int_{\Sigma_t}r^p\left[(L\phi)^2+|\sla\nabla\phi|^2+r^{-2}\phi^2+r^{-2}(\pd_r\phi)^2\right] \\
&+\int_{\tilde\Sigma_t}r^p\left[(L\psi)^2+|\tsla\nabla\psi|^2+r^{-2}\psi^2+r^{-2}(\pd_r\psi)^2\right],
\end{align*}
\begin{align*}
B_p(t)=&\int_{\Sigma_t}r^{p-1}\left[\chi_{trap}(L\phi)^2+\chi_{trap}|\sla\nabla\phi|^2+r^{-2}\phi^2+r^{-2}(\pd_r\phi)^2\right] \\
&+\int_{\tilde\Sigma_t}r^{p-1}\left[\chi_{trap}(L\psi)^2+\chi_{trap}|\tsla\nabla\psi|^2+r^{-2}\psi^2+r^{-2}(\pd_r\psi)^2\right],
\end{align*}
$$(L_1)_p(t)=\int_{\Sigma_t\cap\{r>R\}}r^p(|L\phi|+r^{-1}|\phi|)|\mathcal{L}_\phi|+\int_{\tilde\Sigma_t\cap\{r>R\}}r^p(|L\psi|+r^{-1}|\psi|)|\mathcal{L}_\psi|,$$
$$N(t)=\int_{\Sigma_t}|\pd_t\phi||\Box_g\phi-\mathcal{L}_\phi| +\int_{\tilde{\Sigma}_t}|\pd_t\psi||\Box_{\tilde g}\psi-\mathcal{L}_\psi|,$$
$$N_p(t)=\int_{\Sigma_t}r^p|\mathcal{X}_\phi||\Box_g\phi-\mathcal{L}_\phi| + \int_{\tilde\Sigma_t}r^p|\mathcal{X}_\psi||\Box_{\tilde{g}}\psi-\mathcal{L}_\psi|, $$
where $\chi_H=1-\frac{r_H}r$ and $\chi_{trap}=\left(1-\frac{r_{trap}}r\right)^2$.
\end{theorem}
\begin{proof}
This is a direct application of Propositions \ref{translated_energy_estimate_prop} and \ref{p_estimates_prop}. Note that the linear error term $L_1$ derives from the error term $Err_l$ in Proposition \ref{p_estimates_prop}, which derives from the $r^p$ estimates near $i^0$ that were proved only for the scalar wave equations. This is why the linear error term $(L_1)_p$ is supported far away from the Kerr black hole. 
\end{proof}

We will momentarily improve the prevoius theorem by absorbing the linear error term $(L_1)_p$ into the bulk term $B_p(t)$ on the left side. To do this, we need the following lemma.
\begin{lemma}\label{p_L_a_lem}
If $(L_1)_p(t)$ and $B_p(t)$ are as defined in the previous theorem, then
$$(L_1)_p(t)\lesssim \frac{|a|}{M}B_p(t).$$
\end{lemma}
\begin{proof}
We have
\begin{multline*}
\int_{\Sigma_t\cap\{r>R\}}r^p(|L\phi|+r^{-1}|\phi|)|\mathcal{L}_\phi|  \\
\lesssim \left(\int_{\Sigma_t\cap\{r>R\}}r^{p-1}[(L\phi)^2+r^{-2}\phi^2]\right)^{1/2}\left(\int_{\Sigma_t\cap\{r>R\}}r^{p+1}(\mathcal{L}_\phi)^2\right)^{1/2} \\
\lesssim (B_p(t))^{1/2}\left(\int_{\Sigma_t\cap\{r>R\}}r^{p+1}(\mathcal{L}_\phi)^2\right)^{1/2} \\
\lesssim \frac{|a|}{M}B_p(t)+\left(\frac{|a|}{M}\right)^{-1}\int_{\Sigma_t\cap\{r>R\}}r^{p+1}(\mathcal{L}_\phi)^2.
\end{multline*}
An analogous estimate also shows that
$$\int_{\tilde\Sigma_t\cap\{r>R\}}r^p(|L\psi|+r^{-1}|\psi|)|\mathcal{L}_\psi| \lesssim \frac{|a|}{M}B_p(t)+\left(\frac{|a|}{M}\right)^{-1}\int_{\tilde\Sigma_t\cap\{r>R\}}r^{p+1}(\mathcal{L}_\psi)^2.$$
Therefore, it suffices to establish the following estimate.
\begin{equation}\label{mathcalL_p_estimate}
\int_{\Sigma_t\cap\{r>R\}}r^{p+1}(\mathcal{L}_\phi)^2 + \int_{\tilde\Sigma_t\cap\{r>R\}}r^{p+1}(\mathcal{L}_\psi)^2 \lesssim \frac{a^2}{M^2}B_p(t).
\end{equation}
Let us look at the $\mathcal{L}_\phi$ term first. Recall that
$$\mathcal{L}_\phi=-2\frac{\pd^\alpha B}{A}A\pd_\alpha \psi + 2\frac{\pd^\alpha B\pd_\alpha B}{A^2}\phi-4\frac{\pd^\alpha A\pd_\alpha B}{A^2} A\psi.$$
It follows that
\begin{multline*}
\int_{\Sigma_t\cap\{r>R\}}r^{p+1}(\mathcal{L}_\phi)^2 \\
\lesssim \int_{\Sigma_t\cap\{r>R\}}r^{p-1}\left[\left(\frac{r\pd^\alpha B}{A}\pd_\alpha\psi\right)^2A^2 +\left(\frac{r\pd^\alpha B\pd_\alpha B}{A^2}\phi\right)^2+\left(\frac{r\pd^\alpha A\pd_\alpha B}{A^2}\psi\right)^2A^2\right] \\
\lesssim \int_{\Sigma_t\cap\{r>R\}}r^{p-1}\left(\frac{r\pd^\alpha B\pd_\alpha B}{A^2}\phi\right)^2 +\int_{\tilde\Sigma_t\cap\{r>R\}}r^{p-1}\left[\left(\frac{r\pd^\alpha B}{A}\pd_\alpha\psi\right)^2 +\left(\frac{r\pd^\alpha A\pd_\alpha B}{A^2}\psi\right)^2\right].
\end{multline*}
Using the identities from Lemma \ref{small_a_quantities_lem}, one can easily check each of the following estimates.
\begin{align*}
\left(\frac{r\pd^\alpha B\pd_\alpha B}{A^2}\phi\right)^2 &\lesssim \frac{a^2}{M^2}r^{-2}\phi^2 \\
\left(\frac{r\pd^\alpha B}{A}\pd_\alpha \psi\right)^2 &\lesssim \frac{a^2}{M^2}|\tsla\nabla\psi|^2+\frac{a^2}{M^2}\frac{M^2}{r^2}(\pd_r\psi)^2 \\
\left(\frac{r\pd^\alpha A\pd_\alpha B}{A^2}\psi\right)^2 &\lesssim \frac{a^2}{M^2}r^{-2}\psi^2.
\end{align*}
Thus,
$$\int_{\Sigma_t\cap\{r>R\}}r^{p+1}(\mathcal{L}_\phi)^2 \lesssim \frac{a^2}{M^2} B_p(t).$$
Let us now look at the $\mathcal{L}_\psi$ term. Recall that
$$\mathcal{L}_\psi=-2\frac{\pd^\alpha A_2}{A_2}\pd_\alpha\psi+2A^{-1}\frac{\pd^\alpha B}{A}\pd_\alpha\phi+2\frac{\pd^\alpha B\pd_\alpha B}{A^2}\psi.$$
It follows that
\begin{multline*}
\int_{\tilde\Sigma_t\cap\{r>R\}}r^{p+1}(\mathcal{L}_\psi)^2 \\
\lesssim \int_{\tilde{\Sigma}_t\cap\{r>R\}}r^{p-1}\left[\left(\frac{r\pd^\alpha A_2}{A_2}\pd_\alpha\psi\right)^2+\left(\frac{r\pd^\alpha B}{A}\pd_\alpha \phi\right)^2A^{-2}+\left(\frac{r\pd^\alpha B\pd_\alpha B}{A^2}\psi\right)^2\right] \\
\lesssim \int_{\tilde{\Sigma}_t\cap\{r>R\}}r^{p-1}\left[\left(\frac{r\pd^\alpha A_2}{A_2}\pd_\alpha\psi\right)^2+\left(\frac{r\pd^\alpha B\pd_\alpha B}{A^2}\psi\right)^2\right] +\int_{\Sigma_t\cap\{r>R\}}r^{p-1}\left(\frac{r\pd^\alpha B}{A}\pd_\alpha \phi\right)^2. 
\end{multline*}
Using again the identities from Lemma \ref{small_a_quantities_lem}, one can easily check each of the following estimates.
\begin{align*}
\left(\frac{r\pd^\alpha A_2}{A_2}\pd_\alpha\psi\right)^2 &\lesssim \frac{a^2}{M^2}|\tsla\nabla\psi|^2+\frac{a^2}{M^2}\frac{M^2}{r^2}(\pd_r\psi)^2 \\
\left(\frac{r\pd^\alpha B\pd_\alpha B}{A^2}\psi\right)^2 &\lesssim \frac{a^2}{M^2}r^{-2}\psi^2\\
\left(\frac{r\pd^\alpha B}{A}\pd_\alpha \phi\right)^2 &\lesssim \frac{a^2}{M^2}|\sla\nabla\phi|^2+\frac{a^2}{M^2}\frac{M^2}{r^2}(\pd_r\phi)^2.
\end{align*}
Thus,
$$\int_{\tilde\Sigma_t\cap\{r>R\}}r^{p+1}(\mathcal{L}_\psi)^2 \lesssim \frac{a^2}{M^2} B_p(t).$$
This completes the proof.
\end{proof}

With Lemma \ref{p_L_a_lem}, Theorem \ref{p_L_thm} can be improved so that the linear error term appearing on the right side of the $r^p$ estimate is removed. Also, by assuming the full nonlinear equations, we can simplify the nonlinear error term. This is the purpose of the following theorem.
\begin{theorem}\label{p_thm} (Improved version of Theorem \ref{p_L_thm})
Suppose $|a|/M$ is sufficiently small and fix $\delm,\delp>0$. Suppose furthermore that the pair $(\phi,\psi)$ satisfies the system
\begin{align*}
\Box_g\phi &= \mathcal{L}_\phi+\mathcal{N}_\phi, \\
\Box_{\tilde{g}}\psi &= \mathcal{L}_\psi + \mathcal{N}_\psi.
\end{align*}
Then the following estimates hold for $p\in [\delm,2-\delp]$.
$$E(t_2)\lesssim E(t_1)+\int_{t_1}^{t_2}N(t)dt,$$
$$E_p(t_2)+\int_{t_1}^{t_2}B_p(t)dt\lesssim E_p(t_1)+\int_{t_1}^{t_2}N_p(t)dt,$$
where $E(t)$, $E_p(t)$, and $B_p(t)$ are as defined in Theorem \ref{p_L_thm}, and
$$N(t)=(E(t))^{1/2}\left(||\mathcal{N}_\phi||_{L^2(\Sigma_t)}+||\mathcal{N}_\psi||_{L^2(\tilde{\Sigma}_t)}\right),$$
$$N_p(t)=\int_{\Sigma_t}r^{p+1}(\mathcal{N}_\phi)^2 + \int_{\tilde\Sigma_t}r^{p+1}(\mathcal{N}_\psi)^2 + \int_{\Sigma_t\cap\{r\approx r_{trap}\}}|\pd_t\phi\mathcal{N}_\phi| + \int_{\tilde\Sigma_t\cap\{r\approx r_{trap}\}}|\pd_t\psi\mathcal{N}_\psi|. $$
\end{theorem}

\begin{proof}
By Theorem \ref{p_L_thm}, we have
$$E(t_2)\lesssim E(t_1)+\int_{t_1}^{t_2}N'(t)dt,$$
where
\begin{align*}
N'(t)&=\int_{\Sigma_t}|\pd_t\phi||\Box_g\phi-\mathcal{L}_\phi|+\int_{\tilde{\Sigma}_t}|\pd_t\psi||\Box_{\tilde{g}}\psi-\mathcal{L}_\psi| \\
&=\int_{\Sigma_t}|\pd_t\phi||\mathcal{N}_\phi|+\int_{\tilde{\Sigma}_t}|\pd_t\psi||\mathcal{N}_\psi|.
\end{align*}
Observe that
\begin{align*}
N'(t) &=\int_{\Sigma_t}|\pd_t\phi||\mathcal{N}_\phi|+\int_{\tilde{\Sigma}_t}|\pd_t\psi||\mathcal{N}_\psi| \\
&\lesssim ||\pd_t\phi||_{L^2(\Sigma_t)}||\mathcal{N}_\phi||_{L^2(\Sigma_t)}+||\pd_t\psi||_{L^2(\tilde\Sigma_t)}||\mathcal{N}_\psi||_{L^2(\tilde\Sigma_t)} \\
&\lesssim (E(t))^{1/2}||\mathcal{N}_\phi||_{L^2(\Sigma_t)}+(E(t))^{1/2}||\mathcal{N}_\psi||_{L^2(\tilde\Sigma_t)} \\
&\lesssim N(t).
\end{align*}
This proves the first estimate of the theorem.

By Theorem \ref{p_L_thm} and Lemma \ref{p_L_a_lem}, we have
$$E_p(t_2)+\int_{t_1}^{t_2}B_p(t)dt\lesssim E_p(t_1)+\frac{|a|}{M}\int_{t_1}^{t_2}B_p(t)dt+\int_{t_1}^{t_2}N_p'(t)dt,$$
where
\begin{align*}
N_p'(t)&=\int_{\Sigma_t}r^p|\mathcal{X}_\phi||\Box_g\phi-\mathcal{L}_\phi| + \int_{\tilde\Sigma_t}r^p|\mathcal{X}_\psi||\Box_{\tilde{g}}\psi-\mathcal{L}_\psi| \\
&=\int_{\Sigma_t}r^p|\mathcal{X}_\phi||\mathcal{N}_\phi| + \int_{\tilde\Sigma_t}r^p|\mathcal{X}_\psi||\mathcal{N}_\psi|.
\end{align*}
Observe that
\begin{multline*}
\int_{\Sigma_t}r^p|\mathcal{X}_\phi||\mathcal{N}_\phi| \\
\lesssim \int_{\Sigma_t}r^p(\chi_{trap}|L\phi|+r^{-1}|\phi|+r^{-1}A|\psi|+r^{-2}|\pd_r\phi|)|\mathcal{N}_\phi| +\int_{\Sigma_t\cap\{r\approx r_{trap}\}}|\pd_t\phi\mathcal{N}_\phi| \\
\lesssim (B_p(t))^{1/2}\left(\int_{\Sigma_t}r^{p+1}|\mathcal{N}_\phi|^2\right)^{1/2} +\int_{\Sigma_t\cap\{r\approx r_{trap}\}}|\pd_t\phi||\mathcal{N}_\phi| \\
\lesssim \epsilon B_p(t) +\epsilon^{-1}\int_{\Sigma_t}r^{p+1}|\mathcal{N}_\phi|^2 +\int_{\Sigma_t\cap\{r\approx r_{trap}\}}|\pd_t\phi||\mathcal{N}_\phi| \\
\lesssim \epsilon B_p(t) +\epsilon^{-1}N_p(t).
\end{multline*}
An analogous estimate also shows that
$$\int_{\tilde\Sigma_t}r^p|\mathcal{X}_\psi||\mathcal{N}_\psi|\lesssim \epsilon B_p(t)+\epsilon^{-1}N_p(t).$$
Therefore, we have
$$E_p(t_2)+\int_{t_1}^{t_2}B_p(t)dt\lesssim E_p(t_1)+(|a|/M+\epsilon)\int_{t_1}^{t_2}B_p(t)dt+\epsilon^{-1}\int_{t_1}^{t_2}N_p(t)dt.$$
By taking $|a|/M$ and $\epsilon$ sufficiently small, the bulk term on the right side can be absorbed into the left side. The result is the second estimate of the theorem. 
\end{proof}

\subsection{The dynamic estimates for $(\pd_t^s\phi,\pd_t^s\psi)$}

We now begin to derive higher order estimates analogous to Theorem \ref{p_thm} by commuting with the linear system. The only operator that completely commutes with the linear system is the operator $\pd_t$. Therefore, we immediately have the following homogeneous estimate, which will be important at the highest level of derivatives in the proof of the main theorem.

\begin{theorem}\label{p_o_thm} (Generalization of Theorem \ref{p_thm})
Suppose $|a|/M$ is sufficiently small and fix $\delm,\delp>0$. Suppose furthermore that the pair $(\phi,\psi)$ satisfies the system
\begin{align*}
\Box_g\phi &= \mathcal{L}_\phi+\mathcal{N}_\phi, \\
\Box_{\tilde{g}}\psi &= \mathcal{L}_\psi + \mathcal{N}_\psi.
\end{align*}
Then the following estimates hold for $p\in [\delm,2-\delp]$.
$$\mathring{E}^s(t_2)\lesssim \mathring{E}^s(t_1)+\int_{t_1}^{t_2}\mathring{N}^s(t)dt,$$
$$\mathring{E}_p^s(t_2)+\int_{t_1}^{t_2}\mathring{B}_p^s(t)dt\lesssim \mathring{E}_p^s(t_1)+\int_{t_1}^{t_2}\mathring{N}_p^s(t)dt,$$
where 
$$\mathring{E}^s(t)=\sum_{s'\le s}E[(\pd_t^{s'}\phi,\pd_t^{s'}\psi)](t),$$
$$\mathring{E}_p^s(t)=\sum_{s'\le s}E_p[(\pd_t^{s'}\phi,\pd_t^{s'}\psi)](t),$$
$$\mathring{B}_p^s(t)=\sum_{s'\le s}B_p[(\pd_t^{s'}\phi,\pd_t^{s'}\psi)](t),$$
and
$$\mathring{N}^s(t)=(\mathring{E}^s(t))^{1/2}\sum_{s'\le s}\left(||\pd_t^{s'}\mathcal{N}_\phi||_{L^2(\Sigma_t)}+||\pd_t^{s'}\mathcal{N}_\psi||_{L^2(\tilde{\Sigma}_t)}\right),$$
\begin{multline*}
\mathring{N}_p^s(t)=\sum_{s'\le s}\int_{\Sigma_t}r^{p+1}(\pd_t^{s'}\mathcal{N}_\phi)^2 + \sum_{s'\le s}\int_{\tilde\Sigma_t}r^{p+1}(\pd_t^{s'}\mathcal{N}_\psi)^2 \\
+ \sum_{s'\le s}\int_{\Sigma_t\cap\{r\approx r_{trap}\}}|\pd_t^{s'+1}\phi\pd_t^{s'}\mathcal{N}_\phi| + \sum_{s'\le s}\int_{\tilde\Sigma_t\cap\{r\approx r_{trap}\}}|\pd_t^{s'+1}\psi\pd_t^{s'}\mathcal{N}_\psi|. 
\end{multline*}
\end{theorem}
\begin{proof}
The proof is a direct application of Theorem \ref{p_thm} by making the substitutions
\begin{align*}
\phi &\mapsto \pd_t^{s'}\phi \\
\psi &\mapsto \pd_t^{s'}\psi
\end{align*}
for all values of $s'$ where $s'\le s$, and observing that if
\begin{align*}
\Box_g\phi &= \mathcal{L}_\phi+\mathcal{N}_\phi, \\
\Box_{\tilde{g}}\psi &= \mathcal{L}_\psi + \mathcal{N}_\psi,
\end{align*}
then
\begin{align*}
\Box_g(\pd_t^{s'}\phi) &= \mathcal{L}_{(\pd_t^{s'}\phi)}+\pd_t^{s'}\mathcal{N}_\phi, \\
\Box_{\tilde{g}}(\pd_t^{s'}\psi) &= \mathcal{L}_{(\pd_t^{s'}\psi)} + \pd_t^{s'}\mathcal{N}_\psi,
\end{align*}
where $\mathcal{L}_{(\pd_t^{s'}\phi)}$ and $\mathcal{L}_{(\pd_t^{s'}\psi)}$ are the expressions obtained by replacing $(\phi,\psi)$ with $(\pd_t^{s'}\phi,\pd_t^{s'}\psi)$ in $\mathcal{L}_\phi$ and $\mathcal{L}_\psi$ respectively.
\end{proof}

\subsection{The dynamic estimates for $(\phi^s,\psi^s)$}

Even though $\pd_t$ is the only operator that completely commutes with the linear system, the second order Carter operator $Q$ and its modified version $\tilde{Q}$ commute with the wave operators $q^2\Box_g$ and $q^2\Box_{\tilde{g}}$ respectively.

\begin{definition}
$$Q=\pd_\theta^2+\cot\theta \pd_\theta +a^2\sin^2\theta \pd_t^2,$$
$$\tilde{Q}=\pd_\theta^2+5\cot\theta \pd_\theta +a^2\sin^2\theta \pd_t^2.$$
\end{definition}

\begin{lemma}
$$[Q,q^2\Box_g]=0,$$
$$[\tilde{Q},q^2\Box_{\tilde{g}}]=0.$$
\end{lemma}
\begin{proof}
One can check that the operator $q^2\Box_g-Q$ does not depend on $\theta$ or $t$ and has no $\pd_\theta$ operators. Since $Q$ only depends on $\theta$, and has only $\pd_\theta$ and $\pd_t$ operators, that means
$$0=[Q,q^2\Box_g-Q]=[Q,q^2\Box_g].$$
Similarly,
\begin{align*}
q^2\Box_{\tilde{g}}-\tilde{Q} &= q^2\Box_g+2q^2\frac{\pd^\alpha A_1}{A_1}\pd_\alpha-(Q+4\cot\theta\pd_\theta) \\
&=\Box_g-Q+\left(2q^2\frac{\pd^\alpha A_1}{A_1}\pd_\alpha-4\cot\theta\pd_\theta\right) \\
&=\Box_g-Q+\frac{4r}{r^2+a^2}(q^2g^{\alpha r})\pd_\alpha.
\end{align*}
Since $q^2\Box_g-Q$ and $\frac{4r}{r^2+a^2}(q^2g^{\alpha r})\pd_\alpha$ do not depend on $\theta$ or $t$ and have no $\pd_\theta$ operators, and since $\tilde{Q}$ only depends on $\theta$ and has only $\pd_\theta$ and $\pd_t$ operators, that means
$$0=[\tilde{Q},q^2\Box_{\tilde{g}}-\tilde{Q}]=[\tilde{Q},q^2\Box_{\tilde{g}}].$$
This completes the proof. 
\end{proof}

We define the $s$-order commutators $\Gamma$ and $\tilde{\Gamma}$.
\begin{definition}
$$\Gamma^su = Q^{l}\pd_t^{s-2l}u$$
$$\tilde\Gamma^su = \tilde{Q}^{l}\pd_t^{s-2l}u,$$
where $0\le 2l\le s$.
\end{definition}
We also define the $s$-order dynamic quantities $\phi^s$ and $\psi^s$.
\begin{definition}
$$\phi^s=\Gamma^s\phi,$$
$$\psi^s=\tilde\Gamma^s\psi.$$
\end{definition}
Additionally, we define the $s$-order analogues of $\mathcal{L}_\phi$, $\mathcal{L}_\psi$, $\mathcal{X}_\phi$, and $\mathcal{X}_\psi$.
\begin{definition}
$$\mathcal{L}_{\phi^s}=-2\frac{\pd^\alpha B}{A}A\pd_\alpha \psi^s + 2\frac{\pd^\alpha B\pd_\alpha B}{A^2}\phi^s-4\frac{\pd^\alpha A\pd_\alpha B}{A^2} A\psi^s$$
$$\mathcal{L}_{\psi^s}=-2\frac{\pd^\alpha A_2}{A_2}\pd_\alpha\psi^s+2\frac{\pd^\alpha B\pd_\alpha B}{A^2}\psi^s + 2A^{-1}\frac{\pd^\alpha B}{A}\pd_\alpha\phi^s$$
and
\begin{align*}
\mathcal{X}_{\phi^s} &= 2X(\phi^s)+w\phi^s+w_{(a)}\psi^s \\
\mathcal{X}_{\psi^s} &= 2X(\psi^s)+\tilde{w}\psi^s+\tilde{w}_{(a)}\phi^s,
\end{align*}
where, in each expression, the exact same operator $\Gamma^s$ or $\tilde{\Gamma}^s$ is to be used in each term on the right side, replacing $Q$ with $\tilde{Q}$ where appropriate. For example, the expression
$$-2\frac{\pd^\alpha B}{A}A\pd_\alpha (\tilde{Q}\psi) + 2\frac{\pd^\alpha B\pd_\alpha B}{A^2}(Q\phi)-4\frac{\pd^\alpha A\pd_\alpha B}{A^2} A(\tilde{Q}\psi)$$
belongs to $\mathcal{L}_{\phi^s}$ ($s=2$) while the expression
$$-2\frac{\pd^\alpha B}{A}A\pd_\alpha(\pd_t^2\psi) + 2\frac{\pd^\alpha B\pd_\alpha B}{A^2}(Q\phi)-4\frac{\pd^\alpha A\pd_\alpha B}{A^2} A(\tilde{Q}\psi)$$
does not.
\end{definition}

We take a look again at the equations
$$\Box_g\phi = \mathcal{L}_\phi +\mathcal{N}_\phi,$$
$$\Box_{\tilde{g}}\psi = \mathcal{L}_\psi +\mathcal{N}_\psi.$$
By applying $\Gamma^s$ and $\tilde\Gamma^s$ repsectively, we obtain additional useful equations.
\begin{align*}
\Box_g\phi^s &= q^{-2}\Gamma^s(q^2\mathcal{L}_\phi) + q^{-2}\Gamma^s(q^2\mathcal{N}_\phi) \\
&= \mathcal{L}_{\phi^s} +(q^{-2}\Gamma^s(q^2\mathcal{L}_\phi)-\mathcal{L}_{\phi^s}) +q^{-2}\Gamma^s(q^2\mathcal{N}_\phi),
\end{align*}
\begin{align*}
\Box_{\tilde{g}}\psi^s &= q^{-2}\tilde\Gamma^s(q^2\mathcal{L}_\psi) + q^{-2}\tilde\Gamma^s(q^2\mathcal{N}_\psi) \\
&= \mathcal{L}_{\psi^s} +(q^{-2}\tilde\Gamma^s(q^2\mathcal{L}_\psi)-\mathcal{L}_{\psi^s}) +q^{-2}\tilde\Gamma^s(q^2\mathcal{N}_\psi).
\end{align*}
Therefore, Theorem \ref{p_L_thm} can be generalized to the following theorem. (Note the presence of the additional terms $(L_2)^s(t)$ and $(L_2)_p^s(t)$, which arise from the fact that the $\Gamma$ and $\tilde\Gamma$ operators do not completely commute with the linear system.)
\begin{theorem}\label{p_s_L_thm}
Suppose $|a|/M$ is sufficiently small and fix $\delm,\delp>0$. The following estimates hold for $p\in [\delm,2-\delp]$.
$$E^s(t_2)\lesssim E^s(t_1)+\int_{t_1}^{t_2}(L_2)^s(t)+N^s(t)dt,$$
$$E_p^s(t_2)+\int_{t_1}^{t_2}B_p^s(t)dt\lesssim E_p^s(t_1)+\int_{t_1}^{t_2}(L_1)_p^s(t)+(L_2)_p^s(t)+N_p^s(t)dt,$$
where
\begin{align*}
E^s(t) &= \sum_{s'\le s} E[(\phi^{s'},\psi^{s'})](t), \\
E_p^s(t) &= \sum_{s'\le s} E_p[(\phi^{s'},\psi^{s'})](t), \\
B_p^s(t) &= \sum_{s'\le s} B_p[(\phi^{s'},\psi^{s'})](t),
\end{align*}
and
$$(L_2)^s(t)=\sum_{s'\le s}\int_{\Sigma_t}|\pd_t\phi^{s'} (q^{-2}\Gamma^{s'}(q^2\mathcal{L}_\phi)-\mathcal{L}_{\phi^{s'}})|+\sum_{s'\le s}\int_{\tilde\Sigma_t}|\pd_t\psi^{s'}(q^{-2}\tilde{\Gamma}^{s'}(q^2\mathcal{L}_\psi)-\mathcal{L}_{\psi^{s'}})|,$$
$$N^s(t)=\sum_{s'\le s}\int_{\Sigma_t}|\pd_t\phi^{s'}(\Box_g\phi^{s'}-q^{-2}\Gamma^{s'}(q^2\mathcal{L}_\phi))| +\sum_{s'\le s}\int_{\tilde{\Sigma}_t}|\pd_t\psi^{s'}(\Box_{\tilde g}\psi^{s'}-q^{-2}\tilde\Gamma^{s'}(q^2\mathcal{L}_\psi))|,$$
$$(L_1)_p^s(t) = \sum_{s'\le s}\int_{\Sigma_t\cap\{r>R\}}r^p(|L\phi^{s'}|+r^{-1}|\phi^{s'}|)|\mathcal{L}_{\phi^{s'}}|+\sum_{s'\le s}\int_{\tilde\Sigma_t\cap\{r>R\}}r^p(|L\psi^{s'}|+r^{-1}|\psi^{s'}|)|\mathcal{L}_{\psi^{s'}}|,$$
$$(L_2)_p^s(t) = \sum_{s'\le s}\int_{\Sigma_t}r^p|\mathcal{X}_{\phi^s}||q^{-2}\Gamma^{s'}(q^2\mathcal{L}_\phi)-\mathcal{L}_{\phi^{s'}}| + \sum_{s'\le s}\int_{\tilde\Sigma_t}r^p|\mathcal{X}_{\psi^s}||q^{-2}\tilde\Gamma^{s'}(q^2\mathcal{L}_\psi)-\mathcal{L}_{\psi^{s'}}|,$$
$$N_p^s(t)= \sum_{s'\le s}\int_{\Sigma_t}r^p|\mathcal{X}_{\phi^s}||\Box_g\phi^{s'}-q^{-2}\Gamma^{s'}(q^2\mathcal{L}_\phi)| + \sum_{s'\le s}\int_{\tilde\Sigma_t}r^p|\mathcal{X}_{\psi^s}||\Box_{\tilde{g}}\psi^{s'}-q^{-2}\tilde\Gamma^{s'}(q^2\mathcal{L}_\psi)|.$$
\end{theorem}
\begin{proof}
The proof is a direct application of Theorem \ref{p_L_thm} by making the substitutions
\begin{align*}
\phi&\mapsto\phi^{s'} \\
\psi&\mapsto\psi^{s'}
\end{align*}
for all values of $s'$ (and all commutators represented by $\Gamma^{s'}$ and $\tilde\Gamma^{s'}$) where $s'\le s$. The following estimates are used.
$$|\Box_g\phi^{s'}-\mathcal{L}_{\phi^{s'}}| \le |\Box_g\phi^{s'}-q^{-2}\Gamma^{s'}(q^2\mathcal{L}_\phi)|+|q^{-2}\Gamma^{s'}(q^2\mathcal{L}_\phi)-\mathcal{L}_{\phi^{s'}}|,$$
$$|\Box_{\tilde{g}}\psi^{s'}-\mathcal{L}_{\psi^{s'}}| \le |\Box_{\tilde{g}}\psi^{s'}-q^{-2}\tilde\Gamma^{s'}(q^2\mathcal{L}_\psi)|+|q^{-2}\tilde\Gamma^{s'}(q^2\mathcal{L}_\psi)-\mathcal{L}_{\psi^{s'}}|.$$
The resulting error terms have been grouped into the parts that will either be linear or nonlinear when using the equations for the fully nonlinear system. (See Theorem \ref{p_s_thm} below.) 
\end{proof}

As in the proof of Theorem \ref{p_thm}, we would like to absorb the linear error terms into the bulk. But unfortunately, the terms $(L_2)^s(t)$ and $(L_2)_p^s(t)$ are not as straightforward to eliminate. The term $(L_2)^s(t)$ cannot be absorbed, because it belongs to the classic energy estimate, which has no bulk quantity on the left side. The term $(L_2)_p^s(t)$, which belongs to the $r^p$ estimate, cannot be completely absorbed into the bulk on the left side, because of a complication at the trapping radius.

The strategy is as follows. In Lemma \ref{p_s_L_a_1_lem}, the linear terms $(L_2)^s(t)$, $(L_1)_p^s(t)$, and $(L_2)_p^s(t)$ are estimated by the appropriate bulk norms except near the trapping radius. The trapping radius will need special care, because the factors $\pd_t\phi^s$ and $\pd_t\psi^s$, which appear in each of these linear terms, cannot be estimated by the appropriate bulk norm. Actually, for the particular case where $\pd_t\phi^s$ (resp. $\pd_t\psi^s$) represents $\pd_t^{s+1}\phi$ (resp. $\pd_t^{s+1}\psi$), then this factor can be estimated by the homogeneous bulk norm $\mathring{B}^{s+1}(t)$ with a loss of one derivative. The advantage in this case is that the homogenous bulk norm has already been estimated in Theorem \ref{p_o_thm}. The problem is that $\pd_t\phi^s$ (resp. $\pd_t\psi^s$) also represents terms with the operator $Q$ (resp. $\tilde{Q}$). In Lemma \ref{pd_t_identities_lem}, an approximate identity is given for the factor $\pd_t\phi^s$ (resp. $\pd_t\psi^s$), which expresses it as a sum of $\pd_t^{s+1}\phi$ (resp. $\pd_t^{s+1}\psi$) and other terms. Then in Lemma \ref{p_s_L_a_lem}, this approximate identity is used to refine Lemma \ref{p_s_L_a_1_lem}. The result is that there will be a loss of one derivative (although to a homogeneous norm) as well as nonlinear terms.

\begin{lemma}\label{p_s_L_a_1_lem}
If $(L_2)^s(t)$, $(L_1)_p^s(t)$, $(L_2)_p^s(t)$, and $B_p(t)$ are defined as in the previous theorem (but with the absolute values moved outside the integral--see the remark below), then
$$(L_2)^s(t) \lesssim \frac{|a|}{M}B_1^s(t) +Err_{trap}$$
and
$$(L_1)_p^s(t)+(L_2)_p^s(t) \lesssim \frac{|a|}{M}B_p^s(t) +Err_{trap},$$
where
\begin{multline*}
Err_{trap} =\left|\int_{\Sigma_t\cap\{r\approx r_{trap}\}}\pd_t\phi^{s}(q^{-2}\Gamma^{s}(q^2\mathcal{L}_\phi)-\mathcal{L}_{\phi^{s}})\right| \\
+\left|\int_{\tilde\Sigma_t\cap\{r\approx r_{trap}\}}\pd_t\psi^{s}(q^{-2}\tilde\Gamma^{s}(q^2\mathcal{L}_\psi)-\mathcal{L}_{\psi^{s}})\right|.
\end{multline*}
\end{lemma}
\begin{remark}
The quantities $(L_2)^s(t)$, $(L_1)_p^s(t)$, and $(L_2)_p^s(t)$ have been slightly redefined so that the absolute value is moved outside the integral. This small detail will be important in the proof of Lemma \ref{p_s_L_a_lem}, because it allows integration by parts. The concerned reader can easily check that all of the estimates developed so far are also valid with the absolute value outside the integral.
\end{remark}
\begin{proof}
From Lemma \ref{p_L_a_lem}, by replacing $\phi$ and $\psi$ with $\phi^s$ and $\psi^s$, we have
$$(L_1)_p^s(t)\lesssim \frac{|a|}{M}B_p^s(t).$$
The challenge is to estimate the new terms $(L_2)^s(t)$ and $(L_2)_p^s(t)$, which arise from commuting with the operators $Q$ and $\tilde{Q}$.

Observe that
\begin{multline*}
\left|\int_{\Sigma_t}\pd_t\phi^s(q^{-2}\Gamma^s(q^2\mathcal{L}_\phi)-\mathcal{L}_{\phi^s})\right| \\
\lesssim \frac{|a|}{M}\int_{\Sigma_t\setminus \{r\approx r_{trap}\}}r^{-2}(\pd_t\phi^s)^2 + \left(\frac{|a|}{M}\right)^{-1}\int_{\Sigma_t\setminus \{r\approx r_{trap}\}}r^2(q^{-2}\Gamma^s(q^2\mathcal{L}_\phi)-\mathcal{L}_{\phi^s})^2 + Err_{trap} \\
\lesssim \frac{|a|}{M}B_1^s(t)+ \left(\frac{|a|}{M}\right)^{-1}\int_{\Sigma_t\setminus \{r\approx r_{trap}\}}r^2(q^{-2}\Gamma^s(q^2\mathcal{L}_\phi)-\mathcal{L}_{\phi^s})^2 + Err_{trap},
\end{multline*}
and
\begin{multline*}
\left|\int_{\Sigma_t}r^pL\phi^s(q^{-2}\Gamma^s(q^2\mathcal{L}_\phi)-\mathcal{L}_{\phi^s})\right| \\
\lesssim \frac{|a|}{M}\int_{\Sigma_t\setminus \{r\approx r_{trap}\}}r^{p-1}(L\phi^s)^2 + \left(\frac{|a|}{M}\right)^{-1}\int_{\Sigma_t\setminus \{r\approx r_{trap}\}}r^{p+1}(q^{-2}\Gamma^s(q^2\mathcal{L}_\phi)-\mathcal{L}_{\phi^s})^2 + Err_{trap} \\
\lesssim \frac{|a|}{M}B^s_p(t)+ \left(\frac{|a|}{M}\right)^{-1}\int_{\Sigma_t\setminus \{r\approx r_{trap}\}}r^{p+1}(q^{-2}\Gamma^s(q^2\mathcal{L}_\phi)-\mathcal{L}_{\phi^s})^2 + Err_{trap}.
\end{multline*}
From these two particular example estimates, it should become clear that the lemma reduces to the following estimate
\begin{multline*}
\int_{\Sigma_t\setminus\{r\approx r_{trap}\}} r^{p+1}(q^{-2}\Gamma^s(q^2\mathcal{L}_\phi)-\mathcal{L}_{\phi^s})^2 + \int_{\tilde\Sigma_t\setminus\{r\approx r_{trap}\}} r^{p+1}(q^{-2}\tilde\Gamma^s(q^2\mathcal{L}_\psi)-\mathcal{L}_{\psi^s})^2 \\
\lesssim \frac{a^2}{M^2}B_p^s(t),
\end{multline*}
since one can take $p=1$ to estimate the $(L_2)^s(t)$ term.

Given the estimate (\ref{mathcalL_p_estimate}) for $\mathcal{L}_\phi$ and $\mathcal{L}_\psi$ established in Lemma \ref{p_L_a_lem}, it suffices to show
$$\int_{\Sigma_t\setminus\{r\approx r_{trap}\}} r^{p+1}(q^{-2}\Gamma^s(q^2\mathcal{L}_\phi))^2 + \int_{\tilde\Sigma_t\setminus\{r\approx r_{trap}\}} r^{p+1}(q^{-2}\tilde\Gamma^s(q^2\mathcal{L}_\psi))^2 \lesssim \frac{a^2}{M^2}B_p^s(t).$$
This follows from the formalism developed in \S\ref{regularity_sec}. 
\end{proof}

In a moment, we will estimate the error terms $Err_{trap}$ from the previous lemma. But in order to do so, we need the approximate identities given by the following lemma.

\begin{lemma}\label{pd_t_identities_lem} In a neighborhood of $r_{trap}$, the following identities hold in the sense that each term on the right side is missing a smooth factor.
$$\pd_t\phi^s \approx \pd_t^{s+1}\phi + \pd_r^2\phi^{s-1} +r^{-1}\pd_r\phi^{s-1}+q^{-2}\Gamma^{s-1}(q^2\mathcal{L}_\phi) +q^{-2}\Gamma^{s-1}(q^2\Box_g\phi-q^2\mathcal{L}_\phi),$$
$$\pd_t\psi^s \approx \pd_t^{s+1}\psi + \pd_r^2\psi^{s-1} +r^{-1}\pd_r\psi^{s-1}+q^{-2}\Gamma^{s-1}(q^2\mathcal{L}_\psi) +q^{-2}\Gamma^{s-1}(q^2\Box_{\tilde{g}}\psi-q^2\mathcal{L}_\psi).$$
\end{lemma}
\begin{proof}
By the definition of $\phi^s$ and the fact that $[\pd_t,Q]=0$,
$$\pd_t\phi^s = \pd_t^{s-2i+1}Q^i\phi.$$
We use the approximate identity
$$Q\approx \pd_t^2 +\pd_r^2 +r^{-1}\pd_r +\Box_g.$$
(In reality, there is a factor of $q^2$ missing, but this function is smooth and bounded in a neighborhood of $r_{trap}$.)

We have
\begin{align*}
\pd_t^{s-2i+1}Q^i\phi &\approx \pd_t^{s-2i+1}Q^{i-1}\left(\pd_t^2\phi +\pd_r^2\phi +r^{-1}\pd_r\phi +\Box_g\phi\right) \\
&\approx \pd_t^{s-2(i-1)+1}Q^{i-1}\phi + \pd_r^2\phi^{s-1}+r^{-1}\pd_r\phi^{s-1}+\Gamma^{s-1}(\Box_g\phi-\mathcal{L}_\phi) +\Gamma^{s-1}(\mathcal{L}_\phi)
\end{align*}
Repeating this procedure $i-1$ more times proves the first approximate identity of the lemma. The second approximate identity is proved the same way. 
\end{proof}

Now we improve the estimates from Lemma \ref{p_s_L_a_1_lem} by estimating the error term $Err_{trap}$.
\begin{lemma}\label{p_s_L_a_lem}
If $(L_2)^s(t)$, $(L_1)_p^s(t)$, $(L_2)_p^s(t)$, and $B_p^s(t)$ are as defined in Theorem \ref{p_s_L_thm} (but with the absolute values moved outside of the integral as in Lemma \ref{p_s_L_a_1_lem}), then
$$(L_2)^s(t) \lesssim \frac{|a|}{M}(B_1^s(t)+\mathring{B}_1^{s+1}(t)+N^s(t))$$
and
$$(L_1)_p^s(t)+(L_2)_p^s(t)\lesssim \frac{|a|}{M}(B_p^s(t)+\mathring{B}_p^{s+1}(t)+N_p^s(t)),$$
where
\begin{multline*}
N^s(t)=(E^s(t))^{1/2}\left(\sum_{s'\le s}||q^{-2}\Gamma^{s'}(q^2\Box_g\phi-q^2\mathcal{L}_\phi)||_{L^2(\Sigma_t)}\right. \\
\left.+\sum_{s'\le s}||q^{-2}\tilde\Gamma^{s'}(q^2\Box_{\tilde{g}}\psi-q^2\mathcal{L}_\psi)||_{L^2(\tilde{\Sigma}_t)}\right),
\end{multline*}
$$N_p^s(t) = \sum_{s'\le s}\int_{\Sigma_t}r^{p+1}|q^{-2}\Gamma^{s'}(q^2\Box_g\phi-q^2\mathcal{L}_\phi)|^2 + \sum_{s'\le s}\int_{\tilde\Sigma_t}r^{p+1}|q^{-2}\tilde\Gamma^{s'}(q^2\Box_{\tilde{g}}\psi-q^2\mathcal{L}_\psi)|^2.$$
\end{lemma}
\begin{proof}
From Lemma \ref{p_s_L_a_1_lem}, it suffices to show that
\begin{multline*}
\left|\int_{\Sigma_t\cap\{r\approx r_{trap}\}}\pd_t\phi^{s}(q^{-2}\Gamma^{s}(q^2\mathcal{L}_\phi)-\mathcal{L}_{\phi^{s}})\right| +\left|\int_{\tilde\Sigma_t\cap\{r\approx r_{trap}\}}\pd_t\psi^{s}(q^{-2}\tilde\Gamma^{s}(q^2\mathcal{L}_\psi)-\mathcal{L}_{\psi^{s}})\right| \\
\lesssim \frac{|a|}{M}(B_p^s(t) + \mathring{B}_p^{s+1}(t)+N^s(t))
\end{multline*}
and
\begin{multline*}
\left|\int_{\Sigma_t\cap\{r\approx r_{trap}\}}\pd_t\phi^{s}(q^{-2}\Gamma^{s}(q^2\mathcal{L}_\phi)-\mathcal{L}_{\phi^{s}})\right| +\left|\int_{\tilde\Sigma_t\cap\{r\approx r_{trap}\}}\pd_t\psi^{s}(q^{-2}\tilde\Gamma^{s}(q^2\mathcal{L}_\psi)-\mathcal{L}_{\psi^{s}})\right| \\
\lesssim \frac{|a|}{M}(B_p^s(t) + \mathring{B}_p^{s+1}(t)+N_p^s(t)).
\end{multline*}
We prove these estimates by using the approximate identities from Lemma \ref{pd_t_identities_lem}, ignoring the factor of $|a|/M$, which clearly comes from the factors $(q^{-2}\Gamma^s(q^2\mathcal{L}_\phi)-\mathcal{L}_{\phi^s})$ and $(q^{-2}\tilde{\Gamma}^s(q^2\mathcal{L}_\psi)-\mathcal{L}_{\psi^s})$. The following examples illustrate all of the difficulties.
\begin{multline*}
\left|\int_{\Sigma_t\cap\{r\approx r_{trap}\}}(\pd_t^{s+1}\phi)(q^{-2}\Gamma^{s}(q^2\mathcal{L}_\phi)-\mathcal{L}_{\phi^{s}})\right| \\
\lesssim \int_{\Sigma_t\cap\{r\approx r_{trap}\}}(\pd_t^{s+1}\phi)^2 + \int_{\Sigma_t\cap\{r\approx r_{trap}\}}(q^{-2}\Gamma^{s}(q^2\mathcal{L}_\phi)-\mathcal{L}_{\phi^{s}})^2 \\
\lesssim \mathring{B}_p^{s+1}(t) + B_p^s(t).
\end{multline*}
\begin{multline*}
\left|\int_{\Sigma_t\cap\{r\approx r_{trap}\}}(\pd_r^2\Gamma^{s-1}\phi)(q^{-2}\Gamma^{s}(q^2\mathcal{L}_\phi)-\mathcal{L}_{\phi^{s}})\right| \\
\lesssim \left|\int_{\Sigma_t\cap\{r\approx r_{trap}\}}(\pd_r\Gamma^{s-1}\phi)\pd_r(q^{-2}\Gamma^{s}(q^2\mathcal{L}_\phi)-\mathcal{L}_{\phi^{s}})\right| \\
\lesssim \int_{\Sigma_t\cap\{r\approx r_{trap}\}}(\pd_r\Gamma^{s-1}\phi)^2 + \int_{\Sigma_t\cap\{r\approx r_{trap}\}}(\pd_r(q^{-2}\Gamma^{s}(q^2\mathcal{L}_\phi)-\mathcal{L}_{\phi^{s}}))^2 \\
\lesssim B_p^s(t).
\end{multline*}
In particular, note in the previous estimate the need to integrate by parts. The expression $q^{-2}\Gamma^s(q^2\mathcal{L}_\phi)-\mathcal{L}_{\phi^s}$ is a commutator, so it has at most $s$ derivatives. We can assume that none of these derivatives is $\pd_r$, because otherwise, instead of using the approximate identities from Lemma \ref{pd_t_identities_lem}, we could have simply integrated by parts to move one of the angular derivatives in $\pd_t\phi^s$ to the other factor. We continue with more examples.
\begin{multline*}
\left|\int_{\Sigma_t\cap\{r\approx r_{trap}\}}(q^{-2}\Gamma^{s-1}(q^2\mathcal{L}_\phi))(q^{-2}\Gamma^{s}(q^2\mathcal{L}_\phi)-\mathcal{L}_{\phi^{s}})\right|  \\
\lesssim \int_{\Sigma_t\cap\{r\approx r_{trap}\}}(q^{-2}\Gamma^{s-1}(q^2\mathcal{L}_\phi))^2 + \int_{\Sigma_t\cap\{r\approx r_{trap}\}}(q^{-2}\Gamma^{s}(q^2\mathcal{L}_\phi)-\mathcal{L}_{\phi^{s}})^2 \\
\lesssim B_p^s(t).
\end{multline*}

\begin{multline*}
\left|\int_{\Sigma_t\cap\{r\approx r_{trap}\}}q^{-2}\Gamma^{s-1}(q^2\Box_g\phi-q^2\mathcal{L}_\phi)(q^{-2}\Gamma^{s}(q^2\mathcal{L}_\phi)-\mathcal{L}_{\phi^{s}})\right| \\
\lesssim \int_{\Sigma_t\cap\{r\approx r_{trap}\}}(q^{-2}\Gamma^{s-1}(q^2\Box_g\phi-q^2\mathcal{L}_\phi))^2 + \int_{\Sigma_t\cap\{r\approx r_{trap}\}}(q^{-2}\Gamma^{s}(q^2\mathcal{L}_\phi)-\mathcal{L}_{\phi^{s}})^2 \\
\lesssim N_p^s(t)+B_p^s(t).
\end{multline*}
Also,
\begin{multline*}
\left|\int_{\Sigma_t\cap\{r\approx r_{trap}\}}q^{-2}\Gamma^{s-1}(q^2\Box_g\phi-q^2\mathcal{L}_\phi)(q^{-2}\Gamma^{s}(q^2\mathcal{L}_\phi)-\mathcal{L}_{\phi^{s}})\right| \\
\lesssim \left(\int_{\Sigma_t\cap\{r\approx r_{trap}\}}(q^{-2}\Gamma^{s-1}(q^2\Box_g\phi-q^2\mathcal{L}_\phi))^2\right)^{1/2}\left(\int_{\Sigma_t\cap\{r\approx r_{trap}\}}(q^{-2}\Gamma^{s}(q^2\mathcal{L}_\phi)-\mathcal{L}_{\phi^{s}})^2\right)^{1/2} \\
\lesssim ||q^{-2}\Gamma^{s-1}(q^2\Box_g\phi-q^2\mathcal{L}_\phi)||_{L^2(\Sigma_t)}(E^s(t))^{1/2} \\
\lesssim N^s(t).
\end{multline*}
These example estimates are sufficient to verify the lemma. 
\end{proof}

By using Lemma \ref{p_s_L_a_lem} and assuming the full nonlinear equations, Theorem \ref{p_s_L_thm} can be improved slightly. This is stated in the following theorem.
\begin{theorem}\label{p_s_thm}(Improved version of Theorem \ref{p_s_L_thm}) Suppose $|a|/M$ is sufficiently small and fix $\delm,\delp>0$. Suppose furthermore that the pair $(\phi,\psi)$ satisfies the system
\begin{align*}
\Box_g\phi &= \mathcal{L}_\phi+\mathcal{N}_\phi, \\
\Box_{\tilde{g}}\psi &= \mathcal{L}_\psi + \mathcal{N}_\psi.
\end{align*}
Then the following estimates hold for $p\in [\delm,2-\delp]$.
$$E^s(t_2)\lesssim E^s(t_1)+\int_{t_1}^{t_2}\mathring{B}_1^{s+1}(t)+B_1^s(t)+N^s(t)dt$$
$$E_p^s(t_2)+\int_{t_1}^{t_2}B_p^s(t)dt\lesssim E_p^s(t_1)+\int_{t_1}^{t_2}\mathring{B}_p^{s+1}(t)+N_p^s(t)dt,$$
where $E^s(t)$, $E_p^s(t)$, and $B_p^s(t)$ are as defined in Theorem \ref{p_s_L_thm}, $\mathring{B}^s(t)$ is as defined in Theorem \ref{p_o_thm}, and
$$N^s(t)=(E^s(t))^{1/2}\left(\sum_{s'\le s}||q^{-2}\Gamma^{s'}(q^2\mathcal{N}_\phi)||_{L^2(\Sigma_t)}+\sum_{s'\le s}||q^{-2}\tilde\Gamma^{s'}(q^2\mathcal{N}_\psi)||_{L^2(\tilde{\Sigma}_t)}\right),$$
\begin{align*}
N_p^s(t) =& \sum_{s'\le s}\int_{\Sigma_t}r^{p+1}|q^{-2}\Gamma^{s'}(q^2\mathcal{N}_\phi)|^2 + \sum_{s'\le s}\int_{\tilde\Sigma_t}r^{p+1}|q^{-2}\tilde\Gamma^{s'}(q^2\mathcal{N}_\psi)|^2.
\end{align*}
\end{theorem}
\begin{proof}
By Theorem \ref{p_s_L_thm} and Lemma \ref{p_s_L_a_lem}, we have
$$E^s(t_2)\lesssim E^s(t_1)+\int_{t_1}^{t_2}\mathring{B}_1^{s+1}(t)+B_1^s(t)+N^s(t)+(N')^s(t)dt,$$
where
\begin{align*}
(N')^s(t)&=\sum_{s'\le s}\int_{\Sigma_t}|\pd_t\phi^{s'}(\Box_g\phi^{s'}-q^{-2}\Gamma^{s'}(q^2\mathcal{L}_\phi))| \\
&\hspace{2in}+\sum_{s'\le s}\int_{\tilde{\Sigma}_t}|\pd_t\psi^{s'}(\Box_{\tilde g}\psi^{s'}-q^{-2}\tilde\Gamma^{s'}(q^2\mathcal{L}_\psi))| \\
&=\sum_{s'\le s}\int_{\Sigma_t}|\pd_t\phi^{s'}q^{-2}\Gamma^{s'}(q^2\mathcal{N}_{\phi})| +\sum_{s'\le s}\int_{\tilde{\Sigma}_t}|\pd_t\psi^{s'}q^{-2}\tilde{\Gamma}^{s'}(q^2\mathcal{N}_{\psi})|.
\end{align*}
Observe that
\begin{align*}
(N')^s(t) &=\sum_{s'\le s}\int_{\Sigma_t}|\pd_t\phi^{s'}q^{-2}\Gamma^{s'}(q^2\mathcal{N}_{\phi})| +\sum_{s'\le s}\int_{\tilde{\Sigma}_t}|\pd_t\psi^{s'}q^{-2}\tilde{\Gamma}^{s'}(q^2\mathcal{N}_{\psi})| \\
&\lesssim \sum_{s'\le s}||\pd_t\phi^{s'}||_{L^2(\Sigma_t)}||q^{-2}\Gamma^{s'}(q^2\mathcal{N}_\phi)||_{L^2(\Sigma_t)} \\
&\hspace{2.5in}+ \sum_{s'\le s}||\pd_t\psi^{s'}||_{L^2(\tilde\Sigma_t)}||q^{-2}\tilde\Gamma^{s'}(q^2\mathcal{N}_\psi)||_{L^2(\tilde\Sigma_t)} \\
&\lesssim \sum_{s'\le s}(E^{s'}(t))^{1/2}||q^{-2}\Gamma^{s'}(q^2\mathcal{N}_\phi)||_{L^2(\Sigma_t)}+ \sum_{s'\le s}(E^{s'}(t))^{1/2}||q^{-2}\tilde\Gamma^{s'}(q^2\mathcal{N}_\psi)||_{L^2(\tilde\Sigma_t)} \\
&\lesssim N^s(t).
\end{align*}
This proves the first estimate of the theorem.

By Theorem \ref{p_s_L_thm} and Lemma \ref{p_s_L_a_lem}, we have
$$E_p^s(t_2)+\int_{t_1}^{t_2}B_p^s(t)dt\lesssim E_p^s(t_1)+\int_{t_1}^{t_2}\mathring{B}_p^s(t) +\frac{|a|}{M}B_p^s(t)+N_p^s(t)+(N')_p^s(t)dt,$$
where
\begin{align*}
(N')_p^s(t)=& \sum_{s'\le s}\int_{\Sigma_t}r^p|\mathcal{X}_{\phi^s}||\Box_g\phi^{s'}-q^{-2}\Gamma^{s'}(q^2\mathcal{L}_\phi)| \\
&\hspace{2in}+ \sum_{s'\le s}\int_{\tilde\Sigma_t}r^p|\mathcal{X}_{\psi^s}||\Box_{\tilde{g}}\psi^{s'}-q^{-2}\tilde\Gamma^{s'}(q^2\mathcal{L}_\psi)| \\
=& \sum_{s'\le s}\int_{\Sigma_t}r^p|\mathcal{X}_{\phi^s}||q^{-2}\Gamma^{s'}(q^2\mathcal{N}_\phi)| + \sum_{s'\le s}\int_{\tilde\Sigma_t}r^p|\mathcal{X}_{\psi^s}||q^{-2}\tilde\Gamma^{s'}(q^2\mathcal{N}_\psi)|.
\end{align*}
Observe that
\begin{multline*}
\sum_{s'\le s}\int_{\Sigma_t}r^p|\mathcal{X}_{\phi^s}||q^{-2}\Gamma^{s'}(q^2\mathcal{N}_\phi)| \\
\lesssim \sum_{s'\le s} \left(\int_{\Sigma_t}r^{p-1}(\mathcal{X}_{\phi^s})^2\right)^{1/2}\left(\int_{\Sigma_t}r^{p+1}(q^{-2}\Gamma^{s'}(q^2\mathcal{N}_\phi))^2\right)^{1/2} \\
\lesssim \sum_{s'\le s}\left(\mathring{B}_p^{s'+1}(t)+B_p^{s'}(t)\right)^{1/2}\left(N_p^{s'}(t)\right)^{1/2} \\
\lesssim \epsilon\mathring{B}_p^{s+1}(t)+\epsilon B_p^s(t) +\epsilon^{-1}N_p^s(t).
\end{multline*}
An analogous estimate also shows that
$$\sum_{s'\le s}\int_{\tilde\Sigma_t}r^p|\mathcal{X}_{\psi^s}||q^{-2}\tilde\Gamma^{s'}(q^2\mathcal{N}_\psi)| \lesssim \epsilon\mathring{B}_p^{s+1}(t)+\epsilon B_p^s(t)+\epsilon^{-1}N_p^s(t).$$
Therefore, we have
$$E_p^s(t_2)+\int_{t_1}^{t_2}B_p^s(t)dt\lesssim E_p^s(t_1)+(|a|/M+\epsilon)\int_{t_1}^{t_2}B_p^s(t)dt+\int_{t_1}^{t_2}\mathring{B}_p^s(t)+\epsilon^{-1}N_p^s(t)dt.$$
By taking $|a|/M$ and $\epsilon$ sufficiently small, the bulk term on the right side can be absorbed into the left side. The result is the second estimate of the theorem. 
\end{proof}

\subsection{The dynamic estimates for $(\phi^{s,k},\psi^{s,k})$}

To handle nonlinear terms with a factor of $\pd_r\phi$ or $\pd_r\psi$ near the event horizon, it will be necessary to use the operator $\tg$ introduced in \S\ref{k:additional_commutator_sec} in Chapter \ref{kerr_chap}.

We recall the form of the commutator $[\tg,q^2\Box_g]$ from Lemma \ref{k:dr_comm_2_lem}.
\begin{lemma}\label{dr_comm_2_lem}
$$[\tg^k,q^2\Box_g]u=k\Delta'\pd_r\tg^ku+\{\pd_r\tg^{\le k-1}\Gamma^{\le 1}u,\tg^{\le k-1}\Gamma^{\le 2}u,\tg^{\le k-1}(q^2\Box_gu)\},$$
$$[\tg^k,q^2\Box_{\tilde{g}}]u=k\Delta'\pd_r\tg^ku+\{\pd_r\tg^{\le k-1}\tilde\Gamma^{\le 1}u,\tg^{\le k-1}\tilde\Gamma^{\le 2}u,\tg^{\le k-1}(q^2\Box_{\tilde{g}}u)\},$$
where the $\{...\}$ notation is used to represent terms supported near the event horizon with smooth factors. (See Lemma \ref{k:dr_comm_1_lem}.)
\end{lemma}
\begin{proof}
The first identity is from Lemma \ref{k:dr_comm_2_lem}. The proof for the second identity is similar to the proof of Lemma \ref{k:dr_comm_2_lem}.
\end{proof}

Since $\tg$ doesn't commute with $\Box_g$ or $\Box_{\tilde{g}}$, it is treated seperately from the previous commutators. We define the $s,k$-order dynamic quantities $\phi^{s,k}$ and $\psi^{s,k}$.
\begin{definition}
$$\phi^{s,k} = \tg^k\Gamma^s\phi,$$
$$\psi^{s,k} = \tg^k\tilde\Gamma^s\psi.$$
\end{definition}
We also define the $s,k$-order analogues of $\mathcal{L}_\phi$, $\mathcal{L}_\psi$, $\mathcal{X}_\phi$, and $\mathcal{X}_\psi$.
\begin{definition}
$$\mathcal{L}_{\phi^{s,k}}=-2\frac{\pd^\alpha B}{A}A\pd_\alpha \psi^{s,k} + 2\frac{\pd^\alpha B\pd_\alpha B}{A^2}\phi^{s,k}-4\frac{\pd^\alpha A\pd_\alpha B}{A^2} A\psi^{s,k}$$
$$\mathcal{L}_{\psi^{s,k}}=-2\frac{\pd^\alpha A_2}{A_2}\pd_\alpha\psi^{s,k}+2\frac{\pd^\alpha B\pd_\alpha B}{A^2}\psi^{s,k} + 2A^{-1}\frac{\pd^\alpha B}{A}\pd_\alpha\phi^{s,k}$$
and
\begin{align*}
\mathcal{X}_{\phi^{s,k}} &= 2X(\phi^{s,k})+w\phi^{s,k}+w_{(a)}\psi^{s,k} \\
\mathcal{X}_{\psi^{s,k}} &= 2X(\psi^{s,k})+\tilde{w}\psi^{s,k}+\tilde{w}_{(a)}\phi^{s,k},
\end{align*}
where, in each expression, the exact same operator $\tg^k\Gamma^s$ or $\tg^k\tilde{\Gamma}^s$ is to be used in each term on the right side, replacing $Q$ with $\tilde{Q}$ where appropriate. For example, the expression
$$-2\frac{\pd^\alpha B}{A}A\pd_\alpha (\tg\tilde{Q}\psi) + 2\frac{\pd^\alpha B\pd_\alpha B}{A^2}(\tg Q\phi)-4\frac{\pd^\alpha A\pd_\alpha B}{A^2} A (\tg\tilde{Q}\psi)$$
belongs to $\mathcal{L}_{\phi^{s,k}}$ ($s=2,k=1$) while the expression
$$-2\frac{\pd^\alpha B}{A}A\pd_\alpha(\tg\pd_t^2\psi) + 2\frac{\pd^\alpha B\pd_\alpha B}{A^2}(\tg Q\phi)-4\frac{\pd^\alpha A\pd_\alpha B}{A^2} A (\tg\tilde{Q}\psi)$$
does not.
\end{definition}

We take a look again at the equations
$$\Box_g\phi = \mathcal{L}_\phi +\mathcal{N}_\phi,$$
$$\Box_{\tilde{g}}\psi = \mathcal{L}_\psi +\mathcal{N}_\psi.$$
By applying $\tg^k\Gamma^s$ and $\tg^k\tilde\Gamma^s$ respectively, we obtain additional useful equations.
\begin{align*}
\Box_g\phi^{s,k} &= (\Box_g\phi^{s,k}-q^{-2}\tg^k(q^2\Box_g\phi^s)) +q^{-2}\tg^k\Gamma^s(q^2\mathcal{L}_\phi) + q^{-2}\tg^k\Gamma^s(q^2\mathcal{N}_\phi) \\
&= \mathcal{L}_{\phi^{s,k}} + q^{-2}[q^2\Box_g,\tg^k]\phi^s + (q^{-2}\tg^k\Gamma^s(q^2\mathcal{L}_\phi)-\mathcal{L}_{\phi^{s,k}}) + q^{-2}\tg^k\Gamma^s(q^2\mathcal{N}_\phi),
\end{align*}
\begin{align*}
\Box_{\tilde{g}}\psi^{s,k} &= (\Box_{\tilde{g}}\psi^{s,k}-q^{-2}\tg^k(q^2(\Box_{\tilde{g}}\psi^s)) +q^{-2}\tg^k\tilde{\Gamma}^s(q^2\mathcal{L}_\psi) + q^{-2}\tg^k\tilde{\Gamma}^s(q^2\mathcal{N}_\psi) \\
&= \mathcal{L}_{\psi^{s,k}} + q^{-2}[q^2\Box_{\tilde{g}},\tg^k]\psi^s + (q^{-2}\tg^k\tilde{\Gamma}^s(q^2\mathcal{L}_\psi)-\mathcal{L}_{\psi^{s,k}}) + q^{-2}\tg^k\tilde{\Gamma}^s(q^2\mathcal{N}_\psi).
\end{align*}
Therefore, Theorem \ref{p_s_L_thm} (which generalized Theorem \ref{p_L_thm}) can be generalized to the following theorem. (Note the presence of the additional terms $(L_3)^s(t)$ and $(L_3)_p^s(t)$, which arise from the fact that the $\tg$ operators do not commute with the wave operators.)
\begin{theorem}\label{p_s_k_L_thm}
Suppose $|a|/M$ is sufficiently small and fix $\delm,\delp>0$. The following estimates hold for $p\in [\delm,2-\delp]$.
$$E^{s,k}(t_2)\lesssim E^{s,k}(t_1)+\int_{t_1}^{t_2}(L_2)^{s,k}(t)+(L_3)^{s,k}(t)+N^{s,k}(t)dt,$$
$$E_p^{s,k}(t_2)+\int_{t_1}^{t_2}B_p^{s,k}(t)dt\lesssim E_p^{s,k}(t_1)+\int_{t_1}^{t_2}(L_1)_p^s(t)+(L_2)_p^{s,k}(t)+(L_3)_p^{s,k}(t)+N_p^{s,k}(t)dt,$$
where
\begin{align*}
E^{s,k}(t) &= \sum_{\substack{s'\le s \\ k'\le k}} E[(\phi^{s',k'},\psi^{s',k'})](t), \\
E_p^{s,k}(t) &= \sum_{\substack{s'\le s \\ k'\le k}} E_p[(\phi^{s',k'},\psi^{s',k'})](t), \\
B_p^{s,k}(t) &= \sum_{\substack{s'\le s \\ k'\le k}} B_p[(\phi^{s',k'},\psi^{s',k'})](t),
\end{align*}
and
\begin{align*}
(L_2)^{s,k}(t)=&\sum_{\substack{s'\le s \\ k'\le k}}\int_{\Sigma_t}|\pd_t\phi^{s',k'} (q^{-2}\tg^{k'}\Gamma^{s'}(q^2\mathcal{L}_\phi)-\mathcal{L}_{\phi^{s',k'}})| \\
&+\sum_{\substack{s'\le s \\ k'\le k}}\int_{\tilde\Sigma_t}|\pd_t\psi^{s',k'}(q^{-2}\tg^{k'}\tilde{\Gamma}^{s'}(q^2\mathcal{L}_\psi)-\mathcal{L}_{\psi^{s',k'}})|,
\end{align*}
\begin{align*}
(L_3)^{s,k}(t)=&\sum_{\substack{s'\le s \\ k'\le k}}\int_{\Sigma_t}|\pd_t\phi^{s',k'} (q^{-2}[q^2\Box_g,\tg^{k'}]\phi^{s'})|
+\sum_{\substack{s'\le s \\ k'\le k}}\int_{\tilde\Sigma_t}|\pd_t\psi^{s',k'} (q^{-2}[q^2\Box_{\tilde{g}},\tg^{k'}]\psi^{s'})|,
\end{align*}
\begin{align*}
N^{s,k}(t) =& \sum_{\substack{s'\le s \\ k'\le k}}\int_{\Sigma_t}|\pd_t\phi^{s',k'}(q^{-2}\tg^{k'}(q^2\Box_g\phi^{s'}-\Gamma^{s'}(q^2\mathcal{L}_\phi)))| \\
&+\sum_{\substack{s'\le s \\ k'\le k}}\int_{\tilde{\Sigma}_t}|\pd_t\psi^{s',k'}(q^{-2}\tg^{k'}(q^2\Box_{\tilde g}\psi^{s'}-\tilde\Gamma^{s'}(q^2\mathcal{L}_\psi)))|,
\end{align*}
\begin{align*}
(L_1)_p^s(t) =& \sum_{s'\le s}\int_{\Sigma_t\cap\{r>R\}}r^p(|L\phi^{s'}|+r^{-1}|\phi^{s'}|)|\mathcal{L}_{\phi^{s'}}| \\
&+\sum_{s'\le s}\int_{\tilde\Sigma_t\cap\{r>R\}}r^p(|L\psi^{s'}|+r^{-1}|\psi^{s'}|)|\mathcal{L}_{\psi^{s'}}|,
\end{align*}
\begin{align*}
(L_2)_p^{s,k}(t) =& \sum_{\substack{s'\le s \\ k'\le k}}\int_{\Sigma_t}r^p|\mathcal{X}_{\phi^{s,k}}||q^{-2}\tg^{k'}\Gamma^{s'}(q^2\mathcal{L}_\phi)-\mathcal{L}_{\phi^{s',k'}}|  \\
&+ \sum_{\substack{s'\le s \\ k'\le k}}\int_{\tilde\Sigma_t}r^p|\mathcal{X}_{\psi^{s,k}}||q^{-2}\tg^{k'}\tilde\Gamma^{s'}(q^2\mathcal{L}_\psi)-\mathcal{L}_{\psi^{s',k'}}|,
\end{align*}
\begin{align*}
(L_3)_p^{s,k}(t) =& \sum_{\substack{s'\le s \\ k'\le k}}\int_{\Sigma_t}r^p|\mathcal{X}_{\phi^{s,k}}||q^{-2}[q^2\Box_g,\tg^{k'}]\phi^{s'}|  
+ \sum_{\substack{s'\le s \\ k'\le k}}\int_{\tilde\Sigma_t}r^p|\mathcal{X}_{\psi^{s,k}}||q^{-2}[q^2\Box_{\tilde{g}},\tg^{k'}]\psi^{s'}|,
\end{align*}
\begin{align*}
N_p^{s,k}(t)=& \sum_{\substack{s'\le s \\ k'\le k}}\int_{\Sigma_t}r^p|\mathcal{X}_{\phi^{s,k}}||q^{-2}\tg^{k'}(q^2\Box_g\phi^{s'}-\Gamma^{s'}(q^2\mathcal{L}_\phi))| \\
&+ \sum_{\substack{s'\le s \\ k'\le k}}\int_{\tilde\Sigma_t}r^p|\mathcal{X}_{\psi^{s,k}}||q^{-2}\tg^{k'}(q^2\Box_{\tilde{g}}\psi^{s'}-\tilde\Gamma^{s'}(q^2\mathcal{L}_\psi))|.
\end{align*}
\end{theorem}
\begin{proof}
The proof is a direct application of Theorem \ref{p_L_thm} by making the substitutions
\begin{align*}
\phi&\mapsto\phi^{s',k'} \\
\psi&\mapsto\psi^{s',k'}
\end{align*}
for all values of $s'$ and $k'$ (and all commutators represented by $\tg^k\Gamma^{s'}$ and $\tg^k\tilde\Gamma^{s'}$) where $s'\le s$ and $k'\le k$. The following estimates are used.
\begin{multline*}
|\Box_g\phi^{s',k'}-\mathcal{L}_{\phi^{s',k'}}| \\
\le |q^{-2}[q^2\Box_g,\tg^{k'}]\phi^{s'}|+ |q^{-2}\tg^{k'}(q^2\Box_g\phi^{s'}-\Gamma^{s'}(q^2\mathcal{L}_\phi)|+|q^{-2}\tg^{k'}\Gamma^{s'}(q^2\mathcal{L}_\phi)-\mathcal{L}_{\phi^{s',k'}}|,
\end{multline*}
\begin{multline*}
|\Box_{\tilde{g}}\psi^{s',k'}-\mathcal{L}_{\psi^{s',k'}}| \\
\le |q^{-2}[q^2\Box_{\tilde{g}},\tg^{k'}]\psi^{s'}|+ |q^{-2}\tg^{k'}(q^2\Box_{\tilde{g}}\psi^{s'}-\tilde\Gamma^{s'}(q^2\mathcal{L}_\psi)|+|q^{-2}\tg^{k'}\tilde\Gamma^{s'}(q^2\mathcal{L}_\psi)-\mathcal{L}_{\psi^{s',k'}}|.
\end{multline*}
The resulting error terms have been grouped into the parts that will either be linear or nonlinear when using the equations for the fully nonlinear system. (See Theorem \ref{p_s_k_thm} below.) 
\end{proof}

Once again, the plan is to improve the prevoius theorem after proving a lemma that handles the linear error terms on the right side. In this case, the new linear error terms are $(L_3)^{s,k}(t)$ and $(L_3)_p^{s,k}(t)$. These terms arise from the fact that the operator $\tg$ does not commute with the wave operators. Since $\tg$ is supported near the event horizon, there are no new issues related to trapping. The important observation to make is that one of the terms in the commutator has an appropriate sign on the event horizon. This is the point of the following lemma.
\begin{lemma}\label{p_s_k_L_a_lem}
If $(L_2)^{s,k}(t)$, $(L_3)^{s,k}(t)$, $(L_1)_p^s(t)$, $(L_2)_p^{s,k}(t)$, $(L_3)_p^{s,k}(t)$, and $B_p^{s,k}(t)$ are as defined in Theorem \ref{p_s_k_L_thm}, then
$$(L_2)^{s,k}(t)+(L_3)^{s,k}(t) \lesssim B_{p'}^{s+2,k-1}(t)+\frac{|a|}{M}(B_1^{s,k}(t)+\mathring{B}_1^{s+1}(t))+N^{s,k}(t)$$
and
$$(L_1)_p^s(t)+(L_2)_p^{s,k}(t)+(L_3)_p^{s,k}(t)\lesssim B_{p'}^{s+2,k-1}(t) + \frac{|a|}{M}(B_p^{s,k}(t)+\mathring{B}_p^{s+1}(t))+N_p^{s,k}(t),$$
where
\begin{multline*}
N^{s,k}(t)=(E^{s,k}(t))^{1/2}\left(\sum_{\substack{s'\le s \\ k'\le k}}||q^{-2}\tg^{k'}\Gamma^{s'}(q^2\Box_g\phi-q^2\mathcal{L}_\phi)||_{L^2(\Sigma_t)}\right. \\
+\left.\sum_{\substack{s'\le s \\ k'\le k}}||q^{-2}\tg^{k'}\tilde\Gamma^{s'}(q^2\Box_{\tilde{g}}\psi-q^2\mathcal{L}_\psi)||_{L^2(\tilde{\Sigma}_t)}\right),
\end{multline*}
\begin{multline*}
N_p^{s,k}(t) = \sum_{\substack{s'\le s \\ k'\le k}}\int_{\Sigma_t}r^{p+1}|q^{-2}\tg^{k'}\Gamma^{s'}(q^2\Box_g\phi-q^2\mathcal{L}_\phi)|^2  \\
+ \sum_{\substack{s'\le s \\ k'\le k}}\int_{\tilde\Sigma_t}r^{p+1}|q^{-2}\tg^{k'}\tilde\Gamma^{s'}(q^2\Box_{\tilde{g}}\psi-q^2\mathcal{L}_\psi)|^2,
\end{multline*}
and $p'$ is arbitrary.
\end{lemma}
\begin{remark}
The reason for the arbitrary $p'$ on the right side is that the bulk norms with the arbitrary $p'$ are only used to control the new terms related to commuting with $\tg$. These terms are all supported on a compact radial interval, so the factor $r^{p'-1}$ that appears in $B_{p'}^{s,k}(t)$ can be approximated by a constant.
\end{remark}
\begin{proof}
The proof of Lemma \ref{p_s_L_a_lem} can be adapted to show that
$$(L_2)^{s,k}(t)\lesssim \frac{|a|}{M}(B_1^{s,k}(t)+\mathring{B}_1^{s+1}(t)+N^{s,k}(t)),$$
$$(L_1)_p^{s,k}(t)+(L_2)_p^{s,k}(t)\lesssim \frac{|a|}{M}(B_p^{s,k}(t)+\mathring{B}_p^{s+1}(t)+N_p^{s,k}(t)).$$
It suffices to show that
$$(L_3)^{s,k}(t)\lesssim B_{p'}^{s+2,k-1}(t) + N^{s,k}(t),$$
$$(L_3)_p^{s,k}(t)\lesssim B_{p'}^{s+2,k-1}(t) + N_p^{s,k}(t).$$
The key observation is to recognize that the term represented by $L_3$ actually has a good sign near the event horizon. That is, according to Lemma \ref{dr_comm_2_lem},
\begin{align*}
\int_{\Sigma_t}-2X(\tg^k\phi^s)q^{-2}[q^2\Box_g,\tg^k]\phi^s &= \int_{\Sigma_t}-2(-X^r)\pd_r(\tg^k\phi^s)q^{-2}k\Delta'\pd_r(\tg^k\phi^s)+err \\
&=\int_{\Sigma_t}-2(-X^r)q^{-2}\Delta'(\pd_r\tg^k\phi^s)^2+err
\end{align*}
Since $X^r<0$ near the event horizon and $\Delta'>0$, \textbf{the principal term becomes minus a square}, so it can be ignored, or used to control small error terms.

We now investigate the error terms, which come from the remaining part of the vectorfield $X$ and functions $w$, $\tilde{w}$, $w_{(a)}$, and $\tilde{w}_{(a)}$ in $\mathcal{X}_{\phi^{s,k}}$ and $\mathcal{X}_{\psi^{s,k}}$, as well as the remainder in Lemma \ref{dr_comm_2_lem}.
\begin{multline*}
err = \int_{\Sigma_t\cap \{r_H\le r\le r_H+\delh\}}|\mathcal{X}_{\phi^{s,k}}-X^r\pd_r\phi^{s,k}||q^{-2}k\Delta'\pd_r(\tg^k\phi^s)| \\
+\int_{\Sigma_t\cap\{r_H\le r\le r_H+\delh\}}|\mathcal{X}_{\phi^{s,k}}|(|\pd_r\psi^{k-1,s+1}|+|\psi^{k-1,s+2}|+|q^{-2}\tg^{\le k-1}(q^2\Box_g\phi^s)|).
\end{multline*}
Since there is no product of principal factors (ie. factors of the form $\pd_r\tg^k\phi^s$) each term can be estimated in such a way that at least one factor is estimated by one of $B_{p'}^{s+2,k-1}(t)$ or $N^{s,k}(t)$ or $N_p^{s,k}(t)$, and at most one factor is estimated by $\epsilon (\pd_r\tg^k\phi^s)^2$ after separating the factors. This latter term can be absorbed into the term with the good sign. The same observation and procedure can be repeated for the $\psi$ integral as well.
\end{proof}

Finally, we arrive at the following theorem, which is the most general form required by the proof of the main theorem.
\begin{theorem}\label{p_s_k_thm} (Improved version of Theorem \ref{p_s_k_L_thm}) Suppose $|a|/M$ is sufficiently small and fix $\delm,\delp>0$. Suppose furthermore that the pair $(\phi,\psi)$ satisfies the system
\begin{align*}
\Box_g\phi &= \mathcal{L}_\phi+\mathcal{N}_\phi, \\
\Box_{\tilde{g}}\psi &= \mathcal{L}_\psi + \mathcal{N}_\psi.
\end{align*}
Then the following estimates hold for $p\in [\delm,2-\delp]$, arbitrary $p'$, and integers $s\ge 0$ and $k\ge 1$.
$$E^{s,k}(t_2)\lesssim E^{s,k}(t_1)+\int_{t_1}^{t_2}B_{p'}^{s+2,k-1}(t)+\mathring{B}_1^{s+1}(t)+B_1^{s,k}(t)+N^{s,k}(t)dt,$$
$$E_p^{s,k}(t_2)+\int_{t_1}^{t_2}B_p^{s,k}(t)dt\lesssim E_p^{s,k}(t_1)+\int_{t_1}^{t_2}B_{p'}^{s+2,k-1}(t)+\mathring{B}_p^{s+1}(t)+N_p^{s,k}(t)dt,$$
where $E^{s,k}(t)$, $E_p^{s,k}(t)$, and $B_p^{s,k}(t)$ are as defined in Theorem \ref{p_s_k_L_thm}, and
$$N^{s,k}(t)=(E^{s,k}(t))^{1/2}\left(\sum_{\substack{s'\le s \\ k'\le k}}||q^{-2}\tg^{k'}\Gamma^{s'}(q^2\mathcal{N}_\phi)||_{L^2(\Sigma_t)}+\sum_{\substack{s'\le s \\ k'\le k}}||q^{-2}\tg^{k'}\tilde\Gamma^{s'}(q^2\mathcal{N}_\psi)||_{L^2(\tilde{\Sigma}_t)}\right),$$
\begin{align*}
N_p^{s,k}(t) =& \sum_{\substack{s'\le s \\ k'\le k}}\int_{\Sigma_t}r^{p+1}|q^{-2}\tg^{k'}\Gamma^{s'}(q^2\mathcal{N}_\phi)|^2 + \sum_{\substack{s'\le s \\ k'\le k}}\int_{\tilde\Sigma_t}r^{p+1}|q^{-2}\tg^{k'}\tilde\Gamma^{s'}(q^2\mathcal{N}_\psi)|^2.
\end{align*}
\end{theorem}

\begin{proof}
The proof is very similar to the proof of Theorem \ref{p_s_thm}. An outline of the proof is given here.

By Theorem \ref{p_s_k_L_thm} and Lemma \ref{p_s_k_L_a_lem}, we have
$$
E^{s,k}(t_2) \lesssim E^{s,k}(t_1)+\int_{t_1}^{t_2}B_{p'}^{s+2,k-1}(t)+\mathring{B}_1^{s+1}(t)+B_1^{s,k}(t)+N^{s,k}(t)+(N')^{s,k}(t),$$
where $(N')^{s,k}(t)$ is the quantity $N^{s,k}(t)$ from Theorem \ref{p_s_k_L_thm}. By the same argument as in the proof of Theorem \ref{p_s_thm}, 
$$(N')^{s,k}(t)\lesssim N^{s,k}(t).$$
This proves the first estimate of the theorem.

By Theorem \ref{p_s_k_L_thm} and Lemma \ref{p_s_k_L_a_lem}, we have
\begin{multline*}
E_p^{s,k}(t_2)+\int_{t_1}^{t_2}B_p^{s,k}(t)dt \\
\lesssim E_p^{s,k}(t_1)+\int_{t_1}^{t_2}B_{p'}^{s+2,k-1}(t)+\frac{|a|}{M}B_p^{s,k}(t)+\mathring{B}_p^{s+1}(t)+N_p^{s,k}(t)+(N')_p^{s,k}(t)dt,
\end{multline*}
where $(N')_p^{s,k}(t)$ is the quantity $N_p^{s,k}(t)$ from Theorem \ref{p_s_k_L_thm}. By the same argument as in the proof of Theorem \ref{p_s_thm}, 
$$(N')_p^{s,k}(t)\lesssim \epsilon\mathring{B}_p^{s+1}(t)+\epsilon B_p^{s,k}(t)+\epsilon^{-1}N_p^{s,k}(t).$$
Proceeding as in the proof of Theorem \ref{p_s_thm}, we conclude that if $\epsilon$ and $\frac{|a|}{M}$ are sufficiently small, then the bulk term
$$(|a|/M+\epsilon)\int_{t_1}^{t_2}B_p^{s,k}(t)dt$$
can be absorbed into the left side. The result is the second estimate of the theorem. 
\end{proof}

\section{The $L^\infty$ estimates}\label{wk:pointwise_sec}

In this section, we prove $L^\infty$ estimates for certain derivatives of $\phi$ and $\psi$ that will appear in the nonlinear quantities $q^{-2}\tg^k\Gamma^s(q^2\mathcal{N}_\phi)$ and $q^{-2}\tg^k\tilde{\Gamma}^s(q^2\mathcal{N}_\psi)$. In \S\ref{nl_s_k_structure_sec}, we will see that these quantities are not simply sums of products of single derivatives of $\phi^{s,k}$ and $\psi^{s,k}$. Instead, they will take a slightly more general form consisting of sums of products of single derivatives of $\fd^l\phi^{s-l,k}$ and $\fd^l\psi^{s-l,k}$, where the $\fd^l$ operators are defined in Appendix \ref{regularity_sec}. For this reason, the estimates established in this section will apply to $\fd^l\phi^{s-l,k}$ and $\fd^l\psi^{s-l,k}$ and their first derivatives. However, for simplicity the reader is welcome to think of these quantities as simply $\phi^{s,k}$ and $\psi^{s,k}$ or even more simply as $\phi$ and $\psi$.

We begin with Lemma \ref{infty_base_lem}, which estimates an arbitrary function in $L^\infty(\Sigma_t)$ (which is the same space as $L^\infty(\tilde{\Sigma}_t)$) by certain Sobolev norms on $\Sigma_t$ and $\tilde{\Sigma}_t$. This lemma is then repeatedly applied to single derivatives of $\fd^l\phi^{s-l,k}$ and $\fd^l\psi^{s-l,k}$, resulting in Sobolev norms that can be estimated by the energy norms. (See Lemmas \ref{low_order_infty_lem}-\ref{gh_infty_lem}.)  These estimates are all summarized at the end in Proposition \ref{infinity_prop}.

\subsection{Sobolev-type estimates for $\Sigma_t$ and $\tilde{\Sigma}_t$}

First, we prove the following lemma, which is a Sobolev-type estimate. The weight that is gained in the following estimate depends on the volume form for the associated space, so $L^\infty$ estimates on $\Sigma_t$ gain two factors of $r$ and $L^\infty$ estimates on $\tilde{\Sigma}_t$ gain six factors of $r$. Also, the fact that the volume form for $\tilde{\Sigma}_t$ has additional factors of $\sin\theta$ means that more derivatives are required in the Sobolev estimate.

\begin{lemma}\label{infty_base_lem}
Let $u$ be an arbitrary function decaying sufficiently fast as $r\rightarrow\infty$. Then for any $r_0\ge r_H$,
$$||\fd^lu||^2_{L^\infty(\Sigma_t\cap\{r>r_0\})}\lesssim \int_{\Sigma_t\cap\{r>r_0\}}r^{-2}\left[(\pd_r\Gamma^{\le l+3}u)^2+(\Gamma^{\le l+3}\phi)^2\right]$$
and
$$||\fd^lu||^2_{L^\infty(\tilde\Sigma_t\cap\{r>r_0\})}\lesssim \int_{\tilde{\Sigma}_t\cap\{r>r_0\}}r^{-6}\left[(\pd_r\tilde{\Gamma}^{\le l+5} u)^2+(\tilde{\Gamma}^{\le l+5} u)^2\right].$$
If $l$ is even, the same results hold with only $\Gamma^{\le l+2}$ and $\tilde{\Gamma}^{\le l+4}$ respectively.
\end{lemma}
\begin{proof}
For a fixed $r$, denote by $\bar{u}:S^2(1)\rightarrow\R{}$ the pullback of the function $u:S^2(r)\rightarrow\R{}$ via the canonical map from $S^2(1)$ to $S^2(r)$. It is straightforward to show that
$$\overline{\fd^l u}=\fd^l\bar{u}$$
and
$$\overline{(q^2\sla\triangle)^lu}=\sla\triangle^l \bar{u}.$$
Also, denote by $d\omega$ the measure on $S^2(1)$.

\textbf{If $l$ is even, then}
\begin{multline*}
||\fd^lu||^2_{L^\infty(S^2(r))}=||\overline{\fd^l u}||^2_{L^\infty(S^2(1))}=||\fd^l\bar{u}||^2_{L^\infty(S^2(1))} \lesssim ||\fd^{\le l+2}\bar{u}||_{L^2(S^2(1))}^2 \\
\lesssim ||\sla\triangle^{\le (l+2)/2} \bar{u}||_{L^2(S^2(1))}^2 =\int_{S^2(r)}((q^2\sla\triangle)^{\le (l+2)/2}u)^2d\omega \lesssim \int_{S^2(r)}(\Gamma^{\le l+2}u)^2d\omega.
\end{multline*}
We have used Lemma \ref{spherical_infty_lem}, Theorem \ref{cl_embedding_thm}, and Lemma \ref{main_elliptic_lemma} in the three $\lesssim$ steps.

\textbf{If instead $l$ is odd, then}
$$||\fd^lu||_{L^\infty(S^2(r))}^2\lesssim ... \lesssim ||\fd^{\le l+3}\bar{u}||_{L^2(S^2(1))}^2 \lesssim ... \lesssim \int_{S^2(r)}\left(\Gamma^{\le l+3}u\right)^2d\omega.$$
That is, the calculation is the exact same as for the even case, except in the application of Lemma \ref{spherical_infty_lem}, which requires $\fd^{\le l+3}$ instead of $\fd^{\le l+2}$.

\textbf{Thus, in both cases},
$$||\fd^lu||_{L^\infty(S^2(r))}^2\lesssim \int_{S^2(r)}\left(\Gamma^{\le l+3} u\right)^2d\omega.$$
Now, set $f(r)=\int_{S^2(r)}(\Gamma^{\le l+3}u)^2d\omega$. Note that
\begin{multline*}
|f'(r)|\lesssim \int_{S^2(r)}|\Gamma^{\le l+3}u\pd_r\Gamma^{\le l+3}u|d\omega\lesssim \int_{S^2(r)}\left[(\pd_r \Gamma^{\le l+3}u)^2+(\Gamma^{\le l+3}u)^2\right]d\omega \\
\lesssim\int_{S^2(r)}r^{-2}\left[(\pd_r \Gamma^{\le l+3}u)^2+(\Gamma^{\le l+3}u)^2\right]q^2d\omega.
\end{multline*}
Then, assuming $\lim_{r\rightarrow\infty}f(r)=0$,
\begin{multline*}
|\fd^l u(r_0)|^2\lesssim f(r_0)\lesssim \int_{r_0}^\infty |f'(r)|dr\lesssim \int_{r_0}^\infty \int_{S^2(r)}r^{-2}\left[(\pd_r \Gamma^{\le l+3}u)^2+(\Gamma^{\le l+3}u)^2\right]q^2d\omega dr \\
\lesssim \int_{\Sigma_t\cap\{r\ge r_0\}}r^{-2}\left[(\pd_r \Gamma^{\le l+3}u)^2+(\Gamma^{\le l+3}u)^2\right].
\end{multline*}
The same procedure can be used to prove the analogous estimate on $\tilde\Sigma$, but there are two differences. The first is that there will be a loss of two extra derivatives since the $L^\infty$ estimate is now applied on $S^6$ instead of $S^2$. The second is that there will be a factor of $r^{-6}$ instead of a factor of $r^{-2}$, since the volume form for $\tilde{\Sigma}_t$ is $q^2A^2=O(r^6)$. For a more detailed explanation of this fact, see the proof of Lemma \ref{ws:infty_base_lem}. This explains the differences between the two estimates of the lemma.
\end{proof}

\subsection{Estimating derivatives using the Sobolev-type estimate}

Now, we repeatedly apply Lemma \ref{infty_base_lem} to estimate various derivatives with $r$ weights. We will assume that $\phi$ and $\psi$ decay sufficiently fast as $r\rightarrow\infty$

The following lemma estimates $\phi$ and $\psi$, as well as the higher order analogues $\fd^l\phi^{s-l,k}$ and $\fd^l\psi^{s-l,k}$.
\begin{lemma}\label{low_order_infty_lem}
For $r\ge r_H$,
$$|r^p\fd^l\phi^{s-l,k}|^2+|r^{p+2}\fd^l\psi^{s-l,k}|^2\lesssim E_{2p}^{s+5,k}(t)$$
and for $r\ge r_0>r_H$,
$$|\fd^l\phi^{s-l,k}|^2+|r^2\fd^l\psi^{s-l,k}|^2\lesssim E^{s+5,k}(t).$$
\end{lemma}
\begin{proof}
First, we apply Lemma \ref{infty_base_lem} with $u=r^p\phi^{s-l,k}$.
\begin{align*}
|r^p\fd^l\phi^{s-l,k}|^2 &=|\fd^l(r^p\phi^{s-l,k})|^2 \\
&\lesssim \int_{\Sigma_t}r^{-2}\left[(\pd_r\Gamma^{\le l+3}(r^p\phi^{s-l,k}))^2+(\Gamma^{\le l+3}(r^p\phi^{s-l,k}))^2\right] \\
&\lesssim \int_{\Sigma_t}r^{2p-2}\left[(\pd_r\Gamma^{\le l+3}\phi^{s-l,k})^2+(\Gamma^{\le l+3}\phi^{s-l,k})^2\right] \\
&\lesssim E_{2p}^{s+3}(t).
\end{align*}
Then, we apply Lemma \ref{infty_base_lem} with $u=r^{p+2}\psi^{s-l,k}$.
\begin{align*}
|r^{p+2}\fd^l\psi^{s-l,k}|^2 &=|\fd^l(r^{p+2}\psi^{s-l,k})|^2 \\
&\lesssim \int_{\tilde\Sigma_t}r^{-6}\left[(\pd_r\Gamma^{\le l+5}(r^{p+2}\psi^{s-l,k}))^2+(\Gamma^{\le l+5}(r^{p+2}\psi^{s-l,k}))^2\right] \\
&\lesssim \int_{\tilde\Sigma_t}r^{2p-2}\left[(\pd_r\Gamma^{\le l+5}\psi^{s-l,k})^2+(\Gamma^{\le l+5}\psi^{s-l,k})^2\right] \\
&\lesssim E_{2p}^{s+5}(t).
\end{align*}
Together, these estimates prove the first estimate of the lemma. The second estimate follows from the same exact argument in the special case $p=0$, and the observation that as long as $r\ge r_0>r_H$, then $E^{s,k}(t)$ can be used in place of $E_0^{s,k}(t)$.
\end{proof}

The following lemma estimates $\pd_t\phi$ and $\pd_t\psi$ as well as the higher order analogues $\pd_t\fd^l\phi^{s-l,k}$ and $\pd_t\fd^l\psi^{s-l,k}$.
\begin{lemma}\label{pd_t_infty_lem}
For $r\ge r_H$,
$$|r^p\pd_t\fd^l\phi^{s-l,k}|^2+|r^{p+2}\pd_t\fd^l\psi^{s-l,k}|^2 \lesssim E_{2p}^{s+6,k}(t)$$
and for $r\ge r_0>r_H$,
$$|r\pd_t\fd^l\phi^{s-l,k}|^2+|r^3\pd_t\fd^l\psi^{s-l,k}|^2 \lesssim E^{s+6,k}(t).$$
\end{lemma}
\begin{proof}
The first estimate reduces to Lemma \ref{low_order_infty_lem} by observing that $\pd_t\fd^l\phi^{s-l,k}=\fd^l\pd_t\phi^{s-l,k}=\fd^l\phi^{s+1-l,k}$ and likewise $\pd_t\fd^l\psi^{s-l,k}=\fd^l\psi^{s+1-l,k}$. We now prove the second estimate.

First, we apply Lemma \ref{infty_base_lem} with $u=r\pd_t\phi^{s-l,k}$.
\begin{align*}
|r\pd_t\fd^l\phi^{s-l,k}|^2 &= |\fd^l(r\pd_t\phi^{s-l,k})|^2 \\
&\lesssim \int_{\Sigma_t\cap\{r>r_0\}}r^{-2}\left[(\pd_r\Gamma^{\le l+3}(r\pd_t\phi^{s-l,k}))^2+(\Gamma^{\le l+3}(r\pd_t\phi^{s-l,k}))^2\right] \\
&\lesssim \int_{\Sigma_t\cap\{r>r_0\}}\left[(\pd_r\Gamma^{\le l+3}\pd_t\phi^{s-l,k})^2+(\Gamma^{\le l+3}\pd_t\phi^{s-l,k})^2\right] \\
&\lesssim \int_{\Sigma_t\cap\{r>r_0\}}\left[(\pd_r\Gamma^{\le l+3}\phi^{s+1-l,k})^2+(\pd_t\Gamma^{\le l+3}\phi^{s-l,k})^2\right] \\
&\lesssim E^{s+4,k}(t).
\end{align*}
Next, by applying Lemma \ref{infty_base_lem} with $u=r^3\pd_t\psi^{s-l,k}$ and repeating the same procedure, we arrive at the following estimate.
$$|r^3\pd_t\fd^l\psi^{s-l,k}|^2\lesssim E^{s+6,k}(t).$$
Together, these estimates prove the second estimate of the lemma.
\end{proof}

The following lemma estimates $(r^{-1}\pd_\theta)\phi$ and $(r^{-1}\pd_\theta)\psi$ as well as the higher order analogues $(r^{-1}\pd_\theta)\fd^l\phi^{s-l,k}$ and $(r^{-1}\pd_\theta)\fd^l\psi^{s-l,k}$.
\begin{lemma}\label{pd_theta_infty_lem}
For $r\ge r_H$,
$$|r^{p+1}(r^{-1}\pd_\theta)\fd^l\phi^{s-l,k}|^2+|r^{p+3}(r^{-1}\pd_\theta)\fd^l\psi^{s-l,k}|^2\lesssim E_{2p}^{s+6,k}(t)$$
and for $r\ge r_0>r_H$,
$$|r(r^{-1}\pd_\theta)\fd^l\phi^{s-l,k}|^2+|r^3(r^{-1}\pd_\theta)\fd^l\psi^{s-l,k}|^2\lesssim E^{s+6,k}(t).$$
\end{lemma}
\begin{proof}
This lemma reduces to Lemma \ref{low_order_infty_lem} by observing that
$$r^{p+1}(r^{-1}\pd_\theta)\fd^l\phi^{s-l,k}=r^p\fd\fd^l\phi^{s-l,k}=r^p\fd^{l+1}\phi^{s-l,k}\subset r^p\fd^l\phi^{s+1-l,k},$$
and likewise
$$r^{p+3}(r^{-1}\pd_\theta)\fd^l\psi^{s-l,k}\subset r^{p+2}\fd^l\psi^{s+1-l,k}.$$
\end{proof}

The following lemma estimates $L\phi$ and $L\psi$ as well as the higher order analogues $L\fd^l\phi^{s-l,k}$ and $L\fd^l\psi^{s-l,k}$.
\begin{lemma}\label{L_infty_lem}
Letting $L=\alpha\pd_r+\pd_t$, where $\alpha=\frac{\Delta}{r^2+a^2}$, we have that for $r\ge r_H$,
$$|r^{p+1}L\fd^l\phi^{s-l,k}|^2+|r^{p+3}L\fd^l\psi^{s-l,k}|^2\lesssim E_{2p}^{s+7,k}(t)+\int_{\Sigma_t}r^{2p}(\Box_g\phi^{s+3,k})^2+\int_{\tilde\Sigma_t}r^{2p}(\Box_g\psi^{s+5,k})^2$$
and for $r\ge r_0>r_H$,
$$|rL\fd^l\phi^{s-l,k}|^2+|r^3L\fd^l\psi^{s-l,k}|^2\lesssim E^{s+7,k}(t)+\int_{\Sigma_t\cap\{r>r_0\}}(\Box_g\phi^{s+3,k})^2+\int_{\tilde\Sigma_t\cap\{r>r_0\}}(\Box_g\psi^{s+5,k})^2.$$
\end{lemma}
\begin{proof}
We recall from the proof of Lemma \ref{k:L_pointwise_lem} that
$$(\pd_rLu)^2\lesssim (\Box_gu)^2+(L\pd_tu)^2+r^{-2}(\pd_t^2u)^2+r^{-2}(\pd_r\pd_tu)^2+r^{-2}(\pd_ru)^2+r^{-2}(Qu)^2,$$
By the same technique,
$$(\pd_rLu)^2\lesssim (\Box_{\tilde g}u)^2+(L\pd_tu)^2+r^{-2}(\pd_t^2u)^2+r^{-2}(\pd_r\pd_tu)^2+r^{-2}(\pd_ru)^2+r^{-2}(\tilde{Q}u)^2.$$

We now apply Lemma \ref{infty_base_lem} with $u=r^{p+1}L\phi^{s-l,k}$.
\begin{align*}
|r^{p+1}L\fd^l\phi^{s-l,k}|^2 &= |\fd^l(r^{p+1}L\phi^{s-l,k})|^2 \\
&\lesssim \int_{\Sigma_t}r^{-2}\left[(\pd_r\Gamma^{\le l+3}(r^{p+1}L\phi^{s-l,k}))^2+(\Gamma^{\le l+3}(r^{p+1}L\phi^{s-l,k}))^2\right] \\
&\lesssim \int_{\Sigma_t}r^{2p}\left[(\pd_r\Gamma^{\le l+3}L\phi^{s-l,k})^2+(\Gamma^{\le l+3}L\phi^{s-l,k})^2\right] \\
&\lesssim \int_{\Sigma_t}r^{2p}\left[(\pd_rL\phi^{s+3,k})^2+(L\phi^{s+3,k})^2\right] \\
&\lesssim E_{2p}^{s+3,k}(t) +\int_{\Sigma_t}r^{2p}(\pd_rL\phi^{s+3,k})^2.
\end{align*}
Now,
\begin{align*}
\int_{\Sigma_t}r^{2p}&(\pd_rL\phi^{s+3,k})^2 \\
&\lesssim \int_{\Sigma_t}r^{2p}\left[(\Box_g\phi^{s+3,k})^2+(L\pd_t\phi^{s+3,k})^2+r^{-2}(\pd_t^2\phi^{s+3,k})^2+r^{-2}(\pd_r\pd_t\phi^{s+3,k})^2\right.\\
&\hspace{3in}\left.+r^{-2}(\pd_r\phi^{s+3,k})^2+r^{-2}(Q\phi^{s+3,k})^2\right] \\
&\lesssim \int_{\Sigma_t}r^{2p}\left[(\Box_g\phi^{s+3,k})^2+(L\phi^{s+4,k})^2+r^{-2}(\phi^{s+5,k})^2+r^{-2}(\pd_r\phi^{s+4,k})^2\right].
\end{align*}
It follows that
$$|r^{p+1}L\fd^l\phi^{s-l,k}|^2\lesssim E_{2p}^{s+5,k}(t)+\int_{\Sigma_t}r^{2p}(\Box_g\phi^{s+3,k})^2.$$
By a similar argument,
$$|r^{p+3}L\fd^l\psi^{s-l,k}|^2\lesssim E_{2p}^{s+7,k}(t)+\int_{\tilde{\Sigma}_t}r^{2p}(\Box_{\tilde{g}}\psi^{s+5,k})^2.$$
These two estimates together establish the first estimate of the lemma. The second estimate follows from the same exact argument in the special case $p=0$, and the observation that as long as $r\ge r_0>r_H$, then $E^{s,k}(t)$ can be used in place of $E_0^{s,k}(t)$.
\end{proof}

The following lemma estimates $\tg\phi$ and $\tg\psi$ as well as the higher order analogues $\tg\fd^l\phi^{s-l,k}$ and $\tg\fd^l\psi^{s-l,k}$.
\begin{lemma}\label{gh_infty_lem}
Keeping in mind that $\tg$ is supported in a neighborhood of the event horizon, for arbitrary $p'$, we have for $r\ge r_H$,
$$|\tg \fd^l\phi^{s-l,k}|^2+|\tg \fd^l\psi^{s-l,k}|^2 \lesssim E_{p'}^{s+5,k+1}(t)$$
and for $r\ge r_0>r_H$,
$$|\tg \fd^l\phi^{s-l,k}|^2+|\tg \fd^l\psi^{s-l,k}|^2 \lesssim E^{s+5,k+1}(t).$$
\end{lemma}
\begin{proof}
We apply Lemma \ref{infty_base_lem} with $u=\tg\phi^{s-l,k}$, and freely introduce a factor of $r^{p'}$ since $\tg$ is supported on a compact interval in $r$.
\begin{align*}
|\tg\fd^l\phi^{s-l,k}|^2 &= |\fd^l\phi^{s-l,k+1}|^2 \\
&\lesssim \int_{\Sigma_t}r^{-2}\left[(\pd_r\Gamma^{\le l+3}\phi^{s-l,k+1})^2+(\Gamma^{\le l+3}\phi^{s-l,k+1})^2\right] \\
&\lesssim \int_{\Sigma_t}r^{p'-2}\left[(\pd_r\Gamma^{\le l+3}\phi^{s-l,k+1})^2+(\Gamma^{\le l+3}\phi^{s-l,k+1})^2\right] \\
&\lesssim E_{p'}^{s+3,k+1}(t).
\end{align*}
A similar argument shows that
$$|\tg \fd^l\psi^{s-l,k}|^2\lesssim E_{p'}^{s+5,k+1}(t).$$
Together, these estimates prove the first estimate of the lemma. The second estimate follows from the same argument, and the observation that as long as $r\ge r_0>r_H$, then $E^{s,k}(t)$ can be used in place of $E_{p'}^{s,k}(t)$.
\end{proof}

\subsection{Summarizing the $L^\infty$ estimates}
To conclude this section, we summarize the previous lemmas in a single proposition, making use of the following definition.
\begin{definition} We define two families of operators.
$$\bar{D}=\{L,r^{-1}\pd_\theta\},$$
$$D=\{L,\pd_r,r^{-1}\pd_\theta\}.$$
\end{definition}

\begin{proposition}\label{infinity_prop}
For $r\ge r_H$,
\begin{multline*}
|r^{p+1}\bar{D}\fd^l\phi^{s-l,k}|^2+|r^pD\fd^l\phi^{s-l,k}|^2+|r^p\fd^l\phi^{s-l,k}|^2 \\
+|r^{p+3}\bar{D}\fd^l\psi^{s-l,k}|^2+|r^{p+2}D\fd^l\psi^{s-l,k}|^2+|r^{p+2}\fd^l\psi^{s-l,k}|^2 \\
\lesssim E_{2p}^{s+5,k+1}(t)+E_{2p}^{s+7,k}(t)+\int_{\Sigma_t}r^{2p}(\Box_g\phi^{s+5,k})^2+\int_{\tilde\Sigma_t}r^{2p}(\Box_{\tilde{g}}\psi^{s+5,k})^2
\end{multline*}
and for $r\ge r_0>r_H$,
\begin{multline*}
|rD\fd^l\phi^{s-l,k}|^2+|r^3D\fd^l\psi^{s-l,k}|^2 \\
\lesssim E^{s+5,k+1}(t)+E^{s+7,k}(t)+\int_{\Sigma_t\cap\{r>r_0\}}(\Box_g\phi^{s+5,k})^2+\int_{\tilde\Sigma_t\cap\{r>r_0\}}(\Box_{\tilde{g}}\psi^{s+5,k})^2.
\end{multline*}
\end{proposition}
\begin{proof}
With the exception of the operator $\pd_r$, all of the cases have been proved in Lemmas \ref{low_order_infty_lem}-\ref{gh_infty_lem}. Finally, observe that
$$|r^p\pd_r\fd^l\phi^{s-l,k}|^2\lesssim |r^pL\fd^l\phi^{s-l,k}|^2+|r^p\pd_t\fd^l\phi^{s-l,k}|^2+|\tg\fd^l\phi^{s-l,k}|^2$$
and
$$|r^p\pd_r\fd^l\phi^{s-l,k}|^2\lesssim |r^pL\fd^l\phi^{s-l,k}|^2+|r^p\pd_t\fd^l\phi^{s-l,k}|^2+|\tg\fd^l\phi^{s-l,k}|^2.$$
Thus, even the case of the operator $\pd_r$ can be reduced to Lemmas \ref{low_order_infty_lem}-\ref{gh_infty_lem}.
\end{proof}

\section{The structure of the nonlinear terms}\label{wk:structure_sec}

In this section, we carefully examine the nonlinear terms $\mathcal{N}_\phi$ and $\mathcal{N}_\psi$ as well as their higher order analogues $q^{-2}\tg^k\Gamma^s(q^2\mathcal{N}_\phi)$ and $q^{-2}\tg^k\tilde\Gamma^s(q^2\mathcal{N}_\psi)$ and determine a procedure for estimating them in the proof of the main theorem. In particular, we note two important structural conditions that these terms satisfy. The first is the well-known \textit{null condition} and the second is a condition on the axis that allows us to use the formalism provided in Appendix \ref{regularity_sec}.

We start in \S\ref{wk:phi_psi_derivation_sec} by deriving the equations for $(\phi,\psi)$. Then in \S\ref{nl_structure_sec} we define and prove the precise structures of $\mathcal{N}_\phi$ and $\mathcal{N}_\psi$. These structures will then be used in \S\ref{nl_s_k_structure_sec} to define and prove the precise structures of $q^{-2}\tg^k\Gamma^s(q^2\mathcal{N}_\phi)$ and $q^{-2}\tg^k\tilde\Gamma^s(q^2\mathcal{N}_\psi)$. Finally, in \S\ref{nl_strategy_revisited_sec} we prove a proposition that uses these structures to estimate nonlinear quantities that will show up in the proof of the main theorem (Theorem \ref{main_thm}).

\subsection{Deriving the equations for $(\phi,\psi)$}\label{wk:phi_psi_derivation_sec}

In order to proceed, we first calculate the nonlinear terms $\mathcal{N}_\phi$ and $\mathcal{N}_\psi$.
\begin{proposition}\label{N_identities_prop}
If one makes the substitutions
\begin{align*}
X&= A+A\phi \\
Y&=B+A^2\psi
\end{align*}
and requires that $\phi$ and $\psi$ are axisymmetric functions, then the wave map system (\ref{wk:X_eqn}-\ref{wk:Y_eqn}) reduces to the following system of equations for $\phi$ and $\psi$.
$$\Box_g \phi = \mathcal{L}_\phi + \mathcal{N}_\phi$$
$$\Box_{\tilde{g}} \psi = \mathcal{L}_\psi + \mathcal{N}_\psi,$$
where
$$\mathcal{L}_\phi=-2\frac{\pd^\alpha B}{A}A\pd_\alpha \psi + 2\frac{\pd^\alpha B\pd_\alpha B}{A^2}\phi-4\frac{\pd^\alpha A\pd_\alpha B}{A^2} A\psi$$
$$\mathcal{L}_\psi=-2\frac{\pd^\alpha A_2}{A_2}\pd_\alpha\psi+2\frac{\pd^\alpha B\pd_\alpha B}{A^2}\psi + 2A^{-1}\frac{\pd^\alpha B}{A}\pd_\alpha\phi$$
are the linear terms that first appeared in \S\ref{wk:phi_psi_equations_intro_sec}, and
\begin{align*}
(1+\phi)\mathcal{N}_\phi &= \pd^\alpha \phi \pd_\alpha \phi - A\pd^\alpha \psi A\pd_\alpha \psi +2\frac{\pd^\alpha B}{A} \phi A\pd_\alpha \psi -4\frac{\pd^\alpha A}{A}A\psi A\pd_\alpha \psi \\
&\hspace{.5in}-\frac{\pd^\alpha B\pd_\alpha B}{A^2}\phi^2+4\frac{\pd^\alpha A\pd_\alpha B}{A^2}\phi A\psi -4\frac{\pd^\alpha A\pd_\alpha A}{A^2}(A\psi)^2
\end{align*}
$$(1+\phi)\mathcal{N}_\psi = 2\pd^\alpha \phi \pd_\alpha\psi  +4\frac{\pd^\alpha A}{A}\psi\pd_\alpha \phi -2\frac{\pd^\alpha B}{A}A^{-1}\phi\pd_\alpha\phi$$
are the nonlinear terms.
\end{proposition}
\begin{proof}
The first equation of the wave map system (\ref{wk:X_eqn}) is
$$\Box_gX = \frac{\pd^\alpha X\pd_\alpha X}{X}-\frac{\pd^\alpha Y\pd_\alpha Y}{X}.$$
We substitute $X=A(1+\phi)$ and $Y=B+A^2\psi$.
$$\Box_g(A(1+\phi))=\frac{\pd^\alpha (A(1+\phi))\pd_\alpha (A(1+\phi))}{A(1+\phi)}-\frac{\pd^\alpha (B+A^2\psi)\pd_\alpha (B+A^2\psi)}{A(1+\phi)}.$$
Now, we expand each term as follows.
$$\Box_g(A(1+\phi)) = (1+\phi)\Box_gA +2\pd^\alpha A\pd_\alpha \phi+ A\Box_g\phi$$
$$\frac{\pd^\alpha (A(1+\phi))\pd_\alpha (A(1+\phi))}{A(1+\phi)} = (1+\phi)\frac{\pd^\alpha A\pd_\alpha A}{A} + 2\pd^\alpha A\pd_\alpha\phi+A\frac{\pd^\alpha \phi\pd_\alpha\phi}{(1+\phi)}$$
\begin{multline*}
-\frac{\pd^\alpha (B+A^2\psi)\pd_\alpha (B+A^2\psi)}{A(1+\phi)} = -(1+\phi)\frac{\pd^\alpha B\pd_\alpha B}{A}-\left(\frac{1}{1+\phi}-(1+\phi)\right)\frac{\pd^\alpha B\pd_\alpha B}{A} \\
-2\frac{\pd^\alpha B\pd_\alpha (A^2\psi)}{A(1+\phi)}-\frac{\pd^\alpha (A^2\psi)\pd_\alpha(A^2\psi)}{A(1+\phi)}
\end{multline*}
Using the fact that $(X,Y)=(A,B)$ also solves equation (\ref{wk:X_eqn}), we determine that
$$A\Box_g\phi = A\frac{\pd^\alpha \phi\pd_\alpha\phi}{(1+\phi)} -\left(\frac{1}{1+\phi}-(1+\phi)\right)\frac{\pd^\alpha B\pd_\alpha B}{A} -2\frac{\pd^\alpha B\pd_\alpha (A^2\psi)}{A(1+\phi)}-\frac{\pd^\alpha (A^2\psi)\pd_\alpha(A^2\psi)}{A(1+\phi)}.$$
Note the following identities.
$$-\frac{1}{1+\phi}+1+\phi = \frac{-1+(1+\phi)^2}{1+\phi}=\frac{2\phi+\phi^2}{1+\phi}=\frac{\phi(2(1+\phi)-\phi)}{1+\phi}=2\phi-\frac{\phi^2}{1+\phi}$$
$$\frac{1}{1+\phi} = 1+\left(\frac1{1+\phi}-1\right) = 1+\frac{1-(1+\phi)}{1+\phi} = 1-\frac{\phi}{1+\phi}.$$
We conclude that
\begin{align*}
\mathcal{L}_\phi &=\frac1A\left[ 2\phi\frac{\pd^\alpha B\pd_\alpha B}{A}-2\frac{\pd^\alpha B\pd_\alpha (A^2\psi)}{A}\right] \\
&= 2\frac{\pd^\alpha B\pd_\alpha B}{A^2}\phi -4\frac{\pd^\alpha A\pd_\alpha B}{A^2}A\psi  -2\frac{\pd^\alpha B}{A}A\pd_\alpha \psi
\end{align*}
and
\begin{align*}
(1+\phi)\mathcal{N}_\phi &= \frac{1+\phi}{A}\left[ A\frac{\pd^\alpha\phi\pd_\alpha\phi}{1+\phi}-\frac{\phi^2}{1+\phi}\frac{\pd^\alpha B\pd_\alpha B}{A}+2\frac{\phi}{1+\phi}\frac{\pd^\alpha B\pd_\alpha (A^2\psi)}{A}\right. \\
&\hspace{3.2in}\left.-\frac{\pd^\alpha(A^2\psi)\pd_\alpha (A^2\psi)}{A(1+\phi)} \right]  \\
&= \pd^\alpha\phi\pd_\alpha\phi -\frac{\pd^\alpha B\pd_\alpha B}{A^2}\phi^2 +4\frac{\pd^\alpha A\pd_\alpha B}{A^2}\phi A\psi+2\frac{\pd^\alpha B}{A}\phi A\pd_\alpha\psi \\
&\hspace{1in}-4\frac{\pd^\alpha A\pd_\alpha A}{A^2}(A\psi)^2-4\frac{\pd^\alpha A}{A}A\psi A\pd_\alpha\psi-A\pd^\alpha \psi A\pd_\alpha \psi.
\end{align*}

The second equation of the wave map system (\ref{wk:Y_eqn}) is
$$\Box_gY = 2\frac{\pd^\alpha X\pd_\alpha Y}{X}.$$
We substitute $X=A(1+\phi)$ and $Y=B+A^2\psi$.
$$\Box_g(B+A^2\psi) =2\frac{\pd^\alpha (A(1+\phi))\pd_\alpha (B+A^2\psi)}{A(1+\phi)}.$$
The left side simplifies as follows (again using the equation for $\Box_gA$).
\begin{align*}
\Box_g(B+A^2\psi) &= \Box_g B+2A\psi\Box_gA+2\pd^\alpha A\pd_\alpha A\psi +4A\pd^\alpha A\pd_\alpha \psi+A^2\Box_g\psi \\
&= \Box_g B-2\pd^\alpha B\pd_\alpha B\psi+4\pd^\alpha A\pd_\alpha A\psi +4A\pd^\alpha A\pd_\alpha \psi+A^2\Box_g\psi
\end{align*}
The right side simplifies as follows.
\begin{multline*}
2\frac{\pd^\alpha (A(1+\phi))\pd_\alpha (B+A^2\psi)}{A(1+\phi)} \\
= 2\frac{\pd^\alpha A\pd_\alpha B}{A}+2\frac{\pd^\alpha\phi\pd_\alpha B}{(1+\phi)}+2\frac{\pd^\alpha A\pd_\alpha(A^2\psi)}{A}+2\frac{\pd^\alpha\phi\pd_\alpha(A^2\psi)}{1+\phi} \\
= 2\frac{\pd^\alpha A\pd_\alpha B}{A}+2\frac{\pd^\alpha\phi\pd_\alpha B}{(1+\phi)}+4\pd^\alpha A\pd_\alpha A\psi+2A\pd^\alpha A \pd_\alpha \psi+2\frac{\pd^\alpha\phi\pd_\alpha(A^2\psi)}{1+\phi} 
\end{multline*}
Note that the term $4\pd^\alpha A\pd_\alpha A\psi$ appears on both sides and therefore cancels out. \textbf{This cancellation is the motivation for the choice of linearization $Y=B+A^2\psi$.} Using the fact that $B$ solves equation (\ref{wk:Y_eqn}), we conclude that
$$A^2\Box_g\psi +2A\pd^\alpha A\pd_\alpha \psi = 2\pd^\alpha B\pd_\alpha B \psi +2\frac{\pd^\alpha\phi\pd_\alpha B}{(1+\phi)}+2\frac{\pd^\alpha\phi\pd_\alpha(A^2\psi)}{1+\phi}.$$
Dividing by $A^2$ and using the fact that
$$\Box_g+2\frac{\pd^\alpha A}{A}\pd_\alpha = \Box_g+2\frac{\pd^\alpha A_1}{A_1}\pd_\alpha+2\frac{\pd^\alpha A_2}{A_2}\pd_\alpha = \Box_{\tilde{g}}+2\frac{\pd^\alpha A_2}{A_2}\pd_\alpha,$$
we obtain
$$\Box_{\tilde{g}}\psi+2\frac{\pd^\alpha A_2}{A_2}\pd_\alpha \psi = 2\frac{\pd^\alpha B\pd_\alpha B}{A^2}\psi +\frac{2}{1+\phi}A^{-1}\frac{\pd^\alpha B}{A}\pd_\alpha\phi +\frac{4}{1+\phi} \frac{\pd^\alpha A}{A}\psi \pd_\alpha\phi +\frac{2}{1+\phi}\pd^\alpha\phi \pd_\alpha\psi.$$
Again, since
$$\frac{1}{1+\phi} = 1-\frac{\phi}{1+\phi},$$
we conclude that
$$\mathcal{L}_\psi = -\frac{\pd^\alpha A_2}{A_2}\pd_\alpha \psi +2A^{-1}\frac{\pd^\alpha B}{A}\pd_\alpha\phi+2\frac{\pd^\alpha B\pd_\alpha B}{A^2}\psi$$
and
$$(1+\phi)\mathcal{N}_\psi = -2\frac{\pd^\alpha B}{A}A^{-1}\phi\pd_\alpha\phi +4\frac{\pd^\alpha A}{A}\psi\pd_\alpha \phi +2\pd^\alpha\phi\pd_\alpha\psi.$$
\end{proof}

\subsection{The structures of $\mathcal{N}_\phi$ and $\mathcal{N}_\psi$}\label{nl_structure_sec}

We now categorize each term in $\mathcal{N}_\phi$ and $\mathcal{N}_\psi$, beginning with the definition of the null condition.

\begin{definition} We define two families of terms.
$$\gamma = \{L\phi,r^{-1}\pd_r\phi,r^{-1}\pd_\theta\phi,r^{-1}\phi,r^2L\psi,r\pd_r\psi,r\pd_\theta\psi,r\psi\}$$
$$\beta = \{L\phi,\pd_r\phi,r^{-1}\pd_\theta\phi,r^{-1}\phi,r^2L\psi,r^2\pd_r\psi,r\pd_\theta\psi,r\psi\}$$
The null condition states that any nonlinear term must be a product with at least one $\gamma$ factor.
\end{definition}

We now prove two lemmas (one for $\mathcal{N}_\phi$ and one for $\mathcal{N}_\psi$) which guarantee the null condition as well as a structural condition on the axis.

\begin{lemma}\label{Nphi_structure_lem}(Structure of $\mathcal{N}_\phi$)
The nonlinear term $\mathcal{N}_\phi$ can be expressed as a sum of terms of the form
$$\frac{f\gamma\beta}{1+\phi},$$
satisfying the following additional rules. \\
i) The factor $f$ is smooth and bounded. \\
ii) In the language of Appendix \ref{regularity_sec}, the term belongs to 
$$\tau_{(\le 1)}(P_\theta^\infty\cup \{\phi,\pd_t\phi,\pd_r\phi,(1+\phi)^{-1},\psi,\pd_t\psi,\pd_r\psi\}).$$
iii) If $\psi$ appears at least once in the term, then $f$ has a factor of $\sin^2\theta$.
\end{lemma}
\begin{proof}
In this proof, we use the sign $\approx$ to emphasize that smooth and bounded factors (including $M/r$ and $a/r$) are neglected. One can check (using Lemma \ref{small_a_quantities_lem} as a guide) that
$$ \frac{\pd_\theta A}{A} \approx \frac1{\sin\theta}, \hspace{.5in}
\frac{\pd_r A}{A} \approx r^{-1}, \hspace{.5in}
\frac{\pd_\theta B}{A} \approx \sin\theta, \hspace{.5in}
\frac{\pd_r B}{A} \approx r^{-1}\sin^2\theta $$
$$ \frac{\pd^\alpha A\pd_\alpha A}{A^2} \approx \frac1{r^2\sin^2\theta}, \hspace{.5in}
\frac{\pd^\alpha A\pd_\alpha B}{A^2} \approx r^{-2}, \hspace{.5in}
\frac{\pd^\alpha B\pd_\alpha B}{A^2} \approx r^{-2}\sin^2\theta. $$
Using these, we investigate each term in $(1+\phi)\mathcal{N}_\phi$.
\begin{align*}
\pd^\alpha \phi\pd_\alpha \phi &\approx L\phi\lbar \phi + r^{-2}\pd_\theta\phi\pd_\theta \phi \\
&\approx (L\phi)(\lbar\phi)+(r^{-1}\pd_\theta\phi)(r^{-1}\pd_\theta\phi) \\
&\approx (L\phi)(\pd_r\phi)+(L\phi)(L\phi)+(r^{-1}\pd_\theta\phi)(r^{-1}\pd_\theta\phi).
\end{align*}
\begin{align*}
A\pd^\alpha\psi A\pd_\alpha \psi &\approx r^4\sin^4\theta (L\psi \lbar\psi+r^{-2}\pd_\theta\psi\pd_\theta\psi) \\
&\approx \sin^4\theta (r^2L\psi)(r^2\lbar\psi)+\sin^4\theta(r\pd_\theta\psi)(r\pd_\theta\psi) \\
&\approx \sin^4\theta (r^2L\psi)(r^2\pd_r\psi)+\sin^4\theta (r^2L\psi)(r^2L\psi)+\sin^4\theta(r\pd_\theta\psi)(r\pd_\theta\psi).
\end{align*}
\begin{align*}
\frac{\pd^\alpha B}{A}\phi A\pd_\alpha \psi &\approx r^{-2}\frac{\pd_\theta B}{A}\phi A\pd_\theta\psi +\frac{\pd_r B}{A}\phi A\pd_r\psi \\
&\approx \sin^3\theta \phi\pd_\theta\psi +r^{-1}\sin^4\theta \phi\pd_r\psi \\
&\approx \sin^3\theta (r^{-1}\phi)(r\pd_\theta\psi) + \sin^4\theta (r^{-1}\phi)(\pd_r\psi).
\end{align*}
\begin{align*}
\frac{\pd^\alpha A}{A} A\psi A\pd_\alpha \psi &\approx r^{-2}\frac{\pd_\theta A}{A}A\psi A\pd_\theta\psi +\frac{\pd_rA}{A}A\psi A\pd_r\psi \\
&\approx r^2\sin^3\theta\psi\pd_\theta\psi + r^3\sin^4\theta\psi\pd_r\psi \\
&\approx \sin^3\theta (r\psi)(r\pd_\theta\psi) +\sin^4\theta (r\psi)(r^2\pd_r\psi).
\end{align*}
$$\frac{\pd^\alpha B\pd_\alpha B}{A^2}\phi^2 \approx r^{-2}\sin^2\theta\phi^2 \approx \sin^2\theta (r^{-1}\phi)(r^{-1}\phi).$$
$$\frac{\pd^\alpha A\pd_\alpha A}{A^2}\phi A\psi \approx \sin^2\theta\phi\psi \approx \sin^2\theta(r^{-1}\phi)(r\psi).$$
$$\frac{\pd^\alpha A\pd_\alpha A}{A^2}A\psi A\psi \approx r^2\sin^2\theta\psi^2 \approx \sin^2\theta (r\psi)(r\psi).$$
It is now straightforward to check for each of the above calculations that the first factor in parentheses is a $\gamma$ term, and the second factor in parentheses is a $\beta$ term. It is also straightforward to check that any term containing a $\psi$ factor also has a factor of $\sin^2\theta$. 
\end{proof}

\begin{lemma}\label{wk:structure_Npsi_0_lem}(Structure of $\mathcal{N}_\psi$)
The nonlinear term $\mathcal{N}_\psi$ can be expressed as a sum of terms of the form
$$\frac{r^{-2}f\gamma\beta}{1+\phi}\text{ or }\frac{r^{-2}f\gamma (r^{-1}\fb \phi)}{1+\phi},$$
satisfying the following additional rules. \\
i) The factor $f$ is smooth and bounded. \\
ii) In the language of Appendix \ref{regularity_sec}, the term belongs to 
$$\tau_{(\le 1)}(P_\theta^\infty\cup \{\phi,\pd_t\phi,\pd_r\phi,(1+\phi)^{-1},\psi,\pd_t\psi,\pd_r\psi\}).$$
\end{lemma}
\begin{proof}
We proceed as in the proof of the previous lemma, this time investigating each term in $(1+\phi)\mathcal{N}_\psi$.
\begin{align*}
\pd^\alpha\phi\pd_\alpha\psi &\approx L\phi\lbar\psi +\lbar\phi L\psi +r^{-2}\pd_\theta\phi\pd_\theta\psi \\
&\approx r^{-2}(L\phi)(r^2\lbar\psi)+r^{-2}(r^2L\psi)(\lbar\phi)+r^{-2}(r^{-1}\pd_\theta\phi)(r\pd_\theta\psi) \\
&\approx r^{-2}(L\phi)(r^2\pd_r\psi)+r^{-2}(L\phi)(r^2L\psi)+r^{-2}(r^2L\psi)(\pd_r\phi)+r^{-2}(r^2L\psi)(L\phi) \\
&\hspace{4in}+r^{-2}(r^{-1}\pd_\theta\phi)(r\pd_\theta\psi).
\end{align*}
\begin{align*}
\frac{\pd^\alpha A}{A}\psi\pd_\alpha \phi &\approx r^{-2}\frac{\pd_\theta A}{A}\psi\pd_\theta\phi + \frac{\pd_r A}{A}\psi\pd_r\phi \\
&\approx r^{-2}\psi\fb\phi +r^{-1}\psi\pd_r\phi \\
&\approx r^{-2}(r\psi)(r^{-1}\fb\phi) + r^{-2}(r\psi)(\pd_r\phi).
\end{align*}
\begin{align*}
\frac{\pd^\alpha B}{A^2}\phi\pd_\alpha \phi &\approx r^{-2}\frac{\pd_\theta B}{A^2}\phi\pd_\theta\phi +\frac{\pd_r B}{A^2}\phi\pd_r\phi \\
&\approx r^{-4}\phi \fb\phi +r^{-4}\sin\theta \phi\pd_\theta\phi +r^{-3}\phi\pd_r\phi \\
&\approx r^{-2}(r^{-1}\phi)(r^{-1}\fb\phi)+r^{-2}\sin\theta(r^{-1}\phi)(r^{-1}\pd_\theta\phi)+r^{-2}(r^{-1}\phi)(\pd_r\phi).
\end{align*}
It is now straightforward to check for each of the above calculations that the first factor in parentheses is a $\gamma$ term, and the second factor in parentheses is either a $\beta$ term or $r^{-1}\fb\phi$. Finally, every term has an additional factor of $r^{-2}$.
\end{proof}

\subsection{The structures of $q^{-2}\tg^k\Gamma^s(q^2\mathcal{N}_\phi)$ and $q^{-2}\tg^k\tilde\Gamma^s(q^2\mathcal{N}_\psi)$}\label{nl_s_k_structure_sec}

We now generalize the previous two lemmas by applying the commutators.

\begin{definition}\label{wk:gamma_beta_s_k_def} We generalize the previous families of terms to the following.
\begin{multline*}
\gamma^{s,k} = \{L\fd^l\phi^{s-l,k},r^{-1}\pd_r\fd^l\phi^{s-l,k},r^{-1}\pd_\theta\fd^l\phi^{s-l,k},r^{-1}\fd^l\phi^{s-l,k}, \\
r^2L\fd^l\psi^{s-l,k},r\pd_r\fd^l\psi^{s-l,k},r\pd_\theta\fd^l\psi^{s-l,k},r\fd^l\psi^{s-l,k}\}
\end{multline*}
\begin{multline*}
\beta^{s,k} = \{L\fd^l\phi^{s-l,k},\pd_r\fd^l\phi^{s-l,k},r^{-1}\pd_\theta\fd^l\phi^{s-l,k},r^{-1}\fd^l\phi^{s-l,k}, \\
r^2L\fd^l\psi^{s-l,k},r^2\pd_r\fd^l\psi^{s-l,k},r\pd_\theta\fd^l\psi^{s-l,k},r\fd^l\psi^{s-l,k}\}
\end{multline*}
The null condition still states that any nonlinear term must be a product with at least one $\gamma$ factor.
\end{definition}


The following two lemmas (one for $q^{-2}\tg^k\Gamma^s(q^2\mathcal{N}_\phi)$ and one for $q^{-2}\tg^k\tilde{\Gamma}^s(q^2\mathcal{N}_\psi)$) generalize the lemmas in \S\ref{nl_structure_sec}.

\begin{lemma}\label{Nphi_s_k_structure_lem}
(Structure of $q^{-2}\tg^k\Gamma^s(q^2\mathcal{N}_\phi)$) The nonlinear term $q^{-2}\tg^k\Gamma^s(q^2\mathcal{N}_\phi)$ can be expressed as a sum of terms of the following form. (The integer $j$ represents the number of times a differential operator acts on the denominator.) \\
$$\frac{f\gamma^{s_1,k_1}\beta^{s_2,k_2}(r\beta^{s_3,k_3})...(r\beta^{s_{2+j},k_{2+j}})}{(1+\phi)^{j+1}},$$
where $0\le j\le s+k$, $s_1+...+s_{2+j}\le s+2$, $\max_is_i\le s$, $k_1+...+k_{2+j}\le k$, and the following rules apply. \\
\bp The factor $f$ is smooth and bounded. \\
\bp If $\psi$ appears at least once in the term, then $f$ has a factor of $(\sin\theta)^{\max(0,2-\min_i(s-s_i))}$.
\end{lemma}
\begin{proof}
We start with the case $s=k=0$, which was proved in Lemma \ref{Nphi_structure_lem}. In particular, since $r^{-1}\pd_\theta\phi=\pd_\theta (r^{-1}\phi)$ and $r\pd_\theta\psi=\pd_\theta (r\psi)$, each term in $\mathcal{N}_\phi$ can be written in one of the following forms.
$$\frac{f\gamma_0\beta_0}{1+\phi}\text{ or }
\frac{f_1\pd_\theta f_2\pd_\theta\gamma_0\beta_0}{1+\phi}\text{ or }
\frac{f_1\pd_\theta f_2\gamma_0\pd_\theta\beta_0}{1+\phi}\text{ or }
\frac{f\pd_\theta\gamma_0\pd_\theta\beta_0}{1+\phi},$$
where
$$\gamma_0\in \gamma\setminus \{r^{-1}\pd_\theta\phi,r\pd_\theta\psi\}\text{ and }
\beta_0\in \beta\setminus \{r^{-1}\pd_\theta\phi,r\pd_\theta\psi\}.$$
The effect of applying $\pd_t^{s'}$ to any of these terms is to obtain terms of the form
$$\frac{f\pd_t^{s_1}\gamma_0\pd_t^{s_2}\beta_0\pd_t^{s_3}\phi...\pd_t^{s_{2+j}}\phi}{(1+\phi)^{1+j}}\text{ or }
\frac{f_1\pd_\theta f_2\pd_\theta\pd_t^{s_1}\gamma_0\pd_t^{s_2}\beta_0\pd_t^{s_3}\phi...\pd_t^{s_{2+j}}\phi}{(1+\phi)^{1+j}}\text{ etc.},$$
where $s_1+...+s_{2+j}\le s'$ and $j\le s'$. The additional factors of $\pd_t^{s_i}\phi$ appear each time one of the $\pd_t$ operators acts on the denominator.

Recall that $\tg^k\Gamma^s$ is composed not only of $\pd_t$, but also of $Q$ and $\tg$. Since these operators commute with each other, the order in which they are applied is not important. For the sake of simplicity, we first apply the $\pd_t^{s'}$ operators (which has already been done) and then the $\tg^k$ operators. Since these are both first order operators, it should be clear that the resulting terms are of the form
$$\frac{f\gamma_0^{s_1,k_1}\beta_0^{s_2,k_2}\eta^{s_3,k_3}...\eta^{s_{2+j},k_{2+j}}}{(1+\phi)^{1+j}}\text{ or }
\frac{f_1\pd_\theta f_2\pd_\theta\gamma_0^{s_1,k_1}\beta_0^{s_2,k_2}\eta^{s_3,k_3}...\eta^{s_{2+j},k_{2+j}}}{(1+\phi)^{1+j}}\text{ etc.},$$
where $j\le s'+k$, $s_1+...+s_{2+j}\le s'$, $k_1+...+k_{2+j}\le k$, and $\gamma_0^{s_1,k_1}=\tg^{k_1}\pd_t^{s_1}\gamma_0$, $\beta_0^{s_2,k_2}=\tg^{k_2}\pd_t^{s_2}\beta_0$, and $\eta^{s_i,k_i}=\tg^{k_i}\pd_t^{s_i}\phi$.
This is the point at which we apply the operator $Q$. Recall that
$$Q=\fa+\fb+a^2\sin^2\theta\pd_t^2,$$
where the operators $\fa$ and $\fb$ are defined in Appendix \ref{regularity_sec}. We handle $Q$ as if applying these operators seperately. The operator $a^2\sin^2\theta\pd_t^2$ is equivalent to applying $\pd_t$ twice and multiplying by a smooth bounded function. However, the operators $\fa$ and $\fb$ require Lemma \ref{cl_raise_degree_lem}. After applying $Q^{s''}$, the terms should take the following form.
$$
\frac{f_1\fd^{l_0} f_2 \fd^{l_1}\gamma_0^{s_1,k_1}\fd^{l_2}\beta_0^{s_2,k_2}\fd^{l_3}\eta^{s_3,k_3}...\fd^{l_{2+j}}\eta^{s_{2+j},k_{2+j}}}{(1+\phi)^{1+j}}
$$
where $l_0=1$ or $\fd^{l_0}f_2=1$ so that $l_0+l_1+...l_{2+j}$ is even. We have in addition that $l_0+l_1+...+l_{2+j}\le s''+2$, since there were initially zero or two $\fd$s and each application of $Q$ contributes at most two more. Finally, we observe that the $\eta$ factors are all of the form $r\beta$. From here, the lemma follows.
\end{proof}

\begin{lemma}\label{Npsi_s_k_structure_lem}
(Structure of $q^{-2}\tg^k\tilde\Gamma^s(q^2\mathcal{N}_\psi)$) The nonlinear term $q^{-2}\tg^k\tilde\Gamma^s(q^2\mathcal{N}_\psi)$ can be expressed as a sum of terms, each taking one of the following forms. (The integer $j$ represents the number of times a differential operator acts on the denominator.) \\
$$\frac{r^{-2}f\gamma^{s_1,k_1}\beta^{s_2,k_2}(r\beta^{s_3,k_3})...(r\beta^{s_{2+j},k_{2+j}})}{(1+\phi)^{j+1}}$$
or
$$\frac{r^{-2}f\gamma^{s_1,k_1}(r^{-1}\fb\fc^l\phi^{s_2-2l,k_2})(r\beta^{s_3,k_3})...(r\beta^{s_{2+j},k_{2+j}})}{(1+\phi)^{j+1}},$$
where $0\le j\le s+k$, $s_1+...+s_{2+j}\le s+2$, $\max_is_i\le s$, $k_1+...+k_{2+j}\le k$, and the factor $f$ is smooth and bounded.
\end{lemma}
\begin{proof}
The proof is essentially the same as the proof of the previous lemma, with the one difference being the possible presence of the factor $r^{-1}\fb\phi$. If this is treated like the other factors, the result after applying the commutators is a factor of the form 
$$r^{-1}\fd^i\fb\phi^{s_2-i,k_2}.$$
If $i$ is odd, this term can be rewritten as 
$$r^{-1}\fd^i\fb\phi^{s_2-i,k_2}=r^{-1}\fd^{i+2}\phi^{s_2-i,k_2}=r^{-1}\pd_\theta \fd^{i+1}\phi^{s_2-i,k_2}=\beta^{s_1+1,k_2}.$$
This is at the level $s_2+1$, which is less than or equal to $s$, because $i$ is odd.

If instead $i$ is even, this term can be rewritten as $r^{-1}\fc^l\fb\phi^{s_2-2l,k_2}$, where $i=2l$. By Lemma \ref{fb_fcl_commutation_lem}, this is (up to lower order terms in $s$) equal to $r^{-1}\fb\fc^l\phi^{s_2-2l,k_2}$.

With this exceptional case having been addressed, the proof is complete. 
\end{proof}

\subsection{The general strategy}\label{nl_strategy_revisited_sec}

Now that the structures of $q^{-2}\tg^k\Gamma^s(q^2\mathcal{N}_\phi)$ and $q^{-2}\tg^k\tilde\Gamma^s(q^2\mathcal{N}_\psi)$ have been determined, the strategy to be employed in the proof of the main theorem can be summarized in the following proposition.

\begin{proposition}\label{strategy_revisited_prop}
Suppose the following estimates hold.
$$(1+\phi)^{-1}\lesssim 1,$$
$$\int_R^\infty\int_0^\pi (r^{(p-1)/2}\gamma^s)^2 r^2\sin^5\theta d\theta dr \lesssim E_{p-1}^s(t),$$
$$||r\beta^s||_{L^\infty(\Sigma_t\{r>R\})}^2\lesssim E^{s+7}(t),$$
$$\int_R^\infty\int_0^\pi (\beta^s)^2r^2\sin^5\theta d\theta dr\lesssim E^s(t),$$
$$||r^{(p-1)/2+1}\gamma^s||_{L^\infty(\Sigma_t\cap\{r>R\})}^2\lesssim E_{p-1}^{s+7}(t).$$
Then
\begin{multline*}
\int_{\Sigma_t\cap\{r>R\}}r^{p+1}|q^{-2}\tg^k\Gamma^s(q^2\mathcal{N}_\phi)|^2 + \int_{\tilde{\Sigma}_t\cap\{r>R\}}r^{p+1}|q^{-2}\tg^k\tilde{\Gamma}^s(q^2\mathcal{N}_\psi)|^2 \\
\lesssim \left(E_{p-1}^s(t)E^{s/2+8}(t)+E^s(t)E_{p-1}^{s/2+8}(t)\right)\sum_{j\le s}(E^{s/2+8}(t))^j.
\end{multline*}
\end{proposition}
\begin{remark}
The second and fourth assumptions are automatically true. The first, third, and fifth assumptions will be true in the context of the bootstrap assumptions in the proof of the main theorem (Theorem \ref{main_thm}).
\end{remark}
\begin{proof}
Note that Lemmas \ref{Nphi_s_k_structure_lem} and \ref{Npsi_s_k_structure_lem} both had the requirements $s_1+...+s_{2+j}\le s+2$ and $\max_i s_i\le s$. This means there can be at most one factor (denoted $\gamma^{hi}$ or $\beta^{hi}$) with up to $s$ derivatives and all of the remaining factors (denoted $\gamma^{lo}$ or $\beta^{lo}$) have at most $(s+2)/2=s/2+1$ derivatives. Since the integrals are taken over the range $r>R$, where $\tg=0$, we are free to ignore all $k$ indices. 

The strategy in this proof is rather simple. All of the factors that have at most $s/2+1$ derivatives are estimated in $L^\infty$, while the remaining factor with at most $s$ derivatives is estimated in $L^2$. The procedure then takes one of two possible directions, depending on whether the high derivative term is a $\gamma$ term or a $\beta$ term.

\textbf{First, we estimate the integral over $\Sigma_t$, using the form given in Lemma \ref{Nphi_s_k_structure_lem}.}
\begin{multline*}
\int_{\Sigma_t\cap\{r>R\}}r^{p+1}|q^{-2}\tg^k\Gamma^s(q^2\mathcal{N}_\phi)|^2 \\
\lesssim \sum_{0\le j\le s} \int_{\Sigma_t\cap\{r>R\}}r^{p+1}\left(f\gamma^{hi}\beta^{lo}(r\beta^{lo})^j\right)^2 + \sum_{0\le j\le s} \int_{\Sigma_t\cap\{r>R\}}r^{p+1}\left(f\beta^{hi}\gamma^{lo}(r\beta^{lo})^j\right)^2 .
\end{multline*}
Now,
\begin{align*}
\int_{\Sigma_t\cap\{r>R\}}r^{p+1}\left(f\gamma^{hi}\beta^{lo}(r\beta^{lo})^j\right)^2
&\lesssim \int_{\Sigma_t\cap\{r>R\}}r^{p-1}\left(f\gamma^{hi}(r\beta^{lo})(r\beta^{lo})^j\right)^2 \\
&\lesssim ||r\beta^{lo}||_{L^\infty}^{2(j+1)}\int_{\Sigma_t\cap\{r>R\}}r^{p-1}(f\gamma^{hi})^2 \\
&\lesssim (E^{s/2+8}(t))^{j+1}\int_{\Sigma_t\cap\{r>R\}}r^{p-1}(f\gamma^{hi})^2 \\
&\lesssim E_{p-1}^s(t)(E^{s/2+8}(t))^{j+1}.
\end{align*}
The last step requires further justification, because the $\gamma$ term could possibly be a $\psi$ term. But in that case, Lemma \ref{Nphi_s_k_structure_lem} also states that $f$ has an additional factor of $\sin^2\theta$ (or in the case where $\gamma$ does not have exactly $s$ derivatives, this factor might be either $\sin\theta$ or $1$, but we can apply Lemma \ref{gain_sin_lem}--see the end of the proof of the example Lemma \ref{ws:example_Nphi_term_lem}).

Also,
\begin{align*}
\int_{\Sigma_t\cap\{r>R\}}r^{p+1}\left(f\beta^{hi}\gamma^{lo}(r\beta^{lo})^j\right)^2 
&\lesssim \int_{\Sigma_t\cap\{r>R\}}\left(f\beta^{hi}(r^{(p+1)/2}\gamma^{lo})(r\beta^{lo})^j\right)^2 \\
&\lesssim ||r^{(p+1)/2}\gamma^{lo}||_{L^\infty}^2||r\beta||_{L^\infty}^{2j}\int_{\Sigma_t\cap\{r>R\}}(f\beta^{hi})^2 \\
&\lesssim E_{p-1}^{s/2+8}(t)(E^{s/2+8}(t))^j\int_{\Sigma_t\cap\{r>R\}}(f\beta^{hi})^2 \\
&\lesssim E^s(t)E_{p-1}^{s/2+8}(t)(E^{s/2+8}(t))^j.
\end{align*}
The last step requires the same justification given in the preceeding calculation.

Combining both estimates, we conclude that 
$$\int_{\Sigma_t\cap\{r>R\}}r^{p+1}|q^{-2}\tg^k\Gamma^s(q^2\mathcal{N}_\phi)|^2 \lesssim \left(E_{p-1}^s(t)E^{s/2+8}(t)+E^s(t)E_{p-1}^{s/2+8}(t)\right)\sum_{j\le s}(E^{s/2+8}(t))^j.$$

\textbf{Next, we estimate the integral over $\tilde{\Sigma}_t$, using the forms given in Lemma \ref{Npsi_s_k_structure_lem}.} The first of the two forms given in Lemma \ref{Npsi_s_k_structure_lem} is very similar to the form given in Lemma \ref{Nphi_s_k_structure_lem}, and so we will not repeat the estimates in detail. The only difference is the presence of an additional factor of $r^{-2}$ that comes along with $f$. This is required in the two final steps.
$$\int_{\tilde{\Sigma}_t\cap\{r>R\}}r^{p-1}(r^{-2}f\gamma^{hi})^2\lesssim E_{p-1}^s(t)$$
or
$$\int_{\tilde{\Sigma}_t\cap\{r>R\}}(r^{-2}f\beta^{hi})^2\lesssim E^s(t).$$
The factor of $r^{-2}$ (inside the square) is needed to offset the additional factor of $r^4$ in the volume form for $\tilde{\Sigma}_t$.

The second of the two forms has the factor $r^{-1}\fb\fc^l\phi^{s_2-2l}$, which takes the place of the $\beta^{s_2}$ factor from the first form. We now show that it can be estimated the same way the $\beta^{s_2}$ factor was estimated.

For the case $s_2\le s/2+1$, we must estimate $r$ times the factor $r^{-1}\fb\fc^l\phi^{s_2-2l}$ in $L^\infty$. That is, using the Lemma \ref{low_order_infty_lem},
$$||\fb\fc^l\phi^{s/2+1-2l}||_{L^\infty}^2=||\fd^{2l+2}\phi^{s/2+1-2l}||_{L^\infty}^2 \lesssim E^{s/2+8}(t).$$
For the case where $s_2$ could be at most $s$, we use the fact that there are additional factors of $\sin\theta$ in the volume form for $\tilde{\Sigma}_t$ and again the fact that an additional factor of $r^{-2}$ accompanies $f$. That is,
\begin{align*}
\int_{\tilde{\Sigma}_t\cap\{r>R\}}(r^{-2}fr^{-1}\fb\fc^l\phi^{s-2l})^2
&\lesssim \int_R^\infty\int_0^\pi(r^{-3}\sin^{-1}\theta\pd_\theta\fc^l\phi^{s-2l})^2r^6\sin^5\theta d\theta dr\\
&\lesssim \int_R^\infty\int_0^\pi(r^{-1}\pd_\theta \fc^l\phi^{s-2l})^2r^2\sin^3\theta d\theta dr\\
&\lesssim \int_R^\infty\int_0^\pi(r^{-1}\pd_\theta \phi^s)^2r^2\sin\theta d\theta dr \\
&\lesssim E^s(t).
\end{align*}
(In the second-to-last step, we used an estimate from Theorem \ref{cl_embedding_thm}.) Both of these estimates show that the $r^{-1}\fb\fc^l\phi^{s_2-2l}$ factor can be treated the same way as the $\beta^{s_2}$ factor.

We conclude that
$$\int_{\tilde{\Sigma}_t\cap\{r>R\}}r^{p+1}|q^{-2}\tg^k\tilde{\Gamma}^s(q^2\mathcal{N}_\psi)|^2 \lesssim \left(E_{p-1}^s(t)E^{s/2+8}(t)+E^s(t)E_{p-1}^{s/2+8}(t)\right)\sum_{j\le s}(E^{s/2+8}(t))^j.$$
This completes the proof. 
\end{proof}

\section{Theorem: Global boundedness and decay for axisymmetric perturbations of the nontrivial solutions to the wave map problem from slowly rotating Kerr spacetimes to the hyperbolic plane preserving angular momentum}\label{wk:main_thm_sec}

In this final section, we state and prove the main theorem. We begin in \S\ref{main_thm_statement_sec} with the precise statement of the main theorem, which completely describes the future asymptotic behavior of various energy norms as well as weighted derivatives of $\phi$ and $\psi$. Then the proof follows.

The proof is a bootstrap argument with the main bootstrap assumptions stated in \S\ref{bootstrap_assumptions_sec}. Then in \S\ref{improved_pointwise_estimates_sec}, the bootstrap assumptions are used to improve the $L^\infty$ estimates in Proposition \ref{infinity_prop} by removing the error terms from the right side. The resulting improved $L^\infty$ estimates play a crucial role in the remainder of the proof. Next, in \S\ref{recover_phi_ba_sec} the $L^\infty$ estimates are used to recover the first bootstrap assumption, which gives an upper bound on the factor $(1+\phi)^{-1}$ that appears in the nonlinear terms. The remainder of the proof is split into two parts. 

The first part, which starts with \S\ref{homogeneous_case_sec}, recovers the bootstrap assumptions only for the homogeneous norms $\mathring{E}^s(t)$, $\mathring{E}_p^s(t)$, etc. These norms satisfy slightly simpler estimates, because they are based on the commutator $\pd_t$, which commutes with the entire linear system. As a reminder, these estimates, which are given in Theorem \ref{p_o_thm}, are
$$\mathring{E}^s(t_2)\lesssim \mathring{E}^s(t_1)+\int_{t_1}^{t_2}\mathring{N}^s(t)dt,$$
$$\mathring{E}_p^s(t_2)+\int_{t_1}^{t_2}\mathring{B}_p^s(t)dt\lesssim \mathring{E}_p^s(t_1)+\int_{t_1}^{t_2}\mathring{N}_p^s(t)dt.$$

The second part, which starts with \S\ref{general_case_s_k_sec}, recovers the bootstrap assumptions for the norms $E^{s,k}(t)$, $E_p^{s,k}(t)$, etc. These norms satisfy slightly more complicated estimates, because they use other commutators that do not completely commute with the linear system. As a reminder, these estimates, which are given in Theorem \ref{p_s_k_thm}, are
$$E^{s,k}(t_2)\lesssim E^{s,k}(t_1)+\int_{t_1}^{t_2}B_{p'}^{s+2,k-1}(t)+\mathring{B}_1^{s+1}(t)+B_1^{s,k}(t)+N^{s,k}(t)dt,$$
$$E_p^{s,k}(t_2)+\int_{t_1}^{t_2}B_p^{s,k}(t)dt\lesssim E_p^{s,k}(t_1)+\int_{t_1}^{t_2}B_{p'}^{s+2,k-1}(t)+\mathring{B}_p^{s+1}(t)+N_p^{s,k}(t)dt.$$
In particular, the fact that these estimates depend on $\mathring{B}_p^{s+1}(t)$ means that this second part depends on the results proved in the first part.

\subsection{The main theorem}\label{main_thm_statement_sec}
We begin with the statement of the main theorem.
\begin{theorem}\label{main_thm}
Let
\begin{align*}
X&=A+A\phi, \\
Y&=B+A^2\psi.
\end{align*}
(In particular, the assumption $Y-B=O(\sin^4\theta)$ near the axis excludes any perturbations corresponding to a change in angular momentum in the sense described in \S\ref{ernst_potential_axis_derivation}.) 

Suppose the pair $(X,Y)$ is axisymmetric and satisfies the wave map system
\begin{align*}
X\Box_g X&=\pd^\alpha X\pd_\alpha X-\pd^\alpha Y\pd_\alpha Y \\
X\Box_g Y&=2\pd^\alpha X\pd_\alpha Y,
\end{align*}
where $g$ is a Kerr metric with sufficiently small angular momentum $|a|/M$.

Define the energies
$$E^{\ul{n}}(t)=\mathring{E}^n(t)+\sum_{s+2k=n-1}E^{s,k}(t),$$
$$E_p^{\ul{n}}(t)=\mathring{E}_p^{n}(t)+\sum_{s+2k=n-1}E_p^{s,k}(t).$$
Then for $\delp,\delm>0$ sufficiently small, if the initial data for $(\phi,\psi)$ decay sufficiently fast as $r\rightarrow \infty$ and have size
\begin{equation*}
I_0=E^{\ul{29}}(0)+E^{\ul{29}}_{2-\delp}(0)
\end{equation*}
sufficiently small, then the following estimates hold for $t\ge 0$ (with $T=1+t$).

I) The energies satisfy
$$E^{\ul{29}}(t)\lesssim I_0$$
$$E^{\ul{29}}_{p\in[\delm,2-\delp]}(t)\lesssim I_0$$
$$E^{\ul{27}}_{p\in[1-\delp,2-\delp]}(t)\lesssim T^{p-2+\delp}I_0$$
$$E^{\ul{25}}_{p\in[\delm,2-\delp]}(t)\lesssim T^{p-2+\delp}I_0$$
$$\int_t^{\infty}E^{\ul{23}}_{p\in[\delm-1,\delm]}(\tau)d\tau\lesssim T^{p-2+\delp+1}I_0$$

II) For all $s,k$ such that $s+2k\le 28$, the following $L^\infty$ estimates hold.
\begin{multline*}
|r^{p+1}\bar{D}\fd^l\phi^{s-l,k}|^2+|r^pD\fd^l\phi^{s-l,k}|^2+|r^p\fd^l\phi^{s-l,k}|^2 \\
+|r^{p+3}\bar{D}\fd^l\psi^{s-l,k}|^2+|r^{p+2}D\fd^l\psi^{s-l,k}|^2+|r^{p+2}\fd^l\psi^{s-l,k}|^2 \\
\lesssim E_{2p}^{s+5,k+1}(t)+E_{2p}^{s+7,k}(t)
\end{multline*}

III) Together, (I) and (II) imply that if $s+2k\le 15$, for all $p\in [\delm/2,(2-\delp)/2]$,
\begin{multline*}
|r^{p+1}\bar{D}\fd^l\phi^{s-l,k}|+|r^pD\fd^l\phi^{s-l,k}|+|r^p\fd^l\phi^{s-l,k}| \\
+|r^{p+3}\bar{D}\fd^l\psi^{s-l,k}|+|r^{p+2}D\fd^l\psi^{s-l,k}|+|r^{p+2}\fd^l\psi^{s-l,k}| \\
\lesssim T^{(2p-2+\delp)/2}I_0^{1/2}
\end{multline*}
and additionally for $p\in [(\delm-1)/2,\delm/2]$,
\begin{multline*}
\int_t^\infty |r^{p+1}\bar{D}\fd^l\phi^{s-l,k}|+|r^pD\fd^l\phi^{s-l,k}|+|r^p\fd^l\phi^{s-l,k}| \\
+\int_t^\infty |r^{p+3}\bar{D}\fd^l\psi^{s-l,k}|+|r^{p+2}D\fd^l\psi^{s-l,k}|+|r^{p+2}\fd^l\psi^{s-l,k}| \\
\lesssim T^{(2p-2+\delp)/2+1}I_0^{1/2}.
\end{multline*}
The final estimate should be interpreted as saying that $|r^{(\delm-1)/2+1}\bar{D}\fd^l\phi^{s-l,k}|$, \\
$|r^{(\delm-1)/2}D\fd^l\phi^{s-l,k}|$, $|r^{(\delm-1)/2}\fd^l\phi^{s-l,k}|$, $|r^{(\delm-1)/2+3}\bar{D}\fd^l\psi^{s-l,k}|$, $|r^{(\delm-1)/2+2}D\fd^l\psi^{s-l,k}|$, and $|r^{(\delm-1)/2+2}\fd^l\psi^{s-l,k}|$ decay like $T^{(\delm-3+\delp)/2}$ in a weak sense.

\end{theorem}

\subsection{Bootstrap assumptions}\label{bootstrap_assumptions_sec}

We begin the proof of Theorem \ref{main_thm} by making the following bootstrap assumptions.
\begin{align}
|\phi|&\le 1/2, \label{phi_ba} \\
E^{\ul{29}}(t)&\le C_bI_0, \label{energy_ba} \\
E^{\ul{29}}_{p\in[\delm,2-\delp]}(t)&\le C_b I_0, \label{_29_ba} \\
E^{\ul{25}}_{p\in[\delm,2-\delp]}(t)&\le C_b T^{p-2+\delp}I_0, \label{_25_ba} \\
\int_t^\infty E^{\ul{23}}_{p\in[\delm-1,\delm]}(\tau)d\tau &\le C_b T^{p-2+\delp+1}I_0. \label{_23_ba}
\end{align}
Note that, with the exception of the highest order energies, these bootstrap assumptions are consistent with the general principle that $\mathring{E}_p^s(t)$ and $E_p^{s,k}(t)$ behave like $T^{p-2+\delp}$, which the reader should keep in mind throughout the proof of the main theorem.

\subsection{Improved $L^\infty$ estimates}\label{improved_pointwise_estimates_sec}

The $L^\infty$ estimates from Proposition \ref{infinity_prop} are essential to the argument of the proof of the main theorem. But for the sake of clarity, we first remove the error terms from these estimates and summarize them in the following lemma.
\begin{lemma}\label{simplified_pointwise_lemma}
In the context of the bootstrap assumptions provided in \S\ref{bootstrap_assumptions_sec}, the following $L^\infty$ estimates hold for $s+2k\le 28$ and all $p$ in any bounded range.

For $r\ge r_H$,
\begin{multline*}
|r^{p+1}\bar{D}\fd^l\phi^{s-l,k}|^2+|r^pD\fd^l\phi^{s-l,k}|^2+|r^p\fd^l\phi^{s-l,k}|^2 \\
+|r^{p+3}\bar{D}\fd^l\psi^{s-l,k}|^2+|r^{p+2}D\fd^l\psi^{s-l,k}|^2+|r^{p+2}\fd^l\psi^{s-l,k}|^2 \\
\lesssim E_{2p}^{s+5,k+1}(t)+E_{2p}^{s+7,k}(t)
\end{multline*}
and for $r\ge r_0>r_H$,
\begin{equation*}
|rD\fd^l\phi^{s-l,k}|^2+|r^3D\fd^l\psi^{s-l,k}|^2\lesssim E^{s+5,k+1}(t)+E^{s+7,k}(t)
\end{equation*}
These are the same as the estimates from Proposition \ref{infinity_prop}, except that the error terms have been removed.
\end{lemma}
\begin{proof} 
Recall that, according to Proposition \ref{infinity_prop}, for $r\ge r_H$,
\begin{multline*}
|r^{p+1}\bar{D}\fd^l\phi^{s-l,k}|^2+|r^pD\fd^l\phi^{s-l,k}|^2+|r^p\fd^l\phi^{s-l,k}|^2 \\
+|r^{p+3}\bar{D}\fd^l\psi^{s-l,k}|^2+|r^{p+2}D\fd^l\psi^{s-l,k}|^2+|r^{p+2}\fd^l\psi^{s-l,k}|^2 \\
\lesssim E_{2p}^{s+5,k+1}(t)+E_{2p}^{s+7,k}(t)+\int_{\Sigma_t}r^{2p}(\Box_g\phi^{s+5,k})^2+\int_{\tilde\Sigma_t}r^{2p}(\Box_{\tilde{g}}\psi^{s+5,k})^2
\end{multline*}
and for $r\ge r_0>r_H$,
\begin{multline*}
|rD\fd^l\phi^{s-l,k}|^2+|r^3D\fd^l\psi^{s-l,k}|^2 \\
\lesssim E^{s+5,k+1}(t)+E^{s+7,k}(t)+\int_{\Sigma_t\cap\{r>r_0\}}(\Box_g\phi^{s+5,k})^2+\int_{\tilde\Sigma_t\cap\{r>r_0\}}(\Box_{\tilde{g}}\psi^{s+5,k})^2.
\end{multline*}
To prove the lemma, it suffices to show that
$$\int_{\Sigma_t}r^{2p}(\Box_g\phi^{s+5,k})^2+\int_{\tilde\Sigma_t}r^{2p}(\Box_{\tilde{g}}\psi^{s+5,k})^2 \lesssim E_{2p}^{s+5,k+1}(t)+E_{2p}^{s+7,k}(t)$$
and
$$\int_{\Sigma_t\cap\{r>r_0\}}(\Box_g\phi^{s+5,k})^2+\int_{\tilde\Sigma_t\cap\{r>r_0\}}(\Box_{\tilde{g}}\psi^{s+5,k})^2 \lesssim E^{s+5,k+1}(t)+E^{s+7,k}(t).$$
To do this, we use the equations for $\Box_g\phi^{s+5,k}$ and $\Box_{\tilde{g}}\psi^{s+5,k}$. This will result in a number of linear and nonlinear terms. While the linear terms can be estimated directly by the energy norms, the nonlinear terms will again require the use of $L^\infty$ estimates. For this reason, the proof of the lemma requires a nested bootstrap argument.

The bootstrap assumptions for this argument are that for $r\ge r_H$,
\begin{multline}\label{infinity_inner_Ep_ba}
|r^{p+1}\bar{D}\fd^l\phi^{s-l,k}|^2+|r^pD\fd^l\phi^{s-l,k}|^2+|r^p\fd^l\phi^{s-l,k}|^2 \\
+|r^{p+3}\bar{D}\fd^l\psi^{s-l,k}|^2+|r^{p+2}D\fd^l\psi^{s-l,k}|^2+|r^{p+2}\fd^l\psi^{s-l,k}|^2 \\
\le C_b\left(E_{2p}^{s+5,k+1}(t)+E_{2p}^{s+7,k}(t)\right)
\end{multline}
and for $r\ge r_0>r_H$,
\begin{equation}\label{infinity_inner_E_ba}
|rD\fd^l\phi^{s-l,k}|^2+|r^3D\fd^l\psi^{s-l,k}|^2\le C_b \left(E^{s+5,k+1}(t)+E^{s+7,k}(t) \right).
\end{equation}

We use one representative example to illustrate the bootstrap argument. Consider the term 
$$\frac{L\phi^s\pd_r\phi+L\phi\pd_r\phi^s}{1+\phi},$$
which shows up in the equation for $\Box_g\phi^s$ (see Proposition \ref{N_identities_prop}). We estimate this term as follows.
\begin{multline*}
\int_{\Sigma_t}r^{2p}\left(\frac{L\phi^s\pd_r\phi+L\phi\pd_r\phi^s}{1+\phi}\right)^2 \\
\lesssim \int_{\Sigma_t}r^{2p}(L\phi^s)^2||\pd_r\phi||^2_{L^\infty(\Sigma_t)} +\int_{\Sigma_t}r^{2p-2}(\pd_r\phi^s)^2||rL\phi||^2_{L^\infty(\Sigma_t)}.
\end{multline*}
Then we apply the bootstrap assumption (\ref{infinity_inner_Ep_ba}) with $p=0$ to estimate the $L^\infty$ norms.
\begin{multline*}
\int_{\Sigma_t}r^{2p}\left(\frac{L\phi^s\pd_r\phi+L\phi\pd_r\phi^s}{1+\phi}\right)^2 \\
\lesssim \int_{\Sigma_t}r^{2p}(L\phi^s)^2C_b\left(E_0^{5,1}(t)+E_0^{7,0}(t)\right) +\int_{\Sigma_t}r^{2p-2}(\pd_r\phi^s)^2C_b\left(E_0^{5,1}(t)+E_0^{7,0}(t)\right).
\end{multline*}
(In particular, since we only needed to use the $p=0$ norms, this estimate is not particularly delicate.) Next, we use the bootstrap assumption (\ref{_29_ba}) for the main theorem to estimate the energy norms.
$$\int_{\Sigma_t}r^{2p}\left(\frac{L\phi^s\pd_r\phi+L\phi\pd_r\phi^s}{1+\phi}\right)^2 
\lesssim \int_{\Sigma_t}r^{2p}(L\phi^s)^2C_b(C_bI_0) +\int_{\Sigma_t}r^{2p-2}(\pd_r\phi^s)^2C_b(C_bI_0).$$
And finally, we estimate the remaining integrated quantities by the energy norm $E_{2p}^{s,0}(t)$.
$$\int_{\Sigma_t}r^{2p}\left(\frac{L\phi^s\pd_r\phi+L\phi\pd_r\phi^s}{1+\phi}\right)^2 \lesssim C_b^2I_0E_{2p}^{s,0}(t).$$
As a side note, this procedure will be repeatedly used in the remainder of the proof of the main theorem.

By repeating the procedure in the above example for all the different terms that arise in the equations for $\Box_g\phi^{s+5,k}$ and $\Box_{\tilde{g}}\psi^{s+5,k}$, we eventually conclude that
$$\int_{\Sigma_t}r^{2p}(\Box_g\phi^{s+5,k})^2+\int_{\tilde\Sigma_t}r^{2p}(\Box_{\tilde{g}}\psi^{s+5,k})^2 \lesssim \sum_{j=0}^{s+5+k}(C_b^2I_0)^j\left( E_{2p}^{s+5,k+1}(t)+E_{2p}^{s+7,k}(t)\right)$$
and
$$\int_{\Sigma_t\cap\{r>r_0\}}(\Box_g\phi^{s+5,k})^2+\int_{\tilde\Sigma_t\cap\{r>r_0\}}(\Box_{\tilde{g}}\psi^{s+5,k})^2 \lesssim \sum_{j=0}^{s+5+k}(C_b^2I_0)^j\left(E^{s+5,k+1}(t)+E^{s+7,k}(t)\right).$$
(The exponent $j$ corresponds to the number of times a differential operator acts on the denominator $1+\phi$--see Lemmas \ref{Nphi_s_k_structure_lem} and \ref{Npsi_s_k_structure_lem}. It can also be zero because of the presence of linear terms.) By taking $I_0$ sufficiently small so that $C_b^2I_0\lesssim 1$, the dependence on the bootstrap constant $C_b$ can be removed. The resulting estimates can be used together with Proposition \ref{infinity_prop} to conclude that for $r\ge r_H$,
\begin{multline*}
|r^{p+1}\bar{D}\fd^l\phi^{s-l,k}|^2+|r^pD\fd^l\phi^{s-l,k}|^2+|r^p\fd^l\phi^{s-l,k}|^2 \\
+|r^{p+3}\bar{D}\fd^l\psi^{s-l,k}|^2+|r^{p+2}D\fd^l\psi^{s-l,k}|^2+|r^{p+2}\fd^l\psi^{s-l,k}|^2 \\
\lesssim E_{2p}^{s+5,k+1}(t)+E_{2p}^{s+7,k}(t)
\end{multline*}
and for $r\ge r_0>r_H$,
\begin{equation*}
|rD\fd^l\phi^{s-l,k}|^2+|r^3D\fd^l\psi^{s-l,k}|^2\lesssim E^{s+5,k+1}(t)+E^{s+7,k}(t).
\end{equation*}
This recovers the bootstrap assumptions (\ref{infinity_inner_Ep_ba}-\ref{infinity_inner_E_ba}).
\end{proof}
The conclusion of the above lemma is the same as the statement of part (II) of the main theorem.

\begin{remark}
Lemma \ref{simplified_pointwise_lemma}, which is a simplified version of Proposition \ref{infinity_prop} in the sense that there are no error terms on the right side of any of the estimates in Lemma \ref{simplified_pointwise_lemma}, will be used in the remainder of the proof of the main theorem as a replacement for Proposition \ref{infinity_prop}.
\end{remark}

\subsection{Recovering the bootstrap assumption for $|\phi|$}\label{recover_phi_ba_sec}

With the improved $L^\infty$ estimates, we can now recover the bootstrap assumption (\ref{phi_ba}) for $|\phi|$. The $L^\infty$ estimates imply that
$$|\phi|^2\lesssim E_0^{5,1}(t)+E_0^{7,0}(t).$$
The bootstrap assumption (\ref{_29_ba}) implies that
$$E_0^{5,1}(t)+E_0^{7,0}(t)\lesssim C_b I_0.$$
Thus,
$$|\phi|^2\lesssim C_b I_0.$$
As long as $C_bI_0$ is sufficiently small, this guarantees that $|\phi|$ is also sufficiently small. This fact allows us to recover the bootstrap assumtion (\ref{phi_ba}).

\subsection{The highest order (homogeneous) case}\label{homogeneous_case_sec}

As explained in the introduction of this section, this is the point where we turn to the first of two parts of the remainder of the proof. In particular, we will recover the bootstrap assumptions (\ref{energy_ba}-\ref{_23_ba}) for the homogeneous norms $\mathring{E}^s(t)$, $\mathring{E}_p^s(t)$, etc.

To begin, in \S\ref{part_1_refined_estimates_sec} we prove refined estimates for the nonlinear norms $\mathring{N}^s(t)$ and $\mathring{N}_p^s(t)$. These estimates constitute the crucial step of the proof, because they handle the nonlinear terms near $i^0$ using $L^\infty$ estimates sharply. In \S\ref{part_1_E_boundedness_sec} and \S\ref{part_1_Ep_boundedness_sec}, the bootstrap assumptions (\ref{energy_ba}) and (\ref{_29_ba}) at the highest level $s=29$ are recovered. Then in \S\ref{part_1_decay_sec}, a decay lemma is proved and used to recover the bootstrap assumption (\ref{_25_ba}) at the level $s=25$. Finally, in \S\ref{part_1_weak_decay_sec}, the bootstrap assumption (\ref{_23_ba}) at the level $s=23$, which assumes a weaker form of decay, is recovered.

\subsection{Refined estimates for $\mathring{N}^s(t)$ and $\mathring{N}_p^s(t)$}\label{part_1_refined_estimates_sec}

The $L^\infty$ estimates given in Lemma \ref{simplified_pointwise_lemma} allow us to provide refined estimates for the nonlinear error terms. \textbf{This is the crucial step of the proof.}

\begin{lemma}\label{mr_refined_nl_lem}
In the context of the bootstrap assumptions provided in \S\ref{bootstrap_assumptions_sec}, if $s\le 29$ and $C_bI_0\le 1$, then
$$\mathring{N}^s(t) \lesssim (\mathring{E}^s(t))^{1/2}\left((\mathring{E}^s(t))^{1/2}(E_{\delm-1}^{\ul{22}}(t))^{1/2}+(\mathring{E}^s_{1-\delm}(t))^{1/2}(E^{\ul{22}}_{\delm-1}(t))^{1/2}\right),$$
$$\mathring{N}_p^s(t) \lesssim \mathring{E}^s(t)B_p^{\ul{22}}(t)+\mathring{B}_p^s(t)E^{\ul{22}}(t)+\mathring{E}^s_{p'}(t)(E_{p''}^{\ul{22}}(t))^{1/2}.$$
\end{lemma}
\begin{proof}
First, we recall the definitions of $\mathring{N}^s(t)$ and $\mathring{N}_p^s(t)$ from Theorem \ref{p_o_thm}.
$$\mathring{N}^s(t)=(\mathring{E}^s(t))^{1/2}\sum_{s'\le s}\left(||\pd_t^{s'}\mathcal{N}_\phi||_{L^2(\Sigma_t)}+||\pd_t^{s'}\mathcal{N}_\psi||_{L^2(\tilde{\Sigma}_t)}\right),$$
\begin{multline*}
\mathring{N}_p^s(t)=\sum_{s'\le s}\int_{\Sigma_t}r^{p+1}(\pd_t^{s'}\mathcal{N}_\phi)^2 + \sum_{s'\le s}\int_{\tilde\Sigma_t}r^{p+1}(\pd_t^{s'}\mathcal{N}_\psi)^2 \\
+ \sum_{s'\le s}\int_{\Sigma_t\cap\{r\approx r_{trap}\}}|\pd_t^{s'+1}\phi\pd_t^{s'}\mathcal{N}_\phi| + \sum_{s'\le s}\int_{\tilde\Sigma_t\cap\{r\approx r_{trap}\}}|\pd_t^{s'+1}\psi\pd_t^{s'}\mathcal{N}_\psi|. 
\end{multline*}
Therefore, it suffices to prove the following three estimates.
\begin{equation*}
||\pd_t^s\mathcal{N}_\phi||_{L^2(\Sigma_t)}+||\pd_t^s\mathcal{N}_\psi||_{L^2(\tilde{\Sigma}_t)}
\lesssim (\mathring{E}^s(t))^{1/2}(E_{\delm-1}^{\ul{22}}(t))^{1/2}+(\mathring{E}^s_{1-\delm}(t))^{1/2}(E^{\ul{22}}_{\delm-1}(t))^{1/2}
\end{equation*}
\begin{equation*}
\int_{\Sigma_t\cap\{r>R\}}r^{p+1}(\pd_t^s\mathcal{N}_\phi)^2 + \int_{\tilde\Sigma_t\cap\{r>R\}}r^{p+1}(\pd_t^s\mathcal{N}_\psi)^2 \lesssim \mathring{E}^s(t)B_p^{\ul{22}}(t)+\mathring{B}_p^s(t)E^{\ul{22}}(t)
\end{equation*}
\begin{multline*}
\int_{\Sigma_t\cap\{r<R\}}(\pd_t^s\mathcal{N}_\phi)^2 + \int_{\tilde\Sigma_t\cap\{r<R\}}(\pd_t^s\mathcal{N}_\psi)^2 \\
+ \int_{\Sigma_t\cap\{r\approx r_{trap}\}}|\pd_t^{s+1}\phi\pd_t^s\mathcal{N}_\phi| + \int_{\tilde\Sigma_t\cap\{r\approx r_{trap}\}}|\pd_t^{s+1}\psi\pd_t^s\mathcal{N}_\psi| \\
\lesssim \mathring{E}^s_{p'}(t)(E_{p''}^{\ul{22}}(t))^{1/2} 
\end{multline*}
We begin by observing that the second estimate is a consequence of Proposition \ref{strategy_revisited_prop} (adapted to the homogeneous case) with two unimportant changes. The first is that the norms $\mathring{E}_{p-1}^s(t)$ and $E_{p-1}^{s/2+8}(t)$ used in Proposition \ref{strategy_revisited_prop} have been replaced with the norms $\mathring{B}_p^s(t)$ and $B_p^{\ul{22}}(t)$, since they are equivalent in a region excluding the trapping radius. The second is that the factor of $\sum_{j\le 29}(\mathring{E}^{29}(t))^j$ that appears on the right side of the estimate from Proposition \ref{strategy_revisited_prop} has been eliminated. Given the bootstrap assumption (\ref{energy_ba}) that $\mathring{E}^{29}(t)\lesssim C_bI_0$ and the fact that $C_bI_0\le 1$, this factor is unimportant.

The third estimate is much simpler. Since all of the integrated quantites are defined on a bounded radius, the $r$ factors can be replaced with constants and therefore we have the freedom to choose $p'$ and $p''$. Using again the fact that $C_bI_0\le 1$, we ignore additional factors of any energy norm in this estimate.

The first estimate is a hybrid of the other two. To start,
\begin{multline*}
||\pd_t^s\mathcal{N}_\phi||_{L^2(\Sigma_t)}+||\pd_t^s\mathcal{N}_\psi||_{L^2(\tilde{\Sigma}_t)} \\
\lesssim ||\pd_t^s\mathcal{N}_\phi||_{L^2(\Sigma_t\cap\{r<R\})}+||\pd_t^s\mathcal{N}_\psi||_{L^2(\tilde{\Sigma}_t\cap\{r<R\})} \\
+||\pd_t^s\mathcal{N}_\phi||_{L^2(\Sigma_t\cap\{r>R\})}+||\pd_t^s\mathcal{N}_\psi||_{L^2(\tilde{\Sigma}_t\cap\{r>R\})}
\end{multline*}
The norms on a bounded radius can be estimated as
$$||\pd_t^s\mathcal{N}_\phi||_{L^2(\Sigma_t\cap\{r<R\})}+||\pd_t^s\mathcal{N}_\psi||_{L^2(\tilde{\Sigma}_t\cap\{r<R\})} \lesssim \mathring{E}_{1-\delm}^s(t))^{1/2}(E_{\delm-1}^{\ul{22}}(t))^{1/2}$$
 for the same reason as in the third estimate. The norms over the remaining region $\{r>R\}$ can be estimated as
\begin{multline*}
||\pd_t^s\mathcal{N}_\phi||_{L^2(\Sigma_t\cap\{r>R\})}+||\pd_t^s\mathcal{N}_\psi||_{L^2(\tilde{\Sigma}_t\cap\{r>R\})} \\
\lesssim (\mathring{E}^s(t))^{1/2}(E_{\delm-1}^{\ul{22}}(t))^{1/2}+(\mathring{E}^s_{1-\delm}(t))^{1/2}(E^{\ul{22}}_{\delm-1}(t))^{1/2}
\end{multline*}
using the same strategy that is given in Proposition \ref{strategy_revisited_prop}. One simply has to check that the appropriate powers of $r$ can be assigned to each factor.
\end{proof}

\begin{corollary}\label{mr_NL_absorb_bulk}
In the context of the bootstrap assumptions provided in \S\ref{bootstrap_assumptions_sec}, Theorem \ref{p_o_thm} and Lemma \ref{mr_refined_nl_lem} imply that if $s\le 29$ and $C_bI_0$ is sufficiently small, then
$$\mathring{E}_p^{s}(t_2)+\int_{t_1}^{t_2}\mathring{B}_p^{s}(t)dt \lesssim \mathring{E}_p^{s}(t_1)+\int_{t_1}^{t_2}\mathring{E}^s(t)B_p^{\ul{22}}(t)dt+(C_bI_0)^{1/2}\int_{t_1}^{t_2}\mathring{E}_p^{s}(t)T^{(\delm-3+\delp)/2}dt.$$
\end{corollary}
\begin{proof}
According to Theorem \ref{p_o_thm},
$$\mathring{E}_p^s(t_2)+\int_{t_1}^{t_2}\mathring{B}_p^s(t)dt\lesssim \mathring{E}_p^s(t_1)+\int_{t_1}^{t_2}\mathring{N}_p^s(t)dt.$$
By Lemma \ref{mr_refined_nl_lem},
$$\mathring{N}_p^s(t) \lesssim \mathring{E}^s(t)B_p^{\ul{22}}(t)+\mathring{B}_p^s(t)E^{\ul{22}}(t)+\mathring{E}^s_p(t)(E_{\delm-1}^{\ul{22}}(t))^{1/2}.$$
Using the bootstrap assumptions,
$$\int_{t_1}^{t_2}\mathring{B}_p^s(t)E^{\ul{22}}(t)dt \lesssim \int_{t_1}^{t_2}\mathring{B}_p^{s}(t)C_bI_0dt \lesssim C_bI_0\int_{t_1}^{t_2}\mathring{B}_p^{s}(t)dt.$$
Thus, if $C_bI_0$ is sufficiently small, then this error term can be absorbed into the bulk term on the left side.

Again using the bootstrap assumptions and the weak decay principle,
\begin{multline*}
\int_{t_1}^{t_2}\mathring{E}_{p}^s(t)E_{\delm-1}^{\ul{22}}(t)dt \lesssim \int_{t_1}^{t_2}\mathring{E}_p^s(t)(C_bI_0)^{1/2}T^{(\delm-3+\delp)/2}dt \\
\lesssim (C_bI_0)^{1/2}\int_{t_1}^{t_2}\mathring{E}_p^s(t)T^{(\delm-3+\delp)/2}dt.
\end{multline*}
Thus,
$$\mathring{E}_p^{s}(t_2)+\int_{t_1}^{t_2}\mathring{B}_p^{s}(t)dt \lesssim \mathring{E}_p^{s}(t_1)+\int_{t_1}^{t_2}\mathring{E}^s(t)B_p^{\ul{22}}(t)dt+(C_bI_0)^{1/2}\int_{t_1}^{t_2}\mathring{E}_p^{s}(t)T^{(\delm-3+\delp)/2}dt.$$
\end{proof}

\subsection{Recovering boundedness of $\mathring{E}^{29}(t)$}\label{part_1_E_boundedness_sec}

By Theorem \ref{p_o_thm} and Lemma \ref{mr_refined_nl_lem},
\begin{align*}
\mathring{E}^{29}(t) \lesssim & \mathring{E}^{29}(0)+\int_0^tN^{29}(\tau)d\tau \\
\lesssim & \mathring{E}^{29}(0) \\
&+\int_0^t(\mathring{E}^{29}(\tau))^{1/2}\left((\mathring{E}^{29}(\tau))^{1/2}E_{\delm-1}^{\ul{22}}(\tau))^{1/2}+(\mathring{E}^{29}_{1-\delm}(\tau))^{1/2}(E^{\ul{22}}_{\delm-1}(\tau))^{1/2}\right)d\tau \\
\lesssim & I_0+\int_0^t(C_bI_0)^{1/2}(C_bI_0)^{1/2}(C_bT^{\delm-3+\delp}I_0)^{1/2}d\tau \\
\lesssim & (1+C_b^{3/2}I_0^{1/2})I_0.
\end{align*}
In particular, we used the weak decay principle in the third step. It follows that if $I_0$ is chosen sufficiently small so that $C_b^3I_0\lesssim 1$, then
$$\mathring{E}^{29}(t)\lesssim I_0.$$
This recovers the bootstrap assumption (\ref{energy_ba}) at the highest level of derivatives.

\subsection{Recovering boundedness of $\mathring{E}_{2-\delp}^{29}(t)$}\label{part_1_Ep_boundedness_sec}

By Corollary \ref{mr_NL_absorb_bulk},
\begin{multline*}
\mathring{E}_{2-\delp}^{29}(t)+\int_0^t\mathring{B}_{2-\delp}^{29}(\tau)d\tau \\
\lesssim \mathring{E}_{2-\delp}^{29}(0)+\int_0^t\mathring{E}^{29}(\tau)B_{2-\delp}^{\ul{22}}(\tau)d\tau + (C_bI_0)^{1/2}\int_0^t \mathring{E}_{2-\delp}^{29}(\tau)T^{(\delm-3+\delp)/2}d\tau \\
\lesssim I_0+I_0\int_0^tB_{2-\delp}^{\ul{22}}(\tau)d\tau+(C_bI_0)^{1/2}\int_0^tC_bI_0T^{(\delm-3+\delp)/2}d\tau \\
\lesssim I_0+C_bI_0^2+(C_bI_0)^{3/2}.
\end{multline*}
 It follows that if $C_b^3I_0\lesssim 1$, then
$$\mathring{E}_{2-\delp}^{29}(t)+\int_0^t\mathring{B}_{2-\delp}^{29}(\tau)d\tau \lesssim I_0.$$
This recovers the bootstrap assumption (\ref{_29_ba}) at the highest level of derivatives.

\subsection{Proving decay for $\mathring{E}_{p\in[1-\delp,2-\delp]}^{28}(t)$ and $\mathring{E}_{p\in[\delm,2-\delp]}^{27}(t)$}\label{part_1_decay_sec}

We now use the hierarchy of the dynamic estimates to prove decay in time. For convenience, we prove decay in the following lemma, which will then be repeatedly used for a different values of $p$. Essentially, this lemma states that if a $r^{p+1}$-weighted energy norm is bounded or decays at a certain rate, then the $r^p$-weighted energy norm decays at a rate with an additional factor of $T^{-1}$. This is the source of the tradeoff between a factor of $r$ and a factor of $t$.
\begin{lemma}\label{mr_NL_WE_decay}
Suppose  $p+1,p\in [\delm,2-\delp]$ and
$$\mathring{E}_{p+1}^{s+1}(t)\lesssim T^{(p+1)-2+\delp}I_0,$$
$$\int_t^\infty B_{p+1}^{\ul{22}}(\tau)d\tau \le C_bT^{(p+1)-2+\delp}I_0,$$
$$\mathring{E}_p^s(t)\le C_bT^{p-2+\delp}I_0,$$
$$\int_t^{\infty}B_p^{\ul{22}}(\tau)d\tau \le C_bT^{p-2+\delp}I_0.$$
Then if $I_0$ is sufficiently small,
$$\mathring{E}_p^s(t)\lesssim T^{p-2+\delp}I_0.$$
\end{lemma}
\begin{proof}
Using the mean value theorem, for a given $t$, let $t'\in [t/2,t]$ be the value for which $\mathring{B}_{p+1}^{s+1}(t')=\frac2t\int_{t/2}^t\mathring{B}_{p+1}^{s+1}(\tau)d\tau$. Then using Corollary \ref{mr_NL_absorb_bulk},
\begin{align*}
\mathring{E}_p^{s}(t) &\lesssim \mathring{E}_p^{s}(t')+\int_{t'}^{t}\mathring{E}^s(\tau)B_p^{\ul{22}}(\tau)d\tau+(C_bI_0)^{1/2}\int_{t'}^{t}\mathring{E}_p^{s}(\tau)T^{(\delm-3+\delp)/2}d\tau \\
&\lesssim \mathring{E}_p^s(t')+I_0\int_{t'}^tB_p^{\ul{22}}(\tau)d\tau +(C_bI_0)^{3/2}\int_{t'}^t T^{p-2+\delp}T^{(\delm-3+\delp)/2}d\tau \\
&\lesssim \mathring{E}_p^s(t')+C_bI_0^2T^{p-2+\delp}+(C_bI_0)^{3/2}T^{p-2+\delp}T^{\delm+\delp-1/2} \\
&\lesssim \mathring{E}_p^s(t')+(C_bI_0+C_b^{3/2}I_0^{1/2})T^{p-2+\delp}I_0.
\end{align*}
Now by the choice of $t'$, and another application of Corollary \ref{mr_NL_absorb_bulk},
\begin{multline*}
\mathring{E}_p^s(t')\lesssim \mathring{B}_{p+1}^{s+1}(t')=\frac2t\int_{t/2}^t\mathring{B}_{p+1}^{s+1}(\tau)d\tau \\
\lesssim t^{-1}\mathring{E}_{p+1}^{s+1}(t/2)+t^{-1}\int_{t/2}^t\mathring{E}^s(\tau)B_{p+1}^{\ul{22}}(\tau)d\tau +t^{-1}(C_bI_0)^{1/2}\int_{t/2}^t \mathring{E}_{p+1}^{s+1}(\tau)T^{(\delm-3+\delp)/2}d\tau. \\
\lesssim t^{-1}(T^{(p+1)-2+\delp}I_0+C_bI_0^2T^{(p+1)-2+\delp}+(C_bI_0)^{3/2}T^{(p+1)-2+\delp}).
\end{multline*}
Thus,
$$E_p^{s,k}(t)\lesssim (1+C_bI_0+C_b^{3/2}I_0^{1/2})T^{p-2+\delp}I_0.$$
Provided $C_b^{3/2}I_0^{1/2}$ is sufficiently small, the proof is complete.
\end{proof}

By applying Lemma \ref{mr_NL_WE_decay} for $p=1-\delp$ and $s=28$, we obtain
$$\mathring{E}_{1-\delp}^{28}(t) \lesssim T^{-1}I_0.$$
By interpolation we obtain decay for all $p\in [1-\delp,2-\delp]$.
$$\mathring{E}_{p\in[1-\delp,2-\delp]}^{28}(t)\lesssim T^{p-2+\delp}I_0.$$
Then applying Lemma \ref{mr_NL_WE_decay} again for each $p\in [\delm,1-\delp]$ and $s=27$, we obtain
$$\mathring{E}_{p\in[\delm,2-\delp]}^{27}(t)\lesssim T^{p-2+\delp}I_0.$$
This recovers the bootstrap assumption (\ref{_25_ba}) at the highest level of derivatives.

\subsection{Recovering weak decay for $\mathring{E}_{p\in[\delm-1,\delm]}^{26}(t)$}\label{part_1_weak_decay_sec}

To prove estimates for low $p$ (ie. $p$ in the range $[\delm-1,\delm]$), we first observe that
$$\int_t^\infty \mathring{E}_p^{26}(\tau)d\tau \lesssim \int_t^\infty \mathring{B}_{p+1}^{27}(\tau)d\tau.$$
Then by Corollary \ref{mr_NL_absorb_bulk},
\begin{align*}
\int_t^\infty \mathring{B}_{p+1}^{27}(\tau)d\tau &\lesssim \mathring{E}_{p+1}^{27}(t)+\int_t^\infty \mathring{E}^{27}(\tau)B_{p+1}^{\ul{22}}(\tau)d\tau \\
&\hspace{1.9in}+(C_bI_0)^{1/2}\int_t^\infty\mathring{E}_{p+1}^{27}(\tau)T^{(\delm-3+\delp)/2}d\tau \\
&\lesssim T^{(p+1)-2+\delp}I_0+I_0\int_t^\infty B_{p+1}^{\ul{22}}(\tau)d\tau \\
&\hspace{1.5in}+(C_bI_0)^{1/2}\int_t^\infty T^{(p+1)-2+\delp}I_0T^{(\delm-3+\delp)/2}d\tau \\
&\lesssim (1+C_bI_0+(C_bI_0)^{1/2})T^{(p+1)-2+\delp}I_0.
\end{align*}
It follows that if $I_0$ is sufficiently small so that $C_bI_0\lesssim 1$, then
$$\int_t^\infty \mathring{E}_p^{26}(\tau)d\tau \lesssim T^{(p+1)-2+\delp}I_0 = T^{p-2+\delp +1}I_0.$$
This recovers the bootstrap assumption (\ref{_23_ba}) at the highest level of derivatives.

\subsection{The general case $s+2k\le 28$}\label{general_case_s_k_sec}

As explained in the introduction of this section, this is the point where we turn to the second of two parts of the remainder of the proof. In particular, we will recover the bootstrap assumptions (\ref{energy_ba}-\ref{_23_ba}) for the norms $E^{s,k}(t)$, $E_p^{s,k}(t)$, etc. We will use results from the first part, which started with \S\ref{homogeneous_case_sec}, to handle the homogeneous norm $\mathring{B}_p^s(t)$ appearing on the right side of many estimates in this second part.

The outline for this second part is similar to, but slightly more complicated than, the outline for the first part. To begin, in \S\ref{part_2_refined_estimates_sec}, we prove refined estimates for the nonlinear norms $N^{s,k}(t)$ and $N_p^{s,k}(t)$. Just like the estimates in \S\ref{part_1_refined_estimates_sec}, these estimates constitute the crucial step of the proof. The remainder of this second part uses a finite inductive argument with the inductive assumptions listed in \S\ref{part_2_inductive_assumptions_sec}. In \S\ref{part_2_Ep_boundedness_sec} and \S\ref{part_2_E_boundedness_sec}, the bootstrap assumptions (\ref{_29_ba}) and (\ref{energy_ba}) at the highest level $s+2k=28$ are recovered in that order. Then in \S\ref{part_2_decay_sec}, a decay lemma is proved and used to recover the bootstrap assumption (\ref{_25_ba}) at the level $s+2k=24$. Next, in \S\ref{establish_inductive_assumptions_sec} the inductive assumptions for the next step $k+1$ are established. Finally, in \S\ref{part_2_weak_decay_sec}, the bootstrap assumption (\ref{_23_ba}) at the level $s+2k=22$, which assumes a weaker form of decay, is recovered.

\subsection{Refined estimates for $N^{s,k}(t)$ and $N_p^{s,k}(t)$}\label{part_2_refined_estimates_sec}

The $L^\infty$ estimates given in Lemma \ref{simplified_pointwise_lemma} allow us to provide refined estimates for the nonlinear error terms. \textbf{This is the crucial step of the proof.}

\begin{lemma}\label{sk_refined_nl_lem}
In the context of the bootstrap assumptions provided in \S\ref{bootstrap_assumptions_sec}, if $s+2k\le 28$ and $C_bI_0\le 1$, then
$$N^{s,k}(t)\lesssim (E^{s,k}(t))^{1/2}\left((E^{s,k}(t))^{1/2}(E_{\delm-1}^{\ul{23}}(t))^{1/2}+(E^{s,k}_{1-\delm}(t))^{1/2}(E^{\ul{23}}_{\delm-1}(t))^{1/2}\right),$$
$$N_p^{s,k}(t)\lesssim E^{s,0}(t)B_p^{s/2+8,0}(t)+B_p^{s,0}(t)E^{s/2+8,0}(t)+E^{s,k}_{p'}(t)(E_{p''}^{\ul{23}}(t))^{1/2}.$$
\end{lemma}
\begin{remark}
The careful reader may notice that the low order energy norms in this lemma are at the level $s+2k=23$, whereas the low order energy norms in the analogous lemma for the homogeneous case (Lemma \ref{mr_refined_nl_lem}) are at the level $s+2k=22$. This difference is unimportant, but it is due to the fact that the $\pd_t$ commutators commute with the derivatives in the structures of $\mathcal{N}_\phi$ and $\mathcal{N}_\psi$, while the commutators $Q$ and $\tilde{Q}$ have angular parts that mix with the $\pd_\theta$ derivatives in some terms belonging to $\mathcal{N}_\phi$ and $\mathcal{N}_\psi$. This is manifest in the condition $s_1+...+s_{2+j}\le s+2$ that is given in Lemmas \ref{Nphi_s_k_structure_lem} and \ref{Npsi_s_k_structure_lem}. This means that the lower order factors can have up to $(s+2)/2=s/2+1$ derivatives, whereas in the homogeneous case, they can only have up to $s/2$ derivatives.
\end{remark}
\begin{proof}
First, we recall the definitions of $N^{s,k}(t)$ and $N_p^{s,k}(t)$ from Theorem \ref{p_s_k_thm}.
$$N^{s,k}(t)=(E^{s,k}(t))^{1/2}\left(\sum_{\substack{s'\le s \\ k'\le k}}||q^{-2}\tg^{k'}\Gamma^{s'}(q^2\mathcal{N}_\phi)||_{L^2(\Sigma_t)}+\sum_{\substack{s'\le s \\ k'\le k}}||q^{-2}\tg^{k'}\tilde\Gamma^{s'}(q^2\mathcal{N}_\psi)||_{L^2(\tilde{\Sigma}_t)}\right),$$
\begin{align*}
N_p^{s,k}(t) =& \sum_{\substack{s'\le s \\ k'\le k}}\int_{\Sigma_t}r^{p+1}|q^{-2}\tg^{k'}\Gamma^{s'}(q^2\mathcal{N}_\phi)|^2 + \sum_{\substack{s'\le s \\ k'\le k}}\int_{\tilde\Sigma_t}r^{p+1}|q^{-2}\tg^{k'}\tilde\Gamma^{s'}(q^2\mathcal{N}_\psi)|^2.
\end{align*}
Therefore, it suffices to prove the following three estimates.
\begin{multline*}
\sum_{\substack{s'\le s \\ k'\le k}}||q^{-2}\tg^{k'}\Gamma^{s'}(q^2\mathcal{N}_\phi)||_{L^2(\Sigma_t)}+\sum_{\substack{s'\le s \\ k'\le k}}||q^{-2}\tg^{k'}\tilde\Gamma^{s'}(q^2\mathcal{N}_\psi)||_{L^2(\tilde{\Sigma}_t)} \\
\lesssim (E^{s,k}(t))^{1/2}(E_{\delm-1}^{\ul{23}}(t))^{1/2}+(E^{s,k}_{1-\delm}(t))^{1/2}(E^{\ul{23}}_{\delm-1}(t))^{1/2}
\end{multline*}
\begin{multline*}
\sum_{\substack{s'\le s \\ k'\le k}}\int_{\Sigma_t\cap\{r>R\}}r^{p+1}|q^{-2}\tg^{k'}\Gamma^{s'}(q^2\mathcal{N}_\phi)|^2 + \sum_{\substack{s'\le s \\ k'\le k}}\int_{\tilde\Sigma_t\cap\{r>R\}}r^{p+1}|q^{-2}\tg^{k'}\tilde\Gamma^{s'}(q^2\mathcal{N}_\psi)|^2 \\
\lesssim E^{s,0}(t)B_p^{s/2+8,0}(t)+B_p^{s,0}(t)E^{s/2+8,0}(t)
\end{multline*}
\begin{multline*}
\sum_{\substack{s'\le s \\ k'\le k}}\int_{\Sigma_t\cap\{r<R\}}|q^{-2}\tg^{k'}\Gamma^{s'}(q^2\mathcal{N}_\phi)|^2 + \sum_{\substack{s'\le s \\ k'\le k}}\int_{\tilde\Sigma_t\cap\{r<R\}}|q^{-2}\tg^{k'}\tilde\Gamma^{s'}(q^2\mathcal{N}_\psi)|^2 \\
\lesssim E^{s,k}_{p'}(t)(E_{p''}^{\ul{23}}(t))^{1/2}
\end{multline*}
These estimates are analogous to the estimates in the proof of Lemma \ref{mr_refined_nl_lem}, and they can be proved the same way.
\end{proof}

\begin{corollary}\label{sk_NL_absorb_bulk}
In the context of the bootstrap assumptions provided in \S\ref{bootstrap_assumptions_sec}, Theorem \ref{p_s_k_thm} and Lemma \ref{sk_refined_nl_lem} imply that if $s+2k\le 28$ and $C_bI_0$ is sufficiently small, then
\begin{multline*}
E_p^{s,k}(t_2)+\int_{t_1}^{t_2}B_p^{s,k}(t)dt \\
\lesssim E_p^{s,k}(t_1)+\int_{t_1}^{t_2}\mathring{B}_p^{s+1}(t)dt+(C_bI_0)^{1/2}\int_{t_1}^{t_2}E_p^{s,k}(t)T^{(\delm-3+\delp)/2}dt \\
+\sum_{\substack{s'+2k'\le s+2k \\ k'<k}}\int_{t_1}^{t_2}B_p^{s',k'}(t)dt.
\end{multline*}
whenever $s/2+8\le s+2k$ (which means $16\le s+4k$).
\end{corollary}
\begin{proof}
According to Theorem \ref{p_s_k_thm},
$$E_p^{s,k}(t_2)+\int_{t_1}^{t_2}B_p^{s,k}(t)dt\lesssim E_p^{s,k}(t_1)+\int_{t_1}^{t_2}B_p^{s+2,k-1}(t)+\mathring{B}_p^{s+1}(t)+N_p^{s,k}(t)dt.$$
By Lemma \ref{sk_refined_nl_lem},
$$N_p^{s,k}(t)\lesssim E^{s,0}(t)B_p^{s/2+8,0}(t)+B_p^{s,0}(t)E^{s/2+8,0}(t)+E^{s,k}_p(t)(E_{\delm-1}^{\ul{23}}(t))^{1/2}.$$
We now estimate each of the three terms on the right side. The strategy for the first term depends on whether $k=0$ or $k>0$. If $k=0$, then since $s/2+8\le s$ by assumption,
$$E^{s,0}(t)B_p^{s/2+8,0}(t) \lesssim C_bI_0 B_p^{s,0}(t).$$
It follows that if $C_bI_0$ is sufficiently small, this term can be absorbed into the bulk term on the left side. If instead $k> 0$, then since $s/2+8\le s+2k$ by assumption,
$$E^{s,0}(t)B_p^{s/2+8,0}(t) \lesssim C_bI_0 \sum_{\substack{s'+2k'\le s+2k \\ k'<k}}B_p^{s',k'}(t).$$
It follows that if $I_0$ is sufficiently small so that $C_bI_0\le 1$, then this term can be estimated by the lower order (in $k$) term on the right side.

The strategy for the second term is similar. Since $s+2k\le 28$, it follows that $s/2+8\le 28$, so
$$B_p^{s,0}(t)E^{s/2+8,0}(t)\lesssim C_bI_0B_p^{s,0}(t).$$
Again, if $k=0$, then $C_bI_0$ must be taken sufficiently small so that this term can be absorbed into the bulk term on the left side. If instead $k>0$, then since $s\le s+2k$, as long as $C_bI_0\le 1$, this term can be estimated by the lower order (in $k$) term on the right side.

Finally, by the weak decay lemma,
$$\int_{t_1}^{t_2}E_p^{s,k}(t)(E_{\delm-1}^{\ul{23}}(t))^{1/2}dt \lesssim (C_bI_0)^{1/2}\int_{t_1}^{t_2}E_p^{s,k}(t)T^{(\delm-3+\delp)/2}dt.$$
This completes the proof.
\end{proof}

\subsection{Inductive assumptions}\label{inductive_assumptions_sec}\label{part_2_inductive_assumptions_sec}

The remainder of the proof is a finite induction argument. First, estimates are proved for $k=0$, and then for $k=1$, etc. until $k=14$ (which saturates $s+2k\le 28$). For each $k$, it will be necessary to use estimates established for $k-1$. These inductive assumptions are listed here.

Either
$$k=0,$$
or $s+2k=28$ and
$$\sum_{\substack{s'+2k'\le 28 \\ k'<k}}\int_t^\infty B_{p\in[\delm,2-\delp]}^{s',k'}(\tau)d\tau \lesssim I_0,$$
or $s+2k=26$ and
$$\sum_{\substack{s'+2k'\le 26 \\ k'<k}}\int_t^\infty B_{p\in[1-\delp,2-\delp]}^{s',k'}(\tau)d\tau \lesssim T^{p-2+\delp}I_0,$$
or $s+2k=24$ and
$$\sum_{\substack{s'+2k'\le 24 \\ k'<k}}\int_t^\infty B_{p\in[\delm,2-\delp]}^{s',k'}(\tau)d\tau \lesssim T^{p-2+\delp}I_0.$$

For the remainder of the proof, $k$ should be considered fixed. These estimates will be used for the fixed $k$, and eventually (in \S\ref{establish_inductive_assumptions_sec}) the corresponding estimates obtained by replacing $k$ with $k+1$ will be proved, thus closing the induction argument.

\subsection{Recovering boundedness of $E^{s,k}_{2-\delp}(t)$ ($s+2k=28$)}\label{part_2_Ep_boundedness_sec}

Our first application of Corollary \ref{sk_NL_absorb_bulk} is to prove boundedness of $E^{s,k}_{2-\delp}(t)$ and $\int_0^t B_{2-\delp}^{s,k}(\tau)d\tau$. With $s+2k=28$,
\begin{align*}
E_{2-\delp}^{s,k}(t)&+\int_0^t B_{2-\delp}^{s,k}(\tau)d\tau \\
&\lesssim E_{2-\delp}^{s,k}(0)+\int_0^t\mathring{B}^{29}(\tau)d\tau+(C_bI_0)^{1/2}\int_0^tE_{2-\delp}^{s,k}(\tau)T^{(\delm-3+\delp)/2}d\tau \\
&\hspace{3in}+\sum_{\substack{s'+2k'\le s+2k \\ k'<k}}\int_0^tB_{2-\delp}^{s',k'}(\tau)d\tau \\
&\lesssim I_0+I_0+(C_bI_0)^{1/2}\int_0^tC_bT^{(\delm-3+\delp)/2}I_0d\tau +I_0 \\
&\lesssim (1+C_b^{3/2}I_0^{1/2})I_0.
\end{align*}
It follows that if $I_0$ is sufficiently small so that $C_b^3I_0\lesssim 1$, then
$$E_{2-\delp}^{s,k}(t)+\int_0^tB_{2-\delp}^{s,k}(\tau)d\tau\lesssim I_0.$$
This recovers the bootstrap assumption (\ref{_29_ba}) at the level $k$.

\subsection{Recovering boundedness of $E^{s,k}(t)$ ($s+2k=28$)}\label{part_2_E_boundedness_sec}

Since
$$\int_0^tB_1^{s,k}(\tau)d\tau \lesssim \int_0^tB_{2-\delp}^{s,k}(\tau)d\tau \lesssim I_0,$$
we are now able to prove that $E^{s,k}(t)$ is bounded.

Let $s+2k=28$. By Theorem \ref{p_s_k_thm} and Lemma \ref{sk_refined_nl_lem},
\begin{align*}
E^{s,k}(t) \lesssim & E^{s,k}(0)+\int_0^tB_{2-\delp}^{s+2,k-1}(\tau)+\mathring{B}_1^{s+1}(\tau)+B_1^{s,k}(\tau)+N^{s,k}(\tau)d\tau \\
\lesssim & E^{s,k}(0)+\int_0^tB_{2-\delp}^{s+2,k-1}(\tau) + \mathring{B}_1^{s+1}(\tau) + B_1^{s,k}(\tau)d\tau\\
&+\int_0^t(E^{s,k}(\tau))^{1/2}\left((E^{s,k}(\tau))^{1/2}(E_{\delm-1}^{\ul{22}}(\tau))^{1/2}+(E^{s,k}_{1-\delm}(\tau))^{1/2}(E^{\ul{22}}_{\delm-1}(\tau))^{1/2}\right)d\tau \\
\lesssim & I_0+\int_0^t(C_bI_0)^{1/2}(C_bI_0)^{1/2}(C_bT^{\delm-3+\delp}I_0)^{1/2}d\tau \\
\lesssim & (1+C_b^{3/2}I_0^{1/2})I_0.
\end{align*}
In particular, we used the weak decay principle in the third step. It follows that if $I_0$ is chosen sufficiently small so that $C_b^3I_0\lesssim 1$, then
$$E^{s,k}(t)\lesssim I_0.$$
This recovers the bootstrap assumption (\ref{energy_ba}) at the level $k$.

\subsection{Proving decay for $E^{s,k}_{p\in[1-\delp,2-\delp]}(t)$ ($s+2k=26$) and $E^{s,k}_{p\in[\delm,2-\delp]}(t)$ ($s+2k=24$)}\label{part_2_decay_sec}

Once again, we prove a decay lemma for repeated use.
\begin{lemma}\label{sk_NL_WE_decay}
Suppose  $p+1,p\in [\delm,2-\delp]$ and
$$\int_t^\infty \mathring{B}_{p+1}^{s+3}(\tau)d\tau \lesssim T^{(p+1)-2+\delp}I_0,$$
$$E_{p+1}^{s+2,k}(t)\lesssim T^{(p+1)-2+\delp}I_0,$$
$$\sum_{\substack{s'+2k'\le s+2k+2 \\ k'<k}}\int_t^\infty B_{p+1}^{s',k'}(\tau)d\tau\lesssim T^{(p+1)-2+\delp}I_0.$$
$$\int_t^\infty \mathring{B}_p^{s+1}(\tau)d\tau \lesssim T^{p-2+\delp}I_0,$$
$$E_p^{s,k}(t)\le C_bT^{p-2+\delp}I_0,$$
$$\sum_{\substack{s'+2k'\le s+2k \\ k'<k}}\int_t^\infty B_p^{s',k'}(\tau)d\tau\lesssim T^{p-2+\delp}I_0,$$
Then if $I_0$ is sufficiently small,
$$E_p^{s,k}(t)\lesssim T^{p-2+\delp}I_0.$$
\end{lemma}

\begin{proof}
Using the mean value theorem, for a given $t$, let $t'\in [t/2,t]$ be the value for which $B_{p+1}^{s+2,k}(t')=\frac2t\int_{t/2}^tB_{p+1}^{s+2,k}(\tau)d\tau$. Then using Corollary \ref{sk_NL_absorb_bulk},
\begin{align*}
E_p^{s,k}(t) &\lesssim E_p^{s,k}(t')+\int_{t'}^t\mathring{B}_p^{s+1}(\tau)d\tau+(C_bI_0)^{1/2}\int_{t'}^tE_p^{s,k}(\tau)T^{(\delm-3+\delp)/2}d\tau \\
&\hspace{3in}+\sum_{\substack{s'+2k'\le s+2k \\ k'<k}}\int_{t'}^tB_p^{s',k'}(\tau)d\tau  \\
&\lesssim E_p^{s,k}(t')+C_bI_0\int_{t'}^t(C_bT^{p-2+\delp}I_0)T^{\delm-2+\delp}d\tau+T^{p-2+\delp}I_0 \\
&\lesssim E_p^{s,k}(t')+C_b^2I_0T^{p-2+\delp}I_0+T^{p-2+\delp}I_0.
\end{align*}
Now by the choice of $t'$,
\begin{multline*}
E_p^{s,k}(t')\lesssim B_{p+1}^{s+2,k}(t')=\frac2t\int_{t/2}^tB_{p+1}^{s+2,k}(\tau)d\tau \\
\lesssim t^{-1}E_{p+1}^{s+2,k}(t/2)+t^{-1}\int_{t/2}^t\mathring{B}_{p+1}^{s+3}(\tau)d\tau+t^{-1}C_bI_0\int_{t/2}^tE_{p+1}^{s+2,k}(\tau)T^{\delm-2+\delp}d\tau \\
+t^{-1}\sum_{\substack{s'+2k'\le s+2k+2 \\ k'<k}}\int_{t/2}^t B_{p+1}^{s',k'}(\tau)d\tau.
\end{multline*}
Thus,
$$E_p^{s,k}(t)\lesssim (1+C_b^2I_0)T^{p-2+\delp}I_0.$$
Taking $I_0$ sufficiently small so that $C_b^2I_0\lesssim 1$ completes the proof.
\end{proof}

By applying Lemma \ref{sk_NL_WE_decay} for $p=1-\delp$ and $s+2k=26$, we obtain
$$E_{1-\delp}^{s,k}(t)\lesssim T^{-1}I_0.$$
By interpolation we obtain decay for all $p\in [1-\delp,2-\delp]$.
$$E_{p\in [1-\delp,2-\delp]}^{s,k}(t)\lesssim T^{p-2+\delp}I_0.$$
Then applying Lemma \ref{sk_NL_WE_decay} again for each $p\in[\delm,1-\delp]$ and $s+2k=24$, we obtain
\begin{equation*}
E_{p\in [\delm,2-\delp]}^{s,k}(t)\lesssim T^{p-2+\delp}I_0.
\end{equation*}
This recovers the bootstrap assumption (\ref{_25_ba}) at the level $k$.

\subsection{Establishing inductive assumptions for $k+1$}\label{establish_inductive_assumptions_sec}

By Corollary \ref{sk_NL_absorb_bulk},
\begin{multline*}
\int_t^\infty B_p^{s,k}(\tau)d\tau \lesssim E_p^{s,k}(t)+\int_t^\infty\mathring{B}_p^{s+1}(\tau)d\tau+(C_bI_0)^{1/2}\int_t^\infty E_p^{s,k}(\tau)T^{(\delm-3+\delp)/2}d\tau \\
+\sum_{\substack{s'+2k'\le s+2k \\ k'<k}}\int_t^\infty B_p^{s',k'}(\tau)d\tau.
\end{multline*}
The quantities $E_p^{s,k}(t)$, $\int_t^\infty \mathring{B}_p^{s+1}(\tau)d\tau$, and $\int_t^\infty B_p^{s',k'}(\tau)d\tau$ ($s'+2k'\le s+2k$ and $k'<k$) all have the same proven decay rates. For $s+2k=28$ and $p\in[\delm,2-\delp]$, they are bounded in time by $I_0$, for $s+2k=26$ and $p\in[1-\delp,2-\delp]$, they decay at least as fast as $T^{p-2+\delp}I_0$, and for $s+2k=24$ and $p\in[\delm,2-\delp]$, they decay at least as fast as $T^{p-2+\delp}I_0$. Thus,
$$\sum_{\substack{s'+2k'\le 28 \\ k'<k+1}}\int_t^\infty B_{p\in[\delm,2-\delp]}^{s',k'}(\tau)d\tau \lesssim I_0,$$
$$\sum_{\substack{s'+2k'\le 26 \\ k'<k+1}}\int_t^\infty B_{p\in[1-\delp,2-\delp]}^{s',k'}(\tau)d\tau \lesssim T^{p-2+\delp}I_0,$$
$$\sum_{\substack{s'+2k'\le 24 \\ k'<k+1}}\int_t^\infty B_{p\in[\delm,2-\delp]}^{s',k'}(\tau)d\tau \lesssim T^{p-2+\delp}I_0.$$
These are the inductive assumptions at the next level $k+1$.

\subsection{Recovering weak decay for $E_{p\in[\delm-1,\delm]}^{s,k}(t)$ ($s+2k=22$)}\label{part_2_weak_decay_sec}

Finally, set $s+2k=22$ and observe that for $p\in [\delm-1,\delm]$,
$$\int_t^\infty E_p^{s,k}(\tau)d\tau \lesssim \int_t^\infty B_{p+1}^{s+2,k}(\tau)d\tau\lesssim T^{(p+1)-2+\delp}I_0=T^{p-2+\delp+1}I_0.$$
This recovers the bootstrap assumption (\ref{_23_ba}) at the level $k$.

\noindent This completes the proof of Theorem \ref{main_thm}. \qed

\appendix
\chapter{More about the wave map problem}\label{wave_map_appendix}

\section{The general theory of wave maps and the $\xi_a$ system}\label{wave_map_general_theory_sec}

There are at least two useful ways to linearize the wave map system (\ref{wm_X_eqn}-\ref{wm_Y_eqn}) about the nontrivial solution $(X_0,Y_0)=(A,B)$. One of these ways is to introduce a vector bundle formalism (motivated by the geometric nature of the wave map system) and derive an equation for a section $\xi_a$ of this bundle, $\mathcal{B}$. This approach is useful because it suggests an appropriate grouping of terms when deriving a Morawetz estimate. We derive this system here.

\subsection{The general theory of wave maps}
Let $\Phi$ be a map.
$$\Phi:(M,g)\rightarrow (N,h).$$
Let $\Phi_*=d\Phi$ be its pushforward.
$$\Phi_*:T_pM\rightarrow T_{\Phi(p)}N$$
The wave map equation says that
\begin{equation}\label{wave_map_equation}
div\Phi_*=0.
\end{equation}

Let us take a moment to outline a few conventions for index notation. We will use greek indices to represent tensor quantities on $M$ and lower-case latin indices to represent tensor quantities on $N$. We emphasize that these shall be used for \textit{coordinate invariant} quantities only. When using a particular set of coordinates, we will use primed indices instead. For example, we can represent the pushforward $d\Phi$ by
$$d\Phi=d\Phi^a_\mu,$$
but if we specify the map $\Phi$ in coordinates $\Phi^{a'}(x^{\mu'})$, then we have the coordinate dependent equation
$$d\Phi^{a'}_{\mu'}=\pd_{\mu'}\Phi^{a'}.$$
As an exception to the primed index rule, indices $i,j,k,l$ will be used to correspond to a particular orthonormal frame for the bundle $\mathcal{B}$, which will be introduced later.

We use $\nabla$ to denote the Levi-Civita connections on both $M$ and $N$. The particular connection being used will be clear from context or use of indices. The pushforward $\Phi_*$ and the Levi-Civita connection on $N$ induce a connection $D$ taking a vector $\vec{V}\in T_pM$ to a differential operator $D_{\vec{V}}$ acting on tensors on $N$ by
$$D_{\vec{V}}=\nabla_{\Phi_* \vec{V}}.$$

When the context is clear, we use indices to implicitly push forward contravariant tensors on $M$ or pull back covariant tensors on $N$. Thus, for example, 
$$g^{ab}=d\Phi_\alpha^a d\Phi_{\beta}^bg^{\alpha\beta}$$
and
$$h_{\alpha\beta}=d\Phi_\alpha^a d\Phi_\beta^b h_{ab}.$$
We will never use the inverse of $\Phi$ (which may not exist) to push forward covariant tensors or pull back contravariant tensors. This allows for raising and lowering indices without ambiguity. For example, 
$$h^{\alpha\beta}=g^{\alpha\gamma}g^{\beta\delta}h_{\gamma\delta}=g^{\alpha\gamma}g^{\beta\delta}d\Phi_\gamma^c d\Phi_\delta^d h_{cd}.$$
The one drawback of this convention is that indices no longer clarify whether tensors naturally belong to $M$ or $N$. (For example, is $R_{\alpha\beta\gamma\delta}$ the Riemann curvature tensor for $M$ or the pullback of the Riemann curvature tensor for $N$?) This ambiguity must be resolved explicitly when the tensor is introduced, but is not a serious issue in practice.

Using the above index notation, one can directly calculate the following coordinate-dependent equation.
$$D_{\mu'} d\Phi^{a'}_{\nu'}=\pd_{\mu'}\pd_{\nu'}\Phi^{a'}-\Gamma_{\mu'\nu'}^{\lambda'}\pd_{\lambda'}\Phi^{a'}+\pd_{\mu'}\Phi^{b'}\pd_{\nu'}\Phi^{c'}\Gamma_{b'c'}^{a'}.$$
In particluar, the wave map equation (\ref{wave_map_equation}) in coordinates takes the form
\begin{equation}\label{general_wave_eqn_in_coordinates_eqn}
\Box_g\Phi^{a'}+g^{\mu'\nu'}\pd_{\mu'}\Phi^{b'}\pd_{\nu'}\Phi^{c'}\Gamma_{b'c'}^{a'}=0.
\end{equation}

\subsection{Linearized wave maps and the section $\xi_a$}\label{wave_map_bundle_sec}

To examine the linear stability of a solution $\Phi$ to a wave map, one obtains an equation for a vectorfield 
$$\vec{\psi}:M\rightarrow TN,\hspace{.5in}\vec{\psi}(p)\in T_{\Phi(p)}N.$$
The significance of this vectorfield is that if $\Phi(s)$ is a parametrized family of solutions to the wave map equation with $\Phi(0)=\Phi$, then $\vec{\psi}=\left.\frac{d}{ds}\right|_{s=0}\Phi(s)$.

The equation for $\vec{\psi}$ is
\begin{equation}\label{vec_psi_linear_stability_eqn}
\Box_g\psi^a+R^{\lambda a}{}_{\lambda b}\psi^b=0,
\end{equation}
where $R$ is the curvature tensor of the target manifold
$$R^{\lambda a}{}_{\lambda b}=g^{\gamma\delta}d\Phi_\gamma^c d\Phi_\delta^d R_c{}^a{}_{db}.$$

By comparing the model system (\ref{wm_X_eqn}-\ref{wm_Y_eqn}) to the general wave map equation (\ref{general_wave_eqn_in_coordinates_eqn}) one can read off the Christoffel symbols $\Gamma^{a'}_{b'c'}$ and determine that the target manifold $(N,h)$ for the model problem is the hyperbolic half plane--with negative constant curvature. Letting $\epsilon_{ab}$ be the volume form for the hyperbolic plane, we have
$$R_{abcd}=-\epsilon_{ab}\epsilon_{cd}.$$
Contracting equation (\ref{vec_psi_linear_stability_eqn}) with $\epsilon_{ab}$, and using the facts that $D_\mu \epsilon_{ab}=0$ and $\epsilon_{ac}\epsilon^c{}_b=h_{ab}$, we obtain a new equation.
$$\Box_g(\epsilon_{ab}\psi^b)-g_a{}^b\epsilon_{bc}\psi^c=0.$$
Therefore, we introduce the new dynamic quantity
$$\xi_a=\epsilon_{ab}\psi^b$$
and the potential
$$V^{ab}=g^{ab}=g^{\alpha\beta}d\Phi_\alpha^a d\Phi_\beta^b,$$
and the equation (\ref{vec_psi_linear_stability_eqn}) becomes
\begin{equation}\label{xi_a_linearized_eqn}
\Box_g\xi_a-V_a{}^b\xi_b=0.
\end{equation}
We refer to the vector bundle for this equation as $\mathcal{B}$.

\subsection{The equations for $\xi_a$ in component form}

Assign to the hyperbolic half plane $(N,h)$ the coordinates $X$ and $Y$ with the ranges
\begin{align*}
X &\in (0,\infty) \\
Y &\in (-\infty,\infty).
\end{align*}
The metric for the hyperbolic half plane is given by
$$ds^2= \frac{1}{X^2}(dX^2+dY^2).$$

An orthonormal frame is given by
$$e_1=X\pd_X,\hspace{.5in} e_2=X\pd_Y.$$
We can represent the section $\xi_a$ in terms of the dual frame $\{e^1,e^2\}$.
$$\xi_a=\xi_i(e^i)_a=\xi_1(e^1)_a+\xi_2(e^2)_a.$$
The functions $\xi_1$ and $\xi_2$ are scalar quantities and are the objects of study in \S\ref{wk:xi_estimates_sec}. In particular, these quantities linearize the wave map system (\ref{wm_X_eqn}-\ref{wm_Y_eqn}) in the following way.
\begin{align*}
X &= A-A\xi_2 \\
Y &= B+A\xi_1.
\end{align*}

The pushforward $d\Phi$ is given by
\begin{align*}
d\Phi_\alpha^a = \frac{\pd_\alpha A}{A} (e_1)^a+\frac{\pd_\alpha B}{A} (e_2)^a,
\end{align*}
and the Christoffel symbols for the dual frame can be read from the following relations
\begin{align*}
D_\alpha e^1 &= -\frac{\pd_\alpha B}{A} e^2 \\
D_\alpha e^2 &= \frac{\pd_\alpha B}{A} e^1.
\end{align*}
Equation (\ref{xi_a_linearized_eqn}) in component form is given by the following system of equations for the scalar components $\xi_1$ and $\xi_2$.
\begin{align*}
\Box_g\xi_1 &= - 2\frac{\pd^\alpha B}{A}\pd_\alpha\xi_2+\frac{\pd^\alpha A\pd_\alpha A+\pd^\alpha B\pd_\alpha B}{A^2}\xi_1 ,\\
\Box_g\xi_2 &= 2\frac{\pd^\alpha B}{A}\pd_\alpha\xi_1 +2\frac{\pd^\alpha B\pd_\alpha B}{A^2}\xi_2 + 2\frac{\pd^\alpha A\pd_\alpha B}{A^2}\xi_1.
\end{align*}
In the Schwarzschild case, where $A=r^2\sin^2\theta$ and $B=0$, these equations reduce to
\begin{align*}
\Box_g\xi_1 &= \frac{4}{r^2}\left(1-\frac{2M}r\right)\xi_1+\frac{4\cot^2\theta}{r^2}\xi_1  \\
\Box_g\xi_2 &= 0. 
\end{align*}

\section{A proof of equivalence for the three systems}\label{three_systems_equivalent_sec}

The wave map problem is explained using three different dynamic quantities $(X,Y)$, $(\phi,\psi)$, and $(\xi_1,\xi_2)$. Here, we collect the facts, which are scattered throughout this thesis, that the three systems are equivalent.

\begin{proposition}
Let
\begin{align*}
X &=  A+A\phi = A-A\xi_2 \\
Y &= B+A^2\psi = B+A\xi_1.
\end{align*}
Then the systems for $(X,Y)$, $(\phi,\psi)$, and $(\xi_1,\xi_2)$ are related in the following way.

Let the pair $(X,Y)$ satisfy
\begin{align*}
X\Box_gX &= \pd^\alpha X\pd_\alpha X-\pd^\alpha Y\pd_\alpha Y \\
X\Box_gY &= 2\pd^\alpha X\pd_\alpha Y.
\end{align*}
Then the pair $(\phi,\psi)$ satisfies
\begin{align*}
\Box_g\phi &= \mathcal{L}_\phi + \mathcal{N}_\phi \\
\Box_{\tilde{g}}\psi &= \mathcal{L}_\psi + \mathcal{N}_\psi,
\end{align*}
where $\mathcal{L}_\phi$ and $\mathcal{L}_\psi$ are linear in $(\phi,\psi)$ while $\mathcal{N}_\phi$ and $\mathcal{N}_\psi$ are nonlinear in $(\phi,\psi)$. They are defined in \S\ref{wk:phi_psi_equations_intro_sec}.

Furthermore, the linearized equations
\begin{align*}
\Box_g\phi &= \mathcal{L}_\phi \\
\Box_{\tilde{g}}\psi &= \mathcal{L}_\psi,
\end{align*}
are equivalent to the $\mathcal{B}$-valued equation
$$\Box_g\xi_a-V_a{}^b\xi_b = 0,$$
provided $\xi_a=\xi_1 (e^1)_a+\xi_2 (e^2)_a$.
\end{proposition}
\begin{proof}
By Proposition \ref{N_identities_prop}, the system
\begin{align*}
\Box_g\phi &= \mathcal{L}_\phi+\mathcal{N}_\phi \\
\Box_{\tilde{g}}\psi &= \mathcal{L}_\psi+\mathcal{N}_\psi
\end{align*}
derives directly from the system
\begin{align*}
X\Box_gX &= \pd^\alpha X\pd_\alpha X-\pd^\alpha Y\pd_\alpha Y \\
X\Box_gY &= 2\pd^\alpha X\pd_\alpha Y
\end{align*}
when substituting $X=A+A\phi$ and $Y=B+A^2\psi$.

By Lemma \ref{translate_nl_lem},
\begin{align*}
(e_1)^a(\Box_g\xi_a-V_a{}^b\xi_b) &= A(\Box_{\tilde{g}}\psi -\mathcal{L}_\psi) \\
-(e_2)^a(\Box_g\xi_a-V_a{}^b\xi_b) &= \Box_g\phi-\mathcal{L}_\phi.
\end{align*}
Thus, the linearized equations
\begin{align*}
\Box_g\phi &= \mathcal{L}_\phi \\
\Box_{\tilde{g}}\psi &= \mathcal{L}_\psi,
\end{align*}
are equivalent to the $\mathcal{B}$-valued equation
$$\Box_g\xi_a-V_a{}^b\xi_b = 0.$$
\end{proof}

\section{The behavior of $X$ and $Y$ near the axis}\label{ernst_potential_axis_derivation}

Each of the coordinates $(X,Y)$ of the map $\Phi$ has a geometric interpretation. Consider an axisymmetric vacuum spacetime (ie. a spacetime solving (\ref{eve})) and let $K$ be its axial Killing vectorfield. The coordinate $X$ satisfies $X=K^\mu K_\mu$. Therefore, one expects $X$ to vanish on the axis at the same rate as $A$ in Kerr, and that if $K$ is everywhere spacelike, then $X$ should be non-negative.

The coordinate $Y$ satisfies $\nabla_\mu Y=\epsilon_{\mu\nu\kappa\lambda}K^\nu\nabla^\kappa K^\lambda,$ where $\epsilon$ is the volume form for the spacetime. For this reason, it is called the \textit{Ernst potential}. \cite{weinstein} The geometric significance of $Y$ can be understood the following way. Let $\gamma$ be a curve in the spacetime defined on the interval $[0,1]$. Let $S_\gamma$ be the orbit of $\gamma$ under the $U(1)$ isometry group generated by $K$. Then $S_\gamma$ is a $2$-surface parametrized by the parameter for $\gamma$ and the parameter for the $U(1)$ isometry group. We calculate
\begin{multline*}
\int_{S_\gamma} dA_{\kappa\lambda}\nabla^\kappa K^\lambda=\int_0^1\int_0^{2\pi} (\dot{\gamma}^\mu K^\nu \epsilon_{\mu\nu\kappa\lambda})\nabla^\kappa K^\lambda=2\pi\int_0^1\dot{\gamma}^\mu\nabla_\mu Y \\
=2\pi\left(Y(\gamma(1))-Y(\gamma(0))\right).
\end{multline*}

If $\gamma$ is closed or has both endpoints on the axis, then $S_\gamma$ is the boundary of a domain $\Omega$
$$S_\gamma=\pd\Omega$$
 and by Stokes' theorem,
$$\oint_{S_\gamma}dA_{\kappa\lambda}\nabla^\kappa K^\lambda=2\int_\Omega dV_\mu Ric^\mu{}_\nu K^\nu.$$
For a vacuum spacetime, the integral $\oint_{S_\gamma}dA_{\kappa\lambda}\nabla^\kappa K^\lambda$ computes the angular momentum enclosed by the surface $S_\gamma$. Taking $\gamma$ to be a curve starting and ending on the same half of the axis, we conclude that $Y$ must be constant on the same half of the axis. Taking $\gamma$ to be a curve starting and ending on opposite halves of the axis, we conclude that the values of $Y$ on each half axis differ by a constant related to the angular momentum of the black hole.

Let $(\delta X$, $\delta Y)$ correspond to a perturbation of the map $\Phi$
\begin{align*}
X &= X_0+\delta X \\
Y &= Y_0+\delta Y
\end{align*}
preserving angular momentum. To be consistent with the above geometric properites, both $\delta X$ and $\delta Y$ must both vanish on the axis.\footnote{Actually, $\delta Y$ only needs to be constant on the entire axis, but it is easy to check that $Y$ can be defined up to a constant, so it can be assumed that $\delta Y=0$ on the axis.} Assuming low regularity, since $\delta X$ and $\delta Y$ are axisymmetric, they must vanish at least to order $\sin^2\theta$ on the axis. Thus, both scalars 
\begin{align*}
\xi_1 &=A^{-1}\delta Y \\
\xi_2 &=-A^{-1}\delta X
\end{align*}
should be regular on the axis for regular perturbations of $\Phi$ which preserve angular momentum.

The equation for $\xi_1$ tells us more about $\xi_1$ on the axis. If we write $\xi_1=f(r)+\sin^2\theta g(r,\theta)$ for regular functions $f,g$ and assume $\xi_2$ is regular on the axis, then the equation
$$\sin^2\theta\Box_g\xi_1=\sin^2\theta g_1{}^i\xi_i$$
simplifies when $\sin\theta=0$ to
\begin{align*}
0&=\sin^2\theta g^{\theta\theta}d\Phi_\theta^1d\Phi_\theta^1(f(r)+\sin^2\theta g(r,\theta)) \\
&=\sin^2\theta \frac{1}{q^2}(2\cot\theta)^2(f(r)+\sin^2\theta g(r,\theta)) \\
&=\frac{4}{q^2}f(r).
\end{align*}
This means that we can write
$$\xi_1=r^2\sin^2\theta \psi$$
for some function $\psi$ that is regular on the axis. This is consistent with the class of perturbations defined in the main theorems of Chapters \ref{wm_szd_chap} and \ref{wm_kerr_chap}.

\section{The local theory for the $(\phi,\psi)$ system}

Due to the presence of nonlinear terms with factors that are singular on the axis, it is not obvious that the system (\ref{phi_nonlinear_system_eqn}-\ref{psi_nonlinear_system_eqn}) repeated below admits solutions even locally in time. A proof of this fact is sketched here.

First, we recall the system (\ref{phi_nonlinear_system_eqn}-\ref{psi_nonlinear_system_eqn}).
\begin{align*}
\Box_g\phi &= \mathcal{L}_\phi+\mathcal{N}_\phi \\
\Box_{\tilde{g}}\psi &= \mathcal{L}_\psi+\mathcal{N}_\psi,
\end{align*}
where
$$\mathcal{L}_\phi=-2\frac{\pd^\alpha B}{A}A\pd_\alpha \psi + 2\frac{\pd^\alpha B\pd_\alpha B}{A^2}\phi-4\frac{\pd^\alpha A\pd_\alpha B}{A^2} A\psi$$
$$\mathcal{L}_\psi=-2\frac{\pd^\alpha A_2}{A_2}\pd_\alpha\psi+2\frac{\pd^\alpha B\pd_\alpha B}{A^2}\psi + 2A^{-1}\frac{\pd^\alpha B}{A}\pd_\alpha\phi,$$
and
\begin{align*}
(1+\phi)\mathcal{N}_\phi &= \pd^\alpha \phi \pd_\alpha \phi - A\pd^\alpha \psi A\pd_\alpha \psi +2\frac{\pd^\alpha B}{A} \phi A\pd_\alpha \psi -4\frac{\pd^\alpha A}{A}A\psi A\pd_\alpha \psi \\
&\hspace{.5in}-\frac{\pd^\alpha B\pd_\alpha B}{A^2}\phi^2+4\frac{\pd^\alpha A\pd_\alpha B}{A^2}\phi A\psi -4\frac{\pd^\alpha A\pd_\alpha A}{A^2}(A\psi)^2,
\end{align*}
$$(1+\phi)\mathcal{N}_\psi = 2\pd^\alpha \phi \pd_\alpha\psi  +4\frac{\pd^\alpha A}{A}\psi\pd_\alpha \phi -2\frac{\pd^\alpha B}{A}A^{-1}\phi\pd_\alpha\phi.$$

\subsection{The iteration scheme}

We begin by defining an iteration scheme based on a given set of initial data $(\phi,\psi)|_{t=0}$. The iteration is defined as follows.

\begin{definition}
Let $(\phi_1,\psi_1)$ be a solution to the linear equations
\begin{align*}
\Box_g\phi_1 &= \mathcal{L}_{\phi_1} \\
\Box_{\tilde{g}}\psi_1 &= \mathcal{L}_{\psi_1},
\end{align*}
which is guaranteed by the theory of symmetric hyperbolic systems. Then define a sequence $(\phi_n,\psi_n)$ by the relations
\begin{align*}
\Box_g\phi_n -\mathcal{L}_{\phi_n} &= \mathcal{N}_{\phi_{n-1}} \\
\Box_{\tilde{g}}\psi_n -\mathcal{L}_{\psi_n} &= \mathcal{N}_{\psi_{n-1}},
\end{align*}
where each pair $(\phi_n,\psi_n)$ is constructed from the pair $(\phi_{n-1},\psi_{n-1})$ using Duhamel's method.
\end{definition}

By slightly adapting the proof of Theorem \ref{main_thm}, we obtain the following theorem for each of the pairs $(\phi_n,\psi_n)$.
\begin{theorem}\label{lwp:start_thm}
Define the energies
$$E^{\ul{m}}(t)=\mathring{E}^m[\phi_n,\psi_n](t)+\sum_{s+2k=m-1}E^{s,k}[\phi_n,\psi_n](t),$$
$$E_p^{\ul{m}}(t)=\mathring{E}_p^{m}[\phi_n,\psi_n](t)+\sum_{s+2k=m-1}E_p^{s,k}[\phi_n,\psi_n](t).$$
Then for $\delp,\delm>0$ sufficiently small, if the initial data $(\phi,\psi)|_{t=0}$ decay sufficiently fast as $r\rightarrow \infty$ and have size
\begin{equation*}
I_0=E^{\ul{29}}(0)+E^{\ul{29}}_{2-\delp}(0)
\end{equation*}
sufficiently small, then the following estimates hold for $t\in [0,T]$.

I) The energies satisfy
$$E^{\ul{29}}(t)\lesssim I_0$$
$$E^{\ul{29}}_{p\in[\delm,2-\delp]}(t)\lesssim I_0$$

II) For all $s,k$ such that $s+2k\le 28$, the following $L^\infty$ estimates hold.
\begin{multline*}
|r^{p+1}\bar{D}\fd^l\phi_n^{s-l,k}|^2+|r^pD\fd^l\phi_n^{s-l,k}|^2+|r^p\fd^l\phi_n^{s-l,k}|^2 \\
+|r^{p+3}\bar{D}\fd^l\psi_n^{s-l,k}|^2+|r^{p+2}D\fd^l\psi_n^{s-l,k}|^2+|r^{p+2}\fd^l\psi_n^{s-l,k}|^2 \\
\lesssim E_{2p}^{s+5,k+1}[\phi_n,\psi_n](t)+E_{2p}^{s+7,k}[\phi_n,\psi_n](t)
\end{multline*}

III) Together, (I) and (II) imply that if $s+2k\le 15$, for all $p\in [\delm/2,(2-\delp)/2]$,
\begin{multline*}
|r^{p+1}\bar{D}\fd^l\phi_n^{s-l,k}|+|r^pD\fd^l\phi_n^{s-l,k}|+|r^p\fd^l\phi_n^{s-l,k}| \\
+|r^{p+3}\bar{D}\fd^l\psi_n^{s-l,k}|+|r^{p+2}D\fd^l\psi_n^{s-l,k}|+|r^{p+2}\fd^l\psi_n^{s-l,k}| \\
\lesssim I_0^{1/2}.
\end{multline*}
\end{theorem}
\begin{proof}
The proof is very similar to the proof of Theorem \ref{main_thm}. It is proved by induction on $n$. For the base case $n=1$, since $(\phi_1,\psi_1)$ solves the linear system, the proof is substantially easier. For the inductive step, the inductive hypothesis replaces the bootstrap assumptions.
\end{proof}

\subsection{Establishing the contraction map}

Now, we wish to show that the map
that generates the sequence $(\phi_n,\psi_n)$ is a contraction map with respect to the following norm.

\begin{definition}
Fix $T>0$. Then define
$$|||(\phi',\psi')||| = \sup_{t\in [0,T]} \left(E^{\ul{29}}[\phi',\psi'](t) + E_{2-\delp}^{\ul{29}}[\phi',\psi'](t)\right),$$
where the norms $E^{\ul{m}}(t)$ and $E^{\ul{m}}_p(t)$ are as defined in Theorem \ref{lwp:start_thm}.
\end{definition}

We begin by investigating the differences. Let
\begin{align*}
\Delta\phi_n &= \phi_{n+1}-\phi_n \\
\Delta\psi_n &= \psi_{n+1}-\psi_n.
\end{align*}
Then by the principle of linearity,
\begin{align*}
\Box_g\Delta\phi_n-\mathcal{L}_{\Delta\phi_n} &= \mathcal{N}_{\phi_n}-\mathcal{N}_{\phi_{n-1}} \\
\Box_{\tilde{g}}\Delta\psi_n-\mathcal{L}_{\Delta\psi_n} &= \mathcal{N}_{\psi_n}-\mathcal{N}_{\psi_{n-1}}.
\end{align*}
In light of this fact, we have the following lemma.
\begin{lemma}\label{lwp:contraction_1_lem}
Let $p=2-\delp$. Then
\begin{align*}
(1-CT)&|||(\Delta\phi_n,\Delta\psi_n)||| \lesssim [[N_n-N_{n-1}]],
\end{align*}
where
\begin{align*}
[[N_n-&N_{n-1}]] \\
&= \sum_{s\le 29}\int_0^T||\pd_t^s(\mathcal{N}_{\phi_n}-\mathcal{N}_{\phi_{n-1}})||^2_{L^2(\Sigma_t)}dt 
+ \sum_{s\le 29}\int_0^T||\pd_t^s(\mathcal{N}_{\psi_n}-\mathcal{N}_{\psi_{n-1}})||^2_{L^2(\tilde\Sigma_t)}dt \\
&\hspace{.2in}+ \sum_{s+2k\le 28}\int_0^T||q^{-2}\tg^k\Gamma^s(q^2(\mathcal{N}_{\phi_n}-\mathcal{N}_{\phi_{n-1}}))||^2_{L^2(\Sigma_t)}dt \\
&\hspace{.2in}+ \sum_{s+2k\le 28}\int_0^T||q^{-2}\tg^k\Gamma^s(q^2(\mathcal{N}_{\psi_n}-\mathcal{N}_{\psi_{n-1}}))||^2_{L^2(\tilde\Sigma_t)}dt \\
&\hspace{.2in}+ \sum_{s\le 29}\int_0^T||r^{\frac{p+1}2}\pd_t^s(\mathcal{N}_{\phi_n}-\mathcal{N}_{\phi_{n-1}})||^2_{L^2(\Sigma_t)}dt \\
&\hspace{.2in}+ \sum_{s\le 29}\int_0^T||r^{\frac{p+1}2}\pd_t^s(\mathcal{N}_{\psi_n}-\mathcal{N}_{\psi_{n-1}})||^2_{L^2(\tilde\Sigma_t)}dt \\
&\hspace{.2in}+ \sum_{s+2k\le 28}\int_0^T||r^{\frac{p+1}2}q^{-2}\tg^k\Gamma^s(q^2(\mathcal{N}_{\phi_n}-\mathcal{N}_{\phi_{n-1}}))||^2_{L^2(\Sigma_t)}dt \\
&\hspace{.2in}+ \sum_{s+2k\le 28}\int_0^T||r^{\frac{p+1}2}q^{-2}\tg^k\Gamma^s(q^2(\mathcal{N}_{\psi_n}-\mathcal{N}_{\psi_{n-1}}))||^2_{L^2(\tilde\Sigma_t)}dt.
\end{align*}
\end{lemma}
\begin{proof}
From Theorems \ref{p_o_thm}, \ref{p_s_thm}, and \ref{p_s_k_thm}, we have
\begin{align*}
E^{\ul{29}}&[\Delta\phi_n,\Delta\psi_n](t)  \\
&\lesssim E^{\ul{29}}[\Delta\phi_n,\Delta\psi_n](0) +\int_0^t E^{\ul{29}}[\Delta\phi_n,\Delta\psi_n](\tau)d\tau \\
&\hspace{.2in}+ \sum_{s\le 29}\int_0^t||\pd_t^s(\mathcal{N}_{\phi_n}-\mathcal{N}_{\phi_{n-1}})||^2_{L^2(\Sigma_\tau)}d\tau 
+ \sum_{s\le 29}\int_0^t||\pd_t^s(\mathcal{N}_{\psi_n}-\mathcal{N}_{\psi_{n-1}})||^2_{L^2(\tilde\Sigma_\tau)}d\tau \\
&\hspace{.2in}+ \sum_{s+2k\le 28}\int_0^t||q^{-2}\tg^k\Gamma^s(q^2(\mathcal{N}_{\phi_n}-\mathcal{N}_{\phi_{n-1}}))||^2_{L^2(\Sigma_\tau)}d\tau \\
&\hspace{.2in}+ \sum_{s+2k\le 28}\int_0^t||q^{-2}\tg^k\Gamma^s(q^2(\mathcal{N}_{\psi_n}-\mathcal{N}_{\psi_{n-1}}))||^2_{L^2(\tilde\Sigma_\tau)}d\tau.
\end{align*}
Since $(\phi_{n+1},\psi_{n+1})$ and $(\phi_n,\psi_n)$ agree at $t=0$, it follows that
$$E^{\ul{29}}[\Delta\phi_n,\Delta\psi_n](0)=0.$$
Also,
$$\int_0^t E^{\ul{29}}[\Delta\phi_n,\Delta\psi_n](\tau)d\tau \le t \sup_{\tau\in[0,t]} E^{\ul{29}}[\Delta\phi_n,\Delta\psi_n](\tau) \le T |||(\Delta\phi_n,\Delta\psi_n)|||.$$
It follows that
\begin{align*}
\sup_{t\in[0,T]}&E^{\ul{29}}[\Delta\phi_n,\Delta\psi_n](t) \\
&\lesssim T|||(\Delta\phi_n,\Delta\psi_n)||| \\
&\hspace{.2in}+ \sum_{s\le 29}\int_0^T||\pd_t^s(\mathcal{N}_{\phi_n}-\mathcal{N}_{\phi_{n-1}})||^2_{L^2(\Sigma_t)}dt 
+ \sum_{s\le 29}\int_0^T||\pd_t^s(\mathcal{N}_{\psi_n}-\mathcal{N}_{\psi_{n-1}})||^2_{L^2(\tilde\Sigma_t)}dt \\
&\hspace{.2in}+ \sum_{s+2k\le 28}\int_0^T||q^{-2}\tg^k\Gamma^s(q^2(\mathcal{N}_{\phi_n}-\mathcal{N}_{\phi_{n-1}}))||^2_{L^2(\Sigma_t)}dt \\
&\hspace{.2in}+ \sum_{s+2k\le 28}\int_0^T||q^{-2}\tg^k\Gamma^s(q^2(\mathcal{N}_{\psi_n}-\mathcal{N}_{\psi_{n-1}}))||^2_{L^2(\tilde\Sigma_t)}dt.
\end{align*}
By the same reasoning, for all $p\in [\delm,2-\delp]$,
\begin{align*}
\sup_{t\in[0,T]}&E_p^{\ul{29}}[\Delta\phi_n,\Delta\psi_n](t) \\
&\lesssim T|||(\Delta\phi_n,\Delta\psi_n)||| \\
&\hspace{.2in}+ \sum_{s\le 29}\int_0^T||r^{\frac{p+1}2}\pd_t^s(\mathcal{N}_{\phi_n}-\mathcal{N}_{\phi_{n-1}})||^2_{L^2(\Sigma_t)}dt \\
&\hspace{.2in}+ \sum_{s\le 29}\int_0^T||r^{\frac{p+1}2}\pd_t^s(\mathcal{N}_{\psi_n}-\mathcal{N}_{\psi_{n-1}})||^2_{L^2(\tilde\Sigma_t)}dt \\
&\hspace{.2in}+ \sum_{s+2k\le 28}\int_0^T||r^{\frac{p+1}2}q^{-2}\tg^k\Gamma^s(q^2(\mathcal{N}_{\phi_n}-\mathcal{N}_{\phi_{n-1}}))||^2_{L^2(\Sigma_t)}dt \\
&\hspace{.2in}+ \sum_{s+2k\le 28}\int_0^T||r^{\frac{p+1}2}q^{-2}\tg^k\Gamma^s(q^2(\mathcal{N}_{\psi_n}-\mathcal{N}_{\psi_{n-1}}))||^2_{L^2(\tilde\Sigma_t)}dt.
\end{align*}
Adding both of these estimates and using the expression $[[N_n-N_{n-1}]]$ to represent all terms that contain a difference of nonlinear terms,
$$|||(\Delta\phi_n,\Delta\psi_n)||| \lesssim T|||(\Delta\phi_n,\Delta\psi_n)||| + [[N_n-N_{n-1}]]$$
whence
$$(1-CT)|||\Delta\phi_n,\Delta\psi_n)||| \lesssim [[N_n-N_{n-1}]]$$
\end{proof}

Now we prove the following lemma, which is essential to the contraction argument.
\begin{lemma}\label{lwp:contraction_2_lem}
If $I_0$ is sufficiently small, then
$$[[N_n-N_{n-1}]] \lesssim T |||(\Delta\phi_{n-1},\Delta\psi_{n-1})|||.$$
\end{lemma}
\begin{proof}
Recall that the general structures of $q^{-2}\tg^k\Gamma^s(q^2\mathcal{N}_{\phi})$ and $q^{-2}\tg^k\Gamma^s(q^2\mathcal{N}_{\psi})$ are given in Lemmas \ref{Nphi_s_k_structure_lem} and \ref{Npsi_s_k_structure_lem}. Even with these structures, to establish the estimate of this lemma in full detail would be quite laborious. We will illustrate the full proof by only showing that
\begin{equation}\label{lwp:difference_estimate_1}
||q^{-2}\tg^k\Gamma^s(q^2(\mathcal{N}_{\psi_n}-\mathcal{N}_{\psi_{n-1}}))||_{L^2(\tilde{\Sigma_t})} \lesssim E^{\ul{29}}[\Delta\phi_{n-1},\Delta\psi_{n-1}](t).
\end{equation}
We choose this example estimate because it contains the most general term that has a singular factor on the axis. From Lemma \ref{Npsi_s_k_structure_lem}, such a term can be expressed as
$$\frac{r^{-2}f\gamma^{s_1,k_1}(r^{-1}\fb\fc^l\phi^{s_2-2l,k_2})(r\beta^{s_3,k_3})...(r\beta^{s_{2+j},k_{2+j}})}{(1+\phi)^{j+1}},$$
where $f$ is a bounded function and $\gamma^{s,k}$ and $\beta^{s,k}$ are factors, each depending on derivatives of $\phi$ and $\psi$. (See Definition \ref{wk:gamma_beta_s_k_def}.)

In yet another attempt to simplify the proof, we choose a simple representative term with the above structure. The source for this term (before applying commutators) is
$$\frac{r^{-2}(r\psi)(r^{-1}\fb\phi)}{1+\phi}=\frac{r^{-2}\psi\fb\phi}{1+\phi}.$$
(See the proof of Lemma \ref{wk:structure_Npsi_0_lem} for the origin of this term.)

Now, it is necessary to estimate
$$\frac{r^{-2}\psi_n\fb\phi_n}{1+\phi_n}-\frac{r^{-2}\psi_{n-1}\fb\phi_{n-1}}{1+\phi_{n-1}}$$
in terms of the difference
$$(\Delta\phi_{n-1},\Delta\psi_{n-1})=(\phi_n-\phi_{n-1},\psi_n-\psi_{n-1}).$$

\begin{lemma}
\begin{multline*}
\frac{r^{-2}\psi_n\fb\phi_n}{1+\phi_n}-\frac{r^{-2}\psi_{n-1}\fb\phi_{n-1}}{1+\phi_{n-1}} \\
= \frac{r^{-2}\Delta\psi_{n-1} \fb \phi_n}{1+\phi_n}+\frac{r^{-2}\psi_{n-1}\fb\Delta\phi_{n-1}}{1+\phi_n}-\frac{r^{-2}\Delta\phi_{n-1}\psi_{n-1}\fb\phi_{n-1}}{(1+\phi_n)(1+\phi_{n-1})}.
\end{multline*}
\end{lemma}
\begin{proof}
Note the following calculations.
$$\frac{r^{-2}\psi_n\fb\phi_n}{1+\phi_n} = \frac{r^{-2}(\psi_n-\psi_{n-1})\fb\phi_n}{1+\phi_n}+\frac{r^{-2}\psi_{n-1}\fb\phi_n}{1+\phi_n}$$
$$\frac{r^{-2}\psi_{n-1}\fb\phi_n}{1+\phi_n} = \frac{r^{-2}\psi_{n-1}\fb(\phi_n-\phi_{n-1})}{1+\phi_n} +\frac{r^{-2}\psi_{n-1}\fb\phi_{n-1}}{1+\phi_n}$$
$$\frac{r^{-2}\psi_{n-1}\fb\phi_{n-1}}{1+\phi_n} = r^{-2}\psi_{n-1}\fb\phi_{n-1}\left(\frac1{1+\phi_n}-\frac1{1+\phi_{n-1}}\right)+\frac{r^{-2}\psi_{n-1}\fb\phi_{n-1}}{1+\phi_{n-1}}$$
$$\frac1{1+\phi_n}-\frac1{1+\phi_{n-1}} = \frac{1+\phi_{n-1}-(1+\phi_n)}{(1+\phi_n)(1+\phi_{n-1})} = -\frac{\phi_n-\phi_{n-1}}{(1+\phi_n)(1+\phi_{n-1})}$$
The identity of the lemma now follows by combining all of these calculations.
\end{proof}

We illustrate the proof of (\ref{lwp:difference_estimate_1}) through the following two estimates, both of which rely on the $L^\infty$ estimates in Theorem \ref{lwp:start_thm}.
\begin{align*}
\left|\left| \frac{r^{-2}\psi_{n-1}^{lo}\fb\Delta\phi_{n-1}^{hi}}{1+\phi_n} \right|\right|_{L^2(\tilde{\Sigma}_t)}^2 &= \int_{r_H}^\infty\int_0^\pi \left(\frac{r^{-2}\psi_{n-1}^{lo}\fb\Delta\phi_{n-1}^{hi}}{1+\phi_n}\right)^2r^6\sin^5\theta d\theta dr \\
&\lesssim \int_{r_H}^\infty\int_0^\pi (r^{-1}\pd_\theta\Delta\phi_{n-1}^{hi})^2r^2\sin^3\theta d\theta dr \left|\left|\frac{r\psi^{lo}_{n-1}}{1+\phi_n}\right|\right|_{L^\infty(\tilde\Sigma_t)}^2 \\
&\lesssim E^{\ul{29}}[\Delta\phi_{n-1},\Delta\psi_{n-1}](t) I_0 \\
&\lesssim E^{\ul{29}}[\Delta\phi_{n-1},\Delta\psi_{n-1}](t)
\end{align*}
\begin{align*}
\left|\left| \frac{r^{-2}\psi_{n-1}^{hi}\fb\Delta\phi_{n-1}^{lo}}{1+\phi_n} \right|\right|_{L^2(\tilde{\Sigma}_t)}^2 &= \int_{r_H}^\infty\int_0^\pi \left(\frac{r^{-2}\psi_{n-1}^{hi}\fb\Delta\phi_{n-1}^{lo}}{1+\phi_n}\right)^2r^6\sin^5\theta d\theta dr \\
&\lesssim \int_{r_H}^\infty\int_0^\pi (r^{-1}\psi_{n-1}^{hi})^2r^6\sin^5\theta d\theta dr \left|\left|\frac{r^{-1}\fb\Delta\psi^{lo}_{n-1}}{1+\phi_n}\right|\right|_{L^\infty(\tilde\Sigma_t)}^2 \\
&\lesssim I_0 E^{\ul{29}}[\Delta\phi_{n-1},\Delta\psi_{n-1}](t)  \\
&\lesssim E^{\ul{29}}[\Delta\phi_{n-1},\Delta\psi_{n-1}](t)
\end{align*}
\end{proof}

\begin{proposition}\label{lwp:contraction_prop}
If $T$ and $I_0$ are both sufficiently small, then
$$|||(\phi_{n+1}-\phi_n,\psi_{n+1}-\psi_n)||| \le \lambda |||(\phi_n-\phi_{n-1},\psi_n-\psi_{n-1})|||,$$
where $\lambda<1$.
\end{proposition}
\begin{proof}
By Lemmas \ref{lwp:contraction_1_lem} and \ref{lwp:contraction_2_lem},
$$(1-C_1T)|||(\phi_{n+1}-\phi_n,\psi_{n+1}-\psi_n)||| \le C_2 T |||(\phi_n-\phi_{n-1},\psi_n-\psi_{n-1})|||.$$
Therefore, the estimate of the proposition holds with
$$\lambda = \frac{C_2T}{1-C_1T}.$$
If $T$ is sufficiently small, then $\lambda < 1$.
\end{proof}

\subsection{The limit $(\phi,\psi)=\lim_{n\rightarrow\infty}(\phi_n,\psi_n)$}

We have the following corollary of Proposition \ref{lwp:contraction_prop}.
\begin{corollary}\label{lwp:limit_cor}
There exists a pair $(\phi,\psi)$ such that
$$\lim_{n\rightarrow\infty}|||(\phi-\phi_n,\psi-\psi_n)|||=0.$$
\end{corollary}

Therefore, we conclude the following theorem, which establishes the local theory.
\begin{theorem}
Given initial data $(\phi,\psi)|_{t=0}$ with size
$$I_0=E^{\ul{29}}(0)+E_{2-\delp}^{\ul{29}}(0)$$
sufficiently small, there exists a small time $T$ such that a solution $(\phi,\psi)$ to the system (\ref{phi_nonlinear_system_eqn}-\ref{psi_nonlinear_system_eqn}) exists for all $t\in[0,T]$.
\end{theorem}
\begin{proof}
Let $(\phi,\psi)$ be as defined in Corollary \ref{lwp:limit_cor}. Then by the standard theory of Sobolev spaces,
\begin{align*}
\Box_g\phi-\mathcal{L}_\phi &= \lim_{n\rightarrow\infty} \Box_g\phi_n-\mathcal{L}_{\phi_n} = \lim_{n\rightarrow\infty} \mathcal{N}_{\phi_{n-1}} = \mathcal{N}_\phi \\
\Box_{\tilde{g}}\psi-\mathcal{L}_\psi &= \lim_{n\rightarrow\infty} \Box_{\tilde{g}}\psi_n-\mathcal{L}_{\psi_n} = \lim_{n\rightarrow\infty} \mathcal{N}_{\psi_{n-1}} = \mathcal{N}_\psi.
\end{align*}
\end{proof}

\chapter{General ellipic regularity}\label{elliptic_sec}

In this appendix, we prove Lemma \ref{main_elliptic_lemma}, which shows that $Q$ is a useful commutator for controlling angular derivatives.

\section{Elliptic estimates for $\sla\triangle$ and $\sla\triangle^k$}

Before investigating $Q^k$, we derive elliptic estimates for $\sla\triangle^k$. We begin with a well-known calculation for $\sla\triangle$.

\begin{lemma}\label{standard_elliptic_base_lem}
$$\int_{S^2(1)}|\sla\nabla^2u|^2 + \int_{S^2(1)}|\sla\nabla u|^2 = \int_{S^2(1)}(\sla\triangle u)^2.$$
\end{lemma}
\begin{proof}
We calculate
\begin{align*}
\int_{S^2(1)} |\sla\nabla^{2}u|^2&=\int_{S^2(1)} \sla\nabla^{\alpha_1}\sla\nabla^{\alpha_2}u \sla\nabla_{\alpha_1}\sla\nabla_{\alpha_2} u \\
&=-\int_{S^2(1)} \sla\nabla^{\alpha_2} u \sla\triangle\sla\nabla_{\alpha_2} u \\
&=-\int_{S^2(1)} \sla\nabla^{\alpha_2} u \sla\nabla_{\alpha_2}\sla\triangle u-\int_{S^2(1)} \sla\nabla^{\alpha_2} u [\sla\triangle,\sla\nabla]_{\alpha_2}u \\
&= \int_{S^2(1)} (\sla\triangle u)^2-\int_{S^2(1)} Ric(\sla\nabla u,\sla\nabla u).
\end{align*}
Since on $S^2(1)$, 
$$Ric(\sla\nabla u,\sla\nabla u)=|\sla\nabla u|^2,$$
it follows that
$$\int_{S^2(1)} |\sla\nabla^2u|^2+\int_{S^2(1)} |\sla\nabla u|^2 = \int_{S^2(1)} (\sla\triangle u)^2.$$
\end{proof}

Now we generalize the previous lemma to $\sla\triangle^k$.

\begin{lemma}\label{standard_elliptic_lemma}
$$\int_{S^2(1)}|\sla\nabla^{2k}u|^2 + \int_{S^2(1)}|\sla\nabla^{2k-1}u|^2\lesssim \sum_{k'\le k}\int_{S^2(1)}(\sla\triangle^{k'}u)^2.$$
\end{lemma}

\begin{proof}
This in fact holds for any closed manifold with positive definite Ricci curvature and sufficiently regular Riemann curvature.

Throughout the proof, all integrals are taken over $S^2(1)$ and the tensor indices are dropped as they can be determined by contracting over indices in the order that they appear in each term. For example, 
$$\sla\nabla^3 u\cdot\sla\nabla^3 u=\sla\nabla^{\alpha_1}\sla\nabla^{\alpha_2}\sla\nabla^{\alpha_3}u\sla\nabla_{\alpha_1}\sla\nabla_{\alpha_2}\sla\nabla_{\alpha_3}u.$$

We proceed by induction on $k$. The case $k=1$ was established by Lemma \ref{standard_elliptic_base_lem}. We assume that the statement of the lemma holds for $k$ and prove that it also holds for $k+1$. This is accomplished using a calculation similar to (but significantly more general than) the proof of Lemma \ref{standard_elliptic_base_lem}.
\begin{align*}
\int |\sla\nabla^{2(k+1)}u|^2 &= \int \sla\nabla^{2k+2}u\cdot \sla\nabla^{2k+2}u \\
&= -\int \sla\nabla^{2k+1}u \cdot \sla\triangle \sla\nabla^{2k+1}u \\
&= -\int \sla\nabla^{2k+1}u \cdot \sla\nabla \sla\triangle \sla\nabla^{2k}u -\int \sla\nabla^{2k+1}u\cdot [\sla\triangle,\sla\nabla]\sla\nabla^{2k}u \\
&= \int |\sla\triangle \sla\nabla^{2k}u|^2 -\int \sla\nabla^{2k+1}u\cdot [\sla\triangle,\sla\nabla]\sla\nabla^{2k}u.
\end{align*}
The commutator $[\sla\triangle,\sla\nabla]$ will soon be computed explicitly, but at this point it is important to observe that it is on the order of one derivative. Therefore, terms with a commutator have two fewer derivatives than terms without a commutator. In particular, terms with two commutators can be controlled by $\int |\sla\nabla^{2k}u|^2$ which, by the inductive assumption, can be controlled by $\int |\sla\triangle^{\le k}u|^2$. Thus, we ignore these terms and simply write them as $err$. This means that as soon as a single commutator appears in a term, we are free to commute any of the remaining derivatives without bothering to calculate the error terms explicitly.

Our next task is to commute the Laplacians toward $u$. We calculate
\begin{align*}
\int |&\sla\nabla^j\sla\triangle \sla\nabla^{2k-j}u|^2 \\
&= \int \sla\nabla^j\sla\triangle \sla\nabla^{2k-j}u \cdot \sla\nabla^{j+1}\sla\triangle\sla\nabla^{2k-j-1}u + \int \sla\nabla^j\sla\triangle \sla\nabla^{2k-j}u \cdot \sla\nabla^{j}[\sla\triangle,\sla\nabla]\sla\nabla^{2k-j-1}u \\
&= \int \sla\nabla^j\sla\triangle \sla\nabla^{2k-j}u\cdot  \sla\nabla^{j+1}\sla\triangle\sla\nabla^{2k-j-1}u + \int \sla\triangle \sla\nabla^{2k}u\cdot  \sla\nabla^{j}[\sla\triangle,\sla\nabla]\sla\nabla^{2k-j-1}u +err \\
&= \int \sla\nabla^j\sla\triangle \sla\nabla^{2k-j}u\cdot \sla\nabla^{j+1}\sla\triangle\sla\nabla^{2k-j-1}u - \int \sla\nabla^{2k+1}u \cdot \sla\nabla^{j+1}[\sla\triangle,\sla\nabla]\sla\nabla^{2k-j-1}u +err .
\end{align*}
Repeating the above procedure again, we find that
\begin{equation*}
\int |\sla\nabla^j\sla\triangle \sla\nabla^{2k-j}u|^2 = \int |\sla\nabla^{j+1}\sla\triangle \sla\nabla^{2k-j-1}u|^2 - 2\int \sla\nabla^{2k+1}u\cdot \sla\nabla^{j+1}[\sla\triangle,\sla\nabla]\sla\nabla^{2k-j-1}u +err.
\end{equation*}
Therefore,
\begin{equation*}
\int |\sla\triangle \sla\nabla^{2k}u|^2 = \int |\sla\nabla^{2k}\sla\triangle u|^2 - 2\sum_{j=0}^{2k-1}\int \sla\nabla^{2k+1}u\cdot \sla\nabla^{j+1}[\sla\triangle,\sla\nabla]\sla\nabla^{2k-j-1}u +err.
\end{equation*}
Thus, all of our highest order error terms are of the form
$$\int \sla\nabla^{2k+1}u\cdot \sla\nabla^j[\sla\triangle,\sla\nabla]\sla\nabla^{2k-j}u.$$
We now compute the commutator explicitly.
\begin{align*}
\sla\triangle \sla\nabla_\alpha \sla\nabla_{\beta_1}...\sla\nabla_{\beta_l}u &= g^{\mu\nu} \sla\nabla_\mu\sla\nabla_\nu\sla\nabla_\alpha \sla\nabla_{\beta_1}...\sla\nabla_{\beta_l}u \\
&= g^{\mu\nu} \sla\nabla_\mu\sla\nabla_\alpha \sla\nabla_\nu \sla\nabla_{\beta_1}...\sla\nabla_{\beta_l}u - g^{\mu\nu}\sla\nabla_\mu\sum_i R^\lambda{}_{\beta_i\nu\alpha} \sla\nabla_{\beta_1}...\sla\nabla_{\lambda}...\sla\nabla_{\beta_l}u \\
&= \sla\nabla_\alpha \sla\triangle \sla\nabla_{\beta_1}...\sla\nabla_{\beta_l}u - g^{\mu\nu}\sla\nabla_\mu\sum_i R^\lambda{}_{\beta_i\nu\alpha} \sla\nabla_{\beta_1}...\sla\nabla_{\lambda}...\sla\nabla_{\beta_l}u \\
&\hphantom{==} - g^{\mu\nu}R^\lambda{}_{\nu\mu\alpha}\sla\nabla_{\lambda}\sla\nabla_{\beta_1}...\sla\nabla_{\beta_l}u - g^{\mu\nu}\sum_iR^\lambda{}_{\beta_i\mu\alpha}\sla\nabla_{\nu}\sla\nabla_{\beta_1}...\sla\nabla_{\lambda}...\sla\nabla_{\beta_l}u
\end{align*}
So schematically,\footnote{In the special case of the sphere, $\sla\nabla Riem=0$.}
\begin{equation*}
[\sla\triangle,\sla\nabla] \sla\nabla^l u = Ric\cdot \sla\nabla^{l+1}u +\sla\nabla Riem \sla\nabla^l u+Riem\sla\nabla^{l+1}u.
\end{equation*}
By antisymmetry in certain pairs of indices of $Riem$ and symmetry (up to terms with fewer derivatives of $u$) in indices of $\nabla^{2k+1}u$, after integrating by parts twice, we find that only the term with the Ricci curvature stands out.
\begin{equation*}
\int \sla\nabla^{2k+1}u\cdot\sla\nabla^j[\sla\triangle,\sla\nabla]\sla\nabla^{2k-j}u = \int Ric^{\mu\nu}\sla\nabla_\mu \sla\nabla^{2k}u\cdot\sla\nabla_\nu\sla\nabla^{2k}u + err.
\end{equation*}
Noting the sign of the error terms involving $Ric$, and using the fact that $Ric$ is positive definite, we can finally complete the inductive step.
\begin{equation*}
\int |\sla\nabla^{2(k+1)}u|^2+\int |\sla\nabla^{2(k+1)-1}u|^2 \lesssim \int |\sla\nabla^{2(k+1)}u|^2+\int Ric^{\mu\nu}\sla\nabla_\mu\sla\nabla^{2k}u\cdot\sla\nabla_\nu\sla\nabla^{2k}u,
\end{equation*}
and
\begin{align*}
\int |\sla\nabla^{2(k+1)}u|^2+&(4k+1)\int Ric^{\mu\nu}\sla\nabla_\mu\sla\nabla^{2k}u\cdot\sla\nabla_\nu\sla\nabla^{2k}u \\
&= \int |\sla\nabla^{2k}\sla\triangle u|^2 + err \\
&\lesssim \sum_{k'\le k}\int (\sla\triangle^{k'}\sla\triangle u)^2+\sum_{k'\le k}\int (\sla\triangle^{k'}u)^2 \\
&\lesssim \sum_{k'\le k+1}\int (\sla\triangle^{k'} u)^2.
\end{align*}
\end{proof}

\section{Elliptic estimates for $Q^k$}

Now we prove that $Q$ behaves similarly to $\sla\triangle$. Recall that
$$Q = \sla\triangle + a^2\sin^2\theta\pd_t^2.$$

\begin{lemma}\label{main_elliptic_lemma}
Let $\Gamma^s$ represent $Q^{l}\pd_t^{s-2l}$ for any $l\le s/2$. Then
$$\int_{S^2(1)}|\sla\nabla^{2k}u|^2+\int_{S^2(1)}|\sla\nabla^{2k-1}u|^2\lesssim \sum_{s\le 2k}\int_{S^2(1)}(\Gamma^{s}u)^2.$$
\end{lemma}

\begin{proof}

Given Lemma \ref{standard_elliptic_lemma}, it suffices to show that
$$\int_{S^2(1)}(\sla\triangle^{k} u)^2 \lesssim \sum_{s\le 2k}\int_{S^2(1)} (\Gamma^{s}u)^2.$$
We shall prove by induction (on $k$) that for all $s\ge 0$,
$$\int_{S^2(1)}(\sla\triangle^{k}\Gamma^s u)^2 \lesssim \sum_{s'\le 2k+s}\int_{S^2(1)} (\Gamma^{s'}u)^2.$$
The base case ($k=0$) is trivial. Let us assume the above estimate and prove the corresponding estimate for $k+1$. Assume for the moment that
\begin{equation}\label{trivial_eqn}
\int_{S^2(1)} \left(\sla\triangle^k(a^2\sin^2\theta\pd_t^2 \Gamma^s u)\right)^2\lesssim \sum_{k'\le k}\int_{S^2(1)} (\sla\triangle^{k'} \Gamma^{s+2}u)^2.
\end{equation}
Then we can perform the inductive step.
\begin{align*}
\int_{S^2(1)}(\sla\triangle^{k+1}\Gamma^su)^2&=\int_{S^2(1)}\left(\sla\triangle^k(Q-a^2\sin^2\theta\pd_t^2)\Gamma^s u\right)^2 \\
&\lesssim \int_{S^2(1)} (\sla\triangle^k Q \Gamma^su)^2+\int_{S^2(1)}\left(\sla\triangle^k(a^2\sin^2\theta\pd_t^2 \Gamma^su)\right)^2 \\
&\lesssim \sum_{k'\le k}\int_{S^2(1)} (\sla\triangle^{k'}\Gamma^{s+2}u)^2 \\
&\lesssim \sum_{s'\le 2k+(s+2)}\int_{S^2(1)} (\Gamma^{s'}u)^2 = \sum_{s'\le 2(k+1)+s}\int_{S^2(1)} (\Gamma^{s'}u)^2.
\end{align*}
Finally, we verify (\ref{trivial_eqn}). By a Sobolev estimate and then lemma \ref{standard_elliptic_lemma},
\begin{multline*}
\int_{S^2(1)} \left(\sla\triangle^k(a^2\sin^2\theta\pd_t^2 \Gamma^su)\right)^2 \\
\lesssim ||\sin^2\theta||^2_{H^{2k}(S^2)}||\Gamma^{s+2}u||^2_{H^{2k}(S^2)} \\
\lesssim \left(\sum_{k'\le k}\int_{S^2(1)} (\sla\triangle^{k'}\sin^2\theta)^2\right)\left( \sum_{k'\le k}\int_{S^2(1)} (\sla\triangle^{k'}\Gamma^{s+2}u)^2\right).
\end{multline*}
The regularity of $\sin^2\theta$ can be verified by computing
$$\sla\triangle \sin^2\theta=4\cos^2\theta-2\sin^2\theta,$$
$$\sla\triangle \cos^2\theta=-4\cos^2\theta+2\sin^2\theta.$$
This concludes the proof of Lemma \ref{main_elliptic_lemma}.
\end{proof}

\chapter{Regularity for axisymmetric functions}\label{regularity_sec}

In this appendix, we develop a formalism for understanding terms that represent regular axisymmetric functions using certain differential operators on the unit sphere. The main motivation for this chapter in the context of nonlinear wave equations is summarized in \S\ref{intro_regularity_sec}. The formalism developed here can be used for a variety of geometric problems with a coordinate degeneracy such as that in axisymmetry.

To begin, in \S\ref{gothic_operators_intro_sec}, we define the gothic operators $\fa$, $\fb$, $\fc^l$ and $\fd^l$ and prove a few basic properties. Of particular importance is Lemma \ref{cl_raise_degree_lem}, which explains how these operators act on a particular class of products of functions. In \S\ref{cl_embedding_theorem_sec}, we prove an important embedding theorem that allows us to estimate $L^2$ norms with these operators by $L^2$ norms with the standard spherical laplacian. This is important, because the spherical laplacian is a useful commutator in many geometric problems. Finally, in \S\ref{additional_regularity_lemmas_sec} we prove a few additional related estimates that are needed in Chapters \ref{wm_szd_chap} and \ref{wm_kerr_chap}.

For the remainder of this appendix, we focus on functions defined on a sphere $S^2$ (or $S^6$) which depend only on the angle $\theta$ from the north pole. However, the theory discussed here easily extends to the Kerr spacetimes and other geometries with a coordinate singularity.

\section{The operators $\fa$, $\fb$, $\fc^l$, and $\fd^l$}\label{gothic_operators_intro_sec}

We begin by defining the gothic operators $fa$, $\fb$, $\fc^l$, and $\fd^l$. The most fundamental of these are the operators $\fa$ and $\fb$, which are defined first.
\begin{definition}The operators $\fa$ and $\fb$ are
$$\fa=\pd_\theta^2,$$
$$\fb=\cot\theta\pd_\theta.$$
\end{definition}

In many cases, we will use these operators interchangeably, so we define the following shorthand notation.
\begin{definition}The operator $\fc$ is
$$\fc=\fa\text{ or }\fb.$$
More generally, the family of operators $\fc^l$ is
$$\fc^l=\fc_1...\fc_l\text{ where each }\fc_i\text{ is either }\fa\text{ or }\fb.$$
So, for example, $\fc^2$ represents any of the operators $\fa^2$, $\fa\fb$, $\fb\fa$, or $\fb^2$.
\end{definition}

Finally, we generalize the operator family $\fc^l$ slightly.
\begin{definition}The operator $\fd$ is
$$\fd=\pd_\theta.$$
More generally, the family of operators $\fd^l$ is
$$
\fd^l = \left\{\begin{array}{cc}
\fc^{l/2} & l\in 2\mathbb{Z} \\
\pd_\theta \fc^{(l-1)/2} & l\not\in 2\mathbb{Z}.
\end{array}\right.
$$
So, for example, $\fd^4$ represents any of the operators $\fa^2$, $\fa\fb$, $\fb\fa$, or $\fb^2$, while $\fd^5$ represents any of the operators $\pd_\theta\fa^2$, $\pd_\theta\fa\fb$, $\pd_\theta\fb\fa$, or $\pd_\theta\fb^2$.
\end{definition}

\begin{remark}
One can think of the operator $\fd$ as a single spherical derivative and the operator family $\fd^l$ as $l$ spherical derivatives. This motivates the use of the letter d for this family. In contrast, the operator family $\fc^l$ is more like $2l$ spherical derivatives.
\end{remark}

To gain familiarity with the operator family $\fd^l$, we note the following example, which also serves as a caveat.
\begin{example}
Although
$$\fd\fd^l\subset\fd^{l+1},$$
in general
$$\fd^{l_1}\fd^{l_2}\not\subset\fd^{l_1+l_1}.$$
\end{example}
\begin{proof}
If $l$ is even, then 
$$\fd\fd^l=\pd_\theta\fc^{l/2}=\fd^{l+1}.$$
If $l$ is odd, then 
$$\fd\fd^l=\pd_\theta\pd_\theta\fc^{(l-1)/2}=\fa\fc^{(l-1)/2}\subset\fc^{(l+1)/2}=\fd^{l+1}.$$
A counterexample for the more general case of $\fd^{l_1}\fd^{l_2}$ exists when $l_1=2$ and $l_2=1$. In this case,
$$\fb\pd_\theta\in \fd^2\fd,\text{ but }\fb\pd_\theta=\cot\theta\fa \not\in \fd^3.$$
\end{proof}

\subsection{The commutators $[\fb,\fa]$ and $[\fb,\fc^l]$}

The operators $\fa$ and $\fb$ satisfy a simple commutation relation.
\begin{lemma}
$$[\fb,\fa]=2\fb^2+2\fa.$$
\end{lemma}
\begin{proof}
Let $\alpha = \cot\theta$. Then $\fb=\alpha\pd_\theta$. By direct calculation,
\begin{align*}
\alpha' &= \pd_\theta\left(\frac{\cos\theta}{\sin\theta}\right) \\
&= -1-\frac{\cos^2\theta}{\sin^2\theta} \\
&= -\alpha^2-1
\end{align*}
and
\begin{align*}
\alpha'' &= (-\alpha^2-1)' \\
&= -2\alpha\alpha'.
\end{align*}
Now we compute the commutator.
\begin{align*}
[\fb,\fa](f) &= \fb\fa f -\fa\fb f \\
&= \alpha\pd_\theta^3 f - \pd_\theta^2(\alpha\pd_\theta f) \\
&= -2\alpha'\pd_\theta^2 f-\alpha''\pd_\theta f \\
&= 2\alpha^2\pd_\theta^2f + 2\pd_\theta^2 f+2\alpha\alpha'\pd_\theta f \\
&= 2\alpha\pd_\theta(\alpha\pd_\theta f)+2\pd_\theta^2 f \\
&= 2\fb^2f+2\fa f.
\end{align*}
This verifies the lemma.
\end{proof}

We generalize the previous commutation relation to obtain something that will be essential in Chapter \ref{wm_kerr_chap}. (See Lemma \ref{Npsi_s_k_structure_lem}.)
\begin{lemma}\label{fb_fcl_commutation_lem}
The operators $\fb$ and $\fc^l$ satisfy the commutation relation
$$\fc^l\fb \approx \fb\fc^l+\fc^l,$$
where the $\approx$ sign indicates that the identity holds modulo constant factors for the terms on the right side.
\end{lemma}
\begin{proof}
The proof is a simple induction exercise. The base case is handled by the previous lemma. Now, we assume the statement of the lemma at the level $l$ and prove it for the level $l+1$.
\begin{align*}
\fc^{l+1}\fb &\approx \fc(\fc^l\fb)\\
&\approx \fc (\fb \fc^l + \fc^l) \\
&\approx \fc\fb\fc^l+\fc^{l+1} \\
&\approx (\fb\fc +\fc)\fc^l+\fc^{l+1} \\
&\approx \fb\fc^{l+1}+\fc^{l+1}.
\end{align*}
\end{proof}

\subsection{The graded algebra $\bigoplus_{n\in\mathbb{N}}\tau_{(n)}$}

Now, we investigate products of axisymmetric functions, and in particular we want to understand what types of products are regular on the axis. Certainly, a product of regular axisymmetric functions will also be regular, but the converse is not true. There are regular products of functions that are not necessarily regular themselves. The simplest example is the product $(\sin\theta)^2$. This product is regular on the axis, but the function $\sin\theta$ is not regular, because it behaves like $|\theta|$ in a neighborhood of the half axis $\theta=0$.

We understand regularity for products of functions by understanding how the gothic operators act on these products.
\begin{definition}
Let $\mathcal{F}$ be a fixed set of functions. Define $\tau_{(n)}=\tau_{(n)}(\mathcal{F})$ to be a family of terms of the form
$$\fd^{i_1}f_1...\fd^{i_k}f_k$$
where $f_1,...,f_k\in\mathcal{F}$ and $i_1+...+i_k=2n$. We say a term $\tau\in \tau_{(n)}$ \textit{has degree} $n$.
\end{definition}

The following example should make the meaning of the above definition slightly more clear.
\begin{example}
If $f,g,h\in \mathcal{F}$, then $\fc^i f\fc^j g\in\tau_{(i+j)}$ and $\fc^i f \pd_\theta \fc^j g \pd_\theta \fc^k h\in\tau_{(i+j+k+1)}$, however $\fc^if \fc^j g\pd_\theta \fc^k h$ is not in $\tau_{(n)}$ for any $n$.
\end{example}


The following example should illustrate why these families of terms are useful.
\begin{example}\label{appendix_Pl_example}
Let $\mathcal{F}=\{\cos\theta\}$. Then any term $\tau\in\tau_{(n)}(\mathcal{F})$ is regular on the axis. In particular, $(\sin\theta)^2=(\pd_\theta\cos\theta)^2\in\tau_{(1)}$. In contrast, the term $\cos\theta \pd_\theta\cos\theta = -\cos\theta\sin\theta$ is not in $\tau_{(n)}$ for any $n$ and is also not regular at $\theta=0$ and $\theta=\pi$.
\end{example}

The main purpose for this classification of terms is the following fact.
\begin{lemma}\label{cl_raise_degree_lem}
If $\tau$ is a term of degree $n$, then $\fc\tau$ can be expressed as a sum of terms of degree $n+1$.
\end{lemma}
\begin{proof}
By definition, a term of degree $n$ is of the form
$$\fd^{i_1}f_1...\fd^{i_k}f_k,$$
where $i_1+...+i_k=2n$. We now compute
\begin{align*}
\fa(\fd^{i_1}f_1...\fd^{i_k}f_k) &= \pd_\theta^2(\fd^{i_1}f_1...\fd^{i_k}f_k) \\
&= \sum_{\substack{1\le j_1\le k \\ 1\le j_2\le k}} \fd^{i_1}f_1...\fd\fd^{i_{j_1}}f_{j_1}...\fd\fd^{i_{j_2}}f_{j_2}...\fd^{i_k}f_k \\
&= \sum_{\substack{1\le j_1\le k \\ 1\le j_2\le k}} \fd^{i_1}f_1...\fd^{i_{j_1}+1}f_{j_1}...\fd\fd^{i_{j_2}}f_{j_2}...\fd^{i_k}f_k \\
&= \sum_{\substack{1\le j_1\le k \\ 1\le j_2\le k}} \fd^{i_1}f_1...\fd^{i_{j_1}+1}f_{j_1}...\fd^{i_{j_2}+1}f_{j_2}...\fd^{i_k}f_k .
\end{align*}
The last two steps were kept separate to emphasize that the procedure works even for the terms where $j_1=j_2$. It should be clear that the final sum is a sum of terms of degree $n+1$.

We also compute
\begin{align*}
\fb(\fd^{i_1}f_1...\fd^{i_k}f_k) &= \cot\theta \pd_\theta(\fd^{i_1}f_1...\fd^{i_k}f_k) \\
&= \cot\theta \sum_{1\le j_1\le k} \fd^{i_1}f_1...\pd_\theta\fd^{i_{j_1}}f_{j_1}...\fd^{i_k}f_k \\
&= \cot\theta \sum_{1\le j_1\le k} \fd^{i_1}f_1...\fd^{i_{j_1}+1}f_{j_1}...\fd^{i_k}f_k.
\end{align*}
Now, observe that since $i_1+...+i_k=2n$, then $i_1+...+(i_{j_1}+1)+...+i_k=2n+1$. It follows that either $i_{j_1}+1$ is odd or there is some $j_2$ for which $i_{j_2}$ is odd. If $i_{j_1}+1$ is odd, then 
$$\cot\theta\fd^{i_{j_1}+1}f_{j_1}=\cot\theta\pd_\theta\fc^{i_{j_1}/2}f_{j_1}=\fb\fc^{i_{j_1}/2}f_{j_1}=\fc^{i_{j_1}/2+1}=\fd^{i_{j_1}+2}f_{j_1}.$$
Otherwise, for some $j_2$, $i_{j_2}$ is odd, so
$$\cot\theta\fd^{i_{j_2}}f_{j_2}=\cot\theta\pd_\theta\fc^{(i_{j_2}-1)/2}f_{j_2}=\fb\fc^{(i_{j_2}-1)/2}f_{j_2}=\fc^{(i_{j_2}-1)/2+1}f_{j_2}=\fd^{i_{j_2}+1}f_{j_2}.$$
So all terms can either be expressed as
$$\fd^{i_1}f_1...\fd^{i_{j_1}+2}f_{j_1}...\fd^{i_k}j_k$$
or
$$\fd^{i_1}f_1...\fd^{i_{j_1}+1}f_{j_1}...\fd^{i_{j_2}+1}f_{j_2}...\fd^{i_k}j_k.$$
Both of these are terms of degree $n+1$.
\end{proof}

\subsection{The spaces $P_\theta^l(S^2)$}

\begin{definition}
We define function spaces $P_\theta^l(S^2)$ by the following norm
$$|| f||_{P_\theta^l(S^2)} = \max_{l'\le l}\frac{1}{l'!}||\fd^{l'} f||_{L^2(S^2)}.$$
These spaces make sense for all integers $l$ and also for $l=\infty$.
\end{definition}

We illustrate this space with a few examples.
\begin{example}
The function $\cos\theta$ is in $P_\theta^\infty(S^2)$. This can be verified directly by the fact that
\begin{align*}
\fa \cos\theta &= -\cos\theta \\
\fb \cos\theta &= -\cos\theta.
\end{align*}
\end{example}

\begin{example}
The function $\sin\theta$ is not in $P^2_\theta(S^2)$. Although
$$\fa\sin\theta = -\sin\theta,$$
the problem is that
$$\fb\sin\theta = \frac{\cos^2\theta}{\sin\theta}.$$
\end{example}

\begin{example}
For any positive integer $k$, the function $(\sin\theta)^{2k}$ is in $P_\theta^\infty(S^2)$. However, the function $(\sin\theta)^{2k+1}$ is in $P^{2k}_\theta(S^2)$ but not in $P_\theta^{2k+2}(S^2)$. The latter fact can be verified by applying $\fb^{k+1}$.
\end{example}

The above examples illustrate how the $P_\theta^l(S^2)$ spaces are good candidate spaces for what has often been informally referred to as ``the space of regular axisymmetric functions'' in this thesis.

\section{An important embedding theorem}\label{cl_embedding_theorem_sec}
We now prove the following important embedding theorem.
\begin{theorem}\label{cl_embedding_thm}
$$||\fc^l f||_{L^2(S^2)}\lesssim \sum_{i\le l}||\sla\triangle^i f||_{L^2(S^2)},$$
$$||\pd_\theta \fc^l f||_{L^2(S^2)}\lesssim \sum_{i\le l}||\sla\nabla\sla\triangle^i f||_{L^2(S^2)}.$$
Thus,
$$H^l(S^2)\subset P_\theta^l(S^2).$$
\end{theorem}

\begin{proof}
We prove both estimates by induction on $l$. The base case ($l=0$) is trivial. We assume the estimates hold true at the level $l$, and now we must prove them at the level $l+1$.

First, observe that
$$||\fa^{l+1}f||_{L^2(S^2)}=||\pd_\theta^{2(l+1)}f||_{L^2(S^2)}\lesssim ||\sla\nabla^{2(l+1)}f||_{L^2(S^2)}\lesssim \sum_{i\le l+1} ||\sla\triangle^i \phi||_{L^2(S^2)}$$
and
$$||\pd_\theta \fa^{l+1}f||_{L^2(S^2)}=||\pd_\theta\pd_\theta^{2(l+1)}f||_{L^2(S^2)}\lesssim ||\sla\nabla^{2(l+1)+1}f||_{L^2(S^2)}\lesssim \sum_{i\le l+1} ||\sla\nabla \sla\triangle^i \phi||_{L^2(S^2)}.$$
(In both lines, the second step follows from the identity $\sla\nabla_{e_1}^kf=\pd_\theta^kf$ in the orthonormal frame $e_1=\pd_\theta$ and $e_2=\csc\theta\pd_\phi$, and the third step is a standard ellipticity result--see Appendix \ref{elliptic_sec}.)

We now assert the following lemma.
\begin{lemma}\label{ccc_lem}
$$\fc^k \approx \fa^k +\fc^{\le k-1}\sla\triangle + \fc^{\le k-1},$$
where the $\approx$ sign indicates that the identity holds modulo constant factors for the terms on the right side and the expression $\fc^{\le k-1}$ represents operators obtained by adding together constant multiples of operators of the form $\fc^i$ where $i\le k-1$.
\end{lemma}
From Lemma \ref{ccc_lem}, it is straightforward to prove the estimate at the level $l+1$. Set $k=l+1$. Then any term of the form $\fc^{l+1}f$ can be expressed as a sum of terms of the form $\fa^{l+1}f$ (which we just observed to be bounded by the appropriate norm) or $\fc^{\le l}f$ (which is bounded according to the inductive hypothesis) or $\fc^{\le l}\sla\triangle f$ (which is also bounded according to the inductive hypothesis applied to the function $g=\sla\triangle f$).

It remains to prove Lemma \ref{ccc_lem}.
\begin{proof} (of Lemma \ref{ccc_lem}) To prove Lemma \ref{ccc_lem}, we first assert the following lemma.
\begin{lemma}\label{baa_lem}
For all $k\ge 1$,
$$(2k-1)\fb \fa^{k-1} = -\fa^k +\fc^{\le k-1}\sla\triangle + \fc^{\le k-1}.$$
\end{lemma}

To see how Lemma \ref{ccc_lem} follows from Lemma \ref{baa_lem}, we must represent an arbitrary operator of the form $\fc^k=\fc_1...\fc_k$ in the form given by the right side of the identity in Lemma \ref{ccc_lem}. First, if $\fc_1=...=\fc_k=\fa$, then trivially $\fc_1...\fc_k=\fa^k$ is in the form given in Lemma \ref{ccc_lem}. Alternatively, there is some $j$ for which $\fc_1...\fc_k=\fc_1...\fc_j\fb\fa^{k-j-1}$. According to Lemma \ref{baa_lem}, we can rewrite
\begin{align*}
\fc_1...\fc_j\fb\fa^{k-j-1} &\approx \fc_1...\fc_j(\fa^{k-j}+\fc^{\le k-j-1}\sla\triangle+\fc^{\le k-j-1}) \\
&\approx \fc_1...\fc_j\fa^{k-j} +\fc^{\le k-1}\sla\triangle +\fc^{\le k-1}.
\end{align*}
The latter two terms are in the form given by Lemma \ref{ccc_lem} and the first term $\fc_1...\fc_j\fa^{k-j}$ has \textit{strictly more} $\fa$s at the right end. Thus, by repeating this procedure, we will end up with the term $\fa^k$, which is of the form given in Lemma \ref{ccc_lem}.

It remains to prove Lemma \ref{baa_lem}.

\begin{proof} (of Lemma \ref{baa_lem})
We prove Lemma \ref{baa_lem} by induction on $k$. Suppose $k=1$. By a direct calculation,
$$\fb=-\fa+\sla\triangle.$$
We now assume that the identity of Lemma \ref{baa_lem} is true at the level $k$ and prove the identity for the level $k+1$.
\begin{align*}
(2k+1)\fb\fa^k &= (2k-1)\fb\fa^k +2\fb\fa^k \\
&= (2k-1)\fa\fb\fa^{k-1} +(2k-1)[\fb,\fa]\fa^{k-1} +2\fb\fa^k\\
&= (2k-1)\fa\fb\fa^{k-1} +(2k-1)(2\fb^2+2\fa)\fa^{k-1} +2\fb\fa^k \\
&= (\fa+2\fb)((2k-1)\fb\fa^{k-1}) + 2\fb\fa^k + \fc^{\le k} \\
&= (\fa+2\fb)(-\fa^k+\fc^{\le k-1}\sla\triangle+\fc^{\le k-1}) +2\fb\fa^k + \fc^{\le k} \\
&= -\fa^{k+1} +\fc^{\le k}\sla\triangle + \fc^{\le k}.
\end{align*}
This completes the inductive step and thus proves Lemma \ref{baa_lem}.
\end{proof} 

Since Lemma \ref{ccc_lem} was reduced to Lemma \ref{baa_lem}, this concludes the proof of Lemma \ref{ccc_lem}.
\end{proof} 

Since Theorem \ref{cl_embedding_thm} was reduced to Lemma \ref{ccc_lem}, this concludes the proof of Theorem \ref{cl_embedding_thm}.
\end{proof} 

\section{Additional regularity lemmas}\label{additional_regularity_lemmas_sec}

The $L^\infty$ estimates in Chapter \ref{wm_kerr_chap} are based on the following lemma.
\begin{lemma}\label{spherical_infty_lem}
If $l$ is even, then
$$||\fd^l u||_{L^\infty(S^2)}\lesssim \sum_{i\le 2}||\fd^{l+i}u||_{L^2(S^2)}.$$
If $l$ is odd, then
$$||\fd^l u||_{L^\infty(S^2)}\lesssim \sum_{i\le 3}||\fd^{l+i}u||_{L^2(S^2)}.$$
\end{lemma}
\begin{proof}
We prove this theorem by applying the estimate
$$||u||_{L^\infty(S^2)}\lesssim ||u||_{L^2(S^2)}+||\sla\triangle u||_{L^2(S^2)}.$$
If $l$ is even, then
\begin{align*}
||\fd^lu||_{L^\infty(S^2)} &=||\fc^{l/2}u||_{L^\infty(S^2)} \\
&\lesssim ||\fc^{l/2}u||_{L^2(S^2)}+||\sla\triangle\fc^{l/2}u||_{L^2(S^2)} \\
&\lesssim ||\fc^{l/2}u||_{L^2(S^2)}+||\fc\fc^{l/2}u||_{L^2(S^2)} \\
&\lesssim ||\fd^lu||_{L^2(S^2)}+||\fd^{l+2}u||_{L^2(S^2)} \\
&\lesssim \sum_{i\le 2}||\fd^{l+i}u||_{L^2(S^2)}.
\end{align*}
If $l$ is odd, then
\begin{align*}
||\fd^lu||_{L^\infty(S^2)} &=||\pd_\theta \fc^{(l-1)/2}u||_{L^\infty(S^2)} \\
&\lesssim ||\cot\theta\pd_\theta\fc^{(l-1)/2}u||_{L^\infty(S^2)} + ||\sin\theta \pd_\theta\fc^{(l-1)/2}u||_{L^\infty(S^2)} \\
&\lesssim ||\fb\fc^{(l-1)/2}u||_{L^\infty(S^2)}+ ||\sin\theta \pd_\theta\fc^{(l-1)/2}u||_{L^2(S^2)} +||\sla\triangle(\sin\theta \pd_\theta\fc^{(l-1)/2}u)||_{L^2(S^2)} \\
&\lesssim ||\fd^{l+1}u||_{L^\infty(S^2)}+||\fd^lu||_{L^2(S^2)}+||\fc(\fd\cos\theta \fd^lu)||_{L^2(S^2)} \\
&\lesssim ||\fd^{l+1}u||_{L^2(S^2)}+||\fd^{l+3}u||_{L^2(S^2)}+||\fd^l u||_{L^2(S^2)}+\sum_{i\le 2}||\fd^{1+2-i}\cos\theta \fd^{l+i}u||_{L^2(S^2)} \\
&\lesssim \sum_{i\le 3}||\fd^{l+i}u||_{L^2(S^2)}.
\end{align*}
\end{proof}

Finally, we prove a lemma that handles the cases where commutators act on $\sin\theta$ factors that are necessary to transition from $\Sigma_t$ to $\tilde{\Sigma}_t$. (For example, see Lemma \ref{ws:example_Nphi_term_lem}.)
\begin{lemma}\label{gain_sin_lem}
$$||\fd^i(\sin^2\theta)\fd^{l-i}f||_{L^2(S^2)}\lesssim \sum_{j\le i}||\sin^2\theta \fd^{l-j}f||_{L^2(S^2)}$$
\end{lemma}
\begin{proof}
First, we prove that for a general function $f(\theta)$,
\begin{equation}\label{gain_sin_lem_initial_eqn}
\int_0^\pi f^2\sin\theta d\theta \lesssim \int_0^\pi (\pd_\theta f)^2\sin^3\theta d\theta +\int_0^\pi f^2\sin^3\theta d\theta.
\end{equation}
Indeed,
\begin{align*}
\int_0^\pi f^2\sin\theta d\theta &= \frac12\int_0^\pi f^2\pd_\theta(\sin^2\theta)d\theta +\int_0^\pi f^2(1-\cos\theta)\sin\theta d\theta \\
&= - \int_0^\pi f\pd_\theta f \sin^2\theta d\theta +\int_0^\pi f^2(1-\cos\theta)\sin\theta d\theta \\
&\lesssim \epsilon^{-1}\int_0^\pi (\pd_\theta f)^2\sin^3\theta d\theta +\epsilon\int_0^\pi f^2\sin\theta d\theta +\int_0^\pi f^2\sin^3\theta d\theta.
\end{align*}
Taking $\epsilon$ sufficiently small and subtracting the second term on the right side proves the estimate (\ref{gain_sin_lem_initial_eqn}).

Now, we repeatedly apply the estimate (\ref{gain_sin_lem_initial_eqn}) to prove the lemma for the case $i=2$.
\begin{align*}
\int_0^\pi(\fd^{l-2}f)^2\sin\theta d\theta &\lesssim \int_0^\pi(\pd_\theta \fd^{l-2}f)^2\sin^3\theta d\theta+\int_0^\pi (\fd^{l-2}f)^2\sin^3\theta d\theta \\
&\lesssim \int_0^\pi(\pd_\theta^2 \fd^{l-2}f)^2\sin^5\theta d\theta + \int_0^\pi(\pd_\theta \fd^{l-2}f)^2\sin^5\theta d\theta \\
&\hspace{2.6in}+\int_0^\pi (\fd^{l-2}f)^2\sin^5\theta d\theta \\
&\lesssim \int_0^\pi (\sin^2\theta\fd^l f)^2\sin\theta d\theta +\int_0^\pi (\sin^2\theta\fd^{l-1}f)^2\sin\theta d\theta  \\
&\hspace{2.3in}+ \int_0^\pi (\sin^2\theta\fd^{l-2}f)^2\sin\theta d\theta.
\end{align*}
In the first step, estimate (\ref{gain_sin_lem_initial_eqn}) was applied once, and in the second step, estimate (\ref{gain_sin_lem_initial_eqn}) was applied twice. The final step is simply a rearrangement of terms. It should be clear how this procedure generalizes to any larger $i>2$.
\end{proof}

\bibliographystyle{alpha}
\bibliography{bib}

\end{document}